\documentclass[a4paper,twoside,12pt]{memoir}
\pdfoutput=1

\title{Neighborly and almost neighborly configurations, and their duals}
\author{Arnau Padrol}

\usepackage{etoolbox} 
\newtoggle{bwprint} 
\togglefalse{bwprint}
\newtoggle{print} 
\togglefalse{print}

\usepackage{amsmath, amssymb, amsthm}
\usepackage{graphicx,color}
\usepackage[T1]{fontenc} 
\usepackage{bm}
\usepackage{bbm}
\usepackage{paralist}
\usepackage{dsfont}
\usepackage{mathtools}

\captionnamefont{\bfseries\small}
\captiontitlefont{\small}
\hangcaption

\usepackage[catalan,english]{babel}

\usepackage{nameref} 

\makeatletter
\renewcommand{\p@enumii}{\theenumi.}
\makeatother

\usepackage{import}             
\usepackage{pifont}		
\usepackage{kbordermatrix}      

\usepackage[dvipsnames]{xcolor} 

\usepackage{tikz} 
\usetikzlibrary{arrows,shapes}
\usepackage{array} 
\usepackage{multirow}

\usepackage[
colorlinks]{hyperref}

\usepackage{memhfixc} 
\usepackage
{hypcap} 
\hypersetup{
    bookmarksnumbered,
    pdfstartview={FitH},
    citecolor={darkblue},
    linkcolor={black},
    urlcolor={blue},
    pdfpagemode={UseOutlines}
}
\hypersetup{pdftitle={\thetitle}, pdfauthor={\theauthor}} 
\usepackage{url} 
\makeatletter 
\newcommand\org@hypertarget{}
\let\org@hypertarget\hypertarget
\renewcommand\hypertarget[2]{%
  \Hy@raisedlink{\org@hypertarget{#1}{}}#2%
} 
\makeatother 

\newsubfloat{figure} 
\subcaptionsize{\small}
\subcaptionlabelfont{\bfseries}
\subcaptionfont{\normalfont}

\usepackage{pdfpages}

\definecolor{red}{rgb}{1,0,0}
\definecolor{blue}{rgb}{0,0,1}
\iftoggle{bwprint}{%
\definecolor{darkblue}{gray}{0} 
}{%
\iftoggle{print}{%
\definecolor{darkblue}{gray}{0} 
}{%
\definecolor{darkblue}{rgb}{0,0,0.7} 
}
}

\definecolor{gris}{gray}{0.4} 
\definecolor{green}{rgb}{0,1,0}
\definecolor{gray25}{gray}{0.75}
\definecolor{gray50}{gray}{0.5}
\definecolor{gray75}{gray}{0.25}
\definecolor{black}{gray}{0}

\newcommand{\darkblue}{\color{darkblue}}

\newcommand{\black}{\color{black}}

\newcommand{\ligne}{{\black\hrule}} 

\makeatletter 

\newcommand{\partfigure}{}
\newcommand{\namepart}{}
\renewcommand*{\toclevel@part}{0} 
	\def\@partapp{Part}
	\renewcommand*{\partnumfont}{\black\fontseries{m}\fontshape{sl}\fontsize{15}{15}\selectfont\centering} 
	\renewcommand*{\printpartnum}{\partnumfont \@partapp \space \numtoName{\c@part}} 
\makechapterstyle{chapitre}{%
	\openright 
	\setlength{\beforechapskip}{0pt} 
	\setlength{\midchapskip}{18pt} 
	\setlength{\afterchapskip}{30pt} 
	\renewcommand*{\chapnumfont}{\black\fontseries{m}\fontshape{sl}\fontsize{15}{15}\selectfont\flushright} 
	\renewcommand*{\chaptitlefont}{\darkblue\fontseries{m}\fontshape{sc}\fontsize{25}{25}\selectfont\centering} 
	\renewcommand*{\printchapternum}{\chapnumfont \@chapapp \space \numtoName{\c@chapter}} 

	\renewcommand*{\printchaptertitle}[1]{\chaptitlefont \ligne \begin{center} ##1 \end{center} \ligne%
 										\ifnum\value{compteurlof}>\z@ \addcontentsline{lof}{chapter}{\ifnum\c@chapter>\z@ \thechapter~~~ \fi ##1} \fi%
										\setcounter{figure}{0}%
										\thispagestyle{baspage}} 
}
\chapterstyle{chapitre}

\maxsecnumdepth{subsection} 
\setsecheadstyle{\bfseries\boldmath\fontshape{sc}\fontsize{18}{18}\selectfont \noindent\textcolor{darkblue}} 
\setsubsecheadstyle{\bfseries\boldmath\fontshape{n}\fontsize{15}{15}\selectfont \noindent \textcolor{darkblue}} 
\setsubsubsecheadstyle{\bfseries\boldmath\fontshape{n}\fontsize{13}{13}\selectfont \noindent \textcolor{darkblue}} 
\setparaheadstyle{\bfseries\boldmath\fontshape{n}\fontsize{12}{12}\selectfont \noindent \textcolor{darkblue}} 

\settocdepth{subsection} 
\renewcommand{\@pnumwidth}{3em} 

\makepagestyle{hautpage} 
	\makeevenhead{hautpage}{\textbf \thepage}{\leftmark}{} 
	\makeoddhead{hautpage}{}{\rightmark}{\textbf \thepage} 
	\makeheadrule{hautpage}{\textwidth}{.5pt} 
\makepagestyle{hautpageintro} 
	\makeevenhead{hautpageintro}{\textbf \thepage}{\namepart}{} 
	\makeoddhead{hautpageintro}{}{\namepart}{\textbf \thepage} 
	\makeheadrule{hautpageintro}{\textwidth}{.5pt} 
\makepagestyle{baspage} 
	\makeoddfoot{baspage}{}{\thetitle\quad}{\textbf \thepage} 
	\makefootrule{baspage}{\textwidth}{.4pt}{1pt} 
\renewcommand*{\chaptermark}[1]{\markboth{\@chapapp\ \thechapter\ -\ #1}{\@chapapp\ \thechapter\ -\ #1}} 
\makeindex

\makeatother

\newcommand{\defn}[1]{\emph{\darkblue #1}} 
\newcommand{\imp}[1]{\textbf{#1}} 

\newcommand{\ie}{\textit{i.e.},~} 
\newcommand{\etc}{\textit{etc.}} 
\newcommand{\eg}{\textit{e.g.},~} 

\def\cC{\mathcal{C}}
\def\cD{\mathcal{D}}
\def\cE{\mathcal{E}}
\def\cF{\mathcal{F}}
\def\cG{\mathcal{G}}
\def\cH{\mathcal{H}}

\def\cM{\mathcal{M}}
\def\cN{\mathcal{N}}
\def\cO{\mathcal{O}}
\def\cP{\mathcal{P}}

\def\cS{\mathcal{S}}
\def\cT{\mathcal{T}}

\def\cV{\mathcal{V}}

\def\NN{\mathbbm{N}}

\def\RR{\mathbbm{R}}

\def\ZZ{\mathbbm{Z}}
\def\PP{\mathbbm{P}}

\newcommand{\gs}{\sigma}
\newcommand{\gl}{\lambda}
\newcommand{\gn}{\nu}
\newcommand{\gm}{\mu}
\newcommand{\gk}{\gk}
\newcommand{\ep}{\varepsilon}

\def\e{\mathrm{e}}
\newcommand{\tick}{\ding{51}}
\newcommand{\cross}{\ding{55}}
\newcommand{\matzero}{\mathbf{0}}
\newcommand{\matid}{\mathbf{I}}
\newcommand{\veczero}{\boldsymbol{0}}
\newcommand{\trans}[1]{{#1}^\top}
\newcommand{\ul}[1]{\underline{#1}}
\newcommand{\ol}[1]{\overline{#1}}
\newcommand{\vv}[1]{\boldsymbol{\mathbf{#1}}}
\newcommand{\vvc}[1]{\mathrm{#1}}
\newcommand{\vvh}[1]{\boldsymbol{\mathbf{#1}}}

\newcommand{\conv}{\mathrm{conv}}

\newcommand{\verts}{\mathrm{vert}}
\newcommand{\lin}{\mathrm{lin}}
\newcommand{\aff}{\mathrm{aff}}
\newcommand{\cone}{\mathrm{cone}}
\newcommand{\relint}{\mathrm{relint}}
\newcommand{\intr}{\mathrm{int}}
\DeclareMathOperator{\rank}{rank}

\DeclareMathOperator{\disc}{disc}
\newcommand{\sprod}[2]{\langle {#1} , {#2} \rangle}
\newcommand{\set}[2]{\ensuremath{\left\{#1\,\middle|\,#2\right\}}} 

\newcommand{\ffloor}[2]{\left\lfloor{\frac{#1}{#2}}\right\rfloor}
\newcommand{\fceil}[2]{\left\lceil {\frac{#1}{#2}} \right\rceil}
\newcommand{\bracket}[1]{[#1]}

\newcommand{\rst}[1]{\ensuremath{{\mathbin\upharpoonright}\raise-.5ex\hbox{$#1$}}} 
\newcommand{\card}[1]{\vert {#1} \vert}
\def\restriction#1#2{\mathchoice
              {\setbox1\hbox{${\displaystyle #1}_{\scriptstyle #2}$}
              \restrictionaux{#1}{#2}}
              {\setbox1\hbox{${\textstyle #1}_{\scriptstyle #2}$}
              \restrictionaux{#1}{#2}}
              {\setbox1\hbox{${\scriptstyle #1}_{\scriptscriptstyle #2}$}
              \restrictionaux{#1}{#2}}
              {\setbox1\hbox{${\scriptscriptstyle #1}_{\scriptscriptstyle #2}$}
              \restrictionaux{#1}{#2}}}
\def\restrictionaux#1#2{{#1\,\smash{\vrule height .8\ht1 depth .85\dp1}}_{\,#2}} 
\newcommand{\simp}[1]{\vv{\triangle}_{#1}}
\newcommand{\cube}[1]{\vv{\Box}^{#1}}
\newcommand{\cros}[1]{\vv{\Diamond}^{#1}}

\DeclareMathOperator{\bigjoin}{\text{\huge{$\ast$}}}
\DeclareMathOperator{\join}{\ast}

\newcommand{\cyc}[2]{\vv C_{#1}({#2})} 
\newcommand{\stc}[1]{\cD_{#1}}
\newcommand{\lnei}[2]{\operatorname{nb}_l({#1,#2})}
\newcommand{\lpol}[2]{\operatorname{p}_l({#1,#2})}
\newcommand{\nnei}[2]{\operatorname{nb}({#1,#2})}
\newcommand{\lnr}[2]{\operatorname{nr}_l({#1,#2})}
\newcommand{\lle}[2]{\ell_l({#1,#2})}

\newcommand\fF{\cT} 
\DeclareMathOperator{\sew}{Sew}
\def\rd{s}
\def\rr{r}

\def\l{t}
\newcommand{\Gale}[1]{{#1}^{\star}} 
\newcommand{\pGale}[1]{\left( {#1}\right)^{\star}} 
\newcommand{\bGale}[1]{\bar{#1}^{\star}} 
\DeclareMathOperator{\Val}{Val}
\DeclareMathOperator{\Dep}{Dep}
\DeclareMathOperator{\aVal}{a-Val}
\DeclareMathOperator{\aDep}{a-Dep}

\def\ci{\cC} 
\def\co{\Gale\ci} 
\def\ve{\cV} 
\def\cov{\Gale\ve} 

\DeclareMathOperator{\Cayley}{Cayley}

\DeclareMathOperator{\codeg}{codeg}
\let\deg\relax 
\DeclareMathOperator{\deg}{deg}
\def\degc{\deg}
\def\codegc{\codeg}
\def\degG{\Gale{\deg}}
\DeclareMathOperator{\codegG}{\Gale\codeg}
\def\degZ{\deg_{\ZZ}}
\def\codegZ{\codeg_{\ZZ}}

\def\CayleyG{{Cayley$\Gale{}$} }
\def\G{{$\Gale{}$} }
\def\codegreeG{{codegree$\Gale{}$} }

\def\dd{{\delta}}
\def\kk{{\kappa}}

\DeclareMathOperator{\DD}{\Delta}

\def\suc{\mathrm{suc}}
\def\IG{\mathrm{IG}}
\def\Ce{\cC}
\def\Tv{\cD}

\newtheorem{theorem}{Theorem}[chapter] 
\newtheorem{proposition}[theorem]{Proposition} 
\newtheorem{lemma}[theorem]{Lemma} 
\newtheorem{corollary}[theorem]{Corollary}
\newtheorem{conjecture}[theorem]{Conjecture}
\newtheorem{question}[theorem]{Question}

\newtheorem{assumption}{Assumption}

\theoremstyle{remark}
\newtheorem{remark}[theorem]{Remark}
\newtheorem{observation}[theorem]{Observation}
\newtheorem{example}[theorem]{Example}

\theoremstyle{definition}
\newtheorem{definition}[theorem]{Definition} 
\newtheorem{definitions}[theorem]{Definitions}

\newtheoremstyle{constrstyle}
  {\topsep}   
  {\topsep}   
  {\rmfamily}  
  {0pt}       
  {\bfseries} 
  {.}         
  {5pt plus 1pt minus 1pt} 
  {}          

\theoremstyle{constrstyle}
\newtheorem{constr}{Construction}

\definecolor{shadecolor}{gray}{.9}

\begin{document}
\selectlanguage{english}

\frontmatter

 \newcounter{compteurlof}
 \setcounter{compteurlof}{0}
%

\aliaspagestyle{plain}{empty}

\begin{titlingpage}

\vspace*{\stretch{2}}
{\black\hrule}
{
\scshape
\begin{center}
\Huge\textbf{
\thetitle
}
\end{center}
}

\vspace{\stretch{1}}

\begin{center}
{\bfseries \large Arnau Padrol Sureda}
\end{center}
{\black\hrule}

\vspace{\stretch{0.5}}

\vspace{\stretch{3.5}}
\vspace*{\stretch{1}}
\end{titlingpage}

\thispagestyle{empty}
\begin{center}
{Doctoral program} in {Applied Mathematics}
\end{center}

\vspace*{\stretch{2}}
{
\scshape

\begin{center}
\darkblue\Huge \textbf{
\thetitle
}
\end{center}
}
\vspace{\stretch{.2}}

\begin{center}
Thesis submitted by 
 \\[3ex]
{\bfseries \Large Arnau Padrol Sureda}
 \\[3ex]
for the degree of Doctor of Mathematics in the \\[2ex]
{\bfseries\scshape Universitat Polit\`ecnica de Catalunya} 
\end{center}
\vspace{\stretch{0.2}}

\begin{center}
Thesis Advisor\\[2ex]
{\bfseries \large Julian Pfeifle}
\end{center}

\vspace{\stretch{3}}
\begin{center}
{Barcelona, March 2013}
\end{center}
\begin{center}
{ Departament de Matem\`atica Aplicada II}
\\[1ex] {\scshape Universitat Polit\`ecnica de Catalunya} 
\end{center}
\newpage

\thispagestyle{empty}

{
\vspace*{\stretch{1}}
{\sffamily\noindent
\textbf{Arnau Padrol Sureda}\\
Departament de Matem\`atica Aplicada II\\
Universitat Polit\`ecnica de Catalunya\\
Edifici Omega, Jordi Girona 1-3 \\
08034 Barcelona\\
$<$arnau.padrol@upc.edu$>$
}
}
\newpage

\chapter*{Acknowledgements}

\selectlanguage{english}

First and foremost, I want to thank my advisor, \defn{Julian Pfeifle}, for making this work possible.
You have taught me how to do research, how to write a paper and how to give a talk. I am very grateful for that; and even more for the great time we had together.
During these years I have always felt that the doors of your office, your Inbox and every cafeteria in Gr\`acia were open for me whenever I had yet another question. 
Thank you, Julian, for your patience, advice and encouragement.

My thanks also go to \defn{G\"unter Ziegler}, \defn{Francisco Santos} and \defn{Eran Nevo} --- as well as to \defn{Oriol Serra} and \defn{Vincent Pilaud} --- for having accepted to become part of my thesis committee. 
I have had many interesting mathematical discussions with each of you, without which I am sure that this thesis would not be the same. It is not only an honour that you managed to find a spot for my defense in your busy agendas, but I am personally very happy that I can share this moment with you.

I am also 
indebted to \defn{Benjamin Nill} for coauthoring the paper that is the main source for Part II. Collaborating with you has been a wonderful experience with memorable stages in Leuven, Frankfurt and Kyoto.

I wish to thank \defn{Uli Wagner} for his hospitality during my visit to ETH Z\"urich in 2011. 
Our inspiring conversations on neighborly polytopes motivated the research that has become Part I of this thesis.

Next, I want to thank \defn{Christian Haase} for inviting me to Goethe Universit\"at at Frankfurt am Main. Thanks for sharing your knowledge on lattice polytopes with me, and for disproving my groundless conjecture that a human being can only drink a bounded amount of tea.\\

\selectlanguage{catalan}
Vull agrair tamb\'e als meus companys del \defn{MA II} per fer que al departament m'hi trob\'es com a casa durant aquests anys. En particular, m'agradaria mencionar al \defn{Marc Noy} i al \defn{Ferran Hurtado} per l'esfor\c{c} que dediquen a gestionar els grups de recerca.
I, \'es clar, tampoc puc obviar els dinars. El moment de desconnexi\'o diari durant la tert\'ulia dels migdies ha estat imprescindible per no perdre el cap. Gr\`acies \defn{Maria}, \defn{Rodrigo}, \defn{In\^es}, \defn{Mat\'ias}, \defn{Carlos} i \defn{Elisa}. \selectlanguage{english}
I also want to thank \defn{Aaron}. Every PhD student should have somebody like you in some office nearby. I really enjoyed our (mathematical and non-mathematical) conversations. 

\selectlanguage{catalan}
Sense els meus amics, aix\`o no hagu\'es estat possible. Ja des d'abans de comen\c car, el \defn{Marcel} i el \defn{Juanjo} han compaginat el seu rol d'amics amb el de consellers d'assumptes doctorals, i puc assegurar que les dues tasques les fan molt b\'e. Ells van ser els primers d'una fornada de doctorands que m'han fet sentir sempre acompanyat. Em venen al cap especialment el \defn{V\'ictor}, tant per tots els moments que hem compartit omplint paperassa per demanar beques com per quan el que deman\`avem eren cerveses; i el \defn{Guillem}, perqu\`e no em puc imaginar ning\'u millor amb qui compartir una confer\`encia d'estiu. Tampoc podr\'e oblidar les escapades a la BdM amb l'\defn{Inma}, els tuppers a la gespa del Campus Nord amb l'\defn{Elena} i la \defn{Cris}; ni les (sempre massa espor\`adiques) visites de l'\defn{Ari}, la \defn{Marga} i el \defn{Xavi}.
Tamb\'e vull donar les gr\`acies a l'\defn{Albert}, al \defn{Pucho} i al \defn{Pascal}, per acollir-me a casa seva i fer-me d'amfitrions en una ciutat nova; 
 a la \defn{Maria} i al \defn{Steffen}, pels sopars extraordinaris a Can Mantega; al \defn{D\'idac}, perqu\`e encara ric quan recordo l'\`atic de Tres Senyores; i a molts d'altres com l'\defn{Itziar}, l'\defn{Emilia}, el \defn{Ra\"ul} i un bon grapat que de ben segur estic oblidant de manera imperdonable.

Agraeixo als meus pares, el \defn{Josep Maria} i la \defn{Catalina}, tot el suport que m'han donat durant aquests anys.
Us vull donar les gr\`acies per ser sempre a prop, per tenir sempre a punt el consell adequat, i perqu\`e sempre que 
no n'he fet cas (caparrut com s\'oc, massa sovint), s'ha acabat demostrant que ten\'ieu ra\'o.
Tamb\'e vull donar les gr\`acies a la \defn{Padrina} i a l'\defn{\`Avia}, perqu\`e encara que no acabin d'entendre de qu\`e va la tesi, em diran que no n'havien vista mai cap de tan ben feta, i encara menys cap d'escrita per alg\'u tan ben plantat; al \defn{Jeroni}, que va ser el primer doctor de la fam\'ilia, i a l'\defn{Esperan\c ca}, que ser\`a la seg\"uent. 
I acabo el par\`agraf amb l'\defn{Agn\`es}, de qui no puc deixar d'admirar (i d'envejar, que \'es el que fan els bons germans) la seva energia contagiosa.
\\ 

Les darreres l\'inies s\'on per a la \defn{Laura}, que 
ha estat al meu costat tot aquest temps. Moltes gr\`acies per tot, i m\'es. Al final sembla que, a vegades, ens en sortim.

\selectlanguage{english}



 \pagestyle{hautpageintro}

\selectlanguage{english}

\renewcommand{\namepart}{Contents}
\cleardoublepage  \tableofcontents
\newpage
 \renewcommand{\namepart}{List of figures}
\cleardoublepage  \listoffigures

\vspace{4cm}

 \renewcommand{\namepart}{List of tables}
\listoftables

\pagestyle{hautpage}
\mainmatter
\setcounter{compteurlof}{1}

\renewcommand{\namepart}{Introduction}

\setcounter{chapter}{-1}
\chapter{Summary}
\setlength{\epigraphwidth}{.65\textwidth}
\epigraphtextposition{flushleftright}
\setlength{\epigraphrule}{0pt}
\epigraph{\emph{It seems that our main limits in understanding the combinatorial structure of polytopes still lie in our ability to raise the right questions. Another feature that comes to mind (and is not unique to this area) is the lack of examples, methods of constructing them, and means of classifying them.}}
{Gil Kalai~\cite{Kalai1997}}

This thesis provides new applications of Gale duality to the study of polytopes with extremal combinatorial properties.
It is divided into two parts. The first one explores different methods for constructing \defn{neighborly} polytopes and oriented matroids. The second part is devoted to the study of the \defn{degree} of point configurations, which is a combinatorial invariant closely related to neighborliness. 

All relevant definitions are properly given in Chapter~\ref{ch:intro}; for now we just need to know that a \defn{(convex) polytope} is the convex hull of a finite set of points in the Euclidean space $\RR^d$ and that a \defn{face} of a polytope is its intersection with a supporting hyperplane. For example, Figure~\ref{fig:associahedron} shows an associahedron: a $(d=3)$-dimensional polytope with $14$ vertices ($0$-dimensional faces) and $9$ facets ($(d-1)$-dimensional faces). From a combinatorial point of view, one wants to answer questions about the structure of inclusions among faces of a polytope, its \defn{face lattice}.

\iftoggle{bwprint}{%
\begin{figure}[htpb]
\begin{center}
\includegraphics[width=.65\linewidth]{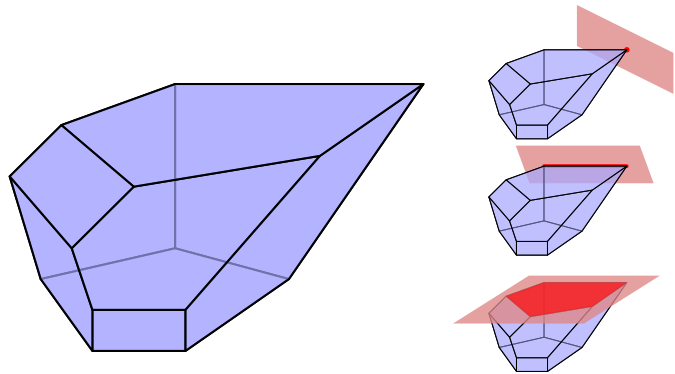}
\end{center}
\caption[An associahedron.]{An associahedron and three of its faces, each with a supporting hyperplane.}
\label{fig:associahedron}
\end{figure} 
}{%
\begin{figure}[htpb]
\begin{center}
\includegraphics[width=.65\linewidth]{Figures/associahedron_col}
\end{center}
\caption[An associahedron.]{An associahedron and three of its faces, each with a supporting hyperplane.}
\label{fig:associahedron}
\end{figure} 
}

A first question could be which $d$-dimensional polytopes with $n$ vertices have the \emph{maximal number of facets}.  
To this, the answer has been known since 1970, when McMullen proved the \defn{Upper Bound Theorem}~\cite{McMullen1970}. It states that the number of facets is maximized by \defn{neighborly polytopes}. These are characterized by the property that every subset of $\ffloor{d}{2}$ vertices forms a face; and are the main characters of \imp{Part~\ref{partI}}. 

Here we encounter for the first time an obstacle that crops up time and again throughout this text: Showing pictures of $d$-dimensional polytopes is difficult for~$d\geq 4$. This means that we will usually have to settle for examples in dimensions $2$ or $3$ to try to get intuition about phenomena happening in arbitrary dimension. However, many interesting properties of polytopes do not appear until higher dimensions. In particular, all $d$-polytopes with $d\leq 3$ are neighborly, and hence there is no non-trivial example until $d=4$.

\begin{figure}[htpb]
\begin{center}
  \includegraphics[width=.3\linewidth]{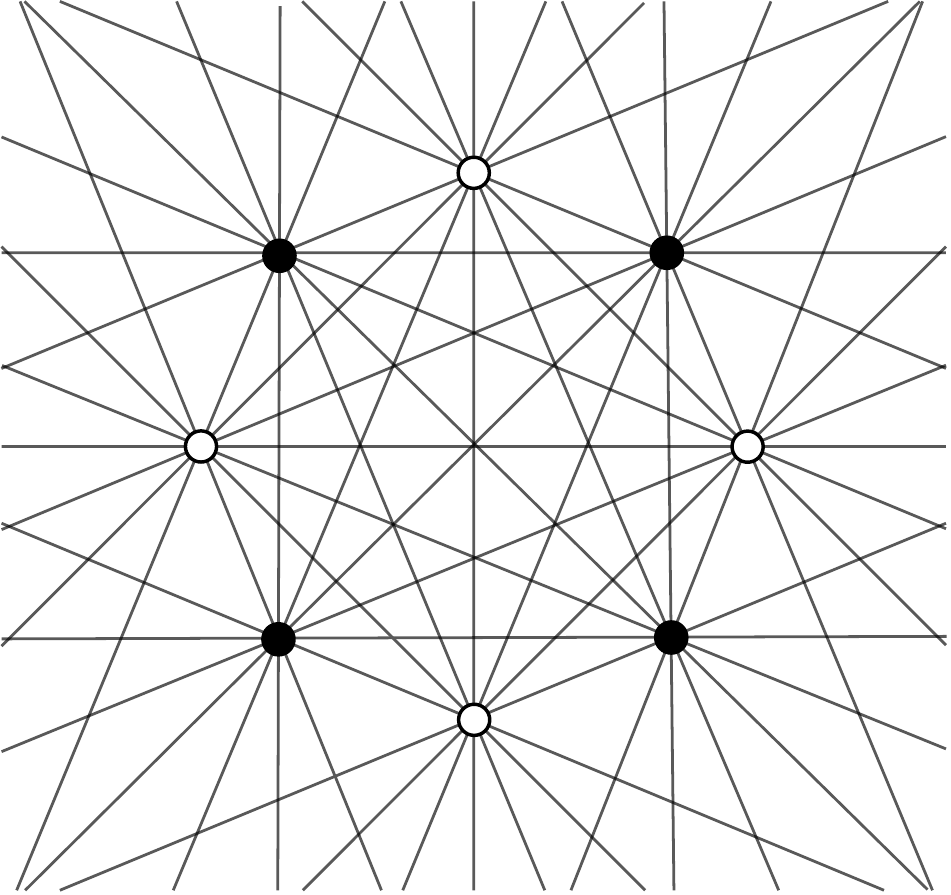}\quad
  \includegraphics[width=.3\linewidth]{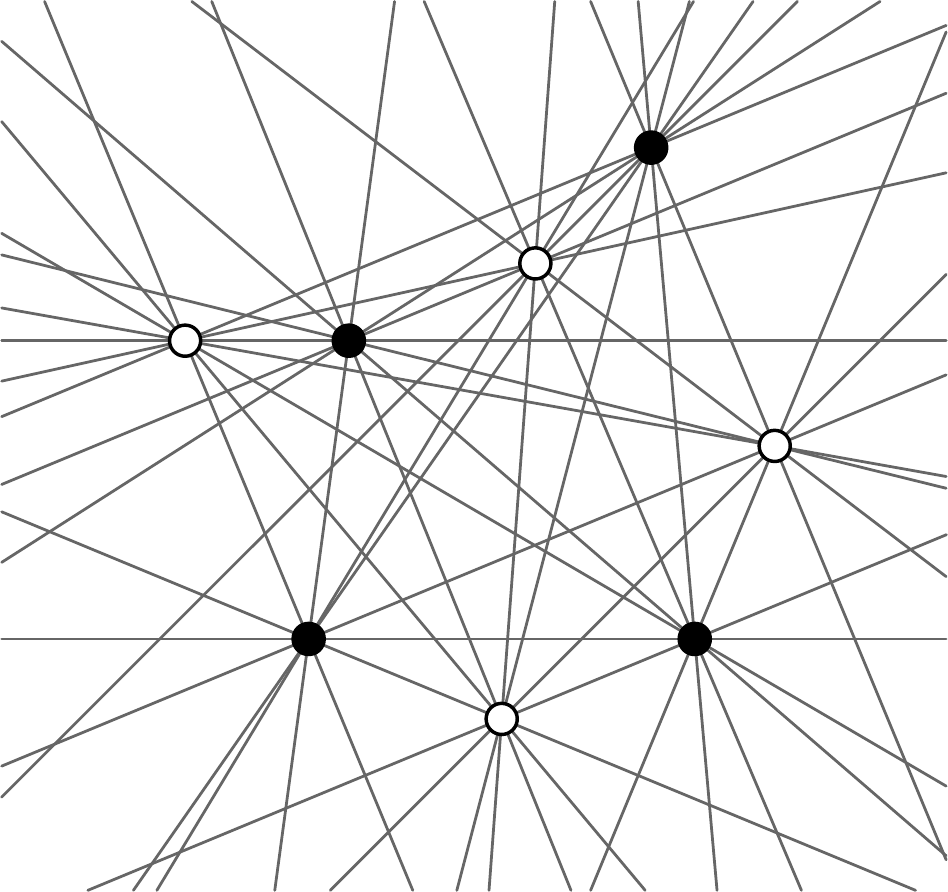}\quad
  \includegraphics[width=.3\linewidth]{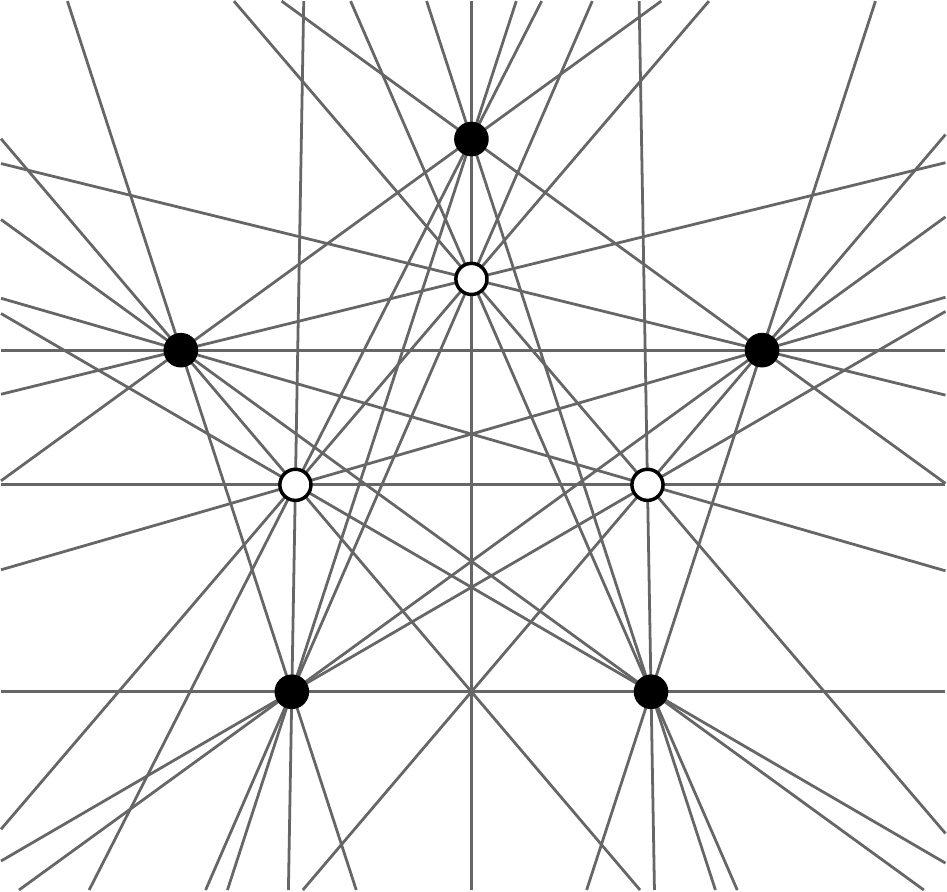}
\end{center}
\caption[Affine Gale diagrams.]{Affine Gale diagrams of the three combinatorial types of neighborly $4$-polytopes with $8$ vertices.}
\label{fig:N48}
\end{figure}

One trick to overcome this obstacle is the use of Gale duality. It allows to represent a configuration of $n$ points in $\RR^d$ by a configuration of~$n$ vectors in $\RR^{n-d-1}$ (or even a colored $n$-point configuration in $\RR^{n-d-2}$). For example, in Figure~\ref{fig:N48}
there are three configurations of $8$ colored points in $\RR^2$, which represent three neighborly $4$-polytopes with $8$~vertices. 
Traditionally, Gale duality has been a very useful tool for obtaining results on polytopes with few vertices, such as $d$-dimensional polytopes with at most~$d+4$ vertices. However, we use it to analyze polytopes with an arbitrary number of vertices, which is a nonstandard application.\\

Back to the questions, it is quite reasonable to ask \emph{how many polytopes there are}. More specifically, if the dimension and the number of vertices are fixed, how many different polytopal face lattices are there? For some time, the best lower bound for this quantity was one found by Shemer, who in the essential paper~\cite{Shemer1982} from 1982 described the \defn{Sewing construction} for neighborly polytopes. 

In Chapter~\ref{ch:shemer} we analyze his construction from the point of view of \defn{oriented matroids}. It turns out that sewing can be easily interpreted using the concept of \defn{lexicographic extensions} of oriented matroids. This provides a simpler proof that sewing works, and can be used to explain related results by Barnette~\cite{Barnette1981}, Bistriczky~\cite{Bisztriczky2000}, Lee \& Menzel~\cite{LeeMenzel2010} and Trelford \& Vigh \cite{TrelfordVigh2011}.
The main result of this chapter is \imp{Theorem~\ref{thm:extshemersewing}}. It is the main ingredient for our \defn{Extended Sewing} technique, which extends Shemer's sewing construction to oriented matroids of any rank, and works with a larger family of flags of faces. Extended Sewing is optimal in the sense that in odd ranks, the flags of faces constructed in this way are the only ones that yield neighborly polytopes (Proposition~\ref{prop:uniqueflags}).

Our second (and most important) new technique is \defn{Gale Sewing}, which is presented in Chapter~\ref{ch:thethm}.
It lexicographically extends {duals} of neighborly polytopes and oriented matroids with the \imp{Double Extension Theorem~\ref{thm:thethm}}. 
Under certain conditions made explicit in Theorem~\ref{thm:cOsubsetcG}, it generalizes the sewing construction. However, both Extended Sewing and Gale Sewing are needed to find many non-realizable neighborly oriented matroids in \imp{Theorem~\ref{thm:nonrealizable}}.

In Chapter~\ref{ch:counting} we estimate the number of polytopes in $\cG$, the family of all polytopes that can be constructed via Gale Sewing. This provides one of the main results of this thesis, which is the following new lower bound for $\lnei{r+d+1}{d}$, the number of combinatorial types of vertex-labeled neighborly $d$-polytopes with $r+d+1$ vertices:
\begin{equation}
\lnei{r+d+1}{d}\geq \frac{\left( r+d \right) ^{\left( \frac{r}{2}+\frac{d}{2} \right) ^{2}}}{{r}^{{(\frac{r}{2})}^{2}}{d}^{{(\frac{d}{2})}^{2}}{{\e}^{3\frac{r}{2}\frac{d}{2}}}}. \end{equation}
This bound is proved in \imp{Theorem~\ref{thm:lblnei}}. Not only it is greater than Shemer's bound, but it even improves the current best lower bound on the number of all polytopes~\cite{Alon1986}. In \imp{Theorem~\ref{thm:nonrealizablebound}} we show that proportional lower bounds also hold for the number of labeled non-realizable neighborly oriented matroids.
\\

The \defn{degree} of a $d$-dimensional point configuration is defined as the maximal codimension of its interior faces; \ie the minimal number $\dd$ such that every subset of $(d-\dd)$ points lies in a common facet. Hence, neighborly $d$-polytopes are precisely those $\ffloor{d}{2}$-simplicial polytopes whose vertex set has degree $\fceil{d}{2}$. In \imp{Part~\ref{partII}} we use a Gale dual interpretation of the degree to prove several results that concern the combinatorial structure of point configurations whose degree is small compared to their dimension; specifically, those whose degree is smaller than $\fceil{d}{2}$, the degree of neighborly polytopes.

The study of the degree comes motivated by a corresponding Ehrhart-theoretic notion for lattice polytopes. In this area, the lattice degree is a measure of complexity of Ehrhart polynomials, and as such, it has been studied quite intensively. We provide new combinatorial interpretations of several recent results on lattice polytopes of small lattice degree.
The analogies between the combinatorial degree and its Ehrhart-theoretic counterpart are explained in Chapter~\ref{ch:intro_degree}, where we also present links to other subjects in geometric combinatorics, such as neighborly polytopes (a relation that goes back at least to Motzkin in 1965~\cite{Motzkin1965}), the Generalized Lower Bound Theorem, and Tverberg theory. In this last setting, our results about small degrees can be translated into statements about point configurations with a non-empty $\kk$-core for some large value of $\kk$.

A first result in this direction is \imp{Corollary~\ref{cor:easybound}}, which states that any $d$-dimensional configuration of $d+1+r$ points of degree $\dd$ with $d \geq r+2 \dd$ must be a pyramid. This should be compared to a theorem in~\cite{Nil08}, according to which any $d$-dimensional lattice polytope with $r + d + 1$ vertices and lattice degree $s$ such that $d > r(2s + 1) + 4s - 2$ must be a lattice pyramid.

In Chapter~\ref{ch:cayley}, we present \defn{weak Cayley configurations} as a tool for understanding the structure of point configurations of small degree. It is shown in~\cite{HNP09} that if the dimension of a lattice polytope~$\vv P$ is greater than a quadratic function on the lattice degree, then $\vv P$~can be projected onto a unimodular simplex, which means that $\vv P$~is a \defn{lattice Cayley polytope}. \imp{Theorem~\ref{thm:d+1-3dd}} is our combinatorial analogue. It states that if the dimension~$d$ of a point configuration~$\vv A$ of degree~$\dd$ is greater than $3\dd$, then it must have a contraction $\vv A/\vv S$ that admits a projection onto the vertex set of a $(d-3\dd)$-simplex, i.e., $\vv A$ is a weak Cayley configuration. We conjecture that this can be further strengthened to $\dd<\frac{d}{2}$, which would be analogous to the classical theorem according to which any polytope that is more than $\ffloor{d}{2}$-neighborly must be a simplex.

The dual version of Theorem~\ref{thm:d+1-3dd} shows that every point in the $\kk$-core of a $r$-dimensional configuration of $n$ points is $({3\kk-2(n-r)})$-divisible, and the conjecture states that it is in fact $({2\kk-(n-r)})$-divisible.

These results can be strengthened for configurations of degree~$1$, which we fully classify in \imp{Theorem~\ref{thm:dd=1}}, the main result of Chapter~\ref{ch:deg1}. These are strongly related to \defn{totally splittable} polytopes, for which we give a combinatorial explanation of their equidecomposability, thus answering a question by Herrmann and Joswig in~\cite{HerrmannJoswig2010}.

The last chapter of the thesis is devoted to a stronger structural conjecture for configurations whose degree is smaller than half the dimension. It is formulated in terms of \defn{codegree decompositions}, a stronger concept than weak Cayley configurations.
If true, \imp{Conjecture~\ref{conj:strongd+1-2dd}} would imply all the previous results on the combinatorial degree. For some special cases, we can prove it. For example, \imp{Corollary~\ref{cor:Lawrence}}
 is a new characterization of Lawrence polytopes in terms of their \defn{vector discrepancy}, the maximal difference between the sizes of a Radon partition of their vertices.
Finally, \imp{Theorem~\ref{thm:2DD}} implies that every configuration of $r+d+1$ points in~$\RR^d$ of degree~$\delta$ smaller than $\frac{2(d+1)-r}{4}$ admits a non-trivial codegree decomposition.

Dualizing, this conjecture states that if $\vv x$ is in the $\kk$-core of a configuration of $n$ points~$\vv A$ in $\RR^r$, $\vv x\in \Ce_\kk(\vv A)$, then there are $m$ disjoint subsets $\vv S_1,\dots,\vv S_m$ of $\vv A$, with $m\geq 2\kk-(n-r)$, such that $\vv x\in \Ce_{\kk_i}(\vv S_i)$ and $\sum_{i=1}^m \kk_i=\kk$. A corollary of \imp{Theorem~\ref{thm:DD=2}} shows that this conjecture holds for $r\leq4$.

\vspace{2cm}

All the results cited in this dissertation are properly attributed to their authors with a reference to the bibliography or stated as ``well known'' in their preceding text. Except for the introductory concepts in Chapter~\ref{ch:intro}, if no authorship is indicated, the result should be understood to be new. Chapters~\ref{ch:intro_degree}, \ref{ch:cayley} and \ref{ch:deg1} are partly joint work with Benjamin Nill.

\chapter{Background and notation}\label{ch:intro}
The goal of this chapter is to fix notation and present some basic concepts. A nice introduction for most of what follows can be found in Chapters 5 and 6 of Matou\v{s}ek's book~\cite{Matousek2002}. For convex polytopes, one of our main references is Gr\"unbaum's classical book~\cite{GruenbaumEtal2003}. The other one is Ziegler's book~\cite{Ziegler1995}, where oriented matroids are also introduced. The monograph~\cite{OrientedMatroids1993} by Bj\"orner et al.\ provides a more comprehensive treatment for oriented matroids. A final recommendation is the book~\cite{DeLoeraRambauSantosBOOK} by De~Loera, Rambau and Santos, which focuses on triangulations.

\section{Polytopes and point configurations}
A \defn{point configuration}\index{point configuration} $\vv A = \{\vv a_1 ,\dots ,\vv a_n\}$ is a finite collection of (labeled) points in the affine space $\RR^d$. We use the word ``configuration'' instead of ``set'' because we do not require the points in $\vv A$ to be different, as long as they have different labels. That is, even if $\vv a_i$ and $\vv a_j$ share the same coordinates, if $i\neq j$ we still consider $\vv a_i$ and $\vv a_j$ to be different elements of~$\vv A$. 
The \defn{affine span}\index{affine!span} of~$\vv A$ is \index{$\aff(\vv A)$}\[\aff(\vv A):=\set{\vv x\in \RR^d}{\vv x=\sum_{i=1}^n\gl_i \vv a_i; \;\;\text{ for }\vv a_i\in \vv A,\; \gl_i\in \RR\text{ and }\sum_{i=1}^n\gl_i=1},\]
and its \defn{dimension}\index{point configuration!dimension} $\dim(\vv A)$ is the dimension of $\aff(\vv A)$. We usually consider  $\vv A$ to be full dimensional, that is, $\dim(\vv A)=d$\index{$\dim(\vv A)$}.

\iftoggle{bwprint}{%
\begin{figure}[htpb]
\centering
 \subbottom[]{\label{sfig:4pointsR2_1}\includegraphics[width=.2\linewidth]{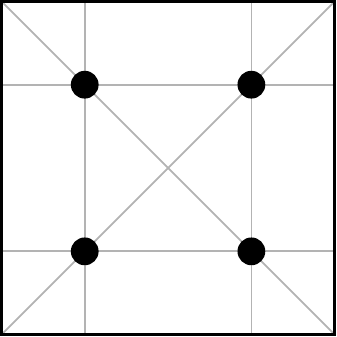}}\qquad
 \subbottom[]{\label{sfig:4pointsR2_2}\includegraphics[width=.2\linewidth]{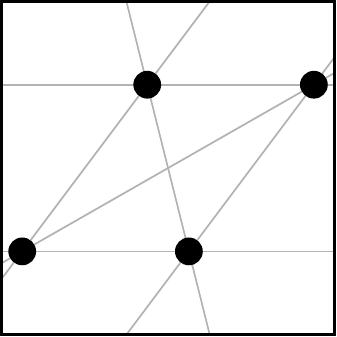}}\qquad
 \subbottom[]{\label{sfig:4pointsR2_3}\includegraphics[width=.2\linewidth]{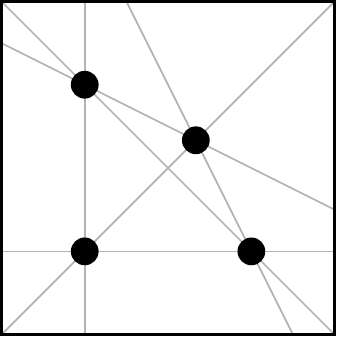}}\qquad
 \subbottom[]{\label{sfig:4pointsR2_4}\includegraphics[width=.2\linewidth]{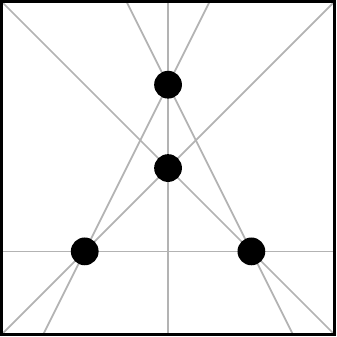}}
\caption{Four point configurations consisting of $4$ points in $\RR^2$}
\label{fig:4pointsR2}
\end{figure}
}{%
\iftoggle{print}{%
\begin{figure}[htpb]
\centering
 \subbottom[]{\label{sfig:4pointsR2_1}\includegraphics[width=.2\linewidth]{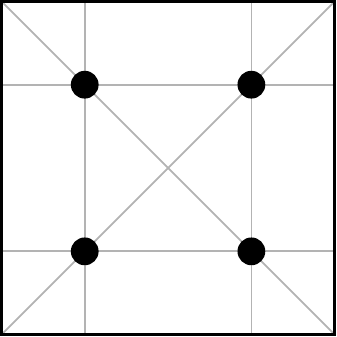}}\qquad
 \subbottom[]{\label{sfig:4pointsR2_2}\includegraphics[width=.2\linewidth]{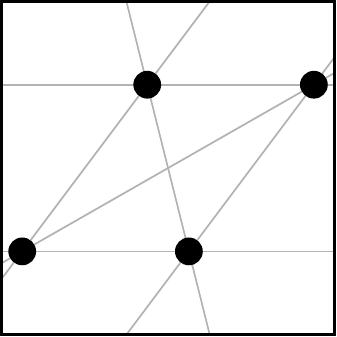}}\qquad
 \subbottom[]{\label{sfig:4pointsR2_3}\includegraphics[width=.2\linewidth]{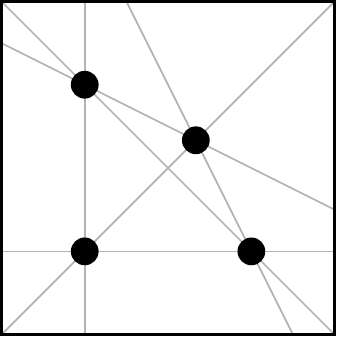}}\qquad
 \subbottom[]{\label{sfig:4pointsR2_4}\includegraphics[width=.2\linewidth]{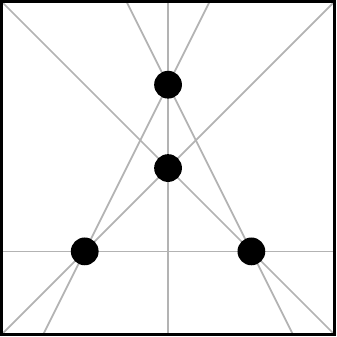}}
\caption{Four point configurations consisting of $4$ points in $\RR^2$}
\label{fig:4pointsR2}
\end{figure}
}{%
\begin{figure}[htpb]
\centering
 \subbottom[]{\label{sfig:4pointsR2_1}\includegraphics[width=.2\linewidth]{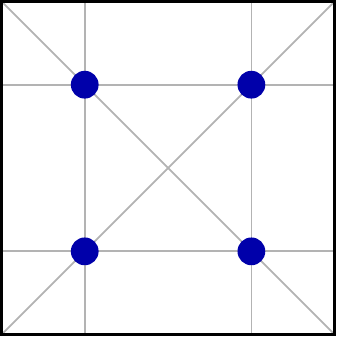}}\quad
 \subbottom[]{\label{sfig:4pointsR2_2}\includegraphics[width=.2\linewidth]{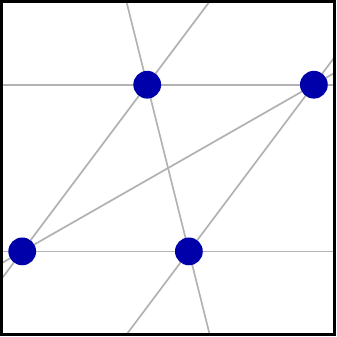}}\quad
 \subbottom[]{\label{sfig:4pointsR2_3}\includegraphics[width=.2\linewidth]{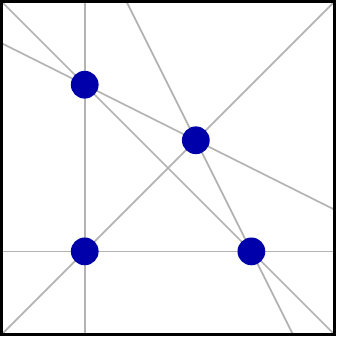}}\quad
 \subbottom[]{\label{sfig:4pointsR2_4}\includegraphics[width=.2\linewidth]{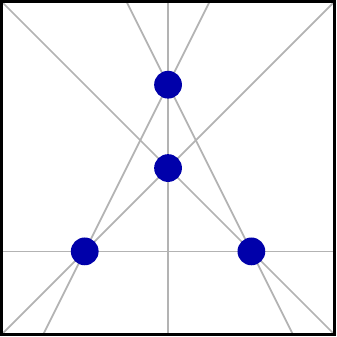}}
\caption{Four point configurations consisting of $4$ points in $\RR^2$}
\label{fig:4pointsR2}
\end{figure}
}
}

Just like here, throughout the whole document we use upright boldface letters such as $\vv A$ to denote subsets of $\RR^d$ (and labeled subsets such as point configurations). Points and vectors in $\RR^d$ are written in boldface lowercase letters such as $\vv a$.
The $i$-th coordinate of~$\vv a$ is denoted by $\vvc a_i$. The standard basis vectors are written as $\vv e_1,\dots,\vv e_d$; the all-zeros vector of appropriate size is represented as~$\veczero$\index{$\veczero$}; and $\matid_n$\index{$\matid_n$} is the $n\times n$ identity matrix.

Figure~\ref{fig:4pointsR2} shows four configurations of four points in the plane. 
The ones in~\ref{sfig:4pointsR2_1} and~\ref{sfig:4pointsR2_2} are strongly related; indeed, they are \defn{affinely isomorphic}\index{affine!isomorphism}. In general $\vv A\subset\RR^d$ and $\vv B\subset\RR^e$ are affinely isomorphic 
 if there is an affine map $f:\RR^d\rightarrow\RR^e$, $\vv x\mapsto M\vv x+ \vv x_0$, that is a bijection between the points of~$\vv A$ and the points of~$\vv B$.
One might also argue that \ref{sfig:4pointsR2_1}~is closer to~\ref{sfig:4pointsR2_3} than to~\ref{sfig:4pointsR2_4}. One would be right, and we will see why in Section~\ref{sec:orientedmatroids}.
 
In the figure we can also see some lines, the \defn{hyperplanes} \index{hyperplane} spanned by subsets of $\vv A$. An (oriented) affine hyperplane\index{hyperplane!affine} $\vvh H$\index{$\vvh H$} --- defined by a normal vector $\vv v\in \RR^d$ and a scalar $c\in\RR$ --- is the set of points $\vvh H:=\set{\vv x}{\sprod{\vv v}{\vv x}=c}$, where $\sprod{\cdot}{\cdot}$ is the standard inner product on~$\RR^d$. Its positive and negative sides are the halfspacess\index{halfspace} $\vvh H^+:=\set{\vv x}{\sprod{\vv v}{\vv x}>c}$\index{$\vvh H$!$\vvh {H}^\pm$} and $\vvh H^-:=\set{\vv x}{\sprod{\vv v}{\vv x}<c}$, respectively. We will also use the closed halfspaces\index{halfspace!closed} $\ol {\vvh H}^+=\vvh H\cup \vvh H^+$\index{$\vvh H$!$\ol{\vvh {H}}^\pm$} and $\ol{\vvh H}^-=\vvh H\cup \vvh H^-$. For a subset $\vv B\subset\vv A$, we say that a hyperplane $\vvh H$ is \defn{spanned} by $\vv B$ if $\vvh H=\aff(\vv B)$. 
If the scalar $c=0$, we say that $\vvh H$ is a \defn{linear hyperplane}\index{hyperplane!linear}. A linear hyperplane $\vvh H$ is \defn{spanned} by a set of vectors $\vv V$ if $\vvh H=\lin(\vv V)$, where $\lin(\vv V)$\index{$\lin(\vv V)$}~is the \defn{linear span}\index{linear!span} of~$\vv V\subset\RR^s$, defined as 
\[\lin(\vv V):=\set{\vv w\in \RR^s}{\vv w=\sum_{i=1}^n\gl_i \vv v_i; \;\;\text{ for }\vv v_i\in \vv V\text{ and } \gl_i\in \RR}.\]

It is sometimes convenient to consider the \defn{homogenization}\index{homogenization} of~$\vv A$. This is the vector configuration $\hom(\vv A)$\index{$\hom(\vv A)$} in the linear space $\RR^{d+1}$ consisting of the vectors $\vv v_i=\binom{\vv a_i}{1}$ obtained by appending a `$1$' to the coordinates of~$\vv a_i$. Affine hyperplanes spanned by subsets of~$\vv A$ are in bijection to linear hyperplanes spanned by subsets of $\hom(\vv A)$. See Figure~\ref{fig:hom} for an example.
\iftoggle{bwprint}{%
\begin{figure}[htpb]
\centering
\includegraphics[width=.6\linewidth]{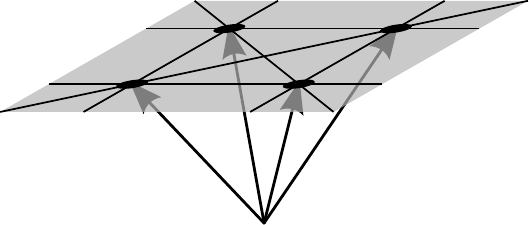}
\caption[Homogenization of a point configuration.]{The homogenization of the point configuration from Figure~\ref{sfig:4pointsR2_1}.}\label{fig:hom}
\end{figure}
}{%
\begin{figure}[htpb]
\centering
\includegraphics[width=.6\linewidth]{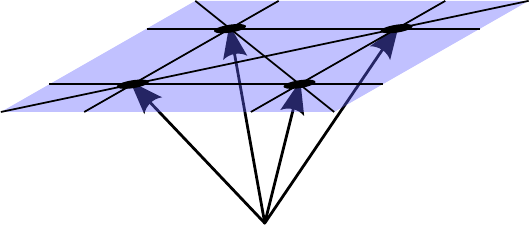}
\caption[Homogenization of a point configuration.]{The homogenization of the point configuration from Figure~\ref{sfig:4pointsR2_1}.}\label{fig:hom}
\end{figure}
}

\medskip

The \defn{convex hull}\index{convex hull} {$\conv(\vv A)$}\index{$\conv(\vv A)$} of~$\vv A$ is the intersection of all closed halfspaces containing $\vv A$. Equivalently,
\[\conv(\vv A)\!:=\!\set{\!\vv x\in \RR^d}{\vv x=\sum_{i=1}^n\!\vvc \gl_i \vv a_i;\text{ for }\vv a_i\in \vv A, \vvc \gl_i\in \RR_{\geq 0}\text{ and } \sum_{i=1}^n\!\vvc \gl_i=1\!}.\]

A (convex) \defn{polytope}\index{polytope} $\vv P$ is the convex hull of a finite set of points in $\RR^d$, or equivalently, a bounded intersection of finitely many closed halfspaces. The \defn{dimension} of~$\vv P$ is the dimension of its affine span, and by a \defn{$d$-polytope} we mean a $d$-dimensional polytope.

\iftoggle{bwprint}{%
\begin{figure}[htpb]
\centering
 \subbottom[$\simp{3}$]{\label{sfig:simplex}\includegraphics[width=.25\linewidth]{Figures/simplex}}\qquad\quad
 \subbottom[$\cube{3}$]{\label{sfig:cube}\includegraphics[width=.25\linewidth]{Figures/cube}}\qquad\quad
 \subbottom[$\cros{3}$]{\label{sfig:cross}\includegraphics[width=.25\linewidth]{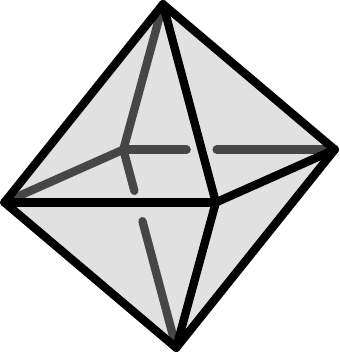}}
\caption{Three $3$-polytopes}
\label{fig:3polytopes}
\end{figure}
}{%
\iftoggle{print}{%
\begin{figure}[htpb]
\centering
 \subbottom[$\simp{3}$]{\label{sfig:simplex}\includegraphics[width=.25\linewidth]{Figures/simplex}}\qquad\quad
 \subbottom[$\cube{3}$]{\label{sfig:cube}\includegraphics[width=.25\linewidth]{Figures/cube}}\qquad\quad
 \subbottom[$\cros{3}$]{\label{sfig:cross}\includegraphics[width=.25\linewidth]{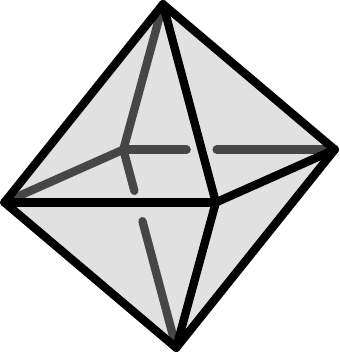}}
\caption{Three $3$-polytopes}
\label{fig:3polytopes}
\end{figure}
}{%
\begin{figure}[htpb]
\centering
 \subbottom[$\simp{3}$]{\label{sfig:simplex}\includegraphics[width=.25\linewidth]{Figures/simplex_col}}\qquad\quad
 \subbottom[$\cube{3}$]{\label{sfig:cube}\includegraphics[width=.25\linewidth]{Figures/cube_col}}\qquad\quad
 \subbottom[$\cros{3}$]{\label{sfig:cross}\includegraphics[width=.25\linewidth]{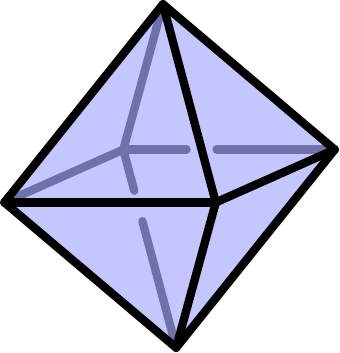}}
\caption{Three $3$-polytopes}
\label{fig:3polytopes}
\end{figure}
}
}

Figure~\ref{fig:3polytopes} shows three $3$-polytopes with~$4$, $8$ and~$6$ vertices and~$4$, $6$ and~$8$ facets respectively. These are the respective $3$-dimensional instances of the 
\defn{$d$-simplex}~$\simp{d}$\index{simplex}\index{$\simp{d}$}, the \defn{$d$-cube}~$\cube{d}$\index{cube}\index{$\cube{d}$} and the \defn{$d$-crosspolytope}~$\cros{d}$\index{crosspolytope}\index{$\cros{d}$}. Their standard versions are given by
\begin{align*}
 \simp{d}&:=\conv\{\vv e_1,\dots,\vv e_{d+1}\}=\set{\vv x\in \RR^{d+1}}{\sum\nolimits_i{x_i}=1,x_i\geq 0},\\
 \cube{d}&:=\conv\left\{\left\{+1,-1\right\}^d\right\}=\set{\vv x\in \RR^{d}}{-1\leq x_i\leq 1},\\
 \cros{d}&:=\conv\{\pm\vv e_1,\dots,\pm \vv e_{d}\}=\set{\vv x\in \RR^{d}}{\sum\nolimits_i{|x_i|}\leq 1}.
\end{align*}

\medskip

A \defn{supporting hyperplane}\index{hyperplane!supporting} of a polytope $\vv P$ is a hyperplane with $\vv P\subset \ol {\vvh H}^+$. The intersection $\vv P\cap {\vvh H}$ of~$\vv P$ with one of its supporting hyperplanes ${\vvh H}$ is a \defn{face}\index{face} of~$\vv P$, which is also a polytope. We consider the whole polytope $\vv P$ to be a face of itself. A $k$-dimensional face is called a \defn{$k$-face}. Faces of a $d$-polytope of dimension $0$, $1$ and $d-1$ are called \defn{vertices}\index{vertex}, \defn{edges}\index{edge} and \defn{facets}\index{facet} respectively. We write $\vv F\leq \vv P$ to denote that $\vv F$ is a face of~$\vv P$. A $k$-face $\vv F$ is \defn{proper}\index{face!proper} if $k<d$, which is denoted $\vv F<\vv P$\index{$\vv F<\vv P$}. If $d\geq1$ the \defn{relative interior}\index{relative interior} of a $d$-polytope ~$\vv P$ is $\relint(\vv P)=\vv P\setminus \bigcup_{\vv F<\vv P} \vv F$\index{$\relint(\vv P)$}, while $\relint\{\vv a\}:=\{\vv a\}$ when $\vv a$ is a single point.

The \defn{face lattice}\index{face lattice} of~$\vv P$ is the set of faces of~$\vv P$, partially ordered by inclusion. Two polytopes $\vv P$ and $\vv Q$ are \defn{combinatorially equivalent}\index{combinatorial equivalence} 
if there is a bijection between their sets of faces that preserves inclusions. Figure~\ref{fig:cubeblown} shows the proper faces of~$\cube{3}$ and its face lattice. The \defn{$f$-vector}\index{$f$-vector} of~$\vv P$ is the vector $\bm f(\vv P):=(f_{-1},f_0,f_1,\dots,f_d)\in \NN^{d+2}$\index{$\bm f(\vv P)$}, where $f_i$~is the number of $i$-dimensional faces of~$\vv P$ (by convention, the empty set is considered a face of dimension $-1$). For example, the $f$-vector of~$\cube{3}$ is $\bm f(\cube{3})=(1,8,12,6,1)$.
Observe how the face lattice of the polytope $\conv (\vv A)$ can be read off from the \defn{polyhedral cone} defined by its homogenization; $\cone(\hom(\vv A)):=\set{\lambda\vv v}{\vv v\in\hom(\vv A),\, 0\leq\lambda\in\RR}$\index{$\cone(\vv V)$}.

\iftoggle{bwprint}{%
\begin{figure}[htpb]
\centering
\raisebox{-0.5\height}{\includegraphics[width=.25\linewidth]{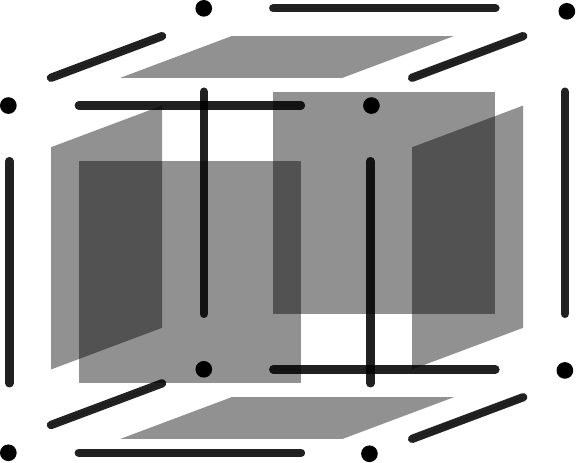}}\qquad\qquad
\raisebox{-0.5\height}{\includegraphics[width=.5\linewidth]{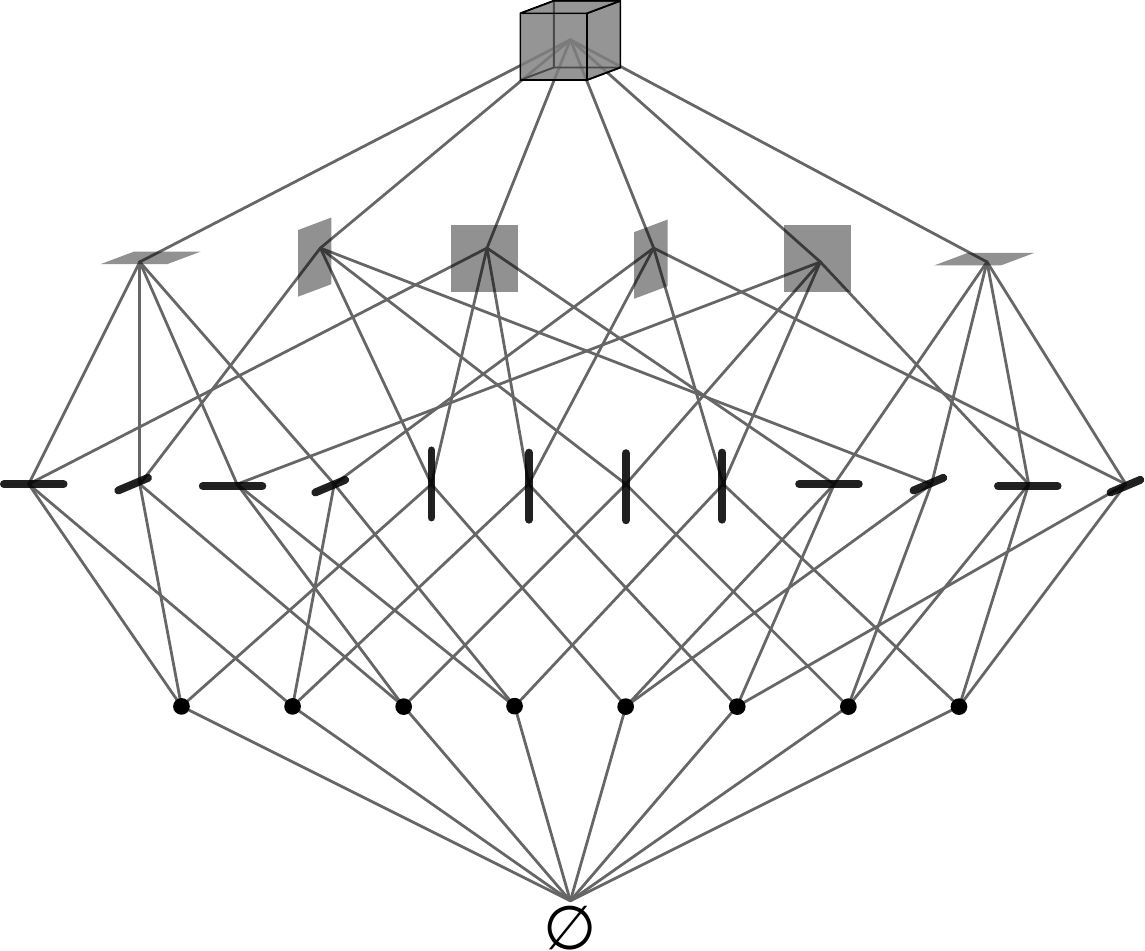}}
\caption{The face lattice of $\cube{3}$.}
\label{fig:cubeblown}
\end{figure}
}{%
\begin{figure}[htpb]
\centering
\raisebox{-0.5\height}{\includegraphics[width=.25\linewidth]{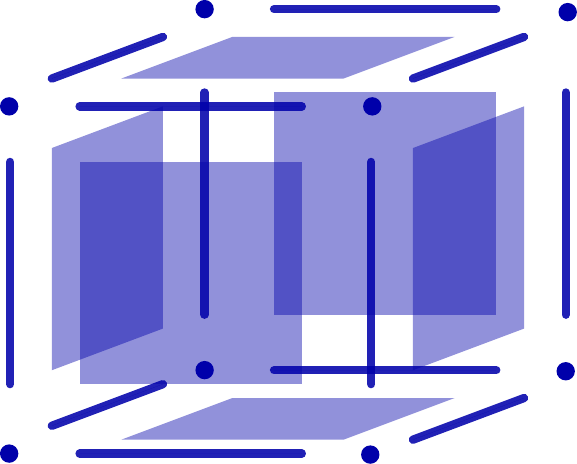}}\qquad\qquad
\raisebox{-0.5\height}{\includegraphics[width=.5\linewidth]{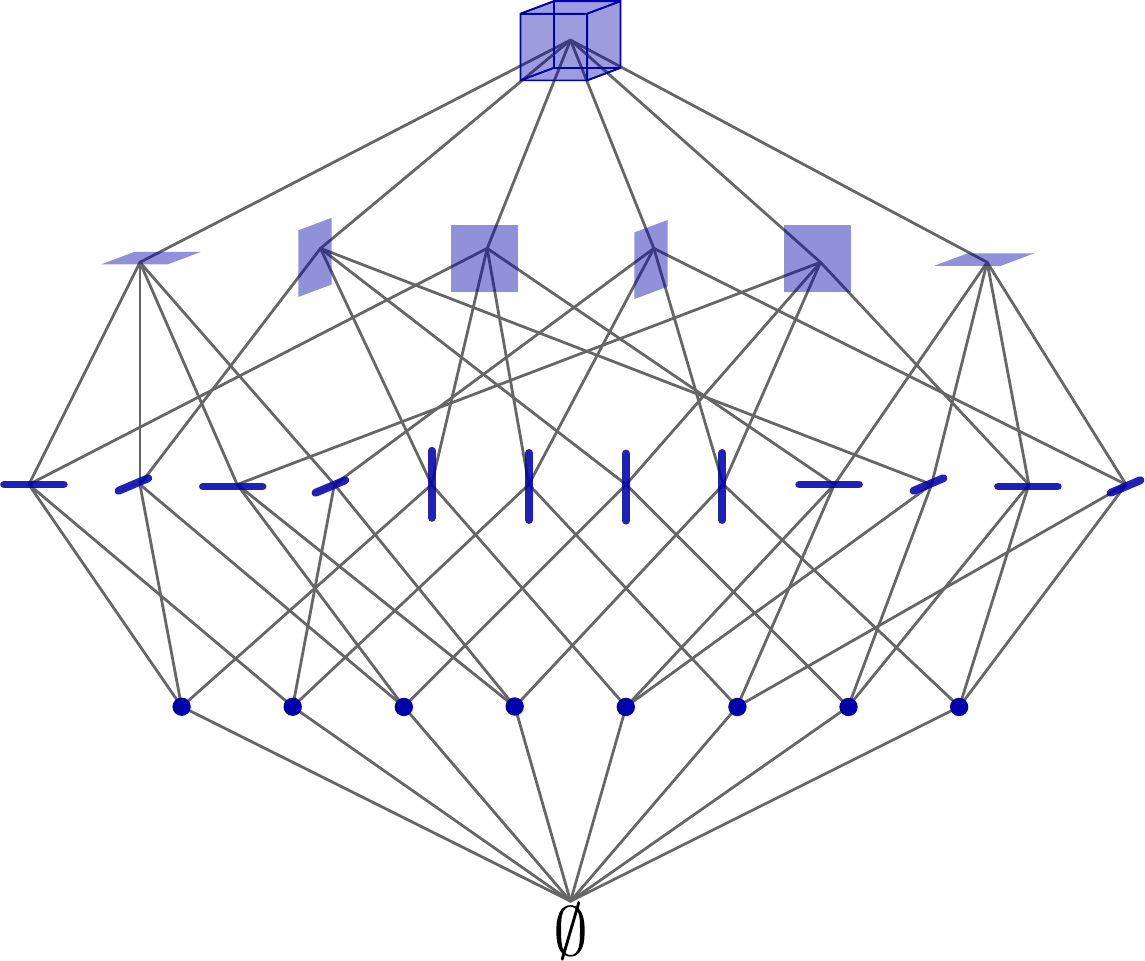}}
\caption{The face lattice of $\cube{3}$.}
\label{fig:cubeblown}
\end{figure}
}

The set of vertices of~$\vv P$ is denoted $\verts(\vv P)$\index{$\verts(\vv P)$}. For each face~$\vv F<\vv P$, its set of vertices is $\verts(\vv F)=\vv F\cap \verts(\vv P)$. This explains why, from a combinatorial point of view, we often identify a face with the corresponding subset of $\verts(\vv P)$. For a general point configuration $\vv A$, we identify a face $\vv F$ of $\conv(\vv A)$ with~$\vv F\cap \vv A$, or even with the set of labels of the points in~$\vv F\cap \vv A$.
\\

A (polyhedral) \defn{subdivision}\index{subdivision} of~$\vv A$ is a collection $\cS$ of subsets of~$\vv A$ that satisfies:
\begin{enumerate}
 \item If $\vv C\in \cS$ and $\vv D=\vv C\cap \vv F$ for some face $\vv F\leq \conv(\vv C)$, then $\vv D\in \cS$.
 \item $\bigcup_{\vv C\in \cS} \conv(\vv C)=\conv(\vv A)$.
 \item If $\vv C\neq \vv D$ are in $\cS$, then $\relint(\vv C)\cap \relint(\vv D)=\emptyset$. 
\end{enumerate}
A \defn{triangulation}\index{triangulation} of~$\vv A$ is a subdivision $\cT$ where the points in each $\vv C\in \cT$ are affinely independent, and thus $\conv (\vv C)$ is a simplex. Some examples are depicted in Figure~\ref{fig:subdivisions}.

\iftoggle{bwprint}{%
\begin{figure}[htpb]
\centering
{\includegraphics[width=.8\linewidth]{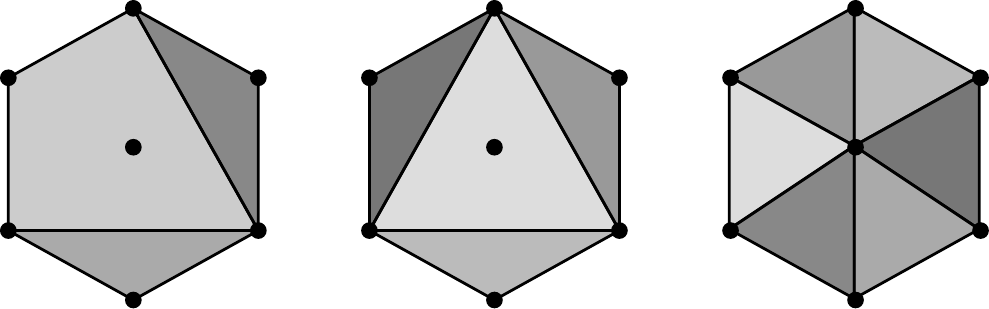}}
\caption{Three subdivisions of a planar configuration of $7$ points, two of which are triangulations.}
\label{fig:subdivisions}
\end{figure}
}{%
\begin{figure}[htpb]
\centering
{\includegraphics[width=.8\linewidth]{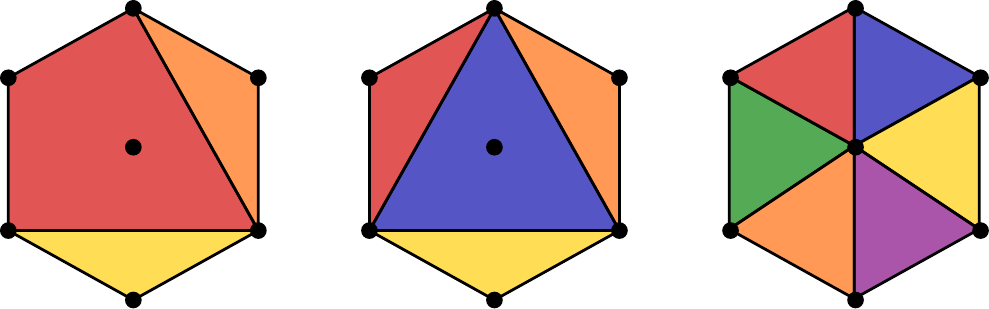}}
\caption{Three subdivisions of a planar configuration of $7$ points, two of which are triangulations.}
\label{fig:subdivisions}
\end{figure}
}

\section{Vector configurations}

A \defn{vector configuration}\index{vector configuration} $\vv V = \{\vv v_1 ,\dots ,\vv v_n\}$ is a finite collection of (labeled) vectors in the linear space $\RR^r$, and its  \defn{rank}\index{rank} is $\rank(\vv V):=\dim(\lin(\vv V))$. We 
define two vector spaces that capture its geometry, namely its \defn{linear evaluations}\index{evaluation!linear} $\Val(\vv V)$\index{$\Val(\vv V)$} and its \defn{linear dependences}\index{dependence!linear} $\Dep({\vv V})$\index{$\Dep({\vv V})$}:
\begin{align*}
\Val(\vv V)&:=\set{\vv \gn\in \RR^n}{\vvc \gn_i= \sprod{\vv c}{\vv v_i},\text{ where }\vv c\in\left(\RR^r\right)^*},\\ 
\Dep(\vv V)&:=\set{\vv \gl \in \RR^n}{ \sum\nolimits_i \vvc \gl_i \vv v_i =0}.
\end{align*}

Given a $d$-dimensional point configuration $\vv A=\{\vv a_1,\dots,\vv a_n\}$, we define its \defn{affine evaluations}\index{evaluation!affine} $\aVal(\vv A)$\index{$\aVal(\vv A)$} and its \defn{affine dependences}\index{dependence!affine} $\aDep(\vv A)$\index{$\aDep(\vv A)$} from its homogenization: $\aVal(\vv A)=\Val(\hom{\vv A})$ and $\aDep=\Dep(\hom{\vv A})$.

Observe that one can tell the face lattice of $\conv(\vv A)$ from its affine functions $\aVal(\vv A)$, since a subset $\vv B$ of points of~$\vv A$ lie in a common supporting hyperplane of $\conv(\vv A)$ if and only if there is an affine evaluation $\vv \gn$ with $\vvc \gn_i=0$ for $\vv a_i\in \vv B$ and $\vvc \gn_i\geq 0$ otherwise.\\

This combinatorial information of $\vv V$ is extracted as follows. A \defn{signed subset}\index{signed set} $X$ of a ground set $E$ is a pair $X=(X^+,X^-)$\index{$(X^+,X^-)$} of disjoint subsets of $E$, one of them called positive and the other negative. Its \defn{support}\index{support} is $\ul X=X^+\cup X^-$\index{$\ul X$}, and the set $E\setminus \ul X$ is abbreviated by $X^0$\index{$X^0$}. Alternatively, we also view a signed subset as a function from $E$ to $\{+,-,0\}$, or even to $\{\pm 1,0\}$. Hence, we say 
$X(e)=+$ or 
$X(e)>0$ for $e\in X^+$\index{$X(e)$}.

For each linear dependence $\vv \gl \in\Dep (\vv V)$, the pair $(U^+,U^-)$ with $U^+:=\set{i}{\vvc \gl_i>0}$ and $U^-:=\set{i}{\vvc \gl_i<0}$ defines a signed subset of the set of labels of~$\vv V$. It is called a \defn{vector}\index{vector} of $\cM(\vv V)$, the \defn{oriented matroid}\index{oriented matroid} of~$\vv V$, which will be formally defined in Section~\ref{sec:orientedmatroids}. The set of vectors of $\cM(\vv V)$\index{$\cM(\vv V)$} is denoted by $\ve(\cM(\vv V))$ or just $\ve(\vv V)$\index{$\ve(\cM)$}.
If $U^-=\emptyset$, we say that $U$ is a \defn{positive vector}\index{vector!positive}. 

We define a partial order for vectors where $U\leq V$ if and only if $U^+\subseteq V^+$ and $U^-\subseteq V^-$. The \defn{circuits}\index{circuit} of $\cM(\vv V)$ are the minimal elements in this partial order, and the set of circuits of $\cM(\vv V)$ is denoted $\ci(\cM(\vv V))$ or~$\ci(\vv V)$\index{$\ci(\cM(\vv V))$}.

We described a vector $U$ of $\cM(\vv V)$ as a signed subset of labels of vectors in $\vv V$. However, we often abuse notation and identify $U$,~$U^+$ and~$U^-$ with the vector subconfigurations $\vv V_U:=\set{\vv v_i\in \vv V}{i\in U}$, $\vv V_{U^+}:=\set{\vv v_i\in \vv V}{i\in U^+}$ and $\vv V_{U^-}:=\set{\vv v_i\in\vv  V}{i\in U^-}$ respectively. Hence, we use $\vv v_i \in U$ and $i \in U$ interchangeably.

In this context, we say that a subconfiguration $\vv W\subseteq\vv  V$ is a positive vector when there is a positive vector $X$ with $\vv V_{X^+}=\vv W$. Observe that $\vv W$~is a positive vector if and only if the origin $\veczero$ is contained in the relative interior of the convex hull of~$\vv W$ (seen as points instead of vectors).

Analogously, we define the set of \defn{covectors}\index{covector} $\cov(\vv V)$\index{$\cov(\cM(\vv V))$} and the set of \defn{cocircuits}\index{cocircuit} $\co(\vv V)$\index{$\co(\cM(\vv V))$} of $\cM(\vv V)$ by extracting a signed subset from each linear evaluation $\vv \gn\in \Val(\vv V)$. Observe that each covector $C$ is defined by a linear hyperplane~${\vvh H}$ where $C^+={\vvh H}^+\cap \vv V$ and $C^-={\vvh H}^-\cap \vv V$. Cocircuits of $\vv V$ are those covectors defined by linear hyperplanes spanned by subsets of~$\vv V$.

In the same way, the oriented matroid $\cM(\vv A)$ of a point configuration $\vv A$ is the oriented matroid of its homogenization $\cM(\hom(\vv A))$ (see Figure~\ref{fig:affinecovector} for an example of an affine covector). 

\iftoggle{bwprint}{%
\begin{figure}[htpb]
\centering
{\includegraphics[width=.6\linewidth]{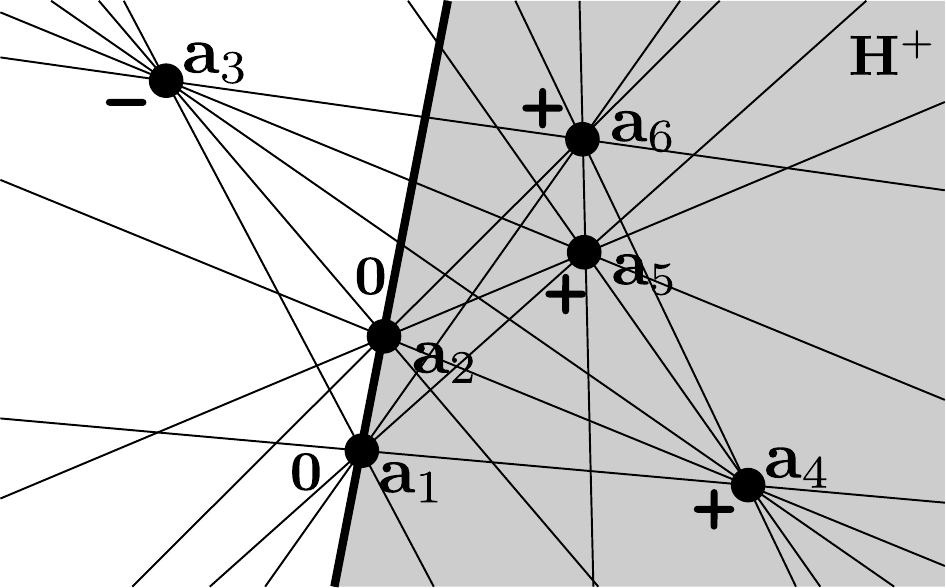}}
\caption[An affine covector defined by a hyperplane ${\vvh H}$.]{An affine covector $C$ defined by a hyperplane ${\vvh H}$. In this example, $C^+=\{\vv a_4,\vv a_5,\vv a_6\}$ and $C^-=\{\vv a_3\}$. Each point is labeled as $+$, $-$ or $0$ accordingly.}
\label{fig:affinecovector}
\end{figure}
}{%
\begin{figure}[htpb]
\centering
{\includegraphics[width=.6\linewidth]{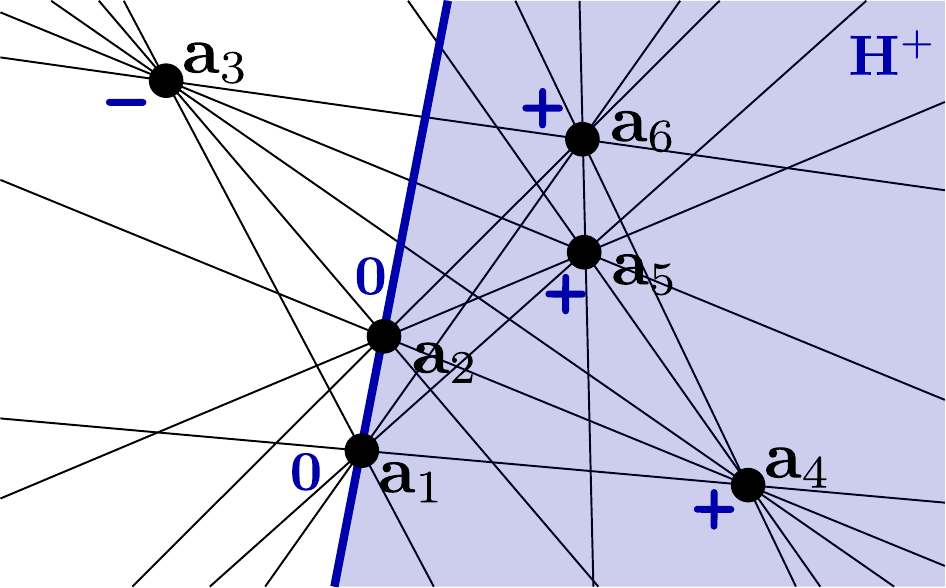}}
\caption[An affine covector defined by a hyperplane ${\vvh H}$.]{An affine covector $C$ defined by a hyperplane ${\vvh H}$. In this example, $C^+=\{\vv a_4,\vv a_5,\vv a_6\}$ and $C^-=\{\vv a_3\}$. Each point is labeled as $+$, $-$ or $0$ accordingly.}
\label{fig:affinecovector}
\end{figure}
}

As we will see later, an oriented matroid $\cM(\vv V)$ can be recovered from any of the sets $\ci(\vv V)$, $\ve(\vv V)$, $\co(\vv V)$ or~$\cov(\vv V)$.

\subsection{Gale duality}

\index{Gale dual}Let $\vv V$ be a vector configuration of size $n$ and rank~$r$; and let $M$ be the $r\times n$ matrix of rank~$r$ whose column vectors contain the coordinates of the vectors of~$\vv V$.  Choose a basis $(\vv b_1,\ldots,\vv b_{n-r})$ of the kernel of $M$, and denote by $\Gale M$ the $n\times (n-r)$ matrix whose column vectors are these~$\vv b_i$'s. In other words, $M\Gale M=\matzero_{r\times (n-r)}$ and $\rank(\Gale M)=n-r$. Finally, define~$\Gale{\vv V}$ to be the configuration of rank $(n-r)$ consisting of the row vectors of $\Gale M$. The configuration $\Gale{\vv V}$\index{$\Gale{\vv V}$} (that shares labels with $\vv V$) is called a \defn{Gale dual} of~$\vv V$.

The Gale dual is well defined up to a linear transformation. With this caveat, it is a duality~---~$\pGale{\Gale {\vv V}}\cong\vv V$~---~whose key property is that it translates linear evaluations into linear dependencies (cf. \cite[Lecture~6]{Ziegler1995}):

\begin{theorem}\label{thm:Galeduality}
 $\Val(\vv V)=\Dep(\Gale{\vv V})$ and  $\Dep(\vv V)=\Val(\Gale{\vv V})$.\qed
\end{theorem}

In this equality, we understand $\vv V$ and $\Gale{\vv V}$ as labeled configurations: the label of the vector in the $i$-th column of $M$ must coincide with the label of the vector in the $i$-th row of~$\Gale M$.\\ 

Usually, $\vv A$ will be a $d$-dimensional configuration of $n$~points, and $\vv V$~the vector configuration of rank~$(n-d-1)$ dual to $\hom(\vv A)$. Hence, $\aVal(\vv A)=\Dep(\vv V)$ and $\aDep(\vv A)=\Val(\vv V)$. Up to a projective transformation, $\vv A$ can be recovered from $\vv V$, and for convenience we write $\Gale{\vv V}=\vv A$ and $\Gale{\vv A}=\vv V$. The vector configuration $\Gale{\vv A}$ is a \defn{Gale diagram}\index{Gale diagram} or the \defn{Gale dual}\index{Gale dual} of~$\vv A$.

\begin{example}\label{ex:Gale}

Consider the point configuration $\vv A$ in $\RR^2$ whose homogenized coordinates are recorded in the following matrix (cf.~Figure~\ref{sfig:GaleA}):
\[
M=\kbordermatrix{&\phantom{.}\vv a_1&\phantom{-}\vv a_2&\phantom{-}\vv a_3&\phantom{-}\vv a_4&\phantom{-}\vv a_5\\
&\phantom{.}0&\phantom{-}0& \phantom{-}0& \phantom{-}1& \phantom{-}2\\
&\phantom{.}0& \phantom{-}1& \phantom{-}2& \phantom{-}0& \phantom{-}0
\\
 &\phantom{.}{1}& \phantom{-}{1}& \phantom{-}{1}& \phantom{-}{1}&\phantom{-}{1}
}.
\]
Its Gale dual is the vector configuration $\vv V=\Gale{\vv A}$ in $\RR^{n-d-1=2}$ defined by the columns of the matrix~$\trans{\left(\Gale M\right)}$ (cf.~Figure~\ref{sfig:GaleV}):
\[\trans{\left(\Gale{M}\right)}=
\kbordermatrix{
&\phantom{-}\vv v_1&\phantom{-}\vv v_2&\phantom{-}\vv v_3&\phantom{-}\vv v_4&\phantom{-}\vv v_5\\
&\phantom{-}1& \phantom{-}0& \phantom{-}0& -2& \phantom{-}1\\
&\phantom{-}1& -2& \phantom{-}1& \phantom{-}0& \phantom{-}0}.\]

\iftoggle{bwprint}{%
\begin{figure}[htpb]
\centering
 \subbottom[ $\vv A$]{\label{sfig:GaleA}\includegraphics[width=.25\textwidth]{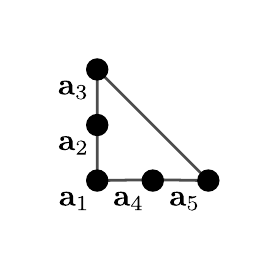}}\qquad\qquad\qquad\qquad
 \subbottom[$\vv V=\Gale{\vv A}$]{\label{sfig:GaleV}\includegraphics[width=.25\textwidth]{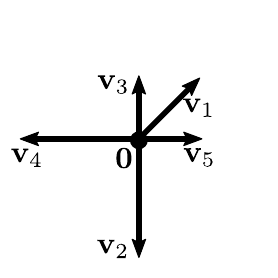}}
 \caption{A point configuration $\vv A$, and its Gale dual $\vv V$.}
\label{fig:exampleGaleDual}
\end{figure}
}{%
\iftoggle{print}{%
\begin{figure}[htpb]
\centering
 \subbottom[ $\vv A$]{\label{sfig:GaleA}\includegraphics[width=.25\textwidth]{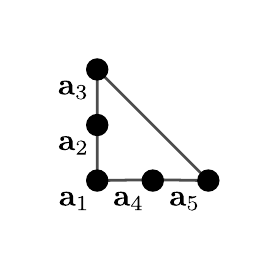}}\qquad\qquad\qquad\qquad
 \subbottom[$\vv V=\Gale{\vv A}$]{\label{sfig:GaleV}\includegraphics[width=.25\textwidth]{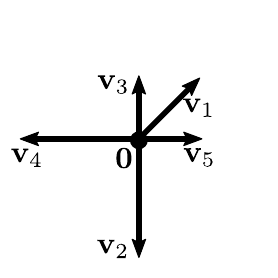}}
 \caption{A point configuration $\vv A$, and its Gale dual $\vv V$.}
\label{fig:exampleGaleDual}
\end{figure}
}{%
\begin{figure}[htpb]
\centering
 \subbottom[ $\vv A$]{\label{sfig:GaleA}\includegraphics[width=.25\textwidth]{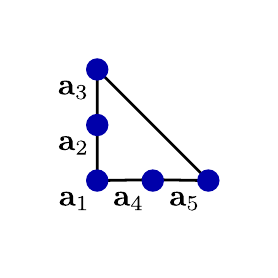}}\qquad\qquad\qquad\qquad
 \subbottom[$\vv V=\Gale{\vv A}$]{\label{sfig:GaleV}\includegraphics[width=.25\textwidth]{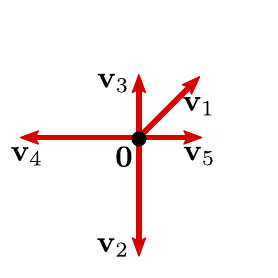}}
 \caption{A point configuration $\vv A$, and its Gale dual $\vv V$.}
\label{fig:exampleGaleDual}
\end{figure}
}
}

\end{example}

Next, we present a lemma that is a particular case of Theorem~\ref{thm:Galeduality} for supporting hyperplanes of point configurations.
\begin{lemma}\label{lem:gale}
Let $\vv A:=\{\vv a_1,\ldots,\vv a_n\}\subset\RR^d$ be a $d$-dimensional point configuration and $\vv V:=\{\vv v_1,\ldots,\vv v_n\}\subset\RR^{n-d-1}$ its Gale dual. For any $I\subseteq [n]$, let $F:=\set{\vv a_i\in \vv A}{i\in I}$ and $\bGale{F}:=\set{\vv v_i\in\vv V}{i\notin I}$\index{$\bGale{F}$}. Then:
\begin{enumerate}[(i)]
 \item\label{it:gale1} $F$ is contained in a supporting hyperplane of~$\conv(\vv A)$ if and only if 
$\bGale{F}$ contains a positive vector of $\cM(\vv V)$.
 \item\label{it:gale2} $F$ are the only points contained in a supporting hyperplane of~$\conv(\vv A)$ if and only if
$\bGale{F}$  is a positive vector of $\cM(\vv V)$.
\end{enumerate}
\end{lemma}

One can check the assertions of this lemma in the configurations~$\vv A$ and~$\vv V$ from Example~\ref{ex:Gale} (cf.~Figure~\ref{fig:exampleGaleDual}). For example, the fact that $\vv a_1$~and~$\vv a_4$ are contained in a common facet of $\conv(\vv A)$ is reflected in the fact that~$\veczero\in\conv\left(\vv v_2,\vv v_3,\vv v_5\right)$. On the other hand, since $\{\vv a_1,\vv a_4\}$ is not a face of~$\vv A$, $\veczero\notin\relint\left(\conv\left(\vv v_2,\vv v_3,\vv v_5\right)\right)$.\\

\begin{observation}
 Observe that multiplying the vectors of~$\vv V$ by positive scalars does not change its oriented matroid.
 Hence, we will often work with normalized versions of the Gale dual,  preserving the combinatorial type of the primal configuration.
\end{observation}
\smallskip

The following theorem characterizes vector configurations that are the Gale dual of a point configuration.

\begin{theorem}
A vector configuration $\vv V$ in $\RR^r$ is a Gale dual of a point configuration in $\RR^{n-r-1}$ (up to positive rescaling) if and only if $\vv V$ \defn{positively spans} $\RR^r$. That is, either $r=0$ or $|{\vvh H}^+ \cap \vv V|\geq 1$ for every oriented linear hyperplane ${\vvh H}$.\qed
\end{theorem}

We can analogously characterize configurations without repeated points:

\begin{lemma}\label{lem:dualnorepeatedpoints}
 A $d$-dimensional configuration $\vv A$ of $n$ points with Gale dual $\vv V:=\Gale{\vv A}$ has no repeated points if and only if either $n=d+1$, or for every linear hyperplane~${\vvh H}$, $|{\vvh H}^+\cap \vv V|\geq 2$ or $|{\vvh H}^- \cap \vv V|\geq 2$. \qed
\end{lemma}

\dots and Gale duals of point configurations in convex position:

\begin{theorem}\label{thm:galepolytope}
A vector configuration $\vv V$ of rank $r$ is a Gale dual of a point configuration in convex position in $\RR^{n-r-1}$ (up to positive rescaling) if and only if $\vv V$ is \defn{positively $2$-spanning}. That is, either $r=0$ or $|{\vvh H}^+ \cap \vv V|\geq 2$ for every oriented linear hyperplane ${\vvh H}$.\qed
\end{theorem}

\iftoggle{bwprint}{%
\begin{figure}[htpb]
\centering
\subbottom[ $\vv V\rightarrow \vv X$]{\label{sfig:affineGaleVA}\includegraphics[width=.53\linewidth]{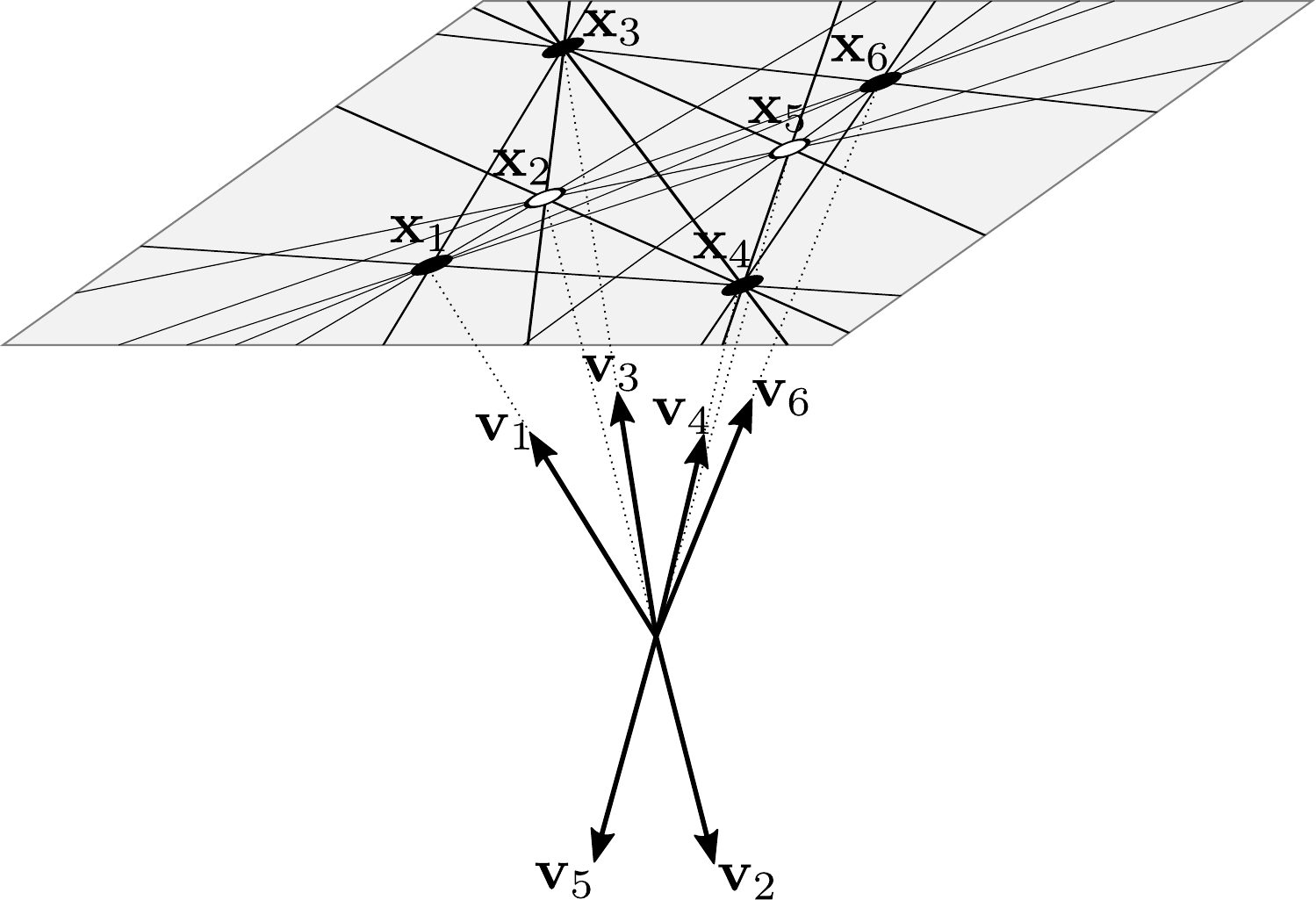}}
\subbottom[ A covector defined by ${\vvh H}$.]{\raisebox{.8cm}{\label{sfig:affineGalecovector}\includegraphics[width=.45\linewidth]{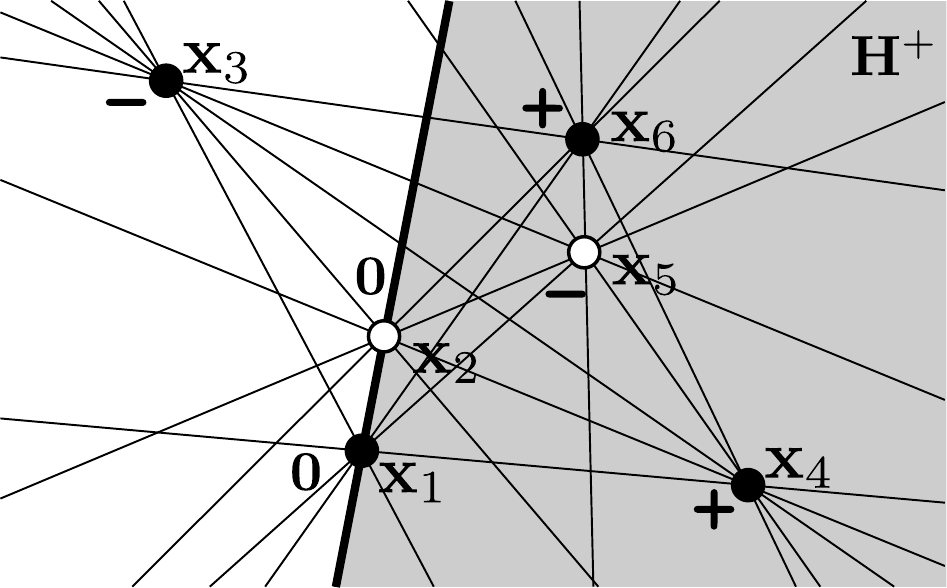}}}
\caption[Example of an affine Gale diagram]{The $3$-dimensional vector configuration $\vv V$ is the Gale dual of the vertex set of an hexagon. Its affine Gale diagram $\vv X$ in $\RR^2$, and an example of a covector $C$ with $C^+=\{\vv x_4,\vv x_3\}$ and $C^-=\{\vv x_3,\vv x_5\}$.}
\label{fig:affineGale}
\end{figure}
}{%
\begin{figure}[htpb]
\centering
\subbottom[ $\vv V\rightarrow \vv X$]{\label{sfig:affineGaleVA}\includegraphics[width=.53\linewidth]{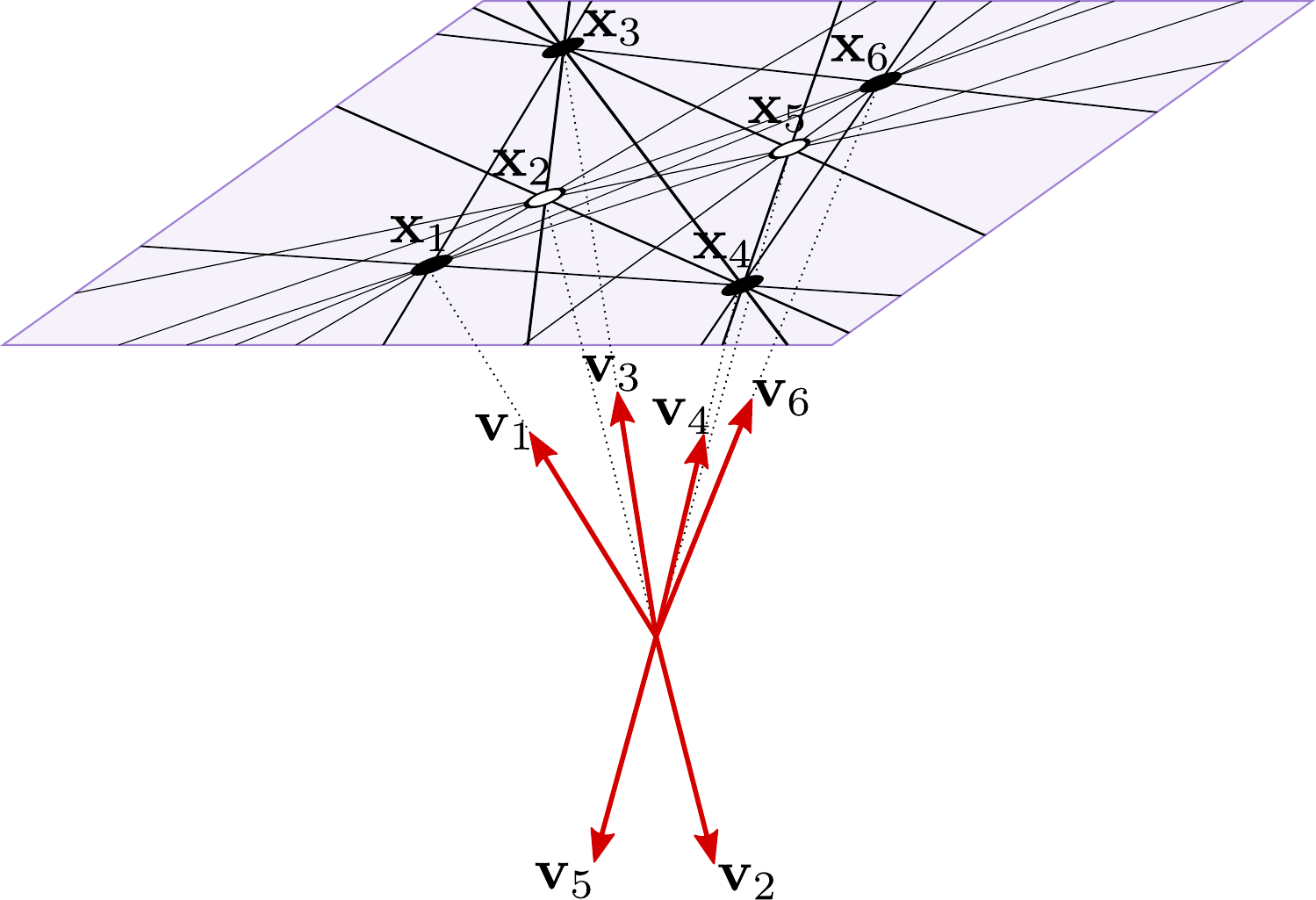}}
\subbottom[ A covector defined by ${\vvh H}$.]{\raisebox{.8cm}{\label{sfig:affineGalecovector}\includegraphics[width=.45\linewidth]{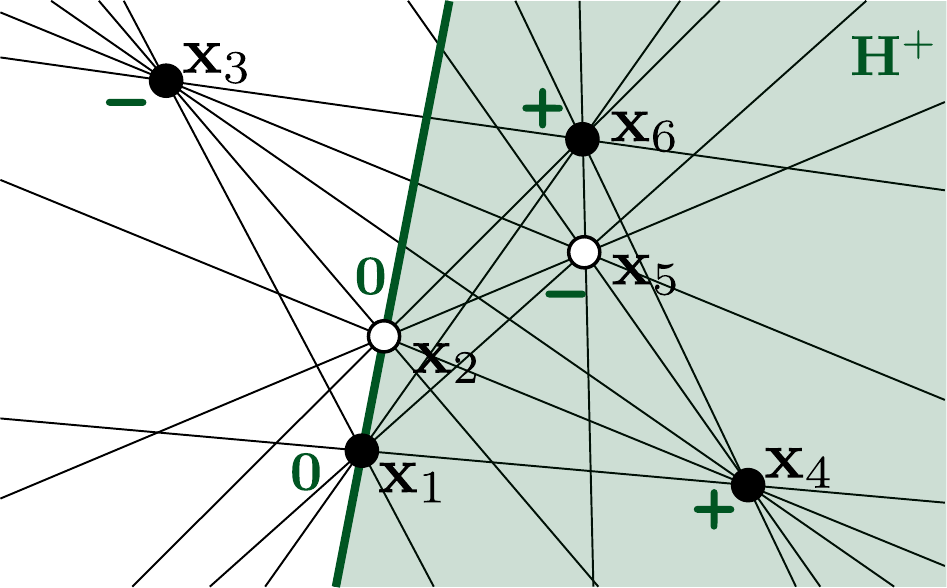}}}
\caption[Example of an affine Gale diagram]{The $3$-dimensional vector configuration $\vv V$ is the Gale dual of the vertex set of an hexagon. Its affine Gale diagram $\vv X$ in $\RR^2$, and an example of a covector $C$ with $C^+=\{\vv x_4,\vv x_3\}$ and $C^-=\{\vv x_3,\vv x_5\}$.}
\label{fig:affineGale}
\end{figure}
}

To visualize a vector configuration we will often use a tool called \defn{affine diagrams}\index{affine diagram}. For a fixed vector configuration $\vv V$ of rank $r$ (for simplicity we assume that $\veczero\notin\vv V$), an affine diagram for $\vv V$ is a 
``colored'' point configuration $\vv X$ in $\RR^{r-1}$. It is constructed using some vector $\vv c\in\RR^{n-d-1}$ such that $\sprod{\vv v_i}{\vv c}\neq 0$ for all $i$, 
and then associating each $\vv v_i\neq \veczero$ with the point $\vv x_i=\frac{\vv v_i}{\sprod{\vv v_i}{\vv c}}$. Then $\vv X$ is a point configuration in the hyperplane $\sprod{\vv x}{\vv c}=1$. We call $\vv x_i$ a \defn{positive point}\index{point!positive} if $\sprod{\vv v_i}{\vv c}>0$, and a \defn{negative point}\index{point!negative} if $\sprod{\vv v_i}{\vv c}<0$.  When $\vv V=\Gale{\vv A}$ we say that $\vv X$ is an \defn{affine Gale diagram} of~$\vv A$. In the concrete example of Figure~\ref{fig:affineGale},  positive points are depicted as full circles and negative points are empty. Observe how the cocircuits of~$\vv V$ can easily be obtained from the affine hyperplanes of $\vv X$, just by counting negative points  negatively. See the example of Figure~\ref{sfig:affineGalecovector} and compare it to Figure~\ref{fig:affinecovector}.

\section{Operations} \label{sec:operations}

\subsection{Deletion and contraction}\label{sec:delcontr}
Two handy operations on a vector configuration $\vv V$ are deletion and contraction of an element. 
The \defn{deletion}\index{deletion} $\vv V\setminus \vv v$\index{$\vv V\setminus \vv v$} of~$\vv v \in \vv V$ is the vector configuration $\vv V\setminus \{\vv v\}$, just like in the example of Figure~\ref{fig:exDeletion}.  
\iftoggle{bwprint}{%
\begin{figure}[htpb]
\centering
 \subbottom[$\vv V$]{\includegraphics[width=.25\textwidth]{Figures/ExampleGale_dual}}
\qquad\raisebox{1.5cm}{\includegraphics[width=.2\textwidth]{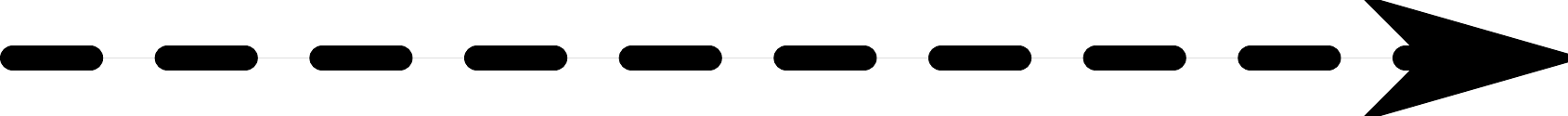}}\qquad
 \subbottom[$\vv V\setminus \vv v_1$]{\label{sfig:Vdeletion}\includegraphics[width=.25\textwidth]{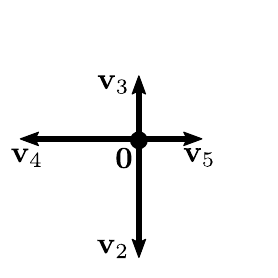}}
\caption{The deletion of $\vv v_1$ from $\vv V$. The labels of the elements of~$\vv V$ are in correspondence with those of $\vv V\setminus\vv v_1$.}
\label{fig:exDeletion}
\end{figure}}{%
\iftoggle{print}{%
\begin{figure}[htpb]
\centering
 \subbottom[$\vv V$]{\includegraphics[width=.25\textwidth]{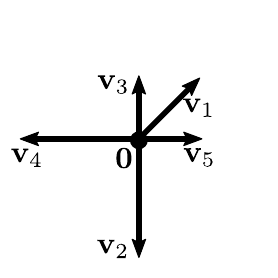}}
\qquad\raisebox{1.5cm}{\includegraphics[width=.2\textwidth]{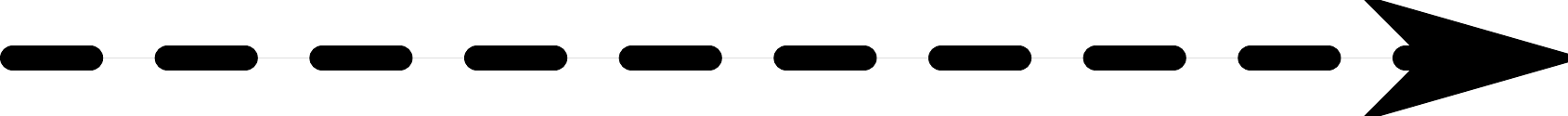}}\qquad
 \subbottom[$\vv V\setminus \vv v_1$]{\label{sfig:Vdeletion}\includegraphics[width=.25\textwidth]{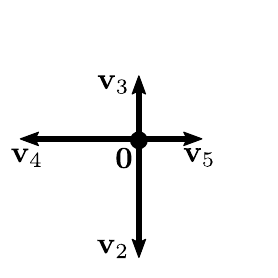}}
\caption{The deletion of $\vv v_1$ from $\vv V$. The labels of the elements of~$\vv V$ are in correspondence with those of $\vv V\setminus\vv v_1$.}
\label{fig:exDeletion}
\end{figure}}{%
\begin{figure}[htpb]
\centering
 \subbottom[$\vv V$]{\includegraphics[width=.25\textwidth]{Figures/ExampleGale_dual_col}}
\qquad\raisebox{1.5cm}{\includegraphics[width=.2\textwidth]{Figures/arrow}}\qquad
 \subbottom[$\vv V\setminus \vv v_1$]{\label{sfig:Vdeletion}\includegraphics[width=.25\textwidth]{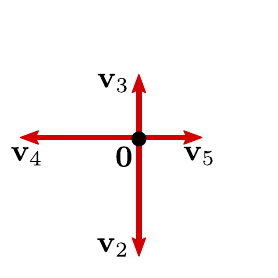}}
\caption{The deletion of $\vv v_1$ from $\vv V$. The labels of the elements of~$\vv V$ are in correspondence with those of $\vv V\setminus\vv v_1$.}
\label{fig:exDeletion}
\end{figure}
}
}

The \defn{contraction}\index{contraction} $\vv V/\vv v$\index{$\vv V/\vv v$} of a non-zero vector $\vv v \in\vv  V$ is given by projecting $\vv V$ parallel to $\vv v$ onto some linear hyperplane that does not contain $\vv v$, and then deleting~$\vv v$. For example, one can use the projection $\vv v_i\mapsto \tilde{\vv  v}_i:={\vv  v}_i -\frac{\sprod{{\vv  v}}{{\vv  v}_i}}{\sprod{{\vv  v}}{{\vv  v}}}{\vv  v}$, and then $\vv V/{\vv  v}=\set{\tilde {\vv  v}_i}{{\vv  v}_i\neq {\vv  v}}$; see Figure~\ref{fig:exContraction} for an example. The contraction of $\veczero$ is just its deletion.
\iftoggle{bwprint}{%
\begin{figure}[htpb]
\centering
 \subbottom[$\vv V$]{\includegraphics[width=.25\textwidth]{Figures/ExampleGale_dual}}
\raisebox{1.5cm}{\includegraphics[width=.075\textwidth]{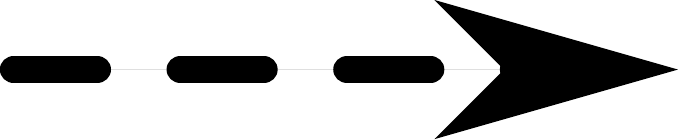}}
\includegraphics[width=.25\textwidth]{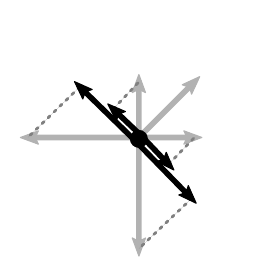}
\raisebox{1.5cm}{\includegraphics[width=.075\textwidth]{Figures/shortarrow}}
 \subbottom[$\vv V/\vv v_1$]{\label{sfig:Vcontraction}\includegraphics[width=.25\textwidth]{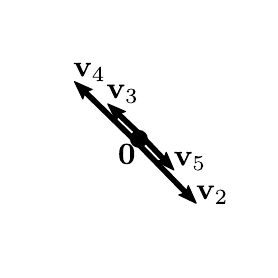}}
\caption{The contraction of $\vv v_1$ from $\vv V$. The labels of the elements of~$\vv V$ are in correspondence with those of $\vv V/\vv v_1$.}
\label{fig:exContraction}
\end{figure}
}{%
\iftoggle{print}{%
\begin{figure}[htpb]
\centering
 \subbottom[$\vv V$]{\includegraphics[width=.25\textwidth]{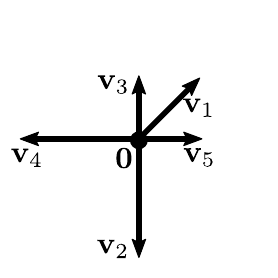}}
\raisebox{1.5cm}{\includegraphics[width=.075\textwidth]{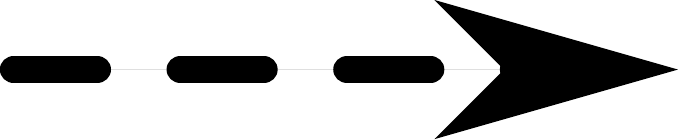}}
\includegraphics[width=.25\textwidth]{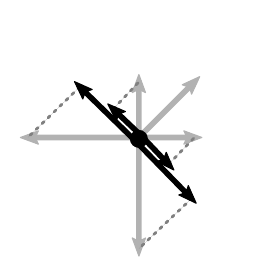}
\raisebox{1.5cm}{\includegraphics[width=.075\textwidth]{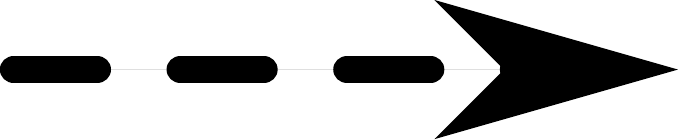}}
 \subbottom[$\vv V/\vv v_1$]{\label{sfig:Vcontraction}\includegraphics[width=.25\textwidth]{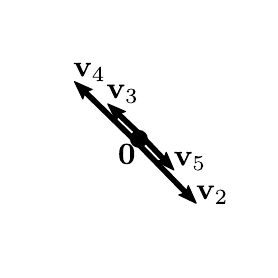}}
\caption{The contraction of $\vv v_1$ from $\vv V$. The labels of the elements of~$\vv V$ are in correspondence with those of $\vv V/\vv v_1$.}
\label{fig:exContraction}
\end{figure}
}{%
\begin{figure}[htpb]
\centering
 \subbottom[$\vv V$]{\includegraphics[width=.25\textwidth]{Figures/ExampleGale_dual_col}}
\raisebox{1.5cm}{\includegraphics[width=.075\textwidth]{Figures/shortarrow}}
\includegraphics[width=.25\textwidth]{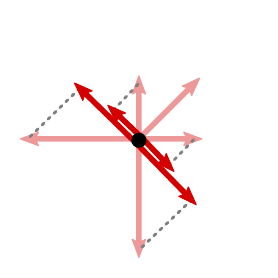}
\raisebox{1.5cm}{\includegraphics[width=.075\textwidth]{Figures/shortarrow}}
 \subbottom[$\vv V/\vv v_1$]{\label{sfig:Vcontraction}\includegraphics[width=.25\textwidth]{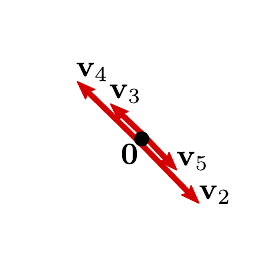}}
\caption{The contraction of $\vv v_1$ from $\vv V$. The labels of the elements of~$\vv V$ are in correspondence with those of $\vv V/\vv v_1$.}
\label{fig:exContraction}
\end{figure}
}
}

In terms of vectors of the oriented matroid $\cM(\vv V)$, we obtain
\begin{align*}
\ve(\vv V\setminus \vv v)&=\set{(U^+,U^-)}{(U^+,U^-)\in \ve(\vv V),\; \vv v \not\in \ul U},\\
\ve(\vv V\, /\, \vv v)&=\set{(U^+\setminus \{\vv v\},U^-\setminus \{\vv v\})}{(U^+,U^-)\in \ve(\vv V)},
\end{align*}
where the equalities of vectors of~$\vv V$ and vectors of~$\vv V/\vv v$ in the previous statement should be understood in the sense that their elements have the same corresponding indices. Indeed, we will often identify the elements of $\vv V/ \vv v$ with the elements of $\vv V\setminus \vv v$ using this natural bijection.

Contraction and deletion are dual operations ---\,  $\pGale{\vv V\setminus \vv v}=\Gale {\vv V}/\vv v$ and $\pGale{\vv V/ \vv v}=\Gale{\vv V}\setminus\vv v$\,--- that commute ---\, $\left(\vv V\setminus \vv v\right)/\vv w= \left(\vv V/\vv w\right)\setminus \vv v$ \,--- and naturally extend to subsets $\vv W\subseteq \vv V$ by iteratively deleting (resp. contracting) every element in $\vv W$. The rank of $\vv V/\vv W$ is \[\rank(\vv V/\vv W)=\rank(\vv V)-\rank(\vv W).\]

These operations can also be defined for point configurations $\vv A$ just by homogenizing them (see Figure~\ref{fig:affineContractionDeletion}). 
However if one wants to stick to the setting of point configurations, then only contractions of faces of~$\conv (\vv A)$ make sense, since otherwise the resulting vector configuration cannot be dehomogenized back into a point configuration. Geometrically, a realization of the contraction $\vv A/\vv a$ can be obtained by projecting $\vv A$ radially from $\vv a$ onto a hyperplane that separates $\vv a$ from $\vv A\setminus \vv a$. The reader can follow this procedure to get the representation of the contraction in Figure~\ref{sfig:Acontraction} from the configuration in Figure~\ref{sfig:Aoriginal}.

\iftoggle{bwprint}{%
\begin{figure}[htpb]
\centering
 \subbottom[$\vv A$]{\label{sfig:Aoriginal}\includegraphics[width=.25\textwidth]{Figures/ExampleGale_primal}}
 \qquad
 \subbottom[$\vv A\setminus \vv a_1$]{\label{sfig:Adeletion}\includegraphics[width=.25\textwidth]{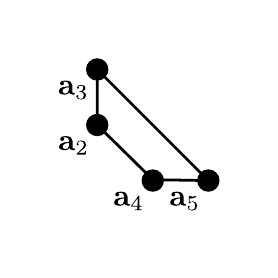}}\qquad
  \subbottom[$\vv A/ \vv a_1$]{\label{sfig:Acontraction}\includegraphics[width=.25\textwidth]{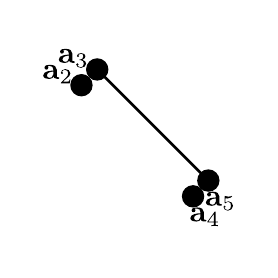}}
\caption{The contraction and deletion of $\vv a_1$ from $\vv A$.}
\label{fig:affineContractionDeletion}
\end{figure}
}{%
\iftoggle{print}{%
\begin{figure}[htpb]
\centering
 \subbottom[$\vv A$]{\label{sfig:Aoriginal}\includegraphics[width=.25\textwidth]{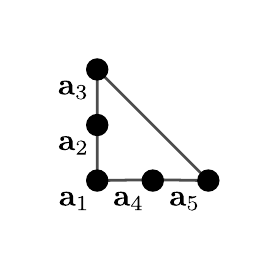}}
 \qquad
 \subbottom[$\vv A\setminus \vv a_1$]{\label{sfig:Adeletion}\includegraphics[width=.25\textwidth]{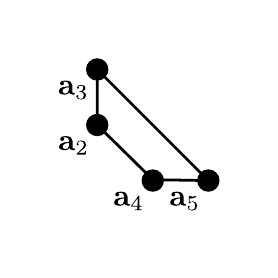}}\qquad
  \subbottom[$\vv A/ \vv a_1$]{\label{sfig:Acontraction}\includegraphics[width=.25\textwidth]{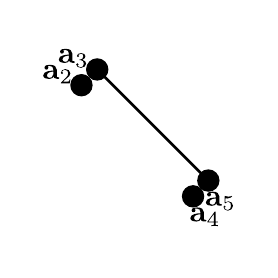}}
\caption{The contraction and deletion of $\vv a_1$ from $\vv A$.}
\label{fig:affineContractionDeletion}
\end{figure}
}{%
\begin{figure}[htpb]
\centering
 \subbottom[$\vv A$]{\label{sfig:Aoriginal}\includegraphics[width=.25\textwidth]{Figures/ExampleGale_primal_col}}
 \qquad
 \subbottom[$\vv A\setminus \vv a_1$]{\label{sfig:Adeletion}\includegraphics[width=.25\textwidth]{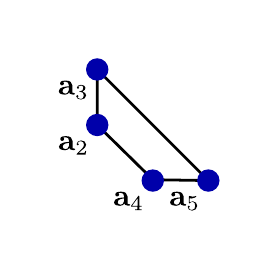}}\qquad
  \subbottom[$\vv A/ \vv a_1$]{\label{sfig:Acontraction}\includegraphics[width=.25\textwidth]{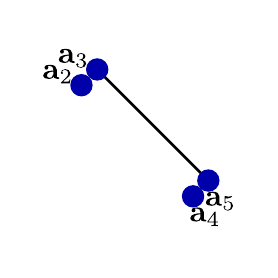}}
\caption{The contraction and deletion of $\vv a_1$ from $\vv A$.}
\label{fig:affineContractionDeletion}
\end{figure}
}
}

When $\vv A$ is the vertex set of a polytope $\vv P=\conv (\vv A)$ and $\vv a\in \vv A$, we say that $\conv (\vv A\setminus \vv a)$ is a \defn{subpolytope}\index{subpolytope} of $\vv A$ and that $\conv(\vv A/\vv a)$ is the \defn{vertex figure} of $\vv a$ in $\vv P$.
Observe that if $\vv F$ is a face of $\vv P$, then the contraction $\vv P/\vv F$ is well defined, and it is a polytope of dimension 
\[\dim(\vv P/\vv F)=\dim(\vv P)-\dim(\vv F)-1.\]\\

\subsection{Projections}\label{sec:proj}
A \defn{(linear) projection}\index{projection!linear} of a vector configuration $\vv V\subset\RR^r$ is its image $\pi(\vv V)\subset\RR^s$ under a linear projection $\pi:\RR^r\rightarrow\RR^s$. Every vector of~$\cM(\vv V)$ is also a vector of $\cM(\pi(\vv V))$, so that $\ve(\vv V)\subset\ve(\pi(\vv V))$, and every covector of $\cM(\pi(\vv V))$ is also a covector of~$\cM(\vv V)$, so that $\cov(\pi(\vv V))\subset\cov(\vv V)$. Besides these relations, the combinatorics of $\pi(\vv V)$ depend strongly on the projection~$\pi$. Its relation to duality is explained in the following lemma.
\begin{lemma}\label{lem:projectionduality}
 If there is a projection $\pi:\RR^r\to\RR^s$ such that $\pi(\vv V)=\vv W$, then there is a projection $\Gale \pi:\RR^{n-s}\rightarrow\RR^{n-r}$ such that $\Gale\pi(\Gale{\vv W})=\Gale{\vv V}$.
\end{lemma}

An \defn{(affine) projection}\index{projection!affine} of a $d$-dimensional point configuration $\vv A$ is the dehomogenization of $\pi(\hom(\vv A))$, where $\pi$ is a linear projection $\pi:\RR^{d+1}\to\RR^{e+1}$ that sends the hyperplane $\set{\vv v\in\RR^{d+1}}{\vvc v_{d+1}=1}$ to the hyperplane $\set{\vv w\in\RR^{e+1}}{\vvc w_{e+1}=1}$.

\subsection{Direct sum, join and pyramid}\label{sec:directsum}
The \defn{direct sum}\index{direct sum} of two vector configurations $\vv V\subset \RR^r$ and $\vv W\subset \RR^s$ is the vector configuration\index{$\vv V\oplus \vv W$}
\[
   \vv V\oplus \vv W 
   \ = \ 
   \set{(\vv v,\veczero)}{\vv v\in \vv V}\cup\set{(\veczero,\vv w)}{\vv w\in \vv W} 
\ \subset \ \RR^{r+s},
\]
see Figure~\ref{fig:directSum} for an example. It is easy to see that $\ci(\vv V\oplus \vv W)=\ci(\vv V)\cup \ci(\vv W )$ and $\co(\vv V\oplus \vv W)=\co(\vv V)\cup \co(\vv W )$, see \cite[Proposition 7.6.1]{OrientedMatroids1993}. 
This shows that:
\begin{lemma}\label{lem:directsumdual}
$\pGale{\vv V\oplus \vv W }=\Gale{\vv V}\oplus \Gale{\vv W }$. 
\end{lemma}

\iftoggle{bwprint}{%
\begin{figure}[htpb]
\centering
 \subbottom[$\vv V$]{\label{sfig:freesumV}\includegraphics[scale=1.2]{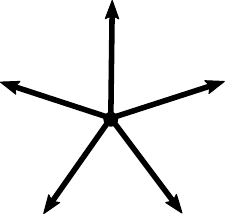}}\qquad\qquad
 \subbottom[$\vv W$]{\label{sfig:freesumW}\quad\qquad\includegraphics[scale=1.2]{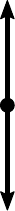}\qquad}\qquad\qquad
 \subbottom[$\vv V\oplus\vv W$]{\label{sfig:freesumVW}\includegraphics[scale=.9]{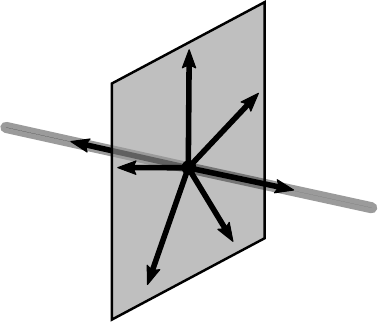}}
\caption[Example of a direct sum.]{Two vector configurations $\vv V$ and $\vv W$ of respective ranks $2$ and $1$, and their direct sum $\vv V\oplus\vv W$, of rank $3$.}
\label{fig:directSum}
\end{figure}
}{%
\begin{figure}[htpb]
\centering
 \subbottom[$\vv V$]{\label{sfig:freesumV}\includegraphics[scale=1.2]{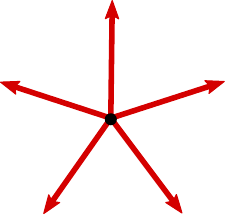}}\qquad\qquad
 \subbottom[$\vv W$]{\label{sfig:freesumW}\quad\qquad\includegraphics[scale=1.2]{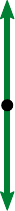}\qquad}\qquad\qquad
 \subbottom[$\vv V\oplus\vv W$]{\label{sfig:freesumVW}\includegraphics[scale=.9]{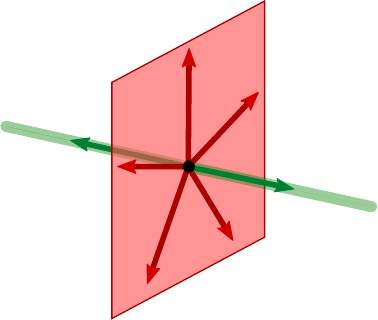}}
\caption[Example of a direct sum.]{Two vector configurations $\vv V$ and $\vv W$ of respective ranks $2$ and $1$, and their direct sum $\vv V\oplus\vv W$, of rank $3$.}
\label{fig:directSum}
\end{figure}
}

When $\vv A\subset \RR^d$ and $\vv B\subset \RR^e$ are point configurations, the dehomogenization of $\hom(\vv A)\oplus \hom(\vv W)$ is the \defn{join}\index{join} of~$\vv A$ and $\vv B$. It is denoted $\vv A\join \vv B$\index{$\vv A\join \vv B$}, and is obtained by embedding the configurations into skew affine subspaces of $\RR^{d+e+1}$. For example,
\[
   \vv A\join \vv B= \set{(\vv a,\veczero,1)\in}{\vv a\in \vv A}\cup\set{(\veczero,\vv b,-1)\in}{\vv b\in \vv B}.
\]
Observe that the dimension of $\vv A\join \vv B$ is \[\dim(\vv A\join \vv B)=\dim(\vv A)+\dim(\vv B)+1.\]

If $\vv B$ is a single point $\{\vv b\}$, then $\vv A\join\{\vv b\}$ is called a \defn{pyramid}\index{pyramid} over $\vv A$ with apex $\vv b$. More generally, if $\vv B$ is the vertex set of a $(k-1)$-simplex, then $\vv A\join\vv B$ is a \defn{$k$-fold pyramid} over $\vv A$. Since the Gale dual of a simplex is $0$-dimensional, a vector configuration $\vv V$ is the Gale dual of a pyramid if and only if $\veczero\in \vv V$.

\section{Oriented Matroids}\label{sec:orientedmatroids}

Given two vectors $U$ and $V$ we define their \defn{composition}\index{composition} $U\circ V$\index{$U\circ V$} to be the vector $W$ with $W^+=U^+\cup(V^+\setminus \ul U^-)$ and $W^-=U^-\cup(V^-\setminus \ul U^+)$. 

In the geometric setting, the composition of covectors is very easy to illustrate. For a fixed vector configuration~$\vv V$, let $\vvh H_1$ and $\vvh H_2$ be linear hyperplanes with normal vectors $\vv v_1$ and $\vv v_2$. We define their composition {$\vvh H_1\circ \vvh H_2$} (with respect to $\vv V$) to be a hyperplane~${\vvh H}$ with normal vector $\vv v_1+\delta \vv v_2$ for some very small $\delta$ (whose precise value depends on $\vv V$).
Let $\vv v\in \vv V$. If $\delta$ is small enough, then $\vv v\in {\vvh H}$ if and only if $\vv v\in {\vvh H}_1\cap {\vvh H}_2$, and $\vv v\in {\vvh H}^\pm$ if and only if either $\vv v\in {\vvh H}_1^\pm$ or $\vv v\in {\vvh H}_1$ and $\vv v\in {\vvh H}_2^\pm$. This procedure can be carried out analogously for affine configurations, as in Figure~\ref{fig:composition}.\\

\iftoggle{bwprint}{%
\begin{figure}[htpb]
\centering
 \subbottom[$\vvh H_1$ and $\vvh H_2$]{\includegraphics[width=.45\linewidth]{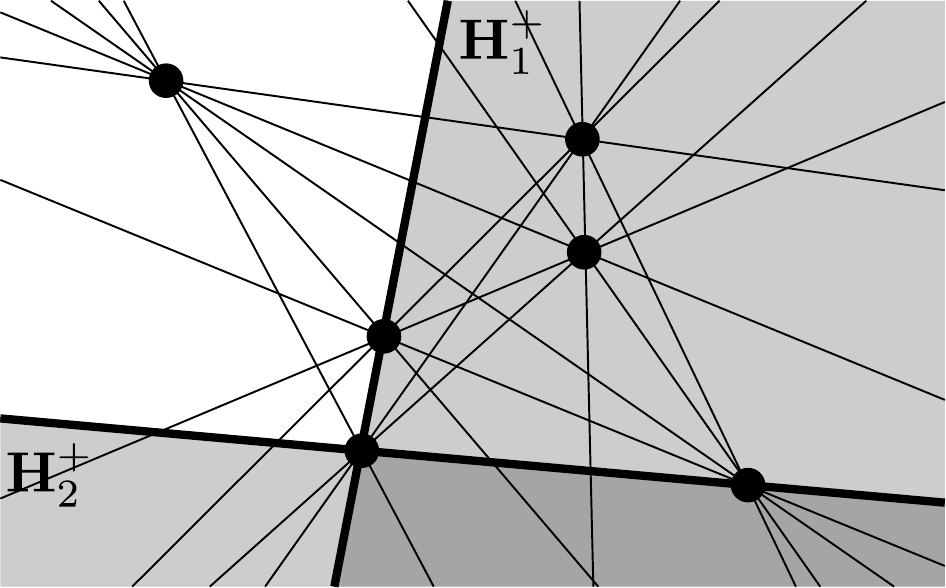}}\quad
 \subbottom[$\vvh H_1\circ \vvh H_2$]{\includegraphics[width=.45\linewidth]{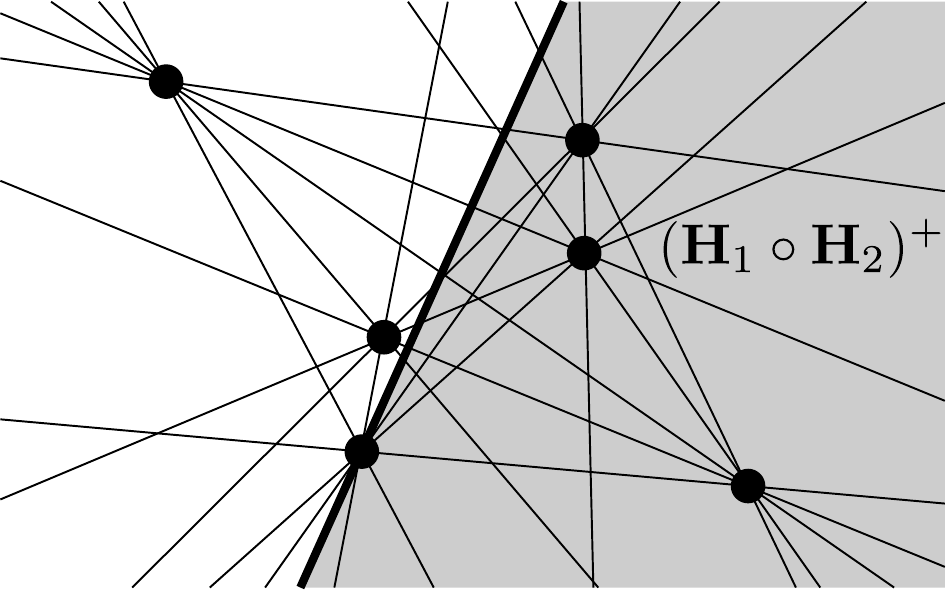}}
 \caption{Example of a composition of affine hyperplanes.}
\label{fig:composition}
\end{figure}
}{%
\begin{figure}[htpb]
\centering
 \subbottom[$\vvh H_1$ and $\vvh H_2$]{\includegraphics[width=.45\linewidth]{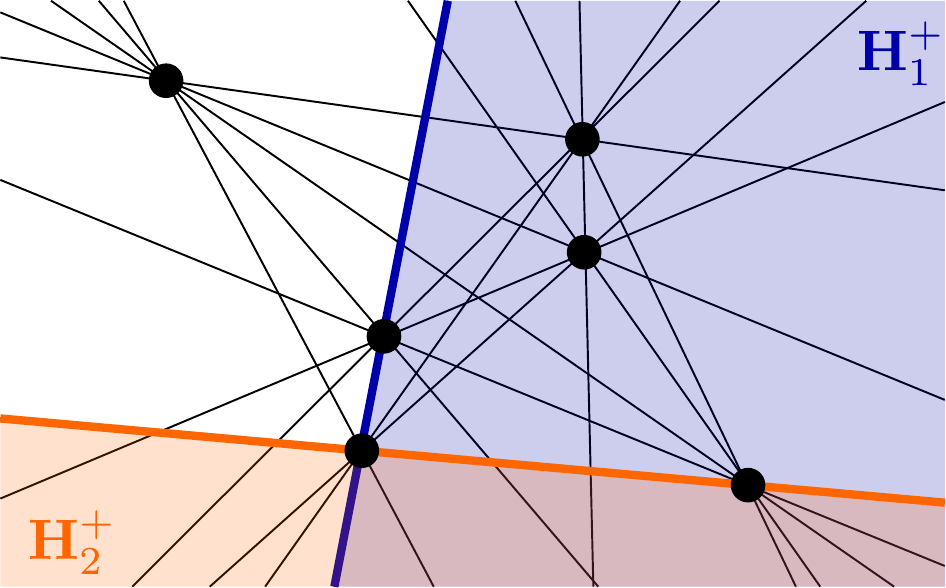}}\quad
 \subbottom[$\vvh H_1\circ \vvh H_2$]{\includegraphics[width=.45\linewidth]{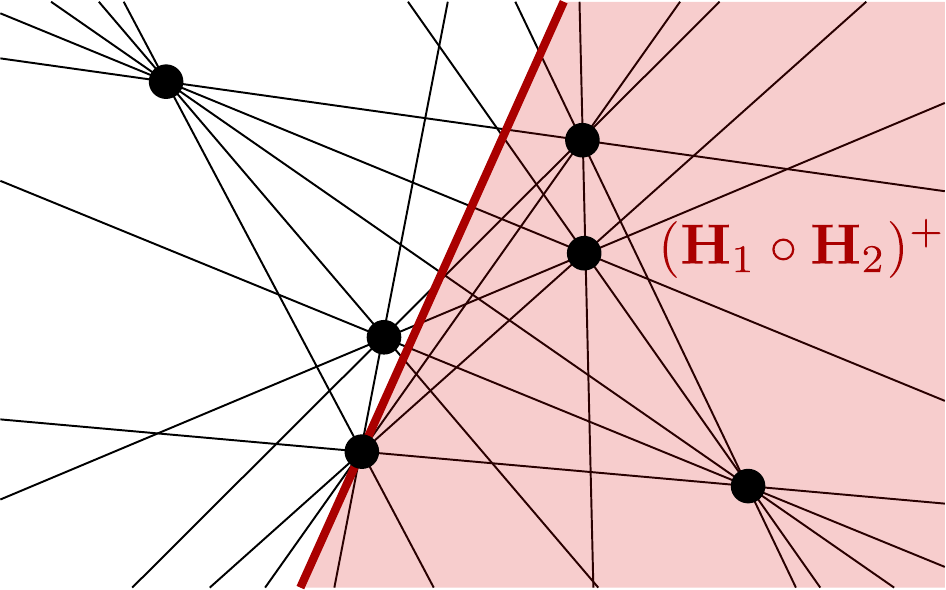}}
 \caption{Example of a composition of affine hyperplanes.}
\label{fig:composition}
\end{figure}
}

With this operation we provide a definition of an \defn{oriented matroid} in terms of vectors. A set $\ve$ of signed subsets of $E$ is the set of vectors of an oriented matroid $\cM$ if it fulfills the following \defn{vector axioms}, where $-U$ is the signed set defined by $(-U)^+=U^-$ and $(-U)^-=U^+$. 

\begin{description}
 \item[(V0)] $\emptyset\in \ve$;
 \item[(V1)] (\defn{vector symmetry}\index{vector!symmetry}) if $U\in\ve$, then $-U\in\ve$;
 \item[(V2)] (\defn{vector composition}\index{vector!composition}) for all $U,V\in \ve$, $U\circ V\in \ve$;
 \item[(V3)] (\defn{vector elimination}\index{vector!elimination}) for all $U,V\in \ve$ and $e\in E$ with $U(e)V(e)<0$, 
there is a $W\in \ve$ such that $W(e)=0$ and $W(f)=(U\circ V)(f)$ for all $f\in E$ with $U(f)V(f)\geq 0$.
\end{description}

For instance, the set of vectors $\ve(\vv V)$ of every vector configuration $\vv V$ fulfills these axioms. Indeed, for \textbf{(V0)} observe that $\veczero$ induces the vector~$\emptyset$, and for \textbf{(V1)} observe that if $\vv \gl$ is a linear dependence, so is $-\vv\gl$. We have seen that \textbf{(V2)}~corresponds to $\vv \gl +\delta \vv \gm$ for some small~$\delta$. Finally, \textbf{(V3)}~corresponds to linear combinations of dependences of the form $\gm_i\vv \gl -\gl_j \vv\gm$.

However, the reciprocal is not true. There are sets of vectors $\ve$ that fulfill \textbf{(V0)}, \textbf{(V1)}, \textbf{(V2)} and \textbf{(V3)} but are not the set of vectors of any vector configuration. The oriented matroids that arise from these sets of vectors are called \defn{non-realizable}\index{oriented matroid!realizable}. They also have many geometric interpretations, for example in terms of \defn{pseudohyperplane arrangements} or \defn{zonotopal tilings}, which we omit here but encourage the reader to discover~\cite{OrientedMatroids1993}.

\iftoggle{bwprint}{%
\begin{example}[The Non-Pappus oriented matroid]
 \begin{figure} [htpb]
\centering
\includegraphics[width=0.5\linewidth]{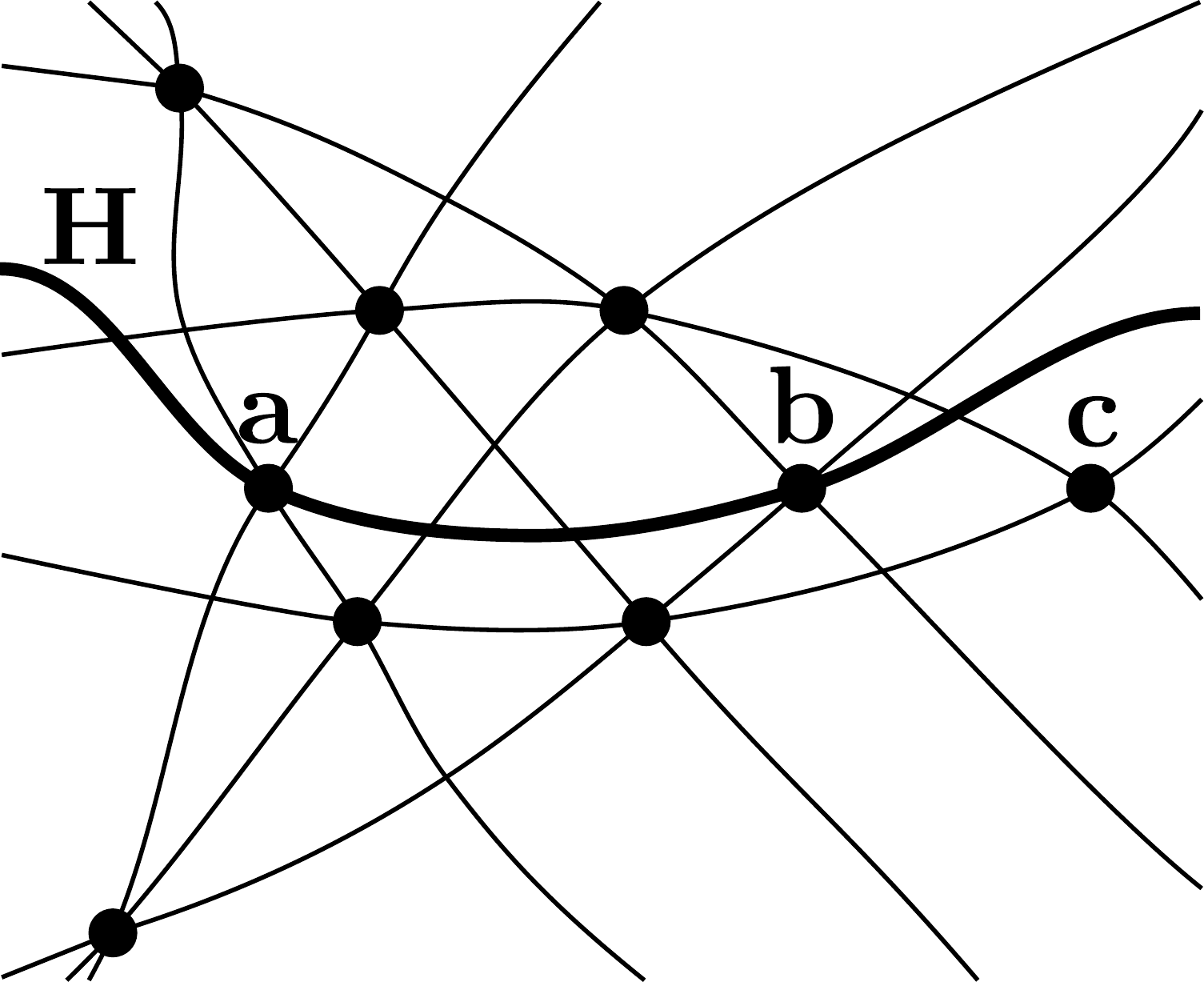}
\caption{A non-realizable oriented matroid.}
\label{fig:nonpappus}
 \end{figure}
}{%
\begin{example}[The Non-Pappus oriented matroid]
 \begin{figure} [htpb]
\centering
\includegraphics[width=0.5\linewidth]{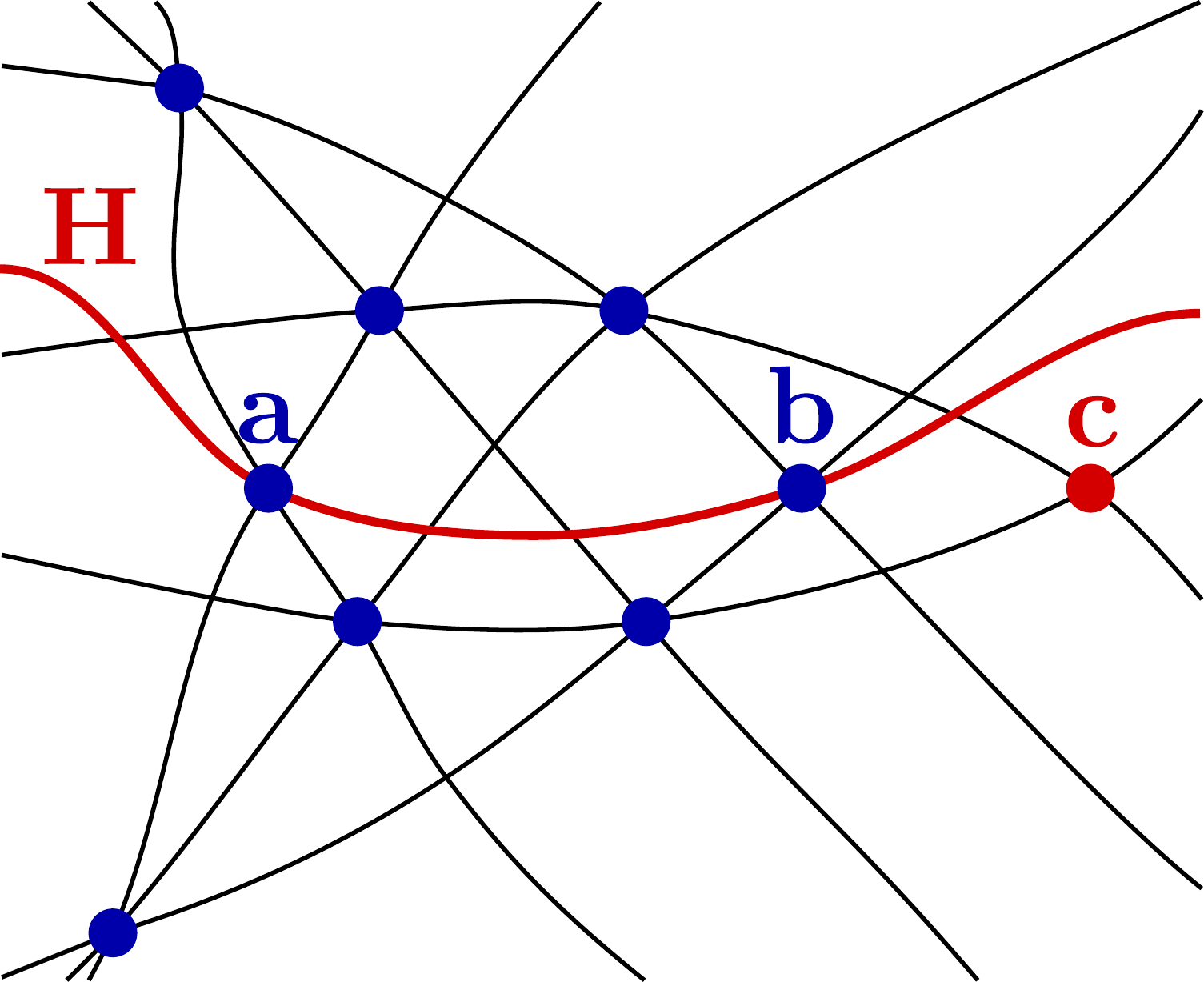}
\caption{A non-realizable oriented matroid.}
\label{fig:nonpappus}
 \end{figure}
}

Consider the point configuration $\vv A$ depicted in Figure~\ref{fig:nonpappus}. We have shown some pseudohyperplanes that determine signed subsets of $\vv A$ just as in Figure~\ref{fig:affinecovector}. These turn out to be the covectors of an oriented matroid. However, this matroid cannot be realizable, since the classical Pappus theorem states that in any affine or projective representation of this point configuration, the points $\vv a$,~$\vv b$ and~$\vv c$ must be collinear, and hence $\vv c$ must belong to ${\vvh H}$.
\end{example}

The \defn{rank}\index{rank} of the oriented matroid $\cM$ with set of vectors $\ve$ is $n-(r+1)$, where~$r$ is the largest length of a chain of vectors $\veczero<U^1<U^2<\dots <U^r$ with $U^i\in \ve$. 
If $\vv V$ is a $\rd$-dimensional vector configuration, $\rank(\cM(\vv V))=\rd$.\\

Two signed sets $C$ and $X$ are \defn{orthogonal}\index{orthogonal}, which we denote $C\perp X$, if either $\ul C\cap \ul X=\emptyset$ or there exist $e,f\in \ul C\cap \ul X$ such that $C(e)X(e)=-C(f)X(f)$. Given a matroid $\cM$ on a ground set $E$ with set of circuits $\ci$, the set
\[\cC^\perp:=\set{C\text{ signed subset of }E}{C\perp X\text{ for all }X\in \ci}\]
is the set of covectors $\cov$ of $\cM$. Moreover, it also fulfills properties \textbf{(V0)}, \textbf{(V1)}, \textbf{(V2)} and \textbf{(V3)}, which means that it is the set of vectors of another matroid $\Gale\cM$, which is called the \defn{dual matroid}\index{oriented matroid!duality} of $\cM$ and has rank $|E|-\rank(\cM)$, which is called the \defn{corank}\index{corank} of $\cM$.
For a realizable matroid $\cM(\vv V)$, its dual matroid is the matroid of the Gale dual, $\pGale{\cM(\vv V)}=\cM(\Gale{\vv V})$.

A matroid $\cM$ of rank $r$ is \defn{uniform} if the underlying matroid $\underline \cM$ is uniform. That is, for every subset $S$ of cardinality $r-1$, there a cocircuit~$C$ with~$C^0=S$, and for every subset $T$ of cardinality $r+1$, there is a circuit~$X$ with $\ul X=T$. A vector configuration $\vv V$ is in \defn{general position}\index{general position} if and only if $\cM(\vv V)$~is uniform.

We say that two oriented matroids $\cM_1$ and $\cM_2$ on respective ground sets $E_1$ and $E_2$ are 
\defn{isomorphic}\index{isomorphic}, which we denote $\cM_1\simeq\cM_2$, when there is a bijection between $E_1$ and $E_2$ that sends circuits of $\cM_1$ to circuits of $\cM_2$ (and equivalently for vectors, cocircuits or covectors). For instance, the point configurations in Figures~\ref{sfig:4pointsR2_1}, \ref{sfig:4pointsR2_2} and \ref{sfig:4pointsR2_3} define isomorphic oriented matroids, which are not isomorphic to the oriented matroid of Figure~\ref{sfig:4pointsR2_4}.\\

When $\cM$ is the oriented matroid of the homogenization of a point configuration $\vv A$, it is said to be \defn{acyclic}\index{oriented matroid!acyclic}. Put differently, the whole ground set is a positive covector. We can define the face lattice of an acyclic matroid that generalizes that of point configurations. Its \defn{facets}\index{facet} are the complements of the supports of its positive cocircuits, and its \defn{faces}\index{face} the complements of its positive covectors. 

Similarly, a matroid is \defn{totally cyclic}\index{oriented matroid!totally cyclic} if it contains the whole ground set 
 as a positive vector. We use the notation $\stc{r}$\index{$\stc{r}$} for the only totally cyclic oriented matroid of rank $r$ with~$r+1$ elements, which is $\cM(\vv V)$ for the vector configuration $\vv V=\{\vv e_1, \dots, \vv e_r, -\sum_{i=1}^r \vv e_i\}$.\\

Of course, the operations presented in Section~\ref{sec:operations} above can also be applied to non-realizable matroids, just by extending their definition in terms of vectors.

\iftoggle{bwprint}{%
\renewcommand{\partfigure}{prelude/Figures/fullsewingexamples}
}{%
\renewcommand{\partfigure}{prelude/Figures/fullsewingexamples_col}
}
\renewcommand{\namepart}{Neighborly}
\part{Neighborly}\label{partI}
\chapter{Introduction}
\label{ch:intro_neigh}
\section{Overview}
Consider the following classical problem in the theory of polytopes. 

\begin{question}[Upper bound problem]
What is the maximal number of $i$-dimensional faces that a $d$-dimensional polytope on $n$ vertices can have?
\end{question}
A first upper bound for $f_i(\vv P)$, the number of $i$-dimensional faces of $\vv P$, in terms of $f_0(\vv P)$, its number of vertices,  is very easy to obtain. 
Since every $i$-face is defined by its set of vertices, which must be of size at least $i+1$, we have the trivial bound \begin{equation}\label{eq:trivialupperbound} f_i(\vv P)\leq \binom{f_0(\vv P)}{i+1}.\end{equation} 

This leads to the question of whether this bound is tight, and in case it is, which polytopes attain it. If \eqref{eq:trivialupperbound} is an equality, we say that $\vv P$ is $(i+1)$-neighborly:

\begin{definition}\label{def:neighborly}
A polytope $\vv P$ is \defn{$k$-neighborly} if every subset of at most~$k$ vertices of $\vv P$ is the set of vertices of a face of $\vv P$.\index{neighborly}
\end{definition}

So the question is now: is there some $k$ such that a $k$-neighborly polytope exists? The following well-known result, see for example~\cite[Chapter~7]{GruenbaumEtal2003}, shows that we will not be able to find non-trivial $k$-neighborly polytopes on $n$ vertices for all values of $k$.

\begin{theorem}\label{thm:kneighd/2}
If a $d$-polytope $\vv P$ is $k$-neighborly for any $k>\ffloor{d}{2}$, then $\vv P$ must be the $d$-dimensional simplex $\Delta_d$.
\end{theorem}

This theorem motivates the definition of a $d$-polytope as \defn{neighborly} if it is $\ffloor{d}{2}$-neighborly. The next step is to show that there exist neighborly polytopes (\ie $k$-neighborly for every $k\leq \ffloor{d}{2}$) different from the simplex.\\

Any set of $n$ points on the \defn{moment curve} in $\RR^d$, $\{(t,t^2,\dots,t^d):t\in \RR\}$, is the set of vertices of a neighborly polytope.
Since it turns out that 
the combinatorial type of this polytope does not depend on the particular choice of points (see~\cite[Section~4.7]{GruenbaumEtal2003}), we denote it as $\cyc{d}{n}$\index{$\cyc{d}{n}$}, the \defn{cyclic polytope} with $n$ vertices in $\RR^d$. 
Cyclic polytopes\index{polytope!cyclic} have been discovered and rediscovered several times (cf.~\cite[Section~7.4]{GruenbaumEtal2003}) among others by Carath\'eodory~\cite{Caratheodory1911}, Gale~\cite{Gale1956,Gale1963} and Motzkin~\cite{Motzkin1957}.\\

In his abstract~\cite{Motzkin1957}, Motzkin claimed that for all $i$, cyclic polytopes $\cyc{d}{n}$ had the maximal number of $i$-faces among $d$-polytopes with $n$ vertices, not only when $i+1\leq \ffloor{d}{2}$, and that cyclic polytopes were unique with this property (in even dimension). The first claim was proven by McMullen in 1970~\cite{McMullen1970}.

\begin{theorem}[Upper Bound Theorem~\cite{McMullen1970}]
 For any $d$-polytope~$\vv P$ with $n$ vertices, and for any $1\leq i \leq d$,
\begin{equation}\label{eq:ubt}
 f_{i}(\vv P)\leq f_i(\cyc{d}{n}).
\end{equation}
Moreover, equality in~\eqref{eq:ubt} for some $i$ with $\ffloor{d}{2}\leq i+1\leq d$ implies that~$\vv P$ is neighborly; and the equality holds for some $i$ with $\fceil{d}{2}\leq i+1\leq d$ if and only if $\vv P$ is simplicial and neighborly.
\end{theorem}

However, the second part of Motzkin's claim is (very) false, as we will soon see. The crucial observation is that all simplicial neighborly $d$-polytopes with $n$ vertices have the same $f$-vector as $\cyc{d}{n}$, and hence are examples of polytopes where~\eqref{eq:ubt} is an equality for all $1\leq i \leq d$.

\begin{remark}\label{rmk:simplicialneighborly}
While it is easy to see that all even-dimensional neighborly polytopes are simplicial using Theorem~\ref{thm:kneighd/2}, there are non-simplicial neighborly polytopes in odd dimension. These do not attain equality in~\eqref{eq:ubt} when $i+1\geq\fceil{d}{2}$. For example, any pyramid over an even dimensional neighborly polytope is neighborly, and (when it is not a simplex) it is non-simplicial.
\end{remark}

\subsection{Many neighborly polytopes}

The first examples of non-cyclic neighborly polytopes were found in 1967 by Gr\"unbaum, who proved in his classical book~\cite[Section 7.2]{GruenbaumEtal2003} that for all $k\geq 2$ there are non-cyclic neighborly $(2k)$-polytopes with $2k+4$ vertices.

In 1981, Barnette introduced the \defn{facet splitting} technique~\cite{Barnette1981}, that allowed him to construct infinitely many neighborly polytopes, and to prove that the number of (combinatorial types of) neighborly $d$-polytopes with $n$ vertices, \defn{$\nnei{n}{d}$}, is bigger than \[\nnei{n}{d}\geq\frac{(2n-4)!}{n!(n-2)!\binom{n}{d-3}}\sim4^n.\]
  
This bound was improved by Shemer in~\cite{Shemer1982}, where he introduced the \defn{Sewing construction} to build an infinite family of neighborly polytopes in any even dimension. 
Given a neighborly $d$-polytope with $n$ vertices and a suitable flag of faces, one can ``sew'' a new vertex onto it to get a new neighborly $d$-polytope with $n+1$ vertices. With this construction, Shemer proved that
$\nnei{n}{d}$
is greater than \[\nnei{n}{d}\geq\frac{1}{2}\left(\left(\frac{d}{2}-1\right)\ffloor{n-2}{d+1}\right)!\sim n^{c_d n},
\] where $c_d\rightarrow \frac{1}{2}$ when $d\rightarrow\infty$.
\\

One of our main results is the following theorem, proved in Chapter~\ref{ch:counting}, that provides a new lower bound for \defn{$\lnei{n}{d}$}\index{$\lnei{n}{d}$}, the number of vertex-labeled combinatorial types of neighborly polytopes with $n$ vertices and dimension~$d$.

\medskip
\noindent\textbf{Theorem \ref{thm:lblnei}}
\emph{The number of labeled neighborly polytopes in even dimension $d$ with $r+d+1$ vertices fulfills}
\begin{equation}\label{eq:thebound}
\lnei{r+d+1}{d}\geq \frac{\left( r+d \right) ^{\left( \frac{r}{2}+\frac{d}{2} \right) ^{2}}}{{r}^{{(\frac{r}{2})}^{2}}{d}^{{(\frac{d}{2})}^{2}}{{\e}^{3\frac{r}{2}\frac{d}{2}}}}. \end{equation}

\medskip
This bound is always greater than 
\[
    \lnei{n}{d}\geq  \left( \frac{n-1}{\e^{3/2}}\right)^{\frac12 (n-d-1)d},
\]
and dividing by $n!$  easily shows this to improve Shemer's bound also in the unlabeled case. Moreover, since pyramids over even-dimensional neighborly polytopes are also neighborly, when $d$ is odd we can use the bound $\lnei{r+d+1}{d}\geq \lnei{r+d}{d-1}$.

Of course, \eqref{eq:thebound} is also a lower bound for \defn{$\lpol{n}{d}$}\index{$\lpol{n}{d}$}, the number of combinatorial types of vertex-labeled $d$-polytopes with $n$ vertices, and even improves 
\[\lpol{n}{d}\geq \left(\frac{n-d}{d}\right)^{\frac{nd}{4}},\]
the current best bound (valid only for $n\geq 2d$), which was given by Alon in 1986~\cite{Alon1986}.
\begin{remark}
To the best of our knowledge, the only known upper bounds for $\lnei{n}{d}$ are the bounds for $\lpol{n}{d}$. Goodman and Pollack proved in~\cite{GoodmanPollack1986} that 
$\lpol{n}{d}\leq n^{d(d+1)n}$, a bound that was improved by Alon in~\cite{Alon1986} to 
\[\lpol{n}{d}\leq \left(\frac{n}{d} \right)^{d^2n(1+o(1))}\text{ when }\tfrac{n}{d}\rightarrow \infty.\]
\end{remark}

\subsection{Constructing neighborly polytopes}

One of the main ingredients for our new results --- in particular for the bound~\eqref{eq:thebound} --- is a new construction for neighborly polytopes. 

We have already mentioned Barnette's facet splitting~\cite{Barnette1981} and Shemer's sewing construction~\cite{Shemer1982}. In both construction techniques, a new vertex is added to an existing neighborly polytope $\vv P$ using a so called \defn{universal flag} of~$\vv P$.

In Section~\ref{sec:shemer}, we show that Shemer's sewing construction can be very transparently explained (and generalized) in terms of \defn{lexicographic extensions} of oriented matroids. In fact, the same framework also explains Lee \& Menzel's related construction of A-sewing for non-simplicial polytopes~\cite{LeeMenzel2010} (Proposition~\ref{prop:Asewing}), and the results in~\cite{TrelfordVigh2011} on faces of sewn polytopes.

Next, we introduce two new construction techniques for polytopes. The first, \defn{Extended Sewing} (Construction~\ref{constr:cE}) is based on our Extended Sewing Theorem~\ref{thm:extshemersewing}. 
It is a generalization of Shemer's sewing to oriented matroids that is valid for any rank, just like Bistriczky's version~\cite{Bisztriczky2000}, and works for a large family of flags of faces, including the ones obtained by Barnette's facet splitting~\cite{Barnette1981}. 
Moreover, Extended Sewing is optimal in the sense that in odd ranks, the flags of faces constructed in this way are the only ones that yield neighborly polytopes (Proposition~\ref{prop:uniqueflags}).

Our second (and most important) new technique is \defn{Gale Sewing} (Construction~\ref{constr:cG}), whose key ingredient is the Double Extension Theorem~\ref{thm:thethm}. It lexicographically extends {duals} of neighborly polytopes and oriented matroids. 
Equivalently, given a neighborly matroid of rank $d$ with $n$ elements, we obtain new neighborly matroids of rank $d+1$ with $n+1$ vertices. 
The bound \eqref{eq:thebound} is obtained by estimating the number of polytopes in $\cG$, the family of all polytopes that can be constructed via Gale Sewing. This family contains all the neighborly polytopes constructed by Devyatov in~\cite{Devyatov2011}, which arise as a special case of Gale Sewing on polytopes of corank~$3$.
\begin{remark}
Using Extended Sewing, we construct three families of neighborly polytopes --- $\cS$, $\cE$ and $\cO$ --- the largest of which is $\cO$. 
In Section~\ref{sec:comparing}, we will see that $\cO\subseteq \cG$, and in this sense, Gale Sewing is a generalization of Extended Sewing. 
However, it is not true that the Double Extension Theorem~\ref{thm:thethm} generalizes the Extended Sewing Theorem~\ref{thm:extshemersewing}, because to build the polytopes in $\cO$ we restrict ourselves to certain universal flags that we find with Proposition~\ref{prop:extnewunifaces} (compare Remark~\ref{rmk:doesnotgeneralize}).
\end{remark}

Both techniques are based on the concept of \defn{lexicographic extension} of an oriented matroid, which turns out to be a very useful tool for constructing polytopes. In Section~\ref{sec:le} we introduce this concept and some of its properties.

\subsection{Neighborly oriented matroids}

We said before that our constructions worked not only with neighborly polytopes but in the more general setting of neighborly oriented matroids. Indeed, neighborliness is a purely combinatorial concept that can be easily defined in terms of oriented matroids. 
\begin{definition}
An oriented matroid $\cM$ of rank $\rd$ on a ground set~$E$ is \defn{neighborly} if
every subset $S\subset E$ of size at most $\ffloor{\rd-1}{2}$ is a face of~$\cM$. That is, if there exists a covector $C$ of $\cM$ with $C^+=E\setminus S$ and $C^-=\emptyset$.\index{neighborly!oriented matroid} 
\end{definition}

With this definition, realizable neighborly oriented matroids correspond to neighborly polytopes. However, not all neighborly oriented matroids are realizable. 
For example, the sphere ``$\cM^{10}_{425}$'' from Altshuler's list~\cite{Altshuler1977} corresponds to a neighborly oriented matroid of rank $5$ with $10$ elements. In~\cite{BokowskiGarms1987}, this matroid is proved to be non-realizable, giving a proof that there exist non-realizable neighborly matroids. 
One can also use a construction presented by Kortenkamp in~\cite{Kortenkamp1997} to build non-realizable neighborly matroids of corank $3$. 

In Theorem~\ref{thm:nonrealizable} we show that there exist non-realizable neighborly oriented matroids with $n$ vertices and rank $\rr$ for any $\rr\geq 5$ and $n\geq r+5$. 
Even more, in Theorem~\ref{thm:nonrealizablebound} we show that lower bounds proportional to~\eqref{eq:thebound} also hold for the number of labeled non-realizable neighborly oriented matroids.\\

Several known results on neighborly polytopes extend to all neighborly oriented matroids with combinatorial proofs that are often simpler than their geometric counterparts. This approach is present in the work of Sturmfels~\cite{Sturmfels1988} and Cordovil and Duchet~\cite{CordovilDuchet2000}.

For example, the Upper Bound Theorem still holds. Indeed, Stanley extended McMullen's Upper Bound Theorem to all simplicial spheres in~\cite{Stanley1975}. 
In particular, for oriented matroids we have that the number of faces of rank~$i$ of an oriented matroid of rank $\rd$ on $n$ elements is maximal for the \defn{alternating matroid} of rank $\rd$ on $n$ elements (the oriented matroid of the cyclic polytope $\cyc{\rd-1}{n}$), and equality is attained for some $i\geq \ffloor{\rd-1}{2}$ only by neighborly oriented matroids.

An important property of neighborly matroids of odd rank (in the realizable case, neighborly polytopes of even dimension) is that they are rigid. 
We call an oriented matroid \defn{rigid}\index{oriented matroid!rigid} if there is no other oriented matroid that has its face lattice; equivalently, if the face lattice (\ie the poset of positive covectors) determines its whole set of covectors. 
This result was first discovered by Shemer for neighborly polytopes~\cite{Shemer1982} and later proven by Sturmfels for all neighborly oriented matroids~\cite{Sturmfels1988}.

\begin{theorem}[{\cite[Theorem 4.2]{Sturmfels1988}}]\label{thm:neigharerigid}
 Every neighborly oriented matroid of odd rank is rigid.
\end{theorem}

\section{Balanced oriented matroids}\label{sec:balvsneigh}

The definition of neighborly matroids that we provided is based on their presentation by cocircuits. Our next goal is to find a characterization in terms of circuits. Said differently, to find a characterization of dual-to-neighborly matroids in terms of cocircuits. These are balanced matroids.

\begin{definition}\label{def:balanced}
An oriented matroid $\cM$ of rank $\rr$ and $n$ elements is \defn{balanced} if every cocircuit~$C$ of $\cM$ is \defn{halving}, \ie 
 \[\ffloor{n-\rr+1}{2}\leq |C^+|\leq\fceil{n-\rr+1}{2}.\]\index{balanced}\index{halving}
\end{definition}

A first observation is that if $n-\rr$ is odd, then $|\ul C|=\rr-1$ for each cocircuit $C\in\co(\cM)$. Hence every balanced matroid of odd corank is uniform. This corresponds to the fact that the vertices of an even dimensional neighborly polytope must be in general position, and implies that all even-dimensional neighborly polytopes are simplicial.

For uniform matroids, balancedness can also be described in terms of discrepancy.\index{discrepancy} 
\begin{definition}\label{def:discrepancy}
Let $\cM$ be an oriented matroid. The \defn{discrepancy} of a cocircuit $C\in\co(\cM)$ is \[\disc(C)=\big|{|C^+|-|C^-|}\big|,\] and the \defn{discrepancy} of an oriented matroid $\cM$ is the maximal discrepancy among its cocircuits, \index{$\disc(C)$}\index{$\disc (\cM)$}\[\disc (\cM)= \max_{C \in \co(\cM)} \disc (C).\]
\end{definition}
In particular, in a uniform oriented matroid, a cocircuit $C$ is halving if and only if it has discrepancy $0$ or $1$, \ie if it has the same number of positive and negative elements ($\pm 1$ if the corank is odd). And a
a uniform oriented matroid $\cM$ is balanced if and only if it has discrepancy $0$ or $1$, \ie all its cocircuits are halving.

\iftoggle{bwprint}{%
\begin{figure}[htpb]
\centering
 \subbottom[$\bullet\bullet\circ\circ\bullet\bullet$]{\label{sfig:gale_3cross}\includegraphics[width=.25\linewidth]{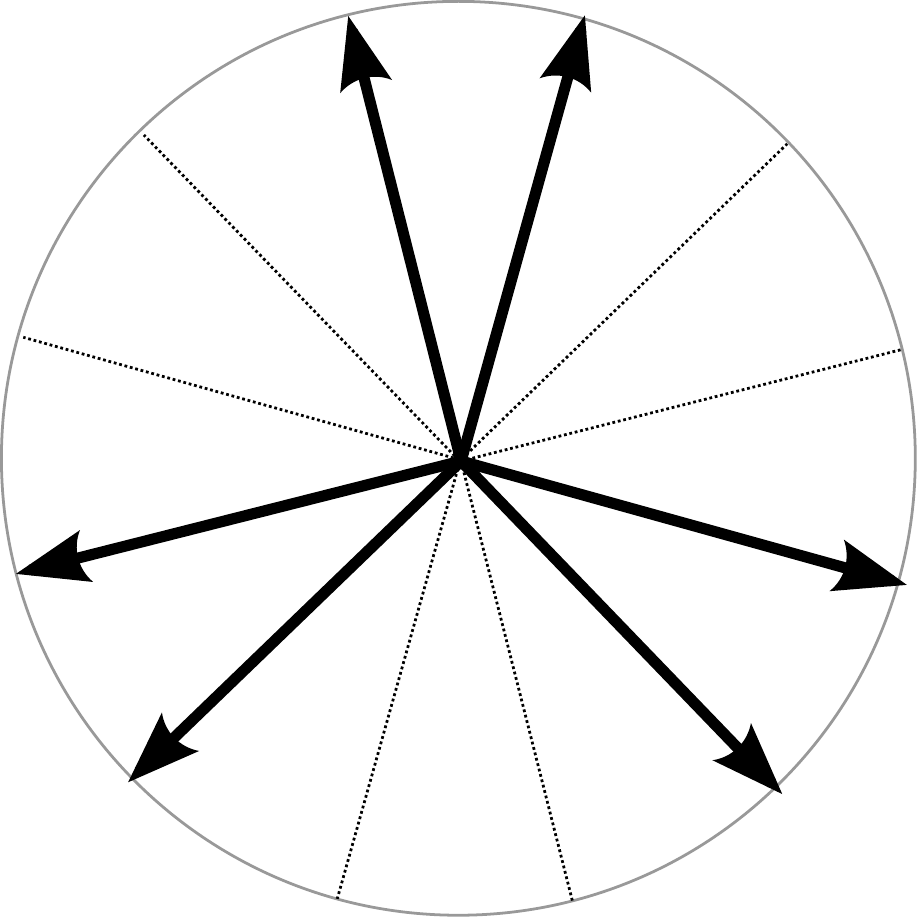}}\qquad\quad
 \subbottom[$\bullet\bullet\circ\bullet\circ\bullet$]{\label{sfig:gale_3cyc}\includegraphics[width=.25\linewidth]{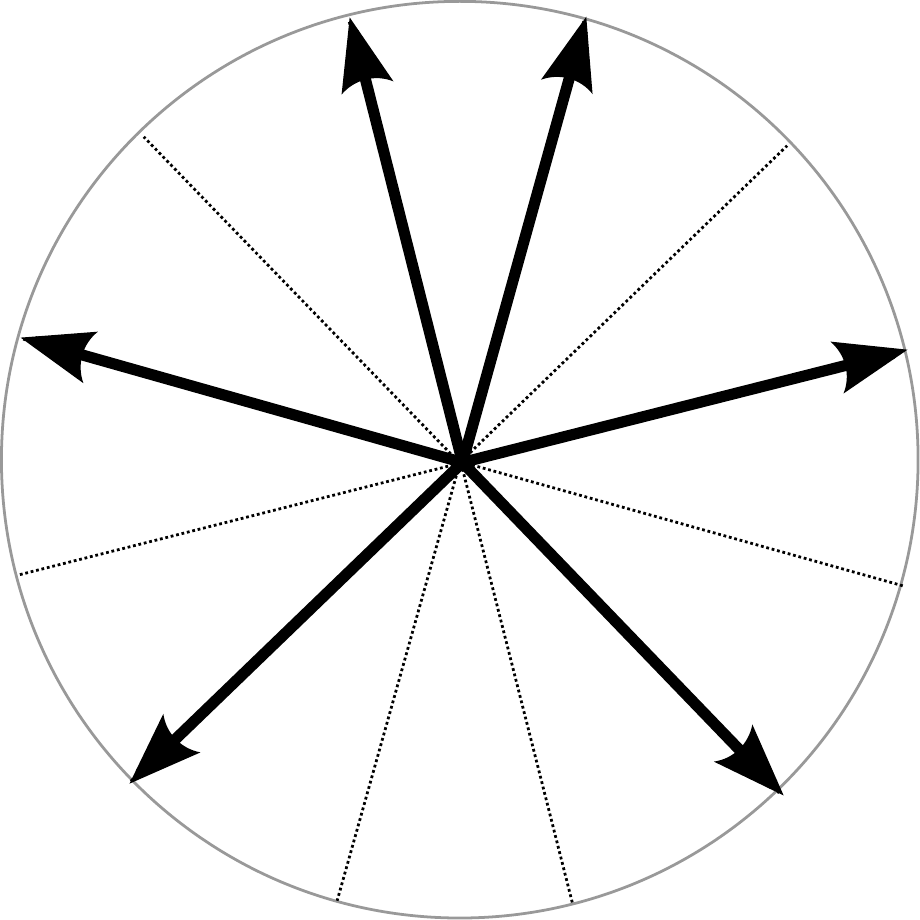}}\qquad\quad
 \subbottom[$\bullet\circ\bullet\circ\bullet\circ\bullet$]{\label{sfig:gale_4cyc}\includegraphics[width=.25\linewidth]{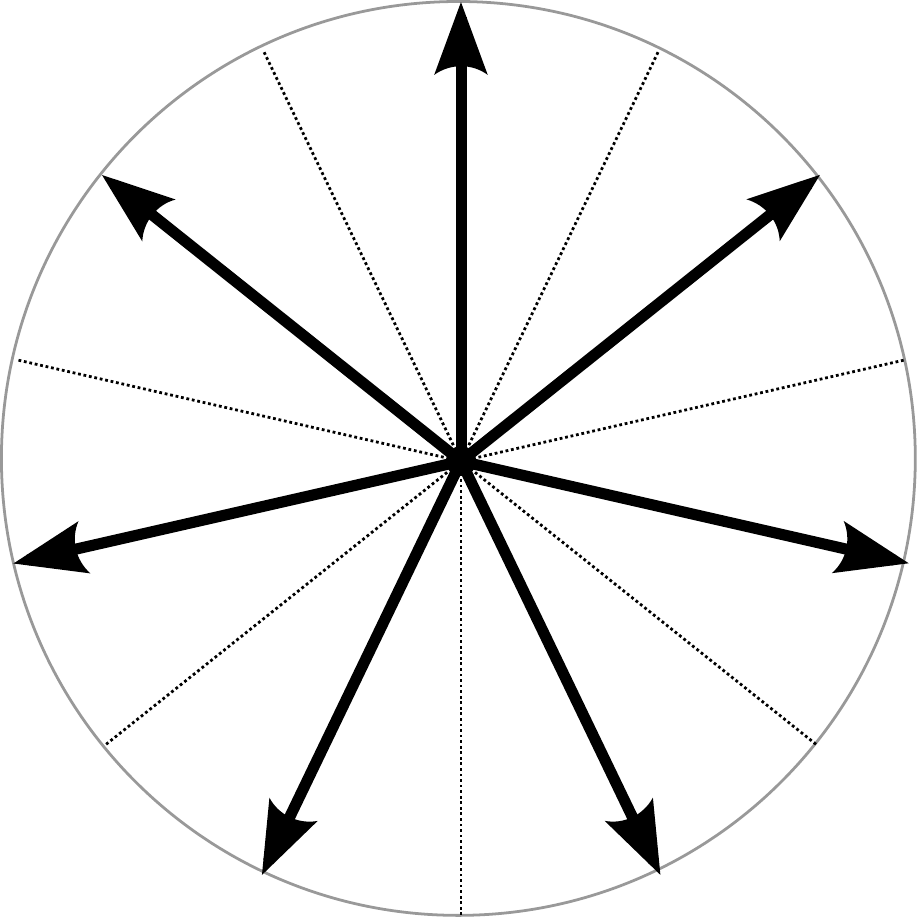}}
\caption[Gale duals of three neighborly polytopes]{Gale duals of three neighborly polytopes (and the corresponding affine diagrams): \subcaptionref{sfig:gale_3cross} a (non-regular) $3$-crosspolytope, \subcaptionref{sfig:gale_3cyc}  a cyclic polytope with $6$ vertices in $\RR^3$ and \subcaptionref{sfig:gale_4cyc} a cyclic polytope with $7$ vertices in~$\RR^4$.}
\label{fig:threegaleduals}
\end{figure}
}{%
\iftoggle{print}{%
\begin{figure}[htpb]
\centering
 \subbottom[$\bullet\bullet\circ\circ\bullet\bullet$]{\label{sfig:gale_3cross}\includegraphics[width=.25\linewidth]{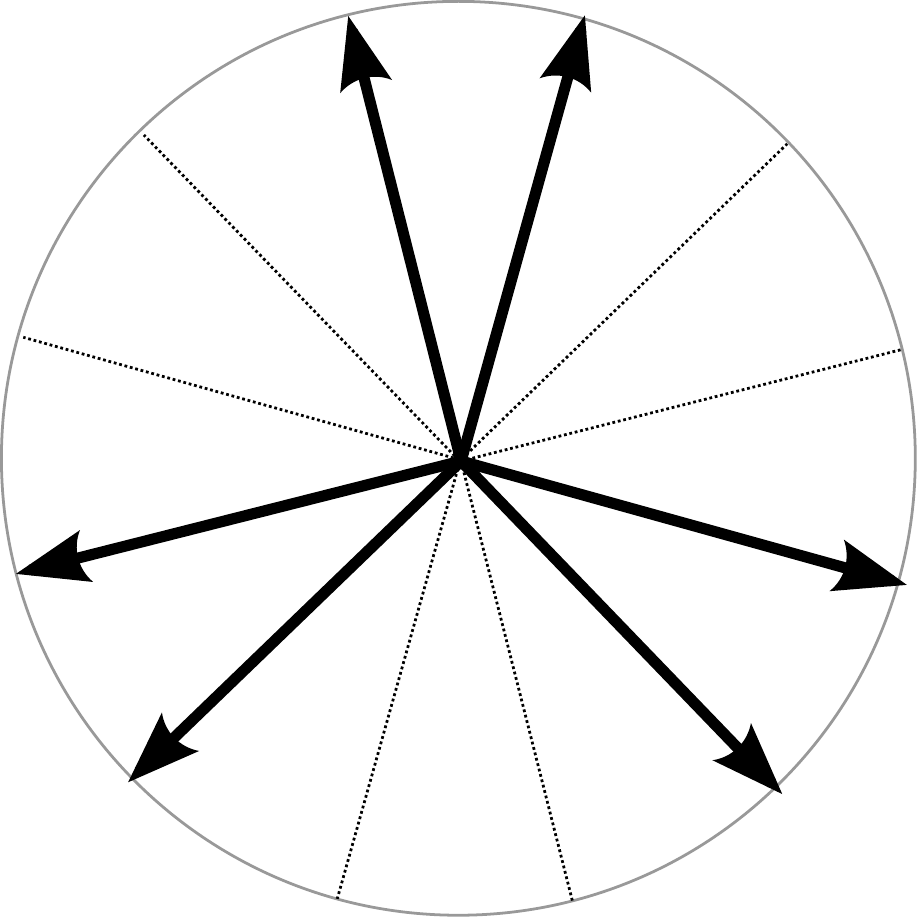}}\qquad\quad
 \subbottom[$\bullet\bullet\circ\bullet\circ\bullet$]{\label{sfig:gale_3cyc}\includegraphics[width=.25\linewidth]{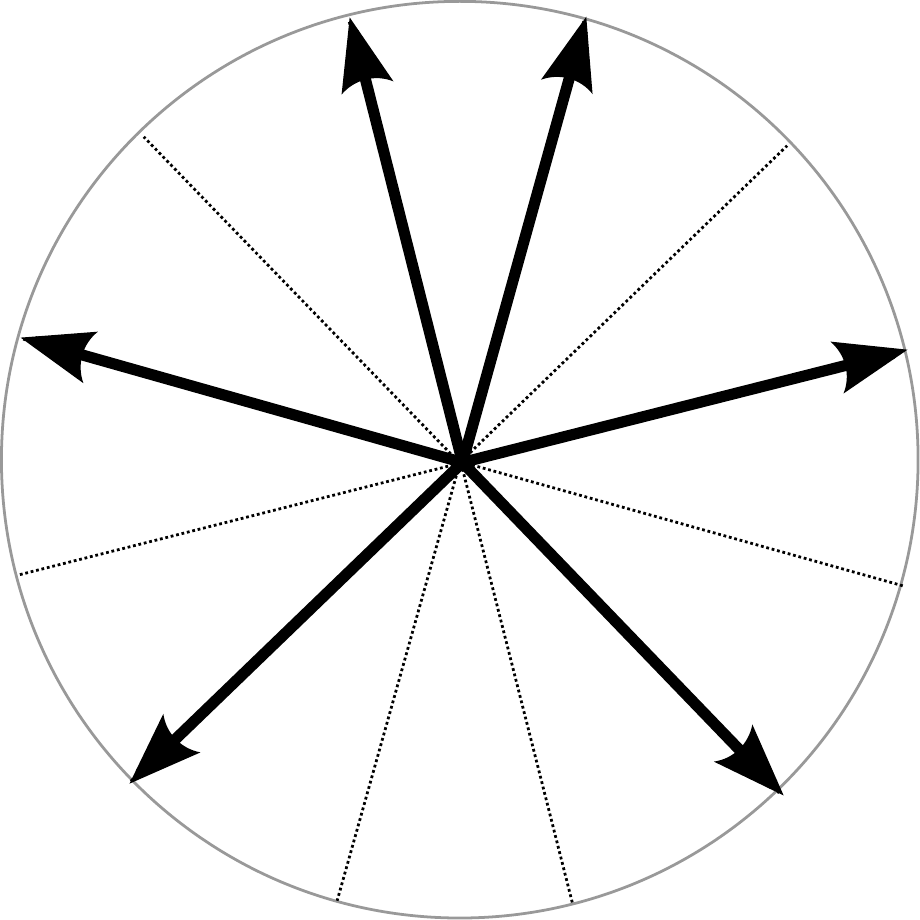}}\qquad\quad
 \subbottom[$\bullet\circ\bullet\circ\bullet\circ\bullet$]{\label{sfig:gale_4cyc}\includegraphics[width=.25\linewidth]{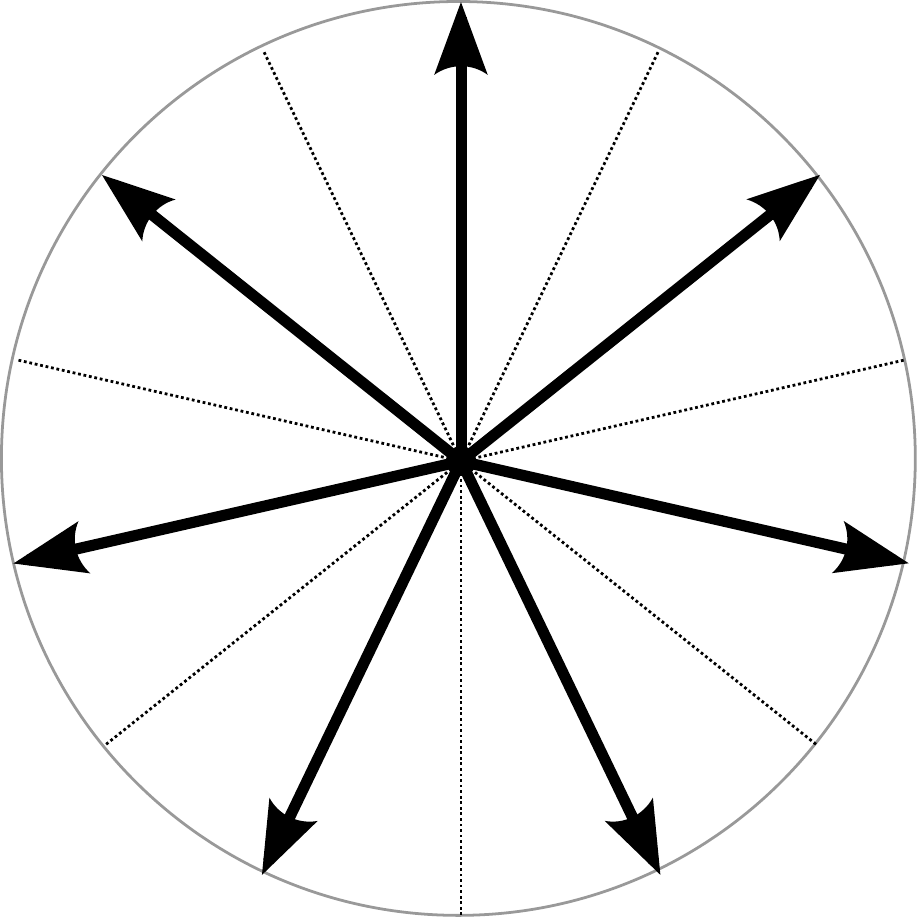}}
\caption[Gale duals of three neighborly polytopes]{Gale duals of three neighborly polytopes (and the corresponding affine diagrams): \subcaptionref{sfig:gale_3cross} a (non-regular) $3$-crosspolytope, \subcaptionref{sfig:gale_3cyc}  a cyclic polytope with $6$ vertices in $\RR^3$ and \subcaptionref{sfig:gale_4cyc} a cyclic polytope with $7$ vertices in~$\RR^4$.}
\label{fig:threegaleduals}
\end{figure}
}{%
\begin{figure}[htpb]
\centering
 \subbottom[$\bullet\bullet\circ\circ\bullet\bullet$]{\label{sfig:gale_3cross}\includegraphics[width=.25\linewidth]{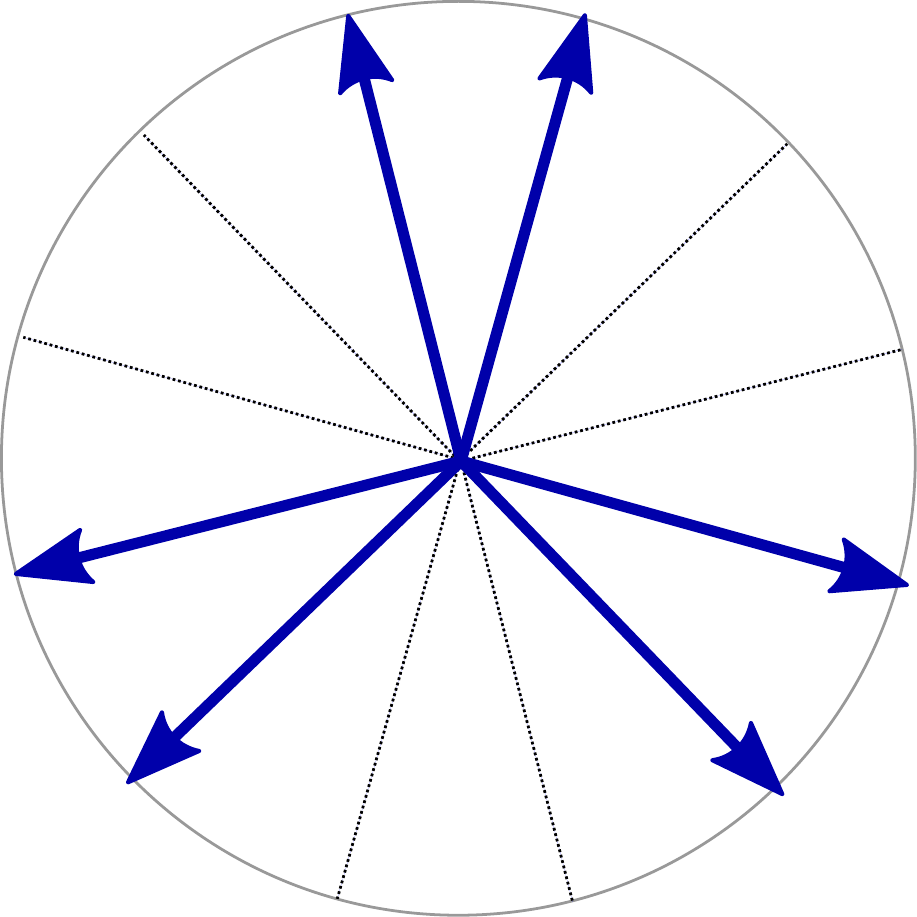}}\qquad\quad
 \subbottom[$\bullet\bullet\circ\bullet\circ\bullet$]{\label{sfig:gale_3cyc}\includegraphics[width=.25\linewidth]{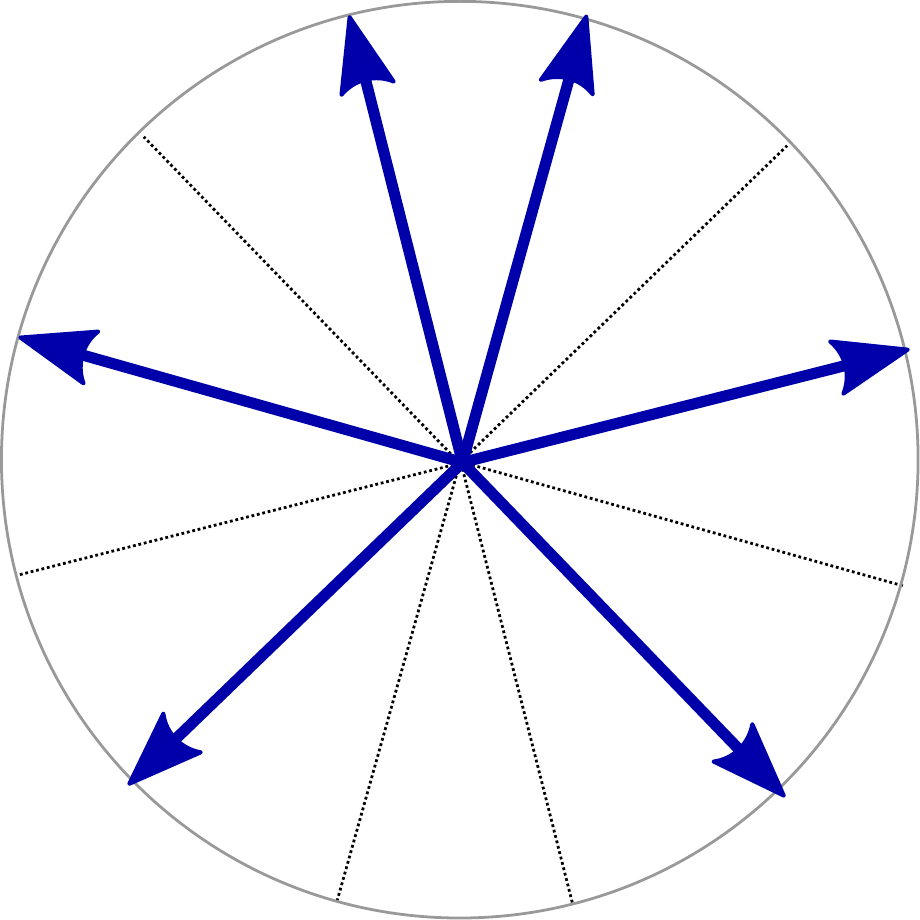}}\qquad\quad
 \subbottom[$\bullet\circ\bullet\circ\bullet\circ\bullet$]{\label{sfig:gale_4cyc}\includegraphics[width=.25\linewidth]{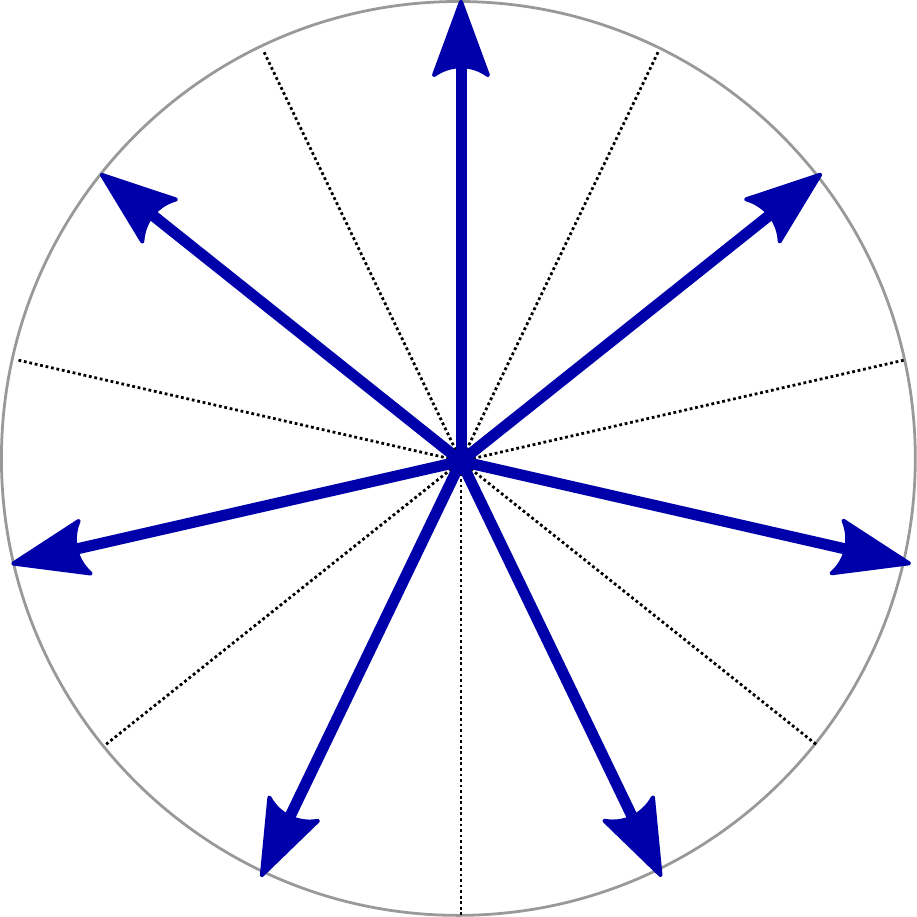}}
\caption[Gale duals of three neighborly polytopes]{Gale duals of three neighborly polytopes (and the corresponding affine diagrams): \subcaptionref{sfig:gale_3cross} a (non-regular) $3$-crosspolytope, \subcaptionref{sfig:gale_3cyc}  a cyclic polytope with $6$ vertices in $\RR^3$ and \subcaptionref{sfig:gale_4cyc} a cyclic polytope with $7$ vertices in~$\RR^4$.}
\label{fig:threegaleduals}
\end{figure}
}
}

 In Figure~\ref{fig:threegaleduals} three balanced vector configurations of rank $2$ are shown. This can be checked observing that for every hyperplane spanned by a vector, the remaining vectors are evenly split between the two sides (with $\pm 1$, if $n$ is even). 
A more detailed discussion on balanced configurations of rank~$2$ is done in Example~\ref{ex:balancedrank2}. The affine Gale diagram that we saw in Figure~\ref{fig:affineGale} is also balanced. The reader is invited to check this claim counting (signed) points at each side of the spanned hyperplanes, where positive points count as $+1$ and negative points as $-1$.
\\

The fact that neighborliness and balancedness are dual concepts is already implicit in the work of Gale~\cite{Gale1963} for polytopes, and one can find a proof for oriented matroids by Sturmfels in~\cite{Sturmfels1988}. Since it is a key result, we provide a proof at the end of the section for the sake of completeness.

\begin{proposition}[{\cite[Proposition 3.2]{Sturmfels1988}}]\label{prop:neighbaldduality}
 An oriented matroid $\cM$ is neighborly if and only if its dual matroid $\Gale\cM$ is balanced.
\end{proposition}

Lemma~\ref{lem:quotientsbalanced} follows directly from the definition of balanced oriented matroid and by duality implies the well known Lemma~\ref{lem:Hellyneighborly}, a Helly type result that is the oriented matroid version of \cite[Exercise 7.3.5.i]{GruenbaumEtal2003}.

\begin{lemma}\label{lem:quotientsbalanced}
An oriented matroid $\cM$ of rank $\rr\geq 1$ on a ground set $E$ is balanced if and only if $\cM/S$ is balanced for every independent subset $S\subset E$ of size $|S|= \rr-1$. \qed
\end{lemma}

\begin{lemma}\label{lem:Hellyneighborly}
An oriented matroid $\cM$ of rank $\rd$ on a ground set $E$ is neighborly if and only if the restriction of $\cM$ to $F$ is neighborly for all $F\subseteq E$ with $|F|=\rd+1$.\qed
\end{lemma}

\begin{example}[Uniform balanced matroids of rank $2$~\cite{AltshulerMcMullen1973}]\label{ex:balancedrank2}
Consider a uniform neighborly matroid $\cP$ of rank~$\rr$ and $\rr+2$ elements. Its dual is a balanced matroid $\cM$ of rank $2$ with $n$ elements. 
Since these are always realizable, we can think of it as a $2$-dimensional vector configuration. Using affine Gale diagrams, we can represent it as a sequence of $n$ black and white points. Figure~\ref{fig:threegaleduals} shows three examples of such a matroid. 

\begin{enumerate}[i)]
 \item If $n=2m+1$, then the affine Gale diagram of $\cP$ must be an alternating sequence of $2m+1$ black and white points. 
In this case, there is only one uniform balanced matroid, with a diagram that can be depicted as
 \(\bullet\left(\circ\bullet\right)^m\), where $\left(\circ\bullet\right)^m$ means $m$ consecutive copies of $\circ\bullet$. This configuration is dual to the alternating matroid, and proves that the only neighborly polytopes with $2m+1$ vertices in $\RR^{2m-2}$ are cyclic.

 \item If $n=2m$, then there are more uniform balanced matroids. 
Their diagrams consist in a sequence that alternates $k+1$ pairs of black points with $k$ pairs of white points, and such that in between a pair of black/white points with the next pair of white/black points there can be an alternating sequence of white-black/black-white points of arbitrary length. That is, 
\[\bullet\bullet(\circ\bullet)^{r_0}\circ\circ(\bullet\circ)^{l_1}\bullet\bullet(\circ\bullet)^{r_1}
\dots\circ\circ(\bullet\circ)^{l_k}\bullet\bullet(\circ\bullet)^{r_k},\]
for some $k$, $r_0,\dots, r_k$ and $l_1,\dots, l_k$ satisfying $1+2k+\sum_{i=1}^k l_i+\sum_{i=0}^k r_i=m$. The dual of the alternating matroid corresponds to the case $k=0$.
\end{enumerate}
\end{example}

To prove Proposition~\ref{prop:neighbaldduality}, we will need this oriented matroid version of Farkas' Lemma. 

\begin{proposition}[3-Painting Lemma~{\cite{BlandLasVergnas1978}}]\label{prop:3coloring}
Let $\cM$ be an oriented matroid on $E$, with vectors $\ve$ and covectors $\cov$. Let $B\cup G\cup W=E$ be a partition of $E$ into disjoint subsets $B, G, W$  and let $e \in B$. Then exactly one of the following statements holds:
\begin{enumerate}[(i)]
\item $\exists C \in \cov$ such that $C(e)>0$, $C(b)\geq 0$ for $b\in B$ and $C(w)= 0$ for $w\in W$;
\item $\exists X \in \ve$ such that $X(e)> 0$, $X(b)\geq 0$ for $b\in B$, $X(g)=0$ for $g\in G$;
 \end{enumerate}
\end{proposition}

\begin{proof}[Proof of Proposition~\ref{prop:neighbaldduality}]
Let $\rd$ be the rank of $\cM$, and $\rr=n-\rd$ the rank of $\Gale\cM$. 

We prove first that if $\cM$ is neighborly, then $\Gale\cM$ is balanced. Let $X=(X^+,X^-)$ be a circuit of $\cM$, we want to see that \[\ffloor{\rd+1}{2}=\ffloor{n-\rr+1}{2}\leq|X^+|\leq\fceil{n-\rr+1}{2}=\fceil{\rd+1}{2}.\]

If $|X^+|\leq\ffloor{\rd+1}{2}-1$, then there must be a covector $C$ of $\cM$ with 
 $C^+=\emptyset$, $C^-=X^-\cup X^0$ and $C^0=X^+$, because $\cM$ is neighborly. 
But then $C\not\perp X$, contradicting circuit-cocircuit orthogonality. Indeed, for all $e\in \ul C\cap \ul X$ we have that $e\in C^-\cap X^-$ and thus that $C(e)X(e)=+$. Hence, we have seen that $|X^+|\geq\ffloor{\rd+1}{2}$. One can analogously prove that $|X^-|\geq\ffloor{\rd+1}{2}$. 
Since $|X^0|\geq n-\rd-1$, this implies that $|X^+|=n-|X^0|-|X^-|\leq n-(n-\rd-1)-\ffloor{\rd+1}{2}=\fceil{\rd+1}{2}$.

To prove the reciprocal, we will need 3-painting. Let $W$ be any subset of $E$ of size $\leq \ffloor{\rd-1}{2}$, and let $B=E\setminus W$. We want to find a covector $C\in\cov(\cM)$ such that $C^+=B$ and $C^0=W$.

If for every $e\in B$ there is a covector $C_e$ with $C_e(e)>0$, $C_e(b)\geq 0$ for $b\in B$ and $C_e(w)=0$ for $w\in W$, then the composition of all these covectors is the desired covector $C=\circ_{e\in B}C_e$.

Assume the contrary and we will reach a contradiction. This means that there is some $e\in B$ such that there does not exist any covector $C \in \cov$ such that $C(e)>0$, $C(b)\geq 0$ for $b\in B$ and $C(w)= 0$ for $w\in W$. 
By the painting lemma, there would be a vector $X$ of $\cM$ such that $X(b)\geq 0$ for all $b\in B$. Hence, $X^-\subset W$, but then $|W|\leq \ffloor{\rd-1}{2}<\ffloor{\rd+1}{2}$, which is a contradiction with the balancedness of $\Gale\cM$.
\end{proof}
\medskip

\begin{remark}
 In general, an oriented matroid $\cM$ is called \defn{$k$-neighborly} when every subset of elements of size $\leq k$ is a face of $\cM$. Proposition~\ref{prop:neighbaldduality} can be generalized to characterize $k$-neighborly oriented matroids in terms of their circuits. Indeed, one can prove that $\cM$ is $k$-neighborly if and only if $|X^+|\geq k+1$ for every circuit $X\in\ci(\cM)$, see~\cite[Proposition 4.1]{CordovilDuchet2000}.
\end{remark}

\section{Inseparability}

Inseparability is an essential (albeit straightforward) tool that will be used extensively in what follows. It is strongly related to the concept of universal edges, which is a basic element of the sewing construction.

\begin{definition}
Given an oriented matroid~$\cM$ on a ground set $E$, and $\alpha\in\{+1,-1\}$, we say that two elements $p,q\in E$ are \defn{$\alpha$-inseparable} in $\cM$ if \begin{equation}\label{eq:definseparable}X(p)=\alpha X(q)\end{equation}
for each circuit $X\in\ci(\cM)$ with $p,q\in \underline X$. 

In the literature, $(+1)$-inseparable elements are also called \defn{covariant} and $(-1)$-inseparable elements \defn{contravariant} (see~\cite[Section 7.8]{OrientedMatroids1993}). \index{inseparable}
\end{definition}

A first useful property is that inseparability is preserved by duality (with a change of sign). This allows to characterize the inseparability of $p$ and~$q$ (originally defined in terms of circuits) with an analogue of \eqref{eq:definseparable} for cocircuits of $\cM$ (again, with a change of sign).

\begin{lemma}[{\cite[Exercise 7.36]{OrientedMatroids1993}}]\label{lem:insep}
 A pair of elements $p$ and $q$ are $\alpha$-inseparable in $\cM$ if and only if they are $(-\alpha)$-inseparable in $\Gale\cM$. \qed
\end{lemma}

The \defn{inseparability graph}\index{inseparability graph} of an oriented matroid $\cM$, \defn{$\IG(\cM)$}\index{$\IG(\cM)$}, is 
the graph that has the elements of $\cM$ as vertices and the pairs of inseparable elements as edges.

\iftoggle{bwprint}{%
\begin{figure}[htpb]
\begin{center}
\includegraphics[width=.6\textwidth]{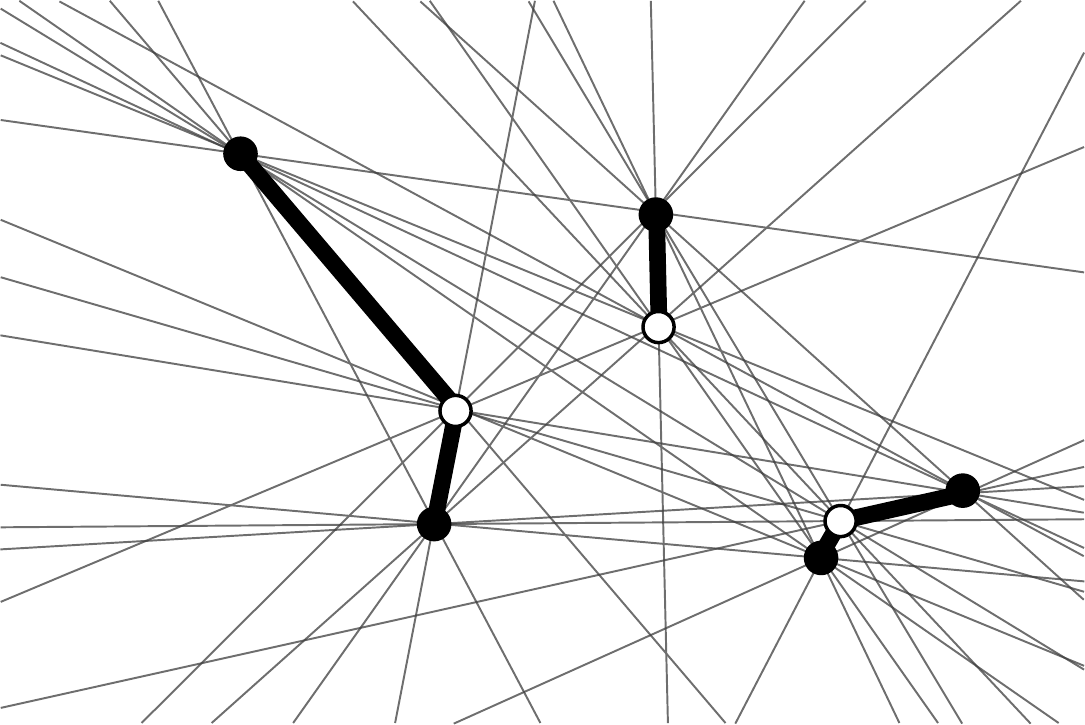}
\end{center}
 \caption{Inseparability graph of a balanced matroid of rank $3$ and $8$ elements. All the inseparable pairs are $(+1)$-inseparable.}
 \label{fig:inseparable}
\end{figure}
}{%
\begin{figure}[htpb]
\begin{center}
\includegraphics[width=.6\textwidth]{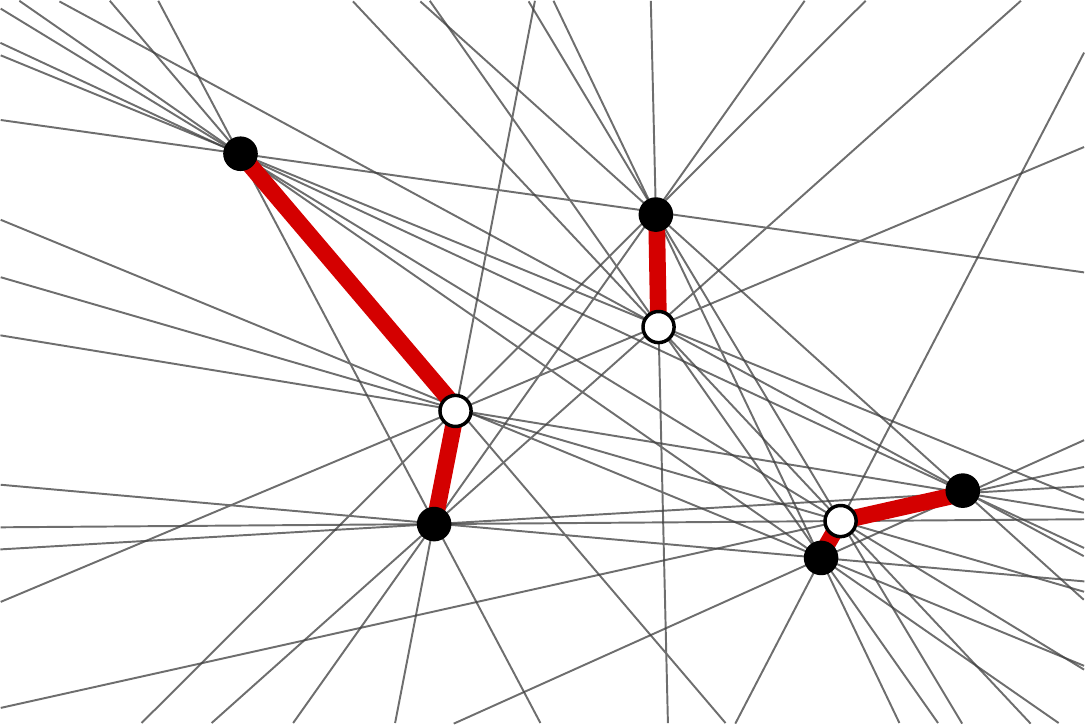}
\end{center}
 \caption{Inseparability graph of a balanced matroid of rank $3$ and $8$ elements. All the inseparable pairs are $(+1)$-inseparable.}
 \label{fig:inseparable}
\end{figure}
}

\begin{example}
Figure~\ref{fig:inseparable} shows the inseparability graph of the balanced matroid of Figure~\ref{fig:affineGale}, which is dual to a neighborly $4$-polytope with $8$ vertices. The pairs joined with a thick segment are $(+1)$-inseparable, since for every hyperplane spanned by the remaining elements, they lie at opposite sides (as signed points, of course).

This neighborly polytope is not combinatorially equivalent to a cyclic polytope, since the inseparability graph of the cyclic polytope is a cycle (cf. Theorem~\ref{thm:IG}).
\end{example}

We proceed to present some results on inseparable elements that will be useful later. The first one concerns inseparable elements and shows the relation between circuits/cocircuits through $x$ and circuits/cocircuits through $y$ when $x$ and $y$ are inseparable.

\begin{lemma}\label{lem:circinseparable}
 Let $\cM$ be a uniform oriented matroid with two $\alpha$-inseparable elements~$x$ and~$y$.
\begin{enumerate}
 \item For every circuit $X\in\ci(\cM)$ with $X(x)=0$ and $X(y)\neq 0$, there is a circuit $X'\in\ci(\cM)$ with $X'(x)=-\alpha X(y)$, $X'(y)=0$ and $X'(e)=X(e)$ for all $e\notin\{x,y\}$;
 \item For every cocircuit $C\in\co(\cM)$ with $C(x)=0$ and $C(y)\neq 0$, there is a cocircuit $C'\in\co(\cM)$ with $C'(x)=\alpha C(y)$, $C'(y)=0$ and $C'(e)=C(e)$ for all $e\notin\{x,y\}$.
\end{enumerate}

\end{lemma}
\begin{proof}
Both statements are equivalent by duality. We prove the first one.

Let $X'\in \ci(\cM)$ be the circuit with support $\underline X'=\underline X \setminus y \cup x$ and such that $X'({x})=-\alpha X(y)$. This circuit exists because $\cM$ is uniform. We will see that $X'({e})= X(e)$ for all $e\in \underline X'\setminus x$.
Let $e\in \underline{X}\setminus y$, and let $C$ be the cocircuit of~$\cM$ with $\underline C=E\setminus(\underline{X}\setminus y\setminus e)$. That makes $\underline C\cap \underline X=\{e,y\}$. Since $y$ and $x$ are $\alpha$-inseparable, $C(x)=-\alpha C(y)$, and by circuit-cocircuit orthogonality,
\begin{displaymath}X(y)X(e)=-C(y)C(e)=\alpha C({x})C(e).\end{displaymath}
But $\underline C\cap \underline X'=\{e,x\}$, and hence, again by orthogonality, \(X'(x)X'(e)=-C({x})C(e).\)
The conclusion now follows from $X'({x})=-\alpha X(y)$.
\end{proof}

In this lemma, uniformity is important, since the result does not hold in general. For example, in the vector configuration $\vv V$ of Figure~\ref{sfig:GaleV}, $\vv v_1$ and $\vv v_3$ are $(-1)$-inseparable. However, the signed set $X$ with $X^+=\{ \vv v_3\}$ and $X^-=\{\vv v_2\}$ is a circuit of $\cM(\vv V)$, while $X'$ with ${X'}^+=\{\vv v_1\}$ and ${X'}^-=\{\vv v_2\}$ is not.

The following lemma concerns inseparable elements of neighborly and balanced oriented matroids and explains why all the inseparable pairs of the previous example are $(+1)$-inseparable.
\begin{lemma}\label{lem:balonlycovar}
 If $\cM$ is a balanced oriented matroid of rank $\rr\geq 2$ with $n$ elements such that $n-\rr-1$ is even, then all inseparable elements in $\cM$ must be $(+1)$-inseparable.

 Analogously, if $\cP$ is a neighborly oriented matroid of odd rank $\rd$ with at least $\rd+2$ elements, then all inseparable elements in $\cP$ must be $(-1)$-inseparable.
\end{lemma}
\begin{proof}
 Both results are equivalent by duality. To prove the second claim, observe that if $p$ and $q$ are $\alpha$-inseparable in $\cP$, then they are also in $\cP\setminus S$ for any $S$ that contains neither $p$ nor $q$. Hence we can remove elements from $\cP$ until we are left with a 
 neighborly matroid with $\rd+2$ elements. The only neighborly matroid of odd rank $\rd$ and corank~$2$ is the alternating matroid (cf. Example~\ref{ex:balancedrank2}), which only has $(-1)$-inseparable pairs.
\end{proof}

We end with an observation about inseparability in balanced matroids, whose proof is fairly easy and left to the reader.

\begin{lemma}\label{lem:insareuni}
  If a pair $x, y$ of elements of a balanced matroid $\cM$ are $(+1)$-inseparable then $\cM\setminus\{x,y\}$ is balanced. If moreover the corank of $\cM$ is odd, the converse is also true; that is, $\cM\setminus\{x,y\}$ is balanced only if $x, y$ are $(+1)$-inseparable.\qed
\end{lemma}

\section{Lexicographic extensions}\label{sec:le}

Lexicographic extensions play a central role in our results. They were introduced by Las Vergnas in 1978~\cite{LasVergnas1978} and have several applications (cf.~\cite{Kortenkamp1997}, \cite{Santos2002}, \cite{SturmfelsZiegler1993}, \cite{Todd1985}). They are also known under the name of \defn{principal extensions} (cf.~\cite{BilleraMunson1984}, \cite{JaggiManiLevitskaSturmfelsWhite1989}).
Even if they are just a particular case of extensions, they are a very versatile tool for constructing many different polytopes and oriented matroids. In this section we present well-known definitions and results on extensions of oriented matroids, mainly from \cite[Chapter~7]{OrientedMatroids1993}, as well as derive some small results that we will user later on.

\subsection{Single element extensions}

Let $\cM$ be an oriented matroid on a ground set $E$. A \defn{single element extension}\index{single element extension} of $\cM$ by an element $p$ is an oriented matroid $\tilde \cM$ on the ground set $\tilde E= E\cup \{p\}$ for some $p\notin E$, such that every circuit of $\cM$ is a circuit in $\tilde \cM$. 
Equivalently, $\tilde \cM$ is a single element extension of $\cM$ if $\cM$ is a restriction of~$\tilde \cM$ by deleting one element. That is, $\tilde\cM\setminus p=\cM$. We will only consider extensions that do not increase the rank, \ie $\rank(\tilde \cM)=\rank(\cM)$.\\

A concept crucial to understanding a single element extension of $\cM$ is its signature, which we define in the following proposition using the formulation in~\cite[Proposition 7.1.4]{OrientedMatroids1993} of a result originally from \cite{LasVergnas1978}.

\begin{proposition}(\cite[Proposition 7.1.4]{OrientedMatroids1993})
Let $\tilde \cM$ be a single element extension of~$\cM$ by $p$. Then, for every cocircuit $C=(C^+,C^-)\in\co(\cM)$, there is a unique way to extend~$C$ to a cocircuit of~$\tilde \cM$: exactly one of $(C^+\cup \{p\},C^-)$, $(C^+,C^-\cup  \{p\})$ or $(C^+,C^-)$ is a cocircuit of $\tilde \cM$. 

That is, there is a unique function $\gs$ from $\co(\cM)\rightarrow \{+,-,0\}$ such that for each $C\in \co(\cM)$ there is a cocircuit $\tilde C\in\co(\tilde\cM)$ with $\tilde C(p)=\gs(C)$ and $\tilde C(e)=C(e)$ for $e\in E$. The function $\gs$ is called the \defn{signature}\index{single element extension!signature} of the extension.

Moreover, the signature $\gs$ uniquely determines the oriented matroid $\tilde\cM$.
\end{proposition}

Not every map from $\co(\cM)$ to $\{0,+,-\}$ corresponds to the signature of an extension, a property that can be checked with the following theorem of Las Vergnas (see also~\cite[Theorem 7.1.8]{OrientedMatroids1993}).

\begin{theorem}[\cite{LasVergnas1978}]\label{thm:localization}
 Let $\cM$ be an oriented matroid and $\gs:\co(\cM)\rightarrow \{+,-,0\}$ a cocircuit signature satisfying $\gs(-C)=-\gs(C)$ for all $C\in \co(\cM)$. Then $\gs$ is a the signature of a single element extension if and only if $\gs$ defines a single element extension on every contraction of $\cM$ of rank~$2$.
\end{theorem}

In the setting of a vector configuration $\vv V$, the signature of the extension by $\vv p$ records at which side of each of the cocircuit-defining hyperplanes does~$\vv p$ lie (see Figure~\ref{fig:localization} for an example). 
Equivalently, we can think of a signature as an orientation of the hyperplane arrangement $\cH$ of all hyperplanes spanned by the vectors in~$\vv V$. If the extension can be realized without modifying $\vv V$, then the signature corresponds to an acyclic orientation of~$\cH$, that is, an orientation where the intersection of all positive halfspaces is not empty (this intersection is the cell of~$\cH$ where we will add our new element~$\vv p$). 

\iftoggle{bwprint}{%
\begin{figure}[htpb]
\begin{center}
\includegraphics[width=.6\linewidth]{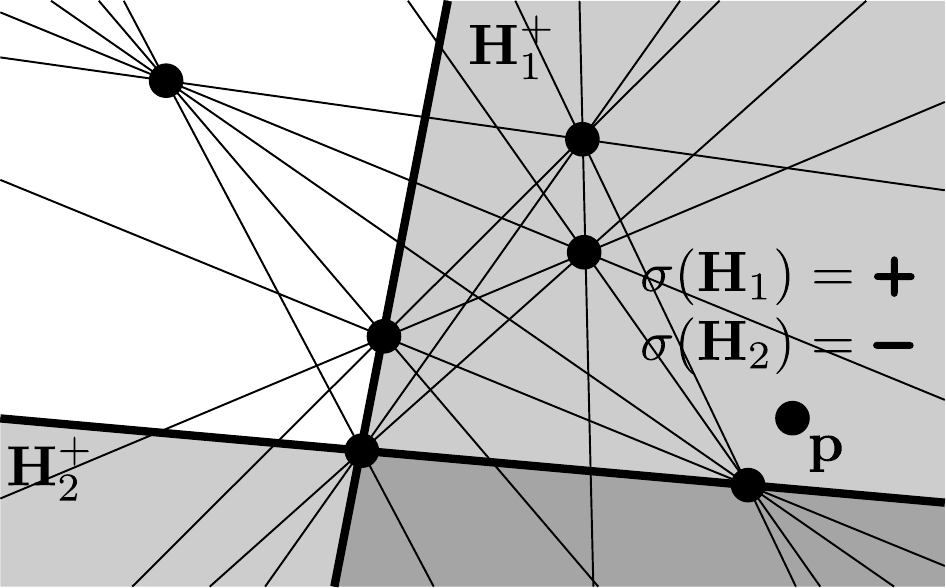}
\end{center}
\caption[Example of a signature of a single element extension.]
{If $C_1$ and $C_2$ are the cocircuits defined by the hyperplanes ${\vvh H}_1$ and ${\vvh H}_2$ then the signature $\gs$ of the extension by $\vv p$ fulfills $\gs(C_1)=+$ and $\gs(C_2)=-$.}\label{fig:localization}
\end{figure}
}{%
\begin{figure}[htpb]
\begin{center}
\includegraphics[width=.6\linewidth]{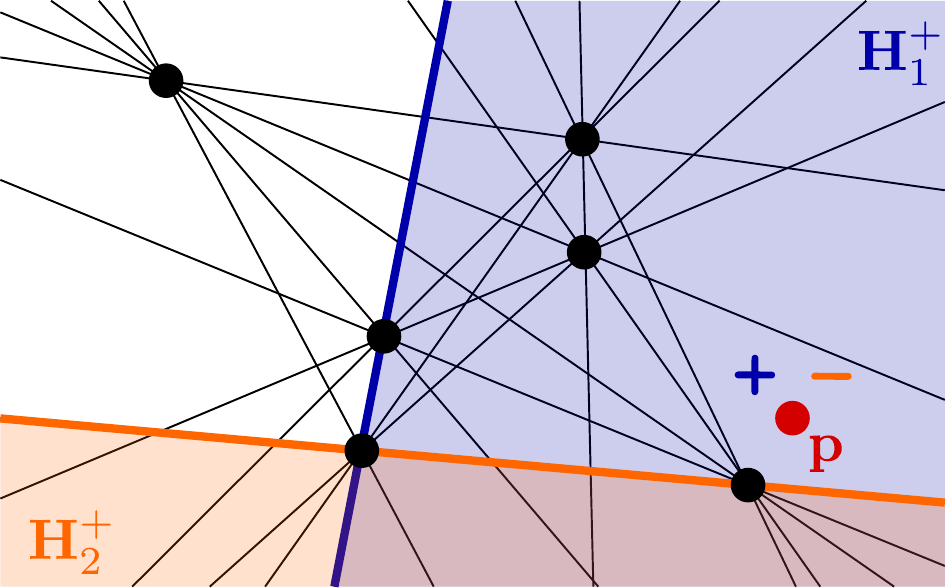}
\end{center}
\caption[Example of a signature of a single element extension.]
{If $C_1$ and $C_2$ are the cocircuits defined by the hyperplanes ${\vvh H}_1$ and ${\vvh H}_2$ then the signature $\gs$ of the extension by $\vv p$ fulfills $\gs(C_1)=+$ and $\gs(C_2)=-$.}\label{fig:localization}
\end{figure}
}

\subsection{Lexicographic extensions}
We are interested in one particular family of single element extensions called lexicographic extensions.\index{lexicographic extension}

\begin{definition}
 Let $\cM$ be a rank~$\rr$ oriented matroid on a ground set~$E$. 
Let $(a_1,a_2,\dots, a_k)$ be an ordered subset of $E$ and let $(\ep_1,\ep_2,\dots,\ep_k)\in\{+,-\}^k$ be a sign vector. 
The \defn{lexicographic extension} $\cM[p]$ of $\cM$ by $p=[a_1^{\ep_1},a_2^{\ep_2},\dots,a_k^{\ep_k}]$ is the oriented matroid on the ground set $E\cup \{p\}$ which is the single element extension of $\cM$ whose signature $\sigma:\co(\cM)\rightarrow \{+,-,0\}$ maps $C\in \co(\cM)$ to
\begin{equation*}
\gs(C)\mapsto
\begin{cases} \ep_iC({a_i})& \text{if $i$ is minimal with $C({a_i})\neq 0$,}
\\
0&\text{if $C({a_i})=0$ for $i=1,\dots,k$.}
\end{cases}
\end{equation*}
We will also use $\cM[a_1^{\ep_1},\dots,a_k^{\ep_k}]$ to denote the lexicographic extension $\cM[p]$ of $\cM$ by $p=[a_1^{\ep_1},\dots,a_k^{\ep_k}]$.
\end{definition}

\begin{remark}
\begin{enumerate}
\item We can always assume that $a_1,\dots,a_k$ are independent. In fact, if $j$ is the first index such that $a_1,\dots,a_j$ are dependent, then \[\cM[a_1^{\ep_1},\dots,a_{j-1}^{\ep_{j-1}},a_{j}^{\ep_{j}},a_{j+1}^{\ep_{j+1}},\dots,a_{k}^{\ep_{k}}]=\cM[a_1^{\ep_1},\dots,a_{j-1}^{\ep_{j-1}},a_{j+1}^{\ep_{j+1}},\dots,a_{k}^{\ep_{k}}].\]
\item If $\cM$ is a uniform matroid of rank $\rr$, then $\cM[a_1^{\ep_1},\dots,a_k^{\ep_k}]$ is uniform if and only if $k=\rr$. This is the most interesting case for us.
\end{enumerate}
\end{remark}

Lexicographic extensions preserve realizability. Indeed, the following observation from \cite[Section 7.2]{OrientedMatroids1993} shows how to realize lexicographic extensions of realizable oriented matroids. Since extensions of a non-realizable matroid must be non-realizable because they have a non-realizable minor, this means that $\cM[p]$ is realizable if and only if $\cM$ is realizable.

\begin{proposition}\label{prop:realizablele}
 If $\cM$ is realizable
then all lexicographic extensions of~$\cM$ are realizable.
\end{proposition}
\begin{proof}[Proof idea]
 Let the vector configuration $\vv V\subset \RR^r$ be a realization of $\cM$. Then a realization of $\cM[\vv v_1^{\ep_1},\vv v_2^{\ep_2},\dots,\vv v_k^{\ep_k}]$ is $\vv V\cup \vv v$, where $\vv v:=\ep_1\vv v_1+\delta \ep_2\vv v_2+\delta^2 \ep_3\vv v_3+\dots+\delta^{k-1} \ep_k\vv v_k$ for some small enough $\delta>0$.
\end{proof}

Actually, in the setting of a vector configuration $\vv V$, the lexicographic extension by $\vv v=[\vv v_1^{\ep_1},\vv v_2^{\ep_2},\dots,\vv v_k^{\ep_k}]$ is very easy to understand.  For every hyperplane ${\vvh H}$ spanned by vectors in $\vv V\setminus\{\vv v_1\}$,
the new vector $\vv v$ must lie on the same side as $\ep_1\vv v_1$; for hyperplanes containing $\vv v_1$ but not $\vv v_2$, $\vv v$ must lie on the same side as $\ep_2\vv v_2$; \etc\  
This is clearly achieved by the vector $\vv v=\ep_1\vv v_1+\delta \ep_2\vv v_2+\delta^2 \ep_3\vv v_3+\dots+\delta^{k-1} \ep_k\vv v_k$, constructed by placing a new vector on top of $\ep_1 \vv v_1$, and perturbing it slightly towards $\ep_2 \vv v_2$, then towards $\ep_3\vv v_3$ and so on. See Figure~\ref{fig:le} for an example of this procedure. 
\\

\iftoggle{bwprint}{%
\begin{figure}[htpb]
\begin{center}
\includegraphics[width=.6\linewidth]{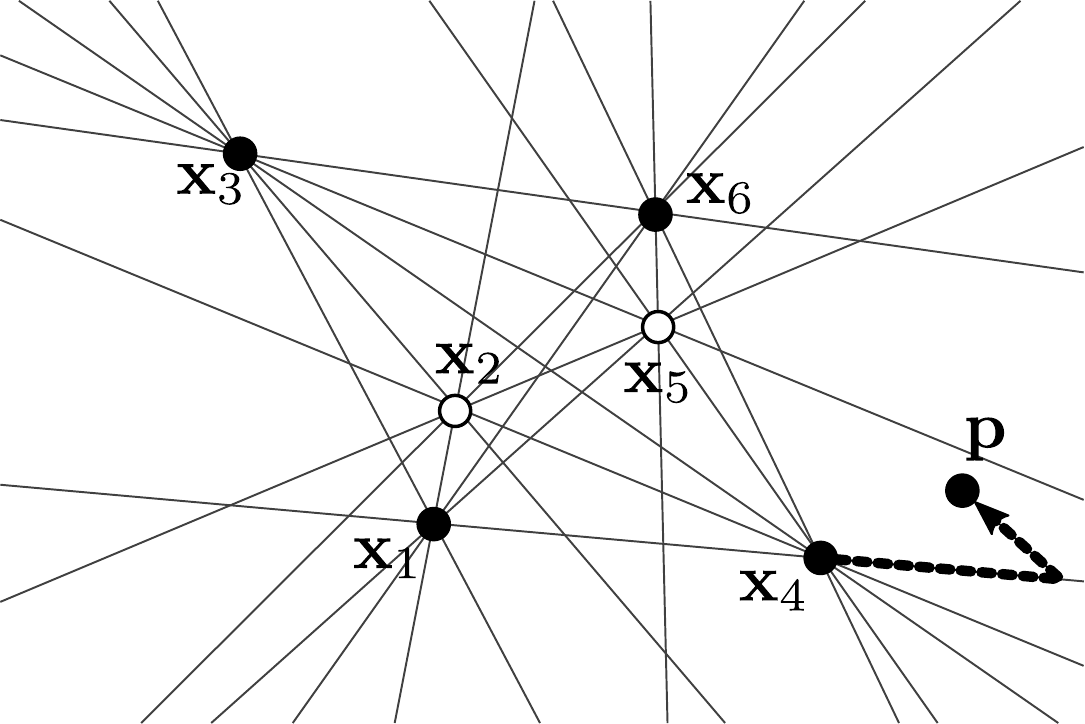}
\end{center}
\caption{An affine Gale diagram, and its lexicographic extension by $\vv p=[\vv x_4^+,\vv x_1^-,\vv x_6^+]$.}\label{fig:le}
\end{figure}
}{%
\begin{figure}[htpb]
\begin{center}
\includegraphics[width=.6\linewidth]{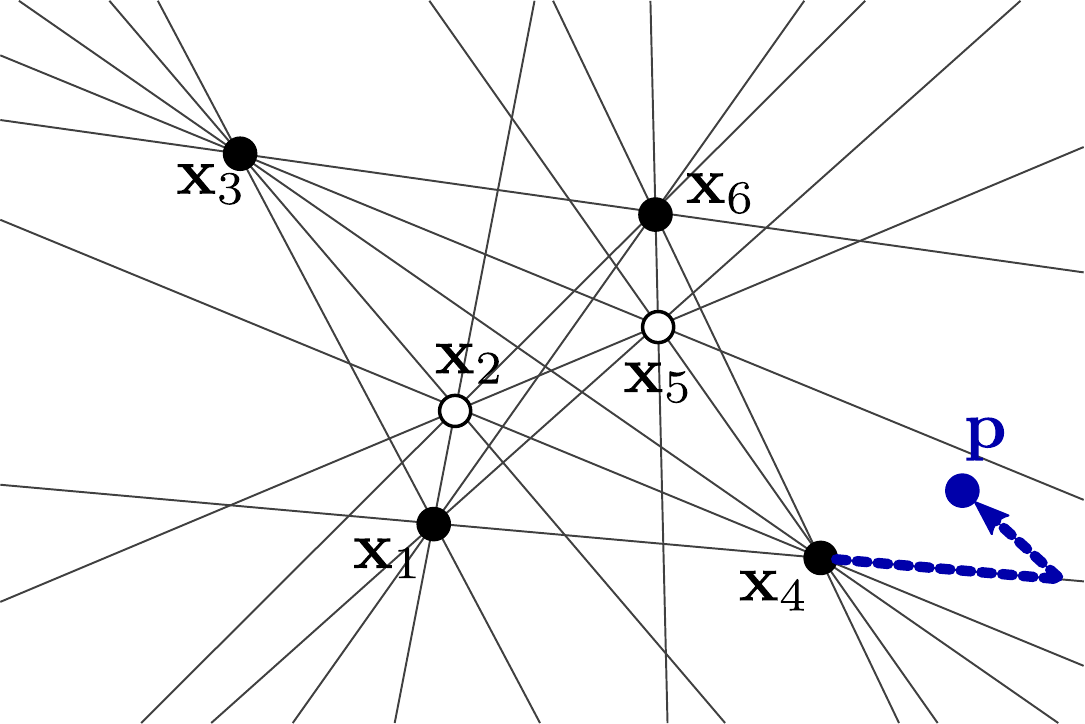}
\end{center}
\caption{An affine Gale diagram, and its lexicographic extension by $\vv p=[\vv x_4^+,\vv x_1^-,\vv x_6^+]$.}\label{fig:le}
\end{figure}
}

Alternatively, lexicographic extensions can be defined in terms of inseparability, as the following result shows. Its proof, which follows easily from the definitions, is omitted.

\begin{lemma}\label{lem:leinseparable}
In a lexicographic extension $\cM [p]$ of $\cM$ by $p=[a_1^{\ep_1},\dots,a_k^{\ep_k}]$, $p$ and~$a_1$ are always $(-\ep_1)$-inseparable. Even more, $p$ and $a_i$ are $(-\ep_i)$-inseparable in $\cM[p]/\{a_1,\dots,a_{i-1}\}$ for $i=1\dots k$, and this property characterizes this single element extension. \qed
\end{lemma}

A key property is that lexicographic extensions on uniform matroids behave well with respect to contractions. 
The upcoming Proposition~\ref{prop:allquotientsofle} can be used to iteratively explain all cocircuits of a lexicographic extension, and hence can be seen as the restriction of~\cite[Proposition 7.1.4]{OrientedMatroids1993} to lexicographic extensions. It is a very useful tool that will be used extensively. Before proving it, it is useful to state the following lemma, which deals with the simultaneous contraction and deletion of $p$ and $a_1$.

\begin{lemma}\label{lem:contractdeletele}
If $\cM$ is uniform and  $p=[a_1^{\ep_1},\dots]$, then
 \[\cM/a_1 = ( \cM[p]\setminus p)/a_1= (\cM[p]\setminus a_1)/p .\]
\end{lemma}
\begin{proof}
The first equality is direct. The second one follows from Lemma~\ref{lem:circinseparable}. Indeed, every cocircuit of $( \cM[p]/a_1)\setminus p$ corresponds to a cocircuit~$C$ of~$\cM[p]$ with $C(a_1)=0$ and $C(p)\neq 0$. By Lemma~\ref{lem:circinseparable}, the
values of $C$ on $e\notin\{a_1,p\}$ coincide with the values of $C'$ on $e\notin\{a_1,p\}$, where $C'$ is a cocircuit of $\cM[p]$ with $C'(a_1)\neq 0$ and $C'(p)=0$. That is, $C'$ corresponds to a cocircuit of $( \cM[p]/p)\setminus a_1$.
\end{proof}

We are now ready to prove the proposition.

\begin{proposition}\label{prop:allquotientsofle}
Let $\cM$ be a uniform oriented matroid of rank $\rr$ on a ground set~$E$, and let 
$\cM[p]$ be the lexicographic extension of $\cM$ by $p=[a_1^{\ep_1},a_2^{\ep_2},\dots,a_\rr^{\ep_\rr}]$.
Then
\begin{align}
\cM[p]/p\ &\stackrel{\varphi}{\simeq}\ (\cM/a_1)[a_2^{-\ep_1\ep_2},\dots,a_\rr^{-\ep_1\ep_\rr}],\label{eq:Mmodp}\\
\cM[p]/a_i\ &=\ (\cM/a_i)[a_1^{\ep_1},\dots,a_{i-1}^{\ep_{i-1}},a_{i+1}^{\ep_{i+1}},\dots,a_\rr^{\ep_\rr}], \text{ and }\label{eq:Mmoda}\\
\cM[p]/e\ &=\ (\cM/e)[a_1^{\ep_1},a_2^{\ep_2},\dots,a_{\rr-1}^{\ep_{\rr-1}}];\label{eq:Mmode}
\end{align}
where $e\in E$ is any element different from $p$ and any $a_i$. The isomorphism $\varphi$ in~\eqref{eq:Mmodp} is $\varphi(e)=e$ for all $e\in E\setminus \{p,a_1\}$ and $\varphi(a_1)=[a_2^{-\ep_1\ep_2},\dots,a_\rr^{-\ep_1\ep_\rr}]$; the latter is the extending element.
\end{proposition}
\begin{proof}
The proof of \eqref{eq:Mmoda} and \eqref{eq:Mmode} is direct 
just by observing the signature of $p$ in~$\cM[p]$.

To prove \eqref{eq:Mmodp}, observe that $(\cM[p]/p)\setminus a_1=\cM/ a_1$ by Lemma~\ref{lem:contractdeletele}. Therefore, we only need to prove that the signature of the extension of $(\cM[p]/p)\setminus a_1$ by $a_1$ coincides with that of the lexicographic extension by $[a_2^{-\ep_1\ep_2},\dots,a_\rr^{-\ep_1\ep_\rr}]$. 
That is, let $C\in \co(\cM[p])$ be a cocircuit of $\cM[p]$ with $C(p)=0$ and $C(a_1)\neq 0$ and let $k>1$ be minimal with $C(a_k)\neq 0$. We want to see that $C(a_1)=-\ep_1\ep_kC(a_k)$. 

Because $a_1$ and $p$ are $(-\ep_1)$-inseparable, Lemma~\ref{lem:circinseparable} yields a cocircuit $C'\in \co(\cM[p])$ with $C'(p)=-\ep_1C(a_1)$ and $C'(a_1)=0$ and such that $k$ is minimal with $C'(a_k)\neq 0$. Moreover $C'(a_k)=C(a_k)$ and by the signature of the lexicographic extension $C'(p)=\ep_kC'(a_k)=\ep_kC(a_k)$. The claim follows from comparing both expressions for $C'(p)$.
\end{proof}

The most interesting case is \eqref{eq:Mmodp}. If $\cM$ is realized by~$\vv V$ and $\vv V\cup \{\vv p\}$ realizes the lexicographic extension of $\cM(\vv V)$ by $\vv p=[\vv v_1^{\ep_1},\vv v_2^{\ep_2},\dots,\vv v_\rr^{\ep_\rr}]$. 
The intuition behind the fact that $\cM(\vv V)[\vv p]/\vv p\simeq \cM(\vv V/\vv v_1)[\vv v_2^{-\ep_1\ep_2},\dots,\vv v_\rr^{-\ep_1\ep_\rr}]$ is that every hyperplane that goes through $\vv p$ and not through $\vv v_1$ looks very much like some hyperplane that goes through~$\vv v_1$ and not through $\vv p$. 
If $\ep_1=+$, then $\vv v_1$ and~$\vv p$ are very close, which means that when we perturb a hyperplane ${\vvh H}$ that goes through $\vv v_1$ with $\vv p$ in ${\vvh H}^+$ to its analogue ${\vvh H}'$ through $\vv p$, then $\vv v_1$ lies in~${{\vvh H}'}^-$ and the remaining elements are on the same side of ${\vvh H}'$ as they were of  ${\vvh H}$. 
On the other hand, if $\ep_1=-$, then $\vv v_1$ and~$-\vv p$ are very close, and to perturb ${\vvh H}$ to ${\vvh H}'$, one must also switch the sign of $\vv v_1$. Hence if $\vv p$ was in ${\vvh H}^+$, then $\vv v_1$ is in ${{\vvh H}'}^-$.

This can be checked in the example of Figure~\ref{fig:le}. For each $i\neq  4$, let ${\vvh H_i}$ be the hyperplane that goes through $\vv x_i$ and $\vv x_4$ oriented with $\vv p\in {\vvh H_i}^+$, and let ${\vvh H_i}'$ be the hyperlane spanned by $\vv x_i$ and $\vv p$ with $\vv x_4\in {\vvh H_i'}^-$. 
Observe how for each $j\notin\{ i, 4\}$, $\vv x_j\in {\vvh H_i'}^+$ if and only if $\vv x_j\in {\vvh H_i}^+$. Moreover, the hyperplane ${\vvh H_4}$ spanned by $\vv x_4$ and $\vv x_1$ with $\vv p\in {\vvh H_4}^+$ can be perturbed to the hyperlane ${\vvh H_4}'$ spanned by $\vv x_4$ and $\vv p$ and fulfills $\vv x_1\in {\vvh H_4'}^+$.
\\

Our next results are not directly necessary for proving any of the results presented later. However, they help to understand lexicographic extensions and they are computationally useful (for example, we used them in the programs that compute the exact number of combinatorial types in Section~\ref{sec:exact}).
First, we present two propositions that give a complete description of $\ci(\cM[p])$ and $\co(\cM[p])$ in terms of $\ci(\cM)$ and $\co(\cM)$, respectively. They can also be deduced from results in \cite{LasVergnas1978} and \cite{Todd1985}.

\begin{proposition}\label{prop:cocircuitslexicographic}
Let $\cM$ be a uniform oriented matroid of rank $\rr$, and let $\cM[p]$ be a the lexicographic extension of $\cM$ by $p=[a_1^{\ep_1},a_2^{\ep_2},\dots,a_\rr^{\ep_\rr}]$.

For each cocircuit $C$ in $\co(\cM)$ let $k_C$ be minimal with $C(a_{k_C})\neq 0$, and define the cocircuits $C'_i$ for $0\leq i < k_{C}$ as follows
\begin{enumerate}[(i)]
 \item For $i=0$
\begin{equation*}
C'_0= 
\begin{cases} (C^+\cup p,C^-)& \text{if $\ep_{k_C}C(a_{k_C})=+$,}
\\
(C^+,C^-\cup p)& \text{if $\ep_{k_C}C(a_{k_C})=-$.}
\end{cases}
\end{equation*}
\item\label{it:lecocirctype2} For all $0<i<{k_C}$
\begin{equation*}
C'_i= 
\begin{cases} (C^+\cup a_i,C^-)& \text{if $\ep_i\ep_{k_C}C(a_{k_C})=-$,}
\\
(C^+,C^-\cup a_i)& \text{if $\ep_i\ep_{k_C}C(a_{k_C})=+$.}
\end{cases}
\end{equation*}
\end{enumerate}
Then $\co(\cM[p])=\set{C'_i}{C\in\co(\cM)\text{ and }0\leq i<k_C}$.
\end{proposition}
\begin{proof}
 The cocircuits of class $C'_0$ are defined by the signature of the extension. The description of the cocircuits $C'_i$ with $0<i<{k_C}$, follows from
\[\cM[a_1^{\ep_1},\dots,a_\rr^{\ep_\rr}]/\{p,a_1,\dots,a_{j-1}\}\simeq(\cM/\{a_1,\dots,a_{j}\})[a_{j+1}^{-\ep_{j}\ep_{j+1}},\dots,a_\rr^{-\ep_{j}\ep_\rr}],\]
which is a corollary of Proposition~\ref{prop:allquotientsofle}.
\end{proof}

\begin{proposition}
Let $\cM$ be a uniform oriented matroid of rank $\rr$, and let $\cM[p]$ be a the lexicographic extension of $\cM$ on $p=[a_1^{\ep_1},a_2^{\ep_2},\dots,a_\rr^{\ep_\rr}]$.

For each circuit $X$ in $\ci(\cM)$ let $l_X$ be minimal with $X(a_{l_X})=0$. We define the circuits $X'_i$ for $0< i < l_{X}$ as follows

\begin{itemize}
 \item $X'_i(a_i)=0$,
 \item $X'_i(p)=\ep_iX(a_i)$,
 \item $X'_i(a_j)=-\ep_j\ep_iX(a_i)$ for $j< i$ and
 \item $X'_i(e)=X(e)$ otherwise.
\end{itemize}

Then $\ci(\cM[p])=\ci(\cM)\cup\set{X'_i}{X\in\ci(\cM)\text{ and }0< i < l_{X}}$.
\end{proposition}
\begin{proof}
Circuits with $X(p)=0$ are clear, since they are circuits of $\cM$. Circuits with $X(p)\neq 0$ and $X(a_1)=0$ are also clear, by Lemma~\ref{lem:circinseparable}. 
Finally, using the definition of contraction we see that for each circuit $X\in\ci(\cM[p])$ with $X(p)\neq 0$ and $X(a_1)\neq 0$ there is a circuit $X'\in \ci(\cM[p]/p)$ such that $X(e)=X'(e)$ for $e\neq p$ and $X(p)=-\ep_1X(a_1)=-\ep_1X'(a_1)$ (because $p$ and $a_1$ are $(-\ep_1)$-inseparable). The result follows by induction using Proposition~\ref{prop:allquotientsofle}.
\end{proof}

We end this section with a characterization of inseparable elements in lexicographic extensions.

\begin{proposition}
 Let $\cM$ be a uniform oriented matroid of rank $\rr$ on a ground set $E$ and let $\cM[p]$ be a lexicographic extension of $\cM$ by $p=[a_1^{\ep_1},a_2^{\ep_2},\dots, a_\rr^{\ep_\rr}]$. Then, 
\begin{itemize}
 \item $p$ and $a_k$ 
are $\alpha$-inseparable if and only if $\alpha=-\ep_k$ and for all $1\leq i\leq  k-1$, $a_k$ and $a_i$ are $(-\ep_k \ep_i)$-inseparable in $\cM/\{a_1,\dots,a_{i-1}\}$. 
 \item $p$ and $y\in E\setminus\{a_1,\dots,a_\rr\}$ are $\alpha$-inseparable if and only if for all $1\leq i\leq \rr$, $y$ and $a_i$ are $(\alpha \ep_i)$-inseparable in $\cM/\{a_1,\dots,a_{i-1}\}$.
 \item $a_k$ and $a_j$ with $k<j$ 
 are $\alpha$-inseparable if and only if $\alpha=\ep_k\ep_{j}$, they are $(\ep_k\ep_{j})$-inseparable in $\cM$ and for all $k+1\leq i\leq j-1$ $a_j$ and $a_i$ are $(-\ep_{j}\ep_i)$-inseparable in $\cM/\{a_1,\dots,a_{i-1}\}$. 
 \item $a_k$ and $y\in E\setminus\{a_1,\dots,a_\rr\}$ are $\alpha$-inseparable if and only if they are $\alpha$-inseparable in $\cM$ and for all $k+1\leq i\leq \rr$, $y$ and $a_i$ are $(-\alpha\ep_k\ep_i)$-inseparable in $\cM/\{a_1,\dots,a_{i-1}\}$. 
 \item $x,y\in E\setminus\{a_1,\dots,a_\rr\}$ are $\alpha$-inseparable if and only if they are $\alpha$-inseparable in~$\cM$. 
\end{itemize}
\end{proposition}
\begin{proof}
By definition, $p$ and $a_1$ are $(-\ep_1)$-inseparable.

The next question is whether there can be any other element $a_k$, $k>1$ that is inseparable with $p$. First of all, observe that any cocircuit $C$ with $C(a_i)=0$ for $i<k$ and $C(a_k)\neq 0$ fulfills $C(p)=\ep_kC(a_k)$, so if~$p$ and~$a_k$ are inseparable, they must be $(-\ep_k)$-inseparable. 
Moreover, let $1\leq i \leq k-1$, and let $C$ be a cocircuit of $\cM[p]$ with $C(p)\neq 0$, $C(a_k)\neq 0$, $C(a_i)\neq 0$ and $C(a_j)= 0$ for $j<i$. Then $C(p)=\ep_iC(a_i)$, so if~$a_k$ and~$p$ are $(-\ep_k)$-inseparable then $\ep_iC(a_i)=\ep_kC(a_k)$. 
Hence, $a_k$ and $a_i$ must be $(-\ep_k \ep_i)$-inseparable in $\cM/\{a_1,\dots,a_{i-1}\}$. Observe that we checked all possible cocircuits $C$ with $C(p)\neq 0$ and $C(a_k)\neq 0$, and therefore these conditions are also sufficient.

The same reasoning works for an element $y\in E\setminus\{a_1,\dots,a_\rr\}$. Observe that in this case, we are treating $y$ as it was $a_{r+1}$.  

Now we see when a pair of elements $x,y\in E$ are $\alpha$-inseparable. Each cocircuit of $\cM[p]$ with $C(p)\neq0$ corresponds to a cocircuit of $\cM$, so $x$ and $y$ must be $\alpha$-inseparable in $\cM$. For circuits $C$ with $C(p)=0$, observe that they are circuits of $\cM[p]/p$, and by Proposition~\ref{prop:allquotientsofle}, $\cM[p]/p\simeq \cM/a_1[a_2^{-\ep_1\ep_2},\dots, a_r^{-\ep_1\ep_r}]$. 
Then the result follows by induction on $r$. The base case $r=1$ is easy and left to the reader.
\end{proof}

\chapter{An update on the Sewing Construction}\label{ch:shemer}

The sewing construction was introduced by Shemer in 1982. In his classical paper~\cite{Shemer1982} he derived several interesting results about neighborly polytopes. In particular he used the \defn{Sewing construction} to give lower bounds for the number of combinatorial types of neighborly polytopes with $n$ vertices in dimension $d$. His bounds were surprisingly high and until Alon's bounds from 1986~\cite{Alon1986} these were the best lower bounds for combinatorial types of polytopes. 
 
His starting point is Gr\"unbaum's \defn{beneath-beyond} technique~\cite[Section 5.2]{GruenbaumEtal2003} (see also the formulation in~\cite{AltshulerShemer1984} and the version for oriented matroids in~\cite[Proposition 9.2.2]{OrientedMatroids1993}). This technique is based on the fact that, given a $d$-polytope $\vv P$ and a point $\vv p\in \RR^d$, the face lattice of $\conv(\vv P\cup \{\vv p\})$ can be read from the relative position (beneath, beyond or on) of $\vv p$ with respect to the hyperplanes defining facets of $\vv P$.
The key for applying this technique is to find subsets of facets of $\vv P$ for which one can prove the existence of certain point $\vv p$ beyond their supporting hyperplanes.

The sewing construction finds such a subset of facets from a flag of faces of $\vv P$. The Sewing Theorem (see Theorem~\ref{thm:shemersewing}) showed that if~$\vv P$ is an even dimensional neighborly polytope and the flag has certain properties (it is a universal flag), then $\conv (\vv P \cup \{\vv p\})$ is still neighborly.
Barnette's facet splitting construction~\cite{Barnette1981} can be seen as a polar version of sewing (with slightly different flags).

In 2000, Bisztriczky proved that the sewing construction can also be used to construct odd dimensional simplicial neighborly polytopes~\cite{Bisztriczky2000}. Additionally, Trelford and Vigh used the structure of vertex figures of neighborly polytopes obtained with the sewing construction to compute their face lattice in \cite{TrelfordVigh2011}. 
A generalization of the sewing construction, called \defn{$A$-sewing}, was used by Lee and Menzel to construct non-simplicial polytopes~\cite{LeeMenzel2010}.
\\

Observing that the sewing construction is in fact a lexicographic extension whose signature is determined by the corresponding face flag, our main result in this chapter is the Extended Sewing Theorem~\ref{thm:extshemersewing}, which is a generalization of Shemer's Sewing Theorem that:
\begin{enumerate}
 \item \emph{Allows to extend any neighborly oriented matroid}: it can also be used to extend non-realizable oriented matroids. This property is exploited by Theorem~\ref{thm:nonrealizable} in Chapter~\ref{ch:thethm} to construct many non-realizable oriented matroids. Despite the matroid formulation, it can still be used to construct polytopes since it preserves realizability. 

 \item \emph{Works for arbitrary rank}: as in Bisztriczky's proof of the Sewing Theorem in~\cite{Bisztriczky2000}, this construction does not put constraints on the parity of the rank.
 \item \emph{Uses a larger family of face flags}: Shemer used the sewing construction on universal flags, which are subflags of the flags used by Barnette with the facet splitting technique~\cite{Barnette1981}. Our construction is extended to any flag that contains a universal subflag (already suggested by Shemer in~\cite[Remark 7.4]{Shemer1982}).
 Proposition~\ref{prop:uniqueflags} shows that these are all the possible flags for constructing neighborly oriented matroids of odd rank (respectively neighborly polytopes of even dimension).
\end{enumerate}

The language of oriented matroids not only allows for working with non-realizable matroids, but makes the results easier to state and prove.
Nevertheless, in this chapter we use the letter $\cP$ for oriented matroids to reinforce the idea that all the following results translate directly to polytopes.

\section{The Sewing Theorem}\label{sec:shemer}
Let $\cP$ be an acyclic oriented matroid on a ground set $E$, and let $F\subset E$ be a facet of~$\cP$. That is, there exists a cocircuit $C_{ F}$ of $\cP$ such that $C_{ F}(e)=0$ if $e\in  F$ and $C_{F}(e)=+$ otherwise. Consider a single element extension of~$\cP$ by~$p$ 
with signature~$\gs_p$. We say that $p$ is 
\defn{beneath}~$ F$ if $\gs_p(C_{ F})=+$, that $p$ is \defn{beyond}~$ F$ when $\gs_p(C_{ F})=-$, and that $p$ is \defn{on} $ F$ if $\gs_p(C_{ F})=0$.
Moreover, we say that $p$ lies \defn{exactly beyond}\index{exactly beyond} a set of facets $\fF$ if it lies beyond all facets in $\fF$ and beneath all facets not in $\fF$. An example is shown in Figure~\ref{fig:beneathbeyond}.\index{beneath-beyond}

\iftoggle{bwprint}{%
\begin{figure}[htpb]
\begin{center}
\includegraphics[width=.4\textwidth]{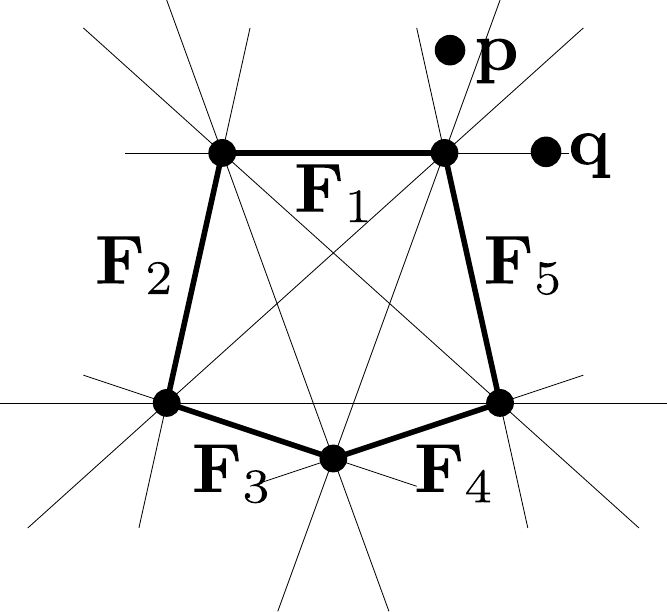}
\end{center}
 \caption[Beneath-beyond example]{The point $\vv p$ is beyond $\vv F_1$ and $\vv F_5$ and beneath $\vv F_2$, $\vv F_3$ and $\vv F_4$. Hence, $\vv p$ is exactly beyond the set $\{\vv F_1,\vv F_5\}$. The point $\vv q$ is on $\vv F_1$, beyond $\vv F_5$ and beneath $\vv F_2$, $\vv F_3$ and $\vv F_4$. }
 \label{fig:beneathbeyond}
\end{figure}
}{%
\begin{figure}[htpb]
\begin{center}
\includegraphics[width=.4\textwidth]{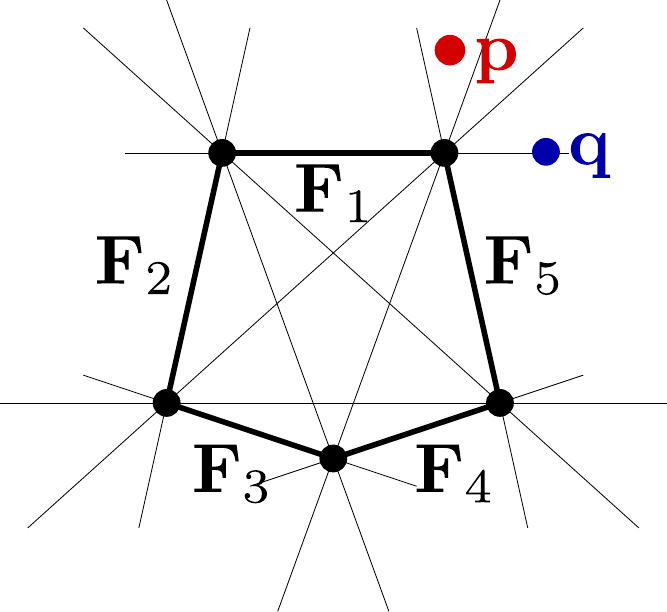}
\end{center}
 \caption[Beneath-beyond example]{The point $\vv p$ is beyond $\vv F_1$ and $\vv F_5$ and beneath $\vv F_2$, $\vv F_3$ and $\vv F_4$. Hence, $\vv p$ is exactly beyond the set $\{\vv F_1,\vv F_5\}$. The point $\vv q$ is on $\vv F_1$, beyond $\vv F_5$ and beneath $\vv F_2$, $\vv F_3$ and $\vv F_4$. }
 \label{fig:beneathbeyond}
\end{figure}
}

\begin{lemma}[{\cite[Proposition 9.2.2]{OrientedMatroids1993}}]\label{lem:benbey} Let $\tilde \cP$ be a single element extension of $\cP$ with signature $\sigma$. Then the values of $\sigma$ on the facet cocircuits of $\cP$ determine the whole face lattice of $\tilde\cP$.
\end{lemma}

A \defn{flag}\index{flag} of an acyclic oriented matroid $\cP$ is a strictly increasing sequence of proper faces $F_1\subset F_2 \subset \dots \subset F_k$. We say that a flag $\cF$ is a \defn{subflag}\index{flag!subflag} of $\cF'$ if for every $F$ that belongs to $\cF$, $F$ also belongs to $\cF'$. 
\\

Fixed a flag $\cF=\{F_j\}_{j=1}^k$ of $\cP$, let $\fF_j$\index{$\fF_j$} be the set of facets of $\cP$ that contain $F_j$, and let $\sew(\cF):=\fF_1\setminus(\fF_2\setminus(\dots\setminus\fF_k)\dots)$\index{$\sew(\cF)$}, so that
\begin{equation*}
\sew(\cF)= 
\begin{cases} (\fF_1\setminus\fF_2) \cup (\fF_3\setminus\fF_4)\cup \dots \cup (\fF_{k-1}\setminus\fF_k)& \text{if $k$ is even,}
\\
(\fF_1\setminus\fF_2) \cup (\fF_3\setminus\fF_4)\cup \dots \cup \fF_k &\text{if $k$ is odd.}
\end{cases}
\end{equation*}

Given a polytope $\vv P$ with a flag $\cF=\vv F_1\subset \vv F_2 \subset \dots \subset \vv F_k$, Shemer proved that there always exists an extension exactly beyond $\sew(\cF)$ (\cite[Lemma 4.4]{Shemer1982}), and called this extension sewing onto the flag. We will show that there is a lexicographic extension that realizes the desired signature.

\begin{definition}[Sewing onto a flag]\label{def:sewing}\index{sewing}
Let $\cF=\{F_j\}_{j=1}^k$ be a flag of an acyclic matroid $\cP$ on a ground set $E$. We extend it with $F_{k+1}=E$ and define $U_j=F_j\setminus F_{j-1}$. 
We say that $p$ is \emph{sewn} onto $\cP$ through~$\cF$, if $\cP[p]$ is a lexicographic extension of $\cP$ by
\begin{displaymath}p=[F_1^+,U_2^-,U_3^+,\dots,U_{k+1}^{(-1)^{k}}],\end{displaymath}
where these sets represent their elements in any order. 
Put differently, the lexicographic extension by $p$ is defined by $p=[a_1^{\ep_1},a_2^{\ep_2},\dots,a_{n}^{\ep_{n}}]$, where $a_1, \dots, a_{n}$ are the elements in $F_{k+1}$ sorted such that 
\begin{itemize}
\item if there is some $m$ such that $a_i\in F_m$ and $a_j\notin F_m$, then $i<j$;
\item if the smallest $m$ such that $a_j\in F_m$ is odd, then $\ep_j=+$; and $\ep_j=-$ otherwise.
\end{itemize}
We use the notation $\cP[\cF]$ to designate the extension $\cP[p]$ when $p$ is sewn onto~$\cP$ through~$\cF$.
\end{definition}

For example, if $\cP$ has rank $4$ and $F_1=\{a_1,a_2\}$ and $F_2=\{a_1,a_2,a_3,a_4\}$ are the elements of two faces of $\cP$, then the lexicographic extensions by $[a_1^+,a_2^+,a_3^-,a_4^-,a_5^+]$, $[a_2^+,a_1^+,a_3^-,a_4^-,a_6^+]$ or $[a_2^+,a_1^+,a_4^-,a_3^-,a_6^+]$ are extensions by an element sewn through~the flag $F_1\subset F_2$.

\iftoggle{bwprint}{%
\begin{figure}[htpb]
\centering
 \subbottom[$\vv P$]{\quad\label{sfig:flagsewing0}\includegraphics[width=.23\linewidth]{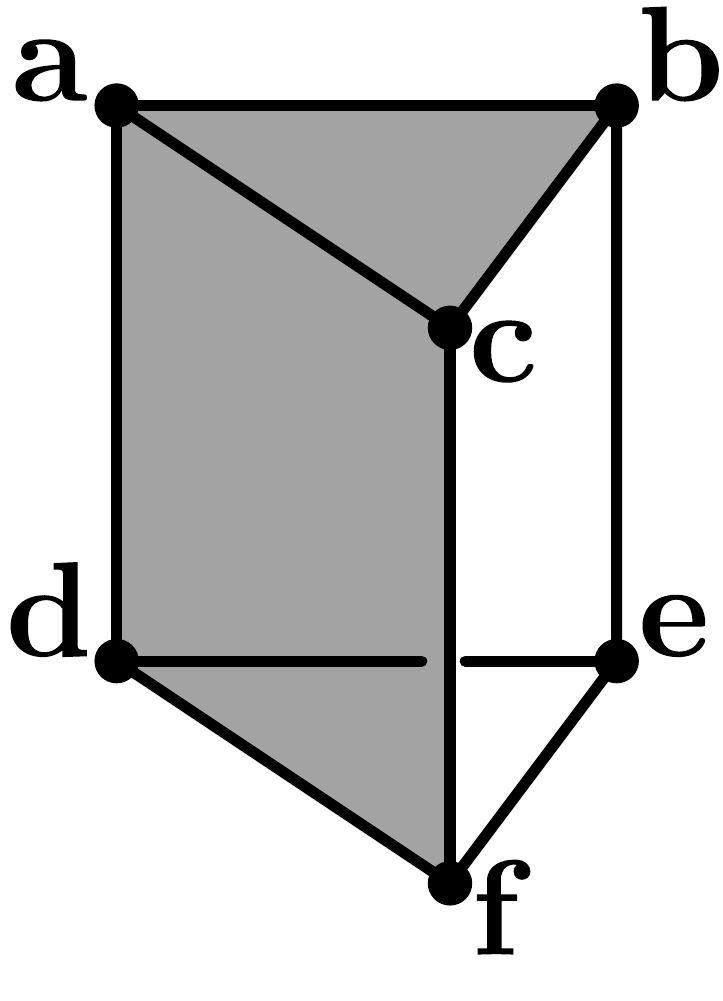}\quad}\quad 
 \subbottom[$\vv P\bracket{\vv c^+}$]{\quad\label{sfig:flagsewing1}\includegraphics[width=.23\linewidth]{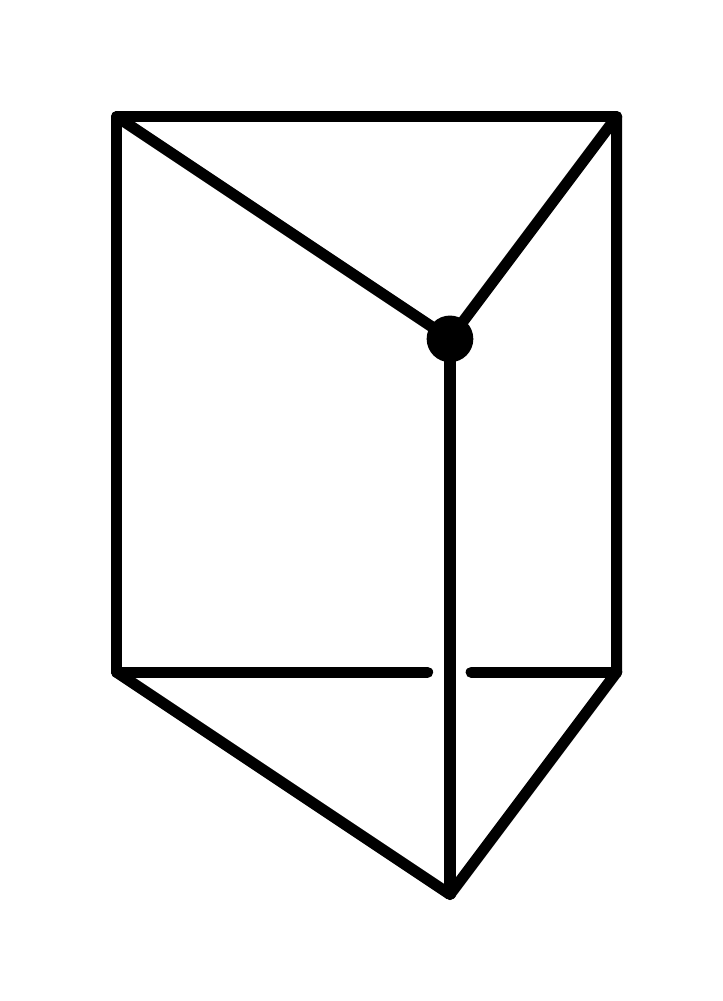}\quad}\quad
 \subbottom[$\vv P\bracket{\vv c^+,\vv b^-}$]{\quad\label{sfig:flagsewing2}\includegraphics[width=.23\linewidth]{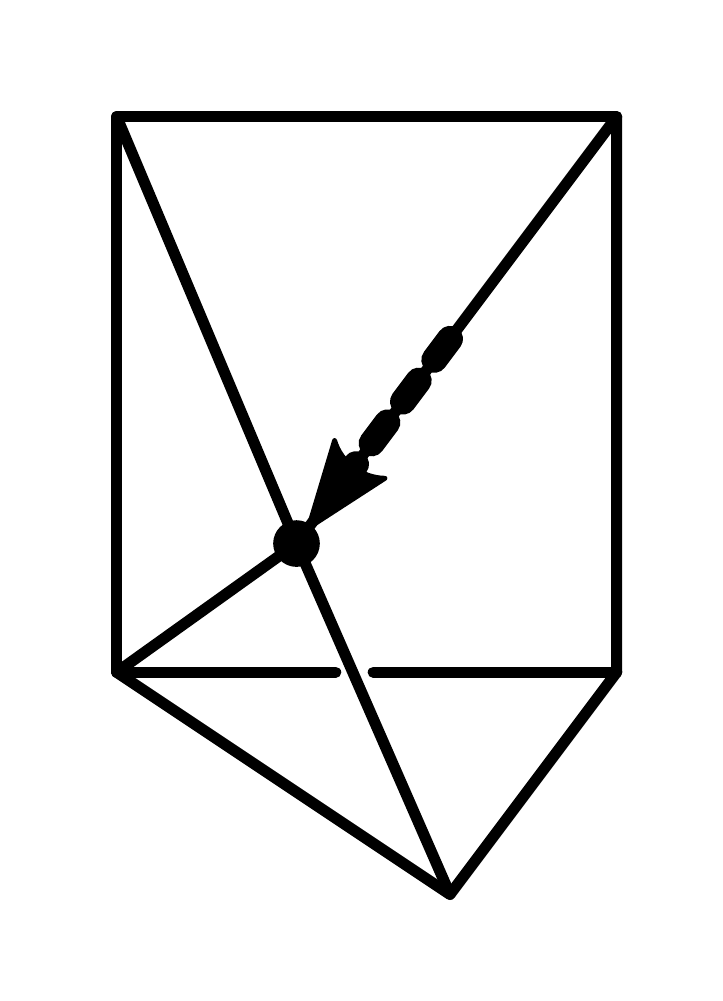}\quad} \\ 
 \subbottom[$\vv P\bracket{\vv c^+,\vv b^-,\vv a^+}$]{\quad\label{sfig:flagsewing3}\includegraphics[width=.23\linewidth]{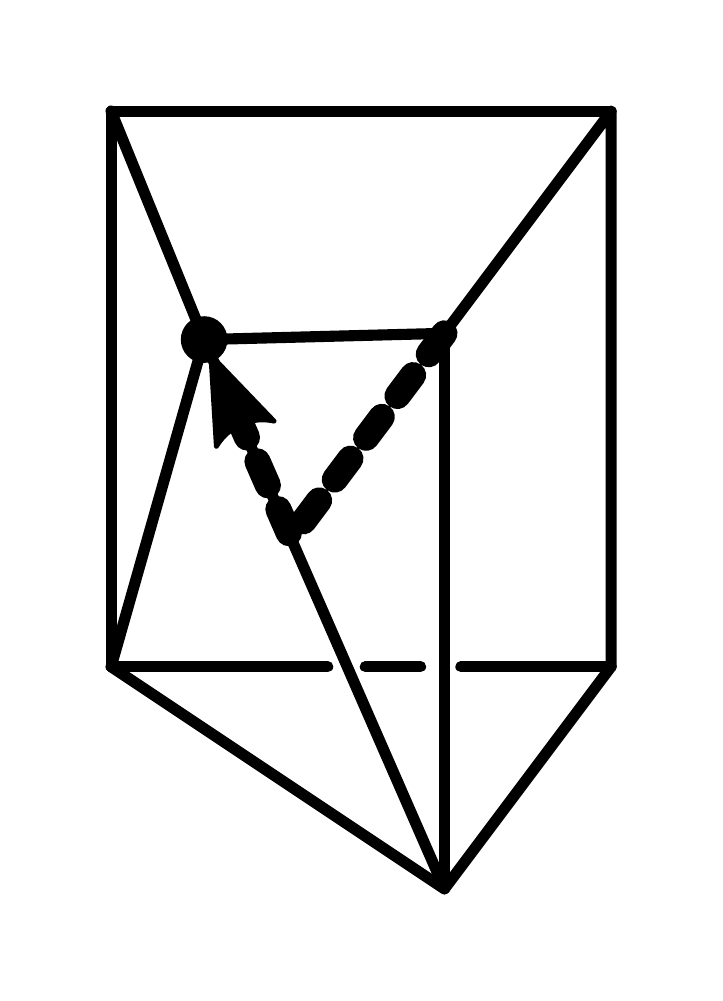}\quad}\quad
 \subbottom[$\vv P\bracket{\vv c^+,\vv b^-,\vv a^+,\vv d^-}$]{\quad\label{sfig:flagsewing4}\includegraphics[width=.23\linewidth]{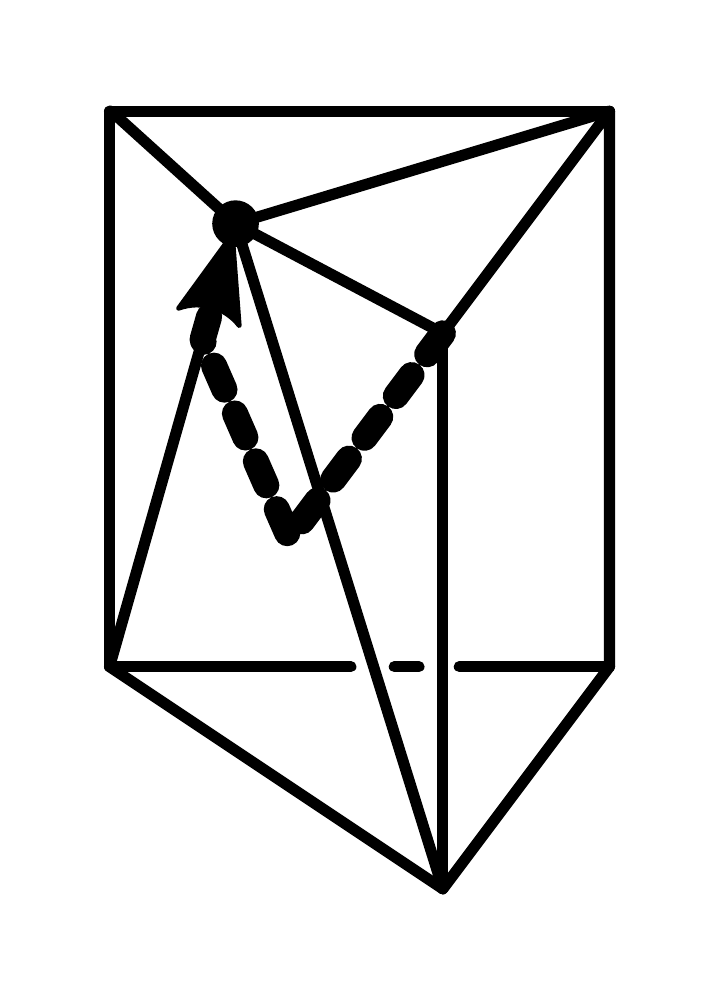}\quad}\quad
 \subbottom[$\vv P\bracket{\cF}$]{\quad\label{sfig:flagsewing5}\includegraphics[width=.23\linewidth]{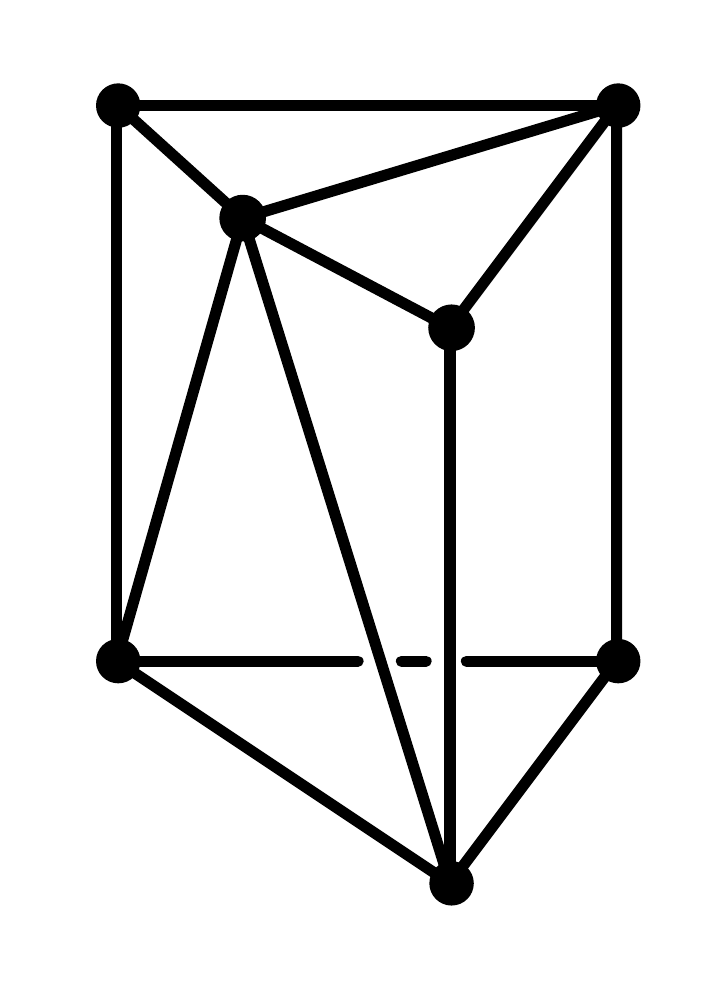}\quad}
 \caption[Sewing example]{A polytope $\vv P:=\conv\{\vv a,\vv b,\vv c,\vv d,\vv e,\vv f\}$. Sewing onto the flag $\cF=\{\vv c\}\subseteq\{\vv c,\vv b\}\subseteq \{\vv c,\vv b,\vv a\}$. Shaded facets in \subcaptionref{sfig:flagsewing0} correspond to~$\sew(\cF)$.}
 \label{fig:sewingonflags}
\end{figure}
}{%
\begin{figure}[htpb]
\centering
 \subbottom[$\vv P$]{\quad\label{sfig:flagsewing0}\includegraphics[width=.23\linewidth]{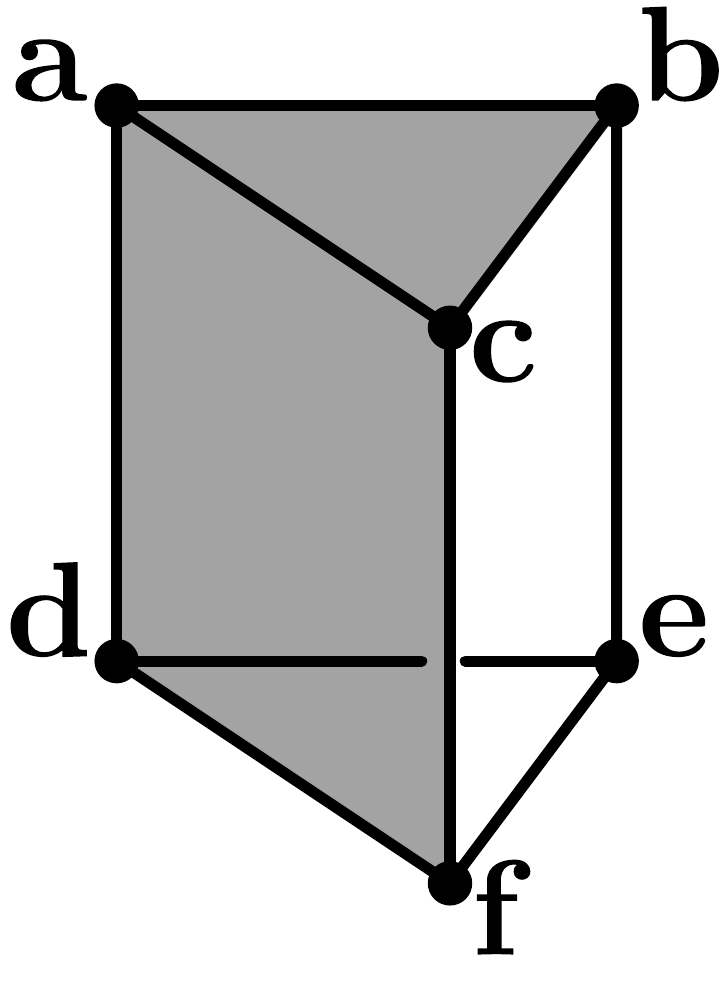}\quad}\quad 
 \subbottom[$\vv P\bracket{\vv c^+}$]{\quad\label{sfig:flagsewing1}\includegraphics[width=.23\linewidth]{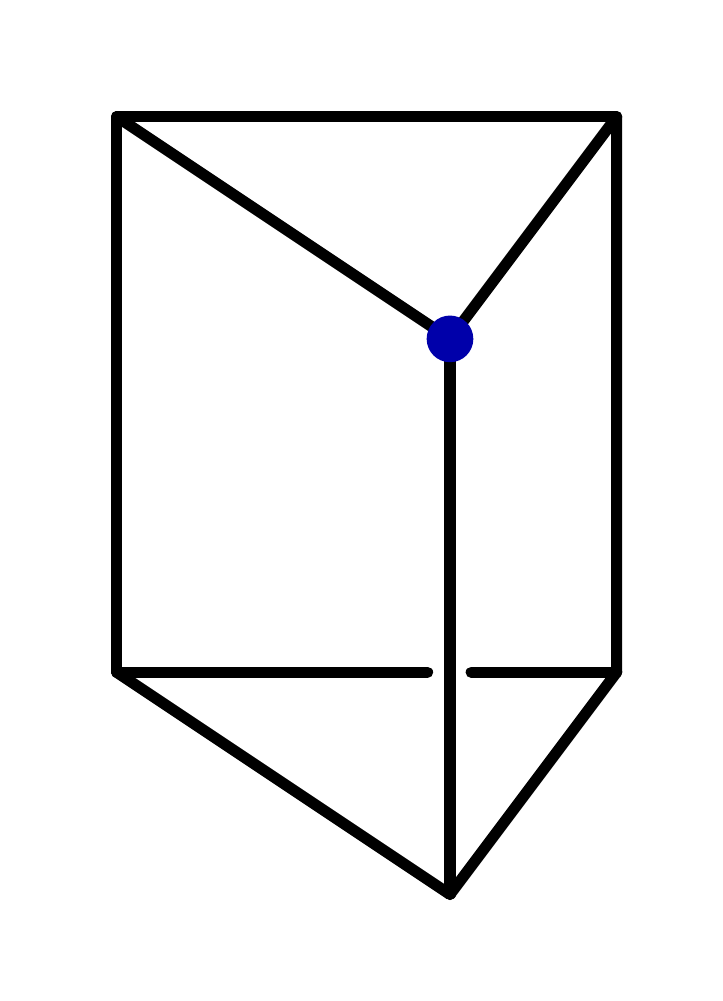}\quad}\quad
 \subbottom[$\vv P\bracket{\vv c^+,\vv b^-}$]{\quad\label{sfig:flagsewing2}\includegraphics[width=.23\linewidth]{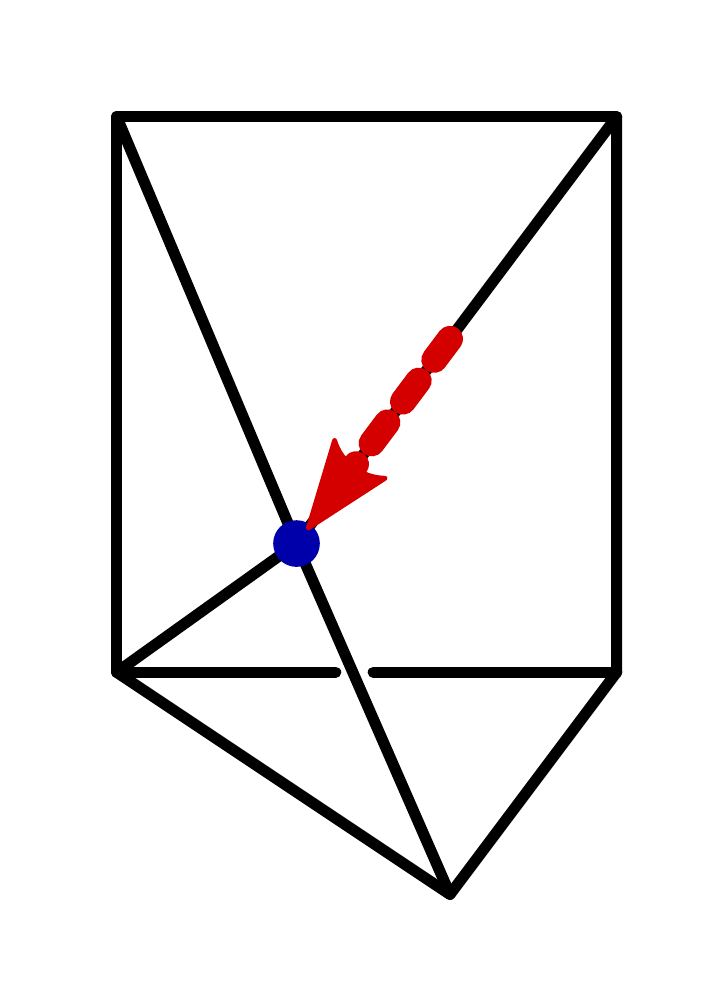}\quad} \\ 
 \subbottom[$\vv P\bracket{\vv c^+,\vv b^-,\vv a^+}$]{\quad\label{sfig:flagsewing3}\includegraphics[width=.23\linewidth]{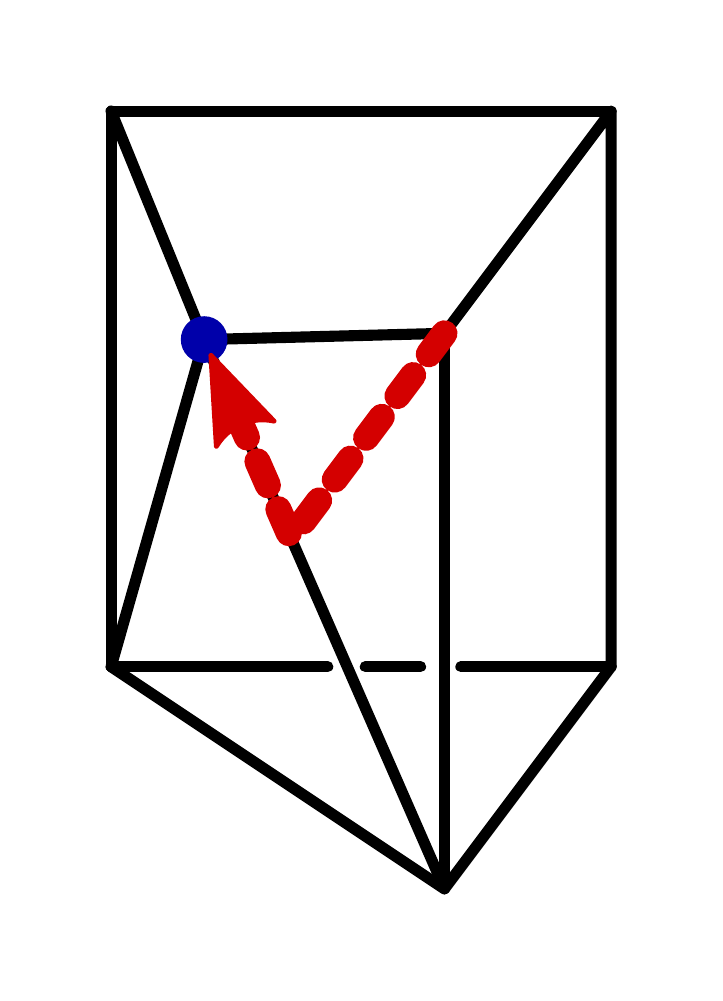}\quad}\quad
 \subbottom[$\vv P\bracket{\vv c^+,\vv b^-,\vv a^+,\vv d^-}$]{\quad\label{sfig:flagsewing4}\includegraphics[width=.23\linewidth]{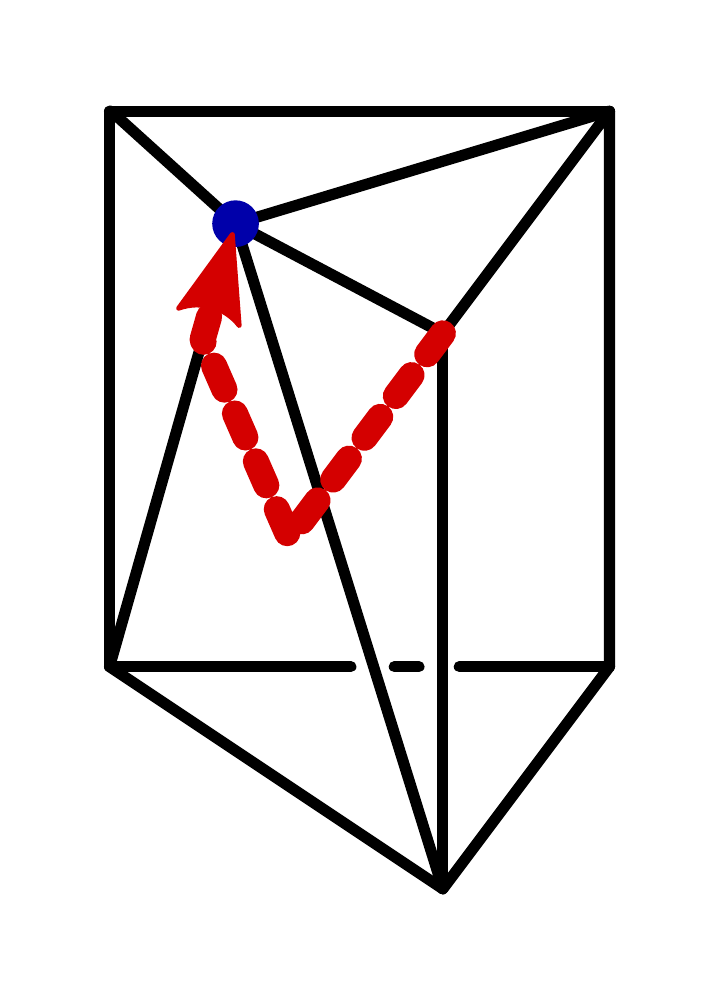}\quad}\quad
 \subbottom[$\vv P\bracket{\cF}$]{\quad\label{sfig:flagsewing5}\includegraphics[width=.23\linewidth]{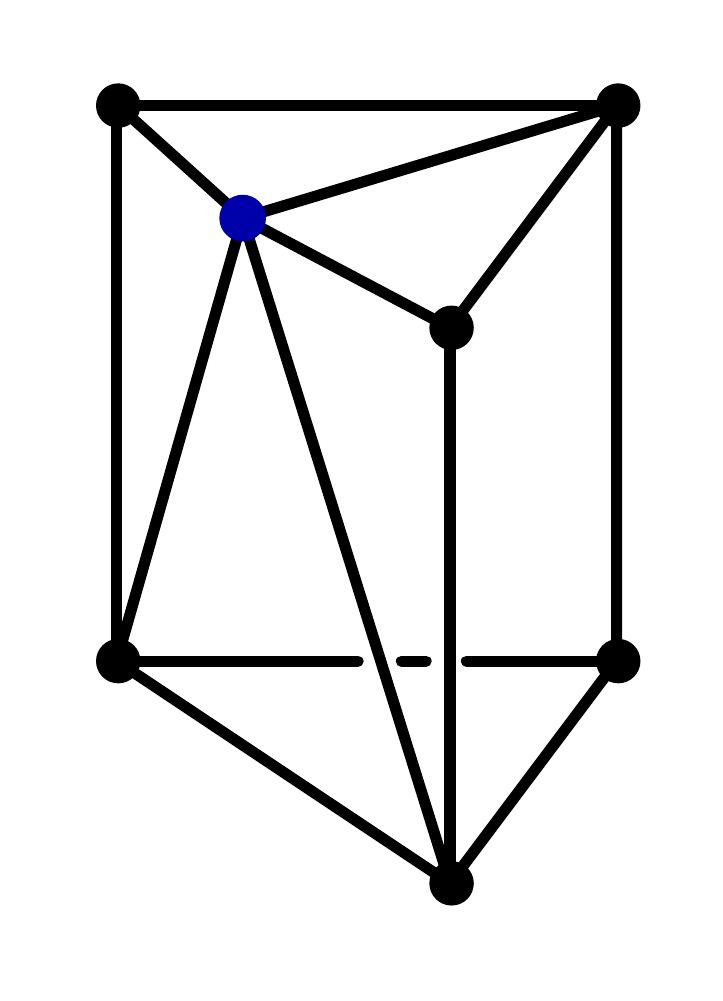}\quad}
 \caption[Sewing example]{A polytope $\vv P:=\conv\{\vv a,\vv b,\vv c,\vv d,\vv e,\vv f\}$. Sewing onto the flag $\cF=\{\vv c\}\subseteq\{\vv c,\vv b\}\subseteq \{\vv c,\vv b,\vv a\}$. Shaded facets in \subcaptionref{sfig:flagsewing0} correspond to~$\sew(\cF)$.}
 \label{fig:sewingonflags}
\end{figure}
}

\begin{example}
Figure~\ref{fig:sewingonflags} shows a point configuration $\vv A=\{\vv a, \vv b, \dots, \vv f\}$ whose convex hull is a polytope $\vv P:=\conv(\vv A)$. Its faces $F_1:=\{\vv c\}$, $F_2:=\{\vv c,\vv b\}$ and $F_3:=\{\vv c,\vv b,\vv a\}$ form a flag $\cF:=F_1\subset F_2\subset F_3$. Then $\fF_1$ contains the facets $\{\vv a,\vv b,\vv c\}$,  $\{\vv a,\vv c,\vv d,\vv f \}$ and $\{\vv b,\vv c,\vv e,\vv f \}$; $\fF_2$ the facets $\{\vv a,\vv b,\vv c\}$ and $\{\vv b,\vv c,\vv e,\vv f \}$; and $\fF_3$ is just the facet $\{\vv a,\vv b,\vv c\}$. Therefore, $\sew(\cF)$ consists of the facets $\{\vv a,\vv b,\vv c\}$,  $\{\vv a,\vv c,\vv d,\vv f \}$. The figure shows how a point is sewn through the flag $\{\vv c\}\subseteq\{\vv c,\vv b\}\subseteq \{\vv c,\vv b,\vv a\}$ in the lexicographic extension $\vv P[\vv c^+,\vv b^-,\vv a^+,\vv d^-]$. Observe that this point lies exactly beyond~$\sew(\cF)$. 
\end{example}

In terms of oriented matroids, the definition of $\cP[\cF]$ is ambiguous, since it can represent different oriented matroids. However, the following proposition (together with Lemma~\ref{lem:benbey}) shows that all the extensions $\cP[\cF]$ have the same face lattice. In particular, this implies that there is no ambiguity when $\cP[\cF]$ is neighborly of odd rank, because these are rigid (Theorem~\ref{thm:neigharerigid}).

\begin{proposition}\label{prop:exactlybeyondle}
Let $\cF=\{F_j\}_{j=1}^k$ be a flag of an acyclic oriented matroid~$\cP$. If $\cP[p]$ is the lexicographic extension $\cP[\cF]$, then $p$ lies exactly beyond $\sew(\cF)$.
\end{proposition}
\begin{proof}
Let the lexicographic extension be by $p=[a_1^{\ep_1},a_2^{\ep_2},\dots,a_{n}^{\ep_{n}}]$ with the elements and signs as in Definition~\ref{def:sewing}.
We have to see that, for $1\leq j\leq k$, $p$ lies beneath any facet in $\fF_{j}\setminus\fF_{j+1}$ if $j$ is even, and beyond any facet in $\fF_{j}\setminus\fF_{j+1}$ if $j$ is odd (with the convention $\fF_{k+1}=\emptyset$).

That is, if $\gs$ is the signature of the lexicographic extension and $F$ a facet of $\cP$ defined by a cocircuit $C_F$, we want to see that 
\[\gs(C_F)=\begin{cases}
            +&\text{ if there is an even $j$ such that $F_{j}\subseteq F$ but $F_{j+1}\not\subseteq F$,}\\
            -&\text{ if there is an odd $j$ such that $F_{j}\subseteq F$ but $F_{j+1}\not\subseteq F$;}
           \end{cases}\]
where we use the convention $F_{k+1}=E$, the ground set of $\cP$.
 
In our case, if $F$ is in $\fF_{j}\setminus\fF_{j+1}$ then the first $a_i$ with $C_F(a_i)\neq 0$ belongs to~$F_{j+1}$ and thus $\ep_i=+$ if $j$ is even and $\ep_i=-$ if $j$ is odd. Therefore, since $\gs(C_F)=\ep_iC_F(a_i)=\ep_i$, $\gs(C_F)=+$ (\ie $p$ is beneath $F$) when $j$ is even while $\gs(C_F)=-$ (\ie $p$ is beyond $F$) when $j$ is odd.
\end{proof}

The proof of the following observation is analogous and left to the reader.
\begin{observation}[$A$-sewing]\label{prop:Asewing}\index{sewing!A-sewing}
In~\cite{LeeMenzel2010}, Lee and Menzel proposed the operation of \emph{$A$-sewing}. Given a flag $\cF=\{F_j\}_{j=1}^k$ of a polytope $P$, it allows to find a point on the facets in $\fF_k$, beyond the facets in $\sew(\cF)\setminus\fF_k$, and beneath the remaining facets. 
In our setting, the process of $A$-sewing corresponds to a lexicographic extension by $[F_1^+,U_2^-,U_3^+,\dots,U_{k}^{(-1)^{k-1}}]$. 
In the example of Figure~\ref{fig:sewingonflags}, the polytopes $\vv P\bracket{\vv c^+,\vv b^-}$ 
and $\vv P\bracket{\vv c^+,\vv b^-,\vv a^+}$ 
correspond to $A$-sewing through the flags $\{\vv c\}\subseteq\{\vv c,\vv b\}$ and $\{\vv c\}\subseteq\{\vv c,\vv b\}\subseteq\{\vv c,\vv b,\vv a\}$ respectively.
\end{observation}

\subsection{Sewing onto universal flags}

Shemer's sewing construction starts with a neighborly oriented matroid $\cP$ of rank~$\rd$ with $n$ elements and gives a neighborly oriented matroid $\tilde \cP$ of rank~$\rd$ with $n+1$ elements, provided that $\cP$ has a \defn{universal flag}.

\begin{definition}
Let $\cP$ be a uniform acyclic oriented matroid of rank~$\rd$, and let $m=\ffloor{\rd-1}{2}$.
\begin{enumerate}[(i)]
 \item A face $F$ of $\cP$ is a \defn{universal face} if the contraction $\cP/F$ is neighborly.
 \item A flag $\cF$ of $\cP$ is a \defn{universal flag} if $\cF=\{F_j\}_{j=1}^m$ where each $F_j$ is a universal face with $2j$ vertices.
\end{enumerate}\index{universal flag}
\end{definition}

\begin{remark}
 According to this definition, Lemma~\ref{lem:insareuni} states that every $(-1)$-inseparable pair of elements of a balanced matroid~$\cM$ forms a universal edge of its dual; and that if additionally the corank of~$\cM$ is odd, then $\{x,y\}$ is a universal edge of~$\Gale\cM$ if and only if $x$ and $y$ are $(-1)$-inseparable in $\cM$.
\end{remark}

The most basic example of neighborly polytopes with universal flags are cyclic polytopes, (cf. {\cite[Theorem 3.4]{Shemer1982}} and \cite[Theorem 1.1]{CordovilDuchet1990}).

\begin{proposition}[{\cite[Theorem 3.4]{Shemer1982}}]\label{prop:universalflagsofcyclic}
Let $\cyc{2m}{n}$ be a cyclic polytope of dimension $2m$, with vertices $\vv a_1,\dots,\vv a_n$ labeled in cyclic order. Then $\{\vv a_i,\vv a_{i+1}\}$ for $1\leq i<n$ and $\{\vv a_1,\vv a_{n}\}$ are universal edges of $\cyc{2m}{n}$. If moreover $n>2m+2$, then these are all the universal edges of~$\cyc{2m}{n}$.
\end{proposition}

\begin{remark}\label{rmk:universalflagsofcyclic}
It is not hard to prove that, for any universal edge $E$ of $\cyc{2m}{n}$,
\[\cyc{2m}{n}/E\simeq \cyc{2m-2}{n-2}\] 
 where the isomorphism is such that the cyclic order is preserved. This observation, combined with Proposition~\ref{prop:universalflagsofcyclic}, provides a recursive method to compute universal flags of $\cyc{2m}{n}$. 

Indeed, let $E$ be a universal edge of~$\cyc{2m}{n}$ which, up to symmetry we can take to be $E=\{\vv a_{n-1},\vv a_n\}$. Let $F_1\subset \dots\subset F_ {m-1}$ be a universal flag of~$\cyc{2m-2}{n-2}$. Then the flag \[E\subset E\cup (E\cup F_1)\subset \dots\subset (E\cup F_ {m-1})\] is a universal flag of~$\cyc{2m}{n}$ (here we identify the vertices of~$\cyc{2m}{n}/E$ with the corresponding vertices of~$\cyc{2m}{n}$). If moreover $n>2m+2$, then all universal flags of~$\cyc{2m}{n}$ arise this way. 
\end{remark}

\begin{example}
Consider the cyclic polytope $\cyc{4}{8}$ with cyclically ordered vertices $\{\vv a_1,\dots,\vv a_8\}$ (its affine Gale diagram is shown in Figure~\ref{fig:cyclicGalesewing}). The edge $E:=\{\vv a_7,\vv a_8\}$ is universal. Then $\cyc{4}{8}/E$ is an hexagon whose six edges are universal. They give rise to six universal flags: 
\begin{align*}
\{\vv a_7,\vv a_8\}&\subset \{\vv a_1,\vv a_2,\vv a_7,\vv a_8\},&
\{\vv a_7,\vv a_8\}&\subset \{\vv a_2,\vv a_3,\vv a_7,\vv a_8\},\\
\{\vv a_7,\vv a_8\}&\subset \{\vv a_3,\vv a_4,\vv a_7,\vv a_8\},&
\{\vv a_7,\vv a_8\}&\subset \{\vv a_4,\vv a_5,\vv a_7,\vv a_8\},\\
\{\vv a_7,\vv a_8\}&\subset \{\vv a_5,\vv a_6,\vv a_7,\vv a_8\},&
\{\vv a_7,\vv a_8\}&\subset \{\vv a_1,\vv a_6,\vv a_7,\vv a_8\}.
\end{align*}
\end{example}

\medskip

With these notions, we are ready to present Shemer's Sewing Theorem.
\begin{theorem}[The Sewing Theorem {\cite[Theorem 4.6]{Shemer1982}}]\label{thm:shemersewing}

Let $\vv P$ be a neighborly $2m$-polytope with a universal flag $\cF=\{\vv F_j\}_{j=1}^m$, for some $m\geq1$ and where $\vv F_j=\bigcup_{i=1}^{j}\{\vv x_i,\vv y_i\}$.
Let $\vv P[\cF]$ be the polytope obtained by sewing $\vv p$ onto $\vv P$ through~$\cF$.
 Then,
\begin{enumerate}
\item\label{it:newisneigh} $\vv P[\cF]$ is a neighborly polytope and $\verts (\vv P[\cF])=\verts (\vv P)\cup \{\vv p\}$.
\item\label{it:newunifaces} For all $1\leq j \leq m$, $\vv F_{j-1}\cup\{\vv x_j,\vv p\}$ and $\vv F_{j-1}\cup\{\vv y_j,\vv p\}$ are universal faces of $\vv P[\cF]$. If moreover $j$ is even, then $\vv F_j$ is also a universal face of $\vv P[\cF]$.
\end{enumerate} 
\end{theorem}
\medskip
Combining Remark~\ref{rmk:universalflagsofcyclic} and the Sewing Theorem~\ref{thm:shemersewing}, one can obtain a large family of neighborly polytopes.

\begin{shaded}\index{sewing}
\begin{constr}[Sewing: the family~$\cS$]\label{constr:cS}
\hspace*{\fill}
  \begin{itemize}[\textbullet]
  \item Let $\vv P_0:=\cyc{d}{n}$ be an even-dimensional cyclic polytope.
  \item Let $\cF_0$ be a universal flag of $\vv P_0$. It can be found using Remark~\ref{rmk:universalflagsofcyclic}.
  \item For $i=1\dots k$:
  \begin{itemize}[\textbullet]
  \item Let $\vv P_i:=\vv P_{i-1}[\cF_{i-1}]$. Then $\vv P_i$ is neighborly by  Theorem~\ref{thm:shemersewing}(\ref{it:newisneigh}).
  \item Theorem~\ref{thm:shemersewing}(\ref{it:newunifaces}) constructs a universal flag $\cF_i$ of $\vv P_i$.
  \end{itemize}
  \item $\vv P:=\vv P_k$ is a neighborly polytope in $\cS$.
 \end{itemize}
\end{constr}
\end{shaded}

This method generates a family of neighborly polytopes that we call \defn{totally sewn} polytopes
 and denote by \defn{$\cS$}. In contrast to Shemer's original definition of totally sewn polytopes, we do not admit arbitrary universal flags of $\vv P[\cF]$ for sewing, but only those that arise from Theorem~\ref{thm:shemersewing}\eqref{it:newunifaces}.

\section{Extended Sewing: flags that contain universal subflags}

We are now almost ready to present our first new construction in the Extended Sewing Theorem~\ref{thm:extshemersewing}. It generalizes the first part of the Sewing Theorem~\ref{thm:shemersewing}, \ie it shows how to extend a neighborly matroid with a universal flag to a new neighborly matroid.
The analogue of the second part is Proposition~\ref{prop:extnewunifaces}, where we find universal faces of the new matroid. 
\\

In order to prove that Extended Sewing works, we need the following lemma, which generalizes~\cite[Theorem 3.1]{TrelfordVigh2011}, and the notation $\cF'/F_i:=\{F_j'/F_i\}_{j=i+1}^{m}$ where $F_j'/F_i$ is the face of $\cP/F_i$ that represents~$F_j'$.

\begin{lemma}\label{lem:quotientsofextShemersewing}
Let $\cP$ be a uniform neighborly matroid of rank $\rd$. Let $\cF'=\{F_k'\}_{k=1}^l$ be a flag of $\cP$ that contains a universal subflag $\cF=\{F_j\}_{j=1}^m$,
 where $m=\ffloor{\rd-1}{2}$ and $F_j=\bigcup_{i=1}^{j}\{x_i,y_i\}$.
 Finally, let $p$ be sewn onto $\cP$ through~$\cF'$.\\

 If $F_{i-1}\cup \{y_{i}\}$ does not belong to $\cF'$, then
\begin{equation*}
\cP[\cF']/\{F_{i-1},x_{i},p\}\simeq (\cP/F_i)[\cF'/F_i].\end{equation*}

This isomorphism sends $y_{i}$ to the vertex sewn through the flag~$[\cF'/F_i]$, while  the remaining vertices are mapped to their natural counterparts.
\end{lemma}
\begin{proof}
By Proposition~\ref{prop:allquotientsofle}, the contraction $\cP[\cF']/F_{i-1}$ is a lexicographic extension of $\cP/F_{i-1}$ whose signature coincides with that of $[\cF']$ by removing the first $2(i-1)$ elements. Hence $\cP[\cF']/F_{i-1}$ must be one of the extensions
\begin{equation*}
\cP[\cF']/F_{i-1} \in
\left.
\begin{cases} 
\quad\cP/F_{i-1}[x_i^+,y_i^+,x_{i+1}^-,\dots],
\\
\quad\cP/F_{i-1}[x_i^-,y_i^-,x_{i+1}^+,\dots],
\\
\quad\cP/F_{i-1}[x_i^+,y_i^-,x_{i+1}^+,\dots],
\\
\quad\cP/F_{i-1}[x_i^-,y_i^+,x_{i+1}^-,\dots]
\end{cases}\right\}.
\end{equation*}

If $F'_{k-1}$ is the face of $\cF'$ corresponding to $F_{i-1}$, and $U'_{k}=F'_{k}\setminus F'_{k-1}$, then the first two cases are possible when $U'_{k}=\{x_i,y_i\}$, and the last two when $U'_{k}= \{x_i\}$ (the case $U'_{k}= \{y_i\}$ is excluded by hypothesis).
We use Proposition~\ref{prop:allquotientsofle} twice on each of these (contracting successively	 $x_i$ and $p$) to get \(\cP[\cF']/\{F_{i-1},x_{i},p\}\simeq (\cP/\{F_{i-1},x_{i},y_i\})[x_{i+1}^+,\dots]=(\cP/F_i)[\cF'/F_i]\).
\end{proof}

We can now state and prove the Extended Sewing Theorem.\index{Extended Sewing}

\begin{theorem}[The Extended Sewing Theorem]\label{thm:extshemersewing}
Let $\cP$ be a uniform neighborly oriented matroid of rank~$\rd$ with a flag $\cF'=\{F_k'\}_{k=1}^l$ that contains a universal subflag $\cF=\{F_j\}_{j=1}^m$, where $F_j=\bigcup_{i=1}^{j}\{x_i,y_i\}$ and $m=\ffloor{\rd-1}{2}$.
Let $p$ be sewn onto $\cP$ through~$\cF'$. Then $\cP[\cF']$ is a uniform neighborly matroid of rank $\rd$.
\end{theorem}
\begin{proof}
The proof is by induction on $\rd$. Observe for the base case that all acyclic uniform matroids of rank $1$ or $2$ are neighborly.

Assign the labels to $x_1$ and $y_1$ in such a way that the extension $\cP[\cF']$ is either the lexicographic extension 
\(\cP\left[x_1^+,y_1^+,\dots\right]\) or \(\cP\left[x_1^+,y_1^-,\dots\right]\) (depending on whether $F'_1=\{x_1,y_1\}$ or $F'_1=\{x_1\}$). 

We check that $\cP[\cF']$ is neighborly by checking that $\pGale{\cP[\cF']}$ is balanced, \ie we check that every circuit $X$ of $\cP[\cF']$ is halving. That is, we want to see that $\ffloor{\rd+1}{2}\leq |X^+|\leq \fceil{\rd+1}{2}$. Let $X\in \ci(\cP)$:

\begin{enumerate}
 \item If $X(p)=0$, then $X$ is halving because it is also a circuit of $\cP$, and $\cP$ is neighborly.
 \item If $X(p)\neq 0$ and $X({x_1})=0$, we use that $p$ and $x_1$ are $(-1)$-inseparable by the definition of lexicographic extension (see Lemma~\ref{lem:leinseparable}).
 By Lemma~\ref{lem:circinseparable}, there is a circuit $X'\in \ci(\tilde \cP)$ with $X'({x_1})=X(p)$, $X'(p)=0$ and $X'(e)=X(e)$ for all $e\notin\{x_1,p\}$. 
 Observe that $|X^+|=|X'^+|$. Since $X'(p)=0$, $X'$ is halving by the previous point, and hence so is $X$.
 \item If $X(p)\neq 0$ and $X({x_1})\neq 0$ then $X(p)=-X({x_1})$ because $p$ and $x_1$ are $(-1)$-inseparable. Observe that the rest of the values of $X$ correspond to a circuit of $\cP[\cF']/\{p,x_1\}$. If $\cP[\cF']/\{p,x_1\}$ is neighborly, we are done.

By Lemma~\ref{lem:quotientsofextShemersewing}, $\cP[\cF']/\{p,x_1\}\simeq \left(\cP/F_1\right)[\cF'/F_1]$. Since the edge $\{x_1,y_1\}$ was universal, the oriented matroid $\cP/F_1$ is  neighborly, and the flag $\cF'/F_1$  contains the universal flag $\cF/F_1$. The result now follows by induction.\qedhere
\end{enumerate}
\end{proof}

\iftoggle{bwprint}{%
\begin{figure}[htpb]
\centering
 \subbottom[$\vv A$]{\label{sfig:sewing0}\includegraphics[width=.28\textwidth]{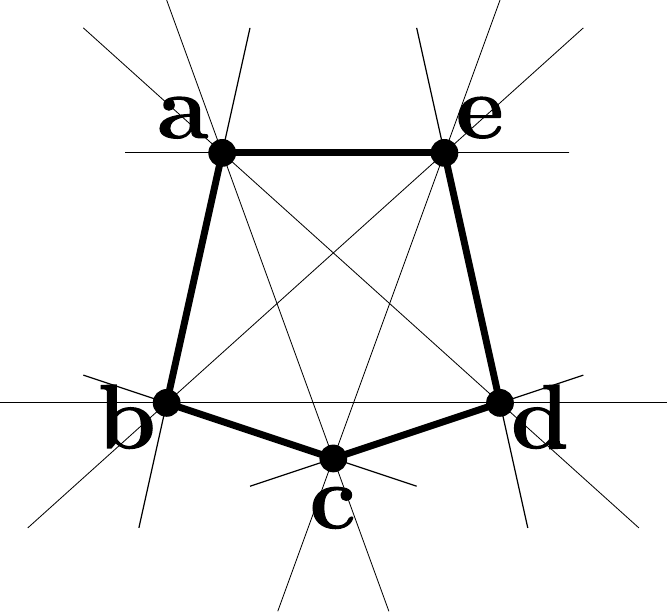}}\qquad 
 \subbottom[$\vv A\bracket{\vv a^+,\vv e^+,\vv c^-}$]{\label{sfig:sewing1}\includegraphics[width=.28\textwidth]{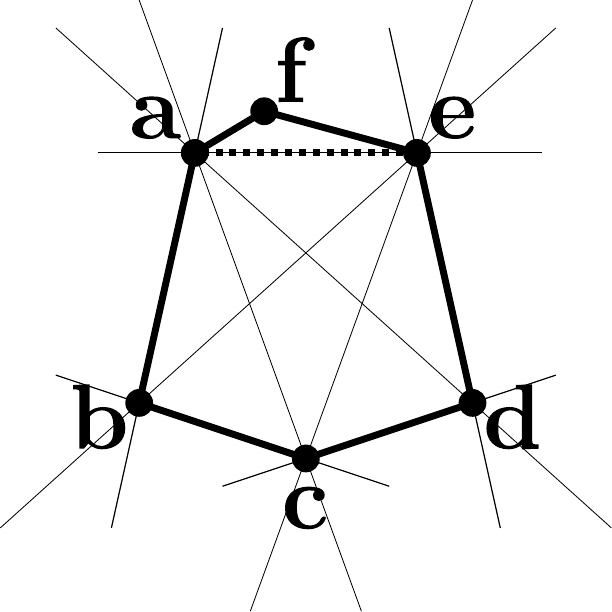}}\qquad 
 \subbottom[$\vv A\bracket{\vv a^+,\vv e^-,\vv c^+}$]{\label{sfig:sewing2}\includegraphics[width=.28\textwidth]{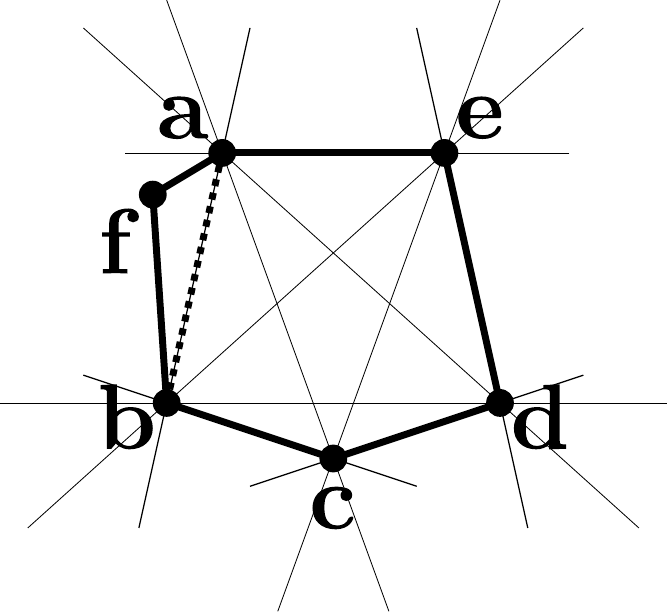}}
 \caption[Extended Sewing example]{Extended Sewing: The point configuration $\vv A$ depicted in Figure~\ref{sfig:sewing0} is the vertex set of a pentagon, and $\{\vv a,\vv e\}$ is one of its universal edges. The sewing of $\vv f$ onto $\{\vv a,\vv e\}$ is shown in Figure~\ref{sfig:sewing1} and the sewing of $\vv f$ onto $\{\vv a\}\subset \{\vv a,\vv e\}$ in Figure~\ref{sfig:sewing2}. In the first case, $\{\vv a,\vv f\}$ and $\{\vv e,\vv f\}$ become universal faces, while $\{\vv a,\vv e\}$ is not a universal face any more. In the second case, $\{\vv a,\vv f\}$ and $\{\vv a,\vv e\}$ are universal faces, while $\{\vv e,\vv f\}$ is not.  }
 \label{fig:sewing}
\end{figure}
}{%
\iftoggle{print}{%
\begin{figure}[htpb]
\centering
 \subbottom[$\vv A$]{\label{sfig:sewing0}\includegraphics[width=.28\textwidth]{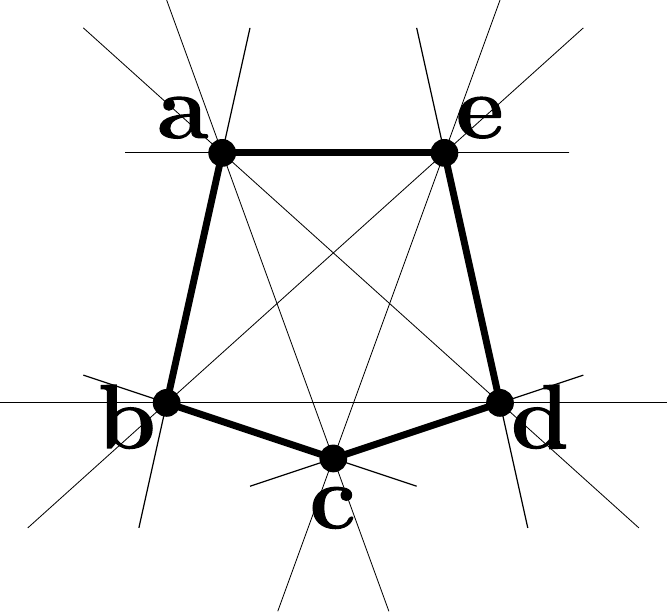}}\qquad 
 \subbottom[$\vv A\bracket{\vv a^+,\vv e^+,\vv c^-}$]{\label{sfig:sewing1}\includegraphics[width=.28\textwidth]{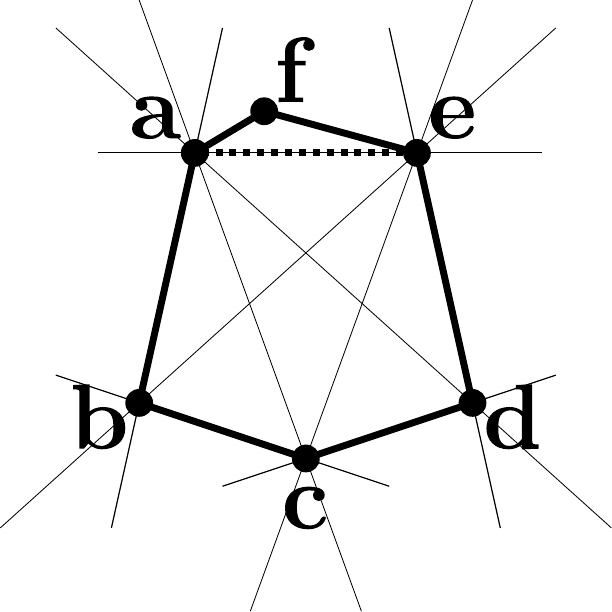}}\qquad 
 \subbottom[$\vv A\bracket{\vv a^+,\vv e^-,\vv c^+}$]{\label{sfig:sewing2}\includegraphics[width=.28\textwidth]{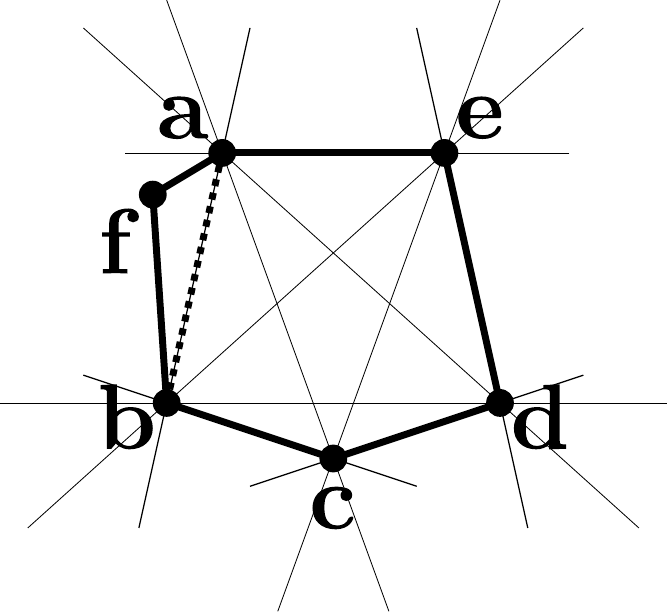}}
 \caption[Extended Sewing example]{Extended Sewing: The point configuration $\vv A$ depicted in Figure~\ref{sfig:sewing0} is the vertex set of a pentagon, and $\{\vv a,\vv e\}$ is one of its universal edges. The sewing of $\vv f$ onto $\{\vv a,\vv e\}$ is shown in Figure~\ref{sfig:sewing1} and the sewing of $\vv f$ onto $\{\vv a\}\subset \{\vv a,\vv e\}$ in Figure~\ref{sfig:sewing2}. In the first case, $\{\vv a,\vv f\}$ and $\{\vv e,\vv f\}$ become universal faces, while $\{\vv a,\vv e\}$ is not a universal face any more. In the second case, $\{\vv a,\vv f\}$ and $\{\vv a,\vv e\}$ are universal faces, while $\{\vv e,\vv f\}$ is not.  }
 \label{fig:sewing}
\end{figure}
}{%
\begin{figure}[htpb]
\centering
 \subbottom[$\vv A$]{\label{sfig:sewing0}\includegraphics[width=.28\textwidth]{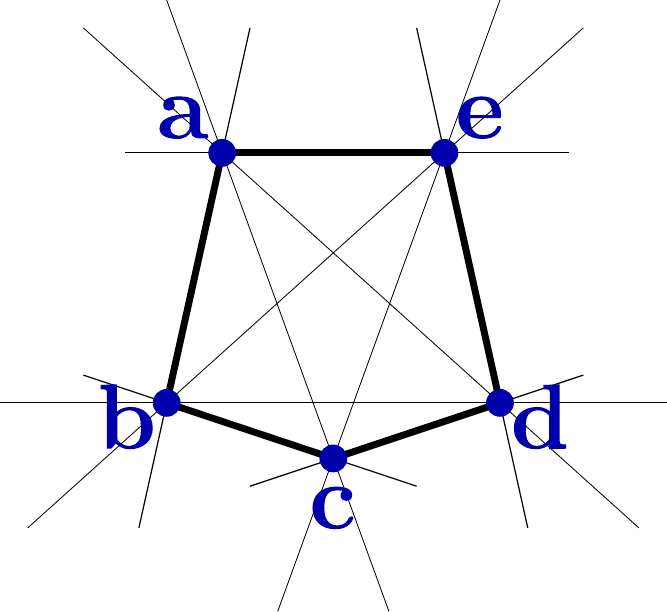}}\qquad 
 \subbottom[$\vv A\bracket{\vv a^+,\vv e^+,\vv c^-}$]{\label{sfig:sewing1}\includegraphics[width=.28\textwidth]{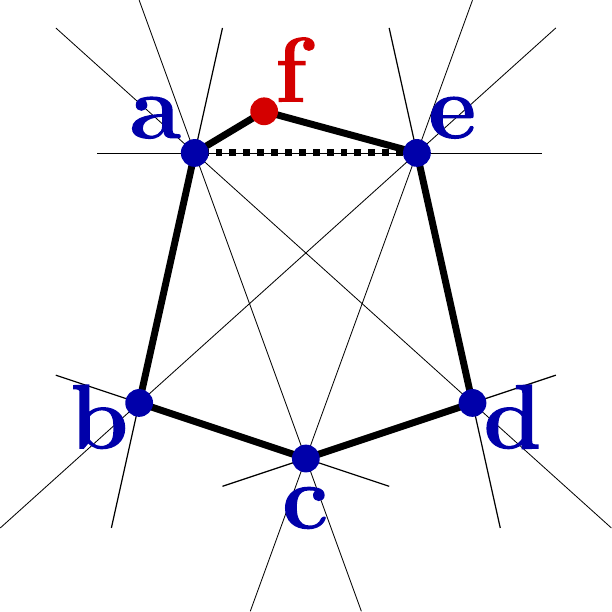}}\qquad 
 \subbottom[$\vv A\bracket{\vv a^+,\vv e^-,\vv c^+}$]{\label{sfig:sewing2}\includegraphics[width=.28\textwidth]{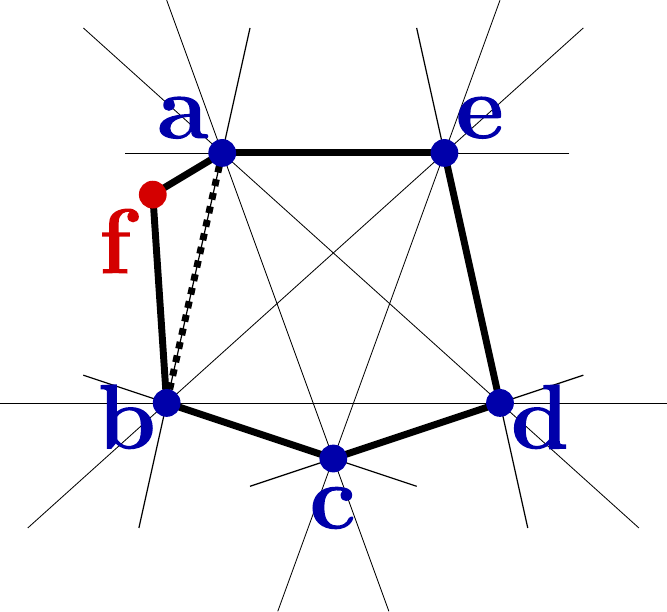}}
 \caption[Extended Sewing example]{Extended Sewing: The point configuration $\vv A$ depicted in Figure~\ref{sfig:sewing0} is the vertex set of a pentagon, and $\{\vv a,\vv e\}$ is one of its universal edges. The sewing of $\vv f$ onto $\{\vv a,\vv e\}$ is shown in Figure~\ref{sfig:sewing1} and the sewing of $\vv f$ onto $\{\vv a\}\subset \{\vv a,\vv e\}$ in Figure~\ref{sfig:sewing2}. In the first case, $\{\vv a,\vv f\}$ and $\{\vv e,\vv f\}$ become universal faces, while $\{\vv a,\vv e\}$ is not a universal face any more. In the second case, $\{\vv a,\vv f\}$ and $\{\vv a,\vv e\}$ are universal faces, while $\{\vv e,\vv f\}$ is not.  }
 \label{fig:sewing}
\end{figure}
}
}

A first application of the Extended Sewing Theorem is the construction of cyclic polytopes.

\begin{proposition}[{\cite[Theorem 5.1]{LeeMenzel2010}}]\label{prop:cyclicareextendedsewn}
Let $\cP$ be the oriented matroid of a cyclic polytope $\cyc{d}{n}$ with elements $a_1,\dots,a_n$ labeled in cyclic order, and let~$\cF$ be the flag $\cF=\{a_n\}\subset\{a_{n-1},a_n\}\subset \dots \subset\{a_{n-d+1},\dots,a_n\}$. Then $\cP[\cF]$ is the oriented matroid of the cyclic polytope $\cyc{d}{n+1}$.
\end{proposition}
\medskip

\subsection{Universal faces created by Extended Sewing}

We can tell many universal faces (and flags) of the neighborly polytopes constructed using the Extended Sewing Theorem~\ref{thm:extshemersewing} thanks to Proposition~\ref{prop:extnewunifaces}, the analogue of the second part of the Sewing Theorem~\ref{thm:shemersewing}.
These faces are best described using the following notation for flags that contain a fixed universal subflag.

\begin{definition}
Let $\cP$ be a neighborly matroid of rank~$2m+1$ and let $\cF'=\{F_k'\}_{k=1}^l$ be a flag of $\cP$ that contains the universal subflag $\cF=\{F_j\}_{j=1}^m$, where $F_j=\bigcup_{i=1}^{j} \{x_i,y_i\}$. Observe that for each $1\leq i\leq j$, either $F_{i-1}\cup\{x_i\}$ or $F_{i-1}\cup\{y_i\}$ or maybe none belongs to~$\cF'$, but not both.
We say that $F_{i}\in \cF$ is \defn{$x_i$-split} (resp. \defn{$y_i$-split}) in $\cF'$ if $F_{i-1}\cup\{x_i\}$ (resp. $F_{i-1}\cup\{y_i\}$) belongs to~$\cF'$, and \defn{non-split} if neither $F_{i-1}\cup\{x_i\}$ nor $F_{i-1}\cup\{y_i\}$ belongs to $\cF'$.
Moreover, we say that $F_i$ is \defn{even} in $\cF'$ if the number of non-split faces~$F_j$ with $j\leq i$ is even, $F_i$ is \defn{odd} otherwise.
\end{definition}

For example, if $l=2$ and $\cF=\{x_1,y_1\}\subset\{x_1,y_1,x_2,y_2\}$ is a universal flag, then $F_1$ is $x_1$-split and $F_2$ is non-split in the flag $\cF'= \{x_1\}\subset\{x_1,y_1\}\subset\{x_1,y_1,x_2,y_2\}$. Moreover, $F_1$ is even in $\cF'$ whereas $F_2$ is odd. In comparison, in the flag $\cF''=\{x_1,y_1\}\subset\{x_1,y_1,y_2\}\subset\{x_1,y_1,x_2,y_2\}$, $F_1$ is non-split and $F_2$ is $y_2$-split; and both $F_1$ and $F_2$ are odd.

\begin{proposition}\label{prop:extnewunifaces}
Let $\cP$ be a uniform neighborly oriented matroid of rank~$\rd$ with a flag $\cF'=\{F_k'\}_{k=1}^l$ that contains a universal subflag $\cF=\{F_j\}_{j=1}^m$, where $F_j=\{x_i,y_i\}_{i=1}^{j}$ and $m=\ffloor{\rd-1}{2}$. 
Let $p$ be sewn onto $\cP$ through~$\cF'$.
Then the following are universal faces of $\cP[\cF']$:
\begin{enumerate}
 \item\label{it:extnewuf1} $F_i$, where $1\leq i \leq m$, if $F_i$ is even.
 \item\label{it:extnewuf2} $(F_{j}\setminus{x_i})\cup p$, where $1\leq i\leq j\leq m$, if
      \begin{enumerate}[(i)]
	\item\label{it:extnewuf21} $F_i$ is not split and $F_j/F_i$ is even in $\cF'/F_i$, or
	\item\label{it:extnewuf22} $F_i$ is $x_i$-split and $F_j/F_i$ is odd in $\cF'/F_i$, or
	\item\label{it:extnewuf23} $F_i$ is $y_i$-split and $F_j/F_i$ is even in $\cF'/F_i$.
      \end{enumerate}
 \item\label{it:extnewuf3} $(F_{j}\setminus{y_i})\cup p$, where $1\leq i\leq j\leq m$, if
      \begin{enumerate}[(i)]
	\item\label{it:extnewuf31} $F_i$ is not split and $F_j/F_i$ is even in $\cF'/F_i$, or
	\item\label{it:extnewuf32} $F_i$ is $x_i$-split and $F_j/F_i$ is even in $\cF'/F_i$, or
	\item\label{it:extnewuf33} $F_i$ is $y_i$-split and $F_j/F_i$ is odd in $\cF'/F_i$.
      \end{enumerate}
\end{enumerate}

\end{proposition}
\begin{proof}
Without loss of generality, we will assume that all split faces are $x_i$-split.
The proof relies on applying, case by case, Proposition~\ref{prop:allquotientsofle} to reduce the contraction to a lexicographic extension that we know to be neighborly because of Theorem~\ref{thm:extshemersewing}.
We also use the following observation about the signature of the extension $\cP[\cF']$: the sign of $x_i$ is $+$ if and only if $F_{i-1}$ is even. This follows from the fact that, by Definition~\ref{def:sewing}, there are some elements $a,b$ and some $\ep=\pm$ such that

\begin{equation*}
\cP[\cF']
=
\begin{cases}
\cP[\dots,a^\ep,x_i^{-\ep},y_i^{\ep},b^{-\ep},\dots]& \text{ if $F_i$ is $x_i$-split,}
\\
\cP[\dots,a^\ep,x_i^{-\ep},y_i^{-\ep},b^{\ep},\dots]&\text{ if it is not split.}
\end{cases}
\end{equation*}
In particular, if $F_i$ is even, then $\cP[\cF']/F_i \simeq (\cP/F_i)[\cF'/F_i]$ and $\cF'/F_i$ is a universal flag of $\cP/F_i$, which is neighborly since $F_i$ is a universal face. This proves point~\ref{it:extnewuf1}.

Moreover, independently of whether $F_i$ is even or odd, 
\begin{equation*}
\cP[\cF']/(F_{i-1}\cup \{p\})\simeq  
\begin{cases}
 \cP/(F_{i-1}\cup \{x_i\})[y_i^{+},x_{i+1}^{-},\dots]& \text{ if $F_i$ is split,}
\\
 \cP/(F_{i-1}\cup \{x_i\})[y_i^{-},x_{i+1}^{+},\dots]& \text{ if it is not.}
\end{cases}
\end{equation*}

Hence, $\cP[\cF']/(F_{i-1}\cup \{x_i,p\})\simeq (\cP/F_i)[\cF'/F_i]$ always. If moreover $F_i$ is not split then $\cP[\cF']/(F_{i-1}\cup \{y_i,p\})\simeq (\cP/F_i)[\cF'/F_i]$. We know that $F_j/F_i$ is a universal face of $(\cP/F_i)[\cF'/F_i]$ when it is even. This proves points~\ref{it:extnewuf21}, \ref{it:extnewuf23}, \ref{it:extnewuf31} and \ref{it:extnewuf32}.

If $F_i$ is split, then $\tilde \cP/(F_{i-1}\cup \{y_i,p\})\simeq (\cP/F_i)[-\cF'/F_i]$, where $[-\cF'/F_i]$ means the extension $[\cF'/F_i]$ with the signs reversed. Using the previous observation, we obtain that $(\cP/F_i[-\cF'/F_i])/(F_j/F_i)\simeq (\cP/F_j)[\cF'/F_j]$ when $F_j/F_i$ is odd, and this proves the remaining points~\ref{it:extnewuf22} and \ref{it:extnewuf33}.
\end{proof}

\begin{remark}\label{rmk:extnewunifaces}
In particular, Proposition~\ref{prop:extnewunifaces} provides a simple way to tell universal flags of~$\cP[\cF']$ (cf. Remark~\ref{rmk:universalflagsofcyclic}). We start with universal edges:

\begin{itemize}
\item If $F_1$~is not split then $\{x_1,p\}$ and $\{y_1,p\}$ are universal edges of $\cP[\cF']$; 
\item if $F_1$~is $x_1$-split, then $\{x_1,p\}$ and $\{x_1,y_1\}$ are universal edges of $\cP[\cF']$; 
\item finally, if $F_1$~is $y_1$-split, then $\{y_1,p\}$ and $\{x_1,y_1\}$ are universal edges of $\cP[\cF']$.
 
\end{itemize}

The contraction of any of these universal edges is a matroid isomorphic to $(\cP/F_1)[\cF'/F_1]$, and we can inductively build a universal flag of $\cP[\cF']$.
\end{remark}

The rank~$3$ example in Figure~\ref{fig:sewing} gives some intuition on why the universal edges of the previous remark appear, but we can see it more clearly in the next example.

\begin{example}\label{ex:universalfaces}
 Let $\cM$ be a neighborly oriented matroid of rank~$5$ with a universal flag $\cF=F_1\subset F_2$, where $F_1=\{a,b\}$ and $F_2=\{a,b,c,d\}$. 
 Consider the lexicographic extensions by the following elements
\begin{align*}
p_1&=[a^+,b^+,c^-,d^-,e^+],\\p_2&=[a^+,b^-,c^+,d^+,e^-],\\p_3&=[a^+,b^+,c^-,d^+,e^-],\text{ and }\\p_4&=[a^+,b^-,c^+,d^-,e^+],
\end{align*}
where $e$ is any element of $\cM$.  For $i=1,2,3,4$, each $p_i$ gives  rise to the oriented matroid $\cM_i=\cM[p_i]$, which corresponds to sewing through the flag~$\cF_i$, with
\begin{align*}
\cF_1&=\{a,b\}\subset\{a,b,c,d\},\\
\cF_2&=\{a\}\subset\{a,b\}\subset\{a,b,c,d\},\\
\cF_3&=\{a,b\}\subset\{a,b,c\}\subset\{a,b,c,d\}\text{ and }\\
\cF_4&=\{a\}\subset\{a,b\}\subset\{a,b,c\}\subset\{a,b,c,d\}.
\end{align*}
Observe that $F_1$ is $a$-split in $\cF_2$ and $\cF_4$, while $F_2$ is $c$-split in $\cF_3$ and $\cF_4$. Moreover, $F_1$ is even in $\cF_2$ and $\cF_4$, and $F_2$ is even in $\cF_1$ and $\cF_4$.
Table~\ref{tb:exunifaces} shows for which~$\cM_i$ each of the following sets of vertices is a universal face.

\begin{table}[htpb]
  \caption{Universal faces in Example~\ref{ex:universalfaces}}	\label{tb:exunifaces}
   \renewcommand{\arraystretch}{2}
\centering
\newcolumntype{C}[1]{>{\centering\let\newline\\\arraybackslash\hspace{0pt}}m{#1}}
\begin{tabular}{ccC{1cm}C{1cm}C{1cm}C{1cm}C{1cm}C{1cm}C{1cm}C{1cm}}
\hline
	&&$ab$&$p_ib$&$ap_i$&$abcd$&$p_ibcd$&$ap_icd$&$abp_id$&$abcp_i$\\\hline
$\cM_1$&&\cross&\tick&\tick&\tick&\cross&\cross&\tick&\tick\\
$\cM_2$&&\tick&\cross&\tick&\cross&\tick&\cross&\tick&\tick\\
$\cM_3$&&\cross&\tick&\tick&\cross&\tick&\tick&\cross&\tick\\
$\cM_4$&&\tick&\cross&\tick&\tick&\cross&\tick&\cross&\tick\\ 
\hline
\end{tabular}
\end{table}

\end{example}
\medskip

\begin{remark}
 Theorem~\ref{thm:extshemersewing} not only generalizes the Sewing Theorem, but also includes Barnette's facet-splitting technique~\cite[Theorem 3]{Barnette1981}, which corresponds to the case where all faces of the universal flag are split.
\end{remark}

\subsection{Extended Sewing (and Omitting)}

Just like in the construction of the family {$\cS$}, we can use the Extended Sewing Theorem~\ref{thm:extshemersewing} and Proposition~\ref{prop:extnewunifaces} to obtain a large family~\defn{$\cE$} of neighborly polytopes that contains $\cS$. In fact, since cyclic polytopes belong to~$\cE$ by Proposition~\ref{prop:cyclicareextendedsewn}, it suffices to start sewing on a simplex (hence, Remark~\ref{rmk:universalflagsofcyclic} is not needed anymore).

\begin{shaded}
\begin{constr}[Extended Sewing: the family $\cE$]\label{constr:cE}
\hspace*{\fill}
  \begin{itemize}[\textbullet]
  \item Let $\vv P_0:=\simp{d}$ be a $d$-dimensional simplex.
  \item Let $\cF_0'$ be a flag of $\vv P_0$ that contains a universal subflag $\cF_0$. $\cF_0$ is built using the fact that all edges of a simplex are universal.
  \item For $i=1\dots k$:
  \begin{itemize}[\textbullet]
  \item Let $\vv P_i:=\vv P_{i-1}[\cF_{i-1}']$, which is neighborly by Theorem~\ref{thm:extshemersewing}.
  \item Use Proposition~\ref{prop:extnewunifaces} (or Remark~\ref{rmk:extnewunifaces}) to find a universal flag $\cF_i$ of $\vv P_i$.
  \item Let $\cF_i'$ be any flag of $\vv P_i$ that contains $\cF_i$ as a subflag.
  \end{itemize}
  \item $\vv P:=\vv P_k$ is a neighborly polytope in $\cE$.
 \end{itemize}
\end{constr}
\end{shaded}

Moreover, since subpolytopes of neighborly polytopes are neighborly, any polytope obtained from a member of~$\cE$ by omitting some vertices is also neighborly. The polytopes that can be obtained in this way via sewing and omitting form a family that we denote \defn{$\cO$}. 

\index{Extended Sewing!and omitting}
\begin{shaded}
\begin{constr}[Extended Sewing and Omitting: the family~$\cO$]\label{constr:cO}
\hspace*{\fill}
  \begin{itemize}[\textbullet]
  \item Let $\vv Q\in \cE$ be a neighborly polytope constructed using Extended Sewing.
  \item Let $\vv S\subseteq\verts(\vv Q)$ be a subset of vertices of $\vv Q$.
  \item $\vv P:=\conv(\vv S)$ is a neighborly polytope in $\cO$.
 \end{itemize}
\end{constr}
\end{shaded}

\subsection{Optimality}

We finish this section by showing that for matroids of odd rank, the flags of the Extended Sewing Theorem~\ref{thm:extshemersewing} are the only ones that yield neighborly polytopes. Therefore, in this sense the sewing construction cannot be further improved.

\begin{proposition}\label{prop:uniqueflags}
Let $\cP$ be a uniform neighborly oriented matroid of odd rank $\rd\geq 3$ with more than $\rd+1$ elements.
Then $\cP[\cF]$ is neighborly if and only if $\cF$ contains a universal subflag.
\end{proposition}
\begin{proof}
By Theorem~\ref{thm:extshemersewing}, this condition is sufficient.
To find necessary conditions, we use that $\cP[\cF]$ is neighborly if and only if every circuit of $\cP[\cF]$ is balanced. 

The proof is by induction on $\rd$. For the base case just observe that neighborly matroids of rank $3$ are polygons, and the only 
flags that yield a polygon with one extra vertex are of the form $\{x\}\subset \{x,y\}$ or just $\{x,y\}$, where $\{x,y\}$ is an edge of the polygon.

By definition, $\cP[\cF]$ is the lexicographic extension $\cP[p]$, with $p$ sewn through~$\cF$. Therefore, 
\(p=\left[a_1^{+},\;a_2^{\ep_2},\dots,\;a_d^{\ep_\rd}\right]\). 
Let $X\in\ci(\cP[\cF])$ be a circuit with $\{p,a_1\}\subset\underline X$. Since $p$ and $a_1$ are $(-1)$-inseparable by Lemma~\ref{lem:leinseparable}, $X(p)=-X({a_1})$. Hence, if $X$ is halving, so is $X\setminus \{p,a_1\}$. Now $X\setminus \{p,a_1\}$ is a circuit of $\cP[\cF]/\{p,a_1\}$, and  all  circuits of $\cP[\cF]/\{p,a_1\}$ arise this way. Hence $\cP[\cF]/\{p,a_1\}$ is neighborly.

By Proposition~\ref{prop:allquotientsofle}, \begin{equation*}\label{eq:Pquopa1}\cP[\cF]/\{p,a_1\}\simeq \cP/\{a_1,a_2\}\left[a_3^{-\ep_2\ep_3},\;a_4^{-\ep_2\ep_4},\;,\dots,\;a_\rd^{-\ep_2\ep_\rd}\right],\end{equation*}
where the second extension is by $a_2$. Hence, by Lemma~\ref{lem:leinseparable}, $a_2$ and $a_3$ are $(\ep_2\ep_3)$-inseparable in $\cP[\cF]/\{p,a_1\}$, which is a neighborly matroid of odd rank and corank at least~$2$. By Lemma~\ref{lem:balonlycovar}, $\ep_2\ep_3=-$.

In particular, either $(\ep_2, \ep_3)=(+,-)$, or $(\ep_2, \ep_3)=(-,+)$. The first option implies that $F_1=\{a_1,a_2\}$, and the second one that $F_1=\{a_1\}$ and $F_2=\{a_1,a_2\}$. Since $(\cP[\cF]/\{p,a_1\})\setminus a_2\simeq\cP/\{a_1,a_2\}$ by Lemma~\ref{lem:contractdeletele}, if $\cP[\cF]/\{p,a_1\}$ is neighborly, then $\cP/\{a_1,a_2\}$ must be neighborly and hence $\{a_1,a_2\}$ must be a universal edge of $\cP$. 
\end{proof}

\chapter{The Gale Sewing Construction}\label{ch:thethm}

In this chapter, we present a different method to construct neighborly matroids. It is also based on lexicographic extensions, but works in the dual. Namely, it extends balanced matroids to new balanced matroids.

A priori, it might seem that the approach taken in the sewing construction (extending neighborly matroids) makes more sense since neighborliness is preserved by deletion of elements while balancedness is not. 
However, observe that if $\cM$ is a balanced matroid such that $\cM \setminus \{x,y\}$ is also balanced, this just means that both $\Gale \cM$ and $\Gale\cM/\{x,y\}$ are neighborly. That is, that $\{x,y\}$ is a universal edge of $\Gale\cM$. 
Recall that all the neighborly polytopes that we have built so far have universal edges and therefore, they can all be constructed by adding elements to its dual matroid. 
In fact, in Section~\ref{sec:comparing} we will see that with our new construction we can build all polytopes in $\cO$. 

Nevertheless, this construction will not be able to construct all neighborly polytopes either (and neither any construction that follows the same approach). 
Indeed, Bokowski and Sturmfels proved in~\cite{BokowskiSturmfels1987} that the sphere ``$ \cM^{10}_{416}$'' is polytopal, providing an instance of a neighborly $4$-polytope with $10$ vertices that has no universal edge. This is a neighborly polytope whose dual cannot be constructed with extensions of a balanced matroid.
\\

Before explaining the construction, we must do some observations about the parity of the corank of $\cM$ because the possibility to extend $\cM$ to a new balanced matroid depends strongly on it.
We omit the proofs of the following pair of lemmas, which are direct observing that one can restrict to matroids of rank $2$ thanks to Lemma~\ref{lem:quotientsbalanced}.

\begin{lemma}
 If $\cM$ is a balanced matroid of rank $r$ with $n$ elements such that $n-r-1$ is even, then any single element extension of $\cM$ is balanced.\qed
\end{lemma}

\begin{lemma}\label{lem:forcedbalancedextension}
 If $\cM$ is a balanced matroid of rank $r$ with $n$ elements such that $n-r-1$ is odd, then there is at most one single element extension of $\cM$ that is balanced.

 This extension exists precisely when the signature $\gs:\co(\cM)\rightarrow \{+,-\}$ with $\gs(C)=+$ if $|C^+|<|C^-|$ and $\gs(C)=-$ if $|C^+|>|C^-|$
 is well defined and fulfills the conditions of Theorem~\ref{thm:localization}. In this case, $\gs$ is the signature of the extension. \qed
\end{lemma}

What these lemmas suggest is that to iteratively construct arbitrarily large balanced matroids, one must do \defn{double element extensions}. 
Starting with a balanced matroid $\cM$ of even corank, construct two consecutive extensions $\cM'$ of $\cM$ and $\cM''$ of $\cM'$. Of course, $\cM''$ is completely determined by~$\cM'$ because of Lemma~\ref{lem:forcedbalancedextension}. 
Hence, one must choose an appropriate $\cM'$ that allows to do this second extension.

In Figure~\ref{fig:unexpanding} there is an example of a balanced matroid of rank~$3$ and $6$ elements. The single element extension with a point in any of the shaded regions does not admit any balanced single element extension. On the other hand, the extension by a point in the white regions does admit such an extension. 

\iftoggle{bwprint}{%
\begin{figure}[htpb]
\begin{center}
\includegraphics[width=.6\textwidth]{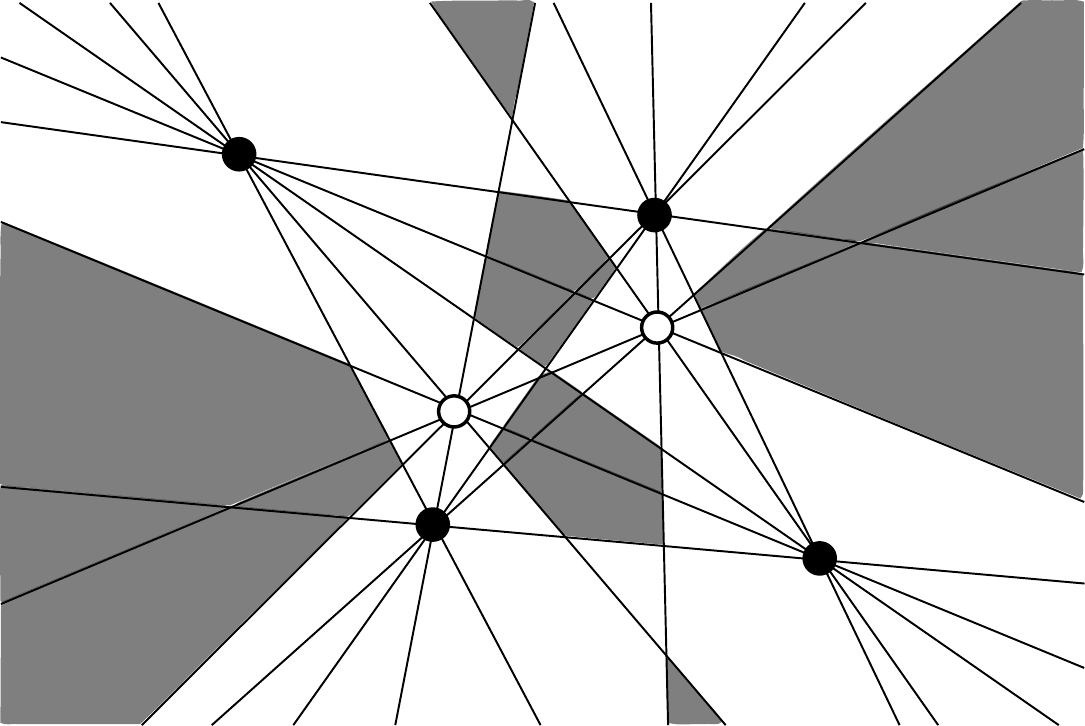}
\end{center}
 \caption[Admissible single element extensions for balanced matroids.]{The single element extension of this matroid with a point in any of the shaded regions is a matroid that does not admit any balanced single element extension.}
 \label{fig:unexpanding}
\end{figure}
}{%
\begin{figure}[htpb]
\begin{center}
\includegraphics[width=.6\textwidth]{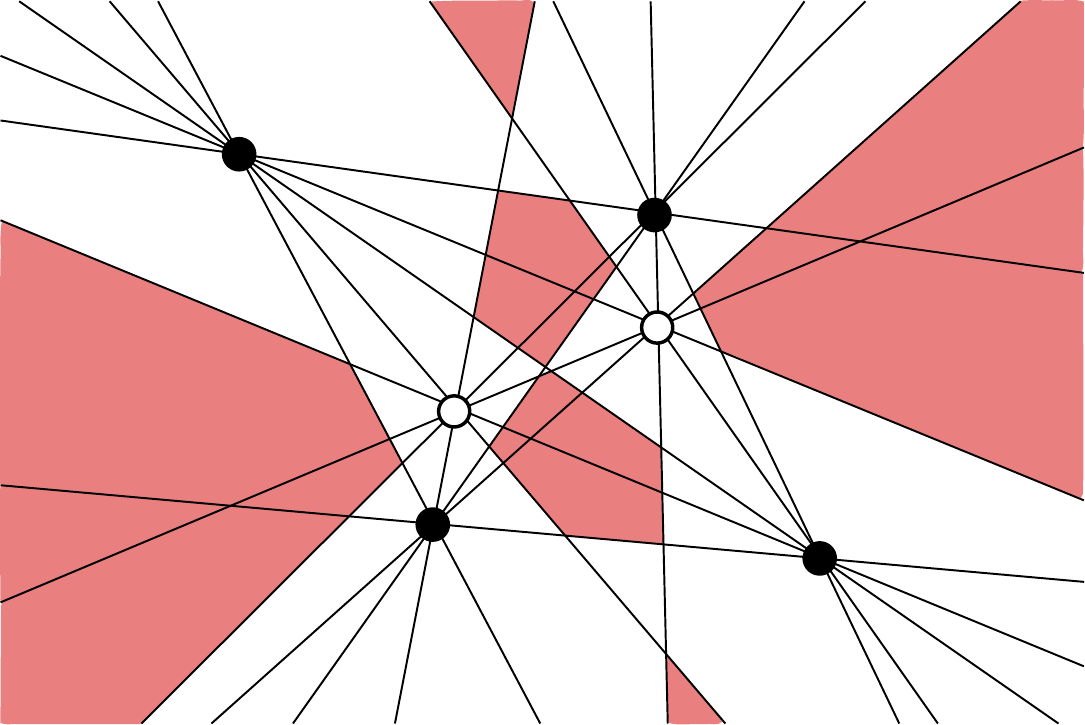}
\end{center}
 \caption[Admissible single element extensions for balanced matroids.]{The single element extension of this matroid with a point in any of the shaded regions is a matroid that does not admit any balanced single element extension.}
 \label{fig:unexpanding}
\end{figure}
}

The observation that in Figure~\ref{fig:unexpanding} all the cells that have one of the original points as a vertex are white will motivate our construction. And the fact that there are white cells without such a vertex hints that we will not be able to make all balanced extensions this way.

\section{The Gale Sewing Construction}

The key ingredient of our construction technique is the Double Extension Theorem~\ref{thm:thethm}, which shows how to perform double element extensions that preserve balancedness:
one must do two consecutive lexicographic extensions by $p=[a_1^{\ep_1},\dots,a_r^{\ep_r}]$ and $q=[p^-,\dots,a_{r-1}^{-}]$. 

Actually, we prove a stronger result, Proposition~\ref{prop:thethm}, that shows that this double element extension preserves the discrepancy. Before proving it, we need a small lemma that explains the contraction of this last element $q$.
\begin{lemma}\label{lem:quotientsofGalesewnareGalesewn}
Let $\cM$ be a uniform oriented matroid of rank $r$, let $a_1\dots a_r$ be elements of $\cM$ and $\ep_1,\dots,\ep_r$ be signs. If  $p$, $q$, $p'$ and $q'$ are defined as
\begin{align*}
p&=[a_1^{\ep_1},a_2^{\ep_2},a_3^{\ep_3},\dots,a_r^{\ep_r}],&q&=[p^-,a_1^{-},a_2^{-},\dots,a_{r-1}^{-}];\\
p'&=[a_2^{-\ep_1\ep_2},a_3^{-\ep_1\ep_3},\dots,a_r^{-\ep_1\ep_r}],& q'&=[p'^{-},a_2^{-},\dots,a_{r-1}^{-}],
\end{align*}
then \[
\left(\cM[p][q]\right)/q\ \simeq\ \left(\cM/a_1\right)[p'][q'].
\]
\end{lemma}
\begin{proof}
Repeatedly applying Proposition~\ref{prop:allquotientsofle},
\begin{align*}
\left(\cM[p][q]\right)/q&= \big(\cM\underbrace{[a_1^{\ep_1},a_2^{\ep_2},\dots,a_r^{\ep_r}]}_{p}\underbrace{[p^-,a_1^{-},\dots,a_{r-1}^{-}]}_{q}\big)/q\\
&\stackrel{\varphi}{\simeq} \big(\cM\underbrace{[a_1^{\ep_1},a_2^{\ep_2},\dots,a_r^{\ep_r}]}_{p}/p\big)\underbrace{[a_1^{-},\dots,a_{r-1}^{-}]}_{\varphi(p)=q'}\\
&\stackrel{\psi}{\simeq} \big(\cM/a_1\big)\underbrace{[a_2^{-\ep_1\ep_2},\dots,a_r^{-\ep_1\ep_r}]}_{\psi(a_1)=p'}\underbrace{[\psi(a_1)^{-},\dots,a_{r-1}^{-}]}_{q'}.\qedhere
\end{align*}
\end{proof}

\begin{proposition}\label{prop:thethm}
Let $\cM$ be a uniform oriented matroid of rank $r$. For any sequence $a_1\dots a_r$ of elements of~$\cM$ and any sequence $\ep_1,\dots,\ep_r$ of signs,
consider the lexicographic extensions 
\begin{itemize}
\item $\cM[p]$ of $\cM$ by $p=[a_1^{\ep_1},a_2^{\ep_2},\dots,a_r^{\ep_r}]$, and 
\item $\cM[p][q]$ of $\cM[p]$ by $q=[p^-,a_1^{-},\dots,a_{r-1}^{-}]$; 
\end{itemize}
then \[\disc(\cM[p][q])=\disc(\cM).\] 
\end{proposition}
\begin{proof}
We prove that $\disc(\cM[p][q])=\disc(\cM)$ by checking that every cocircuit $\tilde C$ of $\cM[p][q]$ has the same discrepancy as some cocircuit $C$ of $\cM$. That is, we prove that $\big||\tilde C^+|-|\tilde C^-|\big|=\big||C^+|-|C^-|\big|$.

If $\tilde C(p)\neq 0$ and $\tilde C(q)\neq 0$ then, by the definition of lexicographic extension, there is a cocircuit $C$ of~$\cM$ such that $\restriction{\tilde C}{\cM}=C$ and $\tilde C(p)=-\tilde C(q)$. Hence $|\tilde C^+|=|C^+|+1$ and $|\tilde C^-|=|C^-|+1$.

The cocircuits $\tilde C$ with $\tilde C(p)=0$ correspond to cocircuits of $(\cM[p][q])/p$, and those with $\tilde C(q)=0$ correspond to cocircuits of $(\cM[p][q])/q$. 
In both cases, there is a cocircuit $C\in \co(\cM)$ with $C(a_1)=0$ that has the same discrepancy as $\tilde C$. Indeed, Proposition~\ref{prop:allquotientsofle} and Lemma~\ref{lem:quotientsofGalesewnareGalesewn} tell us that 
\[(\cM[p][q])/p\simeq (\cM[p][q])/q\simeq(\cM/a_1)\underbrace{[a_2^{-\ep_1\ep_2},\dots,a_r^{-\ep_1\ep_r}]}_{p'}[{p'}^-,a_2^-,\dots,a_{r-1}^-],\] and 
our claim follows by induction on $r$ (it is trivial for $r=1$).
\end{proof}

Our construction theorem is a direct corollary of this result, since a uniform oriented matroid is balanced if and only if its discrepancy is $\leq 1$.

\begin{theorem}[Double Extension]\label{thm:thethm}
Let $\cM$ be a uniform balanced oriented matroid of rank $r$. For any sequence $a_1\dots a_r$ of elements of~$\cM$ and any sequence $\ep_1,\dots,\ep_r$ of signs,
consider the lexicographic extensions 
\begin{itemize}
\item $\cM[p]$ of $\cM$ by $p=[a_1^{\ep_1},a_2^{\ep_2},\dots,a_r^{\ep_r}]$, and 
\item $\cM[p][q]$ of $\cM[p]$ by $q=[p^-,a_1^{-},\dots,a_{r-1}^{-}]$; 
\end{itemize}
then $\cM[p][q]$ is balanced. \qed
\end{theorem}

\iftoggle{bwprint}{%
\begin{figure}[htpb]
\begin{center}
\includegraphics[width=.8\textwidth]{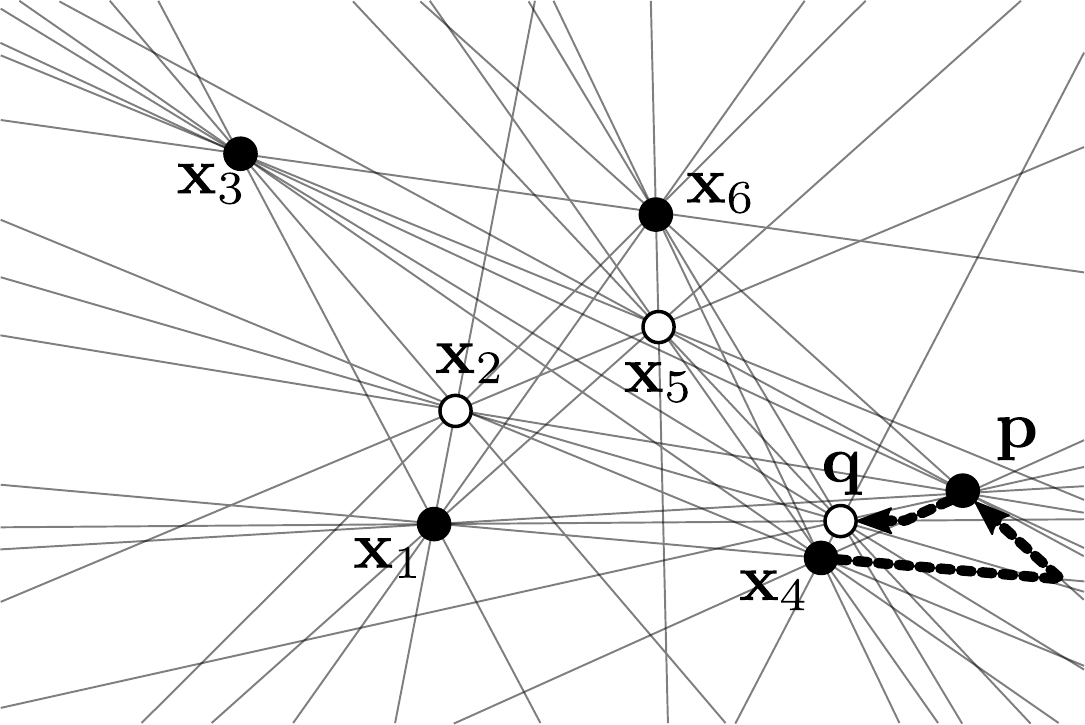}
\end{center}
\caption[Example of Gale Sewing]{The double lexicographic extension of a balanced configuration by $\vv p=[\vv x_4^+,\vv x_1^-,\vv x_6^+]$ and $\vv q=[\vv p^-,\vv x_4^-,\vv x_1^-]$, which is also balanced.}\label{fig:galesewingexample}
\end{figure}
}{%
\begin{figure}[htpb]
\begin{center}
\includegraphics[width=.8\textwidth]{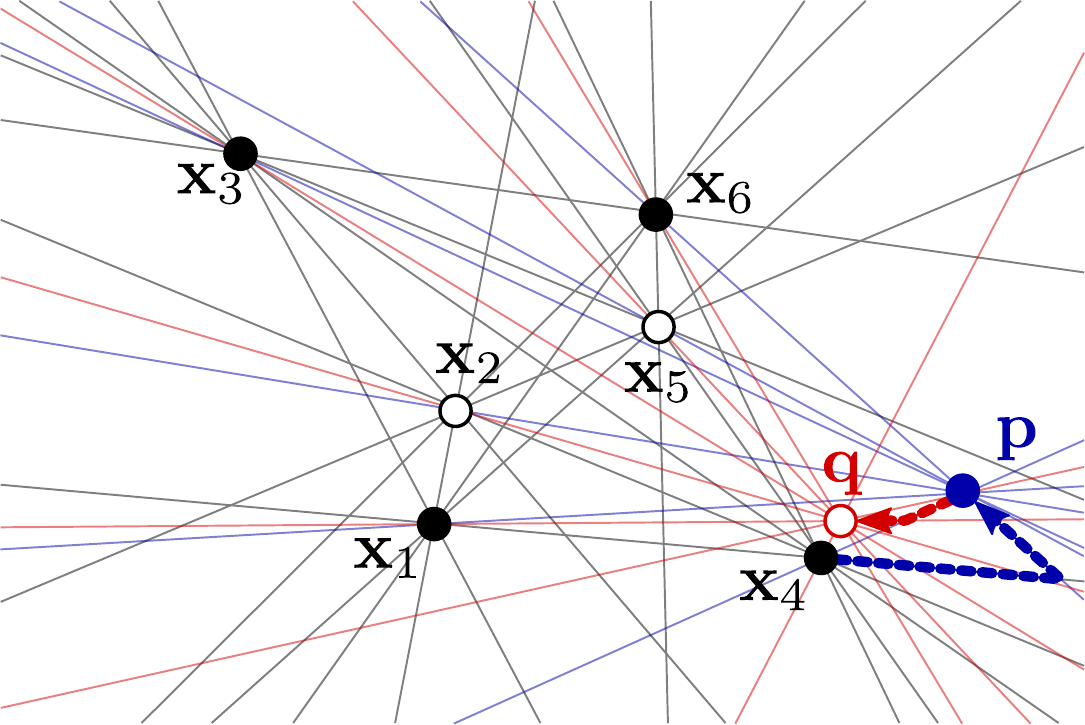}
\end{center}
\caption[Example of Gale Sewing]{The double lexicographic extension of a balanced configuration by $\vv p=[\vv x_4^+,\vv x_1^-,\vv x_6^+]$ and $\vv q=[\vv p^-,\vv x_4^-,\vv x_1^-]$, which is also balanced.}\label{fig:galesewingexample}
\end{figure}
}

\begin{remark}
In a vector configuration $\vv V$, the proof that $\vv V[\vv p][\vv q]$, its lexicographic extension by $\vv p=[\vv a_1^{\ep_1},\dots,\vv a_r^{\ep_r}]$ and $\vv q=[\vv p^-,\dots,\vv a_{r-1}^{-}]$, has the same discrepancy as $\vv V$ is very easy to understand. Every hyperplane $\vvh H$ spanned by a subset of $\vv V$ defines a cocircuit of $\vv V[\vv p][\vv q]$.
The signature of the extension by $\vv q$ implies that if $\vv p\in \vvh H^\pm$ then $\vv q\in \vvh H^\mp$, and hence $\vv q$ balances the discrepancy created by $\vv p$ on this hyperplane.
The other hyperplanes are checked inductively. Indeed, for a hyperplane $\vvh H$ that contains~$\vv p$ but neither~$\vv a_1$ nor~$\vv q$, the fact that $\vv p$ and $\vv a_1$ are inseparable implies that except for~$\vv a_1$, $\vvh H$ looks like a hyperplane spanned by $\vv V$ containing~$\vv a_1$.
Hence $\vv q$ must balance the discrepancy created by $\vv a_1$.
For hyperplanes that go through $\vv p$ and $\vv a_1$ but neither $\vv a_2$ nor $\vv q$, $\vv q$ balances the discrepancy created by~$\vv a_2$; and so on.

Figure~\ref{fig:galesewingexample} displays an example of such a double extension on an affine Gale diagram. 
The reader is invited to follow this justification in the picture (for example, by comparing the hyperplanes spanned by $\{\vv x_4,\vv x_i\}$ with the hyperplanes spanned by $\{\vv p,\vv x_i\}$) and to check how all cocircuits in the diagram are balanced.
\end{remark}

This provides the following method to construct balanced matroids (and hence, by duality, to construct neighborly matroids). 
\begin{shaded}
\begin{constr}[Gale Sewing]\label{constr:cG}
\hspace*{\fill}
  \begin{itemize}[\textbullet]
  \item Let $\cM_0:=\stc{r}$ be the minimal totally cyclic matroid, \ie the oriented matroid of the balanced configuration $\{\vv e_1, \dots, \vv e_r, -\sum_{i=1}^r \vv e_i\}$.
  \item For $k=1\dots m$:
  \begin{itemize}[\textbullet]
  \item Choose different elements $a_{k1},\dots, a_{kr}$ of $\cM_{k-1}$ and choose $\ep_{kj}\in\{+, -\}$ for $j=1\dots r$.
  \item Let $ p_k:=[ a_{k1}^{\ep_{k1}},\dots, a_{kr}^{\ep_{kr}}]$ and $ q_k:=[ p_k^-, a_{k1}^{-},\dots, a_{k(r-1)}^{-}]$.
  \item $\cM_k:=\cM_{k-1}[ p_k][ q_k]$ is balanced because of Theorem~\ref{thm:thethm} and realizable because of Proposition~\ref{prop:realizablele}.
  \end{itemize}
  \item $\cM:=\cM_k$ is a realizable balanced oriented matroid.
  \item $\cP:=\Gale\cM$ is a realizable neighborly oriented matroid. 
  \item Any realization $\vv P$ of $\cP$ is a neighborly polytope in $\cG$.
 \end{itemize}
\end{constr}
\end{shaded}

We call the double extension of Theorem~\ref{thm:thethm} \defn{Gale Sewing}\index{Gale Sewing}, and we denote by \defn{$\cG$} the family of combinatorial types of polytopes whose dual is constructed by repeatedly Gale Sewing from $\{\vv e_1,\dots,\vv e_r,-\sum_{i=1}^r \vv e_i\}$. If $\vv P\in\cG$, we will say that $\vv P$ is \defn{Gale sewn}.

\begin{corollary}\label{cor:primalGaleSewing}
For any neighborly matroid $\cP$ of rank $\rd$ and $n$ elements there is a neighborly matroid $\tilde \cP$ of rank $\rd+2$ with $n+2$ elements, and two distinguished ones $x$ and $y$ among them such that $\tilde\cP/\{x,y\}=\cP$.\qed
\end{corollary}

\iftoggle{bwprint}{%
\begin{figure}[htpb]
\centering
 \subbottom[$\vv X_0=\{\vv x_1,\vv x_2,\vv x_3,\vv x_4\}$,\newline $ \cM_0=\cM(\vv X_0)\simeq\stc{3}$.]{\label{sfig:cyclicGalesewing0}\includegraphics[width=.30\linewidth]{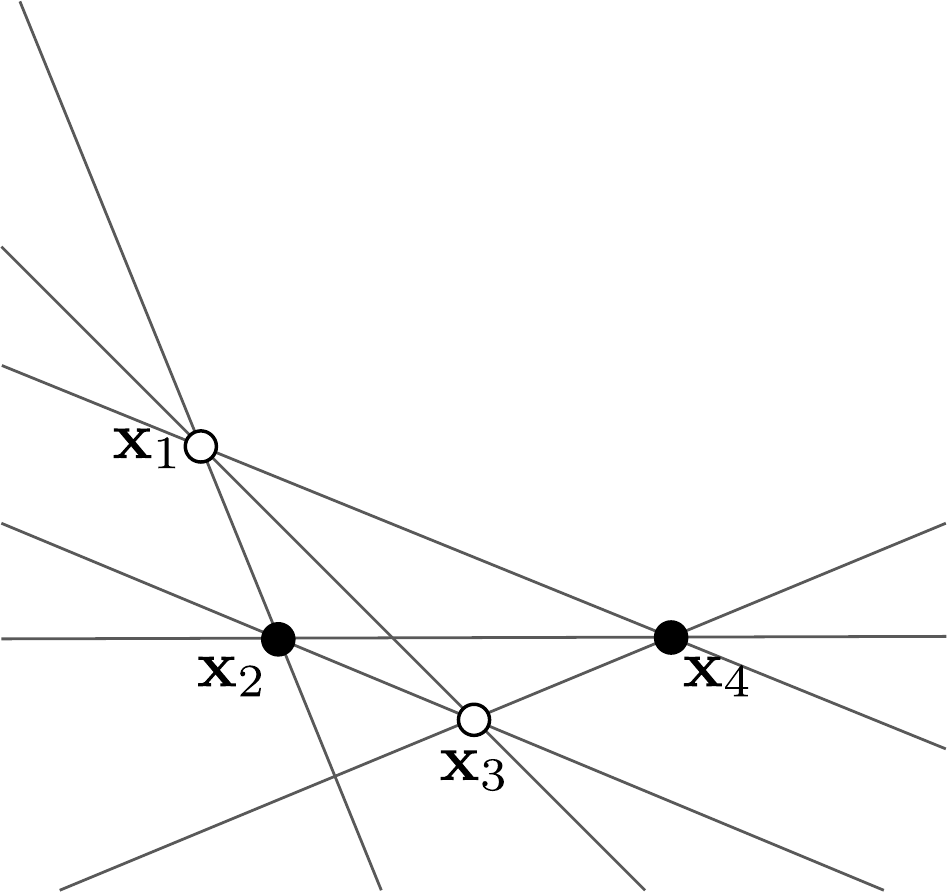}\quad}\, 
 \subbottom[$\cM_1=\cM_0\bracket{\vv x_5}\bracket{\vv x_6}$,\newline $\vv x_5=\bracket{\vv x_4^-,\vv x_3^-,\vv x_2^-}$,\newline $\vv x_6=\bracket{\vv x_5^-,\vv x_4^-,\vv x_3^-}$.]{\label{sfig:cyclicGalesewing1}\includegraphics[width=.30\linewidth]{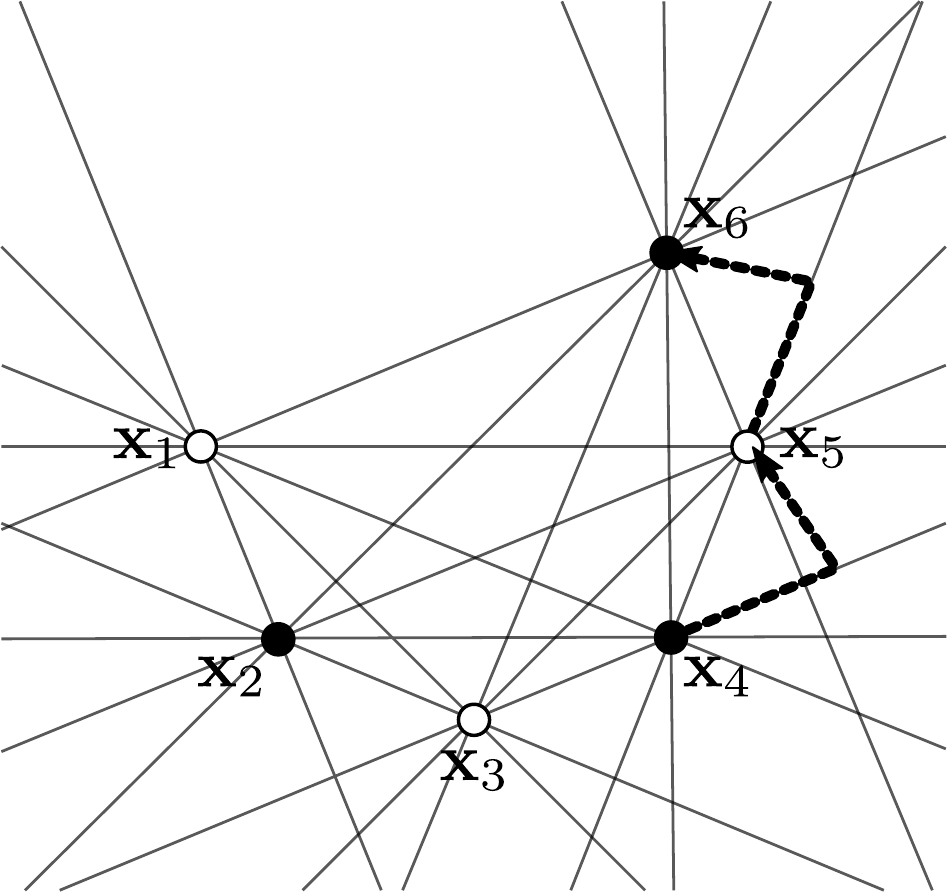}}\, \quad
 \subbottom[$\cM_2=\cM_1\bracket{\vv x_7}\bracket{\vv x_8}$,\newline $\vv x_7=\bracket{\vv x_6^-,\vv x_5^-,\vv x_4^-}$,\newline $\vv x_8=\bracket{\vv x_7^-,\vv x_6^-,\vv x_5^-}$.]{\label{sfig:cyclicGalesewing2}\includegraphics[width=.30\linewidth]{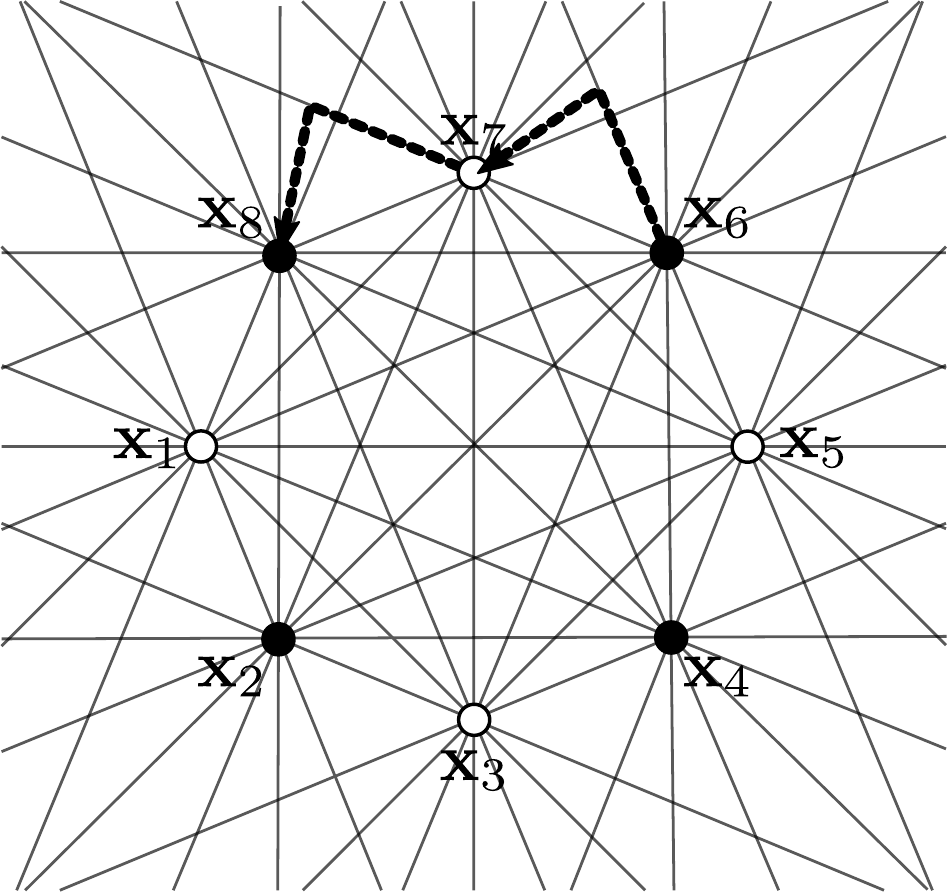}} 
\caption{The process of Gale sewing the dual of $\cyc{4}{8}$. }
 \label{fig:cyclicGalesewing}
\end{figure}}{%
\begin{figure}[htpb]
\centering
 \subbottom[$\vv X_0=\{\vv x_1,\vv x_2,\vv x_3,\vv x_4\}$,\newline $ \cM_0=\cM(\vv X_0)\simeq\stc{3}$.]{\label{sfig:cyclicGalesewing0}\includegraphics[width=.30\linewidth]{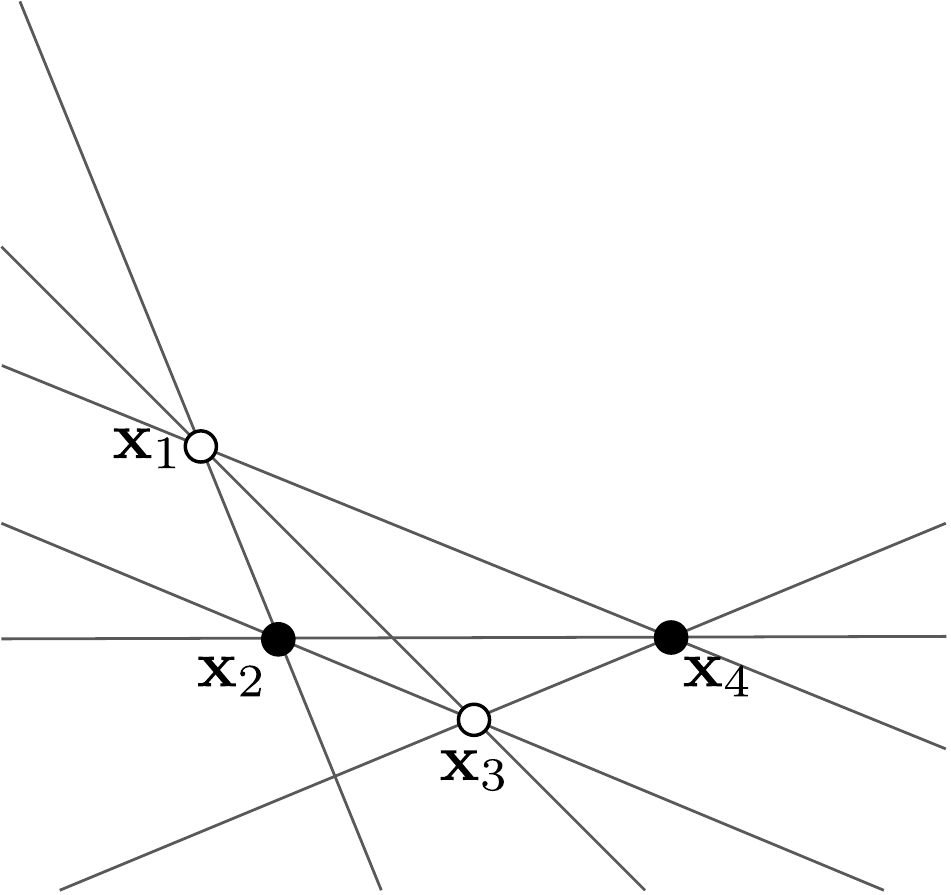}\quad}\, 
 \subbottom[$\cM_1=\cM_0\bracket{\vv x_5}\bracket{\vv x_6}$,\newline $\vv x_5=\bracket{\vv x_4^-,\vv x_3^-,\vv x_2^-}$,\newline $\vv x_6=\bracket{\vv x_5^-,\vv x_4^-,\vv x_3^-}$.]{\label{sfig:cyclicGalesewing1}\includegraphics[width=.30\linewidth]{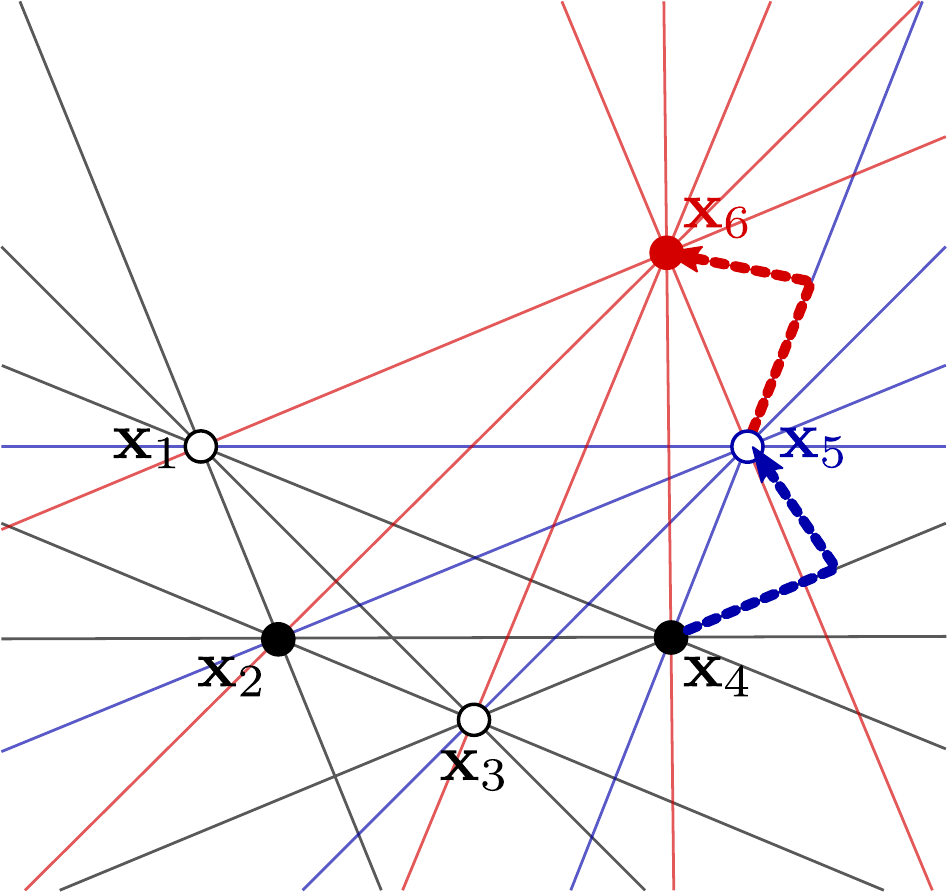}}\, \quad
 \subbottom[$\cM_2=\cM_1\bracket{\vv x_7}\bracket{\vv x_8}$,\newline $\vv x_7=\bracket{\vv x_6^-,\vv x_5^-,\vv x_4^-}$,\newline $\vv x_8=\bracket{\vv x_7^-,\vv x_6^-,\vv x_5^-}$.]{\label{sfig:cyclicGalesewing2}\includegraphics[width=.30\linewidth]{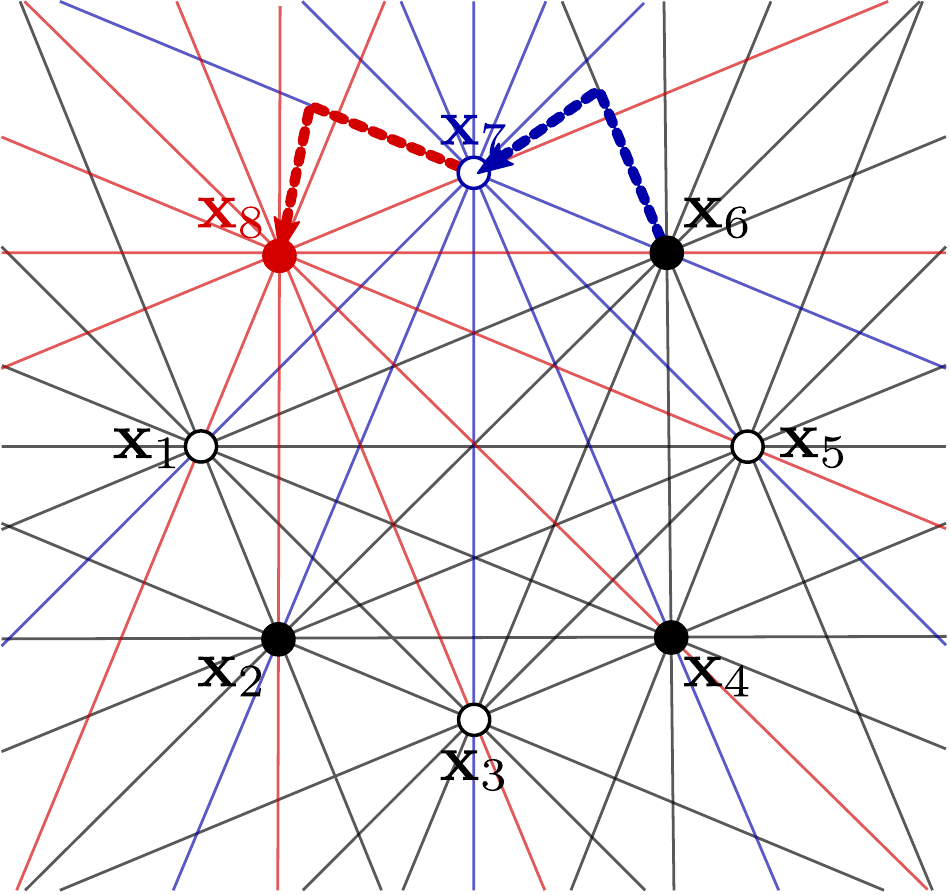}} 
\caption{The process of Gale sewing the dual of $\cyc{4}{8}$. }
 \label{fig:cyclicGalesewing}
\end{figure}}

\begin{example}
In Figure~\ref{fig:cyclicGalesewing} we use the Gale Sewing construction to build the matroid dual to $\cyc{4}{8}$. 
In \ref{sfig:cyclicGalesewing0} we see $\cM_0$, which is isomorphic to~$\stc{3}$, the oriented matroid of $\{\vv e_1,\vv e_2,\vv e_3,-\sum_{i=1}^3 \vv e_i\}$. 
After the first step of Gale Sewing with the extensions by $\vv x_5=\bracket{\vv x_4^-,\vv x_3^-,\vv x_2^-}$ and $\vv x_6=\bracket{\vv x_5^-,\vv x_4^-,\vv x_3^-}$, we obtain~$\cM_1$ in~\ref{sfig:cyclicGalesewing1}. 
This is the Gale dual of an hexagon, and hence it is isomorphic to the configuration of Figure~\ref{fig:affineGale}. 

After doing the lexicographic extensions of~$\cM_1$ by $\vv x_7=\bracket{\vv x_6^-,\vv x_5^-,\vv x_4^-}$ and $\vv x_8=\bracket{\vv x_7^-,\vv x_6^-,\vv x_5^-}$, we obtain~$\cM_2$ in~\ref{sfig:cyclicGalesewing2}. It is the Gale dual of $\cyc{4}{8}$. 
Observe that it is not isomorphic to the configuration constructed with Gale Sewing in Figure~\ref{fig:galesewingexample}, because they do not have the same inseparability graph. Indeed, the inseparability graph of $\cyc{4}{8}$ is the cycle with edges $\vv x_i\sim \vv x_{i+1}$ and $\vv x_1\sim \vv x_8$ (cf. Proposition~\ref{prop:universalflagsofcyclic}), while the inseparability graph of the configuration of Figure~\ref{fig:galesewingexample} is the one depicted in Figure~\ref{fig:inseparable}.

There is a third combinatorial type of neighborly $4$-polytope with $8$~vertices, whose dual was shown in Figure~\ref{fig:N48}. 
We encourage the reader to check that it can also be constructed with the Gale Sewing construction, and that it is not isomorphic to the other two neighborly $4$-polytopes with $8$ vertices that we have presented in this example.
\end{example}

\begin{proposition}[Cyclic polytopes are in $\cG$]\label{prop:cyclicaregalesewn}
 Let $\cM$ be the dual of the alternating matroid of the cyclic polytope $\cyc{d}{n}$, and let $a_1,a_2,\dots, a_n$ be its elements labeled in cyclic order. Then the dual matroid of $\cyc{d+1}{n+1}$ is $\cM[a_{n+1}]$, the single element extension of $\cM$ by $a_{n+1}=[a_{n}^-,a_{n-1}^-,\dots,a_{d}^-]$. 

\end{proposition}
\begin{proof}
We use the following characterization of the circuits of the alternating matroid of rank~$r$
(cf. \cite[Section 9.4]{OrientedMatroids1993}):
the circuits $X$ and $Y$ supported by the $r+1$ elements $x_1<x_2<\dots<x_{r+1}$ (sorted in cyclic order) are those such that $X(x_i)=(-1)^i$ and $Y(x_i)=(-1)^{i+1}$.

If $C$ is a cocircuit of $\cM[a_{n+1}]$ (hence a circuit of its dual) such that $C(a_{n+1})\neq 0$, the signature of the lexicographic extension implies that $C(a_{n+1})$ is opposite to the sign of the largest non-zero element. And thus, by the characterization above, $\cM[a_{n+1}]$ is dual to $\cyc{d+1}{n+1}$. 
\end{proof}

\subsection{Subpolytopes of Gale sewn polytopes}

Our next goal is Proposition~\ref{prop:allquotientsofGalesewnareGalesewn}, that states that subpolytopes (convex hull of subsets of vertices) 
of Gale sewn polytopes are also Gale sewn polytopes.

Its proof uses Proposition~\ref{prop:allquotientsofle}, Lemma~\ref{lem:quotientsofGalesewnareGalesewn} and Lemma~\ref{lem:galesewingorder} below. This lemma shows that when Gale Sewing, the roles of $a_1$, $p$ and $q$ can be exchanged.
Indeed, the isomorphism in~\eqref{eq:changepq} implies that we can switch the roles of $p$ and $q$, while the isomorphism in~\eqref{eq:changepa1} shows how $a_1$ can also be considered as one of the sewn elements.  

\begin{lemma}\label{lem:galesewingorder}
Let $\cM$ be a uniform oriented matroid on a ground set $E$, and consider the lexicographic extensions by
\begin{align*}
 p&=[a_1^{\ep_1},\dots,a_r^{\ep_r}],& q&=[p^-,a_1^{-},\dots,a_{r-1}^{-}];\\
 p'&=[a_1^{-\ep_1},\dots,a_r^{-\ep_r}],& q'&=[p'^-,a_1^{-},\dots,a_{r-1}^{-}];\\
 p''&=[a_1^{+},a_2^{-\ep_1\ep_2},\dots,a_r^{-\ep_1\ep_r}],& q''&=[p''^-,a_1^{-},\dots,a_{r-1}^{-}].
\end{align*}
Then
\vspace{-.5cm}
 \begin{align}
\cM[p][q]&\stackrel{\varphi}{\simeq}\cM[p'][q'],\text{ and }\label{eq:changepq}\\  \cM[p][q]&\stackrel{\psi}{\simeq}\cM[p''][q''],  \label{eq:changepa1} 
 \end{align}
where the bijection $\varphi:E\cup\{p,q\}\rightarrow E\cup\{p',q'\}$ is 
\[\varphi(p)=q',\, \varphi(q)=p'\text{ and }\varphi(e)=e\text{ for }e\in E;\]
when $\ep_1=+$, $\psi:E\cup\{p,q\}\rightarrow E\cup\{p'',q''\}$ is defined as
\[\psi(p)=a_1,\,\,
\psi(q)=q'',\,\,\psi(a_1)=p''\text{ and }\psi(e)=e\text{ for }e\in E\setminus\{a_1\},
\]
 and when $\ep_1=-$, as
\[\psi(p)=q'',\,\,
\psi(q)=a_1,\,\,
\psi(a_1)=p''
\text{ and }
\psi(e)=e\text{ for }e\in E\setminus\{a_1\}.
\]
\end{lemma}
\begin{proof}
We start proving that $\cM[p][q]\stackrel{\varphi}{\simeq}\cM[p'][q']$. For every cocircuit $C\in\co(\cM[p][q])$ we want to find a cocircuit $C'\in\co(\cM[p'][q'])$ with $C'(\varphi(a))=C(a)$ for all $a\in E\cup\{p,q\}$. That is $C'(p')=C(q)$, $C'(q')=C(p)$ and $C'(e)=C(e)$ for $e\in E$. Let $D$ be the restriction of $C$ to $E$.

If $C(q)\neq 0$ and $C(p)\neq 0$, let $i$ be minimal with $D(a_i)\neq 0$. By construction, $C(p)=\ep_iD(a_i)$ and $C(q)=-C(p)=-\ep_iD(a_i)$. By the definition of $\cM[p'][q']$, there is a cocircuit $C'$ that expands $D$, with $C'(p')=-\ep_iD(a_i)=C(q)$ and $C'(q')=-C(p')=\ep_iD(a_i)=C(p)$.

To deal with the case when $C(q)= 0$ or $C(p)= 0$, we use Proposition~\ref{prop:allquotientsofle} to see that $\cM[p][q]/p{\simeq}\cM[p'][q']/q'$ and $\cM[p][q]/q{\simeq}\cM[p'][q']/p'$.
\\

To prove that $\cM[p][q]\stackrel{\psi}{\simeq}\cM[p''][q'']$ we assume that $\ep_1=+$ (otherwise use~\eqref{eq:changepq} to exchange $p$ with $q$).
 In this case we prove that \(\cM[p]\simeq\cM[p'']\), which implies~\eqref{eq:changepa1} because when $p$ and $a_1$ are $(-1)$-inseparable the lexicographic extensions by $[p^-,a_1^{-},\dots,a_{r-1}^{-}]$ and $[a_1^-,p^-,\dots,a_{r-1}^{-}]$ coincide.

For every cocircuit $C\in\co(\cM[p])$ we want to find a cocircuit $C''\in\co(\cM[p''])$ with $C''(p'')=C(a_1)$, $C''(a_1)=C(p'')$ and $C''(e)=C(e)$ for $e\in E\setminus\{a_1\}$. Again, let $D$ be the restriction of $C$ to $E$.

If $C(p)\neq 0$ and $C(a_1)\neq 0$, then $C(p)=C(a_1)=D(a_1)$. Moreover, $D$ is also expanded to a cocircuit $C''$ of $\cM[p'']$ with $C''(p'')=C''(a_1)=D(a_1)$. For circuits with $C(a_1)=0$, observe that $\cM[p]/a_1=\cM[a_2^{\ep_2},\dots,a_r^{\ep_r}]\simeq\cM[p'']/p''$ by Proposition~\ref{prop:allquotientsofle}. Finally, if $C(p)= 0$ then, again by Proposition~\ref{prop:allquotientsofle}, $\cM[p]/p\simeq\cM[a_2^{-\ep_2},\dots,a_r^{-\ep_r}]=\cM[p'']/a_1$.
\end{proof}

 With this lemma we have the last ingredient needed to prove that all the subpolytopes of a Gale sewn polytope are Gale sewn.

\begin{proposition}\label{prop:allquotientsofGalesewnareGalesewn}
 If $\vv P$ is a neighborly polytope in $\cG$, and $\vv a$ is a vertex of~$\vv P$, then $\vv Q=\conv(\verts(\vv P)\setminus \vv a)$ is also a neighborly polytope in $\cG$.
\end{proposition}
\begin{proof}
Let $\cP$ be the oriented matroid of $\vv P$ and $e$ the element of $\cP$ corresponding to the vertex $\vv a$. Observe that $\cP\setminus e$ is the oriented matroid of~$\vv Q$.
The proof is by induction on the rank of $\cP$. When $\cP$ has rank $0$ then $\Gale\cP=\cD_r$ and $\pGale{\cP\setminus e}=\cD_r/e=\cD_{r-1}$.

Otherwise, let $\cM=\Gale\cP$. That $\cP$ belongs to $\cG$ means that there is a matroid $\cN$ whose dual~$\Gale\cN$ is in $\cG$, such that $\cM=\cN[p][q]$ where $p=[a_1^{\ep_1},a_2^{\ep_2},\dots,a_r^{\ep_r}]$ and $q=[p^-,a_1^{-},\dots,a_{r-1}^{-}]$.

We will prove that for every $e\in\cP$, there is some $\tilde e\in \cN$ fulfilling
 \begin{equation}\label{eq:contracte}\pGale{\cP\setminus e}\simeq (\cN/ \tilde e)[\tilde p][\tilde q] =\pGale{\Gale\cN\setminus \tilde e}[\tilde p][\tilde q],\end{equation} 
where $\tilde p=[\tilde a_1^{\tilde \ep_1},\tilde a_2^{\tilde \ep_2},\dots,\tilde a_r^{\tilde \ep_r}]$ and $\tilde q=[\tilde p^-,\tilde a_1^{-},\dots,\tilde a_{r-1}^{-}]$ for some $\tilde a_i$'s and~$\tilde \ep_i$'s. Since $\rank\left( \Gale\cN\right)=\rank\left( \cP\right) -2$, by the induction hypothesis $(\Gale\cN\setminus e')\in \cG$ and our claim follows directly from \eqref{eq:contracte}. 

 If $e=q$, then by Lemma~\ref{lem:quotientsofGalesewnareGalesewn} we know that $\pGale{\cP\setminus e}=\left(\cN[p][q]\right)/q\simeq\left(\cN/a_1\right)[\tilde p][\tilde q]$, where $\tilde p=[a_2^{-\ep_1\ep_2},\dots,a_r^{-\ep_1\ep_r}]$ and $\tilde q=[\tilde p^{-},\dots,a_{r-1}^{-}]$. The case $e=p$ is analogous because of Lemma~\ref{lem:galesewingorder}.
 If $e=a_i$, then $\pGale{\cP\setminus e}\simeq\left(\cN/a_i\right)[\tilde p][\tilde q]$, 
where $\tilde p$ and $\tilde q$ have the same signature as $\tilde p$ and $\tilde q$ but omitting the element $a_i$.
 For the remaining elements $e$, $\pGale{\cP\setminus e}\simeq\left(\cN/e\right)[\tilde p][\tilde q]$.
\end{proof}

\section{\texorpdfstring{Combinatorial description of $\cG$}{Combinatorial description of G}}\label{sec:descriptioncG}

Let $\vv P$ be a simplicial polytope that defines an acyclic uniform oriented matroid $\cP$, and let $\cM:=\Gale\cP$ be its dual matroid. 
The essence of Gale Sewing is to construct a new polytope $\vv{\tilde P}$ whose matroid $\tilde \cP$ is dual to $\tilde\cM:=\cM[p]$, a lexicographic extension of $\cM$ by $p=[a_1^{\ep_1},a_2^{\ep_2},\dots,a_k^{\ep_k}]$. In this section we will see that the combinatorics of $\vv{\tilde P}$ are described by \defn{lexicographic triangulations} of $\vv P$.
 \\

Let $\vv A=\{\vv a_1,\dots,\vv a_n\}$ be the set of vertices of $\vv P\subset \RR^d$. Let $M$ be the $d\times n$ matrix whose columns list the coordinates of the $\vv a_i$'s:

\[
\renewcommand{\arraystretch}{1.5}
       M := 
       \left[
\begin{matrix}
& \kern.4em\vrule height 2ex\kern.2em & \kern.2em\vrule height 2ex\kern.2em& & \vrule height 2ex\kern.6em &\\
&\kern.2em\vv a_1&\vv a_2&\dots&\vv a_n&\\
& \kern.4em\vrule depth 0ex\kern.2em & \kern.2em\vrule depth 0ex\kern.2em& & \vrule depth 0ex\kern.6em &\\
\end{matrix}
 	\right].	
\]

Then there is some $\delta>0$ such that the point configuration ${\tilde {\vv A}}$ defined by the columns of the following $(d+1)\times(n+1)$ matrix $\tilde M$ is a realization of the set of vertices of~${\tilde {\vv P}}$:
\[
\renewcommand{\arraystretch}{1.5}
      \tilde M := 
\kbordermatrix{
&\tilde{\vv a}_1&\tilde{\vv a}_2&\tilde{\vv a}_3&\dots&\tilde{\vv a}_k&\tilde{\vv a}_{k+1}&\dots &\tilde{\vv a}_n& \omit\vrule &\vv p \\
&\kern.4em\vrule height 2ex\kern.4em & \kern.4em\vrule height 2ex\kern.4em & \kern.4em\vrule height 2ex\kern.4em& & \vrule height 2ex\kern.6em & \vrule height 2ex\kern.6em & & \kern.4em\vrule height 2ex\kern.4em & \omit\vrule & \kern.4em\vrule height 2ex\kern.4em \\
&\kern.4em\vv a_1&\vv a_2&\vv a_3&\dots&\vv a_k&\vv a_{k+1}&\dots &\vv a_n& \omit\vrule & \veczero \\
&\kern.4em\vrule depth 1ex\kern.4em & \kern.4em\vrule depth 1ex\kern.4em & \kern.4em\vrule depth 1ex\kern.4em& & \vrule depth 1ex\kern.6em & \vrule depth 1ex\kern.6em & & \kern.4em\vrule depth 1ex\kern.4em& \omit\vrule & \kern.4em\vrule depth 1ex\kern.4em \\
\cline{2-11}
&-\ep_1&-\ep_2 \delta &-\ep_3\delta^2& \dots &-\ep_k\delta^{k-1}&0\kern.4em&\dots&0& \omit\vrule &\kern.2em 1\kern.2em
}.
\]

Geometrically, each of the points $\vv a_i\in \vv A\subset \RR ^d$ is lifted to a point $\tilde {\vv a}_i\in\tilde{\vv A}\subset \RR^{d+1}$ with a height that depends on the signature of the lexicographic extension. Namely, $\tilde{\vv a}_i=\binom{\vv a_i}{ -\ep_i\delta^{i-1}}$ for $i\leq k$ and $\tilde{\vv a}_i=\binom{\vv a_i}{0}$ otherwise. Moreover, $\vv p$ is added to $\tilde{\vv A}$ with coordinates $\binom{\veczero}{1}$.
The vertex figure of $\vv p$ in~${\tilde {\vv P}}$ is combinatorially equivalent to $\vv P$. That is, the faces of ${\tilde{\vv P}}$ that contain $\vv p$ are isomorphic to pyramids over faces of $\vv P$. 
On the other hand, the faces of~${\tilde {\vv P}}$ that do not contain $\vv p$ correspond to faces of a regular subdivision of $\vv P$: the \defn{lexicographic subdivision} of $\vv P$ on $[\vv a_1^{-\ep_1},\vv a_2^{-\ep_2},\dots,\vv a_k^{-\ep_k}]$. 
When $\vv a_1\dots \vv a_k$ form a basis, this subdivision is a triangulation. A concrete example is depicted in~Figure~\ref{fig:lifting}.\\

\iftoggle{bwprint}{%
\begin{figure}[htpb]
\centering
 \subbottom[$\vv P=\conv({\vv A})$]{\includegraphics[width=.25\linewidth]{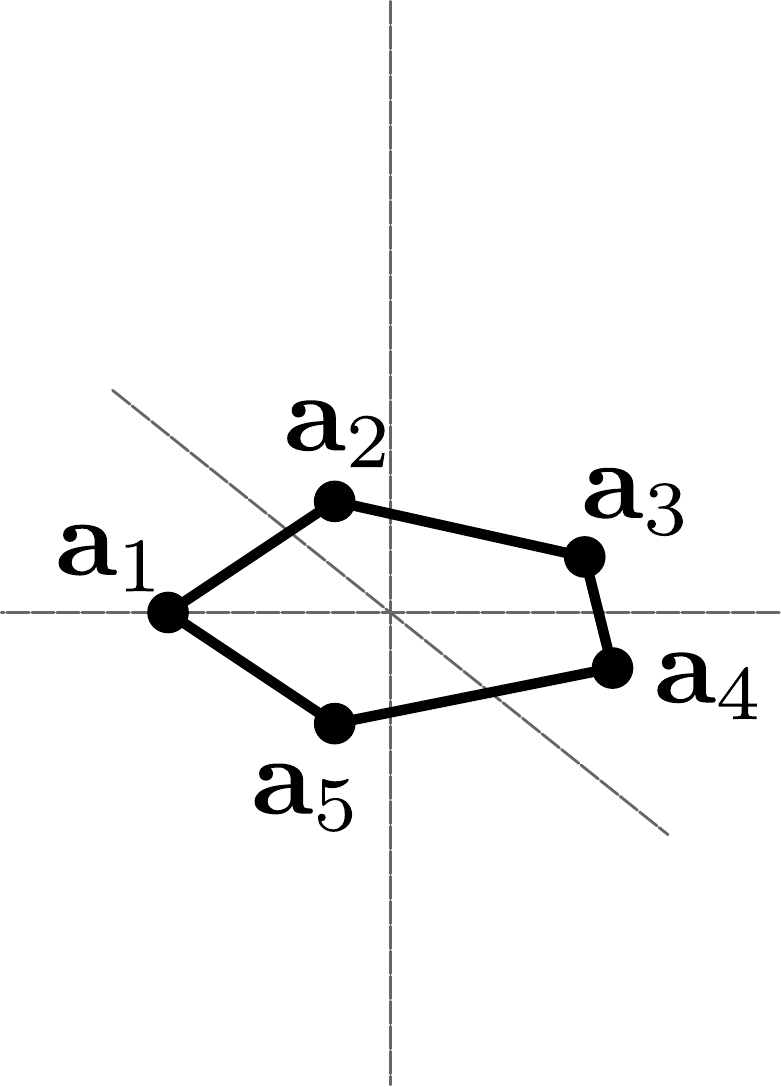}}\qquad\quad 
 \subbottom[Lifting $\vv A$ to $\tilde{\vv A}$]{\includegraphics[width=.25\linewidth]{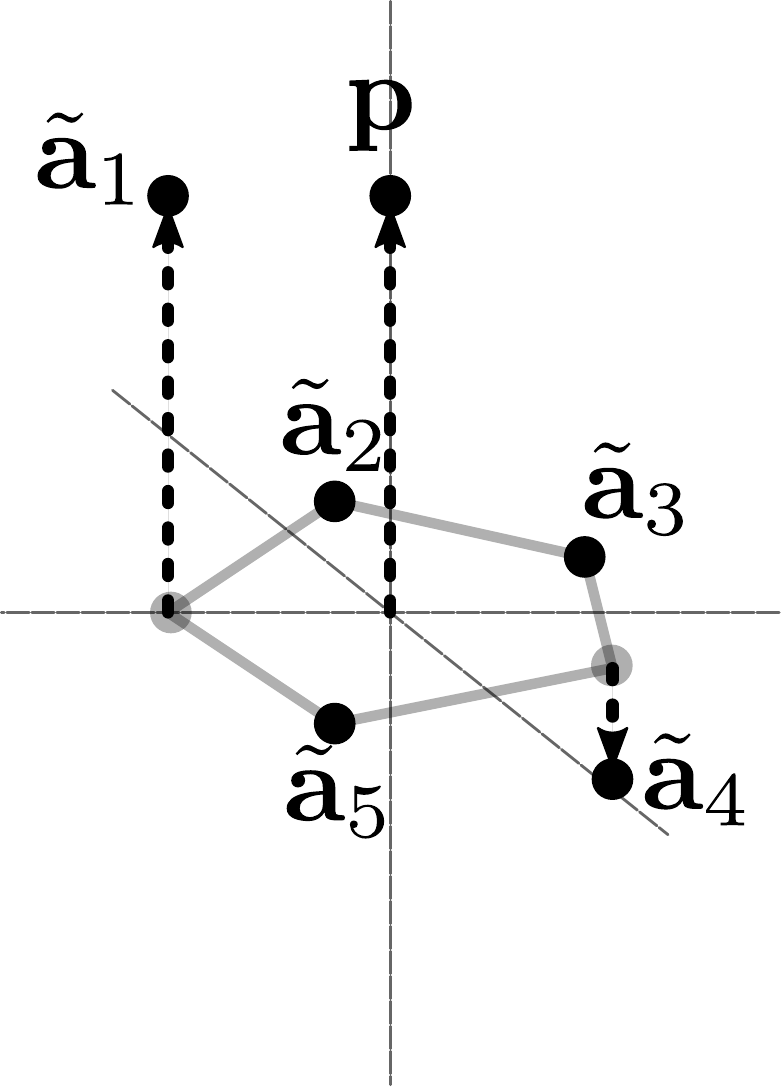}}\qquad\quad 
 \subbottom[$\tilde {\vv P}=\conv(\tilde {\vv A})$]{\includegraphics[width=.25\linewidth]{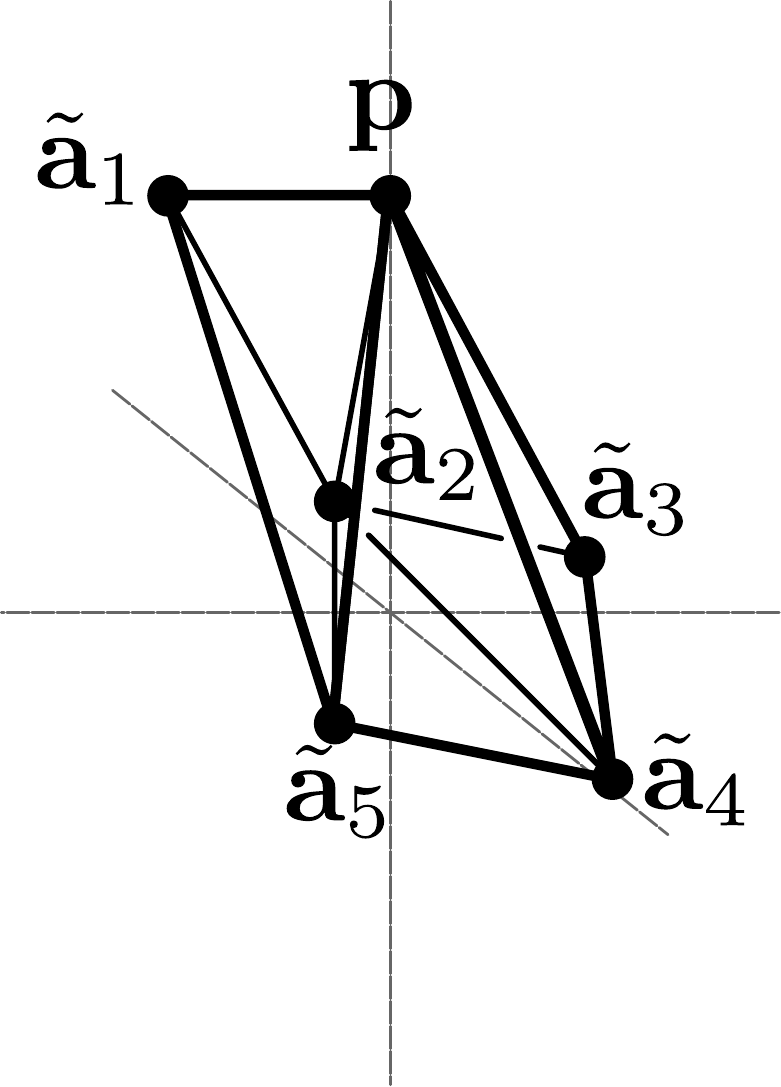}}\\ 
 \subbottom[Upper envelope of~$\tilde{\vv P}$]{\quad\quad\includegraphics[width=.22\linewidth]{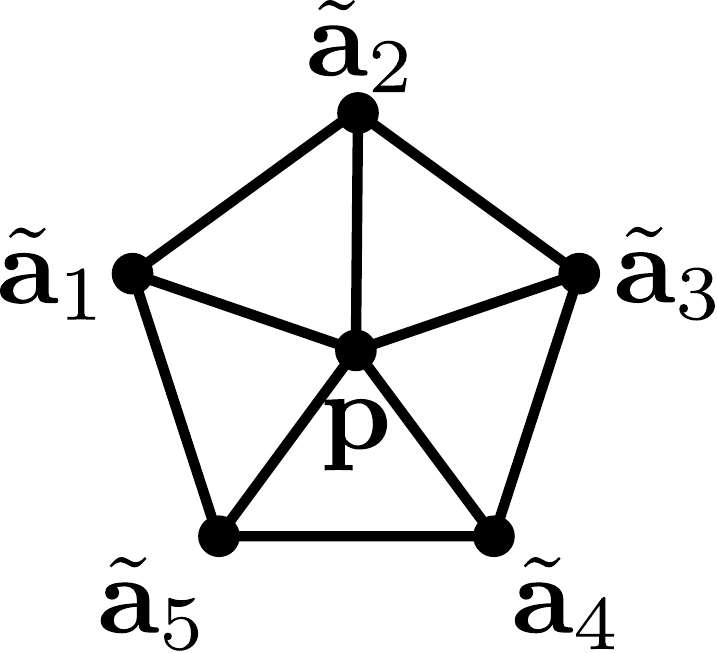}\quad\quad}\qquad 
 \subbottom[Lower envelope of~$\tilde{\vv P}$]{\quad\quad\includegraphics[width=.22\linewidth]{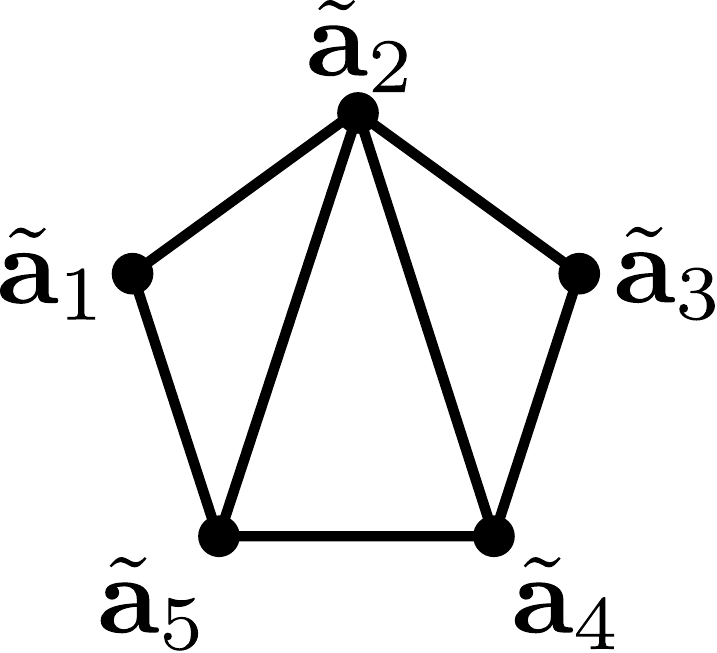}\quad\quad}
 \caption[Lifting triangulation]{Lifting the pentagon $\vv P=\conv({\vv A})$ to $\vv{\tilde P}=\conv(\vv{\tilde A})$ when $\Gale{\vv{\tilde A}}=\Gale{\vv A}[\vv p]$ and $\vv p=[\vv a_1^-,\vv a_4^+]$. Its upper envelope are pyramids over facets of $\vv P$, while the lower envelope is the lexicographic triangulation of $\vv P$ on $[\vv a_1^+,\vv a_4^-]$.}
 \label{fig:lifting}
\end{figure}
}{%
\iftoggle{print}{%
\begin{figure}[htpb]
\centering
 \subbottom[$\vv P=\conv({\vv A})$]{\includegraphics[width=.25\linewidth]{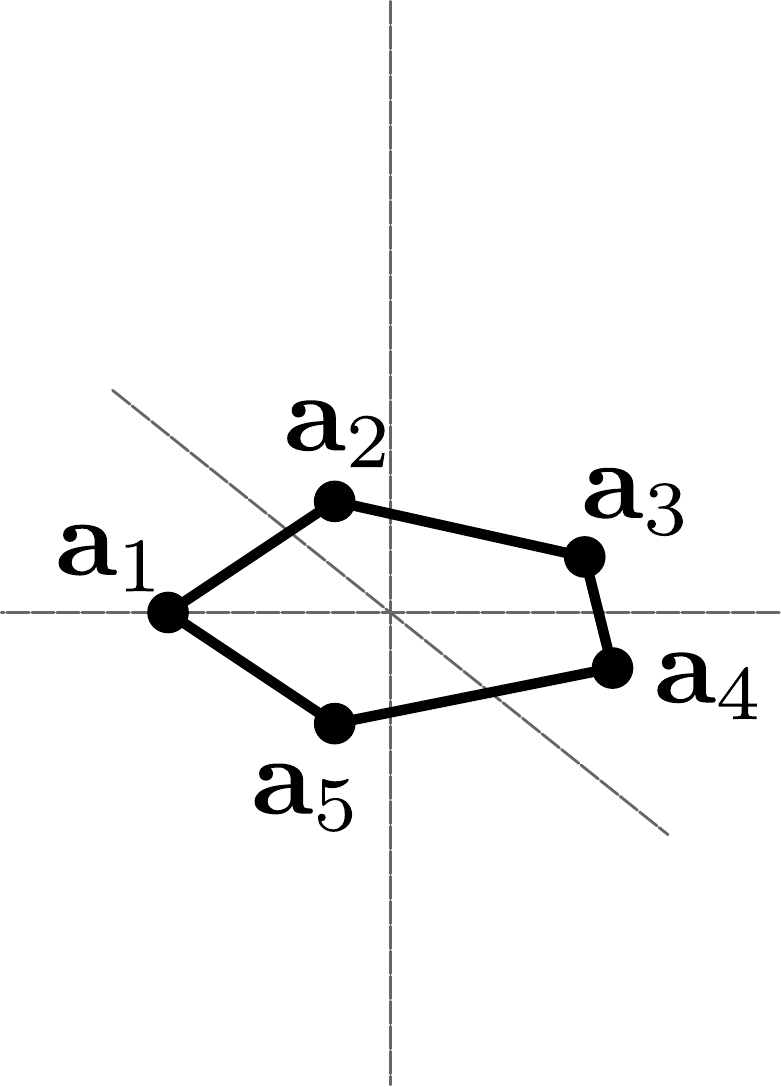}}\qquad\quad 
 \subbottom[Lifting $\vv A$ to $\tilde{\vv A}$]{\includegraphics[width=.25\linewidth]{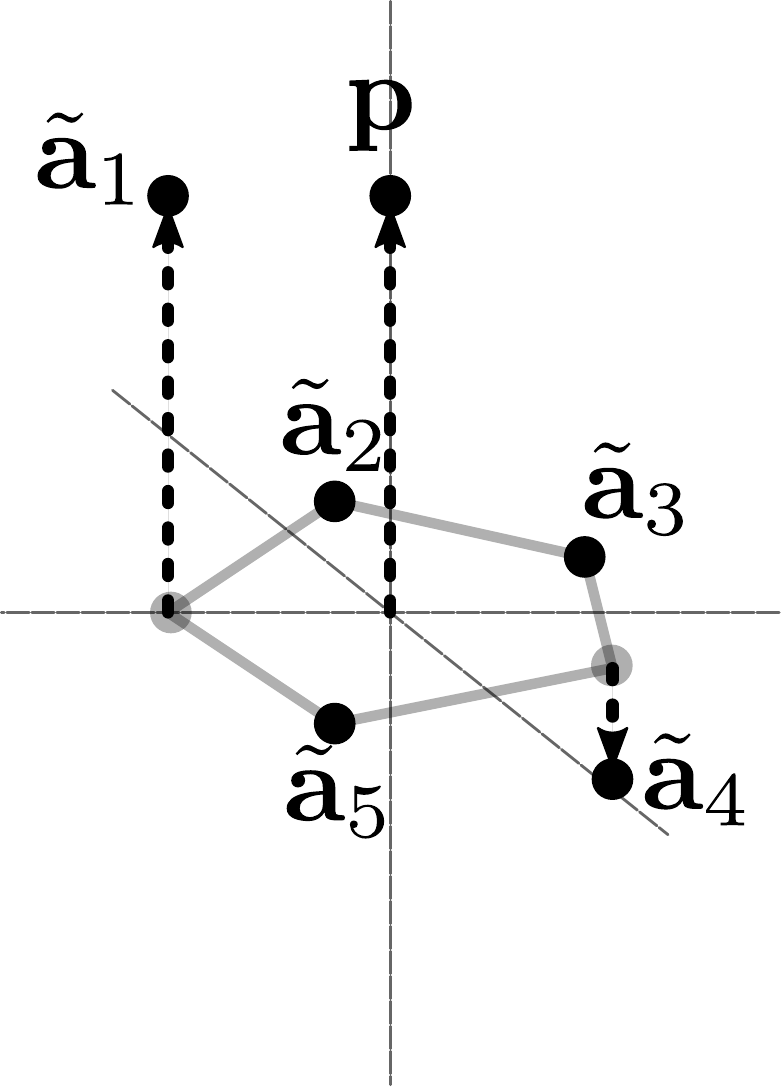}}\qquad\quad 
 \subbottom[$\tilde {\vv P}=\conv(\tilde {\vv A})$]{\includegraphics[width=.25\linewidth]{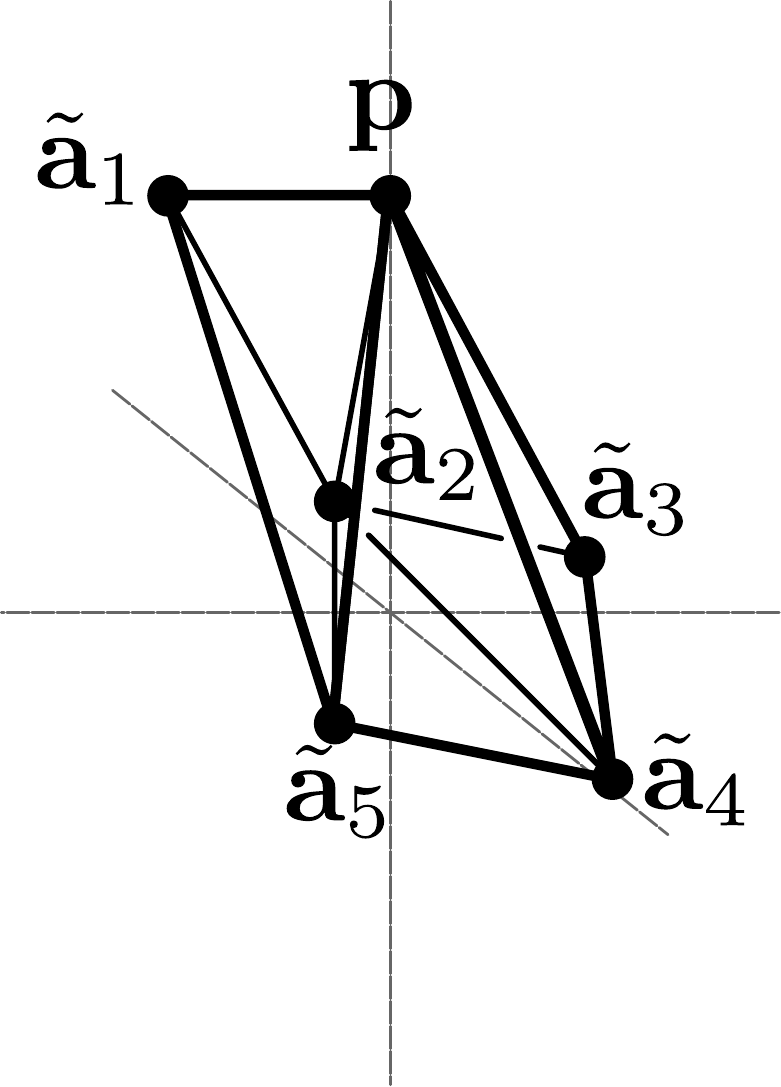}}\\ 
 \subbottom[Upper envelope of~$\tilde{\vv P}$]{\quad\quad\includegraphics[width=.22\linewidth]{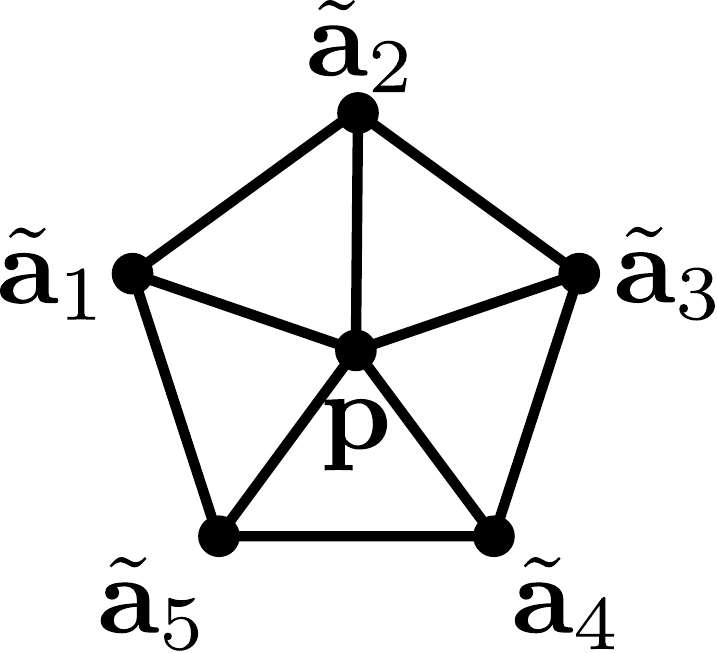}\quad\quad}\qquad 
 \subbottom[Lower envelope of~$\tilde{\vv P}$]{\quad\quad\includegraphics[width=.22\linewidth]{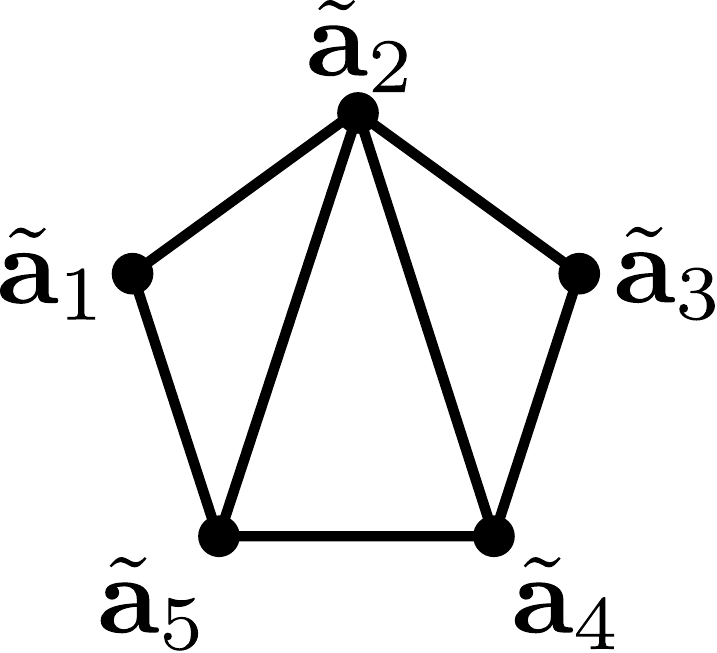}\quad\quad}
 \caption[Lifting triangulation]{Lifting the pentagon $\vv P=\conv({\vv A})$ to $\vv{\tilde P}=\conv(\vv{\tilde A})$ when $\Gale{\vv{\tilde A}}=\Gale{\vv A}[\vv p]$ and $\vv p=[\vv a_1^-,\vv a_4^+]$. Its upper envelope are pyramids over facets of $\vv P$, while the lower envelope is the lexicographic triangulation of $\vv P$ on $[\vv a_1^+,\vv a_4^-]$.}
 \label{fig:lifting}
\end{figure}
}{%
\begin{figure}[htpb]
\centering
 \subbottom[$\vv P=\conv({\vv A})$]{\includegraphics[width=.25\linewidth]{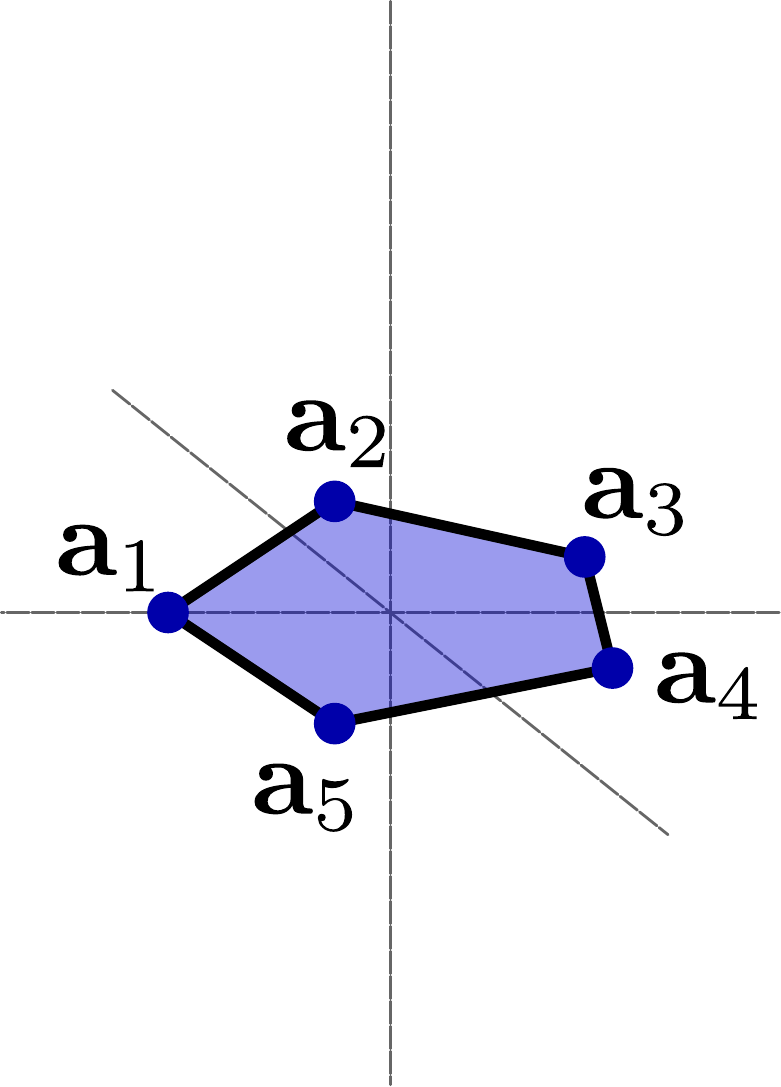}}\qquad\quad 
 \subbottom[Lifting $\vv A$ to $\tilde{\vv A}$]{\includegraphics[width=.25\linewidth]{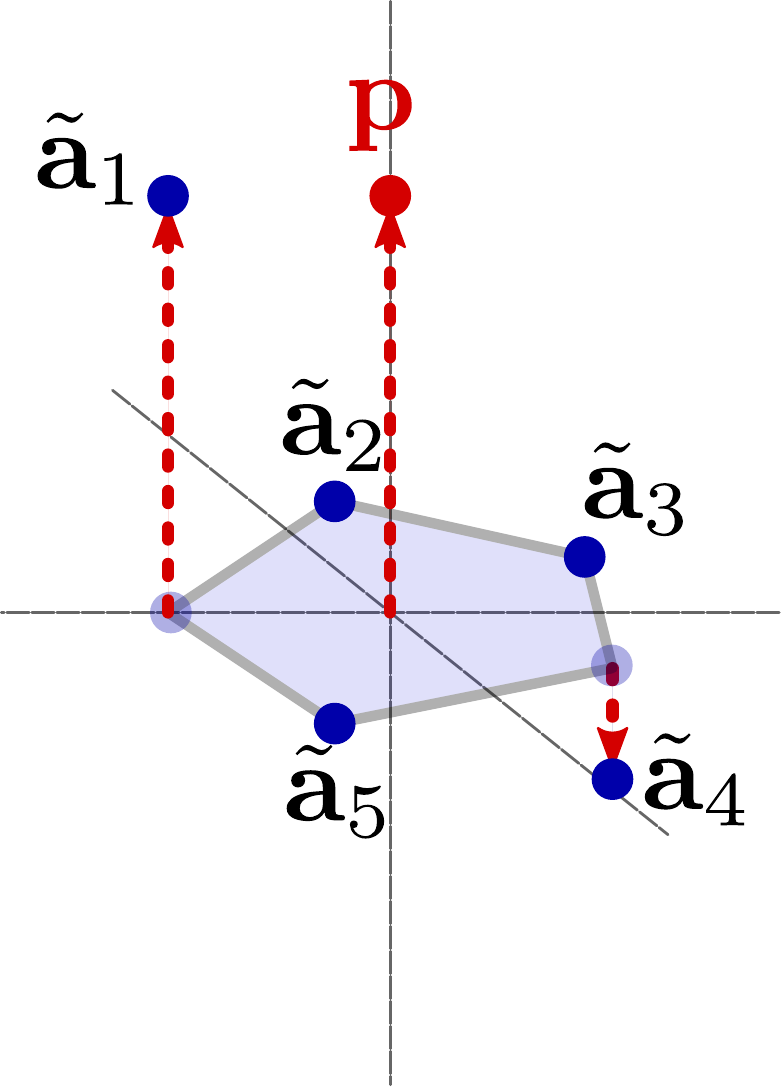}}\qquad\quad 
 \subbottom[$\tilde {\vv P}=\conv(\tilde {\vv A})$]{\includegraphics[width=.25\linewidth]{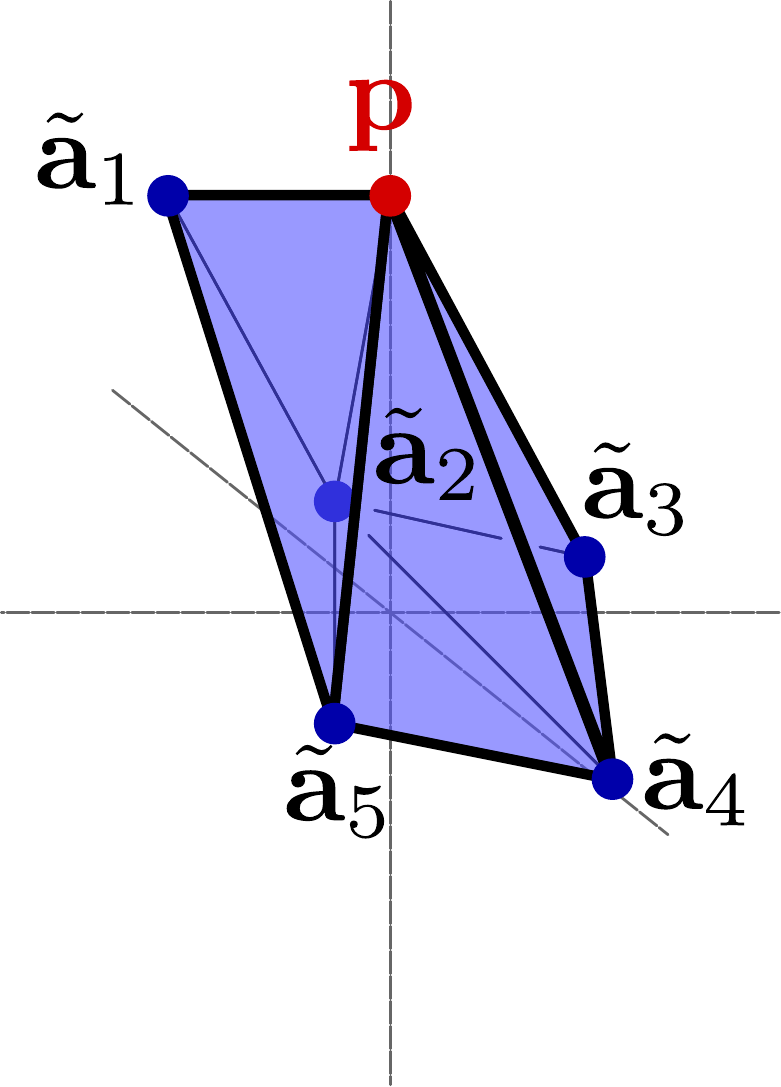}}\\ 
 \subbottom[Upper envelope of~$\tilde{\vv P}$]{\quad\quad\includegraphics[width=.22\linewidth]{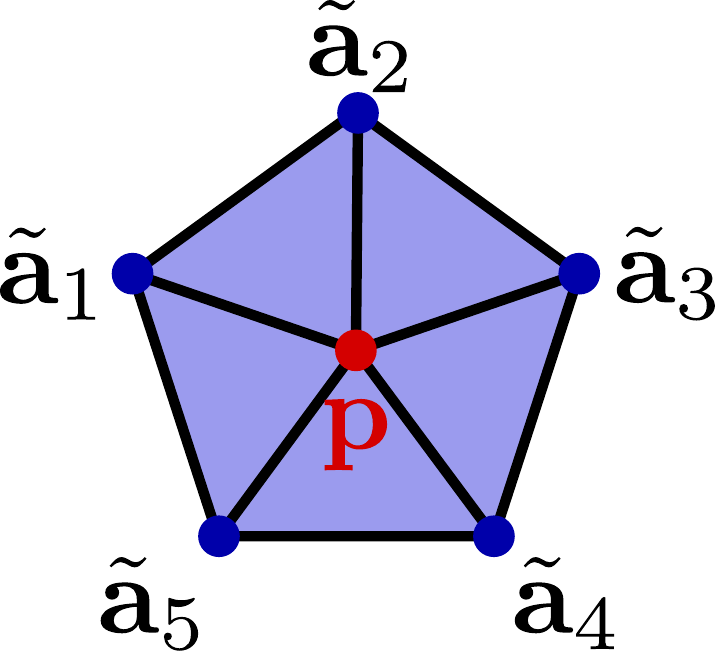}\quad\quad}\qquad 
 \subbottom[Lower envelope of~$\tilde{\vv P}$]{\quad\quad\includegraphics[width=.22\linewidth]{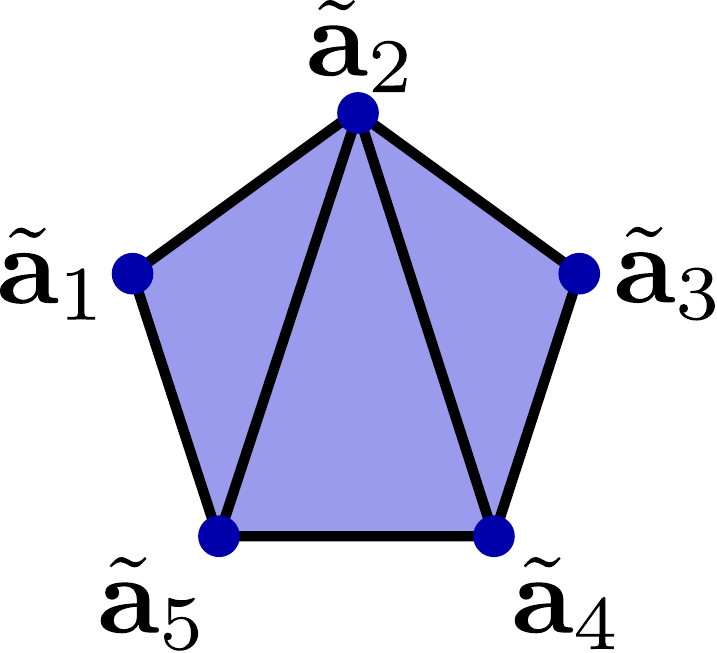}\quad\quad}
 \caption[Lifting triangulation]{Lifting the pentagon $\vv P=\conv({\vv A})$ to $\vv{\tilde P}=\conv(\vv{\tilde A})$ when $\Gale{\vv{\tilde A}}=\Gale{\vv A}[\vv p]$ and $\vv p=[\vv a_1^-,\vv a_4^+]$. Its upper envelope are pyramids over facets of $\vv P$, while the lower envelope is the lexicographic triangulation of $\vv P$ on $[\vv a_1^+,\vv a_4^-]$.}
 \label{fig:lifting}
\end{figure}
}
}

Our formulation of the definition of lexicographic subdivision is based on~\cite{DeLoeraRambauSantosBOOK}. However we use a different ordering, the same as in~\cite{Santos2002}, that mirrors the definition of lexicographic extension (with opposite signs). See also~\cite{Lee1991}.\index{lexicographic subdivision}

\begin{definition}
\label{def:lexsubdivision}
Let $\vv P$ be a $d$-polytope with $n$ vertices $\vv a_1,\dots,\vv a_n$. The \defn{lexicographic subdivision} of $\vv P$ on $[\vv a_1^{\ep_1},\vv a_2^{\ep_2},\dots,\vv a_k^{\ep_k}]$, where $\ep_i=\pm 1$, is defined recursively as follows.
\begin{itemize}
 \item If $\ep_1=+1$ (\defn{pushing}), then the lexicographic subdivision of $\vv P$ is the union of the lexicographic subdivision of $\vv P\setminus \vv a_1$ on $[\vv a_2^{\ep_2},\dots,\vv a_k^{\ep_k}]$, and the simplices joining $\vv a_1$ to the (lexicographically subdivided) faces of $\vv P\setminus \vv a_1$ visible from it.
 \item If $\ep_1=-1$ (\defn{pulling}), then the lexicographic subdivision of $\vv P$ is the unique subdivision in which every maximal cell contains $\vv a_1$ and which, restricted to each proper face $\vv F$ of $\vv P$, coincides with the lexicographic subdivision of that face on $[\vv a_2^{\ep_2},\dots,\vv a_k^{\ep_k}]$.
\end{itemize}
\end{definition}

\begin{remark}
The resemblance with Sanyal and Ziegler's description of the vertex figures of the neighborly cubical polytopes in~\cite{SanyalZiegler2010} is not a coincidence. 
Indeed, all the Gale duals of those vertex figures are lexicographic extensions of the dual of a fixed neighborly polytope. Namely, every vertex of their projected deformed cubes is indexed by a string of signs $\vv \gs\in\{+,-\}^n$. 
In \cite[Theorem 3.7]{SanyalZiegler2010}, they show that the coordinates of the Gale dual $\vv V$ of the vertex figure of the vertex indexed by $\vv \gs$ are given by the columns of the following matrix

\[
\renewcommand{\arraystretch}{1.5}
\kbordermatrix{
&{\vv v}_0&\omit\vrule&{\vv v}_1&{\vv v}_2&\dots&{\vv v}_{n-1}&{\vv v}_{n-d+1}&\dots &{\vv v}_{n-d} \\
&-\omega_1&\omit\vrule& 1& & & & \kern.4em\vrule height 2ex\kern.4em & & \kern.4em\vrule height 2ex\kern.4em\\
&-\omega_2&\omit\vrule& & 1& & & \kern.4em\vrule \kern.4em & & \kern.4em\vrule \kern.4em\\
&\vdots&\omit\vrule & & &\ddots & & {\vv g}_1 &\dots&{\vv g}_{d-1}\\
&-\omega_{n-d}&\omit\vrule& & & & 1& \kern.4em\vrule depth 1ex\kern.4em & & \kern.4em\vrule depth 1ex\kern.4em\\
},
\]
where $\omega_i=(-1)^i\delta^i\prod_{j=1}^i\vvc \gs_i$. Here the vectors $\vv g_i$ are the columns of a matrix~$\ol G$ such that the columns of the matrix $\left[\matid_n\,\, \ol G\right]$ are the Gale dual $\vv W$ of a neighborly polytope. 
Observe how the extension by $\vv v_0$ is a lexicographic extension of $\vv W$ by $\vv v_0=[{\vv v}_1^{(-1)\vvc \gs_1},\dots,{\vv v}_i^{(-1)^i\prod_{j=1}^i\vvc \gs_i},\dots]$.
\\

Additionally, in the same paper it is proved that the number 
of combinatorial types of neighborly simplicial $(d - 2)$-polytopes on $n - 1$ vertices is a lower bound for the number of combinatorial types of $d$-dimensional neighborly cubical polytopes with $2^n$ vertices~\cite[Corollary 3.8]{SanyalZiegler2010}. Hence, the bounds that we obtain in Theorem~\ref{thm:lblnei} automatically yield bounds for the number of neighborly cubical polytopes.
\end{remark}

\section{Comparing the constructions}\label{sec:comparing}
In this section we compare the construction techniques for neighborly polytopes, which are strongly related. Our goal is to prove the following theorem. It states that if a neighborly polytope $\vv P$ is built via Extended Sewing and Omitting (Construction~\ref{constr:cO}), then~$\vv P$~can also be built with Gale Sewing (Construction~\ref{constr:cG}):

\begin{theorem}\label{thm:cOsubsetcG}
 $\cO\subseteq\cG$.
\end{theorem}

\noindent \emph{Start of proof.} By Proposition~\ref{prop:allquotientsofGalesewnareGalesewn}, it suffices to see that $\cE\subseteq \cG$. As a first step, we prove that if we sew on a neighborly polytope in $\cG$ (through~a specific universal flag), we obtain a new polytope in $\cG$. 

Let $\cP$ be a neighborly matroid in $\cG$ whose dual matroid $\cM:=\Gale\cP$ has been constructed using Gale Sewing. 
Specifically, $\cM_0=\stc{r}$ (\ie the oriented matroid of the rank~$r$ configuration $\{\vv e_1, \dots, \vv e_r, -\sum_{i=1}^r \vv e_i\}$) with its elements labeled $\{a_0,\dots,a_r\}$, and $\cM=\cM_m$, where 
\begin{align*}
\cM_k&:=\cM_{k-1}[p_k][q_k], 
\end{align*}
the lexicographic extensions by $p_k$ and $q_k$ are given by
\begin{align*}
p_k&:=[a_{k1}^{\ep_{k1}},\dots,a_{kr}^{\ep_{kr}}],& q_k&:=[p_k^-,a_{k1}^{-},\dots,a_{k(r-1)}^{-}];
\end{align*}
and the $a_{ij}$ are pairwise distinct elements of $\cM_{i-1}$. 
\\

\iftoggle{bwprint}{%
\begin{figure}[htpb]
\centering
 \subbottom[$\cM_0=\cM\{\vv a_0,\vv a_1,\vv a_2\}$.]{\quad\quad\quad\includegraphics[width=.25\linewidth]{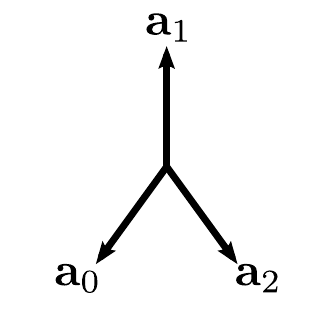}\quad\quad\quad}
\hspace{1cm} 
 \subbottom[$\cM_1=\cM_0\bracket{\vv a_2^-,\vv a_1^-}\bracket{\vv p_1^-,\vv a_2^-}$.]{\quad\quad\quad\includegraphics[width=.25\linewidth]{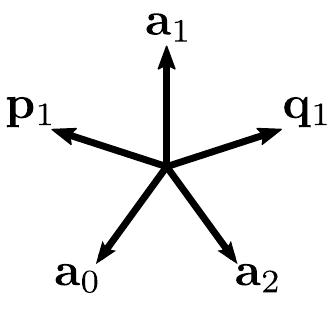}\quad\quad\quad}
\\
 \subbottom[$\cP_0=\cP_1/\{\vv p_1,\vv q_1\}=\Gale\cM_0$.]{\quad\quad\quad\includegraphics[width=.25\linewidth]{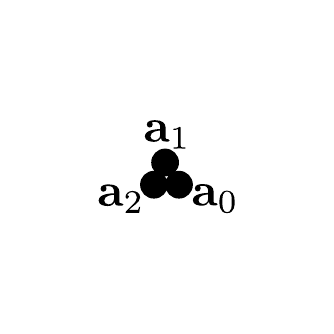}\quad\quad\quad}
\hspace{1cm}  
 \subbottom[$\cP_1=\Gale\cM_1$.]{\quad\quad\quad\includegraphics[width=.22\linewidth]{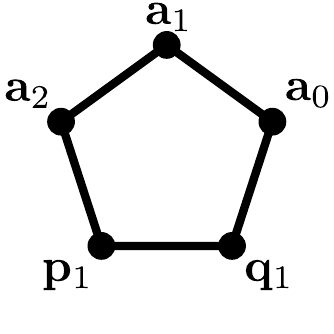}\quad\quad\quad}
 \caption[Primal and dual sewing, part $1$.]{$\cM_1$ is constructed from $\cM_0$ after Gale sewing $\vv p_1=\bracket{\vv a_2^-,\vv a_1^-}$ and $\vv q_1=\bracket{\vv p_1^-,\vv a_2^-}$. The dual of $\cM_1$ is $\cP_1$, whose contraction by $\{\vv p_1,\vv q_1\}$ is isomorphic to $\cM_0$.}
 \label{fig:ShemerVsMe1}
\end{figure}
}{%
\begin{figure}[htpb]
\centering
 \subbottom[$\cM_0=\cM\{\vv a_0,\vv a_1,\vv a_2\}$.]{\quad\quad\quad\includegraphics[width=.25\linewidth]{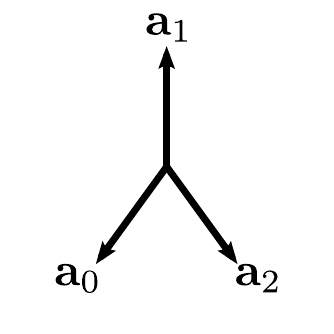}\quad\quad\quad}
\hspace{1cm} 
 \subbottom[$\cM_1=\cM_0\bracket{\vv a_2^-,\vv a_1^-}\bracket{\vv p_1^-,\vv a_2^-}$.]{\quad\quad\quad\includegraphics[width=.25\linewidth]{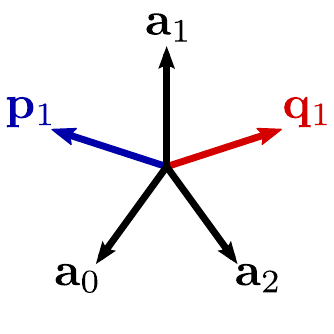}\quad\quad\quad}
\\
 \subbottom[$\cP_0=\cP_1/\{\vv p_1,\vv q_1\}=\Gale\cM_0$.]{\quad\quad\quad\includegraphics[width=.25\linewidth]{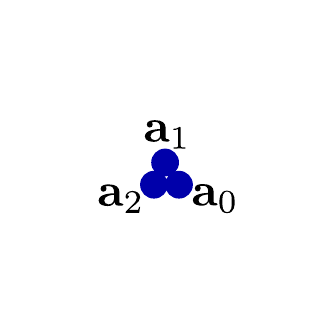}\quad\quad\quad}
\hspace{1cm}  
 \subbottom[$\cP_1=\Gale\cM_1$.]{\quad\quad\quad\includegraphics[width=.22\linewidth]{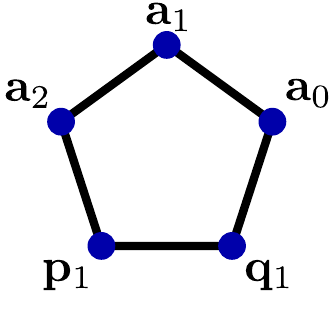}\quad\quad\quad}
 \caption[Primal and dual sewing, part $1$.]{$\cM_1$ is constructed from $\cM_0$ after Gale sewing $\vv p_1=\bracket{\vv a_2^-,\vv a_1^-}$ and $\vv q_1=\bracket{\vv p_1^-,\vv a_2^-}$. The dual of $\cM_1$ is $\cP_1$, whose contraction by $\{\vv p_1,\vv q_1\}$ is isomorphic to $\cM_0$.}
 \label{fig:ShemerVsMe1}
\end{figure}
}

Let $F_j=\bigcup_{i=0}^{j-1} \{p_{m-i},q_{m-i}\}$, so that  $F_{m-k}=\{p_m,q_m,\dots,p_{k+1},q_{k+1}\}$,  and set $\cF=\{F_i\}_{i=1}^{m}$. 
Define $\cP_k:=\cP/F_{m-k}$ for $k=0,\dots,m$, and observe that $\cP_k=\Gale{\cM_k}$ for all~$k$ by deletion-contraction duality. By construction, $\cF$ is a universal flag of $\cP$ because every $\cP_k$ is neighborly.
Let $\cF'$ be some flag that contains $\cF$. Thanks to Lemma~\ref{lem:galesewingorder} we can assume without loss of generality that all split faces in $\cF'$ are $q_i$-split. An example is depicted in Figure~\ref{fig:ShemerVsMe1}.\\

Since $\cP$ is neighborly and $\cF'$ has a universal subflag, we can apply the Extended Sewing Theorem~\ref{thm:extshemersewing} to $\cP$ and $\cF'$. 

Let $\tilde \cP:=\cP[\cF']$, where $p$ is sewn onto~$\cP$ through the flag~$\cF'$. We define $\tilde F_{j+1}=F_{j}\cup q_{m-j}\cup p$	; that is $\tilde F_{m-k}=\{p_m,q_m,\dots,p_{k+2},q_{k+2},q_{k+1},p\}$.	
Then $\tilde \cF=\{\tilde F_i\}_{i=1}^{m}$ is a universal flag of $\tilde \cP$ by Proposition~\ref{prop:extnewunifaces}. We denote $\tilde \cP_k=\tilde \cP/\tilde F_{m-k}$ and observe that by Lemma~\ref{lem:quotientsofextShemersewing}, $\tilde \cP_k\simeq \cP_k[\cF'/F_{m-k}]$, where the sewn vertex is $p_{k+1}$, and thus
$\tilde \cP_k\setminus p_{k+1}\simeq \cP_k$. To provide shorter proofs, we will sometimes refer to $p$ as~$p_{m+1}$.
\\

Finally, let $\tilde \cM:=\tilde \cM_m$ be the oriented matroid constructed by Gale Sewing as follows:
\begin{align*}
\tilde \cM_0&:=\Big\{\vv e_1, \dots, \vv e_{r+1}, -\sum_{i=1}^{r+1} \vv e_i\Big\} \text{ with elements labeled }\{\tilde a_1,\dots,\tilde a_{r+1},\tilde p_1\},\\
\tilde q_k&:=
\begin{cases} [\tilde p_k^+,{(\tilde a_{k1})}^{-\ep_{k1}},\dots,{(\tilde a_{kr})}^{-\ep_{kr}}]& \text{if $F_{m-k+1}$ is not split in $\cF'$,}\\
[\tilde p_k^-,{(\tilde a_{k1})}^{\ep_{k1}},\dots, {(\tilde a_{kr})}^{\ep_{kr}}]& \text{if $F_{m-k+1}$ is split in $\cF'$.}
\end{cases}\\
 \tilde p_{k+1}&:=[\tilde q_k^-,p_k^-,{(\tilde a_{k1})}^{-},\dots,{(\tilde a_{k(r-1)})}^{-}],\\
 \tilde\cM_{k}&:=\tilde\cM_{k-1}[\tilde q_k][\tilde p_{k+1}], 
\end{align*}
where if $a_{ij}$ is the element $x$ of $\cM_{i-1}$ then $\tilde a_{ij}$ denotes the element of $\tilde \cM_{i-1}$ labeled as $\tilde x$. See Figure~\ref{fig:ShemerVsMe2} for an example.

\iftoggle{bwprint}{%
\begin{figure}[htpb]
\centering
 \subbottom[$\tilde{\cP}_0\simeq\Gale{\tilde{\cM}_0}$.]{\qquad\includegraphics[width=.25\linewidth]{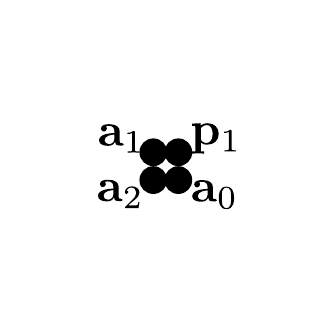}\qquad}\hspace{1cm} 
 \subbottom[$\tilde{\cP}_1=\cP_1\bracket{\vv q_1^+,\vv p_1 ^+, \vv a_1^-}\simeq\Gale{\tilde{\cM}_1}$.]{\qquad\qquad\includegraphics[width=.22\linewidth]{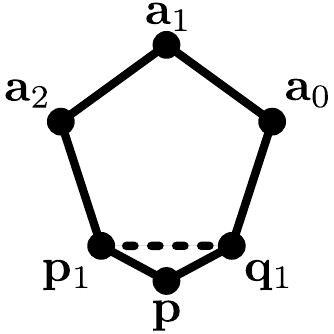}\qquad\qquad}\\ 
 \subbottom[$\tilde{\cM}_0=\cM(\{\tilde{\vv p}_1,\tilde{\vv a}_0,\tilde{\vv a}_1,\tilde{\vv a}_2\})$.]{\quad\qquad\includegraphics[width=.25\linewidth]{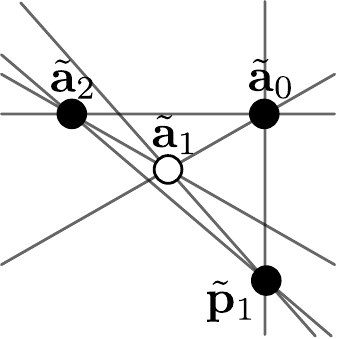}\qquad}\hspace{1cm} 
\subbottom[$\tilde{\cM}_1=\tilde{\cM}_0\bracket{\tilde{\vv p}_1^+,\tilde{\vv a}_2^+,\tilde{\vv a}_1^+}\bracket{\tilde{\vv q}_1^-,\tilde{\vv p}_1^-,\tilde{\vv a}_2^-}$.]{\qquad\qquad\includegraphics[width=.25\linewidth]{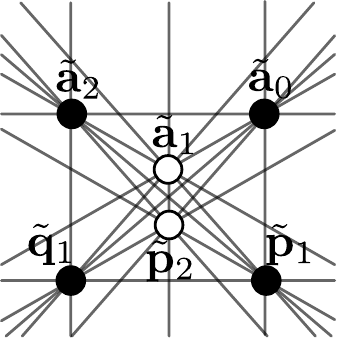}\qquad\qquad} 
 \caption[Primal and dual sewing, part $2$.]{$\tilde{\cP}_1$ is constructed from $\cP_1$ in Figure~\ref{fig:ShemerVsMe1} by sewing $\vv p$ on the flag formed by the universal edge $\{\vv p_1,\vv q_1\}$ (which is not split). Its contraction by $\{\vv p,\vv q_1\}$ is $\tilde{\cP}_0$.
\newline
$\tilde{\cM}_1$ is constructed from $\tilde\cM_0$ by Gale sewing $\tilde{\vv q}_1=\bracket{\tilde{\vv p}_1^+,\tilde{\vv a}_2^+,\tilde{\vv a}_1^+}$ and $\tilde{\vv p}_2=\bracket{\tilde{\vv q}_1^-,\tilde{\vv p}_1^-,\tilde{\vv a}_2^-} $.
\newline
Proposition~\ref{prop:Mk=Pkstar} states that the dual of $\tilde \cP_0$ is isomorphic to $\tilde \cM_0$ and that the dual of $\tilde{\cP}_1$ is isomorphic to $\tilde \cM_1$ with the isomorphism $\tilde{\vv p}_2\mapsto \vv p$, $\tilde{\vv q}_1\mapsto \vv q_1$, $\tilde{\vv p}_1\mapsto \vv p_1$ and $\tilde{\vv a}_i\mapsto \vv a_i$.}
 \label{fig:ShemerVsMe2}
\end{figure}
}{%
\begin{figure}[htpb]
\centering
 \subbottom[$\tilde{\cP}_0\simeq\Gale{\tilde{\cM}_0}$.]{\qquad\includegraphics[width=.25\linewidth]{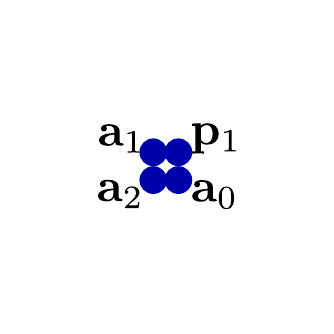}\qquad}\hspace{1cm} 
 \subbottom[$\tilde{\cP}_1=\cP_1\bracket{\vv q_1^+,\vv p_1 ^+, \vv a_1^-}\simeq\Gale{\tilde{\cM}_1}$.]{\qquad\qquad\includegraphics[width=.22\linewidth]{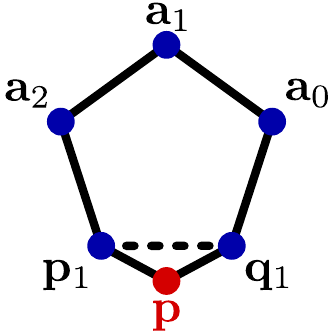}\qquad\qquad}\\ 
 \subbottom[$\tilde{\cM}_0=\cM(\{\tilde{\vv p}_1,\tilde{\vv a}_0,\tilde{\vv a}_1,\tilde{\vv a}_2\})$.]{\quad\qquad\includegraphics[width=.25\linewidth]{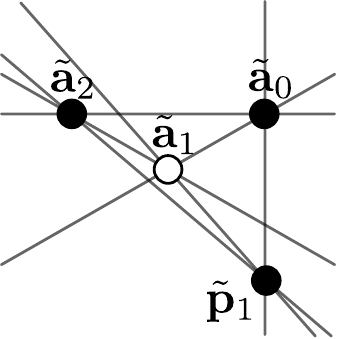}\qquad}\hspace{1cm} 
\subbottom[$\tilde{\cM}_1=\tilde{\cM}_0\bracket{\tilde{\vv p}_1^+,\tilde{\vv a}_2^+,\tilde{\vv a}_1^+}\bracket{\tilde{\vv q}_1^-,\tilde{\vv p}_1^-,\tilde{\vv a}_2^-}$.]{\qquad\qquad\includegraphics[width=.25\linewidth]{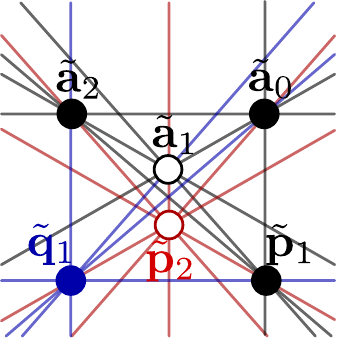}\qquad\qquad} 
 \caption[Primal and dual sewing, part $2$.]{$\tilde{\cP}_1$ is constructed from $\cP_1$ in Figure~\ref{fig:ShemerVsMe1} by sewing $\vv p$ on the flag formed by the universal edge $\{\vv p_1,\vv q_1\}$ (which is not split). Its contraction by $\{\vv p,\vv q_1\}$ is $\tilde{\cP}_0$.
\newline
$\tilde{\cM}_1$ is constructed from $\tilde\cM_0$ by Gale sewing $\tilde{\vv q}_1=\bracket{\tilde{\vv p}_1^+,\tilde{\vv a}_2^+,\tilde{\vv a}_1^+}$ and $\tilde{\vv p}_2=\bracket{\tilde{\vv q}_1^-,\tilde{\vv p}_1^-,\tilde{\vv a}_2^-} $.
\newline
Proposition~\ref{prop:Mk=Pkstar} states that the dual of $\tilde \cP_0$ is isomorphic to $\tilde \cM_0$ and that the dual of $\tilde{\cP}_1$ is isomorphic to $\tilde \cM_1$ with the isomorphism $\tilde{\vv p}_2\mapsto \vv p$, $\tilde{\vv q}_1\mapsto \vv q_1$, $\tilde{\vv p}_1\mapsto \vv p_1$ and $\tilde{\vv a}_i\mapsto \vv a_i$.}
 \label{fig:ShemerVsMe2}
\end{figure}
}

\begin{proposition}\label{prop:Mk=Pkstar}
  In this situation, $\tilde \cM\simeq\Gale{\tilde\cP}$ via the isomorphism $x\mapsto\tilde x$ (cf. Figure~\ref{fig:ShemerVsMe}); hence, $\tilde\cP$ is in $\cG$. 
\end{proposition}

\iftoggle{bwprint}{%
\begin{figure}[htpb]
\begin{center}
\begin{tabular}{m{.18\textwidth}m{.08\textwidth}m{.18\textwidth}m{.48\textwidth}}
\includegraphics[width=.18\textwidth]{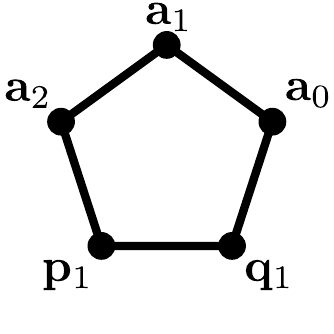}&$ \xleftrightarrow[\cM=\Gale{\cP}]{\text{  duality  }}$&
\includegraphics[width=.18\textwidth]{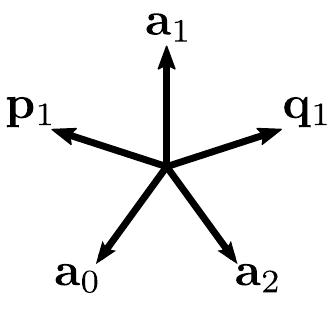}
&
$\cM_0=\{\vv a_0,\vv a_1,\vv a_2\}$\newline $\cM_1=\cM_0[\vv a_2^-,\vv a_1^-][\vv p_1^-,\vv a_2^-]$\newline$\cM=\cM_1$
\\
\centering$\cP$\\
\centering {\rotatebox{90}{$\xleftarrow[{\, \vv p=[\vv q_1^+,\vv p_1 ^+,\vv a_1^-]\,}]{\text{sewing}}$}} &&&
\\
\centering$\tilde \cP=\cP[\vv p]$\\
\includegraphics[width=.18\textwidth]{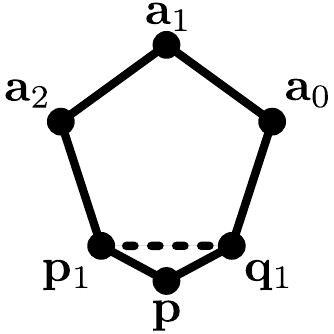}
&$\xleftrightarrow[\Gale{\tilde \cP}\simeq\tilde \cM]{\text{  duality  }}$
&\includegraphics[width=.18\textwidth]{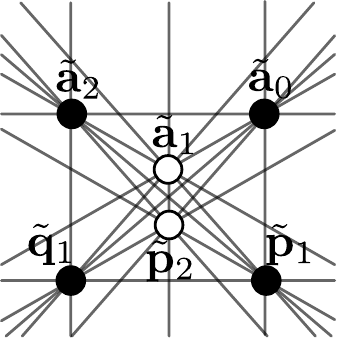}
&$\tilde{\cM}_0=\{\tilde{\vv p}_1,\tilde{\vv a}_0,\tilde{\vv a}_1,\tilde{\vv a}_2\}$\newline 
$\tilde{\cM}_1=\tilde{\cM}_0[\tilde{\vv p}_1^+,\tilde{\vv a}_2^+,\tilde{\vv a}_1^+][\tilde{\vv q}_1^-,\tilde{\vv p}_1^-,\tilde{\vv a}_2^-]$
\newline$\tilde{\cM}=\tilde{\cM}_1$
\end{tabular}
\end{center}
\caption[Primal and dual sewing, comparison.]{We reach the lower left figure by two paths: First (starting in the top right), $\cM$~is the matroid $\cM_1$ constructed by Gale Sewing in Figure~\ref{fig:ShemerVsMe1}. The dual of $\cM$ is $\cP$ (top left). Then $\tilde{\cP}$ (lower left) is constructed from~$\cP$ by sewing $\vv p$ onto the flag formed by the universal edge $\{\vv p_1,\vv q_1\}$ (which is not split). 
\newline
In the second path (lower right), $\tilde{\cM}$~is the matroid $\tilde \cM_1$ constructed by Gale Sewing in Figure~\ref{fig:ShemerVsMe1}; then we dualize to get $\tilde{\cM}^\star$.
\newline
Proposition~\ref{prop:Mk=Pkstar} states that $\Gale{\tilde\cM}\simeq\tilde{\cP}$.
The isomorphism is $\tilde{\vv p}_2\mapsto \vv p$, $\tilde{\vv q}_1\mapsto \vv q_1$, $\tilde{\vv p}_1\mapsto \vv p_1$ and $\tilde{\vv a}_i\mapsto \vv a_i$.}
 \label{fig:ShemerVsMe}
\end{figure}
}{%
\begin{figure}[htpb]
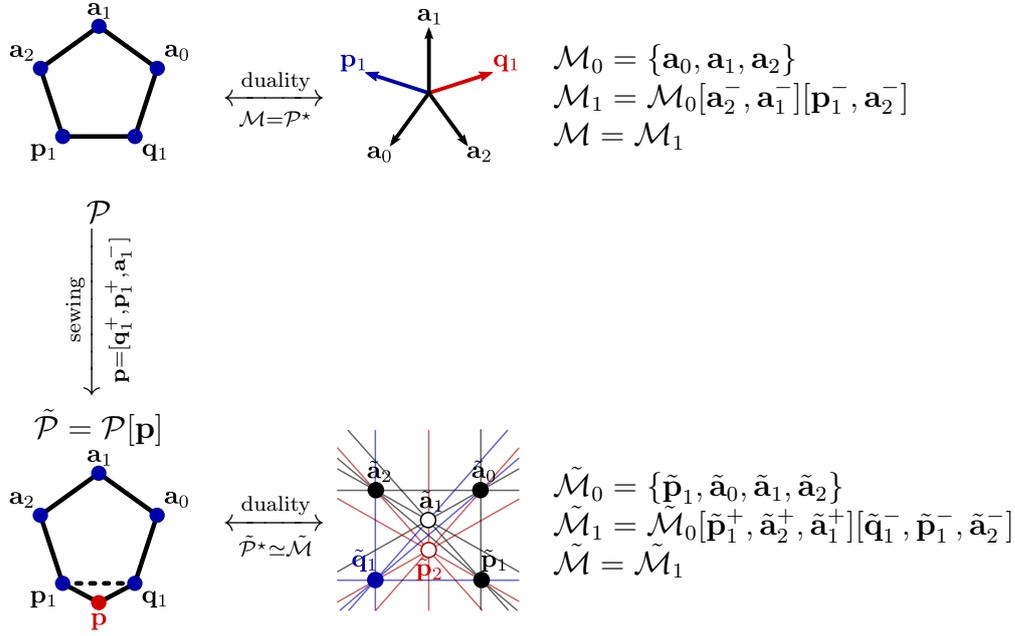

\begin{center}
\begin{tabular}{m{.18\textwidth}m{.08\textwidth}m{.18\textwidth}m{.48\textwidth}}
\includegraphics[width=.18\textwidth]{Figures/SvsM_P1_col}&$ \xleftrightarrow[\cM=\Gale{\cP}]{\text{  duality  }}$&
\includegraphics[width=.18\textwidth]{Figures/SvsM_M1_col}
&
$\cM_0=\{\vv a_0,\vv a_1,\vv a_2\}$\newline $\cM_1=\cM_0[\vv a_2^-,\vv a_1^-][\vv p_1^-,\vv a_2^-]$\newline$\cM=\cM_1$
\\
\centering$\cP$\\
\centering {\rotatebox{90}{$\xleftarrow[{\, \vv p=[\vv q_1^+,\vv p_1 ^+,\vv a_1^-]\,}]{\text{sewing}}$}} &&&
\\
\centering$\tilde \cP=\cP[\vv p]$\\
\includegraphics[width=.18\textwidth]{Figures/SvsM_PP1_col}
&$\xleftrightarrow[\Gale{\tilde \cP}\simeq\tilde \cM]{\text{  duality  }}$
&\includegraphics[width=.18\textwidth]{Figures/SvsM_MM1_col}
&$\tilde{\cM}_0=\{\tilde{\vv p}_1,\tilde{\vv a}_0,\tilde{\vv a}_1,\tilde{\vv a}_2\}$\newline 
$\tilde{\cM}_1=\tilde{\cM}_0[\tilde{\vv p}_1^+,\tilde{\vv a}_2^+,\tilde{\vv a}_1^+][\tilde{\vv q}_1^-,\tilde{\vv p}_1^-,\tilde{\vv a}_2^-]$
\newline$\tilde{\cM}=\tilde{\cM}_1$
\end{tabular}
\end{center}
\caption[Primal and dual sewing, comparison.]{We reach the lower left figure by two paths: First (starting in the top right), $\cM$~is the matroid $\cM_1$ constructed by Gale Sewing in Figure~\ref{fig:ShemerVsMe1}. The dual of $\cM$ is $\cP$ (top left). Then $\tilde{\cP}$ (lower left) is constructed from~$\cP$ by sewing $\vv p$ onto the flag formed by the universal edge $\{\vv p_1,\vv q_1\}$ (which is not split). 
\newline
In the second path (lower right), $\tilde{\cM}$~is the matroid $\tilde \cM_1$ constructed by Gale Sewing in Figure~\ref{fig:ShemerVsMe1}; then we dualize to get $\tilde{\cM}^\star$.
\newline
Proposition~\ref{prop:Mk=Pkstar} states that $\Gale{\tilde\cM}\simeq\tilde{\cP}$.
The isomorphism is $\tilde{\vv p}_2\mapsto \vv p$, $\tilde{\vv q}_1\mapsto \vv q_1$, $\tilde{\vv p}_1\mapsto \vv p_1$ and $\tilde{\vv a}_i\mapsto \vv a_i$.}
 \label{fig:ShemerVsMe}
\end{figure}
}

The proof of Proposition~\ref{prop:Mk=Pkstar} uses the 
following theorem, which states that the contraction and deletion of an element determine the oriented matroid up to the reorientation of that element:

\begin{theorem}[{\cite[Theorem 4.1]{RichterGebertZiegler1994}}]\label{thm:minorsfixmatroid}
Let $\cM'$ and $\cM''$ be two oriented matroids with $n-1$ elements, of respective ranks $\rd$ and $\rd-1$, 
such that $\cov(\cM'')\subseteq\cov(\cM')$. Then there is an oriented matroid $\cM$ with $n$ elements with a special element~$p$
that fulfills $\cM\setminus p =\cM'$ and $\cM/p=\cM''$. The oriented matroid $\cM$ has rank $\rd$ and is unique up to reorientation of $p$.
\end{theorem}

\begin{proof}[Proof of Proposition~\ref{prop:Mk=Pkstar}]
We will prove that $\tilde \cM_k\simeq \Gale{\tilde \cP_k}$, for all $k$. 
The proof will use induction on $k$ and assume that $\tilde \cM_{k-1}\simeq\Gale{\tilde\cP_{k-1}}$. The base case is $\tilde \cM_{0}\simeq\Gale{\tilde\cP_{0}}$, which is straightforward since $\cP_{0}$ is $0$-dimensional.

We will use Theorem~\ref{thm:minorsfixmatroid} twice. Specifically, we will use that if $\cM\setminus p\simeq \cM'\setminus p'$ and $\cM/ p\simeq \cM'/ p'$, then $\cM$ and $\cM'$ are the same oriented matroid up to the reorientation of $p$ and $p'$. 
If additionally there are elements $q\in \cM$ and $q'\in \cM'$ such that $p$ and $q$ are $\alpha$-inseparable in $\cM$ for some $\alpha=\pm 1$ and $p'$ and $q'$ are also $\alpha$-inseparable in $\cM'$, then $\cM\simeq \cM'$.

In particular, we will prove that $\tilde \cM_k/p'_{k+1}\simeq \Gale{\tilde\cP_k}/p_{k+1}$ and that $\tilde \cM_k\setminus p'_{k+1}\simeq \Gale{\tilde\cP_k}\setminus p_{k+1}$. 
Then the claim $\tilde \cM\simeq\Gale{\tilde\cP}$ will follow directly from the fact that $p'_{k+1}$ and $q'_k$ are $(+1)$-inseparable in $\tilde\cM_k$ and that $p_{k+1}$ and $q_k$ are $(+1)$-inseparable in~$\Gale{\tilde\cP_k}$.

\medskip
\noindent \underline{$\tilde \cM_k/\tilde p_{k+1}\simeq \Gale{\tilde\cP_k}/p_{k+1}$}

 Since $\Gale{\tilde\cP_k}/p_{k+1}\simeq\Gale{\cP_k}=\cM_k$ because $\tilde\cP_k\setminus p_{k+1}\simeq\cP_k$, we just need to prove that \begin{equation}\label{eq:quotientsofGaleareShemer}\tilde\cM_k/\tilde p_{k+1}\simeq \cM_k.\end{equation}
By Lemma~\ref{lem:quotientsofGalesewnareGalesewn}, \((\tilde \cM_k/\tilde p_{k+1})\simeq (\tilde \cM_{k-1}/\tilde p_{k})[x'_k][y'_k],\)
where $x'_k=[\tilde a_{k1}^{\ep_{k1}},\dots,\tilde a_{kr}^{\ep_{kr}}]$ and $y'_{k}=[{x'_k}^-,\tilde a_{k1}^{-},\dots,\tilde a_{k(r-1)}^{-}]$.
Using that $\tilde \cM_0/\tilde p_1\simeq \cM_0$ we get the desired result by induction on $k$.

\medskip
\noindent \underline{$\tilde \cM_k\setminus \tilde p_{k+1}\simeq \Gale{\tilde\cP_k}\setminus p_{k+1}$}

The first step is to prove that $(\tilde \cM_k\setminus \tilde p_{k+1})\setminus \tilde q_k \simeq (\Gale{\tilde\cP_{k}}\setminus p_{k+1})\setminus q_k$.
Indeed, using the induction hypothesis,\[(\tilde \cM_k\setminus \tilde p_{k+1})\setminus \tilde q_k =\tilde \cM_{k-1}\simeq\Gale{\tilde\cP_{k-1}}=\pGale{\tilde\cP_{k}/\{q_k,p_{k+1}\}}=(\Gale{\tilde\cP_{k}}\setminus p_{k+1})\setminus q_k.\]

Now we prove that $(\tilde \cM_k\setminus \tilde p_{k+1})/ \tilde q_k \simeq (\Gale{\tilde\cP_{k}}\setminus p_{k+1})/ q_k$. 
First, using Lemma~\ref{lem:contractdeletele} and \eqref{eq:quotientsofGaleareShemer}, we see that
 \begin{displaymath}(\tilde \cM_k\setminus \tilde p_{k+1})/\tilde q_k\simeq(\tilde \cM_k / \tilde p_{k+1})\setminus \tilde q_k\simeq \cM_k\setminus q_k.\end{displaymath}
Now, using again Lemma~\ref{lem:contractdeletele} and that $\tilde \cP_k\setminus p_{k+1}\simeq \cP_k$, we see that
\begin{displaymath} (\Gale{\tilde\cP_k}\setminus p_{k+1})/q_k=\pGale{\tilde\cP_k/p_{k+1}\setminus q_k}\simeq\pGale{\tilde \cP_k\setminus p_{k+1}/ q_k}=\pGale{\cP_k/ q_k}=\Gale{\cP_k}\setminus q_k.\end{displaymath}
Our claim follows since $\cM_k=\Gale{\cP_k}$ by definition.

\smallskip
Because of Theorem~\ref{thm:minorsfixmatroid}, so far we have seen that $\tilde \cM_k\setminus \tilde p_{k+1}\simeq \Gale{\tilde\cP_k}\setminus p_{k+1}$ up to reorientation of $q_k$. We will conclude the proof by seeing that $q_k$ and $\tilde q_k$ have the same orientation in $\cM_k$ and in $\tilde \cM_k$ respectively.

If $F_{m-k+1}$ is not split, then $\tilde q_k$ is $(-1)$-inseparable with $\tilde p_k$ in $\tilde\cM_k$ by construction. Moreover, $q_k$ is $(-1)$-inseparable with $p_k$ in $\pGale{\tilde\cP_{k}/ p_{k+1}}$ since they are $(+1)$-inseparable in the primal: by Proposition~\ref{prop:allquotientsofle} \[\tilde\cP_{k}/ p_{k+1}\simeq\big(\cP_{k}\underbrace{[\cF'/F_{m-k}]}_{p_{k+1}}\big)/p_{k+1}\simeq \big(\cP_{k}/ p_{k}\big)\underbrace{[q_k^-,\dots]}_{p_k},\]
where the last isomorphism sends $p_k$ to the sewn vertex, which is $(+1)$-inseparable from $q_k$.

If $F_{m-k+1}$ is split, then $\tilde q_k$ is $(+1)$-inseparable with $\tilde p_k$ in $\tilde\cM_k$. Moreover, $q_k$ is $(+1)$-inseparable with $p_k$ in $\pGale{\tilde\cP_{k}/ p_{k+1}}$ because $\tilde\cP_{k}/ p_{k+1}\simeq \cP_{k}/ p_{k}[q_k^+,\dots]$. 
\end{proof}

There is a missing detail to conclude that $\cE\subseteq \cG$ from Proposition~\ref{prop:Mk=Pkstar}:
we have to check that after sewing, the universal flags of Proposition~\ref{prop:extnewunifaces} are also of the form $F_j=\cup_{i=0}^{j-1}(p_{m-i},q_{m-i})$ for some Gale Sewing order.
This is a consequence of Lemma~\ref{lem:galesewingorder}, which allows to change the order of the sewings in $\tilde\cM$. This concludes the proof of Theorem~\ref{thm:cOsubsetcG}. \qed

\begin{remark}\label{rmk:doesnotgeneralize}
The fact that $\cE\subsetneq\cG$ implies that in some sense Gale Sewing generalizes ordinary sewing. However, it is not true that the Sewing Theorem~\ref{thm:shemersewing} is a consequence of the Double Extension Theorem~\ref{thm:thethm}, because there are neighborly matroids that have universal flags but are not in~$\cG$. Hence one can sew on them but they 
cannot be treated with Proposition~\ref{prop:Mk=Pkstar}. 
This will become clear in Section~\ref{sec:nonrealizable}, where we work with ``$\cM^{10}_{425}$'', a non-realizable neighborly matroid that has universal flags. Since Gale Sewing (Construction~\ref{constr:cG}) only builds realizable matroids, this matroid is not in $\cG$ and yet one can sew on it. This shows why both constructions are needed.
\end{remark}

\section{Non-realizable neighborly oriented matroids}\label{sec:nonrealizable}
Since the only neighborly matroids of rank~$3$ are cyclic polytopes, there are no non-realizable neighborly matroids of rank~$3$. The sphere ``$\cM^{10}_{425}$'' from Altshuler's list~\cite{Altshuler1977} corresponds to a neighborly matroid of rank~$5$ with $10$ elements. 
In~\cite{BokowskiGarms1987}, this matroid is shown to be non-realizable, thus proving that non-realizable neighborly matroids exist. Kortenkamp's construction~\cite{Kortenkamp1997} can be used to build non-realizable neighborly matroids of corank $3$. 
We combine Theorems~\ref{thm:extshemersewing} and~\ref{thm:thethm} to show that there are many non-realizable neighborly matroids. A lower bound for the cardinality of the number of non-realizable neighborly matroids is derived in Theorem~\ref{thm:nonrealizablebound}.

\begin{table}[htpb]
\caption{List of facets of ``$\cM^{10}_{425}$'' as they appear in~\cite{Altshuler1977}.}\label{tb:facetsM425}
\autorows{c}{5}{c}{
$\{a_0a_1a_2a_3\}$,
$\{a_0a_1a_6a_7\}$,
$\{a_0a_4a_5a_6\}$,
$\{a_1a_2a_6a_9\}$,
$\{a_2a_3a_4a_5\}$,
$\{a_2a_4a_5a_8\}$,
$\{a_4a_5a_6a_7\}$,
$\{a_0a_1a_2a_6\}$,
$\{a_0a_1a_7a_9\}$,
$\{a_0a_4a_6a_7\}$,
$\{a_1a_3a_4a_5\}$,
$\{a_2a_3a_4a_8\}$,
$\{a_2a_5a_7a_8\}$,
$\{a_4a_5a_7a_9\}$,
$\{a_0a_1a_3a_5\}$,
$\{a_0a_1a_8a_9\}$,
$\{a_0a_4a_7a_9\}$,
$\{a_1a_3a_4a_8\}$,
$\{a_2a_3a_5a_7\}$,
$\{a_2a_6a_7a_8\}$,
$\{a_4a_5a_8a_9\}$,
$\{a_0a_1a_4a_5\}$,
$\{a_0a_2a_3a_6\}$,
$\{a_0a_4a_8a_9\}$,
$\{a_1a_3a_8a_9\}$,
$\{a_2a_3a_6a_7\}$,
$\{a_2a_6a_8a_9\}$,
$\{a_5a_7a_8a_9\}$,
$\{a_0a_1a_4a_8\}$,
$\{a_0a_3a_5a_6\}$,
$\{a_1a_2a_3a_9\}$,
$\{a_1a_6a_7a_9\}$,
$\{a_2a_3a_8a_9\}$,
$\{a_3a_5a_6a_7\}$,
$\{a_6a_7a_8a_9\}$
}
\end{table}

\begin{theorem}\label{thm:nonrealizable}
 There exists a non-realizable neighborly matroid of rank $\rd$ with $n$ elements for every $\rd\geq 5$ and $n\geq \rd+5$.
\end{theorem}
\begin{proof}

We start with ``$\cM^{10}_{425}$''. The pairs $\{a_0,a_1\}$, $\{a_2,a_3\}$, $\{a_4,a_5\}$, $\{a_6,a_7\}$ and $\{a_8,a_9\}$ are universal edges of $\cM^{10}_{425}$ because the corresponding contractions are polygons with $8$ vertices (see Table~\ref{tb:facetsM425}). 
In particular, $\{a_0,a_1\}\subset \{a_0,a_1,a_2,a_3\}$ is a universal flag. Hence, applying Theorem~\ref{thm:extshemersewing} we get many non-realizable matroids of rank~$5$ with $n$ vertices for any $n\geq 10$.

Now, applying Corollary~\ref{cor:primalGaleSewing} to these matroids, we get non-realizable oriented matroids of rank $5+2k$ and $n$ vertices for any $k\geq 0$ and any $n\geq 10+2k$.

To get non-realizable matroids of even rank, just observe that any single element extension on the dual of a neighborly matroid of rank $2k+1$ yields the dual of a neighborly matroid of rank $2k+2$.
\end{proof}

\begin{remark}\label{rmk:nonrealizable} 
The result of the previous theorem can be slightly improved. Indeed, in~\cite{RichterGebertSturmfels1991} it is shown that there is a non-realizable cyclic matroid polytope of rank~$4$ with $10$ vertices (cf. \cite[Proposition 9.4.5]{OrientedMatroids1993}). Using Extended Sewing, this implies that there are non-realizable neighborly oriented matroids of rank~$4$ with $n$ elements for any $n\geq 10$.
\end{remark}

All neighborly matroids of rank $2m+1$ that have $n\leq2m+3$ vertices are cyclic polytopes. Moreover, all oriented matroids of rank $5$ with $8$ elements are realizable~\cite[Corollary 8.3.3]{OrientedMatroids1993}. Hence the first case (of odd rank) that Theorem~\ref{thm:nonrealizable} does not deal with are neighborly matroids of rank $5$ with $9$ elements (cf. Figure~\ref{fig:nrn}).

\begin{question}
 Do there exist non-realizable neighborly oriented matroids of rank $5$ on $9$ elements?
\end{question}

\iftoggle{bwprint}{%
\begin{figure}[htpb]
\begin{center}
\includegraphics[width=.8\textwidth]{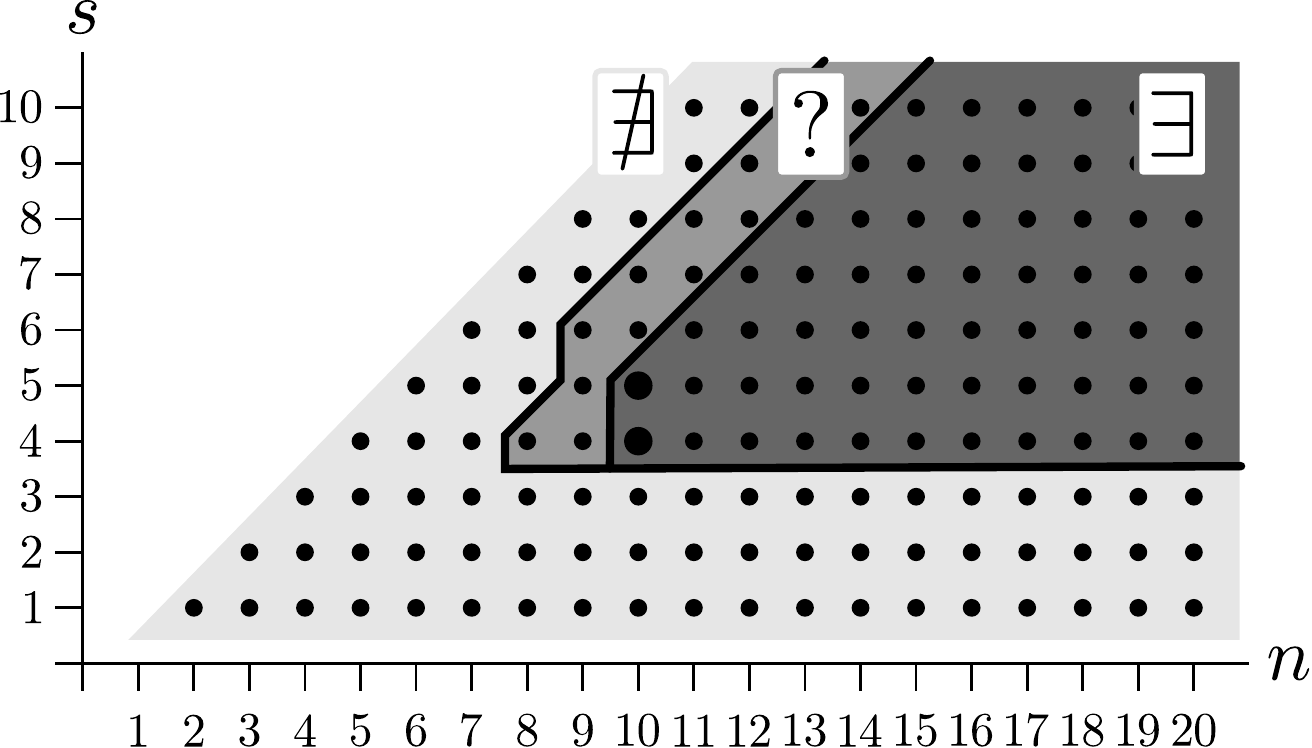}
\end{center}
 \caption[Existence of non-realizable neighborly oriented matroids.]{This table summarizes our knowledge on the existence of non-realizable neighborly oriented matroids of rank $s$ with $n$ elements. It uses~\cite[Corollary 8.3.3]{OrientedMatroids1993}, Theorem~\ref{thm:nonrealizable} and Remark~\ref{rmk:nonrealizable}. Bolder points represent concrete instances from which we sew.} 
 \label{fig:nrn}
\end{figure}
}{%
\begin{figure}[htpb]
\begin{center}
\includegraphics[width=.8\textwidth]{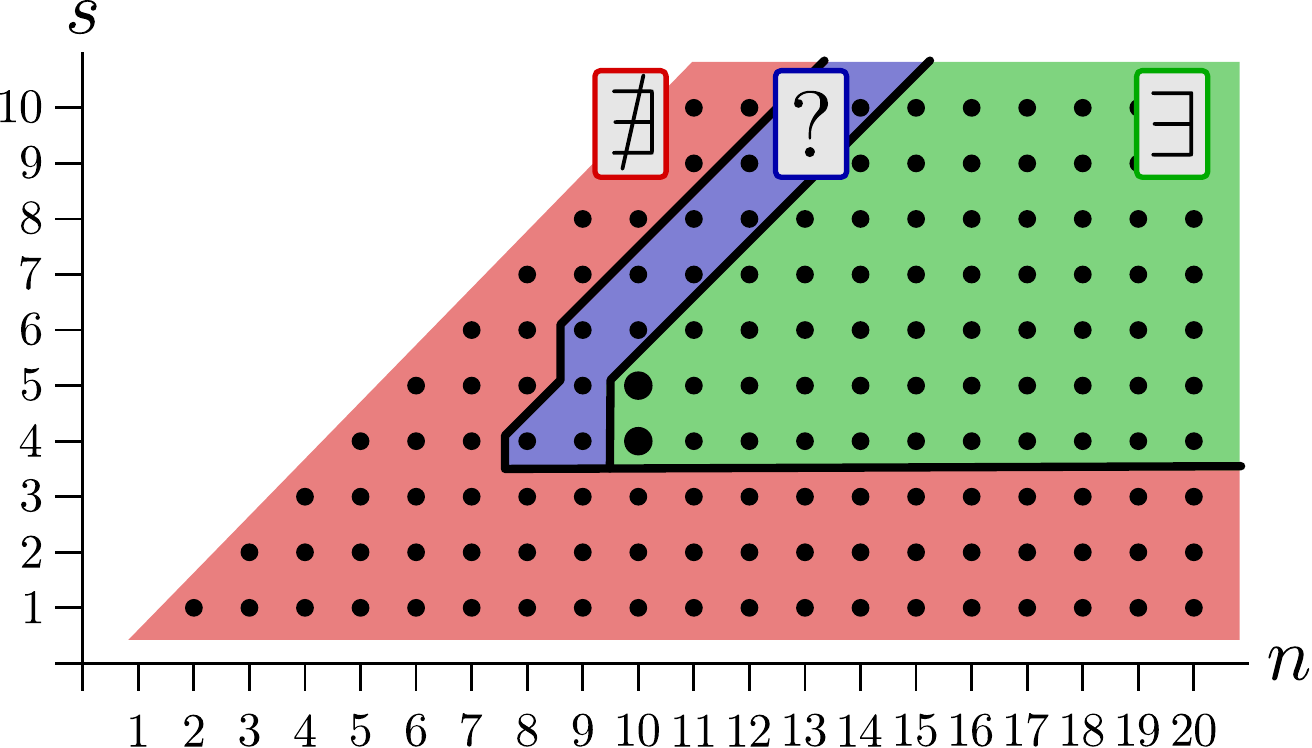}
\end{center}
 \caption[Existence of non-realizable neighborly oriented matroids.]{This table summarizes our knowledge on the existence of non-realizable neighborly oriented matroids of rank $s$ with $n$ elements. It uses~\cite[Corollary 8.3.3]{OrientedMatroids1993}, Theorem~\ref{thm:nonrealizable} and Remark~\ref{rmk:nonrealizable}. Bolder points represent concrete instances from which we sew.} 
 \label{fig:nrn}
\end{figure}
}

\chapter{Many neighborly polytopes and oriented matroids}\label{ch:counting}

In this chapter we approach the question of how many different combinatorial types of neighborly polytopes (and oriented matroids) there are, and in particular, how many of them are in $\cG$. We start with exact numbers for some particular cases of low rank and corank. In Section~\ref{sec:counting} we derive lower bounds for the number of neighborly polytopes by estimating the number of polytopes in $\cG$. Thanks to Theorem~\ref{thm:nonrealizable}, we will be able to use the same techniques to provide lower bounds for the number of non-realizable neighborly oriented matroids.

\section{Some exact numbers}\label{sec:exact}

We have worked with four families of neighborly polytopes:

\begin{description}
 \item[$\cN$:] All neighborly polytopes.
 \item[$\cS$:] Totally sewn neighborly polytopes (Sewing, Construction~\ref{constr:cS}).
 \item[$\cE$:] Neighborly polytopes constructed by Extended Sewing (Construction~\ref{constr:cE}).
 \item[$\cO$:] Neighborly polytopes built by Extended Sewing and Omitting (Construction~\ref{constr:cO}).
 \item[$\cG$:] Gale sewn neighborly polytopes (Gale Sewing, Construction~\ref{constr:cG}).
\end{description}

Table~\ref{tb:numcombtypes} contains the exact number of combinatorial types of neighborly $d$-polytopes with $n$ vertices in each of these families for the cases $d=4$ and $n=8,9$ and for $d=6$ and $n=10$. Exact numbers for $\cN$ come from~\cite{AltshulerSteinberg1973} and~\cite{BokowskiShemer1987}, exact numbers for $\cS$ and $\cO$ come from~\cite{Shemer1982}. Numbers for $\cG$ and $\cE$ have been computed with the help of \texttt{polymake}~\cite{polymake}.

\begin{table}[htpb]
\centering
  \caption{Exact number of combinatorial types}	\label{tb:numcombtypes}
  \setlength{\tabcolsep}{15pt}
   \renewcommand{\arraystretch}{2}
  \begin{tabular}{ c  c c c  c  c  c  c }
    \hline
      $d$& $n$ && $\cS$ & $\cE$ & $\cO$ & $\cG$ & $\cN$\\
    \hline 
      4 & 8 && 3 & 3 & 3 & 3 & 3 \\ 
      4 & 9 && 18 & 18 & 18 & 18 & 23 \\ 
      6 & 10 && 15 & 26 & 28 & 28 & 37 \\ 
    \hline
  \end{tabular}
\end{table}

The known relationships between these families are summarized in the following proposition. The inclusion $\cO\subseteq\cG$ is a consequence of Theorem~\ref{thm:cOsubsetcG}, while the remaining inclusions are by definition. Then Table~\ref{tb:numcombtypes} discards most of the equalities.
\begin{proposition}
\(\cS\subsetneq\cE\subsetneq\cO\subseteq\cG\subsetneq\cN.\) \qed
\end{proposition}

This begs the question:
\begin{question}
Is $\cO=\cG$?
\end{question}

\begin{remark}
 In Table~\ref{tb:numcombtypes} one sees that there are precisely two neighborly $6$-polytopes with $10$ vertices that can be constructed with Gale Sewing but not with Extended Sewing. Labeling their vertices as $\{0,1,\dots,9\}$, we display their facets (as subsets of vertices) in Tables~\ref{tb:notinE1} and~\ref{tb:notinE2}.
\end{remark}

\begin{table}[htpb]
\centering
\caption{List of facets of a neighborly $6$-polytope with $10$ vertices that belongs to $\cG$ and not to $\cE$.}\label{tb:notinE1}

\autorows{c}{6}{c}{
\{$123458\}$,
$\{234589\}$,
$\{123489\}$,
$\{234679\}$,
$\{234569\}$,
$\{345689\}$,
$\{456789\}$,
$\{124589\}$,
$\{124567\}$,
$\{124679\}$,
$\{124569\}$,
$\{145679\}$,
$\{126789\}$,
$\{125689\}$,
$\{156789\}$,
$\{235689\}$,
$\{236789\}$,
$\{135678\}$,
$\{123568\}$,
$\{123678\}$,
$\{012789\}$,
$\{023789\}$,
$\{012378\}$,
$\{013567\}$,
$\{012356\}$,
$\{012367\}$,
$\{012389\}$,
$\{013578\}$,
$\{035678\}$,
$\{036789\}$,
$\{015789\}$,
$\{012567\}$,
$\{014579\}$,
$\{014589\}$,
$\{045789\}$,
$\{046789\}$,
$\{034689\}$,
$\{034679\}$,
$\{045678\}$,
$\{034568\}$,
$\{024567\}$,
$\{023456\}$,
$\{023467\}$,
$\{013489\}$,
$\{013458\}$,
$\{012479\}$,
$\{023479\}$,
$\{012349\}$,
$\{012457\}$,
$\{012345\}$
} 
\end{table}

\begin{table}[htpb]
\centering
\caption{List of facets of a neighborly $6$-polytope with $10$ vertices that belongs to $\cG$ and not to $\cE$.}\label{tb:notinE2}
\autorows{c}{6}{c}{
$\{123467\}$,
$\{123456\}$,
$\{234568\}$,
$\{234589\}$,
$\{123459\}$,
$\{123489\}$,
$\{134567\}$,
$\{345689\}$,
$\{456789\}$,
$\{145678\}$,
$\{124589\}$,
$\{124568\}$,
$\{125679\}$,
$\{126789\}$,
$\{125689\}$,
$\{156789\}$,
$\{235689\}$,
$\{236789\}$,
$\{235679\}$,
$\{123567\}$,
$\{012579\}$,
$\{012789\}$,
$\{023789\}$,
$\{023579\}$,
$\{012357\}$,
$\{023678\}$,
$\{012359\}$,
$\{012389\}$,
$\{035679\}$,
$\{036789\}$,
$\{015789\}$,
$\{012678\}$,
$\{012468\}$,
$\{014678\}$,
$\{012467\}$,
$\{014589\}$,
$\{014578\}$,
$\{045789\}$,
$\{046789\}$,
$\{034689\}$,
$\{045679\}$,
$\{034569\}$,
$\{034567\}$,
$\{013489\}$,
$\{013457\}$,
$\{013459\}$,
$\{012348\}$,
$\{012347\}$,
$\{023468\}$,
$\{023467\}$}
\end{table}

\section{Many neighborly oriented matroids}\label{sec:counting}

The aim of this section is to find lower bounds for \defn{$\lnei{n}{d}$}\index{$\lnei{n}{d}$}, the number of combinatorial types of vertex-labeled neighborly polytopes with $n$ vertices in dimension~$d$. 
Since two neighborly polytopes with the same combinatorial type have the same oriented matroid (Theorem~\ref{thm:neigharerigid}), it suffices to bound the number of labeled realizable neighborly matroids. 

Our strategy will consist in using the Gale Sewing technique of Theorem~\ref{thm:thethm} to construct many neighborly polytopes in~$\cG$ for which we can certify that their oriented matroids are all different.
\\

In the remainder of this chapter, all polytopes and oriented matroids will be labeled. Nevertheless, our bounds will be so large 
as to be almost asymptotically indistinguishable from the naive bounds for unlabeled combinatorial types obtained by dividing by $n!$.

\subsection{Many lexicographic extensions}
 
A first step is to compute lower bounds for \defn{$\lle{n}{r}$}\index{$\lle{n}{r}$}, the smallest number of different labeled lexicographic extensions that a balanced matroid of rank~$r$ with $n$ elements can have. Here, a labeled lexicographic extension of~$\cM$ is a lexicographic extension $\cM[p]$ labeled in such a way that the labels of the elements of~$\cM$ are preserved.

There are $2^{r}\!\frac{n!}{(n-r)!}$ different expressions for lexicographic extensions of a rank~$r$ oriented matroid on $n$ elements, yet not all of them represent different labeled oriented matroids (for example, the extension from Figure~\ref{fig:le} can be represented by both $[\vv x_4^+,\vv x_1^-,\vv x_6^+]$ and $[\vv x_4^+,\vv x_6^+,\vv x_5^+]$ among others). We aim to avoid counting the same extension twice with two different expressions. 

We present two different bounds for the number of lexicographic extensions of a uniform matroid. The bound in Proposition~\ref{prop:lble} is smaller than that in Theorem~\ref{thm:lble2}, but its proof is simpler. Moreover, the factor that we lose in the first bound is asymptotically much smaller than its value. Since all posterior calculations are also simpler with the smaller bound, it is the one that we use. However, we present both for the sake of completeness.

\begin{proposition}\label{prop:lble}
 Let $\cM$ be a rank $r>1$ labeled uniform balanced matroid with $n$ elements. If $n-r-1\ge 2$ is even, then there are at least \begin{equation}\label{eq:bndle1}\lle{n}{r}\geq\frac{2n!}{(n-r+1)!} 
\end{equation}
different uniform labeled lexicographic extensions of $\cM$.
\end{proposition}
\begin{proof}
We focus only on those extensions where $\ep_i=+$ for all $i$, and show that they are unambiguous except for the last element $a_r^{\ep_r}$.

For this, observe that if $r>1$ and the lexicographic extensions by $[a_1^{+},\dots,a_r^{+}]$ and $[{a'_1}^{+},\dots,{a'_r}^{+}]$ yield the same oriented matroid, then either $a_1=a_1'$, or $a_1$ and $a'_1$ are $(-1)$-inseparable. Indeed, for any cocircuit $C\in\co(\cM)$ with $C(a_1)\neq 0$ and $C(a_1')\neq 0$, the signature $\sigma:\co(\cM)\rightarrow\{\pm,0\}$ of the lexicographic extension fulfills $\sigma(C)=C(a_1)$ and $\sigma(C)=C(a_1')$. Thus, if $a_1\neq a_1'$ then $a_1$ and $a_1'$ are $(+1)$-inseparable in $\Gale\cM$ and hence, by Lemma~\ref{lem:insep}, $(-1)$-inseparable in~$\cM$.

But balanced matroids of rank $r\geq 2$ and even corank $\geq 2$ only have $(+1)$-inseparable pairs (see Lemma~\ref{lem:balonlycovar}), which proves that $a_1=a_1'$. Analogously, if $a_i$ and $a_i'$ are the first distinct elements and $i<r$, we can apply the previous argument on the contraction by $\{a_1,\dots,a_{i-1}\}$.

Hence, there are at least $\frac{n!}{(n-r+1)!}$ different choices for the first $r-1$ elements (which give rise to different matroids). For the last element, observe that $\cM/\{a_1,\dots,a_{r-1}\}$ is a matroid of rank~$1$, and that there are exactly two possible different extensions for a matroid of rank~$1$.
\end{proof}

For the stronger bound, we will need the following result by Cordovil and Duchet on the inseparability graph of oriented matroids.\index{inseparability graph}

\begin{theorem}[{\cite[Theorem 1.1]{CordovilDuchet1990}}]\label{thm:IG}
 Let $\cM$ be a rank $r$ uniform oriented matroid with $n$ elements.
\begin{itemize}
 \item If $r\leq 1$ or $r\geq n-1$, then $\IG ( \cM )$ is the complete graph $K_n$.
 \item If $r=2$ or $r=n-2$, then $\IG(\cM)$ is an $n$-cycle.
 \item If $2<r<n-2$, then $\IG(\cM)$ is either a $n$-cycle, or a disjoint union of chains.
\end{itemize}
\end{theorem}

We are ready to prove the strong bound.
\begin{theorem}\label{thm:lble2}
 Let $\cM$ be a rank $r>1$ labeled uniform oriented matroid with $n>\rr+1$ elements. The number of different labeled uniform lexicographic extensions of $\cM$ fulfills \begin{equation}\label{eq:bndle2}\lle{n}{\rr}\geq2^{\rr-1}\frac{n!}{(n-1)(n-\rr)!}.\end{equation}
\end{theorem}
\begin{proof}

To prove this result, we define a family of expressions of lexicographic extensions that guarantees that no two of them define the same oriented matroid, and then give a lower bound on its cardinality. 
The key tool is Observation~\ref{obs:inlble2proof0} below, that is used to show that if $[a_1^{\ep_1},\dots]$ and $[{a_1'}^{\ep_1'},\dots]$ give the same lexicographic extension, then $a_1$ and $a_1'$ must be $(- \ep_1\ep_1')$-inseparable. Therefore, we can use our knowledge on the inseparable elements of~$\cM$ to avoid those expressions that could give repeated lexicographic extensions.\\

Consider the inseparability graph $\IG(\cM)$ of $\cM$, and orient it in such a way that every node has at most one outgoing edge (this is possible because of Theorem~\ref{thm:IG}). Do the same for the inseparability graphs of all contractions of $\cM$ of rank $\geq 2$, $\IG(\cM/S)$.
If $a\rightarrow b$ is a directed edge of $\IG(\cM)$ then we say that $\suc(a)=b$; and we say that $\alpha(a)=\alpha$ when $a$ and~$b$ are $\alpha$-inseparable in $\cM$. In particular, if $a$ and $b$ are inseparable, either $a=\suc(b)$ or $b=\suc(a)$.

Fixed a sequence $a_1,\dots,a_\rr$ of elements of $\cM$ (which will be the elements defining the lexicographic extension), for each $i\leq \rr-2$, let $S_i:=\{a_1,\dots,a_i\}$. We define $\suc_i(a)$ and $\alpha_i(a)$ as $\suc(a)$ and $\alpha(a)$ in the orientation of the inseparability graph of $\cM/S_i$. An example is depicted in Figure~\ref{fig:OrientedIG}.\\

\iftoggle{bwprint}{%
\begin{figure}[htpb]
\centering
 \subbottom[$\cM$]{\includegraphics[width=.4\linewidth]{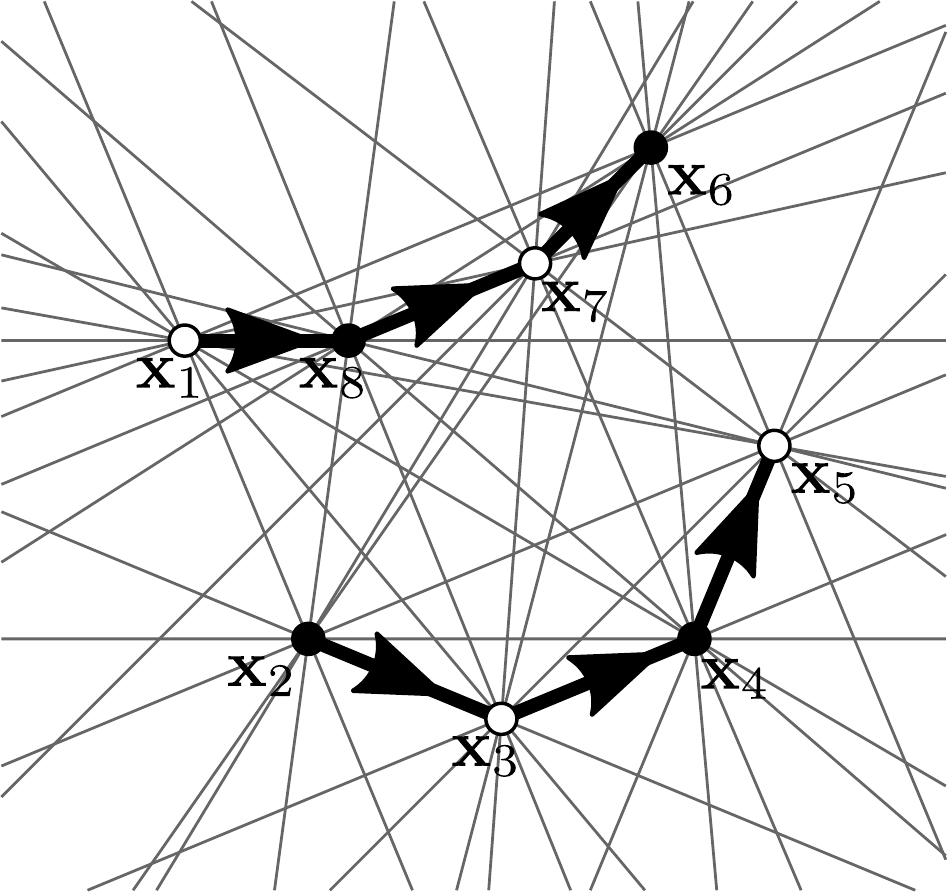}}\hspace{3cm}
 \subbottom[$\cM/\vv x_1$.]{\includegraphics[width=.26\linewidth]{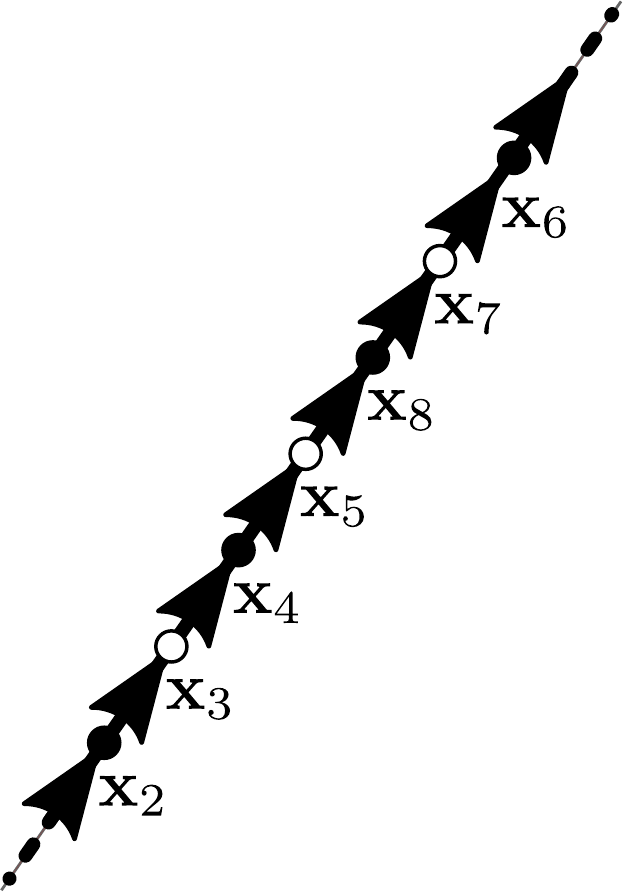}}
 \caption[Orientations of the inseparability graphs.]{An orientation of the inseparability graph of $\cM$ and an orientation of the inseparability graph of $\cM/\vv x_1$. With these orientations, $\suc(\vv x_1)=\vv x_8$ and $\alpha(\vv x_1)=(+1)$; and if $a_1=\vv x_1$ then $\suc_1(\vv x_5)=\vv x_8$ and $\alpha_1(\vv x_5)=(+1)$.}
 \label{fig:OrientedIG}
\end{figure}
}{%
\begin{figure}[htpb]
\centering
 \subbottom[$\cM$]{\includegraphics[width=.4\linewidth]{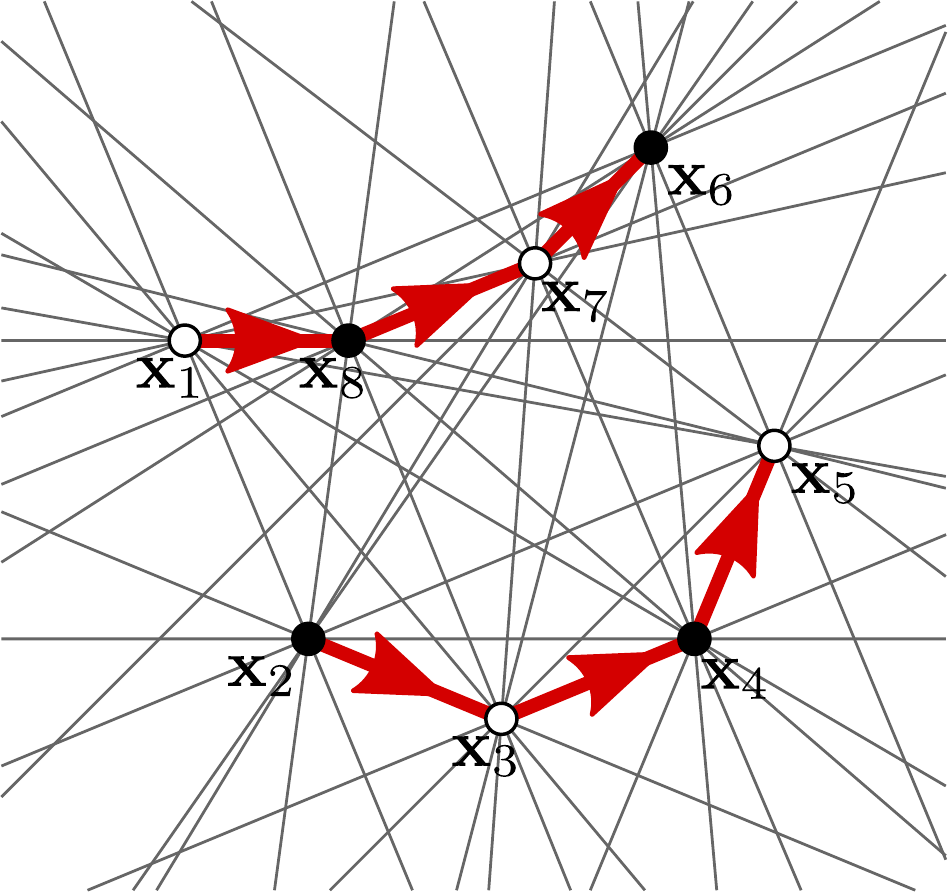}}\hspace{3cm}
 \subbottom[$\cM/\vv x_1$.]{\includegraphics[width=.26\linewidth]{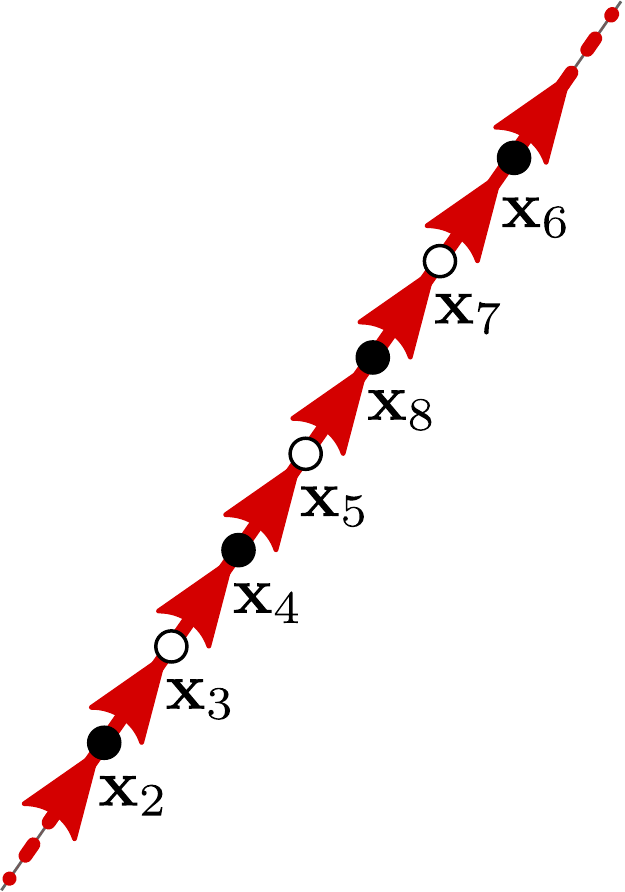}}
 \caption[Orientations of the inseparability graphs.]{An orientation of the inseparability graph of $\cM$ and an orientation of the inseparability graph of $\cM/\vv x_1$. With these orientations, $\suc(\vv x_1)=\vv x_8$ and $\alpha(\vv x_1)=(+1)$; and if $a_1=\vv x_1$ then $\suc_1(\vv x_5)=\vv x_8$ and $\alpha_1(\vv x_5)=(+1)$.}
 \label{fig:OrientedIG}
\end{figure}
}

For these fixed orientations of the inseparability graphs, we consider only the expressions of lexicographic extensions by $p=[a_1^{\ep_1},\dots, a_r^{\ep_r}]$ that fulfill
\begin{enumerate}[(i)]
 \item\label{it:condle1} For $1<i<r$, $a_i^{\ep_i}$ is not $\suc_{i-2}(a_{i-1})^{-\alpha_{i-2}(a_{i-1})\ep_{i-1}}$.
 \item\label{it:condle2} For $1<i<r$, $a_i^{\ep_i}$ is not $\suc_{i-1}(\suc_{i-2}(a_{i-1}))^{\alpha_{i-1}(\suc_{i-2}(a_{i-1}))\alpha_{i-2}(a_{i-1})\ep_{i-1} }$.
 \item\label{it:condle3} $a_r=\suc_{r-2}(a_{r-1})$ and $\ep_r=\alpha_{r-2}(a_{r-1})\ep_{r-1}$.
\end{enumerate}

We prove that if the expression of a lexicographic extension fulfills conditions \eqref{it:condle1}, \eqref{it:condle2} and \eqref{it:condle3}, then $p$ and $\suc_{i-1}(a_i)$ are not $(\alpha_{i-1}(a_i)\ep_i)$-inseparable in $\cM[p]/\{a_1,\dots,a_{i}\}$ for all $i<r$. This condition (together with \eqref{it:condle3}) turns out to be sufficient to certify that the lexicographic extensions are different. \\

Our first step is to prove this claim. That is, there cannot be two different lexicographic extensions fulfilling both that $p$ and $\suc_{i-1}(a_i)$ are not $(\alpha_{i-1}(a_i)\ep_i)$-inseparable in $\cM[p]/\{a_1,\dots,a_{i}\}$ for all $i<r$ and that give rise to the same oriented matroid.

Assume that $\cM[p]$ can be written both as $\cM[a_1^{\ep_1},\dots,a_{i-1}^{\ep_{i-1}},a_i^{\ep_i},\dots]$ and $\cM[a_1^{\ep_1},\dots,a_{i-1}^{\ep_{i-1}},{a_i'}^{\ep_i'},\dots]$ with $a_i\neq a_i'$ for some $i<r$. We will see that then either $p$ and $\suc_{i-1}(a_i)$ are $(\alpha_{i-1}(a_i)\ep_i)$-inseparable in $\cM[p]/\{a_1,\dots,a_{i}\}$ or $p$ and $\suc_{i-1}(a_i')$ are $(\alpha_{i-1}(a_i')\ep_i)$-inseparable in $\cM[p]/\{a_1,\dots,a_{i}\}$.
First, we need one observation.

\begin{observation}\label{obs:inlble2proof0}
If $e$ is $\alpha$-inseparable with $p$ in $\cM[p]/\{a_1,\dots,a_{i-1}\}$, where $p=[a_1^{\ep_1},\dots,a_r^{\ep_r}]$, then either $e=a_i$ or $e$ and $a_i$ are $(\ep_i\alpha)$-inseparable in $\cM/\{a_1,\dots,a_{i-1}\}$.
\end{observation}
\begin{proof}[Proof of the observation]
Any cocircuit $C$ of $\cM/\{a_1,\dots,a_{i-1}\}$ such that $C(a_i)\neq 0$ and $C(e)\neq 0$,  extends to a cocircuit $C'$ of $\cM[p]/\{a_1,\dots,a_{i-1}\}$ with $C'(p)=\ep_iC(a_i)$ and $ C'(p)=-\alpha C(e)$ (by definition of lexicographic extension). Hence $C(a_i)=-\ep_i\alpha C(e)$ for every cocircuit $C$ such that $\{a_1,e\}\subseteq \ul C$.
\end{proof}

As a consequence, since $p$ is $(-\ep_i')$-inseparable with $a_i'$ in $\cM/S_{i-1}$ by Lemma~\ref{lem:leinseparable}, then $a_i$ and $a_i'$ must be $(- \ep_i\ep_i')$-inseparable in $\cM/S_{i-1}$. Hence we can assume without loss of generality that $a_i'=\suc_{i-1}(a_i)$ and $\ep_i'=-\alpha_{i-1}(a_i)\ep_i$. 
But then, our observation of the inseparability of $p$ and $a_i'$ shows that $p$ and $\suc_{i-1}(a_i)=a_i'$ are $(\alpha_{i-1}(a_i)\ep_i=-\ep_i')$-inseparable in $\cM/\{a_1,\dots,a_{i-1}\}$. Therefore, they must also be $(\alpha_{i-1}(a_i)\ep_i)$-inseparable in the contraction $\cM/\{a_1,\dots,a_{i}\}$, which proves our claim.

Of course, it still could happen that there are two different expressions of the form $\cM[a_1^{\ep_1},\dots,a_{r-1}^{\ep_{r-1}},a_r^{\ep_r}]$ and $\cM[a_1^{\ep_1},\dots,a_{r-1}^{\ep_{r-1}},{a_r'}^{\ep_r'}]$ that define the same oriented matroid. However, condition \eqref{it:condle3} implies that in our expressions, $a_r$ and $\ep_r$ are completely determined by $a_{r-1}$.
\\

We proceed to prove that with our conditions \eqref{it:condle1}, \eqref{it:condle2} and \eqref{it:condle3}, $p$ and $\suc_{i-1}(a_i)$ are not $(\alpha_{i-1}(a_i)\ep_i)$-inseparable in $\cM[p]/\{a_1,\dots,a_{i}\}$ for $i<r$. The proof is by induction on $\rr$.

First we deal with $\rr=2$. Because of~\eqref{it:condle3}, the extension must be by $p=[a_1^{\ep_1},\suc(a_1)^{\alpha(a_1)\ep_1}]$. Hence, $\cM[a_1^{\ep_1},s(a_1)^{\alpha(a_1)\ep_1}]/a_1=(\cM/a_1)[s(a_1)^{\alpha(a_1)\ep_1}]$, and thus $p$ and $s(a_1)$ are $(-\alpha(a_1)\ep_1)$-inseparable in $\cM/a_1$.\\

When $r> 2$, we only need to prove that $p$ and $\suc(a_1)$ are not $(\alpha(a_1)\ep_1)$-inseparable in $\cM[p]/a_1$. Indeed, for $a_i$ with $i>1$, the result follows by induction using the contraction $\cM/S_{i-1}$.

Observe that if $\suc(a_1)$ and $p$ are inseparable in $\cM/a_1$, then $a_2$ must be $\suc(a_1)$ or inseparable with $\suc(a_1)$, by Observation~\ref{obs:inlble2proof0}.

\begin{itemize}
 \item 
If $a_2=\suc(a_1)$ then, because of the restrictions on the expression of the lexicographic extension, $\ep_2=\alpha(a_1)\ep_1$ and $p$ and $\suc(a_1)$ are $(-\alpha(a_1)\ep_1)$-inseparable in $\cM/a_1$.

\item Otherwise, if $a_2\neq \suc(a_1)$ and $a_2$ and $\suc(a_1)$ are inseparable in $\cM/a_1$, then either $a_2=\suc_1(\suc(a_1))$ or $\suc_1(a_2)=\suc(a_1)$. 

\begin{itemize}
\item If $a_2=\suc_1(\suc(a_1))$, because of the restrictions on the expression of the lexicographic extension, $\ep_2=-\alpha_1(\suc(a_1))\alpha(a_1)\ep_1$. Moreover, since $a_2$ and $\suc(a_1)$ are $\alpha_1(\suc(a_1))$-inseparable, and again by Observation~\ref{obs:inlble2proof0}, the inseparability of $p$ and $\suc(a_1)$ in $\cM/a_1$ can only be \[\alpha_1(\suc(a_1))\ep_2=-\alpha_1(\suc(a_1))\alpha_1(\suc(a_1))\alpha(a_1)\ep_1=-\alpha(a_1)\ep_1.\]

\item Finally, there is the case $\suc_1(a_2)=\suc(a_1)$. If $p$ and $\suc(a_1)$ were $(\alpha(a_1)\ep_1)$-inseparable in $\cM/a_1$, then $\ep_2=\alpha_1(a_2)\alpha(a_1)\ep_1$. And if $\ep_2=\alpha_1(a_2)\alpha(a_1)\ep_1$, we know by induction hypothesis that the inseparability of $p$ and $\suc_1(a_2)=\suc(a_1)$ cannot be \[\alpha_1(a_2)\ep_2=\alpha_1(a_2)\alpha_1(a_2)\alpha(a_1)\ep_1=\alpha(a_1)\ep_1.\]
\end{itemize}
\end{itemize}

This concludes the proof that there is no pair of different expressions of lexicographic extensions fulfilling \eqref{it:condle1}, \eqref{it:condle2} and \eqref{it:condle3} that define the same oriented matroid.
\\

It only remains to find some bound on the number of lexicographic extensions that fulfill these conditions. For this, observe that if $1<i<r$, then $a_i$ can be any element of $\cM/\{a_1,\dots,a_{i-1}\}$ with $\ep_i$ any sign except for at most two forbidden configurations (if $\suc_{i-2}(a_{i-1})$ and $\suc_{i-1}(\suc_{i-2}(a_{i-1}))$ exist). That means that there are at least $2(n-(i-1))-2=2(n-i)$ choices. For $a_1^{\ep_1}$ there are $2n$ options. And $a_r^{\ep_r}$ is fixed (and well defined, since the inseparability graph of a uniform matroid of rank~$2$ is always a cycle and thus $\suc_{r-2}(a_{r-1})$ always exists). 

Summing up,
\[\lle{n}{r}\geq 2n \prod_{i=2}^{r-1}2(n-i)=
2^{r-1}\frac{n!}{(n-1)(n-r)!}.\qedhere\]
\end{proof}

\subsection{\texorpdfstring{Many neighborly polytopes in $\cG$}{Many neighborly polytopes in G}}

Once we have bounds for $\lle{n}{r}$, we can obtain bounds for $\lnei{n}{d}$ using the Gale Sewing construction. But first we do a case where we know the exact number.

\begin{lemma}\label{lem:nlpolygons}
 The number of labeled balanced matroids of rank $r$ with $r+3$ elements is $\frac12{(r+2)!}$.
\end{lemma}
\begin{proof}
 Balanced matroids of rank $r$ with $r+3$ elements are dual to polygons with $r+3$ vertices in $\RR^2$. There are clearly $\frac12{(r+2)!}$ different combinatorial types of labeled polygons with $r+3$ vertices.
\end{proof}

Using the Double Extension Theorem~\ref{thm:thethm}, the number of different lexicographic extensions can be used to give lower bounds on the number of neighborly polytopes.

\begin{lemma}\label{lem:lbbm}
For $r\geq2$ and $m\geq2$, the number of labeled balanced matroids of rank $r$ with $r+1+2m$ elements is  $\lnei{2m+r+1}{2m}$ and fulfills
\begin{equation}\label{eq:lbbm}
\lnei{2m+r+1}{2m}\geq \lnei{2m+r-1}{2m-2}\frac{r+2m}{2}\lle{r+2m-1}{r}.
\end{equation}
\end{lemma}
\begin{proof}
 The characterization is direct by duality. The bound is a consequence of the Double Extension Theorem~\ref{thm:thethm}. Fix a balanced matroid~$\cM$ of rank $r$ with $r+1+2(m-1)$ elements such that each element receives a label in $\{1,\dots,r+1+2(m-1)\}$. There are at least $\lle{r+2m-1}{r}$ different labeled lexicographic extensions of $\cM$. 
 Consider any such extension $\cM[p]$, the extension of $\cM$ by $p=[a_1^{\ep_1},\dots,a_r^{\ep_r}]$, where $p$ receives the label~$r+2m$.

 Finally let $\cM[p][q]$ be the extension by $q=[p^-,a_1^{-},\dots,a_{r-1}^{-}]$, which is balanced by Theorem~\ref{thm:thethm}. 
 We count the number of ways of labeling $\cM[p][q]$ such that $p$ gets label $r+2m+1$ and the labeling of $\cM[p][q]$ on~$\cM$ preserves the relative order of the original labeling of~$\cM$. There are clearly $r+2m$ ways to do this, which correspond to the possible labels that can be given to~$q$. We claim that each labeled matroid constructed this way is isomorphic to at most one other labeled matroid constructed this way.

 Indeed, observe that $p$ and $q$ are inseparable because of Lemma~\ref{lem:leinseparable}. Moreover, by Theorem~\ref{thm:IG}, $p$ is inseparable from at most two elements in $\cM[p][q]$ because $2\le r$ and $2\le m$. If $p$ is inseparable to only one element in $\cM[p][q]$, then any labeling of $q$ provides a different labeled matroid. Otherwise, if $p$ is inseparable from $q$ and from $a$, it might have been counted as a double extension of $\cM[p][q]\setminus\{p,q\}$ or as a double extension of $\cM[p][q]\setminus\{p,a\}$.
 
 Summing up, there are at least \[(r+2m)\lnei{2m+r-1}{2m-2}\lle{r+2m-1}{r}\] labeled balanced oriented matroids, where each matroid is counted at most twice. This yields the claimed formula.
\end{proof}

This result allows us to give our first explicit lower bound on the number of neighborly polytopes. It is using the simpler bound from~\eqref{eq:bndle1}, because even if it is smaller, with the bound from~\eqref{eq:bndle2} the formulas become more complicated and it adds nothing substantial to the result.

\begin{proposition}\label{prop:brutebound}
 The number of labeled neighborly polytopes in even dimension $d=2m\geq 2$ with $n=r+d+1$ vertices fulfills 
\begin{equation}\label{eq:lblnei1}
  \lnei{2m+r+1}{2m}\geq\prod_{i=1}^{m} {\frac{(r+2i)!}{(2i)!}}
\end{equation}
\end{proposition}
\begin{proof}
Observe that by rigidity (Theorem~\ref{thm:neigharerigid}), counting labeled neighborly polytopes is equivalent to counting labeled neighborly oriented matroids. By duality, this is in turn equivalent to counting balanced oriented matroids. This we do.

Lemma~\ref{lem:nlpolygons} proves the required formula in the initial case $m=1$, and yields $\lnei{2+\rr+1}{2}=\frac12(\rr+2)!$. For $m\ge2$, we observe that by Proposition~\ref{prop:lble},
\[\frac{r+2m}{2}\lle{r+2m-1}{r}\geq 
\frac{(\rr+2m)!}{(2m)!}.\]
Finally, we apply Lemma~\ref{lem:lbbm} to obtain~\eqref{eq:lblnei1}.
\end{proof}

Although Proposition~\ref{prop:brutebound} provides us with the desired bound, it is hard to understand its order of magnitude at first sight. This is the reason why we present the following simplified bound.

\begin{theorem}\label{thm:lblnei}
The number of labeled neighborly polytopes in even dimension $d$ with $n$ vertices fulfills 
\begin{equation*}
\lnei{r+d+1}{d}\geq \frac{\left( r+d \right) ^{\left( \frac{r}{2}+\frac{d}{2} \right) ^{2}}}{{r}^{{(\frac{r}{2})}^{2}}{d}^{{(\frac{d}{2})}^{2}}{\e^{3\frac{r}{2}\frac{d}{2}}}},
\end{equation*}
that is, 
\begin{equation*}
 \lnei{n}{d}\geq \frac{\left(n-1\right) ^{\left( \frac{n-1}{2}\right) ^{2}}}{{(n-d-1)}^{{\left(\frac{n-d-1}{2}\right)}^{2}} d  ^{{(\frac{d}{2})}^{2}}{\e^{\frac{3d(n-d-1)}{4}}}}.
\end{equation*}
\end{theorem}
\begin{proof}
We start from Equation~\eqref{eq:lblnei1}, and approximate the natural logarithm of $\lnei{r+1+2m}{2m}$. Using the fact that $\int_{a-1}^{b} f(s)\  \mathrm{d}s \le \sum_{i=a}^{b} f(i)$ for any increasing function $f$, we obtain
\begin{align*}
\ln\left( \lnei{r+1+2m}{2m}\right)\geq&\, \ln\left(\prod_{i=1}^{m} \frac{(\rr+2i)!}{(2i)!}
\right)=\sum_{i=1}^{m} {\sum_{j=0}^{r-1}{ \ln\left(r+2i-j\right)}}\\
=&\,\sum_{i=1}^{m} {\sum_{j=1}^{r}{ \ln\left(2i+j\right)}}
\geq\int_{i=0}^{m} {\int_{j=0}^{r}{ \ln\left(2i+j\right) \mathrm{d}j} \mathrm{d}i}\\
 =&\,  \frac{1}{4}\left(2m+r\right)^{2}\ln\left( 2m+r\right)-\frac{1}{4}r^{2}\ln \left( r \right) 
\\&\,-{m}^{2} \ln  \left( 2m\right) -\frac{3}{2}mr.
\end{align*}
Hence
\begin{align*}
\lnei{r+1+2m}{2m}& \geq\frac{\left( 2m+r
 \right) ^{\frac14 \left( 2m+r \right) ^{2}}}{{r}^{\frac14{r}^{2}}\left( 2m \right) ^{{m}^{2}}{\e^{\frac32mr}}},
\end{align*}
and we conclude that
\begin{equation}\label{eq:lblnei}
\lnei{r+d+1}{d}\geq \frac{\left( r+d \right) ^{\left( \frac{r}{2}+\frac{d}{2} \right) ^{2}}}{{r}^{{(\frac{r}{2})}^{2}}{d}^{{(\frac{d}{2})}^{2}}{\e^{3\frac{r}{2}\frac{d}{2}}}}.
\end{equation}
\vspace{-.5cm}
\end{proof}

The following corollary is a further simplification of the bound. 
\begin{corollary}
The number of labeled neighborly polytopes in even dimension $d$ with $n$ vertices fulfills 
\[\lnei{n}{d}\geq  \left( \frac{n-1}{\e^{3/2}}\right)^{\frac12(n-d-1)d}.\]
\end{corollary}
\begin{proof}
Since ${r}^{{(\frac{r}{2})}^{2}}{d}^{{(\frac{d}{2})}^{2}}\leq{(r+d)}^{{(\frac{r}{2})}^{2}+{(\frac{d}{2})}^{2}}$
 we obtain
\[
\lnei{r+d+1}{d}\geq \frac{\left( r+d \right) ^{\left( \frac{r}{2}+\frac{d}{2} \right) ^{2}}}{{r}^{{(\frac{r}{2})}^{2}}{d}^{{(\frac{d}{2})}^{2}}{\e^{3\frac{r}{2}\frac{d}{2}}}}\geq  \frac{\left( r+d \right) ^{\frac{rd}{2}}}{{\e^{\frac{3{rd}}{4}}}}.
\qedhere
\]
\end{proof}

Observe that this bound is useful both for neighborly polytopes whose number of vertices is very large with respect to the dimension as well as for neighborly polytopes with fixed corank and large dimension.\\

A final observation is that we can translate these bounds for even dimensional neighborly polytopes to bounds for neighborly polytopes in odd dimension just by taking pyramids, because a pyramid over an even dimensional neighborly polytope is always neighborly. (If simpliciality was needed, any extension in general position of the Gale dual of an even-dimensional neighborly polytope would work too.)

\begin{corollary}\belowdisplayskip=-12pt \label{cor:oddneighpoly}
The number of labeled neighborly polytopes in odd dimension $d$ with $n$ vertices fulfills 
\begin{align*}
\lnei{n}{d}\geq\lnei{n\!-\!1}{d\!-\!1}&\geq
\frac{\left(n-2\right) ^{\left( \frac{n-2}{2}\right) ^{2}}}{{(n-d-1)}^{{\left(\frac{n-d-1}{2}\right)}^{2}} {(d-1)}^{{(\frac{d-1}{2})}^{2}}{\e^{\frac{3(d-1)(n-d-1)}{4}}}}\\
&\geq  \left( \frac{n-2}{\e^{3/2}}\right)^{\frac12(n-d-1)(d-1)}.
\end{align*}\qed
\end{corollary}
\vspace{12pt}

\subsection{Many non-realizable neighborly matroids}

Exactly the same reasoning that leads to the bounds in Theorem~\ref{thm:lblnei} can be applied to give lower bounds for non-realizable neighborly matroids. From now on, let 
\defn{$\lnr{n}{r}$}\index{$\lnr{n}{r}$} represent the number of labeled non-realizable neighborly oriented matroids of rank $r$ with $n$ elements.

\begin{theorem}\label{thm:nonrealizablebound}
 The number of labeled non-realizable neighborly oriented matroids of odd rank~$\rd$ with $n$ elements is at least
\begin{equation*}
\lnr{n}{s}\geq \frac{\left( n-1 \right) ^{\frac12{(\rd-5)(n-\rd)}}}{\left(\frac{n-\rd+4}{2}\right) ^{4}\e^{\frac34(\rd-5)(n-\rd)}}.
\end{equation*}
\end{theorem}
\begin{proof}[Proof sketch]
The principal observation is that an analogue of the inequality~\eqref{eq:lbbm} of Lemma~\ref{lem:lbbm} applies. That is, if $r\geq2$, $m\geq2$ and $n=2m+r+1$, then
\begin{equation*}
\lnr{n}{2m+1}\geq \lnr{n-2}{2m-1}\frac{n-1}{2}\lle{n-2}{r}.
\end{equation*}
This uses the Double Extension Theorem~\ref{thm:thethm} and the fact that all the lexicographic extensions of a non-realizable matroid are non-realizable.\\

Moreover, by Theorem~\ref{thm:nonrealizable}, $\lnr{r+5}{5}\geq 1$ for all $r\geq5$.
Which means that for $m\geq3$ we can mimic the proof of Theorem~\ref{thm:lblnei} to get 
\begin{align*}
\lnr{r+1+2m}{2m+1}\geq& \prod_{i=3}^{m} {\prod_{j=1}^{r}{ \ln\left(2i+j\right)}}\\\geq&\exp\left(\int_{i=2}^{m} {\int_{j=0}^{r}{ \ln\left(2i+j\right) \mathrm{d}j} \mathrm{d}i}\right)\\
=&\frac{2^8\left( 2m+r \right) ^{\frac{\left( 2m+r \right) ^
{2}}{4}}}{ \left( r+4 \right) ^{\frac{\left( r+4 \right) ^{2}}{4}} {(2m)}^{{m}^{2
}}{\e^{\frac{3\left( m-2 \right)r}{2} }}}
\geq\frac{\left( 2m+r \right) ^{(m-2)r}}{\left(\frac{r+4}{2}\right) ^{4}\e^{\frac{3 \left( m-2 \right)r}{2} }}.
\end{align*}
\end{proof}

Again, we can use pyramids to extend these bounds to non-realizable neighborly oriented matroids of even rank.
\begin{corollary}\belowdisplayskip=-12pt \label{cor:evennonreal}
The number of non-realizable neighborly oriented matroids of even rank $\rd$ with $n$ vertices fulfills 
\begin{equation*}
\lnr{n}{s}\geq \frac{\left( n-2 \right) ^{\frac12{(\rd-6)(n-\rd)}}}{\left(\frac{n-\rd+4}{2}\right) ^{4}\e^{\frac34(\rd-6)(n-\rd)}}.
\end{equation*}\qed
\end{corollary}
\vspace{12pt}

\subsection{Many polytopes?}
As we have already said, these bounds for the number of neighborly polytopes are also interesting as bounds for the number of polytopes. A natural question to ask is whether the same tools can be used to get even better bounds for the number of polytopes. 

Of course, one can mimic the Gale Sewing construction to build many not necessarily neighborly polytopes: starting with $\{\vv e_1,\dots,\vv e_r,-\sum_{i=1}^r \vv e_i\}$ and making a sequence of lexicographic extensions one gets the Gale dual of certain point configuration. If this oriented matroid is positively $2$-spanning, the point configuration is in convex position. 

Our main problem is that it is hard to certify how many of these polytopes are different.
Although Theorem~\ref{thm:lble2} can still be used to bound the number of lexicographic extensions, oriented matroids of polytopes are not rigid in general. With neighborly polytopes, we could use that different matroids give rise to different face lattices, which is not longer true in the general case.
\\

In this scenario, the following questions arise naturally.

\begin{question}
Which is the maximal number of realizable oriented matroids that can share the same face lattice?
\end{question}

\begin{question}
Which is the minimal number of (regular) triangulations that a point configuration can have?
\end{question}

\iftoggle{bwprint}{%
\renewcommand{\partfigure}{prelude/Figures/fulldecompexamples}
}{%
\renewcommand{\partfigure}{prelude/Figures/fulldecompexamples_col}
}
\renewcommand{\namepart}{Almost neighborly}
\part{Almost neighborly}\label{partII}
\chapter{Introduction}\label{ch:intro_degree}

\section{Overview}

Consider the following relaxation of the definition of $k$-neighborly point configurations, obtained by changing the condition of ``being the set of vertices of a face'' to ``belonging to a common face'': 

\begin{definition}
A point configuration $\vv A$ is \defn{$k$-almost neighborly}\index{almost neighborly} if every subset of $\vv A$ of size $\leq k$ lies in a common facet of $\conv (\vv A)$. A polytope ${\vv P}$ is \defn{$k$-almost neighborly} if its set of vertices is.
\end{definition}

For example, the $3$-simplex in Figure~\ref{sfig:exdegree1} is $3$-almost neighborly because every subset of $3$ vertices lies in a common face; the prism over a $2$-simplex of Figure~\ref{sfig:exdegree2} is $2$-almost neighborly because every pair of vertices lies in a common face; and the $3$-cube in Figure~\ref{sfig:exdegree3} is just $1$-almost neighborly.

\iftoggle{bwprint}{%
\begin{figure}[htpb]
\centering
 \subbottom[$\simp{3}$]{\label{sfig:exdegree1}\includegraphics[scale=.85]{Figures/simplex}}\qquad
 \subbottom[$\simp{2}\times \simp{1}$]{\label{sfig:exdegree2}\includegraphics[scale=.85]{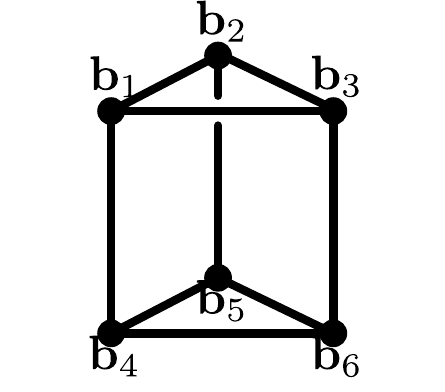}}\qquad
 \subbottom[$\cube{3}$]{\label{sfig:exdegree3}\includegraphics[scale=.85]{Figures/cube}}
 \caption[Polytopes that are $3$-, $2$-, and $1$-almost neighborly.]{Three $3$-polytopes that are $3$-, $2$-, and $1$-almost neighborly, respectively.}
 \label{fig:exdegree}
\end{figure}
}{%
\iftoggle{print}{%
\begin{figure}[htpb]
\centering
 \subbottom[$\simp{3}$]{\label{sfig:exdegree1}\includegraphics[scale=.85]{Figures/simplex}}\qquad
 \subbottom[$\simp{2}\times \simp{1}$]{\label{sfig:exdegree2}\includegraphics[scale=.85]{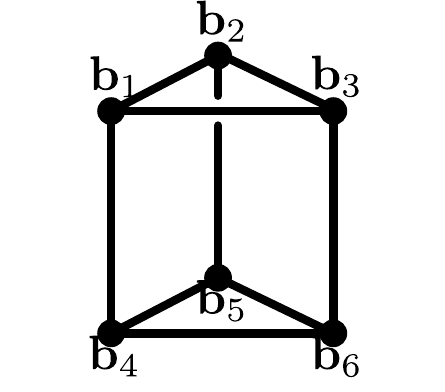}}\qquad
 \subbottom[$\cube{3}$]{\label{sfig:exdegree3}\includegraphics[scale=.85]{Figures/cube}}
 \caption[Polytopes that are $3$-, $2$-, and $1$-almost neighborly.]{Three $3$-polytopes that are $3$-, $2$-, and $1$-almost neighborly, respectively.}
 \label{fig:exdegree}
\end{figure}
}{%
\begin{figure}[htpb]
\centering
 \subbottom[$\simp{3}$]{\label{sfig:exdegree1}\includegraphics[scale=.85]{Figures/simplex_col}}\qquad
 \subbottom[$\simp{2}\times \simp{1}$]{\label{sfig:exdegree2}\includegraphics[scale=.85]{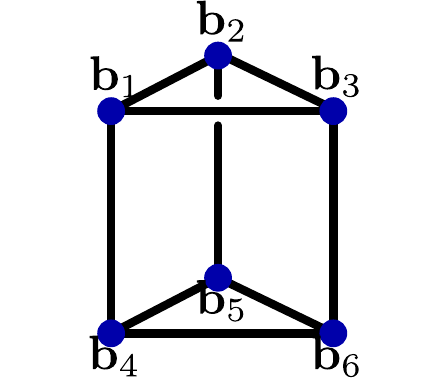}}\qquad
 \subbottom[$\cube{3}$]{\label{sfig:exdegree3}\includegraphics[scale=.85]{Figures/cube_col}}
 \caption[Polytopes that are $3$-, $2$-, and $1$-almost neighborly.]{Three $3$-polytopes that are $3$-, $2$-, and $1$-almost neighborly, respectively.}
 \label{fig:exdegree}
\end{figure}
}
}

This concept is not very different from $k$-neighborliness and, in fact, these notions coincide when ${\vv A}$ is the vertex set of a \defn{$(k-1)$-simplicial}\index{polytope!simplicial} polytope~$\vv P$, which means that each $(k-1)$-dimensional face of~$\vv P$ is a simplex.

The name ``almost neighborly'' was coined by Gr\"unbaum~\cite[Exercices 7.3.5 and 7.3.6]{GruenbaumEtal2003}. According to him, this notion had already considered by Motzkin in 1965 under the name of \defn{$k$-convex sets}~\cite{Motzkin1965}. In 1972, Breen proved that a point configuration is $k$-almost neighborly if and only if all its subconfigurations of size $\leq2d+1$ are~\cite{Breen1972}; a result that should be compared with Lemma~\ref{lem:Hellyneighborly}, its analogue for neighborly polytopes. 

A related concept is that of \defn{weak neighborliness}: $\vv A$ is called \defn{weakly neighborly} if every subset of $k+1$ points is contained in a face $\conv(\vv A)$ of dimension at most $2k$, for all $k$~\cite{Bayer1993}\index{weakly neighborly}. In particular, Bayer~\cite[Theorem 15]{Bayer1993} already classified all $2$-almost neighborly $3$-dimensional polytopes as being triangular prisms and pyramids.\\

Besides the references above, almost neighborliness has been somewhat ignored by researchers in geometric combinatorics. However, it turns out that it is strongly related to an Ehrhart-theoretic invariant of a lattice polytope (called its degree) that has been subject of wide interest over the last years.
Our goal 
is to transfer the intuition acquired in the investigation of the Ehrhart-theoretic degree to this more general and combinatorial setting. 
\\

Let us explain this briefly. 
Given a lattice polytope ${\vv P}$ (a polytope with vertices in the lattice $\ZZ^d$), the Ehrhart polynomial\index{Ehrhart polynomial} counts the number of lattice points in dilates of ${\vv P}$. As it turns out, the complexity of the Ehrhart polynomial is directly related to the largest natural number $k$ such that~$k{\vv P}$ has no interior lattice points. Now, here is our naive observation: for any such~$k$, any set of $k$~vertices of~${\vv P}$ clearly lies in a common facet, since otherwise the sum of those vertices would lie in the interior of $k{\vv P}$. In other words, ${\vv P}$~is $k$-almost neighborly.
\\

We hope to shed some light on the structure of almost neighborly polytopes by mirroring some results on lattice polytopes, even if not all questions can yet be fully answered. 
For instance, Theorem~\ref{thm:kneighd/2} is a basic result on neighborly polytopes to the effect that any $(>\lfloor \frac{d}{2} \rfloor)$-neighborly $d$-polytope is a simplex. What happens for almost neighborly polytopes and point configurations? That is, can we say something about $k$-almost neighborly polytopes and point configurations when $k$ is large in comparison with their dimension $d$?
\\

Our main results can be interpreted as answers to this question. For example, we completely classify $d$- and $(d-1)$-almost neighborly point configurations (Proposition~\ref{prop:dd=0} and Theorem~\ref{thm:dd=1}), and give structural constraints for $k$-almost neighborly point configurations with $k> \lfloor \frac{2d}{3} \rfloor$ (Theorem~\ref{thm:d+1-3dd}) and with $k\geq \lfloor \frac{n+d+1}{4} \rfloor$ (Theorem~\ref{thm:2DD}). In the remainder of this outline, we present these results, show motivations and parallelisms with Ehrhart-theoretic notions, and interpret them in different contexts. 

As will be explained below, Ehrhart theory motivates to focus on the number $d-k$ (which will be called the \defn{degree} of $\vv A$) rather than on~$k$.

\subsection{The degree, the codegree, weak Cayley configurations and codegree decompositions}

Let $\vv A$ be a point configuration in $\RR^d$. 
We say that a non-empty subset $\vv S\subset \vv A$ is an \defn{interior face}\index{face!interior} of $\vv A$ if $\conv(\vv S)$ does not lie on the boundary of $\conv (\vv A)$, so that $\conv (\vv S)\cap\intr \left(\conv (\vv A)\right)$ is non-empty. For example, in Figure~\ref{fig:exdegree}, the faces spanned by $\{\vv a_1,\vv a_2,\vv a_3,\vv a_4\}$, $\{\vv b_1,\vv b_5,\vv b_6\}$ and $\{\vv c_1,\vv c_8\}$ are interior faces of the polytopes in Figures~\ref{sfig:exdegree1}, \ref{sfig:exdegree2} and \ref{sfig:exdegree3} respectively, while $\{\vv a_1,\vv a_2,\vv a_3\}$, $\{\vv b_1,\vv b_6\}$ and $\{\vv c_1,\vv c_7\}$ are not.

\begin{definition}
 \label{def:combdegree}
The \defn{degree}\index{degree} $\degc(\vv A)$\index{$\degc(\vv A)$} of a $d$-dimensional point configuration~$\vv A$ is the maximal codimension of an interior face of~$\vv A$. 

The \defn{codegree}\index{codegree} of $\vv A$ is defined as $\codegc(\vv A) := d+1-\degc(\vv A)$\index{$\codegc(\vv A)$}, 
and equals the maximal positive integer $\kk$ such that every subset of $\vv A$ of size $< \kk$ is contained in a common facet of $\conv(\vv A)$; equivalently, the codegree is the maximal $\kk$ such that $\vv A$ is $(\kk-1)$-almost neighborly.

The  \defn{degree} and \defn{codegree} of a convex polytope ${\vv P}$ are defined to be the degree and codegree of its set of vertices, $\degc(\vv P):=\degc(\verts({\vv P}))$ and $\codegc(\vv P):=\codegc(\verts({\vv P}))$. In particular, the degree of a polytope only depends on its combinatorial type. 
\end{definition}

Observe that $0\leq \degc(\vv A)\leq d$. We consider that every point of a $0$-dimensional point configuration is an interior face, and hence, if $\{\vv a\}$ is a single point, then $\degc(\{\vv a\})=0$ and $\codegc(\{\vv a\})=1$. 
Higher dimensional simplices also have degree~$0$. And Proposition~\ref{prop:dd=0} below shows that they are the only configurations with $\degc(\vv A)=0$.

\medskip\noindent
\textbf{Proposition \ref{prop:dd=0}.}
\emph{
If $\vv A$ is a point configuration with $\degc(\vv A)=0$, then it is the set of vertices of a simplex (possibly with repetitions).}\\

This means that the first interesting configurations have $\degc(\vv A)=1$, such as those depicted in Figure~\ref{fig:exdeg1}. In Theorem~\ref{thm:dd=1} we completely classify these configurations, whose relation to totally splittable point configurations is also discussed in Chapter~\ref{ch:deg1}.

\medskip\noindent
\textbf{Theorem \ref{thm:dd=1}.}
\emph{
Let $\vv A$ be a $d$-dimensional configuration of $n$ points. If $\degc(\vv A)\leq1$, then one of the following holds (up to repeated points)
\begin{enumerate}
\item\label{it:dd=1:d=1} $d\leq 1$; 
or
\item\label{it:dd=1:d=2} $d \ge 2$ and $\vv A$ is a $k$-fold pyramid over a two-dimensional point configuration without interior points; 
or
\item\label{it:dd=1:prism} $d \ge 3$ and $\conv(\vv A)$ is a $k$-fold pyramid over a prism over a simplex with the non-vertex points of $\vv A$ all on the ``vertical'' edges of the prism; or 
\item\label{it:dd=1:simplex} $d \ge 3$ and $\conv(\vv A)$ is a simplex with all non-vertex points of $\vv A$ on the edges adjacent to a vertex $\vv a$ of $\conv(\vv A)$.
\end{enumerate}
}
\iftoggle{bwprint}{%
\begin{figure}[htpb]
\centering
\subbottom[Ex. of type~\eqref{it:dd=1:d=2}.]{\quad\includegraphics[width=.2\linewidth]{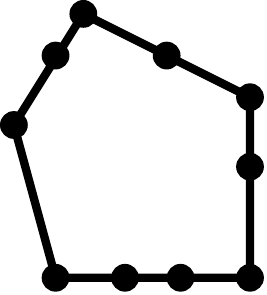}\quad}\qquad  
\subbottom[Ex. of type~\eqref{it:dd=1:prism}.]{\quad\includegraphics[width=.2\linewidth]{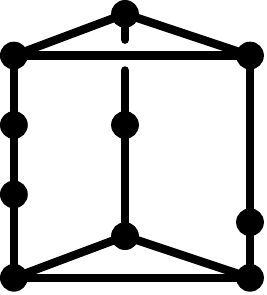}\quad}\qquad
\subbottom[Ex. of type~\eqref{it:dd=1:simplex}.]{\quad\includegraphics[width=.2\linewidth]{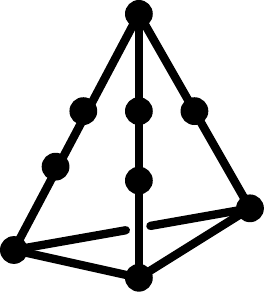}\quad} 
\caption[Some point configurations of degree $1$]{Some point configurations of degree $1$, and their classification according to Theorem~\ref{thm:dd=1}.}\label{fig:exdeg1}
\end{figure}
}{%
\iftoggle{print}{%
\begin{figure}[htpb]
\centering
\subbottom[Ex. of type~\eqref{it:dd=1:d=2}.]{\quad\includegraphics[width=.2\linewidth]{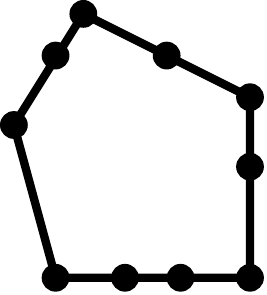}\quad}\qquad  
\subbottom[Ex. of type~\eqref{it:dd=1:prism}.]{\quad\includegraphics[width=.2\linewidth]{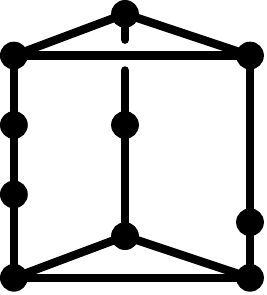}\quad}\qquad
\subbottom[Ex. of type~\eqref{it:dd=1:simplex}.]{\quad\includegraphics[width=.2\linewidth]{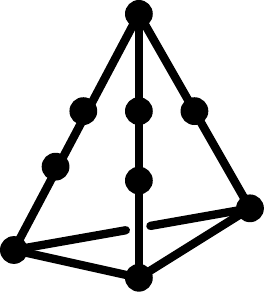}\quad} 
\caption[Some point configurations of degree $1$]{Some point configurations of degree $1$, and their classification according to Theorem~\ref{thm:dd=1}.}\label{fig:exdeg1}
\end{figure}}{%
\begin{figure}[htpb]
\centering
\subbottom[Ex. of type~\eqref{it:dd=1:d=2}.]{\quad\includegraphics[width=.2\linewidth]{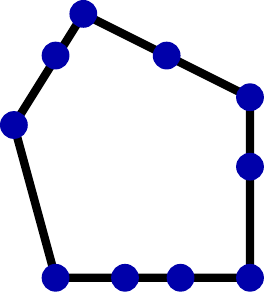}\quad}\qquad  
\subbottom[Ex. of type~\eqref{it:dd=1:prism}.]{\quad\includegraphics[width=.2\linewidth]{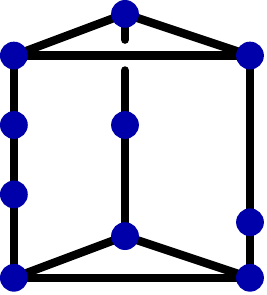}\quad}\qquad
\subbottom[Ex. of type~\eqref{it:dd=1:simplex}.]{\quad\includegraphics[width=.2\linewidth]{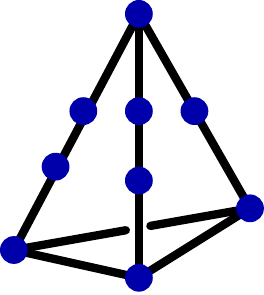}\quad} 
\caption[Some point configurations of degree $1$]{Some point configurations of degree $1$, and their classification according to Theorem~\ref{thm:dd=1}.}\label{fig:exdeg1}
\end{figure}
}
}

In Chapter~\ref{ch:cayley}, we seek more general constraints for configurations whose degree is small with respect to the ambient dimension. Our first result in this direction takes the number of points in the configuration into account. 

\medskip\noindent
\textbf{Corollary \ref{cor:easybound}.}
\emph{Any $d$-dimensional configuration $\vv A$ of $n=r+d+1$ points with
\[d \ge r+2\degc(\vv A)\]
is a pyramid.}

\bigskip

To state our next result, we need the following definition:

\begin{definition}\label{def:strongweakcayley}
 A point configuration $\vv A$ is a \defn{weak Cayley configuration}\index{Cayley configuration!weak} of length $m\ge1$\index{Cayley configuration!weak!length} 
if there exists a partition $\vv A = \vv A_0\uplus \vv A_1 \uplus \cdots \uplus \vv A_m$, such that for any 
$1\leq i\leq m$, $\vv A\setminus \vv A_i$ is the set of points of a proper face of~$\conv (\vv A)$.
The sets $\vv A_1\dots \vv A_m$ are called the \defn{factors}\index{Cayley configuration!weak!factor} of the configuration. 

While we allow~$\vv A_0$ to be the empty set, the factors $\vv A_1, \ldots, \vv A_m$ have to be non-empty, because otherwise $\conv(\vv A\setminus \vv A_i)$ would not be a proper face of~$\conv(\vv A)$.
\end{definition}

\iftoggle{bwprint}{%
\begin{figure}[htpb]
\centering
 \subbottom[$\vv A$]{\label{sfig:weakA}\includegraphics[width=.2\linewidth]{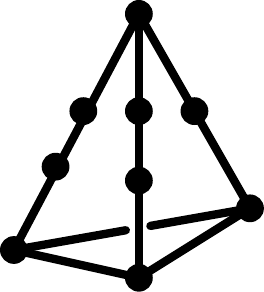}}\quad
 \subbottom[$\vv A_0$]{\label{sfig:weakA0}\includegraphics[width=.15\linewidth]{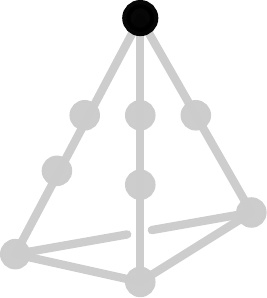}}\quad
 \subbottom[$\vv A_1$]{\label{sfig:weakA1}\includegraphics[width=.15\linewidth]{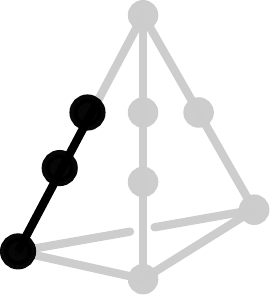}}\quad
 \subbottom[$\vv A_2$]{\label{sfig:weakA2}\includegraphics[width=.15\linewidth]{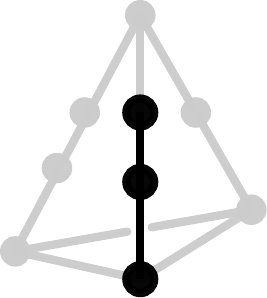}}\quad
 \subbottom[$\vv A_3$]{\label{sfig:weakA3}\includegraphics[width=.15\linewidth]{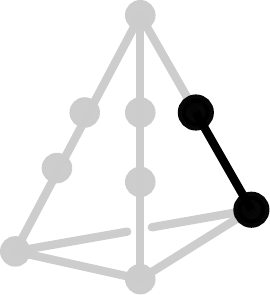}}
  \caption{A weak Cayley configuration of length $3$.}
 \label{fig:exweak}
\end{figure}
}{%
\begin{figure}[htpb]
\centering
 \subbottom[$\vv A$]{\label{sfig:weakA}\includegraphics[width=.2\linewidth]{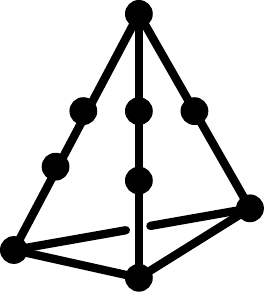}}\quad
 \subbottom[$\vv A_0$]{\label{sfig:weakA0}\includegraphics[width=.15\linewidth]{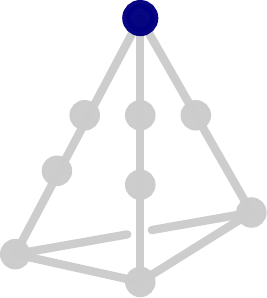}}\quad
 \subbottom[$\vv A_1$]{\label{sfig:weakA1}\includegraphics[width=.15\linewidth]{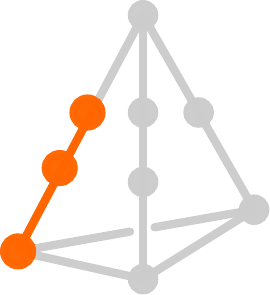}}\quad
 \subbottom[$\vv A_2$]{\label{sfig:weakA2}\includegraphics[width=.15\linewidth]{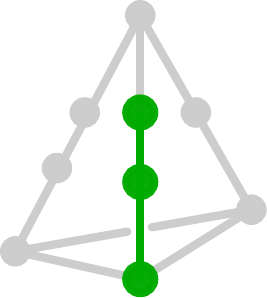}}\quad
 \subbottom[$\vv A_3$]{\label{sfig:weakA3}\includegraphics[width=.15\linewidth]{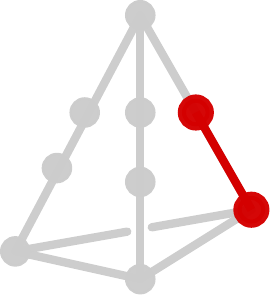}}
  \caption{A weak Cayley configuration of length $3$.}
 \label{fig:exweak}
\end{figure}
}

Figure~\ref{fig:exweak} shows  a $3$-dimensional Cayley configuration of length~$3$ and degree $1$. More generally,  we show in Chapter~\ref{ch:cayley} that ``long'' weak Cayley configurations are examples of point configurations with small degree:

\medskip\noindent
\textbf{Proposition \ref{prop:deg-cayley}.}
\emph{Let $\vv A\subset \RR^d$ be a weak Cayley configuration of length~$m$. Then $\degc(\vv A)\leq d+1-m$.}

\medskip

We are ready to state our main theorem:  every configuration of degree strictly less than~$\frac{d}{3}$ is a non-trivial weak Cayley configuration.

\medskip\noindent
\textbf{Theorem \ref{thm:d+1-3dd}.} \emph{Let $\vv A$ be a $d$-dimensional configuration. If $\degc(\vv A)<\frac{d}{3}$, then $\vv A$ is a weak Cayley configuration of length at least $d-3\degc(\vv A)+1$.}

\medskip

Theorem~\ref{thm:d+1-3dd} should be seen as a converse statement to Proposition~\ref{prop:deg-cayley}. 
The reader may also have noticed that the assumption in Theorem~\ref{thm:d+1-3dd} might be strengthened for $\dd=1$. 
Indeed, if $d > 2$ and $\dd=1$, then $\vv A$ is a weak Cayley configuration of length at least $d-1$. This observation (among others, as will be explained below) motivates our first main conjecture:

\begin{conjecture}\label{conj:d+1-2dd}
Any point configuration $\vv A$ of dimension $d>2\degc(\vv A)$ is a weak Cayley configuration of length at least $d+1-2\degc(\vv A)$. 
\end{conjecture}

The conjectured bound (if correct) is
sharp by Example~\ref{ex:sharpconj}, which shows that an even-dimensional cyclic $d$-polytope with $2d+1$ vertices has degree $\dd=\frac{d}{2}$, but is not a weak Cayley configuration of length~$\geq2$. 
\medskip

A stronger combinatorial constraint is conjectured in Chapter~\ref{ch:conjecture}. It involves the concept of codegree decompositions. 

\begin{definition}
A point configuration $\vv A$ admits a \defn{codegree decomposition}\index{codegree decomposition} if there exist $m$ disjoint subsets $\vv A_1,\dots,\vv A_m$ of~$\vv A$ such that for each $1\leq i\leq m$, $\vv A \setminus \vv A_i$ is the set of points of a proper face of $\conv (\vv A)$, and such that $\codegc(\vv A)=\sum_{i=1}^m\codegc(\vv A/(\vv A \setminus \vv A_i))$. The sets $\vv A_i$ are called the \defn{factors}\index{codegree decomposition!factor} of this decomposition, and $m$ its \defn{length}\index{codegree decomposition!length}. 
\end{definition}

Trivially, if $\vv A$ admits a codegree decomposition of length $m$, then $\vv A$ is a weak Cayley configuration of length~$m$. Therefore, the following conjecture is a strengthening of Conjecture~\ref{conj:d+1-2dd} above.

\medskip\noindent
\textbf{Conjecture \ref{conj:strongd+1-2dd}.}
 \emph{Any point configuration $\vv A$ of dimension $d > 2 \degc(\vv A)$ admits a codegree decomposition of length $m\geq d+1-2\degc(\vv A)$.}\medskip

Pyramids are a first example of configurations that have codegree decompositions. Therefore, a consequence of Corollary~\ref{cor:easybound} is that if $n<2\codegc(\vv A)$, then $\vv A$ admits a codegree decomposition of length at least $2\codegc(\vv A)-n$. The following theorem is a next step in this direction.

\medskip\noindent
\textbf{Theorem \ref{thm:2DD}.}
\emph{
Let $\vv A$ be a $d$-dimensional configuration of $n=r+d+1$ points. If $\vv A$ is not a pyramid, then $\vv A$ admits a codegree decomposition of length at least $2(d+1-2\degc(\vv A))-r$.}\\

All the special cases of Conjecture~\ref{conj:strongd+1-2dd} that we have proved still involve the number of points of the configuration.
For example, we know that the conjecture is true when $r\leq 4$ (Corollary~\ref{cor:rank4}).
A particularly interesting case is the following result, that provides a new characterization of \defn{Lawrence polytopes} in terms of the degree. 

\medskip\noindent
\textbf{Proposition \ref{prop:Lawrence}.}
\emph{
Let $\vv A$ be a configuration of $n$ points. If $\vv A$ is not a pyramid then $\codegc(\vv A)=\frac{n}{2}$ if and only if $\vv A$ is a Lawrence configuration. 
In particular, if $\vv A$ is not a pyramid and $\codegc(\vv A)=\frac{n}{2}$, then $\vv A$ admits a codegree decomposition of length $\frac{n}{2}=d+1-\degc(\vv A)$.}\\

In Chapter~\ref{ch:conjecture} this result is reformulated in terms of the \defn{covector discrepancy} of $\Gale{\vv A}$ (Corollary~\ref{prop:Lawrence}), avoiding the distinction among pyramids and non-pyramids. This concept, which we define in Section~\ref{sec:evidences}, also allows to prove more special cases of the conjecture. For example, Theorem~\ref{thm:DD=2} implies that the conjecture holds when $\vv A$ is not a pyramid and has at most $2\codeg(\vv A)+2$ points.

To finish this thesis, in Section~\ref{sec:consequences} we show how Conjecture~\ref{conj:strongd+1-2dd}, if true, would imply all the results that we have presented above.
\medskip

As in Part~\ref{partI}, Gale duality is an essential tool for our proofs, and all these results are derived using the corresponding Gale dual interpretation of the degree, weak Cayley configurations and codegree decompositions.

Moreover, although we only discuss the degree of point configurations, it can be defined for any oriented matroid and most of the results and proofs presented here can be directly translated in terms of oriented matroids. However, since in this part we do not deal with non-realizable matroids, we restrict our statements and proofs to point/vector configurations for the sake of clarity.

\section{Related concepts}

\subsection{The Generalized Lower Bound Theorem}
\label{sec:triang}

Let $\cT$ be a $(d-1)$-dimensional simplicial complex, and $\bm f(\cT)$ be its $f$-vector. That is, $f_i(\cT)$ denotes the number of $i$-dimensional faces of $\cT$. Then the \defn{$h$-vector}\index{$h$-vector} of $\cT$, $\bm h(\cT)=(h_0(\cT),\dots,h_d(\cT))$ is defined by the polynomial relation \[    \sum_{i=0}^d h_i(\cT)\, t^i= \sum_{i=0}^ d f_{i-1}(\cT)\, t^i \, (1-t)^{d-i}.\]
This polynomial is called the \defn{$h$-polynomial} $h_\cT(t)$ of $\cT$.
 
By the famous $g$-theorem \cite{BilleraLee81,Stanley1980}, $h$-polynomials of the boundary complex of simplicial $d$-polytopes ${\vv P}$ are completely known. In particular, $h_{\partial {\vv P}}(t)$ has degree $d$, it satisfies the Dehn-Sommerville equations $h_i(\partial {\vv P}) = h_{d-i}(\partial {\vv P})$, and it is unimodal (\ie $h_{i}(\partial {\vv P}) \ge h_{i-1}(\partial {\vv P})$ for all $1\le i \le \lfloor d/2\rfloor$). In 1971, McMullen and Walkup~\cite{McMullenWalkup1971} posed the following famous conjecture regarding its unimodality, which is now known as the Generalized Lower Bound Theorem:

\begin{theorem}[Generalized Lower Bound]\label{thm:glbt}
Let ${\vv P}$ be a simplicial $d$-polytope. For $i \in \{1, \ldots, \lfloor d/2\rfloor\}$, 
\begin{enumerate}[(i)]
 \item $h_{i}(\partial {\vv P}) \geq h_{i-1}(\partial {\vv P})$; and
 \item $h_{i}(\partial {\vv P}) = h_{i-1}(\partial {\vv P})$ if and only if ${\vv P}$ can be triangulated without interior faces of dimension $\leq d-i$.
\end{enumerate}
\end{theorem}

The first part of the conjecture was solved by Stanley in 1980, as a part of the proof of the $g$-theorem~\cite{Stanley1980}. The second part of the conjecture had remained open until very recently, when it was proved by Murai and Nevo~\cite{MuraiNevo2012}. 

It is instructive to reformulate the previous theorem. For this, let us consider a triangulation $\cT$ of an arbitrary $d$-polytope ${\vv P}$. 
An \defn{interior face} of $\cT$ is a face of $\cT$ that is not contained in a facet of ${\vv P}$. In this situation, the degree of the $h$-polynomial of $\cT$ is well-known, see \cite[Prop.~2.4]{McM04} or \cite[Corollary~2.6.12]{DeLoeraRambauSantosBOOK}. 

\begin{proposition}
Let $\cT$ be a triangulation of a polytope. Then $\degc(h_\cT(t))$ equals the maximal codimension of an interior face of $\cT$.
\end{proposition}

Considering again a simplicial $d$-polytope ${\vv P}$, one defines $g_0(\partial {\vv P}) := 1$, and $g_i(\partial {\vv P}) = h_{i}(\partial {\vv P})-h_{i-1}(\partial {\vv P})$ for $i=1, \ldots, \lfloor \frac{d}{2} \rfloor$. They form the coefficients of the so-called \defn{\text{$g$-polynomial}} $g_{\partial {\vv P}}(t)$. Therefore, Theorem~\ref{thm:glbt} yields for a simplicial polytope ${\vv P}$ that
\[\degc(g_{\partial {\vv P}}(t)) = \min \left\{\degc(h_\cT(t)) \;:\; \cT \text{ triangulation of }{\vv P}\right\}.\]
In other words, the degree $s$ of the $g$-polynomial of a simplicial polytope certifies the existence of {\em some} triangulation that avoids interior faces of dimension $\le d-1-s$. Equivalently, the simplicial polytope ${\vv P}$ is called \defn{$s$-stacked}\index{polytope!stacked} \cite{McMullenWalkup1971}.

For general polytopes, it is also possible to define (toric) $h$- and $g$-polynomials \cite{Stanley1987}. In this case, by \cite{Stanley1992} any rational polytope ${\vv P}$ (conjecturally any polytope) satisfies 

\[\degc(g_{\partial {\vv P}}(t)) \le \min \left\{\degc(h_\cT(t)) \;:\; \cT \text{ triangulation of }{\vv P}\right\}.\]
It is known that simplices are the only polytopes for which $\degc(g_{\partial {\vv P}}(t))=0$. Note that the previous inequality may not be an equality. For instance, a $3$-polytope $P$ which is a prism over a pentagon satisfies $\degc(h_\cT(t)) = 2$ for any triangulation, while $\degc(g_{\partial {\vv P}}(t)) = 1$. In this general situation, it is a hard, open problem to classify all polytopes with $\degc(g_{\partial {\vv P}}(t))=1$ (these polytopes are called \defn{elementary}, see Section~4.3 in \cite{Kalai1994}).

To describe how our results fit into this framework, let us consider the degree of the vertex set $\verts({\vv P})$ of a $d$-polytope ${\vv P}$. By observing that any interior simplex $\vv S$ of $\verts({\vv P})$ can be extended to a triangulation that uses~$\vv S$ as a face, we see that $\degc(\verts({\vv P}))$ is the maximal codimension of an interior simplex of some triangulation of~${\vv P}$. In other words,
\[\degc(\verts({\vv P})) = \max \left\{\degc(h_\cT(t)) \;:\; \cT \text{ triangulation of }{\vv P}\right\}.\]
Hence, classifying polytopes of degree $\delta$ is equivalent to studying polytopes where {\em all} triangulations avoid interior faces of dimension $\le d-1-\delta$. This problem is more tractable than the one described above, and Theorem~\ref{thm:dd=1}  solves it for $\delta=1$.

Finally, a particular motivation for the study of point configurations of degree $1$ (the main characters of Chapter~\ref{ch:deg1}) comes from the {Lower Bound Theorem} for balls (see~\cite[Theorem 2.6.1]{DeLoeraRambauSantosBOOK}). It states that if $\vv A$ is a $d$-dimensional configuration of $n$ points, then any triangulation using all the points in $\vv A$ uses at least $(n-d)$ full-dimensional simplices, and equality is achieved if and only if every $(d-2)$-face of the triangulation lies on the boundary of $\conv(\vv A)$. Hence, $\degc(\vv A)=1$ holds precisely when \emph{all} triangulations using all the points of~$\vv A$ have size~$(n-d)$. This reflects the fact that all triangulations of $\vv A$ are stacked.

This interpretation of Theorem~\ref{thm:dd=1} is already being used by B\"or\"oczky, Santos and Serra in \cite{BoroczkySantosSerra2013} to derive results in additive combinatorics.

\subsection{Tverberg's Theory}

Let $\vv A$ be a configuration of $n$ points in $\RR^r$.
We say that $\vv x\in \RR ^r$ is in the \defn{$\kk$-core}\index{core} of $\vv A$, denoted by $\Ce_\kk(\vv A)$, if every closed halfspace containing $\vv x$ also contains at least $\kk$ points of $\vv A$. We say that a point $\vv x\in \RR ^r$ is an \defn{$m$-divisible}\index{divisible} point of $\vv A$ if there exist $m$ disjoint non-empty subsets $\vv S_1,\dots,\vv S_m$ of $\vv A$ such that $\vv x\in \conv (\vv S_i)$ for $i=1,\dots m$. We denote by $\Tv_m(\vv A)$ the set of $m$-divisible points of $\vv A$. The well-known Tverberg's Theorem asserts that $\Tv_m(\vv A)\neq \emptyset$ whenever $n\geq (m-1)(r+1)+1$. A good introduction for these concepts can be found in~\cite[Chapter 8]{Matousek2002}.

\iftoggle{bwprint}{%
\begin{figure}[htpb]
	\centering
	\subbottom[][$\Ce_2(\vv A)$]{
	\includegraphics[width=.3\textwidth]{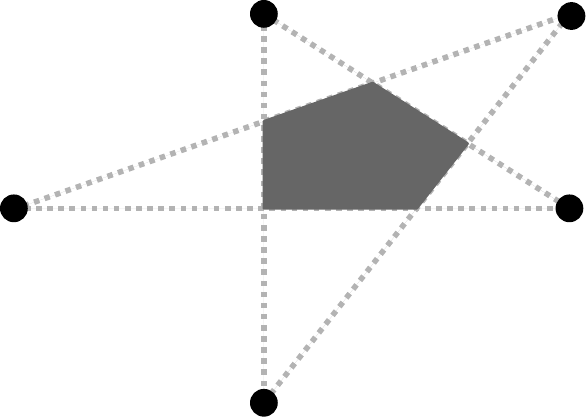}\label{subfig:Core}}
	\qquad\qquad\quad
	\subbottom[][$\Tv_2(\vv A)$]{
	\includegraphics[width=.3\textwidth]{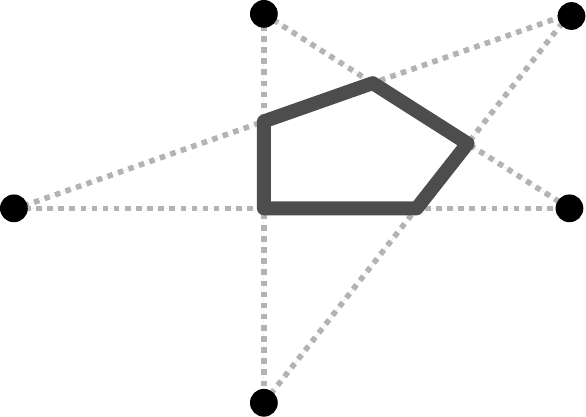}\label{subfig:Divisible}}
\caption[Example of $\kk$-core and $m$-divisible points.]{When $\vv A$ is the vertex set of a pentagon, $\Ce_2(\vv A)$ is the inner pentagon delimited by the interior diagonals, while $\Tv_2(\vv A)$ is only the boundary of this inner pentagon.}\label{fig:CoreAndDivisible}
\end{figure}
}{%
\begin{figure}[htpb]
	\centering
	\subbottom[][$\Ce_2(\vv A)$]{
	\includegraphics[width=.3\textwidth]{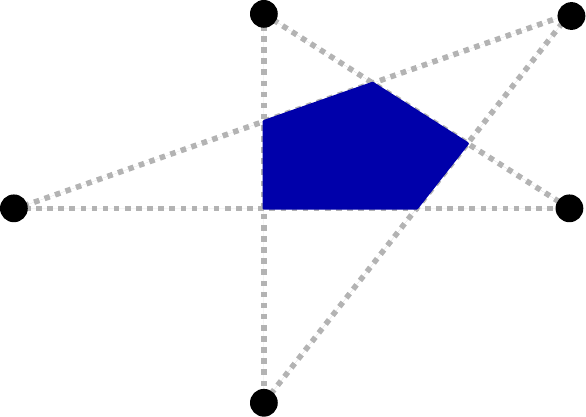}\label{subfig:Core}}
	\qquad\qquad\quad
	\subbottom[][$\Tv_2(\vv A)$]{
	\includegraphics[width=.3\textwidth]{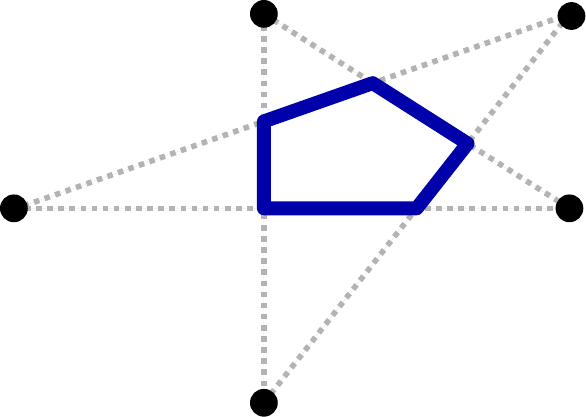}\label{subfig:Divisible}}
\caption[Example of $\kk$-core and $m$-divisible points.]{When $\vv A$ is the vertex set of a pentagon, $\Ce_2(\vv A)$ is the inner pentagon delimited by the interior diagonals, while $\Tv_2(\vv A)$ is only the boundary of this inner pentagon.}\label{fig:CoreAndDivisible}
\end{figure}
}

It is easy to see that $\conv(\Tv_\kk(\vv A))\subset \Ce_\kk(\vv A)$. Equality was conjectured \cite{Reay1982,Sierksma1982}, and actually holds when $r=2$ or $\kk=1$. However, Avis found a counterexample for $n=9$, $r=3$ and $\kk=3$~\cite{Avis1993}, and Onn provided a systematic construction for counterexamples~\cite{Onn2001}.\\

The proof of Theorem~\ref{thm:d+1-3dd} directly yields the following result (see Observations~\ref{obs:depthIsCodeg} and~\ref{obs:divisibleIsCayley}):

\begin{corollary}\label{cor:CoreIsDivisible}
$\Ce_\kk(\vv A)\subset \Tv_{3\kk-2(n-r)}(\vv A)$.
\end{corollary}

Note that this result is only non-trivial if $\kk > \frac23(n-r)$. On the other hand, $C_\kk(\vv A) \not= \emptyset$ implies $\kk \le n-r$. 
Hence, Corollary~\ref{cor:CoreIsDivisible} is of interest for configurations that admit points in some $\kk$-core with a relatively large $\kk$.
In this context, Conjecture~\ref{conj:d+1-2dd}, is equivalent to $\Ce_\kk(\vv A)\subset \Tv_{2\kk-(n-r)}(\vv A)$.

Analogously, Conjecture~\ref{conj:strongd+1-2dd} states that for every $\vv x\in \Ce_\kk(\vv A)$, there are $m$ disjoint (non-empty) subsets $\vv S_1,\dots,\vv S_m$ of $\vv A$, with $m\geq 2\kk-(n-r)$, such that $\vv x\in \Ce_{\kk_i}(\vv S_i)$, $\kk_i\geq 1$ and $\sum_{i=1}^m \kk_i=\kk$. Among some other cases, we can prove this conjecture for all $\vv A\in\RR^r$ with $r\leq4$ (see Corollary~\ref{cor:rank4}). Theorem~\ref{thm:2DD} provides a weaker bound for the number of subsets which, reformulated, states the following:
\begin{theorem}
For any $r$-dimensional configuration of $n$ points $\vv A$, if $\vv x\in \Ce_\kk(\vv A)$ then there are $m=4\kk+r+1-2n$ disjoint (non-empty) subsets $\vv S_1,\dots,\vv S_m$ of $\vv A$ such that $\vv x\in \Ce_{\kk_i}(\vv S_i)$, $\kk_i\geq 1$ and $\sum_{i=1}^m \kk_i=\kk$. 
\end{theorem}

\section{The relation to Ehrhart theory}

\subsection{The lattice degree of a lattice polytope}

Let us consider the situation where ${\vv P} \subset \RR^d$ is a \defn{lattice polytope}\index{lattice polytope}, \ie its vertices are in the lattice $\ZZ^d$. 
We will identify lattice polytopes up to unimodular equivalence, \ie affine isomorphisms of the lattice. 
By a famous result of Ehrhart \cite{Ehr77}, the generating function that enumerates the number of lattice points in 
multiples of ${\vv P}$ is a rational function of the following form:
\[\sum\limits_{k \geq 0} \card{(k {\vv P}) \cap \ZZ^d} \, t^k =  \frac{h^*_{\vv P}(t)}{(1-t)^{d+1}},\]
where the polynomial $h^*_{\vv P}(t) = \sum_{i=0}^d h^*_i t^i$ is called the \defn{$h^*$-polynomial}\index{$h^*$-polynomial} (or $\dd$-polynomial) 
of~${\vv P}$ (see \cite{Bat06,BN07,Hib92,Stanley1993}). Stanley \cite{Stanley1980a,Stanley1986} showed that the coefficients of $h^*_{\vv P}$ are non-negative 
integers. Ehrhart theory can be understood as the study of these coefficients. 

The degree of $h^*_{\vv P}(t)$, i.e., the maximal $i \in \{0, \ldots, d\}$ with $h^*_i\not= 0$, is called the \defn{(lattice) degree}\index{lattice!degree} $\degZ({\vv P})$ of ${\vv P}$ 
\cite{BN07}. The \defn{(lattice) codegree}\index{lattice!codegree} of ${\vv P}$ is defined as $\codegZ({\vv P}) := d+1-\degZ({\vv P})$ and equals the minimal positive integer $k$ 
such that $k {\vv P}$ contains interior lattice points. In recent years these notions and their (algebro-)geometric interpretations 
have been intensively studied \cite{Bat06,BN07,BN08,Nil08,HNP09,DN10,DiRHNP11}. 

The notion of the degree of a lattice polytope was defined in \cite{BN07}, where it was noted that $\degZ({\vv P})$ 
should be considered to be the ``lattice dimension'' of~${\vv P}$. 
This interpretation of the degree was motivated by the following three basic properties:
\begin{enumerate}[(i)]
\item $\degZ({\vv P}) = 0$ if and only if ${\vv P}$ is unimodularly equivalent to the \defn{unimodular simplex} $\conv(\veczero,\vv e_1, \ldots, \vv e_d)$.
\item For a lattice polytope $\vv Q \subset {\vv P}$, we have $\degZ(\vv Q) \le \degZ({\vv P})$ by Stanley's monotoni\-city theorem \cite{Stanley1993}.
\item If ${\vv P}$ is a \defn{lattice pyramid} over $\vv Q$ (i.e., ${\vv P} \cong_\ZZ \conv(\veczero,\vv Q\times\{1\}) \subset \RR^{d+1}$), then 
$\degZ({\vv P}) = \degZ(\vv Q)$.
\end{enumerate}

\subsection{Comparing the degree to the lattice degree}

It was already noted in \cite[Prop.~1.6]{BN07} that a lattice $d$-polytope ${\vv P}$ satisfies \begin{equation}
\label{eq:degcleqdegZ}\degc({\vv P}\cap\ZZ^d) \le \degZ({\vv P}).\end{equation}
Indeed, if $\vv S\subset({\vv P}\cap\ZZ^d)$ is an interior face of $({\vv P}\cap\ZZ^d)$, then $\sum_{\vv a\in\vv S} \vv a$ is an interior point of $k{\vv P}$, where $k=|\vv S|$.

If ${\vv P}$ is a \defn{normal} lattice polytope\index{lattice polytope!normal} (i.e., any lattice point in $k {\vv P}$ is the sum of $k$ lattice points in ${\vv P}$), then $\degc({\vv P}\cap \ZZ^d) = \degZ({\vv P})$. 
However, \eqref{eq:degcleqdegZ} is not an equality in general, as the following example in $3$-space shows: ${\vv P} = \conv(\veczero,\vv e_1,\vv e_2, \vv e_1+\vv e_2+2\vv e_3)$. 
This instance of a {Reeve simplex} \cite{Ree57} satisfies $\degc({\vv P}\cap\ZZ)=0$, but $\degZ({\vv P}) = 2$.

There are also examples of non-normal polytopes where~\eqref{eq:degcleqdegZ} is an equality. The point configuration \[\vv A\!=\!\{(0,0,0,0)\!,(-4, -5, -2,-4)\!,(1, 0, 0,1)\!, (0, 0 ,1,1)\!,(0, 0, 0,1)\!,(3,5,1,1)\}\] is the set of lattice points in a non-normal lattice polytope ${\vv P}=\conv(\vv A)$ with $\degc({\vv P}) = \degZ({\vv P})=3$. This example is due to Aaron Dall (personal communication).
\medskip

In the setting of Section~\ref{sec:triang}, Equation~\eqref{eq:degcleqdegZ} can also be deduced directly from a stronger result by Betke and McMullen \cite{Betke-McM}: 
$h_\cT(t) \le h^*_{\vv P}(t)$, coefficientwise, for any triangulation $\cT$ of ${\vv P\cap \ZZ}$.

From our point of view, the degree may be seen as a natural generalization of the Ehrhart-theoretic lattice degree. 
In particular, all properties of the lattice degree mentioned above also hold in the setting of point configurations $\vv A \subset \RR^d$:

\begin{enumerate}[(i)]
\item $\degc(\vv A) = 0$ if and only if $\vv A$ is the vertex set of a $d$-simplex (Proposition~\ref{prop:dd=0}).
\item For $\vv A' \subset \vv A$, we have $\degc(\vv A') \le \degc(\vv A)$ (Corollary~\ref{cor:degpointdeletioncontraction}).
\item If $\vv A$ is a \defn{pyramid} over $\vv A'$
, then $\degc(\vv A) = \degc(\vv A')$ (Lemma~\ref{lem:pyr}).
\end{enumerate}

\subsection{Cayley configurations}\label{sec:PrimalCayley}

Our main results are motivated by analogous statements in Ehrhart theory. 
In particular, the notion of a weak Cayley configuration originates in the widely used construction of Cayley polytopes. 
Cayley polytopes play a very important role in the study of the degree of lattice polytopes~\cite{BN07,HNP09} and, more generaly, in the study of mixed subdivisions of Minkowski sums via the Cayley trick (cf. Observation~\ref{obs:Cayleytrick}). Let us carefully state some natural generalizations:

\begin{definitions}
 Let $\vv A$ be a point configuration in $\RR^d$. We say that $\vv A$ is 
\begin{itemize}
 \item a \defn{lattice Cayley configuration}\index{Cayley configuration!lattice} of length $m$, if $\vv A\subset \ZZ^d$ and $\vv A$ maps onto $\conv(\veczero,\vv e_1, \ldots, \vv e_{m-1})$ via a lattice projection $\ZZ^d \to \ZZ^{m-1}$.

 \item an \defn{affine Cayley configuration}\index{Cayley configuration!affine} of length $m$, if $\vv A$ maps onto the vertex set of a $(m-1)$-simplex via an affine projection $\RR^d \to \RR^{m-1}$.

 \item a \defn{combinatorial Cayley configuration}\index{Cayley configuration!combinatorial} of length $m$, 
if there exists a partition $\vv A = \vv A_1 \uplus \cdots \uplus \vv A_m$, such that for any $1\leq i\leq m$, $\vv A\setminus \vv A_i$ is the set of points of a proper face of~$\conv (\vv A)$.
\end{itemize}
The sets $\vv A_1,\dots, \vv A_m$ are called the \defn{factors} of the configuration. They must be non-empty, because $\vv A_i=\emptyset$ implies that $\conv(\vv A\setminus\vv A_i)$ is not a proper face of $\conv(\vv A)$. Moreover, observe that this definition implies that for any $\emptyset \not= I \subsetneq \{1, \ldots, m\}$, \(\conv \left( \bigcup\nolimits_{i \in I} \vv A_i\right)\) is a proper face of~$\conv(\vv A)$.\\

We say that a polytope ${\vv P} \subset \RR^d$ is a lattice, affine or combinatorial \defn{Cayley polytope}, if its vertex set is a Cayley configuration of the respective type.
\end{definitions}

Obviously, ``lattice'' implies ``affine'' implies ``combinatorial''. Of course, there are affine Cayley configurations that are not lattice, and there are combinatorial Cayley configurations that are not affine (\eg the vertices of a deformed prism in $\RR^3$).\\

Let us point out that the term ``combinatorial Cayley configuration'' is not ambiguous, but indeed means \emph{combinatorially equivalent to an affine Cayley configuration}, as the following result (proved in Chapter~\ref{ch:cayley}) shows:

\medskip\noindent
\textbf{Proposition \ref{prop:combinatorialisaffine}.}
Every combinatorial Cayley configuration of length~$m$ is combinatorially equivalent to an affine Cayley configuration of length~$m$.

\medskip

\begin{observation}[{The Cayley trick}]\label{obs:Cayleytrick}
 The origin of the definition of Cayley polytopes\index{Cayley polytope} comes from the Cayley trick. If $\vv A$ and $\vv B$ are two point configurations in $\RR^d$, their Cayley embedding is the $(d+1)$-dimensional point configuration $\Cayley(\vv A,\vv B)=\set{(\vv a,1,0)}{\vv a\in \vv A}\cup\set{(\vv b,0,1)}{\vv b\in \vv B}$. The \defn{Cayley trick} states that there is a bijection between triangulations of $\Cayley(\vv A,\vv B)$ and mixed subdivisions of the Minkowski sum $\vv A + \vv B$ (see~\cite[Section 9.2]{DeLoeraRambauSantosBOOK} and references therein). The fact that Cayley polytopes have small degree reinforces the idea that the degree should be considered as the true dimension of a polytope: the triangulations of a Cayley polytope are in bijection with mixed subdivisions of a lower dimensional polytope.
\end{observation}

\subsection{Motivating results}\label{sec:mainmotivation}

Proposition~\ref{prop:dd=0}, Theorem~\ref{thm:dd=1}, Theorem~\ref{thm:d+1-3dd} and Corollary~\ref{cor:easybound}
 should be seen as a combinatorial generalization of known results in the geometry of lattice polytopes.
Here are the original formulations of these statements in the context of lattice polytopes. 
In the following let ${\vv P}$ be a $d$-dimensional lattice polytope with $r + d + 1$ vertices and lattice degree $\degZ({\vv P})=s$.

\begin{enumerate}[(i)]
\item As noted above, ${\vv P}$ has lattice degree $s=0$ if and only if ${\vv P}$ is a unimodular simplex.
\item Lattice $d$-polytopes ${\vv P}$ of lattice degree $s=1$ were classified in \cite{BN07}: either ${\vv P}$
is a $(d-2)$-fold lattice pyramid over the lattice triangle $\conv((0,0),(2,0),(0,2))$, or ${\vv P}$ is a lattice Cayley polytope of length $d$.
\item\label{it:latticepyr} The following result was shown in \cite{Nil08}: If 
\[d > r(2s + 1) + 4s - 2,\]
then ${\vv P}$ is a lattice pyramid over an $(d-1)$-dimensional lattice polytope.
\item\label{it:latticeCayley} And in~\cite{HNP09}: If $d >f(s) := (s^2+19s-4)/2$, then ${\vv P}$ is a lattice Cayley polytope of length $d+1-f(s)$.
\end{enumerate}

The reader is invited to compare these results with the combinatorial statements for arbitrary polytopes or point configurations in Proposition~\ref{prop:dd=0}, Theorem~\ref{thm:dd=1}, Theorem~\ref{thm:d+1-3dd} and Corollary~\ref{cor:easybound}. As is to be expected, the assumptions in combinatorial setting are more general, while the conclusions are weaker. Nevertheless, the bounds in Corollary~\ref{cor:easybound} and Theorem~\ref{thm:d+1-3dd} are better than their lattice analogues. 

It was noted in \cite{BN07} that lattice Cayley polytopes of length $m$ have lattice degree at most $d+1-m$. In particular, lattice Cayley polytopes of 
large length have small lattice degree. It was asked in \cite{BN07} whether there might be a converse to this, and the  above statement \eqref{it:latticeCayley} answered this question affirmatively. The assumption in \eqref{it:latticeCayley} is surely not sharp, and it is conjectured that $f(s)=2s$ should suffice, see \cite{DiRHNP11,DN10}. 

\begin{conjecture} \textup{\cite{DiRHNP11,DN10}}\label{conj:latticeconj}
 Let $\vv P$ be a lattice $d$-polytope. If $d >2\degZ(\vv P)$, then ${\vv P}$ is a lattice Cayley polytope of length $d+1-2\degZ(\vv P)$.
\end{conjecture}

Therefore, it seems at first very tempting to also conjecture the analogous statement for combinatorial types of polytopes: Namely, 
for a $d$-dimensional polytope ${\vv P}$, $d > 2\degc(\verts({\vv P}))$ should imply that ${\vv P}$~is a combinatorial Cayley polytope of length $d+1-2\degc(\verts({\vv P}))$. This statement indeed holds for $\degc(\verts({\vv P}))=1$ by Theorem~\ref{thm:dd=1}; however, 
the following example proves the general guess wrong:

\begin{example}\label{ex:liftedexceptionalsimplex}
Consider the $(d+1)$-dimensional point configuration \[\vv A=\{\veczero,2 \vv e_1, \dots , 2 \vv e_d, \vv e_1 + \vv e_{d+1},\vv e_1-\vv e_{d+1}, \dots, \vv e_d +\vv  e_{d+1},\vv e_d-\vv e_{d+1}\}.\] It is in convex position, i.e., $\vv A$ is the vertex set of $\conv(\vv A)$, and has degree~$2$. However, $\vv A$ is not a combinatorial Cayley configuration of length $>1$.

Indeed, let $\vv U=\{\vv e_1,\dots, \vv e_d\}$, $\vv U^+=\vv U+\vv e_d$ and $\vv U^-=\vv U-\vv e_d$. Observe that if $\vv A$ were a Cayley configuration with factors $\vv A_1, \ldots, \vv A_m$, then every subconfiguration $\vv B$ of $\vv A$ would admit a decomposition with factors $\vv B_i=\vv A_i\cap \vv B$ (of course, some of them might be empty). In particular, $\vv A^+=\conv\{2\vv U,\vv U^+\}$ and $\vv A^-=\conv\{2\vv U,\vv U^-\}$ would be Cayley configurations. But each of $\vv A^+$ and $\vv A^-$ is the vertex set of a prism over a simplex, which only admits the two Cayley decompositions $\vv A^\pm=2\vv U\uplus \vv U^\pm$ and $\vv A^\pm=\biguplus_i \{2\vv e_i, \vv e_i\pm\vv e_{d+1}\}$. None of these decompositions of $\vv A^+$ and~$\vv A^-$ can be extended to a Cayley decomposition of $\vv A$ that includes $\{\veczero\}$. The decomposition of $\vv A$ into $\vv A^+$ and $\vv A^-$ for $d=2$ is shown in Figure~\ref{fig:liftedprism}.

\iftoggle{bwprint}{%
\begin{figure}[htpb]
\centering
 \subbottom[$\vv A$]{\label{sfig:liftedprism1}\includegraphics[width=.27\linewidth]{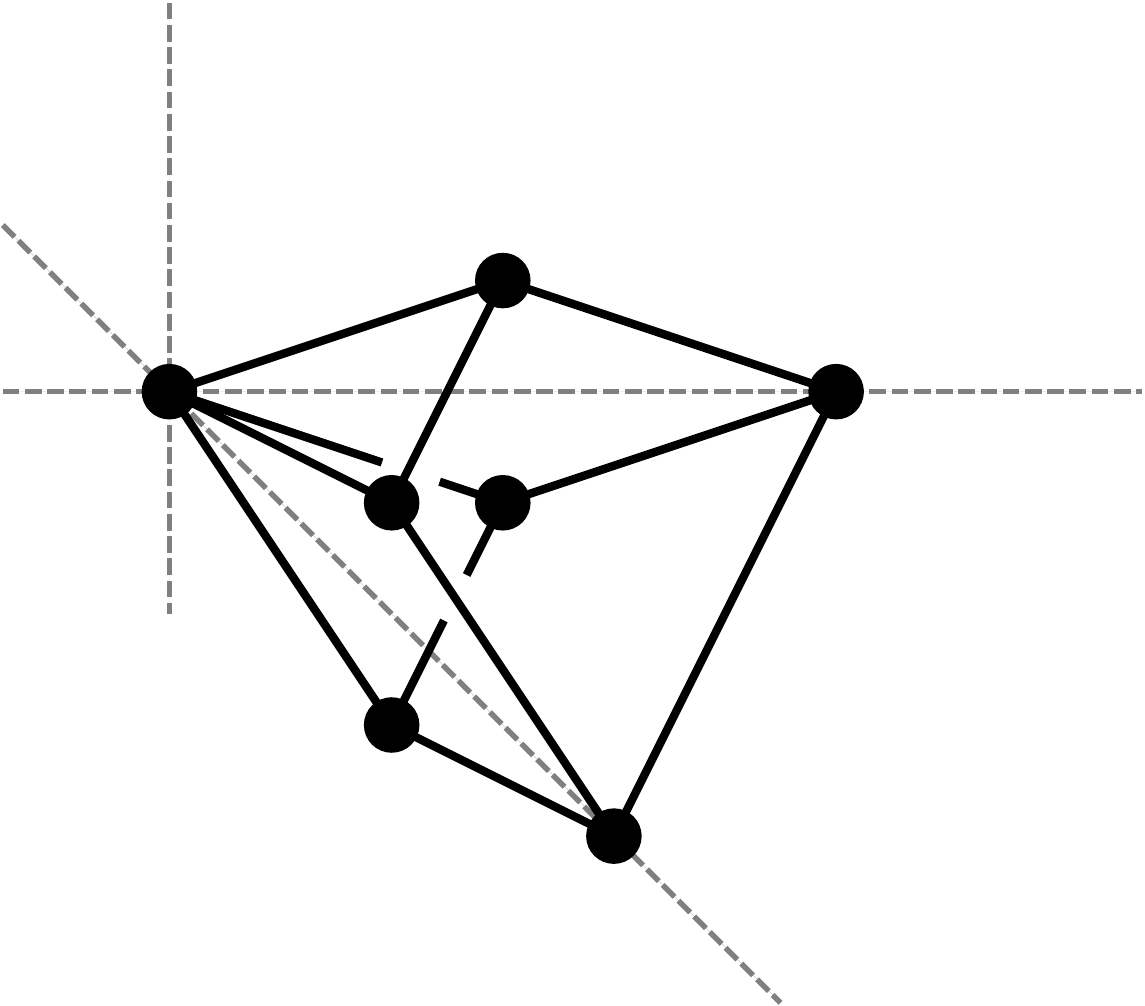}}\qquad
 \subbottom[$\vv A^+$]{\label{sfig:liftedprism2}\includegraphics[width=.27\linewidth]{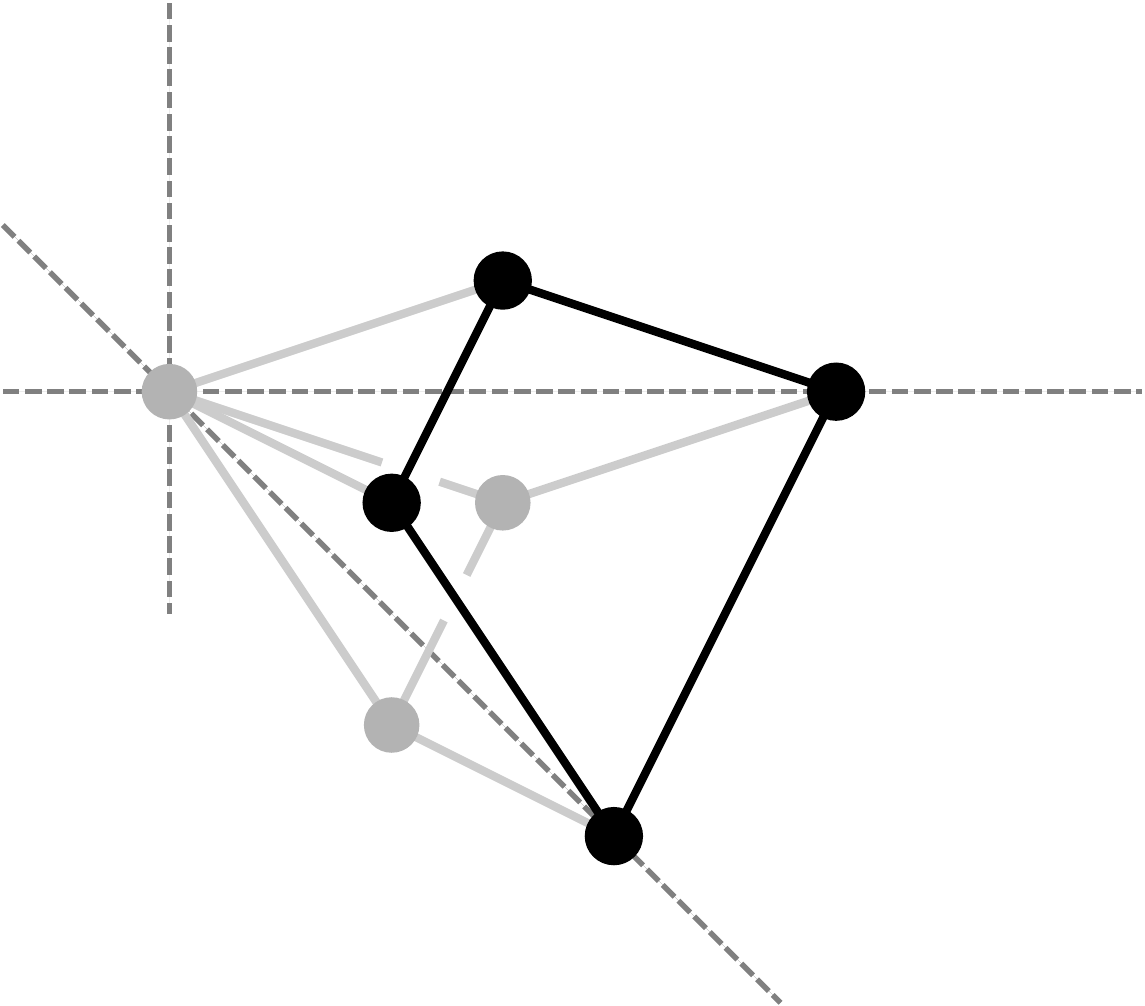}}\qquad
 \subbottom[$\vv A^-$]{\label{sfig:liftedprism3}\includegraphics[width=.27\linewidth]{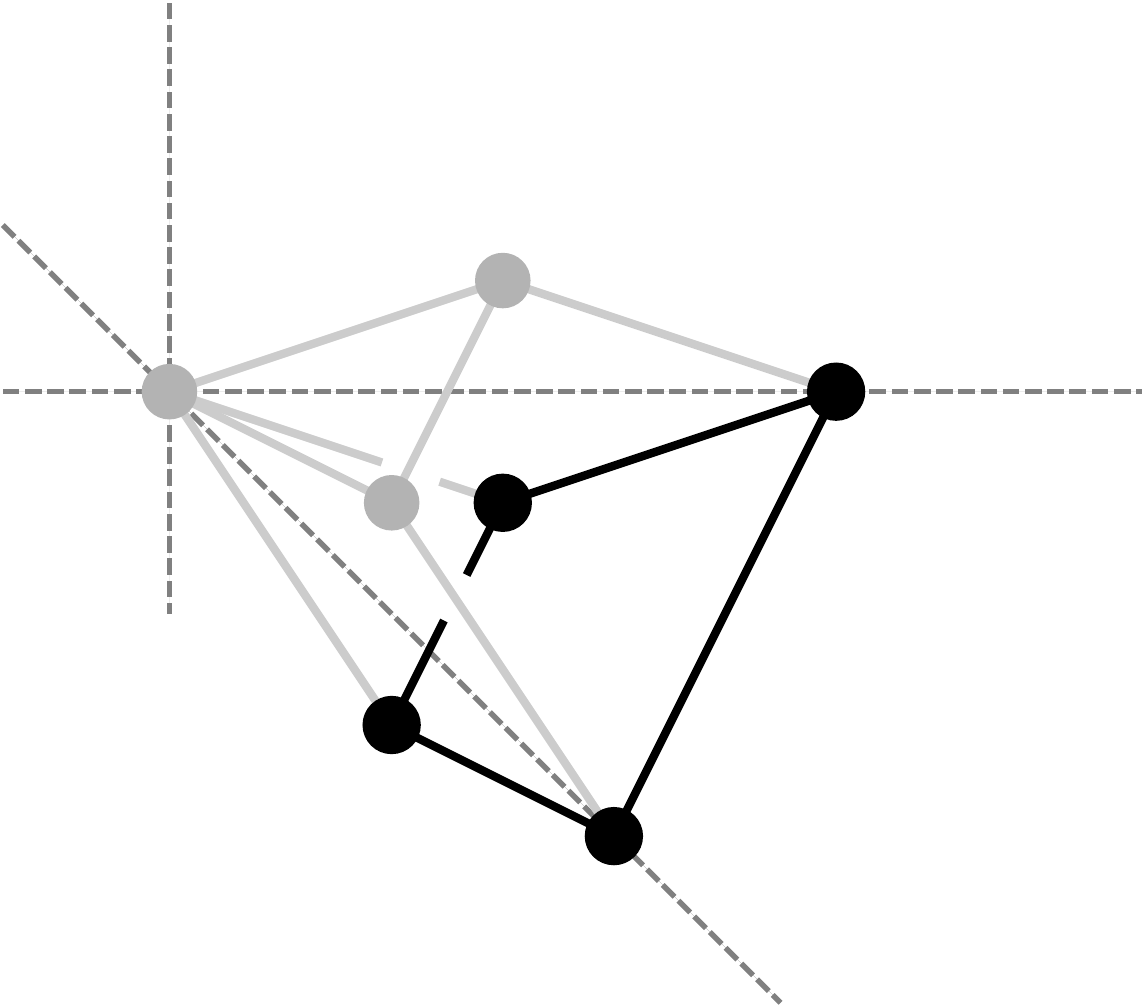}}
 \caption{The point configuration of Example~\ref{ex:liftedexceptionalsimplex}.}
 \label{fig:liftedprism}
\end{figure}
}{%
\begin{figure}[htpb]
\centering
 \subbottom[$\vv A$]{\label{sfig:liftedprism1}\includegraphics[width=.27\linewidth]{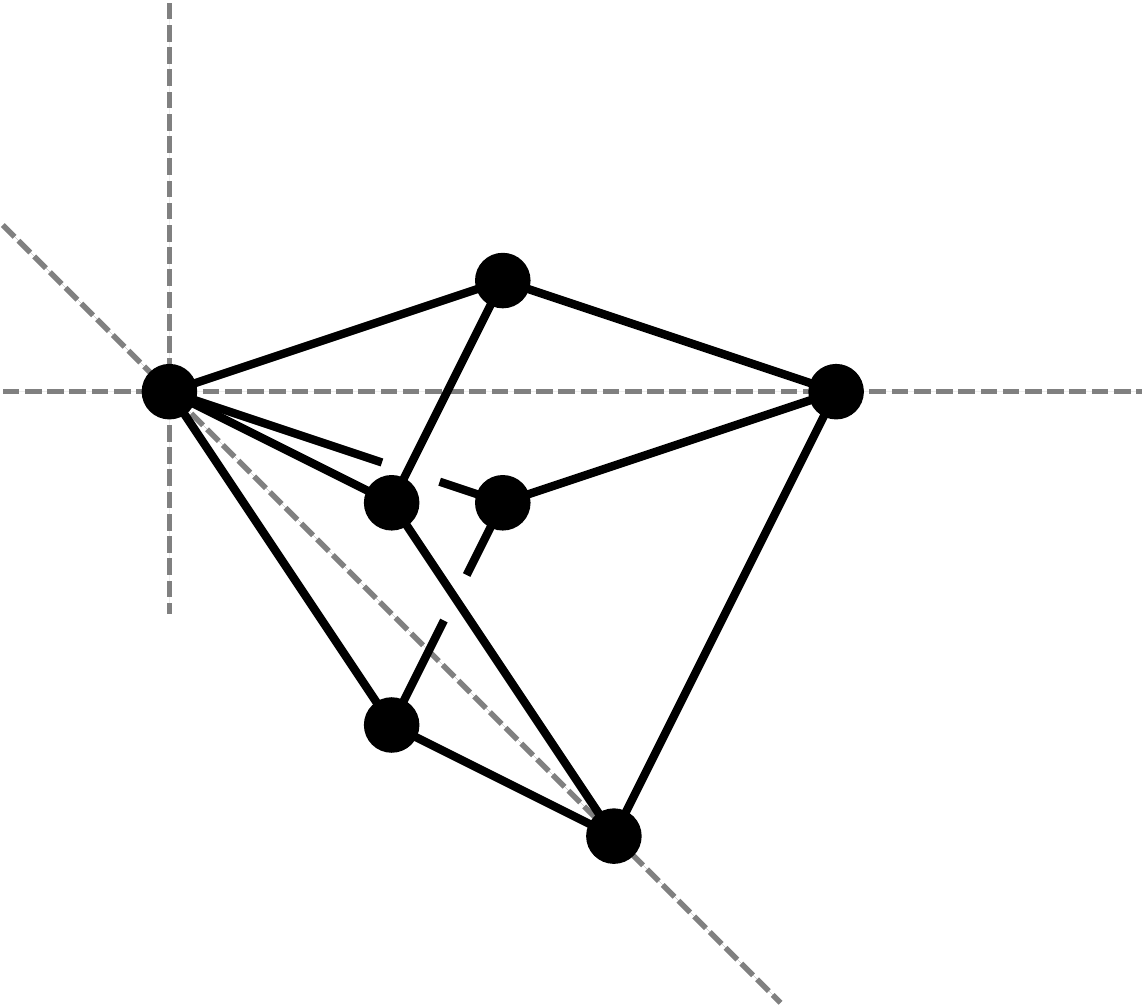}}\qquad
 \subbottom[$\vv A^+$]{\label{sfig:liftedprism2}\includegraphics[width=.27\linewidth]{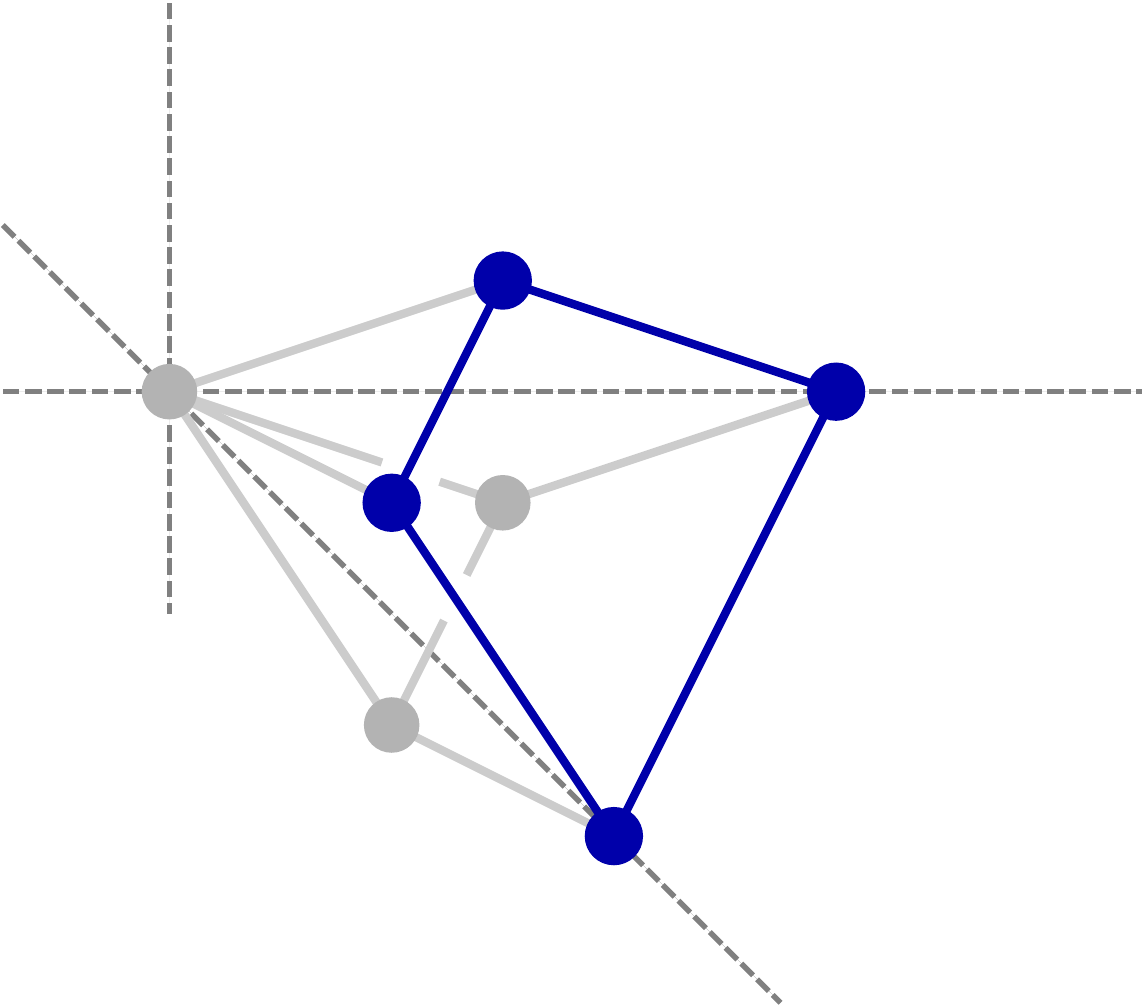}}\qquad
 \subbottom[$\vv A^-$]{\label{sfig:liftedprism3}\includegraphics[width=.27\linewidth]{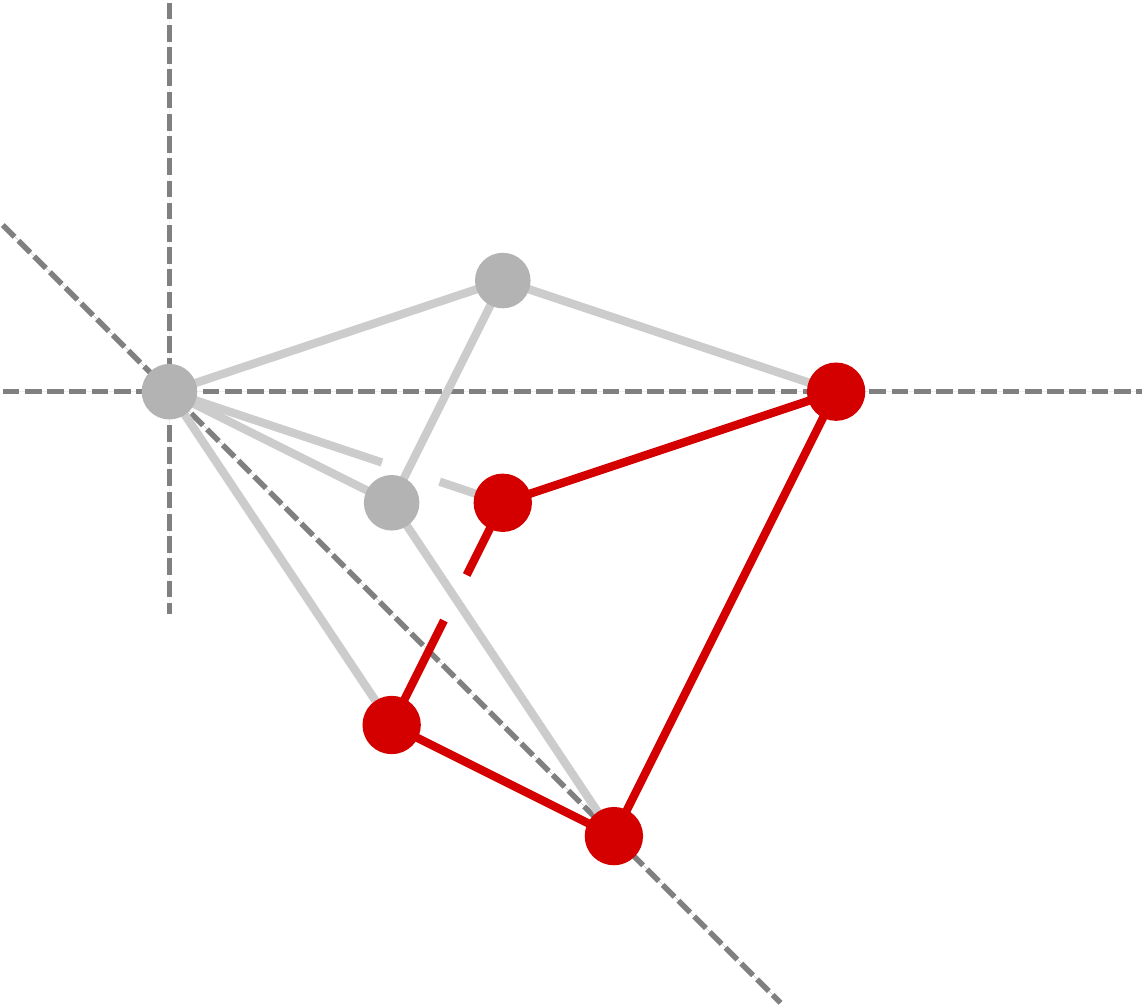}}
 \caption{The point configuration of Example~\ref{ex:liftedexceptionalsimplex}.}
 \label{fig:liftedprism}
\end{figure}
}
\end{example}

\subsection{Weak Cayley configurations}

Even though the point configuration of Example~\ref{ex:liftedexceptionalsimplex} is not a combinatorial Cayley configuration, the subsets $\vv B_i=\{\veczero,2 \vv e_i, \vv e_i + \vv e_{d+1},\vv e_i - \vv e_{d+1}\}$ fulfill all the necessary conditions for this, except for {disjointness}. This motivates our original definition of weak Cayley configuration.
 
\begin{definition}
 A point configuration $\vv A$ is a \defn{weak Cayley configuration} of length $m$, 
if~$\vv A$ can be covered by subsets $\vv A = \vv B_1 \cup \cdots \cup \vv B_m$, such that for any $1\leq i\leq m$, 
$\emptyset \neq I \subsetneq \{1, \ldots, m\}$, \(\bigcup\nolimits_{i \in I} \vv B_i\) is the set of points of a proper face of~$\conv (\vv A)$.
\end{definition}

Setting $\vv A_i=\vv B_i\setminus \bigcup\nolimits_{j\neq i} \vv B_j$ for $1\leq i\leq m$ and $\vv A_0=\vv A\setminus\bigcup\nolimits_{i=1}^d \vv A_i$, it is not hard to prove 
that this definition is indeed equivalent to the stronger Definition~\ref{def:strongweakcayley} stated above in the introduction. We prefer Definition~\ref{def:strongweakcayley} since it is more restrictive and allows for the following formulation:

\begin{observation} 
  $\vv A$ is a weak Cayley configuration of length $m$ if and only if there is a (possibly empty) face $\vv F$ of $\conv(\vv A)$ such that the contraction $\vv A/(\vv F\cap \vv A)$ is a Cayley configuration of length~$m$. Then the factors of~$\vv A$ are defined as the factors of~$\vv A/(\vv F\cap \vv A)$.
\end{observation}

Note that every combinatorial Cayley configuration is a weak Cayley configuration. Example~\ref{ex:liftedexceptionalsimplex} motivates why even for polytopes (instead of more general point configurations) it is necessary to consider weak Cayley configurations. The point configuration in this example is a weak Cayley configuration of length~$d$, with $\vv A_0=\{\veczero\}$ and $\vv A_i=\{2 \vv e_i, \vv e_i + \vv e_{d+1},\vv e_i - \vv e_{d+1}\}$.

Summing up, Theorem~\ref{thm:d+1-3dd} should be seen as the correct combinatorial analogue of the statement \eqref{it:latticeCayley} for lattice polytopes in the previous section. 
Moreover, Conjecture~\ref{conj:latticeconj} for lattice polytopes precisely matches our Conjecture~\ref{conj:d+1-2dd}.

\subsection{Codegree decompositions}

We have presented some results on geometric combinatorics that are inspired in analogue Ehrhart-theoretic results. But this is a two-way path, and we can use our understanding of the degree of point configurations to motivate research on lattice polytopes. In particular, the evidences for Conjecture~\ref{conj:strongd+1-2dd} suggest that the following conjecture might hold.

\begin{conjecture}
 If $\vv P$ is lattice $d$-polytope~$\vv P$ with $d>2\degZ(\vv P)$, then there is a lattice projection that maps $\vv P$ onto a lattice join of lattice polytopes that has the same lattice codegree as $\vv P$.
\end{conjecture}
It is true at least when $\degZ(\vv P)\leq1$ by~\cite[Theorem~2.5]{BN07}.
\chapter{Weak Cayley configurations}\label{ch:cayley}

In this chapter, we use Gale duality to relate the degree to weak Cayley configurations. In particular, we present the proofs of Proposition~\ref{prop:dd=0}, Corollary~\ref{cor:easybound} and Theorem~\ref{thm:d+1-3dd} that were announced in the introduction.

\section{The dual degree}
\label{sec:dualdeg}
 One can mirror the duality between neighborliness and balancedness (see Section~\ref{sec:balvsneigh}) to give a dual interpretation of the degree.

\begin{definition}
Let $\vv V$ be a full-dimensional vector configuration in $\RR^r$. Its \defn{dual degree}\index{degree!dual} is 
\[\degG(\vv V):=\max_{\vvh H} |\vvh H^+\cap \vv V|-r,\]
\index{$\degG(\vv V)$}where $\vvh H$ runs through all linear hyperplanes of $\RR^r$. 

That is, $\degG(\vv V)=\dd$ if and only if $\dd$ is the minimal integer such that for every linear hyperplane $\vvh H$, there are at most $r+\dd$ vectors of $\vv V$ in $\vvh H^+$. 
\end{definition}

This definition is coherent with its primal counterpart:

\begin{proposition}\label{prop:degstarastar}
\(\degc({\vv A})=\degG(\Gale {\vv A}).\)
\end{proposition}
\begin{proof}
Let ${\vv A}$ be a $d$-dimensional configuration of $n$ points. By definition, $\degc({\vv A})=\dd$ if every subset $\vv S$ of ${\vv A}$ of size $d-\dd$ is contained in a supporting hyperplane. Equivalently, if $\vv W$ contains the origin in its convex hull for every $\vv W\subset\Gale{\vv A}$ of size $n-d+\dd=r+\dd+1$ (see Lemma~\ref{lem:gale}).
Therefore, if~$\degc({\vv A})=\dd$ there cannot be a hyperplane $\vvh H$ in $\RR^r$ through the origin that contains more than $r+\dd$ vectors of $\Gale{\vv A}$ in $\vvh H^+$ (by the Farkas Lemma, see \cite[Section 1.4]{Ziegler1995}). This proves that $\degG(\Gale{\vv A})\leq \deg(\vv A)$.
Conversely, if there is a set of $r+\dd$ vectors whose convex hull does not contain the origin, which by Lemma~\ref{lem:gale} means that there is an interior face of $\vv A$ of cardinality $\leq d+1-\dd$, then we can separate this set from the origin by a hyperplane~$\vv H$, again by the Farkas Lemma. This proves that $\degG(\Gale{\vv A})\geq \deg(\vv A)$.
\end{proof}

\begin{definition}\label{def:dualcodegree}
The \defn{dual codegree}\index{codegree!dual} of a vector configuration $\vv V$ is 
\[
   \codegG(\vv V):=\min_{\vvh H} |\ol {\vvh H}^-\cap \vv V|=\min_{\vvh H} |\ol {\vvh H}^+\cap \vv V|,
\]
\index{$\codegG(\vv V)$}where ${\vvh H}$ runs through all linear hyperplanes.
\end{definition}

\begin{remark}
If $|\vv V|=r+d+1$, it is easy to see that \[\codegG(\vv V)=d+1-\degG(\vv V),\] which is consistent with the primal definition (cf. Definition~\ref{def:combdegree}).
\end{remark}

\begin{observation}\label{obs:depthIsCodeg}
These definitions bring the first connection to Tverberg theory mentioned in the introduction. Fix a vector configuration~$\vv V$, and let~$\ol {\vv V}$ be the point configuration consisting of the endpoints of the vectors in $\vv V$. Observe that $\codegG(\vv V) \geq \kk$ if and only if the origin has depth~$\kk $ in~$\ol {\vv V}$, so that $\veczero\in \Ce_\kk(\ol {\vv V})$.
\end{observation}

\begin{example}
The vector configurations of Figure~\ref{fig:dualdegree} have dual degree~$1$, because for each linear hyperplane~$\vvh H$, there are at most $3$~vectors of the configuration in $\vvh H^+$. Their dual codegrees are $3$~and~$2$, respectively.
\end{example}

\iftoggle{bwprint}{%
\begin{figure}[htpb]
\centering
 \subbottom[Gale dual of $\simp{1}\times\simp{2}$.]{\qquad\qquad\label{sfig:dualprism}\includegraphics[width=.25\linewidth]{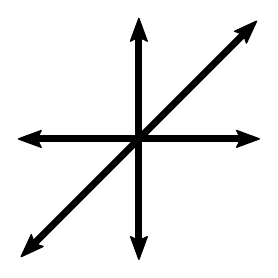}\qquad\qquad}
 \subbottom[Gale dual of a pentagon.]{\qquad\qquad\label{sfig:dualpentagon}\includegraphics[width=.25\linewidth]{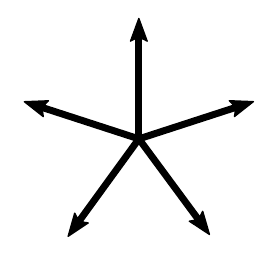}\qquad\qquad}
 \caption{Two vector configurations with dual degree $1$.}
 \label{fig:dualdegree}
\end{figure}
}{%
\iftoggle{print}{%
\begin{figure}[htpb]
\centering
 \subbottom[Gale dual of $\simp{1}\times\simp{2}$.]{\qquad\qquad\label{sfig:dualprism}\includegraphics[width=.25\linewidth]{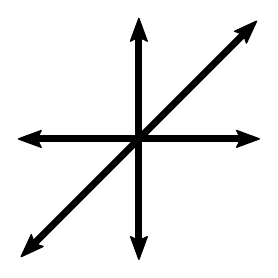}\qquad\qquad}
 \subbottom[Gale dual of a pentagon.]{\qquad\qquad\label{sfig:dualpentagon}\includegraphics[width=.25\linewidth]{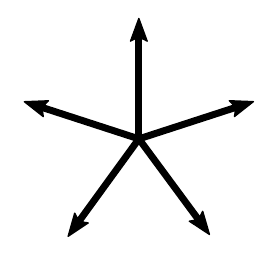}\qquad\qquad}
 \caption{Two vector configurations with dual degree $1$.}
 \label{fig:dualdegree}
\end{figure}
}{%
\begin{figure}[htpb]
\centering
 \subbottom[Gale dual of $\simp{1}\times\simp{2}$.]{\qquad\qquad\label{sfig:dualprism}\includegraphics[width=.25\linewidth]{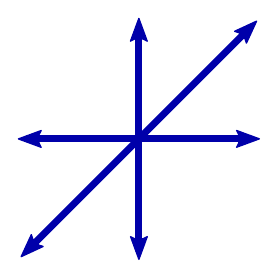}\qquad\qquad}
 \subbottom[Gale dual of a pentagon.]{\qquad\qquad\label{sfig:dualpentagon}\includegraphics[width=.25\linewidth]{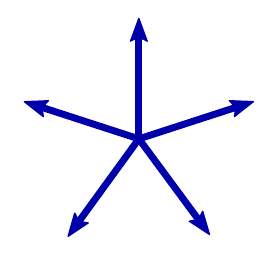}\qquad\qquad}
 \caption{Two vector configurations with dual degree $1$.}
 \label{fig:dualdegree}
\end{figure}
}
}

\begin{example}\label{ex:prism}
 A \defn{prism over a $d$-simplex} is the cartesian product $\simp{1} \times \simp{d}$ of a $1$-simplex $\simp{1}$ with a $d$-simplex $\simp{d}$. Its Gale dual can be represented by the vector configuration $\vv V:=\{\pm \vv e_1,\dots,\pm \vv e_{d}, \pm \sum_{i=1}^d \vv e_i\}\subset\RR^{d}$. 
 Hence, every hyperplane~$\vvh H$ in general position in~$\RR^d$ contains exactly $d+1$ points of $\vv V$ in $\vvh H^+$, and thus its degree is~$1$. Figure~\ref{sfig:dualprism} illustrates the case $d=2$.
\end{example}

The dual degree of a vector configuration (resp., the degree of a point configuration) can only decrease under deletions and contractions.

\begin{proposition}\label{prop:subconfigurations}
For any vector configuration $\vv V$ and for any vector $\vv v\in \vv V$, $\degG(\vv V\setminus \vv v)\leq \degG(\vv V)$ and $\degG(\vv V/\vv v)\leq \degG(\vv V)$.
\end{proposition}
\begin{proof}
The first statement is a consequence of the fact that, for any hyperplane ${\vvh H}$, $|{\vvh H}^+\cap (\vv V\setminus \vv v)|\leq |{\vvh H}^+\cap \vv V|$. This already proves our claim if the deletion of $\vv v$ does not reduce the rank of $\vv V$. On the other hand, if $\rank(\vv V\setminus \vv v)= \rank(\vv V)-1$, the claim also follows easily after observing that $\vv v\in {\vvh H}^+$ for any hyperplane $\vv H$ such that $|{\vvh H}^+\cap \vv V|=\rank(\vv V)+\degG(\vv V)$.

For the second statement we can assume that $\vv v \neq\veczero$. We have to show that no hyperplane ${\vvh H}$ through $\vv v$ can have more than $r+\dd-1$ elements on its positive side, where $r$ is the rank of $\vv V$ and $\dd=\degG(\vv V)$. 
Indeed, suppose that $\degG(\vv V/\vv v)>\degG(\vv V)$, so that ${\vvh H}^+\cap \vv V$ has at least  $(r-1)+(\dd+1)=r+\dd$ elements on its positive side. 
And let ${\vvh H}_{\vv v}$ be the hyperplane whose normal vector is $\vv v$. Then ${\vvh H}\circ {\vvh H}_{\vv v}$ has at least~$\vv v$ and~$\vv V\cap\vvh H$ in its positive side, and this contradicts $\degG(\vv V)=\dd$.
\end{proof}

\begin{corollary}\label{cor:degpointdeletioncontraction}
For any point configuration $\vv A$ and for any point $\vv a\in {\vv A}$, $\degc({\vv A}\setminus \vv a)\leq \degc({\vv A})$ and $\degc({\vv A}/\vv a)\leq \degc({\vv A})$.\qed
\end{corollary}

\subsection{Pure vector configurations}

Corollary~\ref{cor:degpointdeletioncontraction} explains one of the reasons why it is natural to allow configurations that admit repeated points: even if ${\vv A}$ has no repeated points, ${\vv A}/\vv a$ might contain some (see the example of Figure~\ref{fig:affineContractionDeletion}). 
However, it is straightforward to see that deleting repeated points from~${\vv A}$ changes neither the degree nor the property of being a weak Cayley configuration: 

\begin{lemma}\label{lem:removerepeated}
If the point configuration $\vv A'$ is obtained from $\vv A$ after deleting all repeated points, then 
$\degc({\vv A})=\degc({\vv A}')$. Moreover, ${\vv A}$ is a (weak) Cayley configuration of length $m$ if and only if ${\vv A}'$ is. 
\qed
\end{lemma}

For this reason, we usually only consider point configurations without repeated points. Dually, we say that a vector configuration~$\vv V$ is \defn{pure}\index{vector configuration!pure} if~$\Gale{\vv V}$ does not have repeated points. Using the characterization from Lemma~\ref{lem:dualnorepeatedpoints}, this definition can be presented as follows:
\begin{definition}
A vector configuration $\vv V\subset\RR^\rr$ is \defn{pure} if and only if either $r=0$, or for every linear hyperplane~${\vvh H}$, $|{\vvh H}^+\cap \vv V|\geq 2$ or $|{\vvh H}^- \cap \vv V|\geq 2$. 
\end{definition}

It is obvious that if $\vv A$ has no repeated points then $\vv A\setminus \vv a$ neither. Dualizing, this reads:
\begin{lemma}\label{lem:quotientsofpurearepure}
 If $\vv V$ is a pure vector configuration, then $\vv V/\vv v$ is pure for each $\vv v\in\vv V$.\qed
\end{lemma}

A first interesting consequence of this characterization is the following lemma, which will allow us to classify point configurations of degree~$0$.

\begin{lemma}~\label{lem:puredeggeq1}
 If ${\vv V}$ is a pure vector configuration with $\rank(\vv V)\geq 1$, then $\degG({\vv V})\geq 1$.
\end{lemma}
\begin{proof}
Let ${\vvh H}$ be a hyperplane spanned by some subconfiguration $\vv W \subset {\vv V}$. By Lemma~\ref{lem:dualnorepeatedpoints}, we can assume that $|{\vvh H}^+\cap {\vv V}|\geq 2$. Then the contraction~${\vv V}/\vv W$ is a pure configuration of rank $1$ that satisfies $\degG({\vv V}/\vv W)\geq 1$. The result now follows from Proposition~\ref{prop:subconfigurations}.
\end{proof}

\begin{proposition}\label{prop:dd=0}
The degree of a point configuration ${\vv A}$ is~$0$ if and only if ${\vv A}$ is the set of vertices of a simplex (possibly with repetitions).
\end{proposition}
\begin{proof}
Because of Corollary~\ref{cor:degpointdeletioncontraction} and Lemma~\ref{lem:removerepeated}, it is enough to see that there are no $d$-dimensional point configurations of degree~$0$ with $d+2$ points, none of which are repeated; this follows from Lemma~\ref{lem:puredeggeq1}.
\end{proof}

\begin{corollary}\label{cor:deg0}
The dual degree of a vector configuration ${\vv V}$ is~$0$ if and only if $\vv V$ is a direct sum of positive circuits.\qed
\end{corollary}

\subsection{Irreducible vector configurations}
 In this dual setting, some results mentioned in the introduction have a very easy interpretation. For example, recall that Gale duals of pyramids are very easy to deal with. Indeed, if ${\vv A}'$ is a pyramid over ${\vv A}$, then $\Gale {({\vv A}')}=\Gale {\vv A}\cup \{\veczero\}$, adding the origin to $\Gale {\vv A}$.

\begin{lemma}\label{lem:pyr}
 If ${\vv A}'$ is a pyramid over ${\vv A}$ then $\degc ({\vv A}')=\degc ({\vv A})$.
\end{lemma}
\begin{proof}
For every linear hyperplane ${\vvh H}$, we have ${\vvh H}^+\cap \pGale{\vv A'}={\vvh H}^+\cap\Gale{\vv A}$; hence, $\degG (\Gale{\vv A})=\degG (\pGale{\vv A'})$.
\end{proof}

This motivates the following definition.

\begin{definition}
 We say that a vector configuration ${\vv V}$ is \defn{irreducible}\index{vector configuration!irreducible} if it does not contain the origin, that is, if $\Gale {\vv V}$ is not a pyramid.
\end{definition}

Here is a simple observation about irreducible vector configurations.

\begin{proposition}\label{prop:easybound}
An irreducible vector configuration ${\vv V}\in \RR^r$ of dual degree $\dd$ cannot contain more than $2r+2\dd$ vectors.
\end{proposition}
\begin{proof}
 Take any generic linear hyperplane ${\vvh H}$, so that $\vv V \cap \vv H = \emptyset$. By the definition of $\degG$, there are at most $r+\dd$ vectors in~${\vvh H}^+$ and in~${\vvh H}^-$.
\end{proof}

Phrasing this in terms of  the primal setting proves a result we alluded to before:

\begin{corollary}\label{cor:easybound}
Any $d$-dimensional configuration $\vv A$ of $n=r+d+1$ points with
$d \ge r+2\degc(\vv A)$
is a pyramid. \qed
\end{corollary}

Finally, an observation that we will use later.

\begin{lemma}\label{lem:quotientirreducible}
 Let $\vv V$ be a vector configuration, and $\vv W$ a subconfiguration fulfilling $\lin(\vv W)\cap \vv V=\vv W$. Then $\vv V/\vv W$ is irreducible.
\end{lemma}
\begin{proof}
Inded, if $\pi$ is the projection that sends each $\vv v\in \vv V\setminus \vv W$ to a vector in $\vv V/\vv W$ (recall the definition of contraction for vector configurations in Section~\ref{sec:operations}), then $\pi(\vv v)=\veczero$ if and only if $\vv v\in\lin(\vv W)$.
\end{proof}

\section{\texorpdfstring{\CayleyG vector configurations}{Cayley* vector configurations}}\label{sec:DualCayley}

The concepts of Section~\ref{sec:PrimalCayley} can be formulated in the Gale dual setting.

\begin{definitions}
 A vector configuration ${\vv V}$ is 
\begin{itemize}
 \item an \defn{affine \CayleyG configuration}\index{\CayleyG configuration!affine} of length $m$, if there exists a partition ${\vv V}=\vv V_1\uplus \cdots \uplus \vv V_m$ such that $\sum_{\vv v_j\in \vv V_i} \vv v_j=\veczero$ for $i=1\ldots m$.

 \item a \defn{combinatorial \CayleyG configuration}\index{\CayleyG configuration!combinatorial} of length $m$, 
if there exists a partition  ${\vv V}=\vv V_1\uplus \cdots \uplus \vv V_m$ such that  $\vv V_i$ is a positive vector of $\cM({\vv V})$ for each~$i$. 
That is, for each $\vv V_i$ there is a positive vector ${\vv \gl^{(i)}}\in \RR^{|\vv V_i|}$ such that $\sum_{\vv v_j\in \vv V_i} \gl^{(i)}_j \vv v_j=\veczero$.
\end{itemize}
The sets $\vv V_i$ are called the \defn{factors} of the \CayleyG configuration.
\end{definitions}

These concepts coincide with their primal counterparts:

\begin{proposition}\label{cayley-dual}
 ${\vv A}$ is an affine (resp. combinatorial) Cayley configuration of length~$m$ if and only if its Gale dual ${\vv V}$ is an affine (resp. combinatorial) \CayleyG configuration of length~$m$.
\end{proposition}
\begin{proof}
Let ${\vv A}$ be a combinatorial Cayley configuration of length $m$. Then there is a partition ${\vv A} = \vv A_1 \uplus \cdots \uplus \vv A_m$ such that $\vv A\setminus \vv A_i$ is the set of points in a proper face of~$\conv({\vv A})$, for any $1\leq i\leq m$. 
Let $\vv V=\Gale{\vv A}$ be the dual of~$\vv A$, and for $1\leq i\leq m$ let~$\vv V_i:=\set{\vv v_j\in \vv V}{\vv a_j\in \vv A_i}$  be the disjoint subsets of~$\vv V$ corresponding to the respective~$\vv A_i$. 
Then, by Lemma~\ref{lem:gale}, every $\vv V_i$ is a positive vector of $\cM({\vv V})$. Thus, ${\vv V}$ is a combinatorial \CayleyG configuration of length~$m$. The converse is direct.

In the affine case, ${\vv A}$ is an affine Cayley configuration of length $m$ when there is an affine projection $\pi:\RR^d \to \RR^{m-1}$ that maps ${\vv A}$ onto the vertex set of a $(m-1)$-simplex $\Delta_{m-1}$ with vertices $\{\vv w_1,\dots,\vv w_m\}$. 
Let $\vv A_i:={\vv A}\cap \pi^{-1}(\vv w_i)$ and observe that there is an affine function~$f_i$ such that $f_i(\vv a)=1$ if~$\vv a\in\vv A_i$ and~$f_i(\vv a)=0$ otherwise. 
Let ${\vv V} = \{\vv v_1, \ldots, \vv v_n\}$ be the Gale dual of ${\vv A} = \{\vv a_1, \ldots, \vv a_n\}$, where $\vv v_i$ corresponds to $\vv a_i$. 
For $1\leq i\leq m$ we define $\vv V_i := \{\vv v_j \;:\; \vv a_j \in \vv A_i\}$. By duality, affine valuations on ${\vv A}$ correspond to linear dependences of ${\vv V}$. Hence,  we obtain $\sum_{j=1}^n f_i(\vv a_j)\vv v_j=\sum_{\vv v_j\in\vv  V_i} \vv v_j=\veczero$, so that ${\vv V}$ is an affine \CayleyG configuration with factors $\vv V_i$. Again, the converse follows similarly.
\end{proof}

We can now make good on our promise from Section~\ref{sec:PrimalCayley} to show that all combinatorial Cayley configurations can be realized by affine Cayley configurations:

\begin{proposition}\label{prop:combinatorialisaffine}
Every combinatorial Cayley configuration of length~$m$ is combinatorially equivalent (as an oriented matroid) to an affine Cayley configuration of length~$m$.
\end{proposition}

\begin{proof}
Let ${\vv V} = \vv V_1 \uplus \cdots \uplus \vv V_m$ be the combinatorial \CayleyG configuration whose Gale dual is ${\vv A}$. Then, for each factor $\vv V_i$, there is a vector ${\vv \gl^{(i)}}$ in $\RR^{|\vv V_i|}$ such that $\sum_{\vv v_j\in\vv  V_i}\gl^{(i)}_{j} \vv v_j=0$. 
Let $\vv \gl$ be the vector in $\RR^{|{\vv V}|}$ with entries $\gl_{j}= \gl^{(i)}_{j}$ if $\vv v_j\in\vv V_i$. Scaling the vectors of a vector configuration does not affect its oriented matroid. Hence, $\{\gl_{j}\vv v_j\}_{\vv v_j\in {\vv V}}$ is an affine \CayleyG configuration combinatorially equivalent to ${\vv V}$.
\end{proof}

This result motivates the use of the generic term \defn{Cayley configuration} for combinatorial Cayley configurations.

\subsection{\texorpdfstring{Weak \CayleyG configurations}{Weak Cayley* configurations}}

Weak Cayley configurations also have a dual version.

\begin{definition}
 A vector configuration ${\vv V}$ is a \defn{weak \CayleyG configuration}\index{\CayleyG configuration!weak} of length $m$ if it contains $m$ disjoint positive circuits of $\cM({\vv V})$. 
\end{definition}

\iftoggle{bwprint}{%
\begin{figure}[htpb]
\centering
 \subbottom[$\vv V$]{\label{sfig:dualweakV}\includegraphics[width=.2\linewidth]{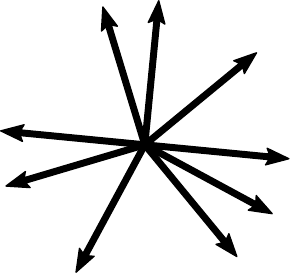}}\quad
 \subbottom[$\vv V_0$]{\label{sfig:dualweakV0}\includegraphics[width=.15\linewidth]{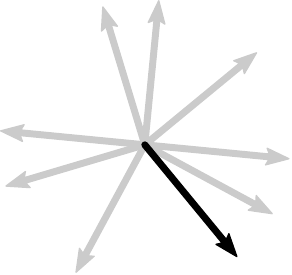}}\quad
 \subbottom[$\vv V_1$]{\label{sfig:dualweakV1}\includegraphics[width=.15\linewidth]{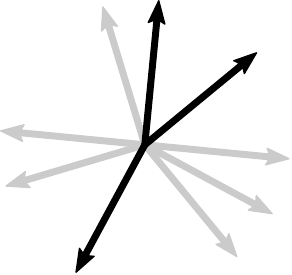}}\quad
 \subbottom[$\vv V_2$]{\label{sfig:dualweakV2}\includegraphics[width=.15\linewidth]{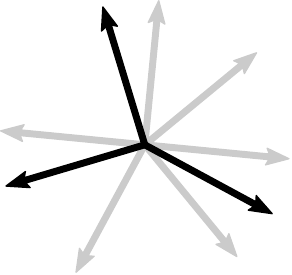}}\quad
 \subbottom[$\vv V_3$]{\label{sfig:dualweakV3}\includegraphics[width=.15\linewidth]{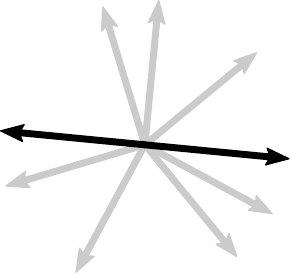}}
  \caption{A weak \CayleyG configuration of length $3$.}
 \label{fig:exdualweak}
\end{figure}
}{%
\begin{figure}[htpb]
\centering
 \subbottom[$\vv V$]{\label{sfig:dualweakV}\includegraphics[width=.2\linewidth]{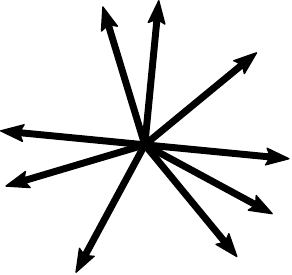}}\quad
 \subbottom[$\vv V_0$]{\label{sfig:dualweakV0}\includegraphics[width=.15\linewidth]{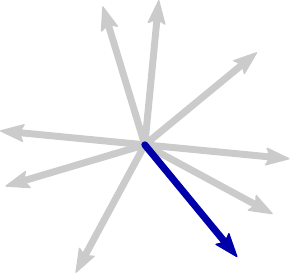}}\quad
 \subbottom[$\vv V_1$]{\label{sfig:dualweakV1}\includegraphics[width=.15\linewidth]{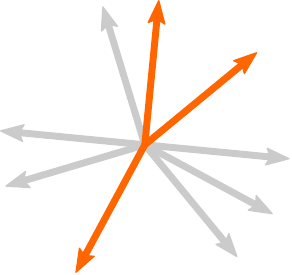}}\quad
 \subbottom[$\vv V_2$]{\label{sfig:dualweakV2}\includegraphics[width=.15\linewidth]{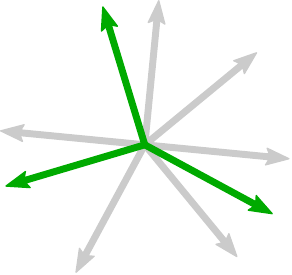}}\quad
 \subbottom[$\vv V_3$]{\label{sfig:dualweakV3}\includegraphics[width=.15\linewidth]{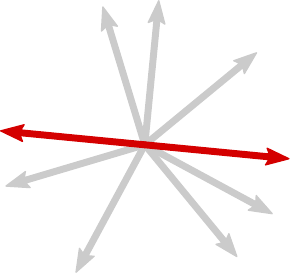}}
  \caption{A weak \CayleyG configuration of length $3$.}
 \label{fig:exdualweak}
\end{figure}
}

In Figure~\ref{fig:exdualweak} there is an example of a weak \CayleyG configuration of length $3$ in $\RR^2$. Observe that the decomposition into factors is not unique; for example, the vector in $\vv V_0$ could be swapped with a vector in $\vv V_1$ or $\vv V_2$.
\\

While ${\vv A}$ is a weak Cayley configuration if and only if it contains a (possibly empty) subset $\vv A_0$ such that the contraction ${\vv A}/\vv A_0$ is a combinatorial Cayley configuration of length $m$, ${\vv V}$ is a weak \CayleyG configuration if and only if it contains a subset $\vv V_0$ such that the deletion ${\vv V}\setminus\vv  V_0$ is a combinatorial \CayleyG configuration of length $m$.
Using this observation and Proposition~\ref{cayley-dual}, it follows directly from the duality of deletion and contraction that this definition is consistent with the primal version:

\begin{proposition}
 ${\vv A}$ is a weak Cayley configuration of length $m$ if and only if $\Gale {\vv A}$ is a weak \CayleyG vector configuration of length $m$.\qed
\end{proposition}

We can now prove the promised estimate on the combinatorial degree of weak Cayley configurations:

\begin{proposition}\label{prop:deg-cayley}
  If $\vv A$ is a $d$-dimensional weak Cayley configuration of length~$m$ then $\degc(\vv A)\leq d+1-m$. 
\end{proposition}

\begin{proof}
 If ${\vv V}$ is a weak \CayleyG vector configuration in $\RR^r$ whose factors are $\vv V_1,\dots,\vv V_m$, then every linear hyperplane ${\vvh H}$ contains at least one element of every factor in $\ol {\vvh H}^-$. Therefore $|{\vvh H}^+\cap {\vv V}|\leq n-m$ for any ${\vvh H}$, which proves that $\degG({\vv V})\leq n-r-m=d+1-m$.
\end{proof}

Moreover, we can easily see why Conjecture~\ref{conj:d+1-2dd} is sharp:

\begin{example}\label{ex:sharpconj}
 Let $\vv P$ be a neighborly polytope in even dimension $d=2e$ with $n\geq 2d+1$ vertices. Recall that $\vv P$ must be simplicial by Remark~\ref{rmk:simplicialneighborly}. 
 Therefore, its Gale dual ${\vv V}$ is a vector configuration in general position in $\RR^r$ (no $r$ vectors in the same linear hyperplane), where $2r+1\geq n$ because $r=n-d-1$ and $n\geq 2d+1$. Since $\vv P$ is neighborly, $\dd:=\degG({\vv V})=\degc(\verts(\vv P))=e$. 
 Regarding Conjecture~\ref{conj:d+1-2dd} we see that ${\vv V}$ is a weak \CayleyG configuration of length $d+1-2\dd=1$, but it cannot be a weak \CayleyG configuration of length $2$. 
 Indeed, since the vectors in ${\vv V}$ are in general position, each circuit~$C$ of ${\vv V}$ has cardinality $r+1$; and since $n<2r+2$, ${\vv V}$ cannot contain two disjoint circuits.
\end{example}

\begin{observation}\label{obs:divisibleIsCayley}
Observe that if ${\vv V}$ is a vector configuration and $\ol {\vv V}$ is its set of endpoints, then ${\vv V}$ is a weak \CayleyG configuration of length $m$ if and only if  $\veczero\in \Tv_m(\ol {\vv V})$, \ie the origin is a $m$-divisible point of $\ol {\vv V}$.
Together with Observation~\ref{obs:depthIsCodeg}, this explains why Theorem~\ref{thm:d+1-3dd} is equivalent to Corollary~\ref{cor:CoreIsDivisible}.
\end{observation}

\section{Small degree implies weak Cayley}\label{sec:degk}
The following proposition relates the degree of the restriction of a vector configuration to a subspace to the degree of its contraction. It is the main ingredient of the proof of Theorem~\ref{thm:d+1-3dd}. In Section~\ref{sec:extremalconfigurations}, the subconfigurations that attain equality in~\eqref{eq:decompositioninequality} are further investigated.

\begin{proposition}\label{prop:subspacequotient}
Let ${\vv V}$ be a vector configuration and let ${\vv W}\subset {\vv V}$ be a subconfiguration of ${\vv V}$ such that $\lin ({\vv W})\cap {\vv V}={\vv W}$.
If we use the notation
\begin{itemize}
 \item $\rank ({\vv V})=r$, $|{\vv V}|=r+d+1$ and $\degG({\vv V})=\dd$;
 \item $\rank ({\vv W})=r_{\vv W}$, $|{\vv W}|=r_{\vv W}+d_{\vv W}+1$ and $\dd_{\vv W}=\degG({\vv W})$ (in $\RR^{r_{\vv W}}$); and
 \item $\rank ({\vv V}/{\vv W})=r_{/{\vv W}}$, $|{\vv V}/{\vv W}|=r_{/{\vv W}}+d_{/{\vv W}}+1$ and $\dd_{/{\vv W}}=\degG({\vv V}/{\vv W})$,
\end{itemize}
then
\begin{align}
 r&=r_{\vv W}+r_{/{\vv W}},\notag\\
 d&=d_{\vv W}+d_{/{\vv W}}+1,\notag\\
 \dd&\geq \dd_{\vv W}+\dd_{/{\vv W}}.\label{eq:decompositioninequality}
 \end{align}
\end{proposition}
\begin{proof}
By construction, $r=r_{\vv W}+r_{/{\vv W}}$. Moreover, counting the number of elements in ${\vv V}$ we get $r+d+1=r_{\vv W}+d_{\vv W}+1+r_{/{\vv W}}+d_{/{\vv W}}+1$, which implies that $d=d_{\vv W}+d_{/{\vv W}}+1$. 

Since the degree of ${\vv W}$ is $\dd_{\vv W}$, there is an oriented hyperplane ${\vvh H}_{\vv W}$ of $\lin ({\vv W})$ that contains $r_{\vv W}+\dd_{\vv W}$ elements of ${\vv W}$ in ${\vvh H}_{\vv W}^+$. 
Let ${\vvh H}_{\vv W}'$ be a hyperplane of $\RR^r$ such that ${\vvh H}'_{\vv W}\cap \lin {\vv W}={\vvh H}_{\vv W}$. Note that such a hyperplane always exists, for example the only hyperplane that contains ${\vvh H}_{\vv W}$ and the orthogonal complement of $\lin ({\vv W})$. 
Since ${\vv V}/{\vv W}$ has degree $\dd_{/{\vv W}}$, there is an oriented hyperplane ${\vvh H}_{/{\vv W}}$ of the quotient $\vv V/\vv W$ that has $r_{/{\vv W}}+\dd_{/{\vv W}}$ elements of~${\vv V/\vv W}$ at ${\vvh H}_{/{\vv W}}^+$. 
By definition of quotient, there is a hyperplane~${\vvh H}_{/{\vv W}}'$ of~$\RR^r$ that contains $\lin ({\vv W})$ and such that ${\vvh H}_{/{\vv W}}'^+\cap \vv V={\vvh H}_{/{\vv W}}^+\cap \vv V/\vv W$ (identifying elements of $\vv V/\vv W$ with the corresponding elements of $\vv V$). Then 
\begin{align*}
r+\dd&\geq |({\vvh H}_{/{\vv W}}'\circ {\vvh H}_{\vv W}')\cap {\vv V}|\\&=|{\vvh H}_{/{\vv W}}^+\cap \vv V/\vv W|+|{\vvh H}_{{\vv W}}^+\cap \vv W|=r_{/{\vv W}}+\dd_{/{\vv W}}+r_{\vv W}+\dd_{\vv W}.
\end{align*} 
therefore, $\dd_{{\vv W}}+\dd_{/{\vv W}}\leq \dd$.
\end{proof}

Observe that we took the ``worst'' hyperplane in $\RR^r$ containing $\lin ({\vv W})$ (worst in terms of $|{\vvh H}^+\cap \vv V|$), and slightly perturbed it so that it cut $\lin ({\vv W})$ in its worst hyperplane. The proposition states that this perturbed hyperplane cannot be worse than the worst hyperplane that cuts $\vv V$.

\iftoggle{bwprint}{%
\begin{figure}[htpb]
\centering
 \subbottom[$\vv V$]{\label{sfig:sq_V}\,\includegraphics[width=.26\linewidth]{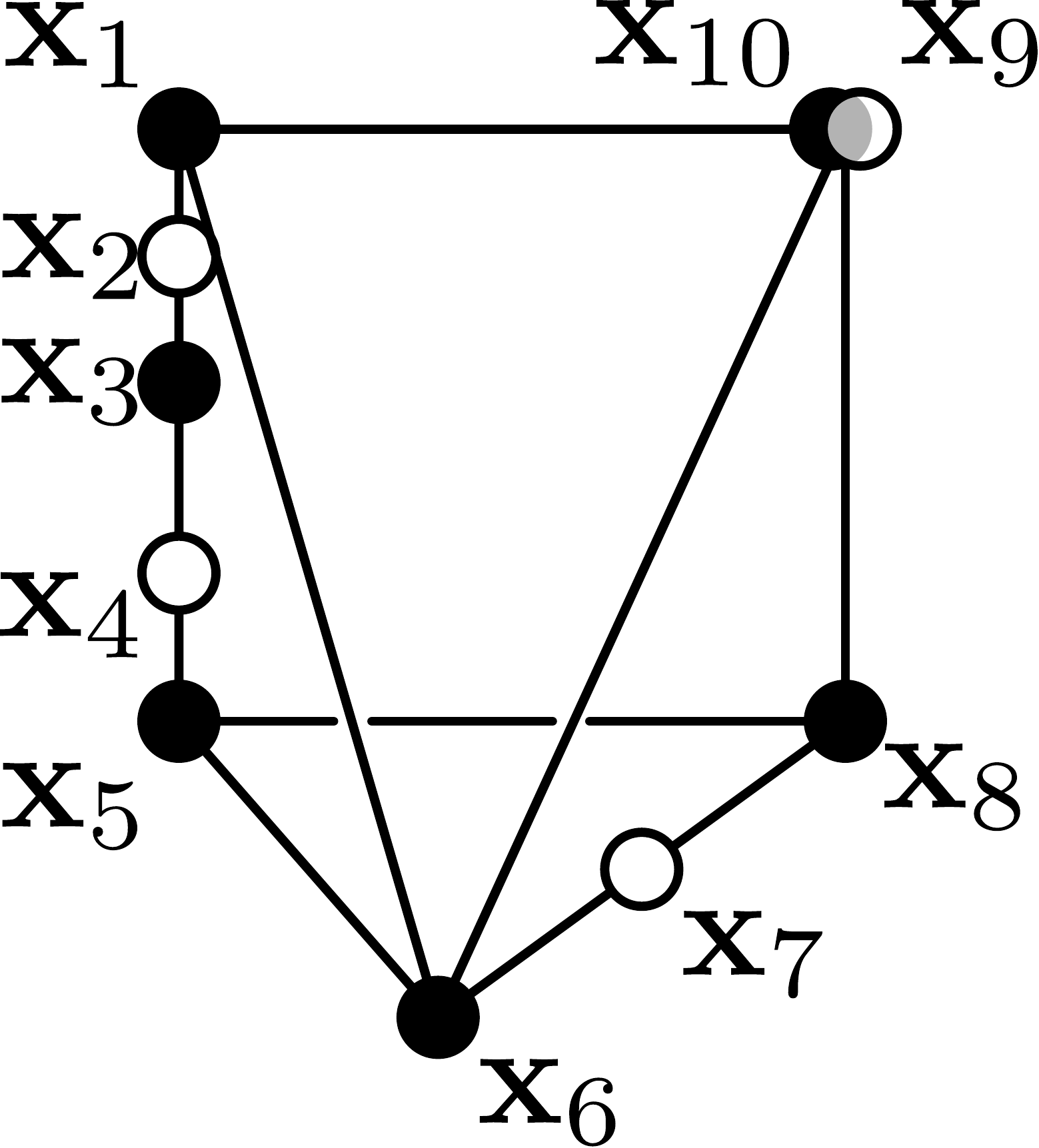}\,}\,\,\qquad
 \subbottom[$\vv W$]{\label{sfig:sq_W}\,\includegraphics[width=.26\linewidth]{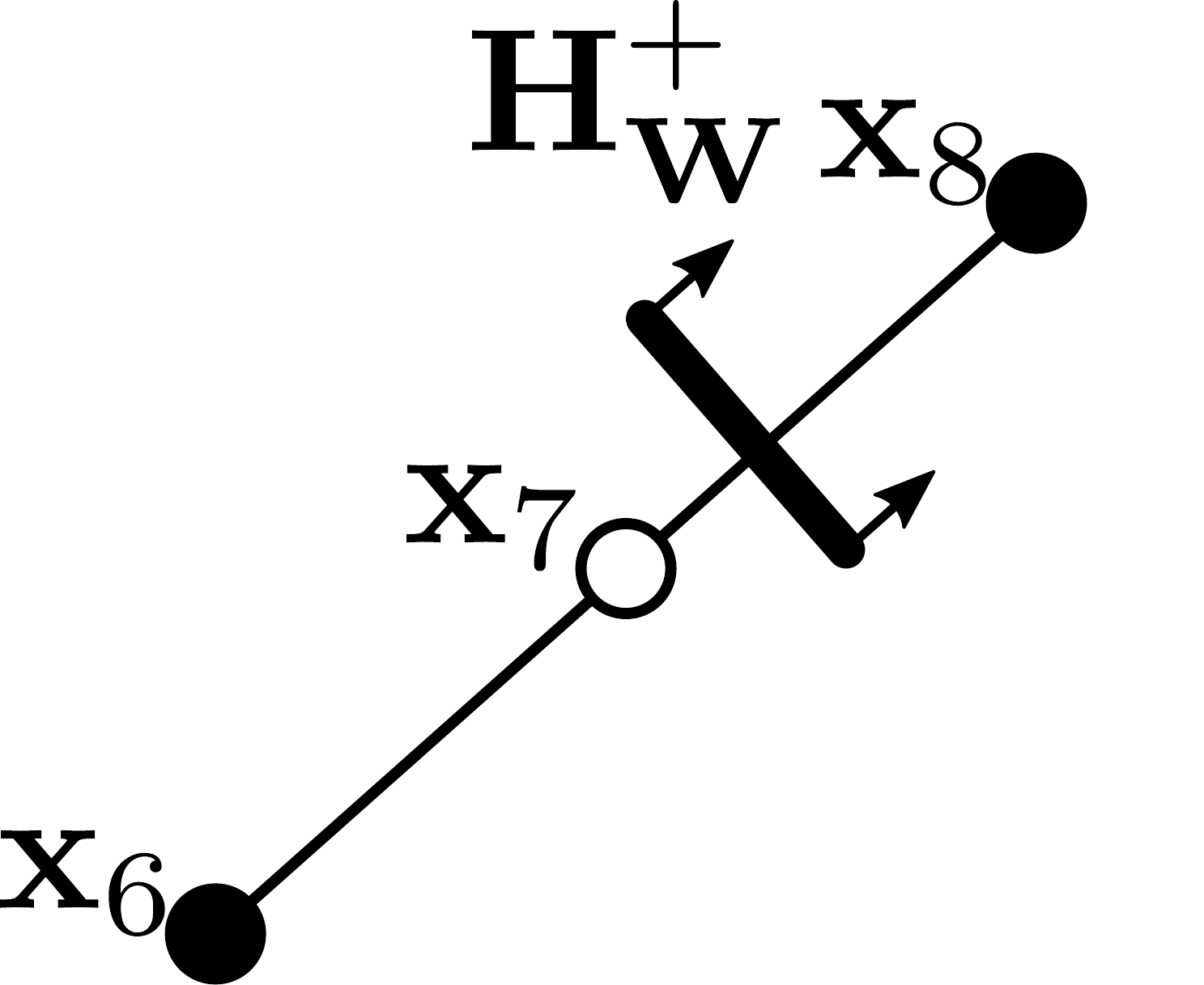}\,}\,\,\qquad
 \subbottom[$\vv V/\vv W$]{\label{sfig:sq_modW}\,\includegraphics[width=.26\linewidth]{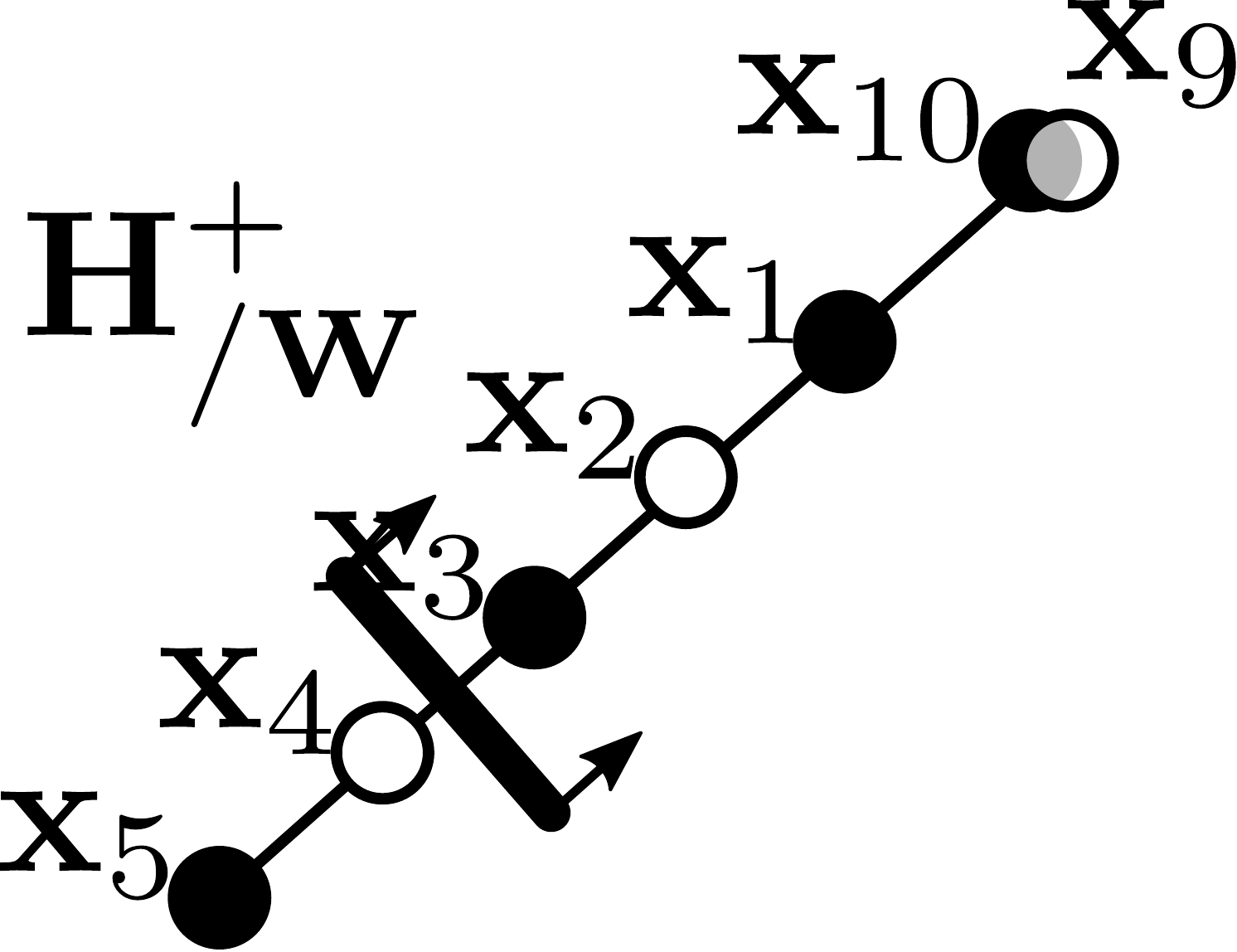}\,}
 \subbottom[$\vvh H_{\vv W}'$]{\label{sfig:sq_W'}\,\includegraphics[width=.26\linewidth]{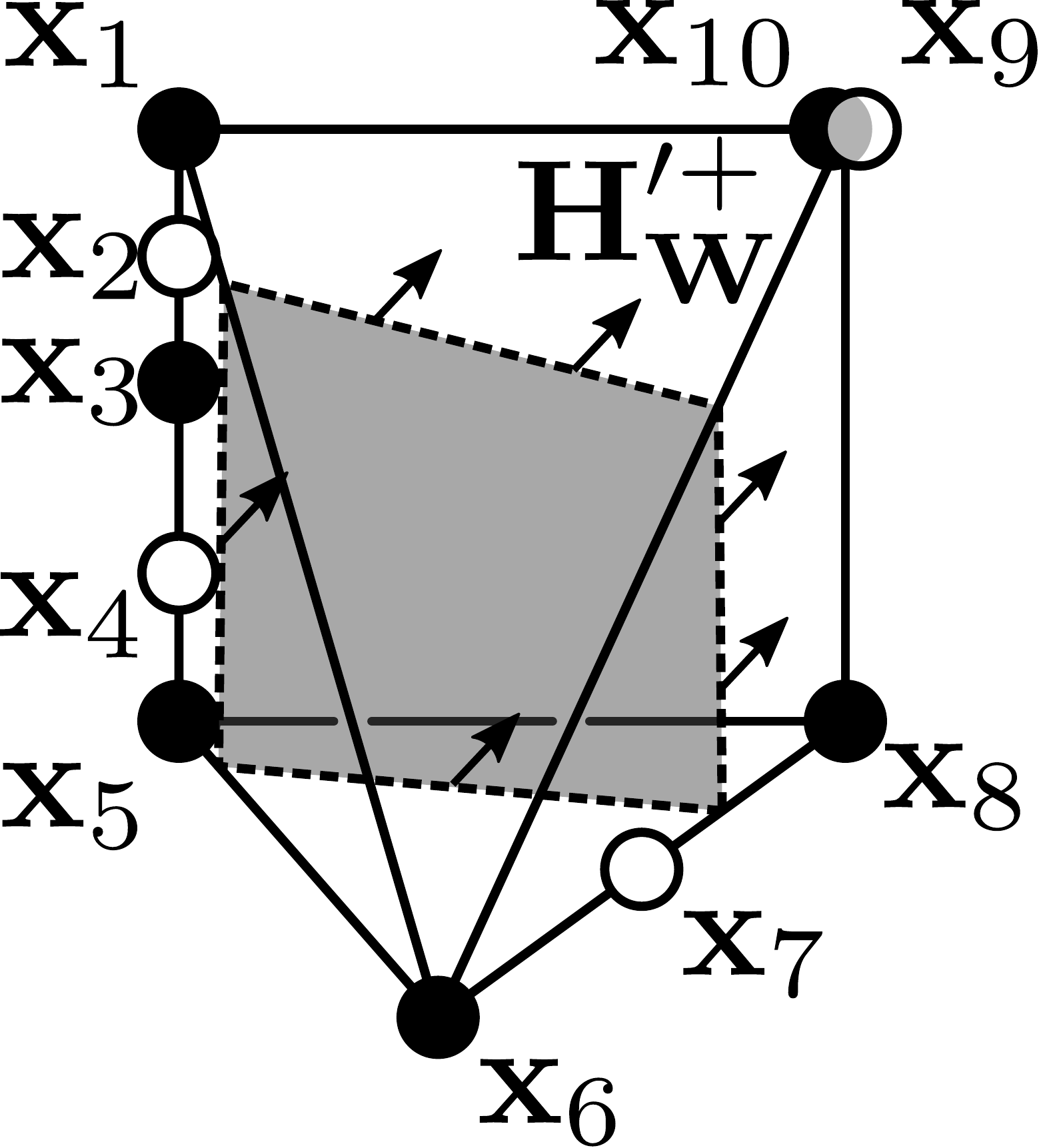}\,}\,\,\qquad
 \subbottom[$\vvh H_{/\vv W}'$]{\label{sfig:sq_modW'}\,\includegraphics[width=.26\linewidth]{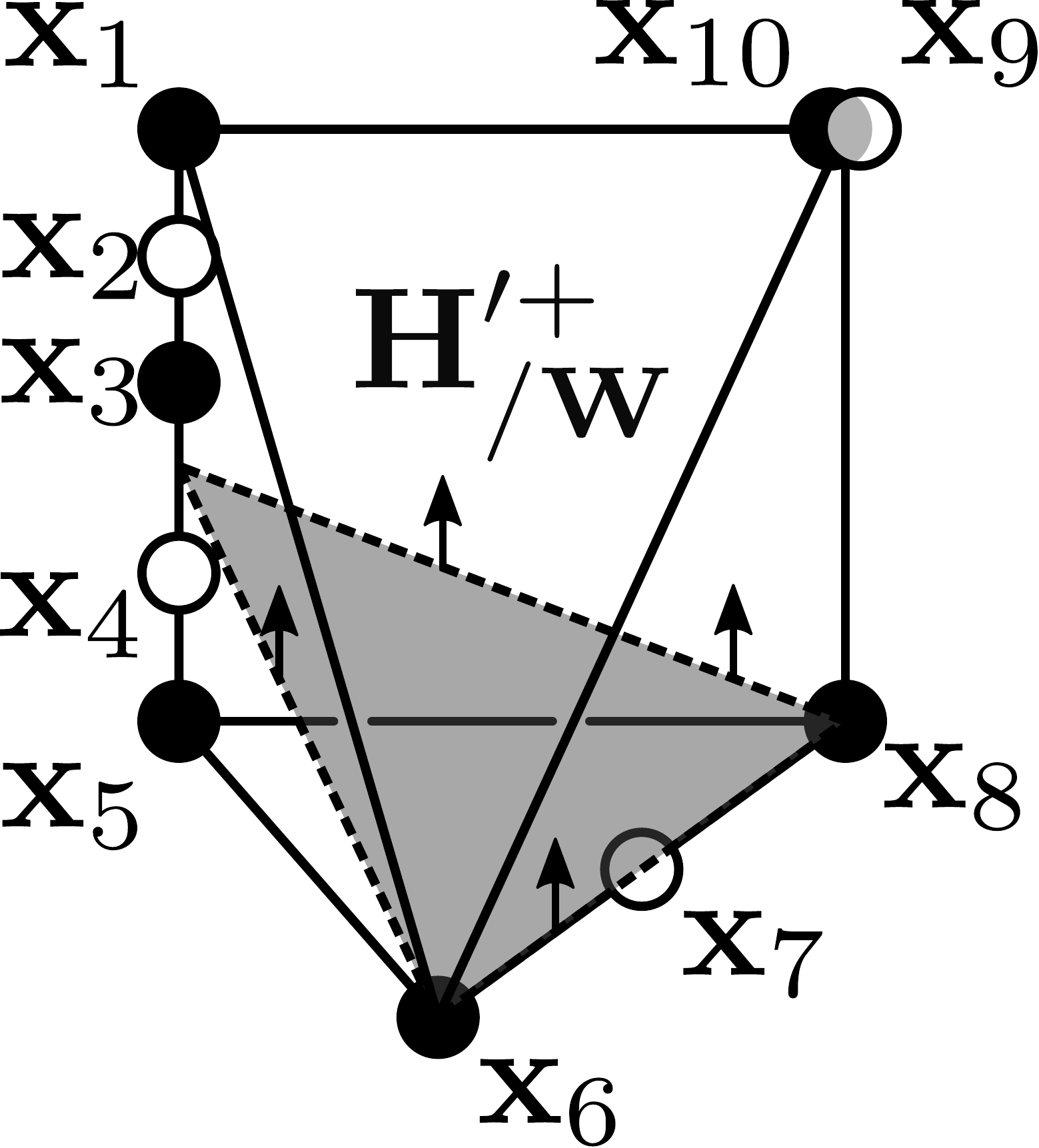}\,}\,\,\qquad
 \subbottom[$\vvh H={\vvh H'}_{/{\vv W}}\circ {\vvh H}_{\vv W}'$]{\label{sfig:sq_H}\,\includegraphics[width=.26\linewidth]{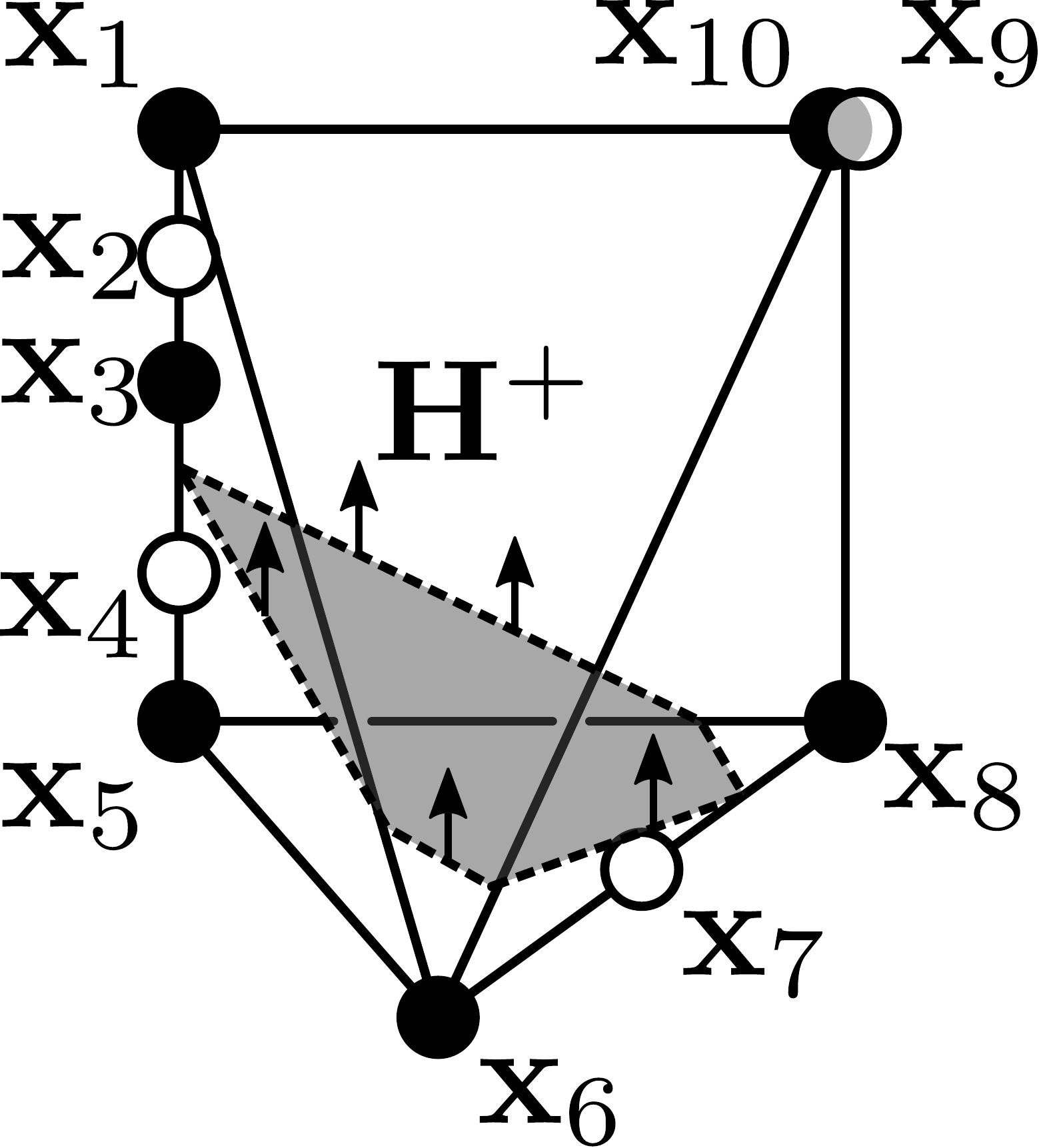}\,} 
\caption[Illustrating Proposition~\ref{prop:subspacequotient}.]{Illustrating Proposition~\ref{prop:subspacequotient}. Overlapping circles represent points that have the same coordinates.}
 \label{fig:subspacequotient}
\end{figure}
}{%
\begin{figure}[htpb]
\centering
 \subbottom[$\vv V$]{\label{sfig:sq_V}\,\includegraphics[width=.26\linewidth]{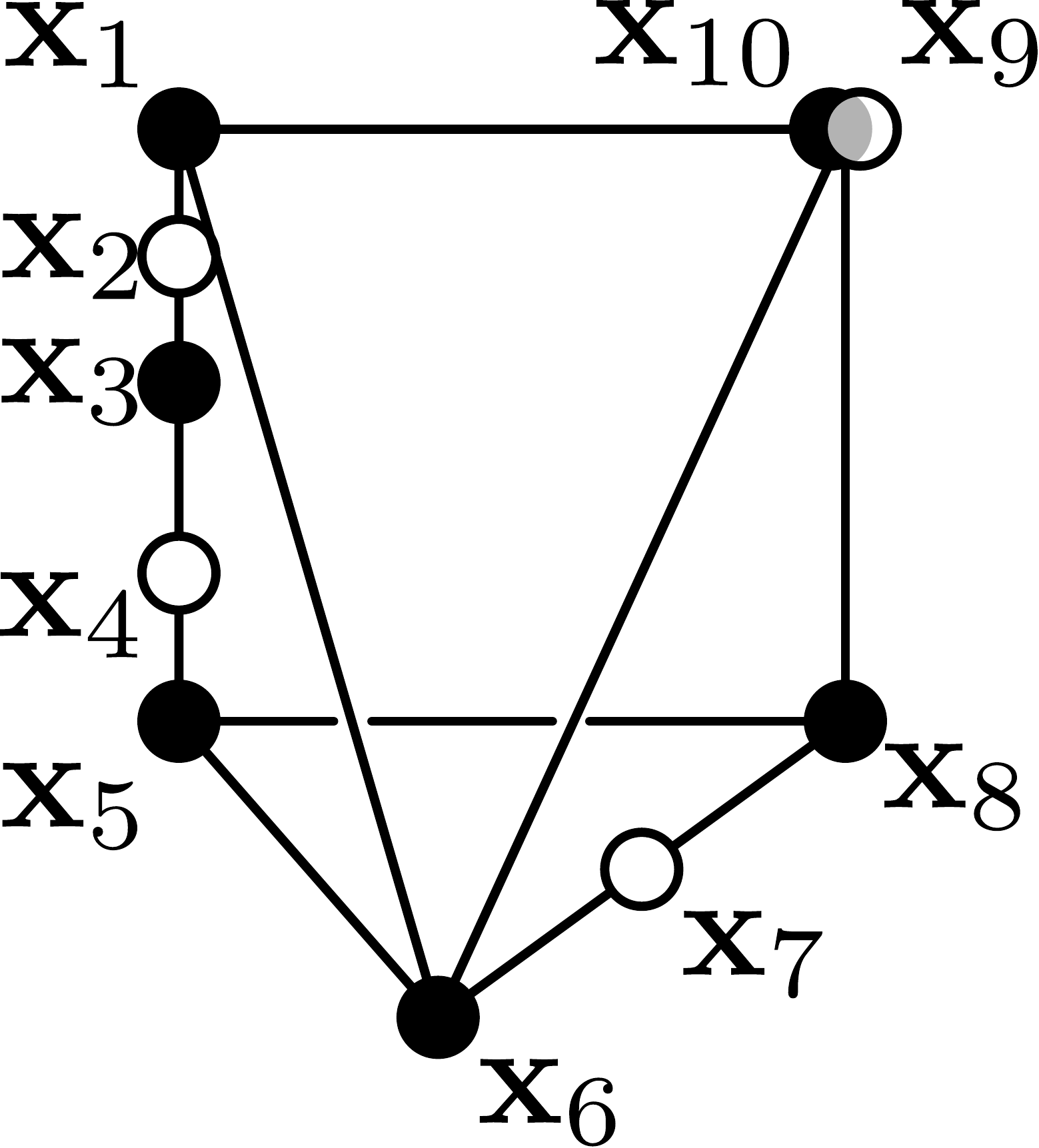}\,}\,\,\qquad
 \subbottom[$\vv W$]{\label{sfig:sq_W}\,\includegraphics[width=.26\linewidth]{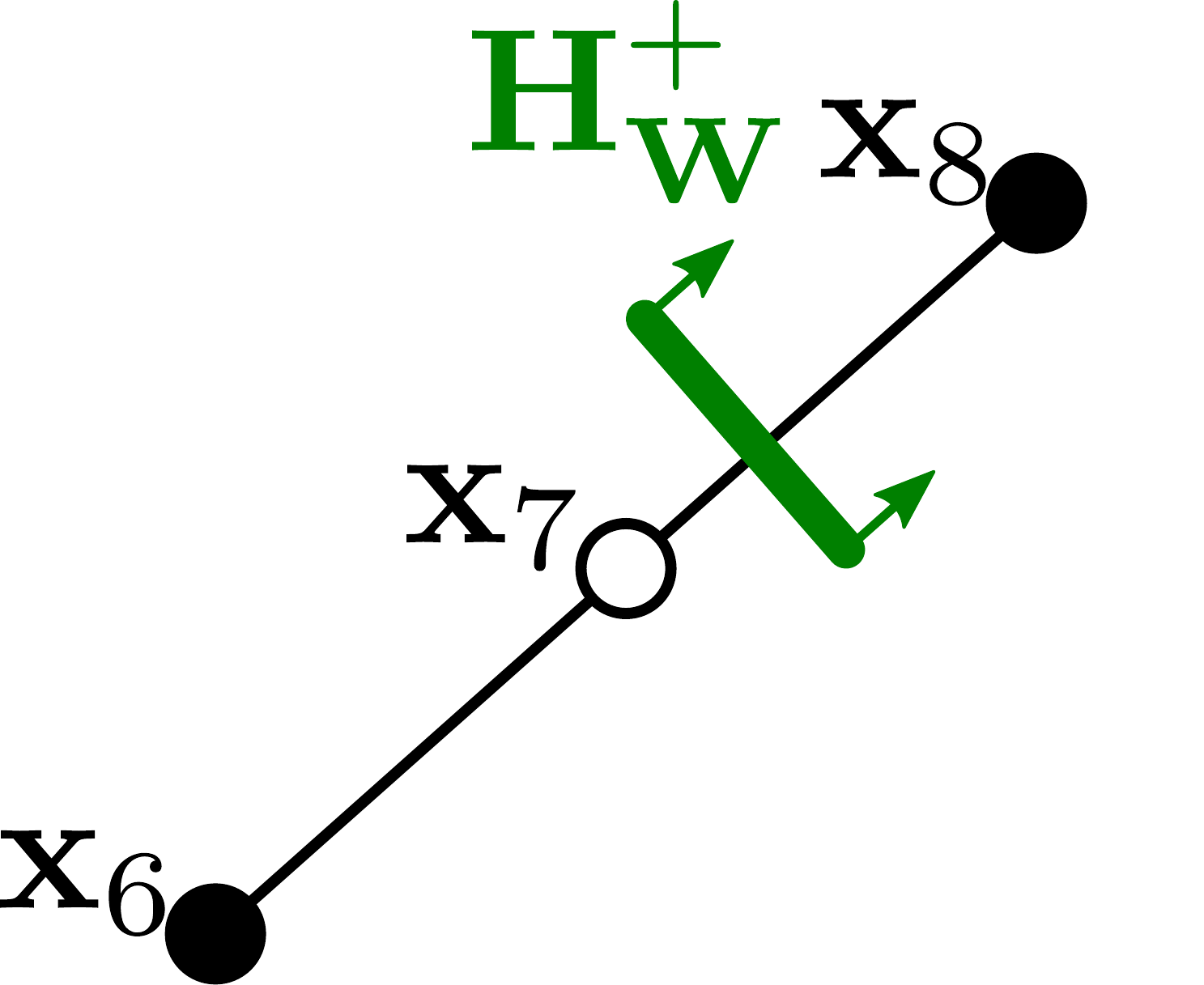}\,}\,\,\qquad
 \subbottom[$\vv V/\vv W$]{\label{sfig:sq_modW}\,\includegraphics[width=.26\linewidth]{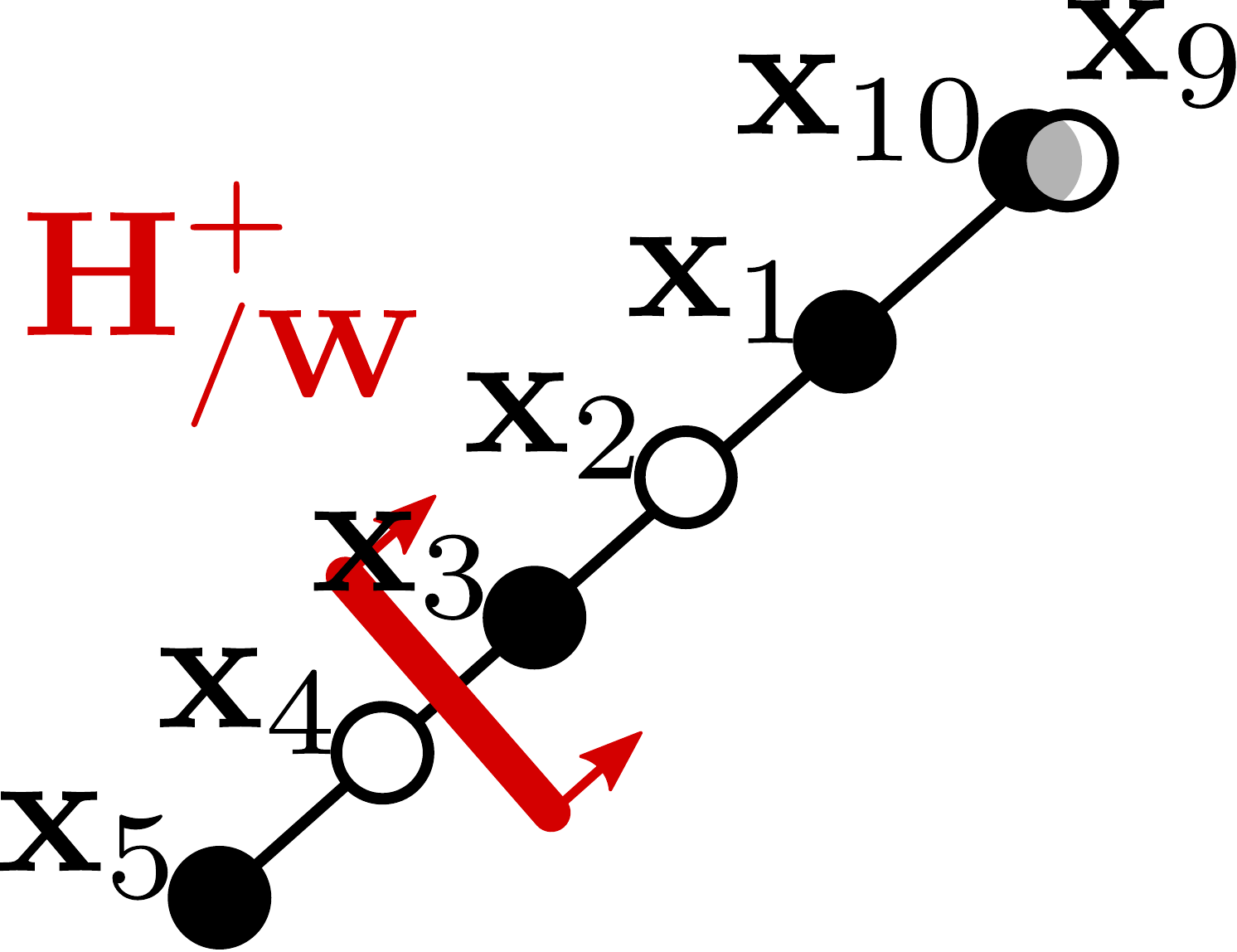}\,}
 \subbottom[$\vvh H_{\vv W}'$]{\label{sfig:sq_W'}\,\includegraphics[width=.26\linewidth]{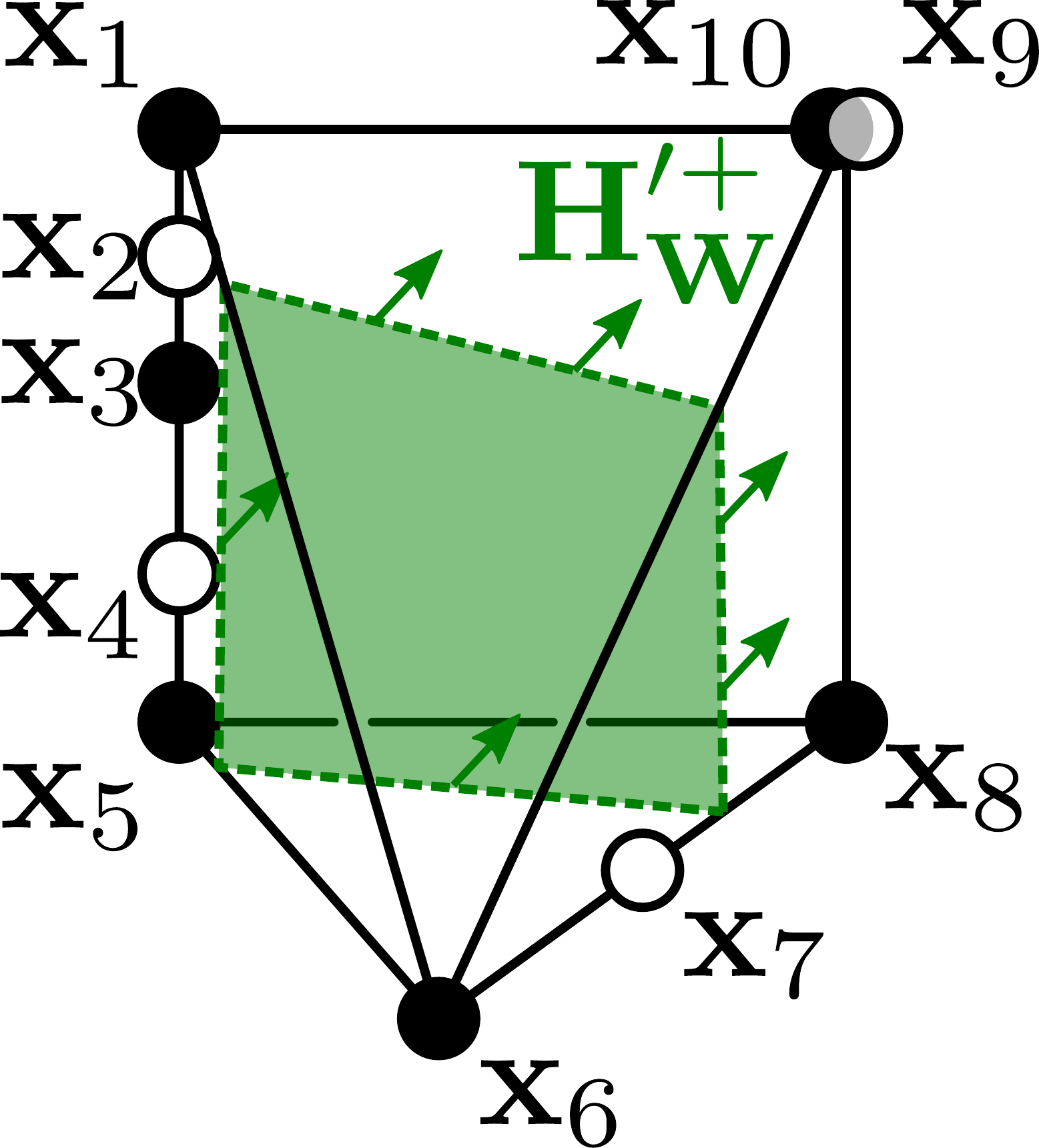}\,}\,\,\qquad
 \subbottom[$\vvh H_{/\vv W}'$]{\label{sfig:sq_modW'}\,\includegraphics[width=.26\linewidth]{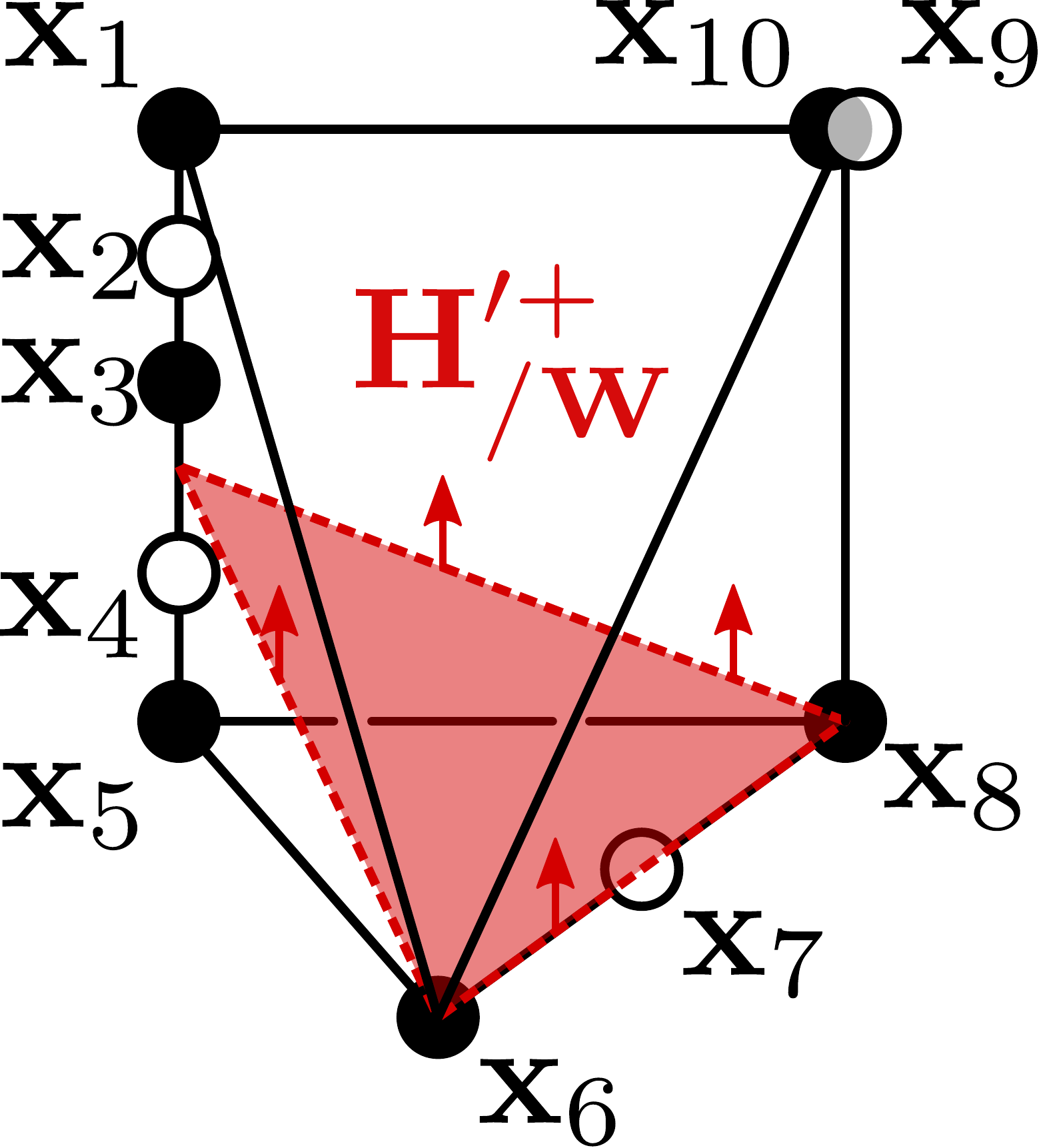}\,}\,\,\qquad
 \subbottom[$\vvh H={\vvh H'}_{/{\vv W}}\circ {\vvh H}_{\vv W}'$]{\label{sfig:sq_H}\,\includegraphics[width=.26\linewidth]{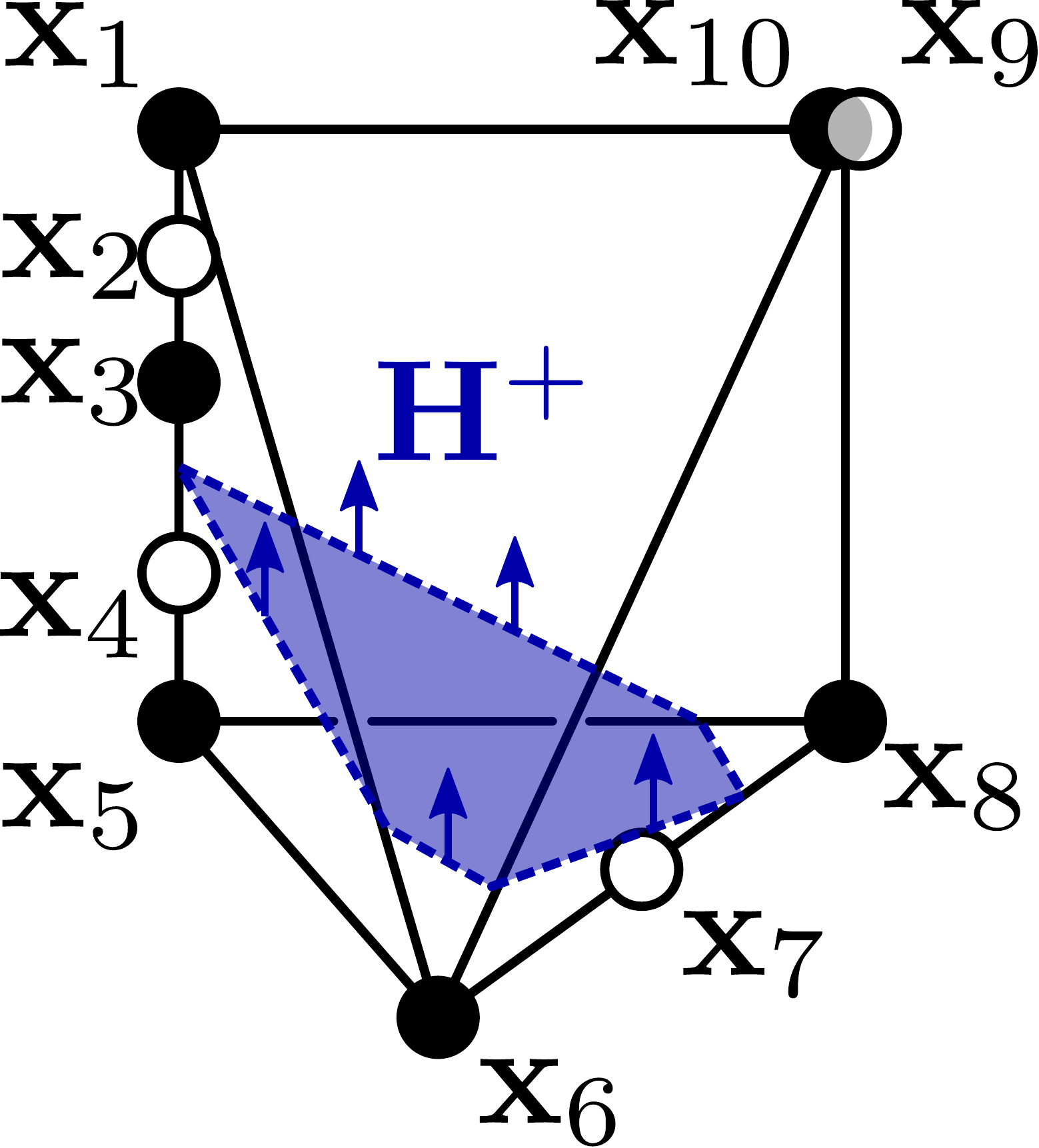}\,} 
\caption[Illustrating Proposition~\ref{prop:subspacequotient}.]{Illustrating Proposition~\ref{prop:subspacequotient}. Overlapping circles represent points that have the same coordinates.}
 \label{fig:subspacequotient}
\end{figure}
}

\begin{example}
 In Figure~\ref{sfig:sq_V} there is an affine diagram representing a vector configuration $\vv V$ of rank~$4$ with $10$ elements and $\degG(\vv V)=2$. In~\ref{sfig:sq_W} we can see $\vv W=\{\vv x_6,\vv x_7,\vv x_8\}$, a subconfiguration of $\vv V$ of rank~$2$ with $3$ elements and $\degG(\vv W)=0$. 
The hyperplane $\vvh H_{\vv W}$ fulfills $|\vvh H_{\vv W}\cap \vv W|=2$. The quotient $\vv V/\vv W$ is shown in~\ref{sfig:sq_W}. It is a configuration of rank~$2$ with~$7$ elements and $\degG(\vv V/\vv W)=2$. The hyperplane $\vvh H_{/\vv W}$ fulfills $|\vvh H_{/\vv W}\cap \vv V/\vv W|=2$. 
In~\ref{sfig:sq_W'} a hyperplane ${\vvh H}_{\vv W}'$ fulfilling ${\vvh H}'_{\vv W}\cap \lin {\vv W}={\vvh H}_{\vv W}$ is depicted. And in~\ref{sfig:sq_modW'} we show ${\vvh H}'_{/{\vv W}}$, a hyperplane that contains $\lin ({\vv W})$ and has $r_{/{\vv W}}+\dd_{/{\vv W}}$ elements of~${\vv V}$ in ${\vvh H}_{/{\vv W}}^+$. 
Now $\vvh H={\vvh H'}_{/{\vv W}}\circ {\vvh H}_{\vv W}'$ contains $r+\degG({\vv W})+\degG({\vv V/{\vv W}})=6$ elements in $\vvh H^+$, which can be verified in~\ref{sfig:sq_H}.
\end{example}

We will use Proposition~\ref{prop:subspacequotient} to prove Theorem~\ref{thm:d+1-3dd}. Recall that in the dual setting our goal is to find many disjoint positive circuits. In our proof, we will iteratively find a subconfiguration ${\vv W}$ of ${\vv V}$ of lower rank that has smaller dual degree. 
Eventually we will find a configuration of degree $0$, and Corollary~\ref{cor:deg0} will certify that in this subconfiguration there are already many disjoint positive circuits.

\begin{theorem}\label{thm:d+1-3dd}
 Let $\vv A$ be a $d$-dimensional configuration. If $\degc(\vv A)<\frac{d}{3}$, then $\vv A$ is a weak Cayley configuration of length at least $d-3\degc(\vv A)+1$.
\end{theorem}

\begin{proof}
We will prove the dual statement, which says that any vector configuration~${\vv V}$ with $r+d+1$ elements and degree~$\dd$ is a weak \CayleyG configuration of length at least $d-3\dd+1$.

By Lemma~\ref{lem:removerepeated}, we can assume that ${\vv V}$ is pure. The proof will be by induction on $\dd$. The base case is $\dd=0$, which we know to hold because of Corollary~\ref{cor:deg0}.

Let ${\vvh H}$ be any hyperplane spanned by elements of ${\vv V}$. Let ${\vv W}={\vv V}\cap {\vvh H}$. Then ${\vv V}/{\vv W}$ is irreducible by Lemma~\ref{lem:quotientirreducible}, pure by Lemma~\ref{lem:quotientsofpurearepure} and has rank~$r_{/{\vv W}}=1$ with $d_{/{\vv W}}+2$ elements and degree $\dd_{/{\vv W}}:=\degG({\vv V}/{\vv W})$. By Lemma~\ref{lem:puredeggeq1}, we know that
\begin{equation}\label{eq:ddgeq1}
 \dd_{/{\vv W}}\geq 1.
\end{equation}
From Proposition~\ref{prop:easybound} we can deduce that $(d_{/{\vv W}} -2\dd_{/{\vv W}})\leq r_{/{\vv W}}-1=0$. Therefore the previous equation~\eqref{eq:ddgeq1} implies that 
\begin{equation}
\label{eq:dmodhleq-1}
 (d_{/{\vv W}} -3\dd_{/{\vv W}})= (d_{/{\vv W}} -2\dd_{/{\vv W}})-\dd_{/{\vv W}}\leq -\dd_{/{\vv W}}\leq -1
\end{equation}

On the other hand, ${\vv W}$ is a vector configuration of rank $r-1$ with $r+d_{{\vv W}}$ elements and degree $\dd_{\vv W}:=\degG({\vv W})$. By Proposition~\ref{prop:subspacequotient}, 
\begin{equation}\label{eq:dwleqd-1}
\dd_{\vv W}\stackrel{\eqref{eq:decompositioninequality}}{\leq} \dd-\dd_{/{\vv W}}\stackrel{\eqref{eq:ddgeq1}}{\leq} \dd -1.                                                                                                                    
\end{equation}
Moreover, again by Proposition~\ref{prop:subspacequotient} and \eqref{eq:dmodhleq-1},
\[
 d_{\vv W}-3\dd_{\vv W}\stackrel{\eqref{eq:decompositioninequality}}{\geq} (d-3\dd)-(d_{/{\vv W}} -3\dd_{/{\vv W}})-1\stackrel{\eqref{eq:dmodhleq-1}}{\geq} d-3\dd.
\]

Since $\dd_{\vv W}\leq \dd-1$ by~\eqref{eq:dwleqd-1}, we can apply induction on ${\vv W}$, which certifies that~$\vv W$ contains at least $d_{\vv W}-3\dd_{\vv W}+1\geq d-3\dd+1$ disjoint positive circuits, and hence so does~${\vv V}$.
\end{proof}

Of course, this theorem is just a first step. It only proves that there is some subspace that contains many disjoint circuits, but ignores the vectors outside of this subspace, which could form more disjoint circuits. Some of the results presented in the next chapter suggest that it should be possible to improve on this in future work. 
Note that there is not even yet a linear bound for the Ehrhart-theoretical counterpart of this theorem (see statement \eqref{it:latticeCayley} in Section~\ref{sec:mainmotivation}).

\chapter{Configurations of degree 1}\label{ch:deg1}

For point configurations of degree~$1$, we can strengthen Theorem~\ref{thm:d+1-3dd} and provide their full classification. The main goal of this chapter is to prove the following result.

\begin{theorem}\label{thm:dd=1}
 For any $d$-dimensional point configuration~${\vv A}$, $\degc({\vv A}) \leq 1$ if and only if one of the following two conditions hold:
\begin{enumerate}
 \item ${\vv A}$ is a $k$-fold pyramid over a $2$-dimensional point configuration without interior points (up to repeated points); or
 \item ${\vv A}$ is a weak Cayley configuration of length $d$.
\end{enumerate}
\end{theorem}
This result, together with Corollary~\ref{cor:deg0}, provides yet more motivation for Conjecture~\ref{conj:d+1-2dd}:

\begin{corollary}\label{cor:deg1}
 Any $d$-dimensional point configuration ${\vv A}$ with $\degc({\vv A}) \leq 1$ is a weak Cayley configuration of length at least $d+1-2\degc(\vv A)$.\qed
\end{corollary}
 
Let us point out that the dimension of each factor of a weak Cayley configuration of length $d$ cannot be greater than $1$, since all factors are included in a flag of faces of length $d-1$. 
As a result:

\begin{proposition}\label{prop:dd=1:d=3}
If $d\geq 3$, the only weak Cayley configurations of length~$d$ are (up to repeated points)
\begin{itemize}
\item either $k$-fold pyramids over prisms over simplices with extra
  points on the ``vertical'' edges (in which case $\vv A_0=\emptyset$,
  so ${\vv A}$ is a Cayley configuration of length $d$);
\item or simplices with a vertex $\vv a$ and points on the edges adjacent to $\vv a$
  (here, $\vv A_0 = \{\aff(\vv a)\cap \vv A\}$ and ${\vv A}/\vv A_0$ is the vertex set of a simplex, a Cayley
  configuration of length~$d$).
\end{itemize}
\end{proposition}
Here, by a \defn{vertical}\index{vertical edge} edge of $\simp{d-1}\times\simp{1}$ we mean an edge of the form $\{\vv a\}\times\simp{1}$ for some vertex $\vv a$ of $\simp{d-1}$.

This proposition yields the final ingredient to recover the formulation of Theorem~\ref{thm:dd=1} presented in Chapter~\ref{ch:intro_degree}, we only need to observe that a $2$-dimensional point configuration $\vv A$ has degree $\deg(\vv A)\leq 1$ if and only if it does not have interior points.

\section{Lawrence polytopes}\label{sec:Lawrence}

Lawrence polytopes\index{Lawrence polytope} form a very interesting family of polytopes (cf.~\cite{BayerSturmfels1990}, \cite[Chapter~9]{OrientedMatroids1993}, \cite{Santos2002} or \cite[Chapter~6]{Ziegler1995}). A \defn{Lawrence polytope} is a polytope $\vv P$ with a centrally symmetric Gale dual~$\vv V=\Gale{\vv P}$. That is, $-\vv V=\vv V$ (as a multiset).

In Example~\ref{ex:prism} we introduced a particular instance of a Lawrence polytope. Namely, we saw that the Gale dual of a prism over a $d$-simplex is the vector configuration $\{\pm \vv e_1,\dots,\pm \vv e_{d}, \pm \sum_{i=1}^d \vv e_i\}\subset\RR^{d}$, which is centrally symmetric.

More generally, pairs of vectors of the form $\{\vv v,-\gl \vv v\}$ for some $\gl>0$ are called \defn{antipodal}\index{antipodal}. They have a nice behavior with respect to the degree because each hyperplane in general position contains exactly one vector of such a pair in its positive side. 

Observe how an irreducible vector configuration is combinatorially equivalent to the Gale dual of a Lawrence polytope precisely when it is a union of antipodal pairs of vectors. In this direction, the following proposition shows that irreducible Lawrence polytopes can be also characterized in terms of their extremal degree. Recall that Proposition~\ref{prop:easybound} stated that every irreducible vector configuration of rank $r$, $r+d+1$ elements and degree $\dd$ fulfills $r\geq d+1-2\dd$; Lawrence polytopes are those that attain the equality.    

\begin{proposition}\label{prop:Lawrence}
 An irreducible vector configuration $\vv V$ of rank $r$, $n=r+d+1$ elements and degree $\dd$ fulfills $r=d+1-2\dd$ if and only if $\vv V$ is centrally symmetric (up to rescaling).
\end{proposition}
\begin{proof}
Observe that if $\vv V$ is centrally symmetric and $\veczero\notin\vv V$, then every hyperplane ${\vvh H}$ in general position contains exactly $\frac{n}{2}$ elements of $\vv V$ in ${\vvh H}^+$. Therefore $\dd=\frac{n}{2}-r=\frac{d+1-r}{2}$.

To prove the converse, we will see that $\vv W:=\lin (\vv v)\cap \vv V$ is centrally symmetric for each $\vv v\in \vv V$.
Let $d_{\vv W}+2$ be the number of elements of $\vv W$ and $\dd_{\vv W}$ its degree. And let $\dd_{/\vv W}$ be the degree of $\vv V/\vv W$, and $(r-1)+d_{/\vv W}+1$ its number of elements. 

By Proposition~\ref{prop:easybound}, $(d_{/\vv  W}-2\dd_{/\vv W})\leq r-2$, and applying Proposition~\ref{prop:subspacequotient} we get that
 $d_{\vv  W}-2\dd_{\vv  W}\geq (d-2\dd)-(d_{/\vv  W}-2\dd_{/\vv W})-1\geq (d-2\dd)-(r-1)=0$.
Moreover, again by Proposition~\ref{prop:easybound}, $d_{\vv W}-2\dd_{\vv  W}\leq 0$. Therefore $d_{\vv  W}=2\dd_{\vv  W}$, and it is easy to check that any configuration of rank~$1$ fulfilling $d_{\vv  W}=2\dd_{\vv  W}$ must be centrally symmetric (again, up to rescaling).
\end{proof}
Corollary~\ref{cor:Lawrence} in Chapter~\ref{ch:conjecture} presents an analogous characterization of Lawrence polytopes in terms of the related concept of covector discrepancy. This allows for a reformulation of this proposition that avoids the hypothesis of irreducibility.

\section{A complete classification}\label{sec:deg1}

\subsection{\texorpdfstring{Circuits in configurations of degree~$1$}{Circuits in configurations of degree 1}}

In order to prove Theorem~\ref{thm:dd=1}, we need the following crucial result about circuits in vector configurations of dual degree~$1$. It states that in a pure vector configuration of dual degree~$1$ all small circuits are positive (or negative).

\begin{proposition}\label{prop:nosmallcircuits}
 Let ${\vv V}$ be a pure vector configuration of rank~$r$ with $\degG({\vv V})=1$. If $C$ is a circuit of $\cM({\vv V})$ with $|C^+|>0$ and $|C^-|>0$, then $|C|=r+1$. 
\end{proposition}

\begin{proof}
Consider ${\vv W}={\vv V}\cap \lin (C)$. By construction, $\rank(\vv W)=|C|-1$. If $|C^+|>0$ and $|C^-|>0$, there is a hyperplane ${\vvh H}$ in $\lin(C)$ with $C\subset {\vvh H}^+$. Indeed, by the Farkas Lemma (see \cite[Section 1.4]{Ziegler1995}), if there is no such  hyperplane, then $C$ must be a positive circuit. Therefore $\degG({\vv W})\geq 1$ because $|{\vvh H}^+\cap\vv W|\geq |{\vvh H}^+\cap C|=\rank(\vv W)+1$. Since ${\vv V}$ is pure, ${\vv V}/{\vv W}$ is also pure by Lemma~\ref{lem:quotientsofpurearepure}. If moreover $|C|\leq r$, then $\rank({\vv V}/{\vv W})\geq 1$, and by Lemma~\ref{lem:puredeggeq1}, $\degG({\vv V}/{\vv W})\geq 1$. Now, Proposition~\ref{prop:subspacequotient} implies that $\degG({\vv V})\geq \degG({\vv W})+\degG({\vv V}/{\vv W})\ge2$, which contradicts the hypothesis that $\degG({\vv V})=1$.
\end{proof}

We deduce some useful corollaries:

\begin{corollary}\label{cor:norepeatedpuredeg1}
If $\vv V$ is a pure vector configuration of rank~$r\ge2$ and $\degG({\vv V})=1$, then it has no repeated vectors except for, perhaps, the zero vector. 
\end{corollary}
\begin{proof}
  By Proposition~\ref{prop:nosmallcircuits}, any circuit with non-empty positive and negative part has size $r+1\geq 3$.
\end{proof}

\begin{corollary}\label{cor:disjointcircuits}
 Let ${\vv V}$ be a pure $r$-dimensional vector configuration with $\degG({\vv V})=1$. If $C\ne D$ are circuits of $\cM({\vv V})$ with $|C\cup D|\leq r+1$, then $C \cap D=\emptyset$.
\end{corollary}
\begin{proof}
Since $C$~and~$D$ are minimal by definition, there must exist $c\in C\setminus D$ and $d\in D\setminus C$. Therefore, $|C|\leq r$ and $|D|\leq r$ and, by Proposition~\ref{prop:nosmallcircuits}, both $C$ and $D$ may be assumed to be positive circuits. 
Suppose there also exists some $p\in C \cap D$. 
Eliminating $p$ on $C$ and $-D$ by oriented matroid circuit elimination, we find a circuit $E$ with $c\in E^+$, $d\in E^-$ of size $|E|\leq |C \cup D|-1 \leq r$. This contradicts Proposition~\ref{prop:nosmallcircuits}.
\end{proof}

Another interesting  consequence is that the factors of a weak \CayleyG configuration of dual degree~$1$ are its only small circuits.
\begin{lemma}\label{lem:onlysmallcircuits}
Let ${\vv V}$ be a pure vector configuration of rank~$r$ with $r+d+1$ elements, $d\ge2$ and $\degG({\vv V})=1$. If ${\vv V}$ is a weak \CayleyG configuration of length $d$ with factors $C_1,\dots,C_d$, and $D$ is a circuit of $\cM({\vv V})$ with $|D|\leq r$, then $D=C_i$ for some $1\leq i\leq d$.
\end{lemma}
\begin{proof}
Assume that $D\neq C_i$ for all $1\leq i\leq d$. If there is some $C_i$ with $|C_i\cap D|= |C_i|-1$, then $|C_i\cup D|\leq r+1$ and we get a contradiction to Corollary~\ref{cor:disjointcircuits}.
Otherwise, if $|C_i\cap D|\leq |C_i|-2$ for all $i$ and $|C_j\cap D|\neq \emptyset$ for some $j$, then 
\begin{align*}
  |C_j\cup D|&\leq n-\left|\bigcup_{i\neq j} \left(C_i\setminus D\right)\right|= n-\sum_{i\neq j}  \left| C_i\right|-\left| C_i\setminus D\right|\leq
 n-(d-1)2\\&=r+d+1-2d+2=r-d+3\leq r+1,
\end{align*}
and we again get a contradiction to Corollary~\ref{cor:disjointcircuits}. Hence, $D$ does not intersect any $C_1, \ldots, C_d$. 
By Proposition~\ref{prop:nosmallcircuits}, $D$ can be assumed to be a positive circuit. Therefore, ${\vv V}$ is a weak \CayleyG configuration of length $d+1$, so Proposition~\ref{prop:deg-cayley} implies that ${\vv V}$ has dual degree $0$, a contradiction.
\end{proof}

In particular, in the situation of the previous lemma any subset ${\vv W}\subset {\vv V}$ with $|{\vv W}|\leq r$ that does not contain any $C_i$ must be linearly independent.\\

Finally, we state another easy consequence of the Farkas Lemma (see \cite[Section 1.4]{Ziegler1995}) whose proof we leave to the reader.
\begin{lemma}\label{lem:poscircuithyperplane}
 Let $C$ be a positive circuit of a vector configuration $\vv V$, and let $\vvh H$ be a hyperplane. If $C\not\subset\vvh H$, then $|\vvh H^+\cap C|\geq 1$ and $|\vvh H^-\cap C|\geq 1$.\qed
\end{lemma}

\subsection{The classification}
We finally have all the tools needed to prove the following proposition, which directly implies Theorem~\ref{thm:dd=1}.

\begin{proposition}\label{prop:deg1}
 Let ${\vv V}$ be an irreducible pure vector configuration in $\RR^r$ with $n=r+d+1$ elements and $d\ge3$. If $\degG({\vv V})=1$, then ${\vv V}$ is a weak \CayleyG configuration of length $d$.
\end{proposition}

\begin{proof}
We fix $d\ge3$ and use induction on $r$. By Proposition~\ref{prop:easybound}, $r\geq d-1$, and our base case is $r=d-1$. Proposition~\ref{prop:Lawrence} tells us that $r=d-1$ if and only if ${\vv V}$ is centrally symmetric (up to rescaling). Observe that each of the pairs of antipodal vectors forms a circuit, and hence ${\vv V}$ is a \CayleyG configuration of length~$d$.

If $r>d-1$ $\vv V$ cannot be centrally symmetric by Proposition~\ref{prop:Lawrence}. Hence, there is some $\vv v\in {\vv V}$ such that ${\vv V}\cap \lin (\vv v)$ is not centrally symmetric. Since ${\vv V}$ does not have multiple vectors by Corollary~\ref{cor:norepeatedpuredeg1}, $\lin(\vv v)\cap \vv V=\vv v$, which is a configuration of rank~$1$, $1$ element and $\degG(\vv v)=0$. 
By Proposition~\ref{prop:subspacequotient} we know that $\degG(\vv V/\vv v)\leq \degG(\vv V)-\degG(\vv v)=1$, and by Lemma~\ref{lem:puredeggeq1} that $\degG(\vv V/\vv v)\ge1$. Combining these inequalities we see that $\degG({\vv V}/\vv v)=1$.
Therefore, ${\vv V}/\vv v$ is a vector configuration of dual degree~$1$ that is irreducible (Lemma~\ref{lem:quotientirreducible}), pure (Lemma~\ref{lem:quotientsofpurearepure}), and has rank~$r-1$ and $(r-1)+d+1$ elements. By the induction hypothesis, ${\vv V}/\vv v$ is therefore a weak \CayleyG configuration with factors $\tilde C_1, \dots, \tilde C_d$, say.
For convenience, we define $\tilde C_0:=\left({\vv V}/\vv v\right)\setminus\bigcup\nolimits_{i=1}^d \tilde C_i$. 

By counting the number of elements in $|{\vv V}/\vv v|$, we see that 
\begin{equation}\label{eq:countingcircuitelements} \sum_{i=0}^d |\tilde C_i|=|\vv V /\vv v|=|\vv V|-1=r+d.
\end{equation} 
After subtracting $2d$ from both sides, $|\tilde C_0|\ge0$ implies that 
\[\sum_{i=1}^d{\big( |\tilde C_i|-2\big)}\leq r-d;
\]
in particular, $|\tilde C_i|\leq r-1$ for all $1\leq i\leq d$ because $d\ge3$ and $|\tilde C_j|\geq 2$ for all $j$. 

For each $1 \le i \le d$, $\tilde C_i$ is a positive circuit of $\cM({\vv V}/\vv v)$ that expands to a circuit $C_i$ of $\cM({\vv V})$ 
(see Section~\ref{sec:operations}). From now on, we consider subsets of~${\vv V}/\vv v$ as subsets of~${\vv V} \setminus \vv v$ by identifying corresponding elements, so that $\tilde C_i = C_i \setminus \{\vv v\}$. 

Since $|C_i|\leq |\tilde C_i|+1\leq r$, Proposition~\ref{prop:nosmallcircuits} shows that either $\vv v\notin C_i$ or $\vv v\in C_i^+$. Hence, $C_i$ is again a positive circuit with either $C_i^+=\tilde C_i^+$ or $C_i^+=\tilde C_i^+\cup \{\vv v\}$. We will show that if some $C_i$ contains $\vv v$, no other $C_j$ can. This will prove our claim because then $C_1,\dots,C_i,\dots,C_d$ are disjoint positive circuits that instate ${\vv V}$ as a weak \CayleyG configuration. 

For this, we assume that $\vv v\in C_1\cap C_2$ and reach a contradiction. We start with some definitions.
For $1\leq i\leq d$, let ${D}_i$ be a subset of $|\tilde C_i|-2$ elements of~$\tilde C_i$, and set ${D}:=\tilde C_0 \cup \bigcup\nolimits_{i=1}^d {D}_i$. Next, choose $\vv v_1\in \tilde C_1 \setminus {D}_1$ and $\vv v_2\in \tilde C_2 \setminus {D}_2$ (so that, in particular, $\vv v\notin\{\vv v_1,\vv v_2\}$) and define ${D}':={D}\cup\{\vv v_1,\vv v_2\}$. A diagram to follow the definition of all these sets is depicted in Figure~\ref{fig:setjungle}.

\iftoggle{bwprint}{%
\begin{figure}[htpb]
\centering
\includegraphics[width=.6\linewidth]{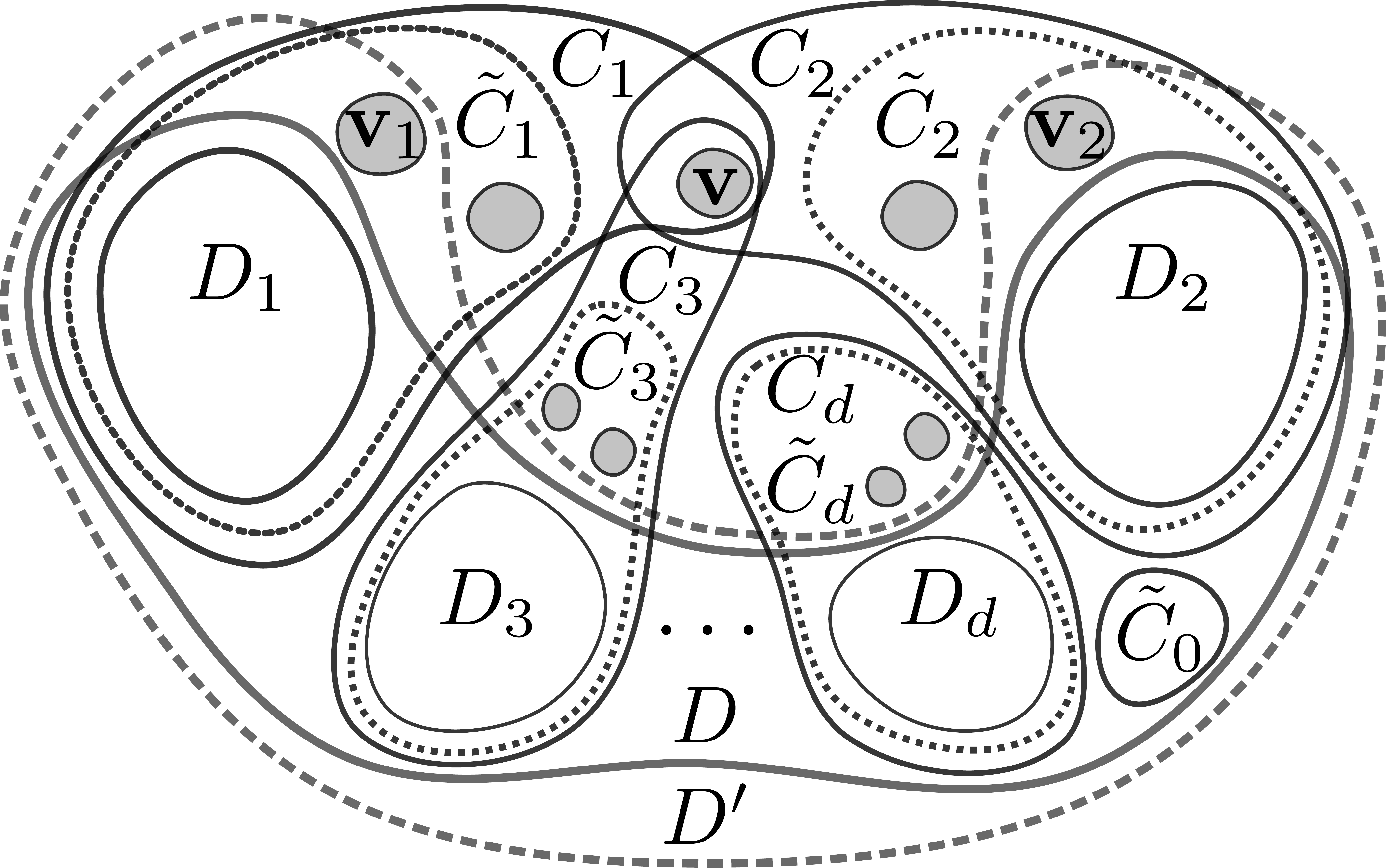}
  \caption[Sketch of the sets involved in the proof of Proposition~\ref{prop:deg1}.]{Sketch of the sets involved in the proof of Proposition~\ref{prop:deg1}. Solid dots represent elements (as the ones labeled as $\vv v$, $\vv v_1$ and $\vv v_2$). For $i>2$ some $C_i$ might contain $\vv v$ and some not; in this example, $\vv v\in C_3$ but $\vv v\notin C_d$.}
 \label{fig:setjungle}
\end{figure}
}{%
\begin{figure}[htpb]
\centering
\includegraphics[width=.6\linewidth]{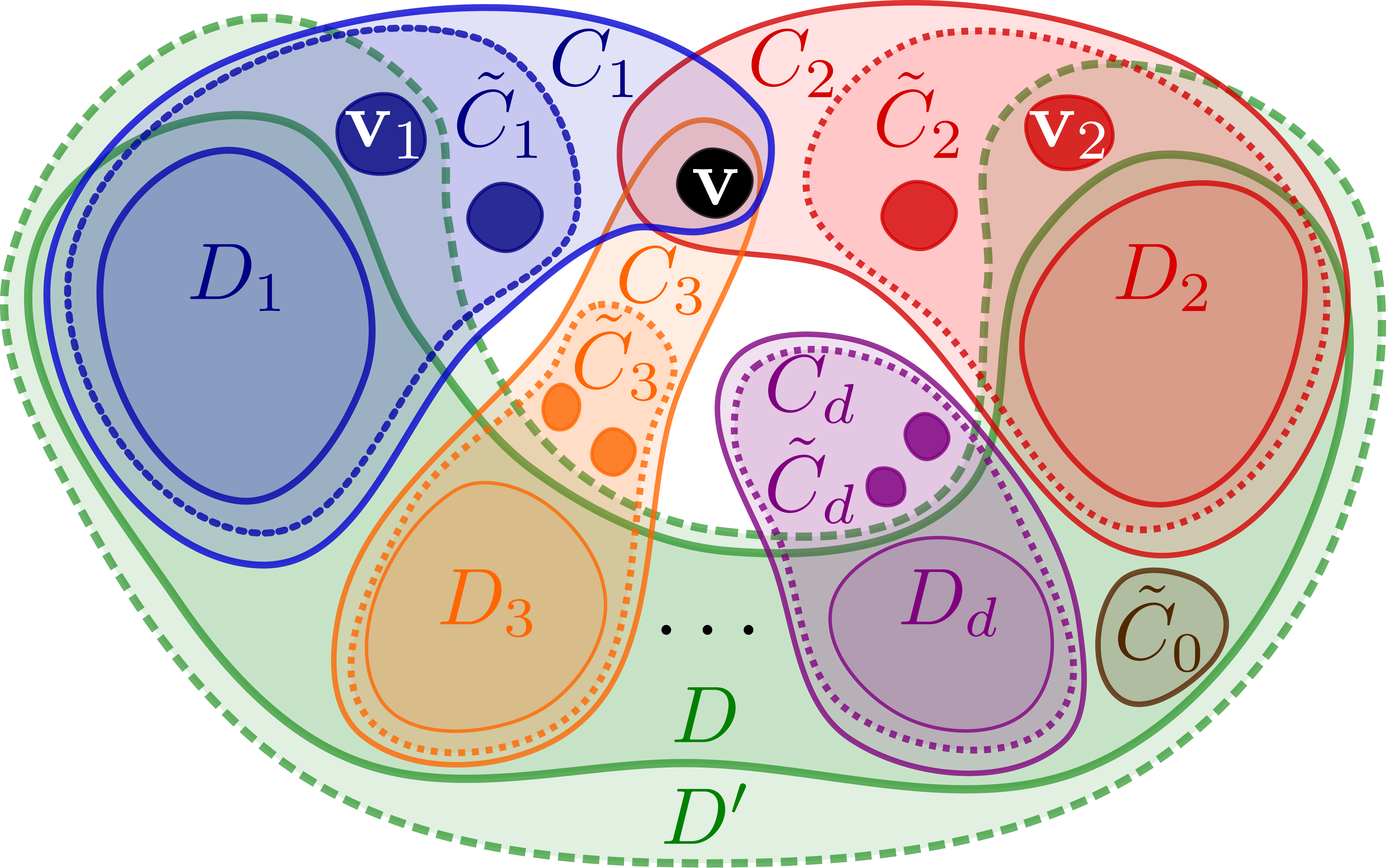}
  \caption[Sketch of the sets involved in the proof of Proposition~\ref{prop:deg1}.]{Sketch of the sets involved in the proof of Proposition~\ref{prop:deg1}. Solid dots represent elements (as the ones labeled as $\vv v$, $\vv v_1$ and $\vv v_2$). For $i>2$ some $C_i$ might contain $\vv v$ and some not; in this example, $\vv v\in C_3$ but $\vv v\notin C_d$.}
 \label{fig:setjungle}
\end{figure}
}

A first observation is that the elements in ${D}'$ must be linearly independent. Indeed, since 
\begin{eqnarray*}
|{D}'|
&=&
2+|\tilde C_0|+\sum_{i=1}^d \big(|\tilde C_i|-2\big)
\\ &\stackrel{\eqref{eq:countingcircuitelements}}{=}&
2 + (r+d) - 2d 
\ = \ 
r+2-d
\ \leq \ 
r-1,
\end{eqnarray*}
already their projections to ${\vv V}/\vv v$ are linearly independent. The reason for this is that if the elements in $D'/\vv v$ were not linearly independent then they would contain a circuit. But this contradicts Lemma~\ref{lem:onlysmallcircuits} because $D'\not\supseteq \tilde C_i$ for all $i$, since by construction $|\tilde C_i \setminus D'|\ge1$ for all~$i$. Now, let ${\vvh H}'$ be a hyperplane through $\lin({D}')$ that is otherwise in general position with respect to~${\vv V}$. This is possible because the rank of~$\vv V$ is~$r$, and $D'$~has at most $r-1$ elements. Observe that $\vv v\notin\lin({D}')$, because otherwise the vectors in ${D}'$ would form a circuit in ${\vv V}/\vv v$. Therefore, $\vv v\notin {\vvh H}'$, and we can orient~${\vvh H}'$ so that $\vv v\in {{\vvh H}'}^-$. Then $|{{\vvh H}'}^+\cap \tilde C_i|=|{{\vvh H}'}^+\cap C_i|=1$ for $i=1,2$ because of Lemma~\ref{lem:poscircuithyperplane} and our assumption that $\vv v\in C_1\cap C_2$.
 Moreover, since the elements in ${D}'$ are linearly independent, we can perturb~${\vvh H}'$ to a hyperplane~${\vvh H}$ through $\lin({D})$ such that $\vv v_1,\vv v_2\in {\vvh H}^+$. This yields 
\[
\big|{{\vvh H}}^+\cap \tilde C_i\big|
=
\big|({{\vvh H}'}^+\cap \tilde C_i)\cup \vv v_i\big|
= 2
\qquad\text{for } i=1,2.
\]

Furthermore, we claim that $|{\vvh H}^+\cap \tilde C_j|\geq 1$ for all $j\ge3$. 
If, on the contrary, there existed some $j\ge3$ with $|{\vvh H}^+\cap \tilde C_j|=0$,  Lemma~\ref{lem:poscircuithyperplane} would yield~$\vv v\notin C_j$ (\ie $C_j=\tilde C_j$), and moreover $C_j$ would be completely contained in~${\vvh H}$. Hence, by construction, $C_j$ would be completely contained in $\lin({D})$. In particular, some $\vv v_j\in \tilde C_j \setminus {D}_j$ would satisfy $\vv v_j\notin {D}$ but $\vv v_j\in \lin ({D})$. Therefore, this element would be part of a circuit in~$\{\vv v_j\}\cup {D}$, distinct from~$C_j$ since $|C_j\cap {D}|=C_j-2$. However, $|C_j\cup {D}|\leq |{D}|+3\leq r$, which would contradict Corollary~\ref{cor:disjointcircuits}.

Finally, let ${\vvh H}''$ be any hyperplane such that ${D}\subset {{\vvh H}''}^+$. Now
\begin{itemize}
\item $\big|({\vvh H}\circ {\vvh H}'')^+ \cap \tilde C_0\big|=|\tilde C_0|$;
\item $\big|({\vvh H}\circ {\vvh H}'')^+ \cap \tilde C_i\big|=|\tilde C_i|$ for $i=1,2$; and
\item $\big|({\vvh H}\circ {\vvh H}'')^+ \cap \tilde C_j\big|\geq |\tilde C_j|-1$ for $3\le j\leq d$. 
\end{itemize}

Therefore, using \eqref{eq:countingcircuitelements} we see that
\[
 \big|({\vvh H}\circ {\vvh H}'')^+ \cap {\vv V}\big|\ \geq \  \sum_{i=0}^d |\tilde C_i|-(d-2)=r+2,
\]
which contradicts $\degG({\vv V})=1$.
\end{proof} 

\section{Totally splittable polytopes}

A \defn{split}\index{point configuration!split} of a point configuration $\vv A$ is a (necessarily regular) polyhedral subdivision of $\vv A$ with exactly two maximal cells, which are separated by a \defn{split hyperplane}. A point configuration $\vv A$ is called \defn{totally splittable}\index{point configuration!totally splittable} if every triangulation of $\vv A$ is a common refinement of splits (see~\cite{HerrmannJoswig2010}). 

A polytope $\vv P$ is totally splittable if $\verts(\vv P)$ is a totally splittable point configuration. In~\cite[Theorem 9]{HerrmannJoswig2010}, Herrmann and Joswig establish a complete classification of totally splittable polytopes: simplices, polygons, prisms over simplices, crosspolytopes and  (possibly multiple) joins of these.

From this classification and Theorem~\ref{thm:dd=1}, it follows that every polytope of degree $\leq 1$ is totally splittable. This section aims to shed more light on this relationship. We also answer an open question from~\cite{HerrmannJoswig2010} by giving a combinatorial proof of the fact that totally splittable polytopes are equidecomposable~\cite[Corollary 30]{HerrmannJoswig2010}.

Two splits of $\vv A$ are called \defn{compatible} if their split hyperplanes do not intersect in the interior of $\vv A$. We will say that $\vv A$ is \defn{strongly totally splittable} if every triangulation of~$\vv A$ is a common refinement of compatible splits. 

\begin{proposition}
 $\vv A$ is strongly totally splittable if and only if $\degc(\vv A)\leq 1$.
\end{proposition}
\begin{proof}
 Let $\vv A$ be a $d$-dimensional point configuration of degree $\degc(\vv A)\leq 1$, and let $\cT$ be a triangulation of $\vv A$. Since $\vv A$ has no interior $(d-2)$-faces, the boundary of every $(d-1)$-face of $\cT$ lies in the boundary of $\vv A$. It is easy to see that such a $(d-1)$-face always defines a split of $\vv A$, for example by using \cite[Observation 3.1]{HerrmannJoswig2008}. Since these faces do not intersect in the interior of~$\vv A$, the splits are compatible, and $\cT$ is their common refinement.

 Analogously, if every triangulation of $\vv A$ is a common refinement of compatible splits, $\vv A$ cannot have interior $(d-2)$-faces. Otherwise a triangulation using such a face must involve splits that intersect in that face, and hence, not compatible.
\end{proof}

As a corollary, every polytope of degree~$1$ is totally splittable. In particular, by analyzing the cases of Herrmann and Joswig's result one could deduce an alternative proof of Theorem~\ref{thm:dd=1} for the case that the points in~${\vv A}$ are in convex position.\\

A point configuration $\vv A$ is \defn{equidecomposable}\index{point configuration!equidecomposable} if all its triangulations share the same $f$-vector. Looking at the classification of totally splittable polytopes, Sanyal observed that they are all equidecomposable~\cite[Corollary 30]{HerrmannJoswig2010}. In that paper, the authors ask for a classification-free proof of this result. Here we provide a proof based on the following characterization, proved in the book by De Loera, Rambau and Santos~\cite{DeLoeraRambauSantosBOOK}; the necessity part had already been found by Bayer in~\cite{Bayer1993}.

\begin{theorem}[{\cite[Theorem 8.5.19]{DeLoeraRambauSantosBOOK}}]\label{thm:equidecomposable}
 $\vv A$ is equidecomposable if and only if $|X^+|=|X^-|$ for every circuit $X$ of $\vv A$.\qed
\end{theorem}

The following property of totally splittable configurations is new.

\begin{proposition}
 If a point configuration $\vv A$ is totally splittable, then $|X^+|\leq 2$ for every circuit $X$ of $\vv A$.
\end{proposition}
\begin{proof}
 The proof is by induction on $r$, the rank of $\Gale{\vv A}$. Observe that $r=0$ if and only if $\vv A$ is the vertex set of a simplex and that the result is easy to prove for $r=1$, i.e., for totally splittable $d$-dimensional configurations of $d+2$ points.

 If $r>1$, since every circuit $X$ of $\vv A$ involves at most $d+2$ points, there is some $\vv a\in \vv A$ with $\vv a\notin \underline X$. Therefore $X$ is a circuit of $\vv A\setminus \vv a$.  It is easy to see that if $\vv A$ is totally splittable then $\vv A\setminus \vv a$ must also be (see~\cite[Proposition 13]{HerrmannJoswig2010}). Hence, by induction $|X^+|\leq 2$.
\end{proof}

This proposition explains the equidecomposability of totally splittable polytopes.

\begin{corollary}
 Every totally splittable polytope is equidecomposable.
\end{corollary}
\begin{proof}
 Let $\vv P$ be a totally splittable polytope. By Theorem~\ref{thm:galepolytope}, every circuit~$X$ of $\vv P$ satisfies $|X^+|\geq 2$, because $\vv P$ is a polytope. Since $|X^+|\leq 2$ by the previous proposition, every circuit $X$ fulfills $|X^+|=|X^-|=2$, and~$\vv P$ is equidecomposable by Theorem~\ref{thm:equidecomposable}.
\end{proof}

And raises the following question.

\begin{question}
For a point configuration~$\vv A$, is it enough to check that $|X^+|\leq 2$ for every circuit~$X$ of $\cM(\vv A)$ to prove that $\vv A$ is totally splittable?
\end{question}

\chapter{Codegree decompositions} \label{ch:conjecture}

In Theorem~\ref{thm:d+1-3dd} we proved that any $d$-dimensional point configuration of degree $\leq\dd$ is a weak Cayley configuration of length at least $d+1-3\dd$, and our Conjecture~\ref{conj:d+1-2dd} suggests that this bound could be improved to $d+1-2\dd$. Furthermore, in Example~\ref{ex:sharpconj} we saw that the conjectured bound, if true, is tight.

However, being a weak Cayley configuration of length $d+1-2\dd$ can only guarantee a degree $\leq 2\dd$ (Proposition~\ref{prop:deg-cayley}). This bound on the degree is also tight, in the sense that it is not hard to find examples of weak Cayley configurations of maximal length $k$ and degree precisely $d+1-k$: For example, the vertices of a $d$-dimensional crosspolytope $\cros{d}$ form a weak Cayley configuration of maximal length $2$ and have degree $d-1$.

This shows that the concept of weak Cayley configuration is not strong enough to characterize point configurations of degree $<\fceil{d}{2}$. However, there are other strategies that, for some fixed $k<\fceil{d}{2}$, allow to construct point configurations of degree $\leq k$. We present three such strategies, which in combination give rise to the concept of \defn{codegree decompositions}. 

\section{Codegree decompositions}
The first construction is the \defn{join} of point configurations (see Section~\ref{sec:operations}).
It is easy to see that the join of an interior face of $\vv A$ with an interior face of~$\vv B$ is an interior face of $\vv A\join \vv B$, and that all interior faces of $\vv A\join \vv B$ arise this way. Thus,
\begin{lemma}\label{lem:codegjoin}
 Let $\vv A$ and $\vv B$ be point configurations, then $\degc(\vv A \join \vv B)=\degc(\vv A)+ \degc(\vv B)$ and $\codegc(\vv A \join \vv B)=\codegc(\vv A)+ \codegc(\vv B)$.\qed
\end{lemma}

The second construction are \defn{liftings}, understood as the inverse of affine projections. The key observation here is that if $\vv S$ is an interior face of $\vv A$, and $\pi$ is an affine projection, then $\pi(\vv S)$ is an interior face of $\pi(\vv A)$.

\begin{lemma}\label{lem:codegproj}
  Let $\vv A$ and $\vv B$ be point configurations. If $\pi(\vv A)=\vv B$ for some affine projection $\pi$, then $\codegc(\vv B)\leq \codegc(\vv A)$.\qed
\end{lemma}

Our last construction, the inverse operation of contraction, is strongly related to lifting. Given a point configuration~$\vv B$, we want to find a point configuration~$\vv A$ with an element~$\vv a$ such that the contraction~$\vv A/\vv a$ coincides with~$\vv B$. To do so, we just need to choose a lifting of~$\vv B$ and place a new point $\vv a$ high enough in the direction of the lifting (see Section~\ref{sec:descriptioncG} for an explicit example).  
The following result complements Corollary~\ref{cor:degpointdeletioncontraction}.

\begin{lemma}
If $\vv a\in \vv A$ then $\codegc(\vv A/\vv a)\leq \codegc(\vv A)$.
\end{lemma}
\begin{proof}
By Gale duality, this is equivalent to seeing that $\codegG(\vv V\setminus \vv v)\leq \codegG(\vv V)$ for each $\vv v\in \vv V$, which is direct by Definition~\ref{def:dualcodegree}:
\[
\codegG(\vv V\setminus \vv v)=\min_{\vvh H} \big|\ol {\vvh H}^-\cap (\vv V\setminus \vv v)\big|
\leq 
\min_{\vvh H} \big|\ol {\vvh H}^-\cap \vv V\big|
=\codegG(\vv V).\qedhere
\]
\end{proof}

The combination of these constructions allows us to build point configurations of large codegree and, hence, of small degree:

\begin{proposition}\label{prop:strongCayleycodeg}
  Let $\vv A$ be a point configuration in $\RR^d$. If there exists a subset~$\vv A_0\subseteq\vv A$ and a projection $\pi$ such that $\pi(\vv A/\vv A_0)=\vv B_1\join \vv B_2 \join \dots \join \vv B_m$, then $\codegc(\vv A)\geq\sum_{i=1}^m \codegc(\vv B_i)$.\qed
\end{proposition}
\medskip

In view of this proposition, we make the following definitions:

\begin{definition}\label{def:decomposable}
 A point configuration $\vv A$ admits an \defn{affine codegree decomposition}\index{codegree decomposition!affine} if it has a codegree preserving contraction $\vv A/\vv A_0$ such that there is a codegree preserving projection $\pi$ from $\vv A/\vv A_0$ to a join of $m$ polytopes:
\begin{itemize}
\item $\pi(\vv A/\vv A_0)=\vv B_1\join \vv B_2 \join \dots \join \vv B_m$, \ and 
\item $\codegc(\vv A)=\codegc(\vv A/\vv A_0)=\sum_{i=1}^m\codegc(\vv B_i)$. 
\end{itemize}
The sets $\vv A_i=\pi^{-1}(\vv B_i)$ are called the \defn{factors} of the decomposition, and $m$ its \defn{length}. 

Analogously, $\vv A$ admits a \defn{(combinatorial) codegree decomposition}\index{codegree decomposition} if there exist disjoint non-empty subsets $\vv A_1,\dots,\vv A_m$ of $\vv A$ such that
\begin{itemize} 
\item For $1\leq i\leq m$, $\vv A \setminus \vv A_i=\vv F\cap \vv A$ for some face $\vv F$ of $\conv (\vv A)$;  and 
\item $\codegc(\vv A)=\sum_{i=1}^m\codegc\big(\vv A/(\vv A \setminus \vv A_i)\big)$. 
\end{itemize}
The sets $\vv A_i$ are called the \defn{factors}\index{codegree decomposition!factors} of this decomposition and $m$ its \defn{length}\index{codegree decomposition!length}. 
\end{definition}

As an example, the reader is invited to find an affine codegree decomposition of length~$3$ of the weak Cayley configuration presented in Figure~\ref{fig:exweak}.

The definition of combinatorial codegree decomposition is motivated by Corollary~\ref{cor:definitionscombdecomposable} (which we prove later). It shows that a point configuration admits a combinatorial codegree decomposition if and only if it is combinatorially equivalent to a configuration that admits an affine codegree decomposition. 
Consequently, the generic term ``codegree decomposition'' is used to refer to combinatorial codegree decompositions.\\

We are finally able to state our main conjecture.

\begin{conjecture}\label{conj:strongd+1-2dd}
 Any point configuration~$\vv A$ of dimension $d>2\deg(\vv A)$ admits a codegree decomposition of length at least $d+1-2\deg(\vv A)$.
\end{conjecture}

Observe first that this conjecture is a converse to Proposition~\ref{prop:strongCayleycodeg}, since it states that, combinatorially, all point configurations of small degree are constructed as liftings of joins.

Moreover, it is very easy to check that if a projection of $\vv A/\vv A_0$ factors as $\pi(\vv A/\vv A_0)=\vv B_1\join \vv B_2 \join \dots \join \vv B_m$, then $\vv A$ is a weak Cayley configuration of length $m$ with factors $\pi^{-1}(\vv B_i)$. Just observe that each $\vv B_i$ can be independently projected to a $0$-dimensional point configuration. This proves the following lemma.
\begin{lemma}\label{lem:decompiscayley}
If $\vv A$~admits a codegree decomposition of length $m$ then $\vv A$ is a weak Cayley configuration of length $m$.\qed
\end{lemma}

Therefore, Conjecture~\ref{conj:d+1-2dd}, which stated that every $d$-dimensional point configuration $\vv A$ is a weak Cayley configuration of length at least $d+1-2\degc(\vv A)$,
is a direct consequence of Conjecture~\ref{conj:strongd+1-2dd}. The difference is that Conjecture~\ref{conj:d+1-2dd} only uses the fact that $\codegc(\vv B_i)>1$, while Conjecture~\ref{conj:strongd+1-2dd} adds structure to the factors by using the entire codegree.

\iftoggle{bwprint}{%
\begin{figure}[htpb]
\centering
 \subbottom[$\vv A$]{\label{sfig:pyramid1}\includegraphics[width=.2\linewidth]{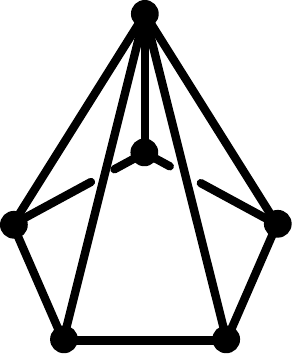}}\qquad\qquad\qquad
 \subbottom[$\vv B$]{\label{sfig:pyramid2}\includegraphics[width=.2\linewidth]{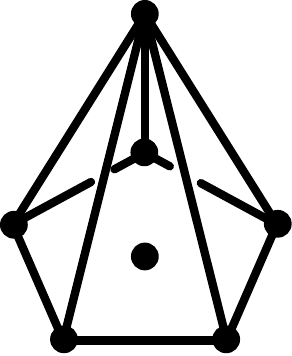}}
 \caption{Two configurations with a codegree decomposition of length $2$.}
\label{fig:exampleCodegreeDecomposition}
\end{figure}
}{%
\iftoggle{print}{%
\begin{figure}[htpb]
\centering
 \subbottom[$\vv A$]{\label{sfig:pyramid1}\includegraphics[width=.2\linewidth]{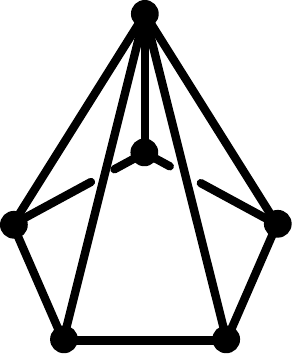}}\qquad\qquad\qquad
 \subbottom[$\vv B$]{\label{sfig:pyramid2}\includegraphics[width=.2\linewidth]{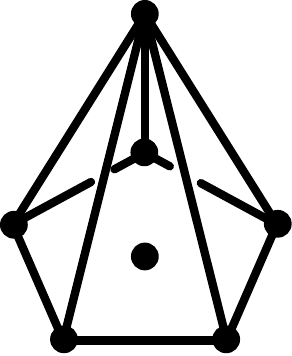}}
 \caption{Two configurations with a codegree decomposition of length $2$.}
\label{fig:exampleCodegreeDecomposition}
\end{figure}
}{%
\begin{figure}[htpb]
\centering
 \subbottom[$\vv A$]{\label{sfig:pyramid1}\includegraphics[width=.2\linewidth]{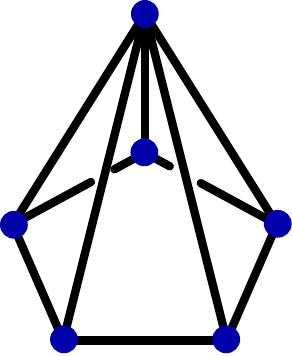}}\qquad\qquad\qquad
 \subbottom[$\vv B$]{\label{sfig:pyramid2}\includegraphics[width=.2\linewidth]{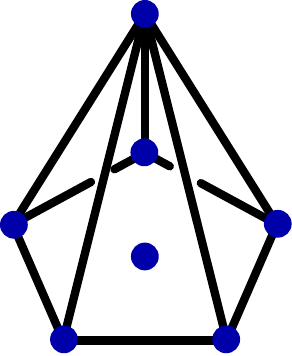}}
 \caption{Two configurations with a codegree decomposition of length $2$.}
\label{fig:exampleCodegreeDecomposition}
\end{figure}
}
}

\begin{example}
Figure~\ref{sfig:pyramid1} depicts a pyramid~$\vv A$ over the vertex set of a pentagon, while $\vv B$ adds an interior point to the base (Figure~\ref{sfig:pyramid2}).

Their properties compare as follows: 
\begin{center}\small
 \renewcommand{\arraystretch}{1.1}
  \begin{tabular}[c]{m{.55\textwidth}ccc}
\hline
    \centering Property && $\vv A$  & $\vv B$ \\
\hline
\quad    degree  && $1$ & $2$ \\
\quad    weak Cayley decomposition with 2 factors && apex $\uplus$ base & apex $\uplus$ base\\
\quad    weak Cayley decomposition with 3 factors && none & none \\
\quad    codegree decomposition with 2 factors && apex $\uplus$ base & apex $\uplus$ base\\
\quad    codegree decomposition with 3 factors && none & none \\
\quad    codegree of apex && $1$ & $1$ \\
\quad    codegree of base && $2$ & $1$ \\
\quad    codegree && $3$ & $2$\\
\hline
\end{tabular}
\end{center}

Observe that weak Cayley decompositions cannot explain why the
codegrees of $\vv A$ and $\vv B$ are different. However, this can easily be read off from the codegrees of the factors of the codegree decomposition.
\end{example}
\smallskip

We next provide a Gale dual interpretation of codegree decompositions. It is used in Section~\ref{sec:evidences} to 
gather evidence for Conjecture~\ref{conj:strongd+1-2dd} by showing, among other results, that any $d$-dimensional irreducible configuration of degree~$\dd$ with $r+d+1$ elements admits a codegree decomposition of length at least $2(d+1-2\dd)-(r-1)$ (Theorem~\ref{thm:2DD}). 

While we do not yet have any general bound that is independent of~$r$, we know that the conjectured bound holds for configurations of degree~$0$ or~$1$ (Propositions~\ref{prop:conjdd=0} and~\ref{prop:conjdd=1}) and for irreducible configurations with $r+d+1$ elements fulfilling $r-2\leq d+1-2\dd$ (Theorem~\ref{thm:DD=2}); which in particular settles the conjecture for configurations of corank $\leq 4$ (Corollary~\ref{cor:rank4}). These results are best explained using the concept of \defn{covector discrepancy}, which is introduced in Section~\ref{sec:evidences}.
 
As a final remark, in Section~\ref{sec:consequences} we show the strength of Conjecture~\ref{conj:strongd+1-2dd} by proving how, if it is was proved to be true, many of our results from previous chapters would follow directly.
  
\section{\texorpdfstring{Codegree\G decompositions}{Codegree* decompositions}}

Just as we did for Cayley configurations in Section~\ref{sec:DualCayley}, we also need a Gale dual interpretation of codegree decompositions. For this, recall from Definition~\ref{def:dualcodegree} that the dual codegree of a vector configuration is $\codegG(\vv V):=\min_{\vvh H} |\ol {\vvh H}^-\cap \vv V|$.

\begin{definition}\label{def:dualcodegreedecompositions}
 A vector configuration $\vv V$ admits a \defn{\codegreeG decomposition}\index{\codegreeG decomposition} of length $m$ if it admits a partition $\vv V=\vv V_0\uplus \vv V_1\uplus \dots \uplus \vv V_m$ such that \begin{equation}\label{eq:defcodegreecodecompositions}
\codegG(\vv V)=\sum_{i=1}^m \codegG(\vv V_i),
\end{equation}  
and where $\codegG(\vv V_i)\geq 1$ for $i=1,\dots,m$. The sets $\vv V_1,\dots, \vv V_m$ are the \defn{factors}\index{\codegreeG decomposition!factor} of the decomposition. 
The decomposition is \defn{affine}\index{\codegreeG decomposition!affine} if additionally \begin{equation}\label{eq:defaffinecodegreecodecompositions}
\sum_{\vv v_j\in\vv V_i}\vv v_j=\veczero \qquad\text{for each }1\leq i \leq m.
\end{equation}
\end{definition}

Figure~\ref{fig:exdualcodegree} shows a \codegreeG decomposition of length~$3$ of a vector configuration~$\vv V$ of dual codegree~$4$. It decomposes $\vv V$ into 
$\vv V_0$ (of dual codegree~$0$), $\vv V_1$ and $\vv V_2$ (of dual codegree~$1$), and $\vv V_3$ (of dual codegree~$2$). This decomposition is affine because the barycenter of each factor is the origin.
It is not the unique \codegreeG decomposition of $\vv V$, since the vector in~$\vv V_0$ could be swapped with a vector in $\vv V_3$ (though the decomposition would not be affine anymore).

Observe how this decomposition trivially implies that $\vv V$ is a weak \CayleyG configuration of length~$3$ (with factors $\vv V_1$, $\vv V_2$ and $\vv V_3$). By contrast, any of the decompositions that result from removing a vector from $\vv V_3$ and adding it to $\vv V_0$ are still weak \CayleyG decompositions of $\vv V$ but not \codegreeG decompositions. 

\iftoggle{bwprint}{%
\begin{figure}[htpb]
\centering
 \subbottom[$\vv V$]{\label{sfig:dualcodegreeV}\includegraphics[width=.2\linewidth]{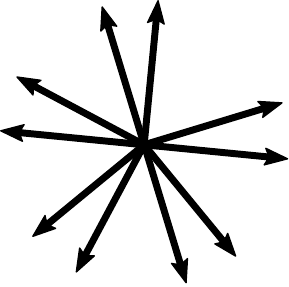}}\quad
 \subbottom[$\vv V_0$]{\label{sfig:dualcodegreeV0}\includegraphics[width=.15\linewidth]{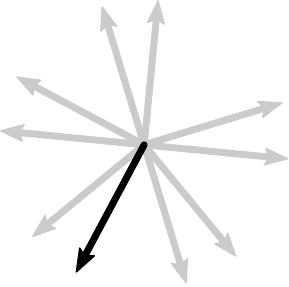}}\quad
 \subbottom[$\vv V_1$]{\label{sfig:dualcodegreeV1}\includegraphics[width=.15\linewidth]{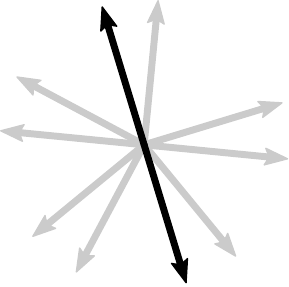}}\quad
 \subbottom[$\vv V_2$]{\label{sfig:dualcodegreeV2}\includegraphics[width=.15\linewidth]{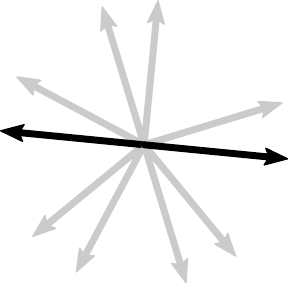}}\quad
 \subbottom[$\vv V_3$]{\label{sfig:dualcodegreeV3}\includegraphics[width=.15\linewidth]{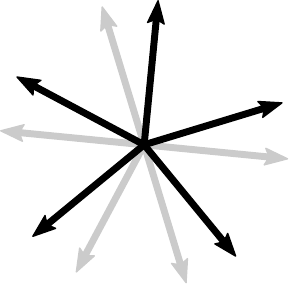}}
  \caption{A \codegreeG decomposition of $\vv V$ of length~$3$ with two factors of dual codegree~$1$ and one of dual~codegree $2$.}
 \label{fig:exdualcodegree}
\end{figure}
}{%
\begin{figure}[htpb]
\centering
 \subbottom[$\vv V$]{\label{sfig:dualcodegreeV}\includegraphics[width=.2\linewidth]{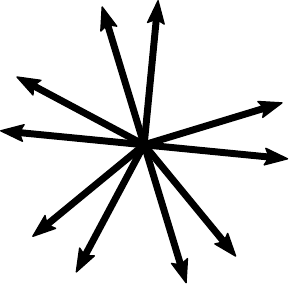}}\quad
 \subbottom[$\vv V_0$]{\label{sfig:dualcodegreeV0}\includegraphics[width=.15\linewidth]{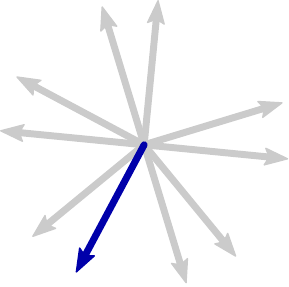}}\quad
 \subbottom[$\vv V_1$]{\label{sfig:dualcodegreeV1}\includegraphics[width=.15\linewidth]{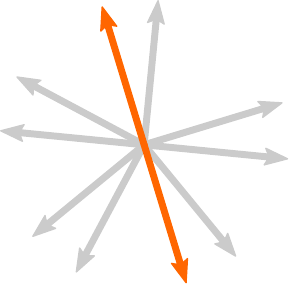}}\quad
 \subbottom[$\vv V_2$]{\label{sfig:dualcodegreeV2}\includegraphics[width=.15\linewidth]{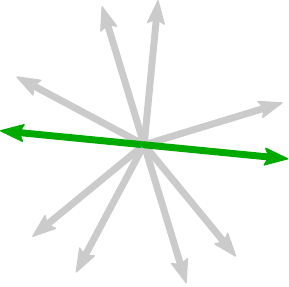}}\quad
 \subbottom[$\vv V_3$]{\label{sfig:dualcodegreeV3}\includegraphics[width=.15\linewidth]{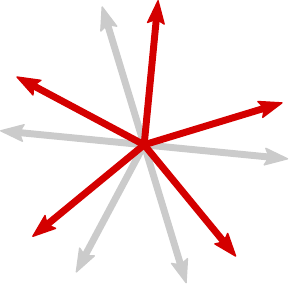}}
  \caption{A \codegreeG decomposition of $\vv V$ of length~$3$ with two factors of dual codegree~$1$ and one of dual~codegree $2$.}
 \label{fig:exdualcodegree}
\end{figure}
}

\begin{remark}
In Definition~\ref{def:dualcodegreedecompositions} we impose no condition on~$\vv V_0$, but in fact $\codegG(\vv V_0)=0$ by the upcoming Proposition~\ref{prop:sumcodegrees}, because $\codegG(\vv V)\geq \sum_{i=0}^m \codegG(\vv V_i)$.
\end{remark}

\begin{remark}
We do not consider $\vv V_0$ to be a factor of a \codegreeG decomposition in order to be consistent with the primal definitions of codegree decomposition and weak Cayley configuration.

On the one hand, $\vv A_0$, the primal counterpart of $\vv V_0$, behaves slightly differently from the other factors in the definition of codegree decomposition. If $\vv A_0=\emptyset$, one can define a codegree decomposition using liftings and joins (and a combinatorial equivalence). 
However, $\vv A_0$ must be treated differently and imposes the introduction of contractions into the definition. The underlying reason for this is that the codegree of a point configuration is always positive. In the dual picture this is much neater, because there are vector configurations of dual codegree~$0$ (they are precisely the acyclic ones), and hence there would be no need to explicitly differentiate $\vv V_0$ from the other factors in the dual definition.

On the other hand, we only require $\vv A\setminus\vv A_i$ to be the set of points on a face of $\conv(\vv A)$ when $i>0$. This implies that if $\vv A$ admits a codegree decomposition with $m$~factors then it must be a weak Cayley configuration of length $m$. However, if $\vv A_0$ were considered to be a factor, it should be treated differently from the factors of positive codegree when computing the length of $\vv A$ as a weak Cayley configuration.
\end{remark}
\smallskip

The definition of \codegreeG decomposition can also be reformulated in terms of projections, following this observation:
\begin{observation}\label{obs:codegreeprojection}
 For any vector configuration $\vv V$ that
decomposes into $\vv V_0\uplus \vv V_1\uplus \dots \uplus \vv V_m$, there is a projection $\pi:\RR^{\sum_i \rank \vv V_i}\to\RR^{\rank \vv V}$ that maps $\vv V_0\oplus \vv V_1\oplus \dots \oplus \vv V_m$ onto $\vv V$. \qed
\end{observation}

Next, we prove that Definition~\ref{def:dualcodegreedecompositions} is coherent with Definition~\ref{def:decomposable}. That is, that every \codegreeG decomposition is in fact dual to an codegree decomposition (and that affine \codegreeG decompositions are dual to affine codegree decompositions). For this, we need the following simple lemma.

\begin{lemma}\label{lem:totallycyclicsubconf}
For any vector configuration  $\vv V$ with $\codegG(\vv V)\geq 1$, there is a totally cyclic subconfiguration $\vv W\subseteq \vv V$ with $\codegG(\vv W)=\codegG(\vv V)$.
\end{lemma}

\begin{proof}
 The proof is by induction on the rank $r$ of $\vv V$, and trivial if $r=0$. 
 If $\vv V$ is not totally cyclic, there must be a hyperplane ${\vvh H}$ with ${\vvh H}^-\cap \vv V=\emptyset$. We can assume that ${\vvh H}$ is spanned by vectors in $\vv V$. Let $\vv W=\vv V\cap {\vvh H}$, and observe that $\codegG(\vv V)\geq \codegG(\vv W)$. Moreover, $\codegG(\vv V/\vv W)=0$ because ${\vvh H}^-\cap \vv V=\emptyset$. Finally, since $\codegG(\vv V)\leq \codegG(\vv W)+\codegG(\vv V/\vv W)$, we see that $\codegG(\vv V)=\codegG(\vv W)$, and the result follows by induction. 
\end{proof}

\begin{proposition}\label{prop:decompequivalence}
 A point configuration $\vv A$ admits a codegree decomposition of length $m$ if and only if $\vv V:=\Gale{\vv A}$ admits a \codegreeG decomposition of length $m$. The decomposition of $\vv A$ is affine if and only if the corresponding decomposition of $\vv V$ is.
\end{proposition}
\begin{proof}
 Thanks to contraction-deletion duality, we can assume that $\vv A_0=\emptyset$ and $\vv V_0=\emptyset$. We prove first the ``only if'' parts.

To prove that each affine codegree decomposition of $\vv A$ induces an affine \codegreeG decomposition of $\vv V$, observe that if $\vv A$ projects onto $\vv B_1\join\dots\join\vv B_m$, then there is a projection $\Gale\pi$ from $\pGale{\vv B_1\join\dots\join\vv B_m}=\Gale{\vv B_1}\oplus\dots\oplus\Gale{\vv B_m}$ onto $\vv V=\Gale{\vv A}$. 
Here we used Lemma~\ref{lem:directsumdual} to commute between the Gale dual of a join and the direct sum of the Gale duals, and the existence of the projection is certified by Lemma~\ref{lem:projectionduality}. 
Let $\vv V_i=\Gale\pi(\Gale{\vv B_i})$. Since $\Gale{\vv B_i}$ is the Gale dual of a (homogenized) affine point configuration,  
$$\sum_{\vv v_j\in\vv V_i}\vv v_j=\sum_{\vv w_j\in\Gale{\vv B_i}}\Gale\pi(\vv w_j)=\Gale\pi\bigg(\sum_{\vv w_j\in\Gale{\vv B_i}}\vv w_j\bigg)=\veczero \qquad\text{for all } i,$$
which proves condition \eqref{eq:defaffinecodegreecodecompositions}.
Next, Proposition~\ref{prop:degstarastar} and our assumption $\codegc(\vv A)=\sum_{i=1}^m\codegc(\vv B_i)$ imply $\codegG(\vv V)=\sum_{i=1}^m \codegG(\Gale{\vv B_i})$.
Moreover $\codegG(\vv V_i)\geq \codegG(\Gale{\vv B_i})$ because $\Gale{\vv B_i}$ projects onto $\vv V_i$. Hence	, 
$$\codegG(\vv V)=\sum_{i=1}^m \codegG(\Gale{\vv B_i})\leq \sum_{i=1}^m \codegG(\vv V_i)\leq \codegG(\vv V),$$
so equality holds throughout, and we have found an affine \codegreeG decomposition of~$\vv V$ of length~$m$. 

In the combinatorial setting, defining $\vv V_i$ as the subset of $\vv V$ that corresponds to $\vv A_i$ for each $1\leq i\leq m$ directly yields the desired combinatorial \codegreeG decomposition of~$\vv V$. Indeed, $\vv V_i=\vv V\setminus(\vv V\setminus \vv V_i)=\pGale{\vv A/(\vv A \setminus \vv A_i)}$, and the justification that these~$\vv V_i$ fulfill~\eqref{eq:defcodegreecodecompositions} is analogous to the affine case.\\

For the reciprocal statements, observe that if $\vv V$ has a \codegreeG decomposition into $\vv V_1\uplus\dots\uplus\vv V_m$, then $\vv V_1\oplus\dots\oplus\vv V_m$ projects onto $\vv V$ by Observation~\ref{obs:codegreeprojection}. This is enough for the affine case by Lemma~\ref{lem:projectionduality}, because each $\vv V_i$ is dual to a homogenized point configuration and this directly implies~\eqref{eq:defaffinecodegreecodecompositions} by Gale duality.
For the combinatorial case, Lemma~\ref{lem:totallycyclicsubconf} finds us some totally cyclic subconfiguration $\vv W_i\subset \vv V_i$ with $\codegG(\vv W_i)=\codegG(\vv V_i)$. These~$\vv W_i$ make $\vv W_1\join\dots\join\vv W_m$ project onto $\bigcup_{i=1}^m\vv W_i$, again by Observation~\ref{obs:codegreeprojection}. 
Hence, if $\vv A_i$ corresponds to~$\vv W_i$ then $\vv A/(\vv A\setminus \bigcup_{i=1}^m \vv A_i)$ projects onto $\bigjoin\nolimits_{i=1}^m\Gale{\vv W_i}$. Moreover, by Gale duality the condition that $\vv A \setminus \vv A_i$ is a face also holds, because $\vv W_i$ is totally cyclic (see Lemma~\ref{lem:gale}).
\end{proof}

A first consequence of this result is that our combinatorial characterization in Definition~\ref{def:decomposable} is well defined.

\begin{corollary}\label{cor:definitionscombdecomposable}
A point configuration $\vv A$ admits a codegree decomposition of length~$m$ if and only if it is combinatorially equivalent (as an oriented matroid) to a configuration that admits an affine codegree decomposition of length $m$.
\end{corollary}

\begin{proof}
By Proposition~\ref{prop:decompequivalence} it is enough to prove that if $\vv A$ admits a codegree decomposition of length~$m$, then its Gale dual $\vv V:=\Gale{\vv A}$ is combinatorially equivalent to a vector configuration that admits an affine \codegreeG decomposition of length $m$.

Let $\vv A=\vv A_0\uplus \vv A_1\uplus\dots\uplus \vv A_m$ be a codegree decomposition of $\vv A$ with factors $\vv A_1,\dots,\vv A_m$. And let~$\vv V_i$ with $0\leq i\leq m$ be the disjoint subsets of its Gale dual $\vv V$ corresponding to the respective $\vv A_i$. 
As in the proof of Proposition~\ref{prop:decompequivalence}, we can assume that $\vv A_0=\emptyset$ and $\vv V_0=\emptyset$ by contraction-deletion duality.

 Since $\vv A \setminus \vv A_i$ is a face, $\vv V_i$ is the support of a positive vector. Let $\vv\lambda^{(i)}\in\RR^{|\vv V_i|}_{>0}$ be such that $\sum_{\vv v_j\in \vv V_i}\lambda^{(i)}_{j}\vv v_j=\veczero$, and let $\vv \lambda\in\RR^{|{\vv V}|}_{>0}$ have entries $\lambda_{j}= \lambda^{(i)}_{j}$ if ${\vv v_j}\in \vv V_i$.
For each $\vv v_j\in \vv V$, set $\vv v_j'= \lambda_{j} \vv v_j$. These vectors form the vector configuration $\vv V':=\{\vv v_j'\}_{\vv v_j\in {\vv V}}$, which is combinatorially equivalent to ${\vv V}$ by construction. We claim that $\biguplus_i\vv V_i'$, with $\vv V_i':=\{\vv v_j'\}_{\vv v_j\in {\vv V_i}}$, is a \codegreeG decomposition of~$\vv V'$ of length~$m$. Notice that each $\codegG(\vv V_i')\ge1$ because $\vv V_i$~is the support of a positive vector. And moreover, $\codegG(\vv V_i')=\codegG(\vv V_i)$ by construction.

To check that this decomposition is affine, we have to check condition~\eqref{eq:defaffinecodegreecodecompositions}. Indeed, by construction
\begin{align*}\sum_{\vv v_j'\in\vv V_i'}\vv v_j'=\sum_{\vv v_j\in \vv V_i}\lambda_{j}\vv v_j=\sum_{\vv v_j\in \vv V_i}\lambda^{(i)}_{j}\vv v_j=\veczero.\end{align*}

The reciprocal is direct.
\end{proof}

The following proposition is the dual version of Proposition~\ref{prop:strongCayleycodeg}, which is simpler to state and prove in this setup.
\begin{proposition}
\label{prop:sumcodegrees}
Let $\vv V=\vv V_0\uplus \vv V_1\uplus \dots \uplus \vv V_m$ be a partition of a vector configuration $\vv V$. Then $\codegG(\vv V)\geq \sum_{i=0}^m \codegG(\vv V_i)$.
\end{proposition}
\begin{proof}
By definition,
 \begin{align*}\codegG(\vv V)&=\min_{\vvh H} |\ol {\vvh H}^-\cap \vv V|= \min_{\vvh H}\sum_{i=0}^m |\ol {\vvh H}^-\cap \vv V_i|\\
&\geq\sum_{i=0}^m \left(\min_{\vvh H} |\ol {\vvh H}^-\cap \vv V_i|\right) =\sum_{i=0}^m \codegG(\vv V_i). \qedhere\end{align*}
\end{proof}

\section{\texorpdfstring{Finding \codegreeG decompositions}{Finding codegree* decompositions}}
\label{sec:evidences}

It is not hard to see that Proposition~\ref{prop:dd=0} and Theorem~\ref{thm:dd=1} directly imply Conjecture~\ref{conj:strongd+1-2dd} for point configurations whose degree is at most $1$ (cf.~Propositions~\ref{prop:conjdd=0} and \ref{prop:conjdd=0}). 
In this section we provide more evidence for the conjecture by proving the Theorems~\ref{thm:2DD} and~\ref{thm:DD=2}  announced in the introduction.
Their statements and proofs involve the concept of covector discrepancy, which is presented below. For convenience, some results from previous chapters are also reformulated using this new notation. 

\subsection{The covector discrepancy}

It turns out that a useful parameter is $\DD(\vv V)$\index{$\DD(\vv V)$}, the maximal discrepancy of a covector of $\vv V$. When $\vv V$ is irreducible, $\DD(\vv V)$ is just a linear combination of its number of elements and its codegree.

\begin{definition}
The \defn{covector discrepancy}\index{covector discrepancy} of a vector configuration $\vv V$ is
\[
   \DD(\vv V)=\max_{C\in\cov(\vv V)}\left||C^+|-|C^-|\right|,
\]
 the maximal discrepancy of a covector of $\vv V$.
\end{definition}

Recall that the neighborliness of a point configuration $\vv A$ can be derived from the maximal discrepancy of a circuit of $\vv A$ (see Definitions~\ref{def:discrepancy} and~\ref{def:balanced} and Proposition~\ref{prop:neighbaldduality}). The following result shows that the almost neighborliness of~$\vv A$ is similarly associated to the maximal discrepancy of a vector of the oriented matroid of~$\vv A$.
  
\begin{lemma}\label{lem:DDrelatestodd}
Any irreducible vector configuration $\vv V$ of rank $r$ with $n=r+d+1$ elements and dual degree $\dd$ satisfies
\begin{equation}\label{eq:DDrelatestodd}
 \DD(\vv V) = r-(d+1-2\dd).
\end{equation}
In particular, plugging in definitions yields
\begin{equation}\label{eq:r+dd=DD+kk}
\rank(\vv V)+\degG(\vv V)=\DD(\vv V)+\codegG(\vv V),
\end{equation} 
and
\begin{equation}
\label{eq:achievecodegstar}
   |\vv V| \ = \ \DD(\vv V) + 2\codegG(\vv V).
\end{equation}
As a consequence, any hyperplane~${\vvh H}$ that achieves the dual codegree, so that $|\ol{\vv H}^-\cap\vv V| = \codegG(\vv V)$, satisfies $|\vv H^+\cap\vv V| =\DD(\vv V)+\codegG(\vv V)$.
\end{lemma}

\begin{proof}
It suffices to prove \eqref{eq:DDrelatestodd} and then use the identities $n=r+d+1$ and $\codegG(\vv V)=d+1-\dd$. By definition,
\begin{align*}
  \DD(\vv V)
  &=
  \max_{C\in\cov(\vv V)}\left||C^+|-|C^-|\right|= \max_{\vvh H}\left(\left|\vvh H^+\cap \vv V\right|-\left|{\vvh H}^-\cap \vv V\right|\right).
\end{align*}
Observe how the hyperplane ${\vvh H}$ that attains the maximum must be in general position with respect to $\vv V$, since otherwise ${\vvh H}$ could be perturbed to increase the discrepancy (using, for example, the strategy followed in the proof of Proposition~\ref{prop:subspacequotient}). Moreover, if $\vv V$ is irreducible and ${\vvh H}$ is in general position, then $\left|{\vvh H}^-\cap \vv V\right|=n-\left|{\vvh H}^+\cap \vv V\right|$. 
Therefore, if $\vv V$ is irreducible then $\max_{\vvh H}\left(\left|\vvh H^+\cap \vv V\right|-\left|{\vvh H}^-\cap \vv V\right|\right)=\max_{\vvh H}\left(2\left|\vvh H^+\cap \vv V\right|-n\right)$. Thus,
\begin{eqnarray}
  \DD(\vv V)
  &=&
  \max_{\vvh H}\left(2\left|\vvh H^+\cap \vv V\right|-n\right)\notag
  \\ &=&
  2\max_{\vvh H}\big|\vvh H^+\cap \vv V\big|-n \notag
  \\ &=& 
  2(r+\dd)-n=r-(d+1-2\dd)=n-2(d+1-\dd), \label{eq:DDrelatestodd2}
\end{eqnarray}
where we use that, by definition, $\dd=\degG(\vv V)=\max_{\vvh H}\big|\vvh H^+\cap \vv V\big|-r$. 
For the last claim in the statement, observe that for any hyperplane $\vvh H$ achieving the dual codegree,  
 \(|\vv H^+\cap\vv V| = |\vv V|-|\ol{\vv H}^-\cap\vv V|= |\vv V|-\codegG(\vv V)\stackrel{\eqref{eq:achievecodegstar}}{=}\DD(\vv V)+\codegG(\vv V).\)
\end{proof}

In terms of covector discrepancy, Conjecture~\ref{conj:strongd+1-2dd} reads as:

\medskip
\noindent\textbf{Conjecture~\ref{conj:strongd+1-2dd} (reformulated)}
\emph{Any irreducible vector configuration $\vv V$ of rank $r$ admits a \codegreeG decomposition of length at least $r-\DD(\vv V)$.}

\medskip

Since Lemma~\ref{lem:DDrelatestodd} only relates $\DD(\vv V)$ to $\degG(\vv V)$ when $\vv V$ is irreducible, from now on we focus on irreducible configurations; that is, vector configurations whose Gale dual is not a pyramid. 
The following observation shows that this assumption is safe, in the sense that Conjecture~\ref{conj:strongd+1-2dd} is true if and only if it is for irreducible configurations. This is because a pyramid over~$\vv B$ admits a codegree decomposition of length $d+1-2\dd$ if and only if the base~$\vv B$ admits a codegree decomposition of length $(d-1)+1-2\dd$:

\begin{observation}\label{obs:pyramids}
 Let $\vv V$ be a vector configuration of rank $r$ with $r+d+1$ elements, dual degree~$\delta$, and exactly $k$~copies $\vv v_1,\dots,\vv v_k$ of~$\veczero$ (\ie $\Gale {\vv V}$ is a $k$-fold pyramid). 
 If $\vv V'=\vv V\setminus \{\vv v_1,\dots,\vv v_k\}$ admits a \codegreeG decomposition of length $m$, then $\vv V$ admits a \codegreeG decomposition of length $m+k$ just by including the sets $\vv V_{m+i}=\{\vv v_i\}$ as factors for $1\leq i \leq k$.
Now $\vv V'$~is a configuration of rank $r'=r$ with $r'+d'+1=r+d+1-k$ elements and dual degree~$\dd'=\dd$, so that $(d+1-2\dd)=(d'+1'-2\dd')+k$. Hence, if $\vv V'$~admits a \codegreeG decomposition of length at least $d'+1'-2\dd'$ as promised in Conjecture~\ref{conj:strongd+1-2dd}, then $\vv V$ also admits a \codegreeG decomposition of length at least $d+1-2\dd$.
\end{observation}

Another important (yet straightforward) remark concerning pyramids that we use later is that adding or removing the origin from a vector configuration does not change the covector discrepancy:
  
\begin{lemma}\label{lem:DDirreducible}
  $\DD(\vv V)=\DD(\vv V\cup\veczero)$ for any vector configuration $\vv V$.\qed
\end{lemma}

Reformulating Corollary~\ref{cor:easybound} in this language makes its proof trivial:

\begin{corollary}\label{cor:easyboundDD}
$\DD(\vv V)\geq 0$ for any vector configuration $\vv V$. In particular, $r\geq d+1-2\dd$ for any irreducible vector configuration of rank $r$ with $r+d+1$ elements and degree~$\dd$.
\end{corollary}

\begin{proof}
The first part is trivial, since $\DD(\vv V)$ is an absolute value. The second part follows from Lemma~\ref{lem:DDrelatestodd}.
\end{proof}

We conclude that \emph{configurations of vanishing discrepancy are extremal}.
We studied those that were irreducible in Proposition~\ref{prop:Lawrence}, which we reformulate (and expand to non-irreducible) as:

\begin{corollary}\label{cor:Lawrence}
 A vector configuration $\vv V$ fulfills $\DD(\vv V)=0$ if and only if it is centrally symmetric (up to rescaling).
\end{corollary}
\begin{proof}
 For irreducible configurations this is proven in Proposition~\ref{prop:Lawrence}, and adding the origin does not change neither central symmetry nor covector discrepancy.
\end{proof}

This result can be interpreted in terms of \codegreeG decompositions and weak Cayley decompositions, as the following corollary shows. Since $d+1-\dd\geq d+1-2\dd$ because $\dd\ge0$, this already proves Conjecture~\ref{conj:strongd+1-2dd} for vector configurations with $\DD(\vv V)=0$.

\begin{corollary}\label{cor:Lawrencedecomposition}
 Let $\vv V$ be a vector configuration of rank~$r$ and dual degree~$\dd$ with $n=r+d+1$ elements. If $\DD(\vv V)=0$, then $\vv V$~admits a \codegreeG decomposition of length 
\[
   d+1-\dd=\codegG(\vv V)\geq \rank(\vv V)
\] 
into pairs of antipodal vectors and copies of the origin. In particular, it is a weak Cayley configuration of length $d+1-\dd$.
\end{corollary}
\begin{proof}
By Corollary~\ref{cor:Lawrence}, $\vv V$ is centrally symmetric. If $\vv V$ is irreducible, then the fact that $\DD(\vv V)=0$ implies that $n=2(d+1-\dd)$ by~\eqref{eq:achievecodegstar}. Hence, the decomposition of~$\vv V$ into pairs of antipodal vectors is a \codegreeG decomposition of length $\frac{n}{2}=\frac{2(d+1-\dd)}{2}=d+1-\dd$. Moreover, when $\DD(\vv V)=0$ then $\codegG(\vv V)=\rank(\vv V)+\degG(\vv V)\geq \rank(\vv V)$ by~\eqref{eq:r+dd=DD+kk}.

On the other hand, if $\vv V$ is not irreducible and contains exactly $k$ copies of the origin, then each of these copies contributes to $1$ to the dual codegree. Hence removing all $k$ copies of the origin creates $k$ factors of rank~$0$ and leaves an irreducible configuration of codegree $\codegG(\vv V)-k$. Then the claim follows from the irreducible case (cf.  Observation~\ref{obs:pyramids}).
\end{proof}

One final observation that concerns addition and deletion of antipodal pairs of vectors.

\begin{lemma}\label{lem:removecs}
 For any vector configuration $\vv V\subset\RR^r$, any vector $\vv v\in\RR^r$ and any real number $\gl>0$, $\DD(\vv V)=\DD(\vv V\cup\{\vv v,-\gl\vv v\})$.
\end{lemma} 
\begin{proof}
Observe that for every hyperplane $\vvh H$, either $\{\vv v,-\gl\vv v\}\subset \vvh H$ or both $|\vvh H^\pm\cap\{\vv v,-\gl\vv v\}|=1$. Hence, for every $\vvh H$
 \[\left|\vvh H^+\cap \vv V\right|-\left|{\vvh H}^-\cap \vv V\right|=\left|\vvh H^+\cap \left( \vv V\cup\{\vv v,-\gl\vv v\}\right)\right|-\left|{\vvh H}^-\cap \left(\vv V\cup\{\vv v,-\gl\vv v\}\right)\right|.\qedhere\]
\end{proof}

\subsection{Extremal subconfigurations} \label{sec:extremalconfigurations}
It is also convenient to rewrite Proposition~\ref{prop:subspacequotient}, which relates the degree of a configuration to the degree of a subconfiguration and the degree of its quotient, in terms of the covector discrepancy.

\begin{corollary}\label{cor:subspacequotient}
Let ${\vv V}$ be a vector configuration, and let ${\vv W}\subset {\vv V}$ be a subconfiguration such that $\lin ({\vv W})\cap {\vv V}={\vv W}$. Then, 
\begin{equation}\label{eq:DDineq}
 \DD(\vv V)\geq \DD(\vv W)+\DD(\vv V/\vv W).
\end{equation}
\end{corollary}
\begin{proof}
If $\vv V$ is irreducible, this is a direct consequence of Proposition~\ref{prop:subspacequotient}, using \eqref{eq:achievecodegstar} and that $|\vv V|=|\vv W|+|\vv V/\vv W|$. On the other hand, if $\veczero\in \vv V$ then $\veczero\in \vv W$ since $\lin ({\vv W})\cap {\vv V}={\vv W}$. The last ingredient to complete the proof is Lemma~\ref{lem:DDirreducible}.
\end{proof}

The same argument easily proves the following equivalences.
\begin{lemma}\label{lem:subspacequotientequality}
Let ${\vv V}$ be a vector configuration, and let ${\vv W}\subset {\vv V}$ be a subconfiguration such that $\lin ({\vv W})\cap {\vv V}={\vv W}$. Then, the following three equalities are equivalent:
\begin{align}
\DD(\vv V)&=\DD(\vv W)+\DD(\vv V/\vv W);\label{eq:idDD}\\ \degG(\vv V)&=\degG(\vv W)+\degG(\vv V/\vv W);\text{ and }\label{eq:iddd}\\\codegG(\vv V)&=\codegG(\vv W)+\codegG(\vv V/\vv W).\label{eq:idkk} 
\end{align}\qed
\end{lemma}

We are particularly interested in the subconfigurations where equality in \eqref{eq:DDineq}--\eqref{eq:idkk} holds, and will now show how to find such subconfigurations using one extremal hyperplane. This procedure, which will be refined in subsequent results, is also explained in Example~\ref{ex:extremalconfigurations} and Figure~\ref{fig:extremalconfigurations}.

\begin{lemma}\label{lem:extremalconfigurations}
For any vector configuration $\vv V$ of rank~$r$, and for any~$s$ with $1\le s\le r-1$, there exists a  subconfiguration~$\vv W$ of rank~$s$ such that 
\begin{equation}\label{eq:extremalconfigurations}
  \DD(\vv V)=\DD(\vv W)+\DD(\vv V/\vv W).
\end{equation}
\end{lemma}

\begin{proof}
Let ${\vvh H}$ be any hyperplane that achieves the dual codegree of~$\vv V$, so that $|\ol{\vvh H}^-\cap \vv V|=\codegG(\vv V)$. Moreover, let $\vvh C_{\vvh H}$ be the polyhedral cone\index{$\vvh C_{\vvh H}$} 
\[
\vvh C_{\vvh H}=\cone\left({\big({\vvh H}^+\cap \vv V\big)\cup \big({\vvh H}^+\cap -\vv V\big)}\right).
\] 
We will show that the subconfiguration $\vv W= \vv V\cap\vv {\vvh H}'$ induced by any supporting hyperplane~${\vvh H}'$ of~$\vvh C_{\vvh H}$ satisfies~\eqref{eq:extremalconfigurations}; choosing one supporting a face of~$\vvh C_{\vvh H}$ of the appropriate dimension then concludes the proof.

For this, observe that $\vvh C_{\vvh H}\subset\ol{\vvh H}'^+$ and that ${\vvh H}'$ fulfills ${\vvh H}'^+\cap \vv V\subseteq {\vvh H}^+\cap \vv V$ and ${\vvh H}'^-\cap \vv V\subseteq {\vvh H}^-\cap \vv V$. Moreover, by Lemma~\ref{lem:subspacequotientequality}, \eqref{eq:extremalconfigurations}~is equivalent to 
\begin{equation}\label{eq:codeg2}
  \codegG(\vv V)=\codegG(\vv W)+\codegG(\vv V/\vv W).
\end{equation}

To prove~\eqref{eq:codeg2}, we use that, by definition, $\codegG(\vv W)\leq|\ol{\vvh H}^-\cap \vv W|$ because~${\vvh H}$ cuts $\lin(\vv W)$ in a hyperplane of~$\lin(\vv W)$. Moreover, ${\vvh H}'$ can be considered as a hyperplane of~$\vv V/\vv W$ because $\vv W\subset{\vvh H}'$. This proves that $\codegG(\vv V/\vv W)\leq |\ol{\vvh H}'^-\cap (\vv V/\vv W)|=|\ol{\vvh H}^-\cap (\vv V\setminus \vv W)|$. Summing up,
 \begin{align*}
\codegG(\vv W)+\codegG(\vv V/\vv W)&\leq \big|\ol{\vvh H}^-\cap \vv W\big|+ \big|\ol{\vvh H}^-\cap (\vv V\setminus \vv W)\big|\\&= |\ol{\vvh H}^-\cap \vv V|=\codegG(\vv V).
\end{align*}
The fact that $\codegG(\vv V)\leq \codegG(\vv W)+\codegG(\vv V/\vv W)$ by Proposition~\ref{prop:subspacequotient} concludes the proof.
\end{proof}

\iftoggle{bwprint}{%
\begin{figure}[htpb]
\centering
 \subbottom[$\vv V$ and $\vvh H$]{\label{sfig:epc_1}\includegraphics[width=.4\linewidth]{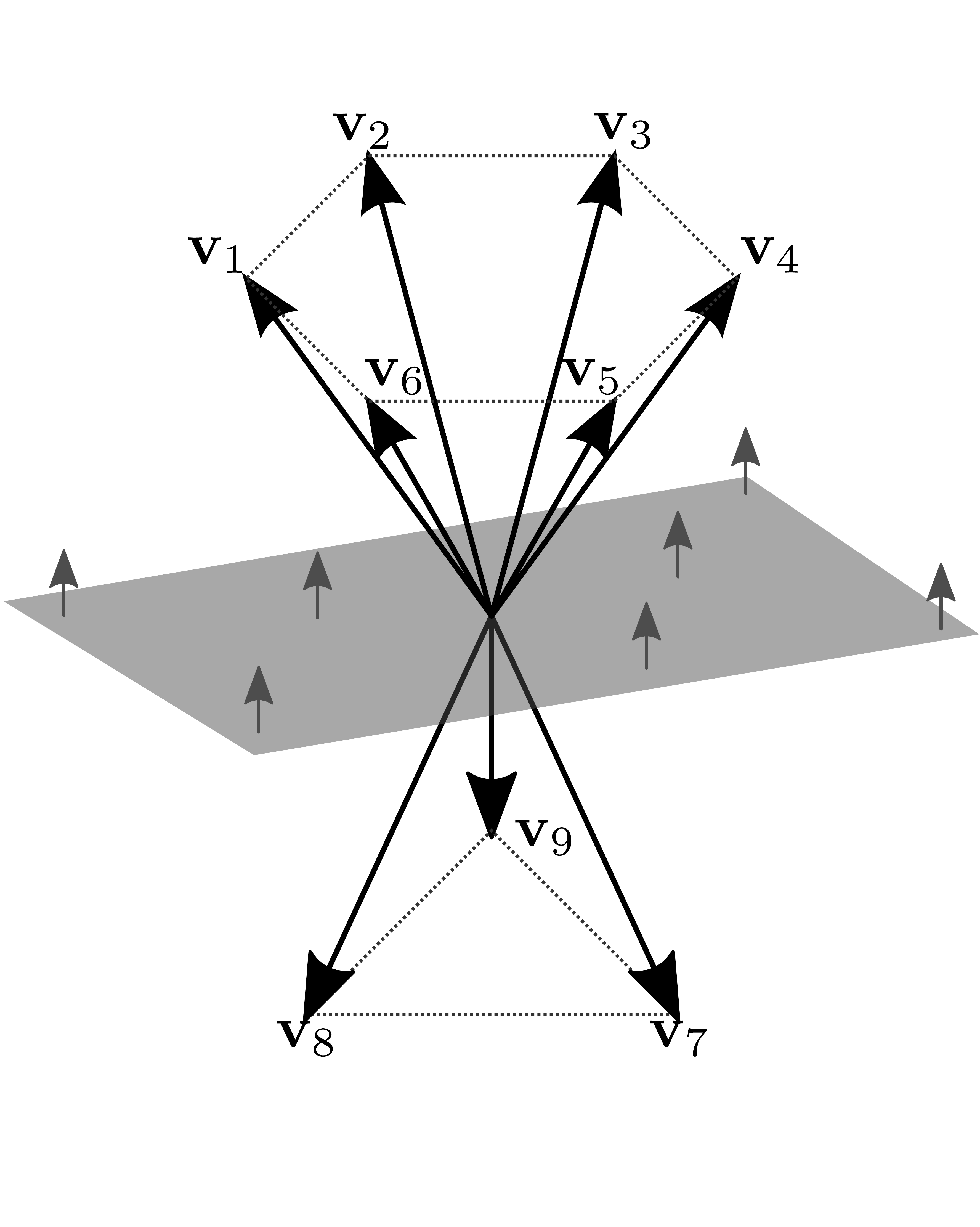}}\quad
 \subbottom[$\vvh C_{\vvh H}$ and ${\vvh H}'$]{\label{sfig:epc_2}\includegraphics[width=.4\linewidth]{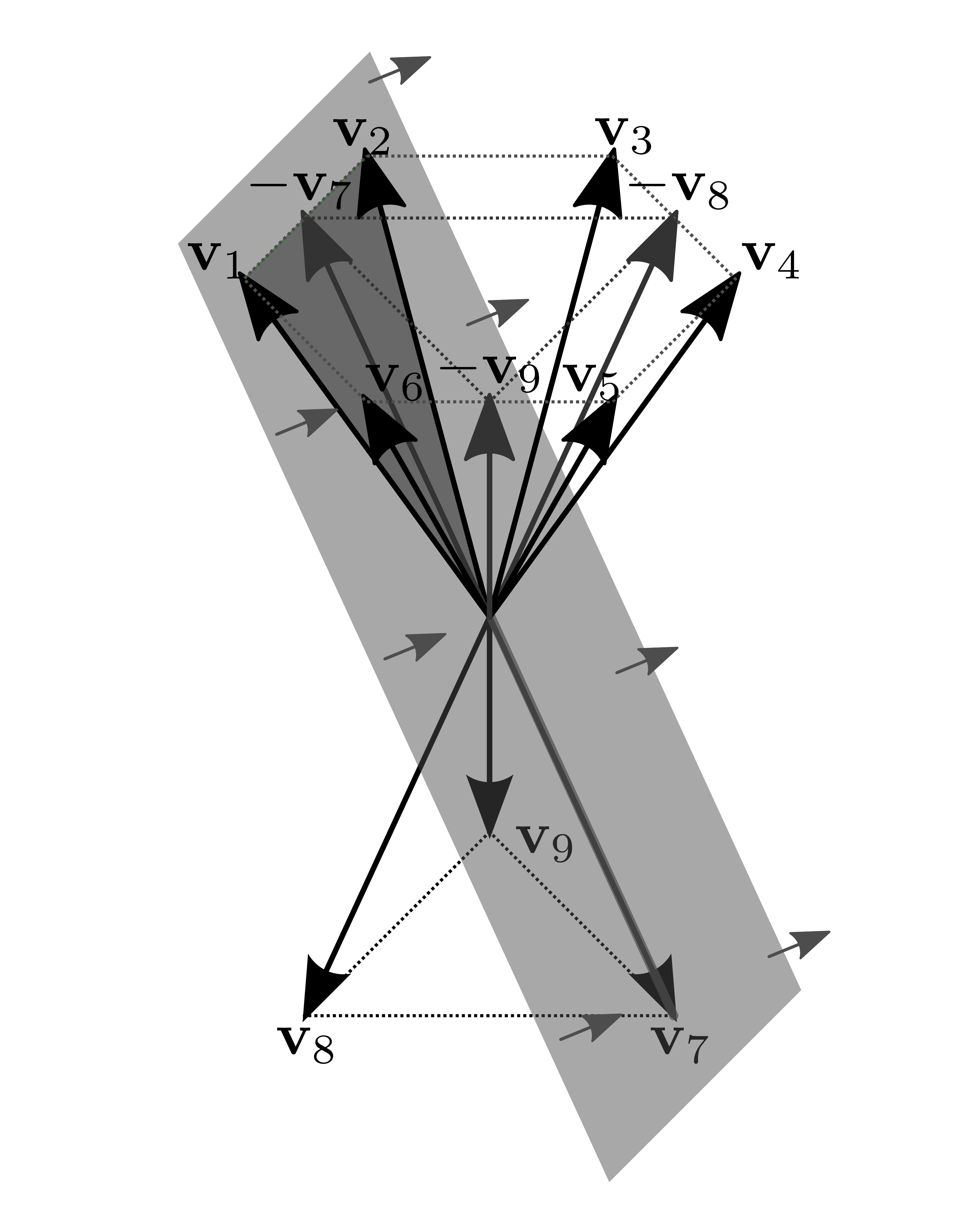}}\quad
 \subbottom[$\vv W$ and ${\vvh H}$]{\quad\label{sfig:epc_4}\includegraphics[width=.4\linewidth]{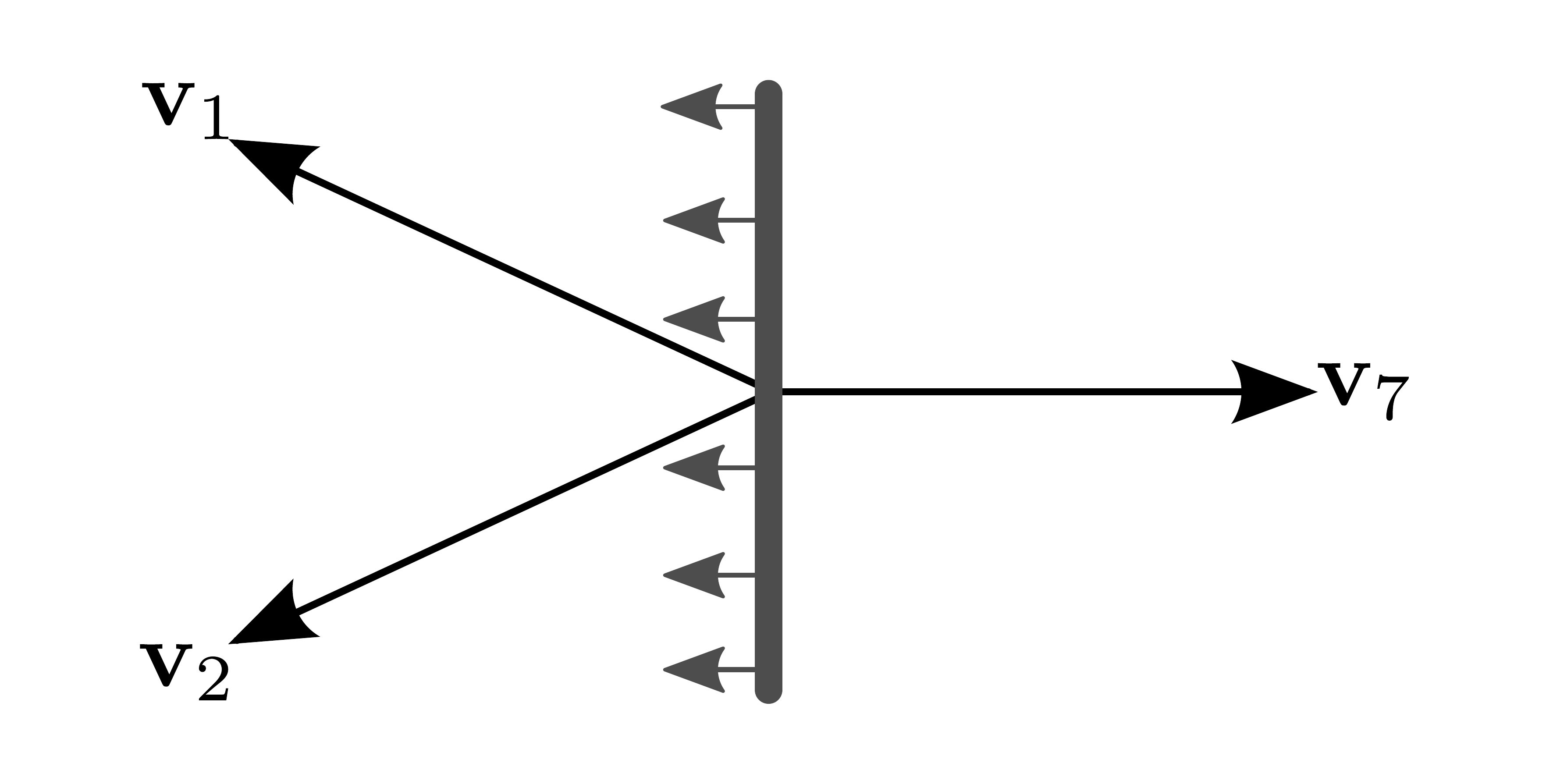}\,\,}
 \subbottom[$\vv V/\vv W$ and ${\vvh H}'$]{\quad\label{sfig:epc_3}\includegraphics[width=.4\linewidth]{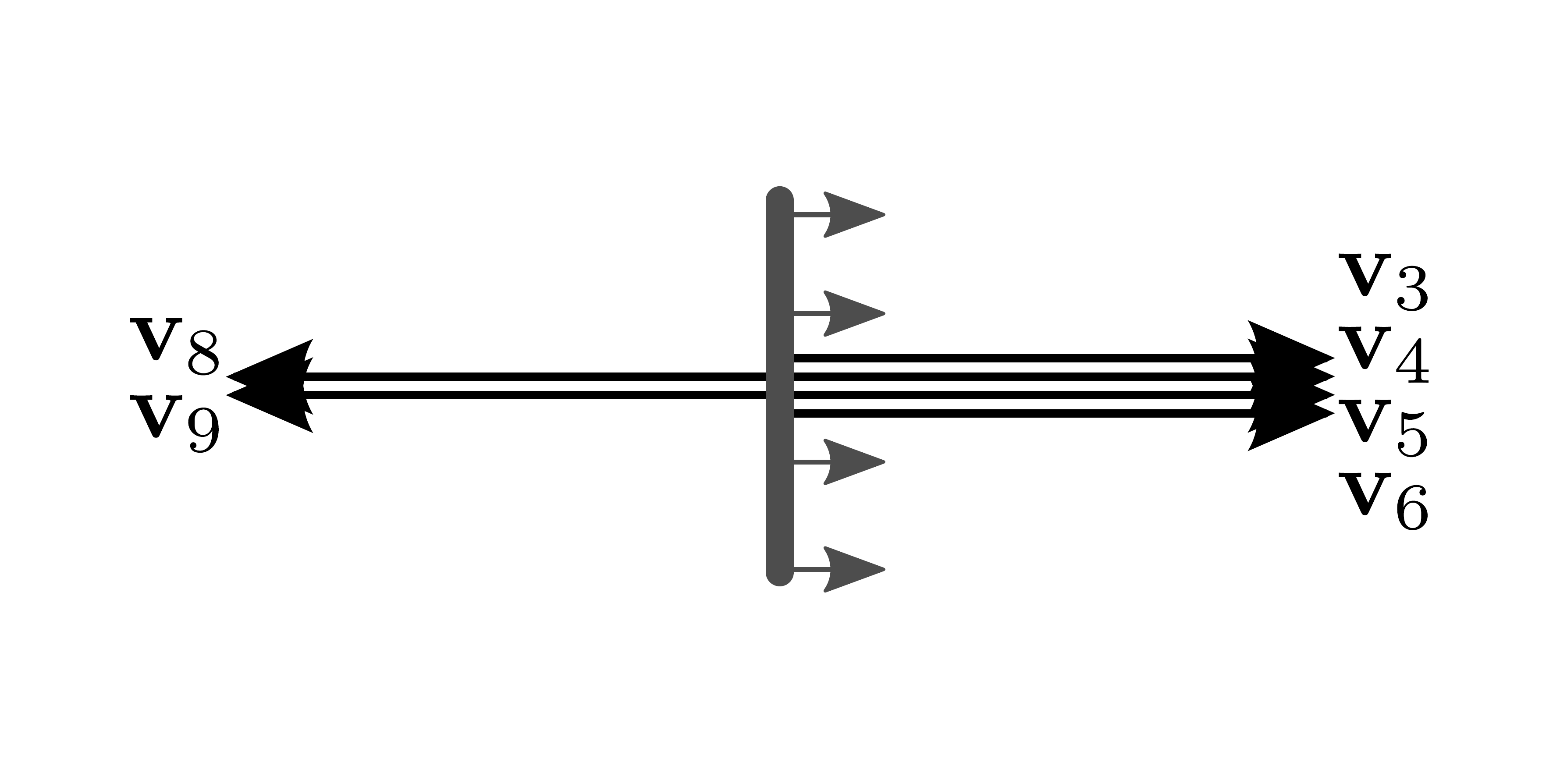}\qquad}
\caption[Illustrating Lemma~\ref{lem:extremalconfigurations}.]{Illustrating Lemma~\ref{lem:extremalconfigurations}. 
}
 \label{fig:extremalconfigurations}
\end{figure}
}{%
\begin{figure}[htpb]
\centering
 \subbottom[$\vv V$ and $\vvh H$]{\label{sfig:epc_1}\includegraphics[width=.4\linewidth]{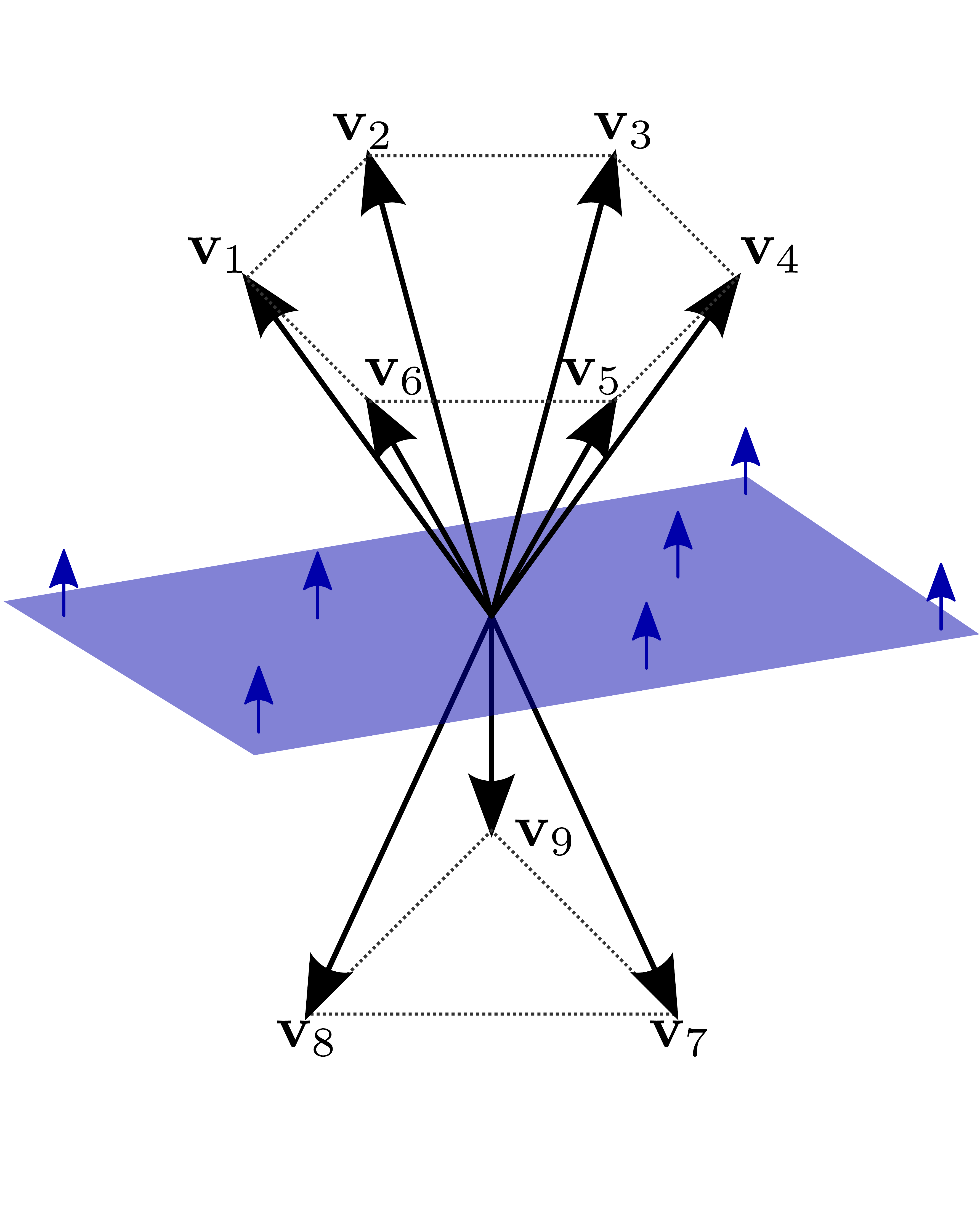}}\quad
 \subbottom[$\vvh C_{\vvh H}$ and ${\vvh H}'$]{\label{sfig:epc_2}\includegraphics[width=.4\linewidth]{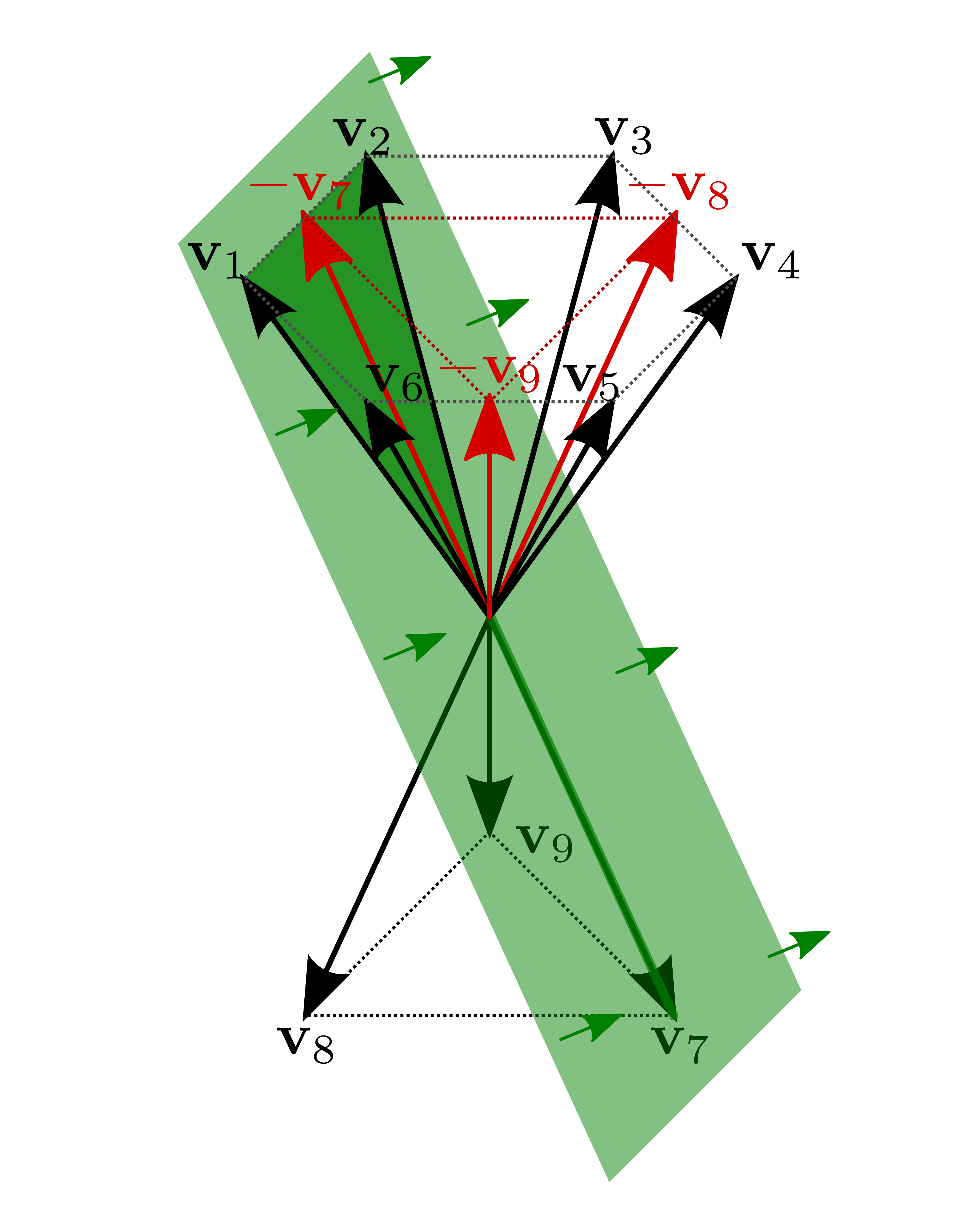}}\quad
 \subbottom[$\vv W$ and ${\vvh H}$]{\quad\label{sfig:epc_4}\includegraphics[width=.4\linewidth]{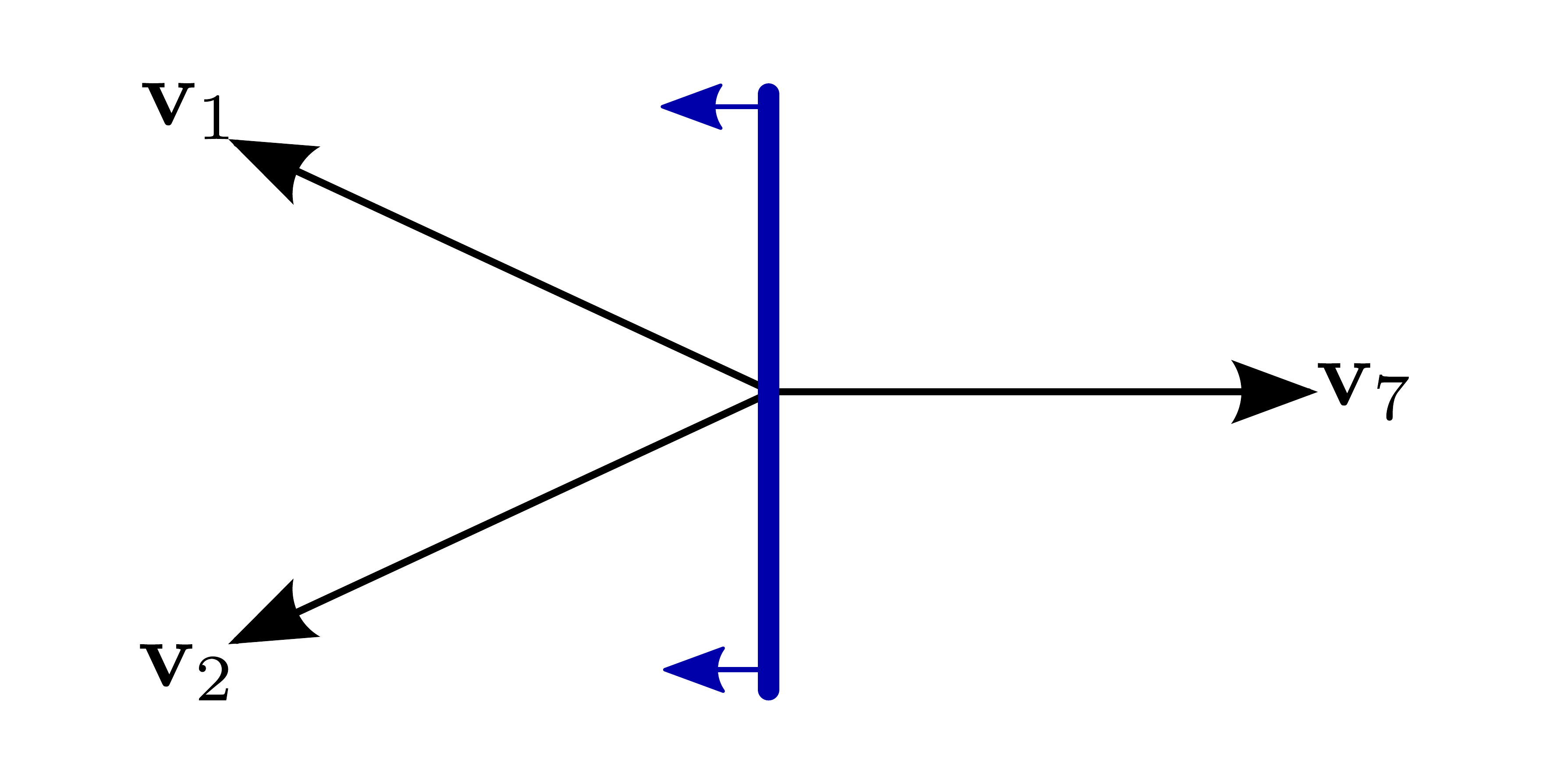}\,\,}
 \subbottom[$\vv V/\vv W$ and ${\vvh H}'$]{\quad\label{sfig:epc_3}\includegraphics[width=.4\linewidth]{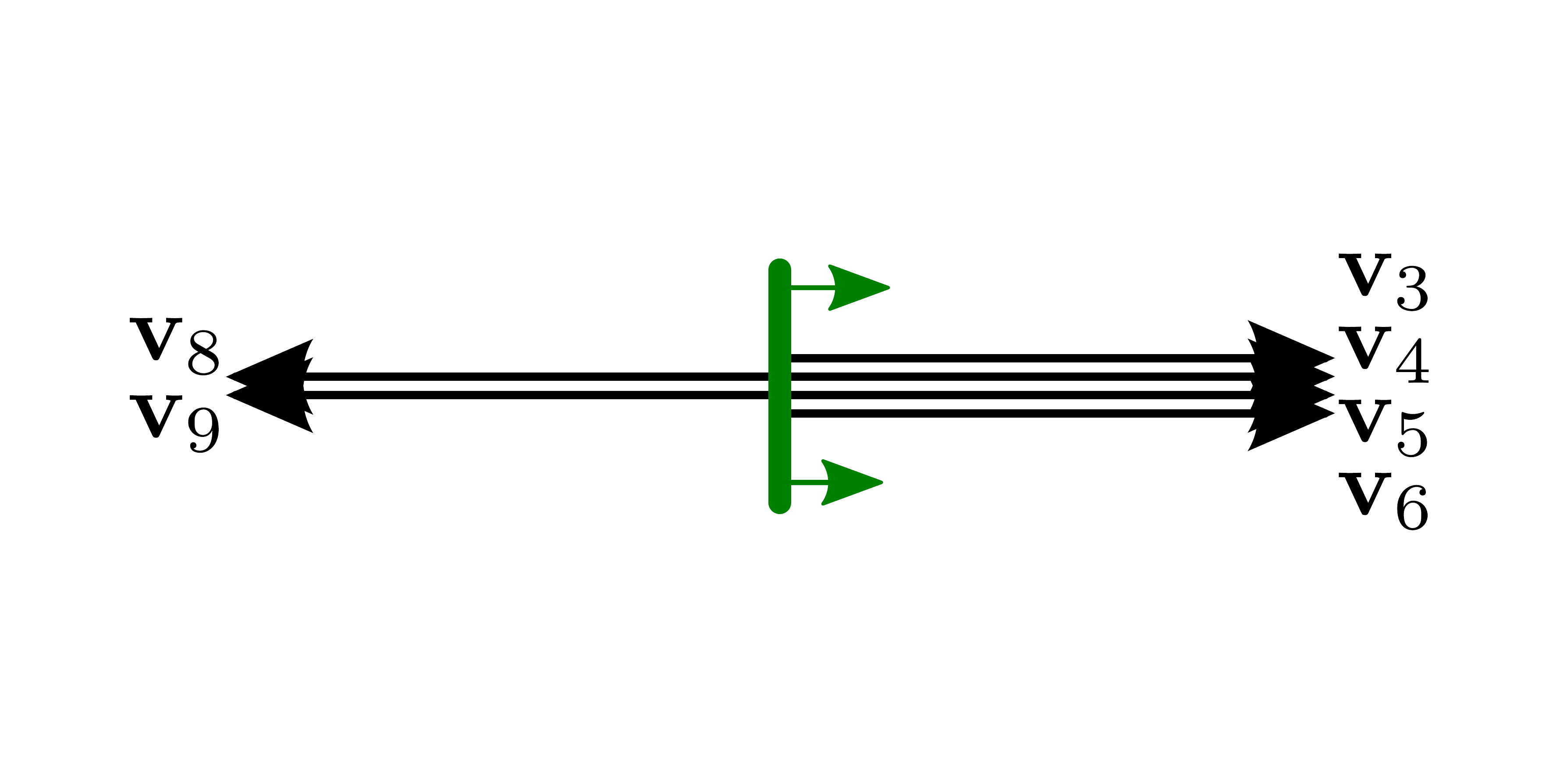}\qquad}
\caption[Illustrating Lemma~\ref{lem:extremalconfigurations}.]{Illustrating Lemma~\ref{lem:extremalconfigurations}. 
}
 \label{fig:extremalconfigurations}
\end{figure}
}

\begin{example}\label{ex:extremalconfigurations}
Consider the configuration $\vv V$ of $9$ vectors in $\RR^3$ displayed in Figure~\ref{fig:extremalconfigurations}.
Since $\degG(\vv V)=3$ and $\codegG(\vv V)=3$, the hyperplane~$\vvh H$ in~\ref{sfig:epc_1} fulfills $|\ol {\vvh H}^-\cap \vv V|=\codegG(\vv V)$. In~\ref{sfig:epc_2} one can see $\vvh C_{\vvh H}$ together with one of its supporting hyperplanes~${\vvh H}'$. ${\vvh H}'\cap \vv V$ is the subconfiguration $\vv W=\{\vv v_1,\vv v_2,\vv v_ 7\}$.
Observe in~\ref{sfig:epc_4} how just by intersecting with $\lin(\vv W)$, one can consider ${\vvh H}$ as a hyperplane of $\lin(\vv W)$. This shows that $|\ol {\vvh H}^-\cap \vv W|\ge\codegG(\vv W)$. Finally, in~\ref{sfig:epc_4} we see how to translate~$\vvh H'$ into a hyperplane of $\vv V/\vv W$. 
This hyperplane fulfills $|\ol {\vvh H'}^-\cap \vv V/\vv W|=|\ol {\vvh H}^-\cap \vv V\setminus \vv W|$, which proves $|\ol {\vvh H}^-\cap \vv V\setminus\vv W|\ge\codegG(\vv V/\vv W)$. This implies that $\codegG(\vv V)=\codegG(\vv W)+\codegG(\vv V/\vv W)$, which we can verify by observing that $\codegG(\vv W)=1$ and $\codegG(\vv V/\vv W)=2$.
\end{example}

\begin{remark}
 In the proof of Lemma~\ref{lem:extremalconfigurations}, we could have defined $\vvh C_{\vvh H}$ directly as $\vvh C_{\vvh H}=\cone\big({\vvh H}^+\cap \vv V\big)$. Indeed, it is not hard to see that each extremal ray of $\vvh C_{\vvh H}$ must be a extremal ray of $\cone\big({\vvh H}^+\cap \vv V\big)$, because otherwise we could perturb $\vvh H$ until it contains less elements in $\ol{\vvh H}^-$.
\end{remark}

In particular, every face of $\vvh C_{\vvh H}$ provides a configuration fulfilling \eqref{eq:idDD}, and specifically, the extremal rays of $\vvh C_{\vvh H}$ give such configurations of rank~$1$. Lemma~\ref{lem:nocs} below shows that if $\DD(\vv V)>0$, we can find them with positive discrepancy. Further conditions are established in Lemma~\ref{lem:strongnocs}.

\begin{lemma}\label{lem:nocs}
Let $\vv V$ be a vector configuration with $\DD(\vv V)>0$. Then there is a subconfiguration $\vv W\subset \vv V$ of rank $1$ with $\lin(\vv W)\cap \vv V=\vv W$ and $\DD(\vv W)>0$, such that $\DD(\vv V)=\DD(\vv W)+\DD(\vv V/\vv W)$.
\end{lemma}
\begin{proof}
 Let ${\vvh H}$ be a hyperplane in general position with respect to $\vv V$ fulfilling $|\ol{\vvh H}^-\cap V|=\codegG(\vv V)$. Perturb ${\vvh H}$ until it hits the first vector $\vv v$ in $\vv V$ such that $\lin (\vv v)\cap \vv V$ is not centrally symmetric. This must happen eventually because $\DD(\vv V)>0$ (see Corollary~\ref{cor:Lawrence}).
 Observe that, whenever $\vvh H$ sweeps across a pair of centrally symmetric vectors, $|\vvh H^+ \cap \vv V|$ remains constant. Therefore, when $\vvh H$ hits the first non-centrally symmetric vector, it is a supporting hyperplane of a cone such as that in the proof of Lemma~\ref{lem:extremalconfigurations}.
Now setting $\vv W = \lin (\vv v)\cap \vv V$ concludes the proof. 
\end{proof}

\begin{lemma}\label{lem:strongnocs}
Let $\vv V$ be a vector configuration and let $\vv W\subset \vv V$ be a subconfiguration such that $\lin(\vv W)\cap \vv V=\vv W$. Then, either $\DD(\vv V\setminus \vv W)=0$, or there is a different subconfiguration $\vv U\subset (\vv V\setminus \vv W)$ of rank $1$ with $\lin(\vv U)\cap \vv V=\vv U$ and $\DD(\vv U)>0$ such that $\DD(\vv V)=\DD(\vv U)+\DD(\vv V/\vv U)$.

In particular, for any linear hyperplane $\vvh H$, either all vectors of $\vv V\setminus\vvh H$ come in antipodal pairs, or we can find such an $\vv U$ in $\vv V\setminus\vvh H$.
\end{lemma}

\begin{proof}
By Lemma~\ref{lem:removecs}, we can assume that $\vv V$ contains no pair of antipodal vectors. If $\vv W=\vv V$ we are trivially done, because then $\vv V\setminus\vv W=\emptyset$. 

 Otherwise, let ${\vvh H}$ be a hyperplane such that $|\ol{\vvh H}^-\cap \vv V|=\codegG(\vv V)$ and $\vvh C_{\vvh H}=\cone\left({({\vvh H}^+\cap \vv V)\cup ({\vvh H}^+\cap (-\vv V))}\right)$. By Lemma~\ref{lem:extremalconfigurations}, our claim reduces to finding an extremal vector of $\vvh C_{\vvh H}$ that does not lie in $\lin (\vv W)$, and this can be done since $\rank(\vvh C_{\vvh H})=\rank(\vv V)>\rank(\vv W)$.
\end{proof}

While our Lemma~\ref{lem:extremalconfigurations} can already be used to find such extremal subconfigurations of rank larger than~$1$, we will eventually need the upcoming Corollary~\ref{cor:nocs}, which is stronger since it adds extra constraints on the discrepancy of the subconfigurations. Its proof is a combination Lemma~\ref{lem:nocs} with Lemma~\ref{lem:iterativeextremal} below.

\begin{lemma}\label{lem:iterativeextremal}
 Let $\vv V$ be a vector configuration, let $\vv W$ be a subconfiguration of~$\vv V$ such that 
\begin{equation}
  \label{eq:DDV}\DD(\vv V)=\DD(\vv W)+\DD(\vv V/\vv W),
\end{equation} 
and let $\vv U$ be a subconfiguration of~$\vv V/\vv W$ such that \begin{equation}\label{eq:DDV/W}
 \DD(\vv V/\vv W)=\DD\big((\vv V/\vv W)/\vv U\big)+\DD(\vv U).
\end{equation}

Then $(\vv U\cup \vv W)\subset\vv V$ fulfills 
\begin{equation}\label{eq:iterativeextremalsum}
\DD(\vv U\cup \vv W)=\DD(\vv U)+\DD(\vv W), 
\end{equation}
and hence 
\[
\DD(\vv V)=\DD(\vv U\cup \vv W)+\DD\big(\vv V/(\vv U\cup\vv W)\big).
\]
Here we identify $\vv U\subset(\vv V/\vv W)$ with the corresponding $\vv U\subset(\vv V\setminus\vv W)$.
\end{lemma}

\begin{proof}
By Corollary~\ref{cor:subspacequotient}, $\DD(\vv U\cup \vv W)\le  \DD(\vv V)-\DD({\vv V/(\vv U\cup \vv W)})$. Since $\vv V/(\vv U\cup \vv W) = (\vv V/\vv W)/\vv U$, using equation~\eqref{eq:DDV/W} and later \eqref{eq:DDV} yields
\begin{eqnarray*}
  \DD(\vv U\cup \vv W)
  &\stackrel{\ref{cor:subspacequotient}}{\leq}&
  \DD(\vv V)-\DD({\vv V/(\vv U\cup \vv W)})
  \\ &=&
  \DD(\vv V)-\DD\big((\vv V/\vv W)/\vv U\big)
  \\ 
  &\stackrel{\eqref{eq:DDV/W}}{=}&\DD(\vv V)-(\DD(\vv V/\vv W)-\DD(\vv U))\\
  &=&(\DD(\vv V)-\DD(\vv V/\vv W))+\DD(\vv U)\\
  &\stackrel{\eqref{eq:DDV}}{=}&
  \DD(\vv W)+\DD(\vv U). 
\end{eqnarray*}
Moreover, again by Corollary~\ref{cor:subspacequotient},
\begin{align*}
\DD(\vv U\cup \vv W)&\geq\DD(\vv W)+\DD((\vv U\cup \vv W)/\vv W)=\DD(\vv W)+\DD(\vv U),
\end{align*}
and we are done.
\end{proof}

\begin{corollary}\label{cor:nocs}
For any vector configuration $\vv V$ of rank $r$ and any~$s$ with $1\leq s\leq r-1$, there is a subconfiguration $\vv W\subset\vv V$ of rank $s$ that fulfills $\DD(\vv V)=\DD(\vv W) +\DD(\vv V/\vv W)$, $\DD(\vv W)\geq \min\{\DD(\vv V),s\}$ and $\vv W=\lin(\vv W)\cap \vv V$.
\end{corollary}
\begin{proof}
If $\DD(\vv V)=0$, the claim is trivially true since $\DD(\vv W)\geq 0$ and $\DD(\vv V/\vv W)\geq 0$ for any $\vv W$.
For $\DD(\vv V)>0$, first observe that the trivial solution $\vv W=\emptyset$ is excluded by the inequality $\DD(\vv W)\geq \min\{\DD(\vv V),s\}$. Our proof proceeds by induction on~$s$. 
The case $s=1$ follows from Lemma~\ref{lem:nocs}. 
Otherwise, if $s\geq 2$, use Lemma~\ref{lem:nocs} to find some $\vv W_0$ of rank~$1$ with $\DD(\vv W_0)\geq\min\{\DD(\vv V),1\}\geq 1$ such that \begin{equation}\label{eq:cornocs1}\DD(\vv V)=\DD(\vv W_0)+\DD(\vv V/\vv W_0).\end{equation}
 By induction hypothesis, there is some $\vv W_1$ of rank~$s-1$ in $\vv V/\vv W_0$ such that $\DD(\vv W_1)\geq \min\{\DD(\vv V/\vv W_0),s-1\}$ and 
\begin{equation}\label{eq:cornocs2}
\DD(\vv V/\vv W_0)=\DD(\vv W_1)+\DD\!\big((\vv V/\vv W_0)/\vv W_1\big).
\end{equation}
 The subconfiguration we are looking for is $\vv W=\vv W_0\cup\vv W_1$. Observe that $\DD(\vv W)=\DD(\vv W_1)+\DD(\vv W_0)$ and $\DD(\vv V)=\DD(\vv W)+\DD(\vv V/\vv W)$ by Lemma~\ref{lem:iterativeextremal}.
We now claim that 
\[
   \DD(\vv W)=\DD(\vv W_1)+\DD(\vv W_0)\geq \min\{\DD(\vv V),s\}.
\]
Indeed, recall that $\DD(\vv W_0)\geq 1$ and $\DD(\vv W_1)\geq \min\{\DD(\vv V/\vv W_0),s-1\}$. 
\begin{itemize}
\item If $\DD(\vv V)<s$, then $\DD(\vv V/\vv W_0)=\DD(\vv V)-\DD(\vv W_0)< s-1$. This means that $\min\{\DD(\vv V/\vv W_0),s-1\}=\DD(\vv V/\vv W_0)$ and hence 
\begin{eqnarray*}
\DD(\vv W_1)+\DD(\vv W_0)
&\geq& 
\DD(\vv V/\vv W_0)+\DD(\vv W_0)
\\ &\stackrel{\eqref{eq:cornocs2}}{=}&
\DD(\vv V)\ =\ \min\{\DD(\vv V),s\}.
\end{eqnarray*}

\item Else, $\DD(\vv V)\geq s$ and $\DD(\vv V/\vv W_0)=\DD(\vv V)-\DD(\vv W_0)\geq s-\DD(\vv W_0)$. Hence $\min\{\DD(\vv V/\vv W_0),s-1\}\geq s-\DD(\vv W_0)$ and
\begin{eqnarray*}
\DD(\vv W_1)+\DD(\vv W_0)
&\geq& 
s- \DD(\vv W_0)+\DD(\vv W_0)
\\ &=&
s\ = \ \min\{\DD(\vv V),s\}.
\end{eqnarray*}
\end{itemize}

Finally,  $\DD(\vv V)\geq\DD(\vv W)$ is a direct consequence of Corollary~\ref{cor:subspacequotient}.
\end{proof}

\subsection{\texorpdfstring{Decompositions of length $>r-2\DD(\vv V)$}{Decompositions of length >r-2 Disc(V)} }

Thanks to these results, we are able to prove Theorem~\ref{thm:2DD}, that finds non-trivial \codegreeG decompositions when $r>2\DD(\vv V)$. The core of its proof is encapsulated in the following proposition.

\begin{proposition}\label{prop:2DD}
Every vector configuration $\vv V$ has a subconfiguration~$\vv W$ of rank at most $2\DD(\vv V)$ such that $\DD(\vv W)=\DD(\vv V)$ and $\DD(\vv V\setminus \vv W)=0$. 
\end{proposition}

Since the proof is quite long (but not actually complicated), we first outline our strategy. We start by finding a subconfiguration~$\vv W_0$ such that the contraction $\vv W/\vv W_0$ is centrally symmetric. This contraction is divided into several centrally symmetric pieces $\vv U_1,\dots,\vv U_m$ of rank~$1$. 
Then we show that the preimage in~$\vv V$ of most of these pieces (that is, before contracting~$\vv W_0$) is still centrally symmetric, where ``most'' means all except for a subconfiguration of rank at most~$2\DD(\vv V)$. See also Example~\ref{ex:2DD} and Figure~\ref{fig:2DD}.

\begin{proof}
We put $\l=\DD(\vv V)$, and make the following simplifying assumptions: 
\begin{itemize}
\item $\rank(\vv V)>2\l$. Otherwise, $\vv W=\vv V$ is the claimed subconfiguration. 
\item $\DD(\vv V)>0$. Otherwise, the result follows from Corollary~\ref{cor:Lawrence}.
\item $\vv V$ is irreducible. Otherwise, we remove all copies of $\veczero$, then apply this proposition for irreducible configurations to find a suitable $\vv W$, and finally plug the copies of $\veczero$ back into $\vv W$. By Lemma~\ref{lem:DDirreducible}, this does not change any of $\DD(\vv W)$, $\DD(\vv V\setminus \vv W)$ or $\DD(\vv V)$.
\end{itemize}

Setting $s=\l$ in Corollary~\ref{cor:nocs}, we know that there is a subconfiguration~$\vv W_0$ of~$\vv V$ such that
\begin{eqnarray}
 \vv W_0 &=& \lin(\vv W_0)\cap \vv V, \label{eq:linW0}\\
 \rank(\vv W_0) &=& \l, \label{eq:rankW0}\\
 \l &=& \DD(\vv W_0)+\DD(\vv V/\vv W_0) \quad\text{ and } \label{eq:DDVprop} \\
 \DD(\vv W_0) &=& \l\, . \label{eq:DDW0}
\end{eqnarray}
Here, the last equation follows because $\DD(\vv W_0)\geq\l$ by the last relation of Corollary~\ref{cor:nocs}, and $\DD(\vv W_0)\leq\l$ trivially because of Corollary~\ref{cor:subspacequotient}.

The combination of \eqref{eq:DDVprop} and \eqref{eq:DDW0} implies that $\DD(\vv V/\vv W_0)=0$. 
Hence, $\vv V/\vv W_0$ is centrally symmetric by Corollary~\ref{cor:Lawrence}. 
Corollary~\ref{cor:Lawrencedecomposition} then finds a \codegreeG decomposition of~$\vv V/\vv W_0$
into $\codegG(\vv V/\vv W_0)$ many pairs of antipodal vectors, where we use that $\vv V/\vv W_0$ is irreducible by~\eqref{eq:linW0} and Lemma~\ref{lem:quotientirreducible}. 
As long as possible, we group collinear pairs among these factors together, and end up with a decomposition of $\vv V/\vv W_0$ into $m$ centrally symmetric factors $\vv U_1,\dots,\vv U_m$ of rank~$1$ that fulfill $\lin(\vv U_i)\cap (\vv V/\vv W_0)=\vv U_i$. By construction, $\rank(\vv V/\vv W_0)\le m\le\codegG(\vv V/\vv W_0)$ and
\begin{equation}\label{eq:decompV/W0}\codegG(\vv V/\vv W_0)=\sum_{i=1}^m\codegG(\vv U_i).\end{equation}
Moreover, $\DD(\vv U_i)=0$ and $\DD\!\big((\vv V/\vv W_0)/\vv U_i\big)=0$ for all~$i$, because all these configurations are centrally symmetric. Thus, trivially,
\begin{equation}\label{eq:decompV/W}0=\DD(\vv V/\vv W_0)=\DD(\vv U_i)+\DD\!\big((\vv V/\vv W_0)/\vv U_i\big).\end{equation}
Next, for $1\leq i\leq m$ we define $\vv W_i\subset \vv V$ to be the subconfiguration of~$\vv V$ that fulfills $\vv W_0\subset\vv W_i$ and $\vv W_i/\vv W_0=\vv U_i$. That is, $\vv W_i=\vv W_0\cup \vv U_i$ after identifying the elements in $\vv U_i\subset \vv V/\vv W_0$ with the corresponding elements in $\vv V\setminus\vv W_0$.
Our next step is to use~\eqref{eq:DDVprop}~and~\eqref{eq:decompV/W} to apply Lemma~\ref{lem:iterativeextremal}. 
Specifically, the equation \eqref{eq:iterativeextremalsum} of that lemma with $\vv U=\vv U_i$, $\vv W=\vv W_0$ and $\vv U\cup \vv W=\vv W_i$ tells us that \begin{equation}\label{eq:DDWi}
\DD(\vv W_i)=\DD(\vv U_i)+\DD(\vv W_0).
\end{equation}
Consequently, $\DD(\vv W_i)=\l$ since $\DD(\vv U_i)=0$ and $\DD(\vv W_0)=\l$ by~\eqref{eq:DDW0}.

Now set $I=\{i \mid 1\le i\le m, \; \DD(\vv W_i\setminus \vv W_0)>0\}$. 
We claim that 
\begin{equation}
\label{eq:rankcondition}
 \rank{\Big({\bigcup\nolimits_{i\in I}{\vv U_i}}\Big)} \le \l\ .
\end{equation}
Observe that we will be done once we have seen that~\eqref{eq:rankcondition} is true. 
Indeed, if we define $\vv W:=\vv W_0\cup\bigcup_{i\in I}\vv W_i$ then,
\begin{eqnarray}
  \rank(\vv W) 
  & \leq & 
  \rank(\vv W_0)+\rank\Big(\bigcup\nolimits_{i\in I}\vv U_i\Big) \notag\\ 
  & \stackrel{\eqref{eq:rankW0}}{=} &
  \l+\rank\Big(\bigcup\nolimits_{i\in I}\vv U_i\Big)\ \stackrel{\eqref{eq:rankcondition}}{\leq} \ 2\l\ .\label{eq:rankinequality}
\end{eqnarray}
Hence, \eqref{eq:rankcondition} implies that $\rank(\vv W)\leq 2\l=2\DD(\vv V)$. Moreover, the definition of~$I$ says that, if $j\notin I$ then $\DD(\vv W_j\setminus \vv W_0)=0$ and hence $\vv W_j\setminus \vv W_0$~is centrally symmetric. Then 
observing that $\vv V\setminus \vv W=\bigcup_{j\notin I}\vv W_j\setminus\vv W_0$ shows that $\vv V\setminus \vv W$ is centrally symmetric. This implies $\DD(\vv V\setminus \vv W)=0$, and moreover $\DD(\vv W)=\DD\big(\vv V\setminus(\vv V\setminus\vv W)\big) = \DD(\vv V)$, because removing a centrally symmetric subconfiguration does not change the covector discrepancy by Lemma~\ref{lem:removecs}. This concludes our proof modulo~\eqref{eq:rankcondition}.\\

We prove \eqref{eq:rankcondition} by contradiction, assuming that $\{1,\dots,{\l+1}\}\subseteq I$ (and in particular, that $m\ge \l+1$),  and that $\vv U_1,\dots,\vv U_{\l+1}$ are linearly independent. Put differently, we assume that $\rank\bigcup_{i=1}^{\l+1}\vv U_i=\l+1$ and that $\DD(\vv W_i\setminus \vv W_0)>0$ for $i\leq \l+1$. 
This last property and~\eqref{eq:linW0} let us use Lemma~\ref{lem:strongnocs} (with $\vv V=\vv W_i$ and $\vv W=\vv W_0$) for each $1\le i\leq \l+1$ to find some $\vv T_i\subset \vv W_i\setminus \vv W_0$ of rank~$1$ with 
\begin{align}\label{eq:DDTige1}\DD(\vv T_i)&\geq 1\qquad\text{and}\\ 
\label{eq:DDWiTi}\DD(\vv W_i)&=\DD(\vv T_i) + \DD(\vv W_i/\vv T_i).
\end{align}
Moreover, $\vv T_i\not\subset\lin(\vv W_0)$ because $\lin(\vv W_0)\cap \vv V=\vv W_0$.
Additionally, we let Lemma~\ref{lem:extremalconfigurations} choose some $\vv T_0\subset \vv W_0$ with $\rank(\vv T_0)=\rank(\vv W_0)-1$ such that 
\[
  \DD(\vv T_0)=\DD(\vv W_0)-\DD(\vv W_0/\vv T_0).
\]
(If $\rank(\vv W_0)=1$, then $\vv T_0$ is empty.) 
By Lemma~\ref{lem:subspacequotientequality}, the displayed equation is equivalent to
\begin{equation}\label{eq:DDT0}
\codegG(\vv W_0/\vv T_0) + \codegG(\vv T_0)=\codegG(\vv W_0).
\end{equation} 
Observe that if $\rank(\vv T_0)\ge1$, then $\vv T_0$ is also linearly independent of the remaining $\vv T_i$, since these were already independent in $\vv V/\vv W_0$.

Now, let ${\vvh H}'$ be a hyperplane such that $\vv T_i\subset {\vvh H}'$ for $0\leq i\leq \l+1$ and that is in general position with respect to the remainder of~$\vv V$. (If $\vv T_0$~is empty, we only require that $\vvh H'$ does not contain $\vv W_0$.) Such an $\vvh H'$ exists because $\vv T_1,\dots,\vv T_m$ (and perhaps also~$\vv T_0$) are linearly independent, and $\rank(\vv V)>2\l$. Since $\vvh H'\cap\vv W_0=\vv T_0$, we can orient ${\vvh H}'$ in such a way that 
\begin{equation}\label{eq:H'capW0}|{{\vvh H}'}^+\cap \vv W_0|=\codegG(\vv W_0/\vv T_0).\end{equation} 
The reason for this is that the contraction $\vv W_0/\vv T_0$ is of rank~$1$ because $\rank(\vv T_0)=\rank(\vv W_0)-1$, and one of the two possible orientations of~${\vvh H}'$ must attain the dual codegree of $\vv W_0/\vv T_0$. 

Our next ingredient are hyperplanes $\vv H_0,\dots,\vv H_{m}$  with $\bigcup_{j\ne i}\vv T_j\subset \vv H_i$, such that
\begin{eqnarray}\label{eq:HicapTi}
  |{\vvh H}_i^+ \cap \vv T_i|
  &=&
  \rank(\vv T_i)+\degG(\vv T_i) \notag
  \\ &\stackrel{\eqref{eq:r+dd=DD+kk}}{=}&
  \codegG(\vv T_i)+\DD(\vv T_i) \qquad\text{for }i\geq 1,
\end{eqnarray}
and
\begin{equation}\label{eq:H0capT0}
  |{\vvh H}_0^+ \cap \vv T_0|=\codegG(\vv T_0).
\end{equation} 
To find them, start from hyperplanes $\vvh H_i'$ in $\lin(\vv T_i)$ fulfilling \eqref{eq:HicapTi} and \eqref{eq:H0capT0} respectively. They achieve the codegree of $\vv W_i$ in ${\vvh H}_i^-$ for $1\leq i \leq m$, and the codegree of $\vv W_0$ in ${\vvh H}_0^+$ (cf.~Lemma~\ref{lem:DDrelatestodd}).
 Then choose $\vvh H_i$ to be any hyperplane that does not contain $\lin(\vv T_i)$ and goes through $\vvh H_i'$ and the~$\vv T_j$ with $j\neq i$; this succeeds by the linear independence of the $\vv T_i$.

Now let ${\vvh H}={\vvh H}'\circ{\vvh H_0}\circ{\vvh H_1}\dots\circ{\vvh H}_{\l+1}$. Perturbing $\vv H$ even further if necessary, we can assume that $\vv V\cap {\vvh H}=\emptyset$ because $\vv V$ is irreducible. 
Observe that
${\vvh H}$ is chosen in such a way that ${\vvh H}^+$ contains as many vectors of~$\vv T_i$ as possible for $1\leq i\leq \l+1$, and as few as possible from~$\vv W$; cf.~Figure~\ref{fig:2DD}.
In particular, since $\vv T_0\subset {\vvh H}'$, by \eqref{eq:H'capW0}, \eqref{eq:H0capT0} and \eqref{eq:DDT0} we see that
\begin{eqnarray}
  |{\vvh H}^+\cap \vv W_0|
  &=&
  |{{\vvh H}'}^+\cap \vv W_0|+|{\vvh H}_0^+\cap \vv T_0|\notag
  \\ &\stackrel{\eqref{eq:H'capW0},\eqref{eq:H0capT0}}{=}&
  \codegG(\vv W_0/\vv T_0)+\codegG(\vv T_0)
  \\ &\stackrel{\eqref{eq:DDT0}}{=}&
  \codegG(\vv W_0).\label{eq:HcapW0}
\end{eqnarray}
Moreover, observe that 
\begin{equation}
 \label{eq:H'capWi}
|{\vvh H'}^+\cap \vv W_i|\geq  \codegG(\vv W_i/\vv T_i)
\qquad\text{for each } 1\leq i\leq \l+1
\end{equation}
since $\vv T_i\subset {\vvh H}'$ and by the definition of the dual codegree. Combining this equation with \eqref{eq:HicapTi} and \eqref{eq:HcapW0}, and then \eqref{eq:DDWiTi} and \eqref{eq:DDWi} in their \codegreeG formulation (using Lemma~\ref{lem:subspacequotientequality}), and finally \eqref{eq:DDTige1} yields, for $i\le\l+1$,
\begin{eqnarray}
  |{\vvh H}^+\cap \vv W_i\setminus\vv  W_0|
  &=& 
  |{\vvh H'}^+\cap \vv W_i|+|{\vvh H}_i^+\cap \vv T_i|-|{\vvh H}^+\cap \vv W_0|\notag\\
  &\stackrel{\eqref{eq:H'capWi},\eqref{eq:HicapTi},\eqref{eq:HcapW0}}{\geq}&
  \codegG(\vv W_i/\vv T_i)+\codegG(\vv T_i)+ \notag
  \\ && \qquad{}+\DD(\vv T_i)-\codegG(\vv W_0)\notag\\
  &\stackrel{ \eqref{eq:DDWiTi}}{=}&
  \codegG(\vv W_i)-\codegG(\vv W_0)+\DD(\vv T_i)\notag\\
  &\stackrel{\eqref{eq:DDWi}}{=}& 
  \codegG(\vv U_i)+\DD(\vv T_i) \notag
  \\ &\stackrel{\eqref{eq:DDTige1}}{\geq}&
  \codegG(\vv U_i)+1.\label{eq:HcapWi1}
\end{eqnarray}
Finally, for $i>\l+1$, we use again the definition of \codegreeG to obtain
\begin{equation}
 \label{eq:HcapWi}
|{\vvh H}^+\cap \vv W_i|\geq  \codegG(\vv W_i).
\end{equation}
Now \eqref{eq:HcapW0} and the \codegreeG version of \eqref{eq:DDWi} yield, for $i>\l+1$,
\begin{eqnarray}
  |{\vvh H}^+\cap \vv W_i\setminus\vv  W_0|
  &=&
  |{\vvh H}^+\cap \vv W_i|-|{\vvh H}^+\cap \vv W_0|\notag\\
  &\stackrel{\eqref{eq:HcapWi},\eqref{eq:HcapW0}}{\geq}& 
  \codegG(\vv W_i)-\codegG(\vv W_0) \notag
  \\ &\stackrel{\eqref{eq:DDWi}}{=}&
  \codegG(\vv U_i).\label{eq:HcapWi2}
\end{eqnarray}

Summing up, \eqref{eq:HcapW0}, \eqref{eq:HcapWi1} and \eqref{eq:HcapWi2}, combined with \eqref{eq:decompV/W0} and \eqref{eq:DDVprop}, and then \eqref{eq:r+dd=DD+kk} from Lemma~\ref{lem:DDrelatestodd} (which we can apply since we assumed $\vv V$ to be irreducible) show that   
\begin{eqnarray}
  \big|{\vvh H}^+\cap \vv V\big|
  &=&
  \big|{\vvh H}^+\cap \vv W_0\big|
   + \sum_{i=1}^{\l+1} \big|{\vvh H}^+\cap \vv W_i\setminus \vv W_0\big|+ 
  \notag\\ &&\qquad
  {}+\sum_{i=\l+2}^{m} |{\vvh H}^+\cap \vv W_i\setminus \vv W_0|\notag\\
  &\stackrel{\eqref{eq:HcapW0},\eqref{eq:HcapWi1},\eqref{eq:HcapWi2}}{\geq}&
  \codegG(\vv W_0)+ \sum_{i=1}^{\l+1} (\codegG(\vv U_i)+1)+
  \notag\\ &&\qquad
  {}+\sum_{i=\l+2}^{m} \codegG(\vv U_i)\notag\\
  &=&
  \codegG(\vv W_0)+\l+1+\sum_{i=1}^{m} \codegG(\vv U_i)\notag\\
  &\stackrel{\eqref{eq:decompV/W0}}{=}&
  \codegG(\vv W_0)+\codegG(\vv V/\vv W_0)+\l+1\notag\\
  &\stackrel{\eqref{eq:DDVprop}}{=}&
  \codegG(\vv V)+\l+1=
  \codegG(\vv V)+\DD(\vv V)+1
  \notag\\ &\stackrel{\eqref{eq:r+dd=DD+kk}}{=}&
  r+\degG(\vv V)+1,\label{eq:lasteq}
\end{eqnarray}
which contradicts the fact that by definition $\max_{\vvh H}\big|{\vvh H}^+\cap\vv V\big|=r+\degG(\vv V)$, and therefore proves~\eqref{eq:rankcondition}.
\end{proof}

\iftoggle{bwprint}{%
\begin{figure}[htpb]
\centering
 \subbottom[$\vv V$]{\label{sfig:2DD_V}\includegraphics[width=.30\linewidth]{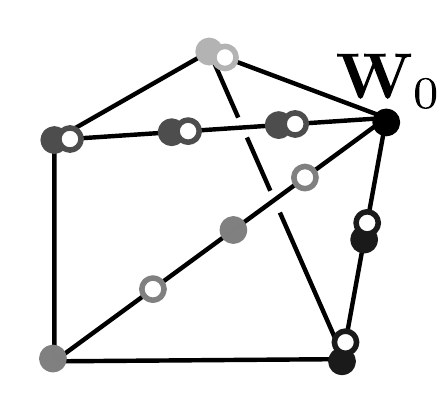}}\hfill
 \subbottom[$\vv V/\vv W_0$]{\label{sfig:2DD_VmodW0}\includegraphics[width=.30\linewidth]{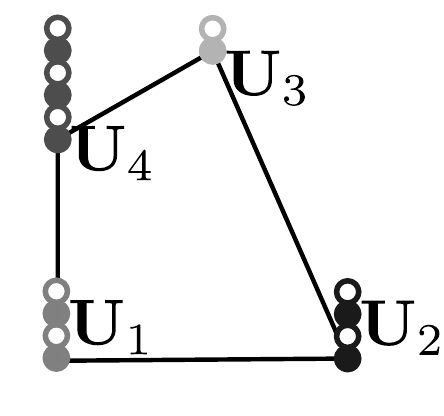}}\hfill
 \subbottom[$\vv W_1,\dots, \vv W_4$]{\label{sfig:2DD_Wi}\includegraphics[width=.30\linewidth]{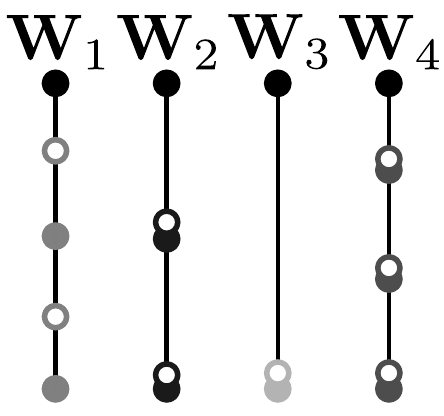}}\newline\newline
 \subbottom[$\vv V'$]{\label{sfig:2DD_VV}\includegraphics[width=.30\linewidth]{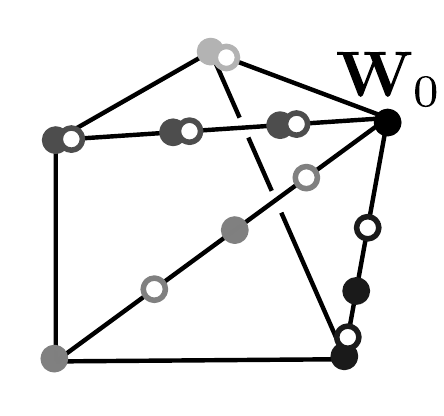}}\hfill
 \subbottom[$\vv T_1,\vv T_2$]{\label{sfig:2DD_Ti}\includegraphics[width=.15\linewidth]{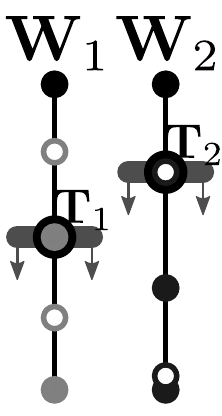}}\hfill
 \subbottom[$\vvh H'$]{\label{sfig:2DD_HH}\includegraphics[width=.30\linewidth]{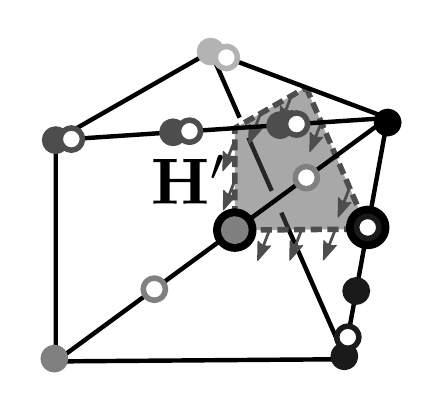}}\hfill
 \subbottom[$\vvh H_1',\vvh H_2'$]{\label{sfig:2DD_Hi}\includegraphics[width=.15\linewidth]{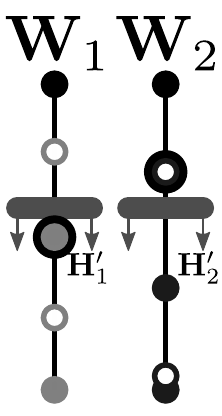}}\newline\newline
 \subbottom[$\vvh H_1$]{\label{sfig:2DD_H1}\includegraphics[width=.30\linewidth]{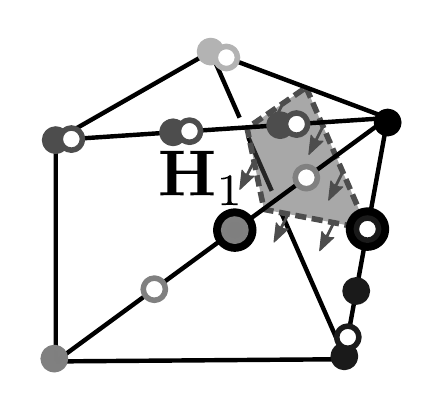}}\hfill
 \subbottom[$\vvh H_2$]{\label{sfig:2DD_H2}\includegraphics[width=.30\linewidth]{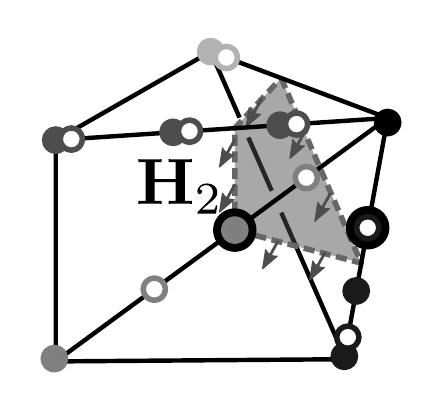}}\hfill
\subbottom[$\vvh H$]{\label{sfig:2DD_H}\includegraphics[width=.30\linewidth]{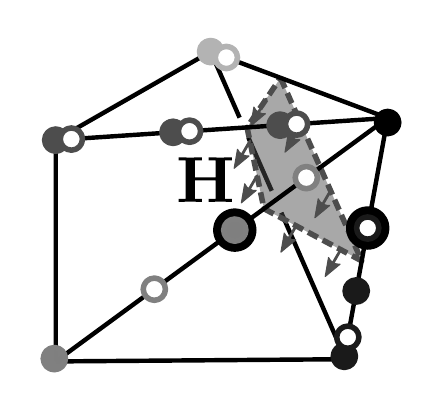}}
 \caption[Illustrating the proof of Proposition~\ref{prop:2DD}.]{Illustrating the proof of Proposition~\ref{prop:2DD}. Overlapping circles represent points that have the same coordinates.}
\label{fig:2DD}
\end{figure}
}{%
\begin{figure}[htpb]
\centering
 \subbottom[$\vv V$]{\label{sfig:2DD_V}\includegraphics[width=.30\linewidth]{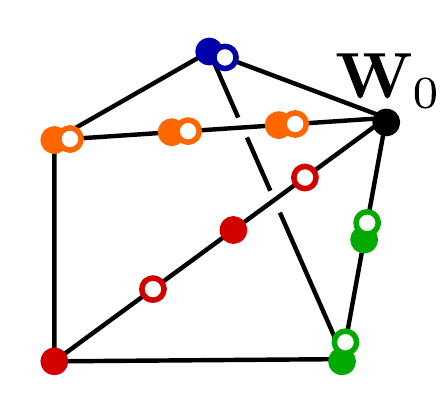}}\hfill
 \subbottom[$\vv V/\vv W_0$]{\label{sfig:2DD_VmodW0}\includegraphics[width=.30\linewidth]{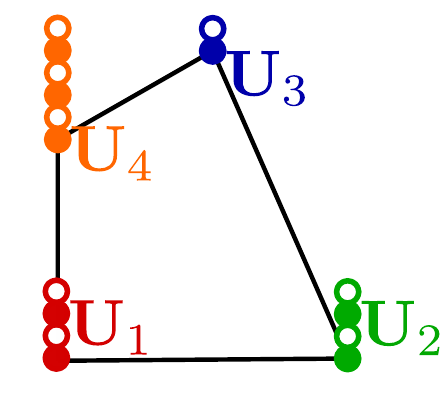}}\hfill
 \subbottom[$\vv W_1,\dots, \vv W_4$]{\label{sfig:2DD_Wi}\includegraphics[width=.30\linewidth]{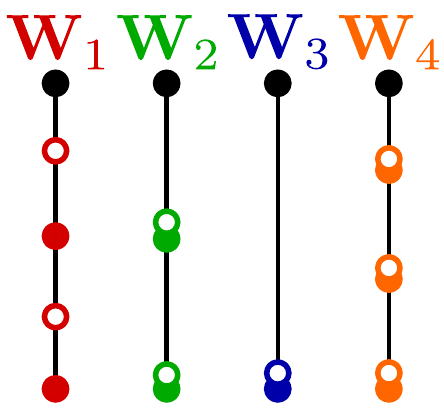}}\newline\newline
 \subbottom[$\vv V'$]{\label{sfig:2DD_VV}\includegraphics[width=.30\linewidth]{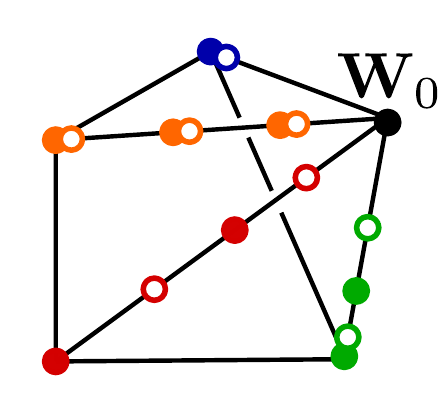}}\hfill
 \subbottom[$\vv T_1,\vv T_2$]{\label{sfig:2DD_Ti}\includegraphics[width=.15\linewidth]{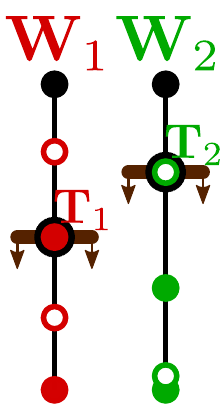}}\hfill
 \subbottom[$\vvh H'$]{\label{sfig:2DD_HH}\includegraphics[width=.30\linewidth]{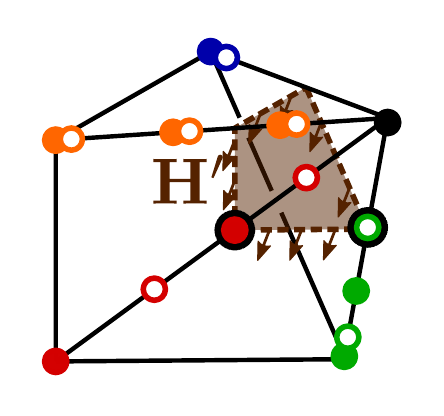}}\hfill
 \subbottom[$\vvh H_1',\vvh H_2'$]{\label{sfig:2DD_Hi}\includegraphics[width=.15\linewidth]{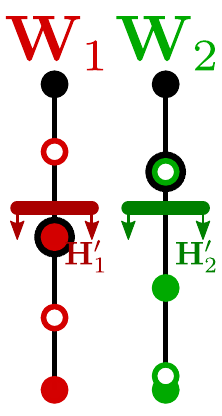}}\newline\newline
 \subbottom[$\vvh H_1$]{\label{sfig:2DD_H1}\includegraphics[width=.30\linewidth]{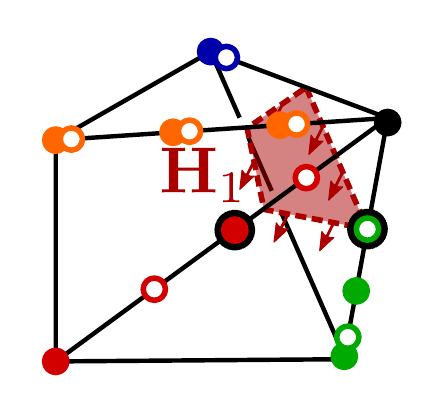}}\hfill
 \subbottom[$\vvh H_2$]{\label{sfig:2DD_H2}\includegraphics[width=.30\linewidth]{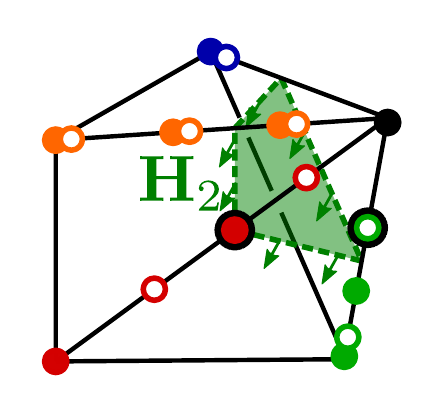}}\hfill
\subbottom[$\vvh H$]{\label{sfig:2DD_H}\includegraphics[width=.30\linewidth]{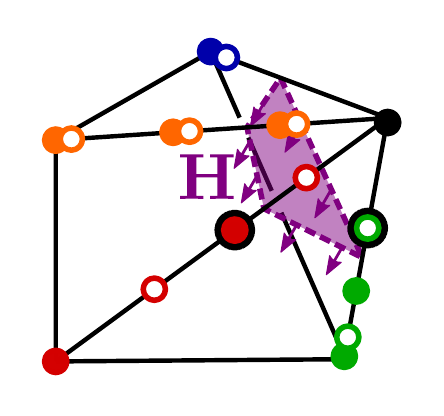}}
 \caption[Illustrating the proof of Proposition~\ref{prop:2DD}.]{Illustrating the proof of Proposition~\ref{prop:2DD}. Overlapping circles represent points that have the same coordinates.}
\label{fig:2DD}
\end{figure}
}

\begin{example}\label{ex:2DD}
 Proposition~\ref{prop:2DD} is illustrated in Figure~\ref{fig:2DD}. In~\ref{sfig:2DD_V}, there is an affine diagram of a vector configuration $\vv V$ that has $17$~elements, rank~$4$, dual degree~$5$, and discrepancy $\DD(\vv V)=1$. A subconfiguration $\vv W_0$ of rank~$1$ with one element is also depicted in~\ref{sfig:2DD_V} (the ``apex''). 
 The contraction $\vv V/\vv W_0$, in \ref{sfig:2DD_VmodW0}, is centrally symmetric and factors into~$\vv U_1$, $\vv U_2$, $\vv U_3$ and~$\vv U_4$. The corresponding sets~$\vv W_i$, which are all subsets of $\vv V$ that contain $\vv W_0$, are shown independently in~\ref{sfig:2DD_Wi}. Observe that $I=\{1\}$, since $\DD(\vv W_1\setminus \vv W_0)=2$ but $\DD(\vv W_j\setminus \vv W_0)=0$ for $j=2,3,4$.

 Figure~\ref{sfig:2DD_VV} shows a slightly modified configuration $\vv V'$ where one of the elements of $\vv W_2$ has been perturbed in such a way that $\{1,2\}\subset I$, \ie $\DD(\vv W_i\setminus \vv W_0)>0$ for $i=1,2$. If $\DD(\vv V')=1$, this would contradict~\eqref{eq:rankcondition}. We will follow the proof of Proposition~\ref{prop:2DD} to show that $\degG(\vv V')\geq 6$ and hence $\DD(\vv V')\geq 3$.

 In \ref{sfig:2DD_Ti}, two subspaces $\vv T_1$ and $\vv T_2$ of rank~$1$ have been chosen so that $\DD(\vv T_1)>0$ and $\DD(\vv T_2)>0$. They are highlighted with double circles. 
 A hyperplane $\vvh H'$ that goes through~$\vv T_1$ and~$\vv T_2$ is shown in \ref{sfig:2DD_HH} (its intersection with $\vv W_1$ and $\vv W_2$ is displayed in~\ref{sfig:2DD_Ti}). The hyperplanes $\vvh H_1$ and $\vvh H_2$ are shown in \ref{sfig:2DD_H1} and \ref{sfig:2DD_H2}, while their respective intersections with $\lin(\vv W_1)$ and $\lin(\vv W_2)$ (the \codegreeG achieving hyperplanes $\vvh H_1'$ and $\vvh H_2'$) can be seen in \ref{sfig:2DD_Hi}. Finally, the hyperplane $\vvh H$ is shown in \ref{sfig:2DD_H}. It fulfills $|\vvh H^+\cap \vv V'|=10$, which shows that $\degG(\vv V')\geq 6$.
\end{example}

\smallskip
As a corollary of Proposition~\ref{prop:2DD} we get the following theorem, which certifies the existence of non-trivial \codegreeG decompositions when the rank is large with respect to the discrepancy. In the primal setting, it finds a codegree decomposition when the dimension is large with respect to the degree and the number of elements.
\begin{theorem}\label{thm:2DD}
Any irreducible vector configuration $\vv V$ with $\DD(\vv V)>0$ admits a \codegreeG decomposition of length at least 
\[
   r+1-2\DD(\vv V) \ = \ 2(d+1-2\dd)-r+1 \ = \ 3d+4-4\dd -n.
\]
where $r:=\rank(\vv V)$, $\dd:=\degG(\vv V)$ and $d:=|\vv V|-r-1$.
\end{theorem}
\begin{proof}
By Proposition~\ref{prop:2DD}, there is a subspace $\vv W$ of rank $\leq 2\DD(\vv V)$ such that $\DD(\vv V\setminus \vv W)=0$. Therefore, $\vv V\setminus \vv W$ admits a \codegreeG decomposition of length 
\begin{eqnarray}
  \rank(\vv V\setminus \vv W)+\degG(\vv V\setminus\vv W)
  &\geq&
  \rank(\vv V\setminus \vv W)\notag
  \\ &\geq& 
  r-\rank(\vv W)\notag
  \\ &\geq& r-2\DD(\vv V),\label{eq:ineq3}
\end{eqnarray}
by Corollary~\ref{cor:Lawrencedecomposition}. Observe that if $\rank(\vv V\setminus \vv W)>r-2\DD(\vv V)$, then we are done, since the decomposition of $\vv V\setminus \vv W$ would already be of the desired length.

The situation $\rank(\vv V\setminus \vv W)=r-2\DD(\vv V)$ is only possible when equality holds in \eqref{eq:ineq3}, which means that $\rank(\vv W)=2\DD(\vv V)$.

Since $\DD(\vv W)=\DD(\vv V)>0$, $\vv W$ is not empty, and we can extend the decomposition of $\vv V\setminus \vv W$ to a \codegreeG decomposition of $\vv V$ by adding~$\vv W$ as the last factor. Observe that $\codegG(\vv W)\geq 1$ by \eqref{eq:achievecodegstar}, since $|\vv W|\geq \rank(\vv W)=2\DD(\vv W)$.

Finally, we use~\eqref{eq:DDrelatestodd} to get the desired expression.
\end{proof}

Corollary~\ref{cor:Lawrence} dealt with configurations of zero discrepancy. In view of Proposition~\ref{prop:easybound}, when these are irreducible of dual degree $\dd$, they have the maximal possible number of elements, $2r+2\dd$. It turns out that Theorem~\ref{thm:2DD} implies Conjecture~\ref{conj:strongd+1-2dd} for irreducible configurations with one element less:
\begin{corollary}\label{cor:minimal+1}
 Let $\vv V$ be an irreducible vector configuration of dual degree~$\dd$ with $n=r+d+1=2r+2\dd-1$ elements. Then $\vv V$ admits a \codegreeG decomposition of length at least $r-1=d+1-2\dd$. 
\end{corollary}

\begin{proof}
 $\DD(\vv V) = 2r+2\dd - n = 1$ by~\eqref{eq:DDrelatestodd2}. By Theorem~\ref{thm:2DD}, $\vv V$ admits a \codegreeG decomposition of length at least $r+1-2\DD(\vv V)=r-1$.
\end{proof}

\subsection{\texorpdfstring{Configurations of rank $4$}{Configurations of rank 4}}

Finally, we can go one step further and prove the conjecture for configurations with $\DD(\vv V)=2$. Observe that for such configurations, Theorem~\ref{thm:2DD} only certifies decompositions of length~$r-3$, but Conjecture~\ref{conj:strongd+1-2dd} predicts configurations of length~$r-2$. In Theorem~\ref{thm:DD=2}, we find the missing factor.

To prove it, we need a lemma that is based on the following variation of the Sylvester-Gallai Theorem.

\begin{lemma}[{Sharp Dual of Sylvester's Theorem, \cite[Theorem 42]{Lenchner2008}}]
In any arrangement of $n\ge3$ lines in the affine plane, not all of them parallel and not all of them passing through a common point, there must be at least one (finite) point contained in exactly two lines of the arrangement.
\end{lemma}

We reformulate this in the following way:

\begin{lemma}\label{lem:strongSylGal}
Let $\vv a\in\vv A$ be member of a configuration of at least $3$~points in $\RR^2$ without repetition. If $\vv A\setminus \vv a$ is not contained in a common line, there is a line~$\vv L$ such that $|\vv L\cap \vv A|=2$ and $\vv a\notin \vv L$.
\end{lemma}

\begin{proof}
 Using projective duality, we translate $\vv A$ into an arrangement~$\vv A^*$ (not to be confused with $\vv A$'s Gale dual $\Gale{\vv A}$) of lines in the projective plane~$\PP^2$, and dehomogenize by taking~$\vv a^*$, the line corresponding to~$\vv a$, to be the line at infinity. The condition that no line contains $\vv A\setminus \vv a$ implies that the lines in the affine arrangement $\vv A^*\setminus \vv a^*\subset\RR^2$ are not all parallel, and do not all pass through a common point.

Therefore, we can apply the Sharp Dual of Sylvester's Theorem to find a (finite) point contained in exactly two lines of the arrangement. The dual of this point is the desired line $\vv L$; it does not contain $\vv a$, since the intersection point was taken to be finite.
\end{proof}

\begin{lemma}\label{lem:chooseinrank3}
 Let $\vv V$ be a vector configuration of rank~$3$, let $\vv V'\subset \vv V$ be any subconfiguration obtained from~$\vv V$ by deleting pairs 
 of antipodal vectors, 
 $\{\vv v, -\gl\vv v\}$ with $\gl>0$, 
 as long as possible; and let $\vv W'\subseteq \vv V'$ be a subconfiguration of rank~$1$. If $\vv W=\lin(\vv W')\cap \vv V$, one of the following holds:
\begin{enumerate}[\upshape(a)]
 \item\label{it0:lem:chooseinrank3} Either $\rank(\vv V')=2$; 
 \item\label{it1:lem:chooseinrank3} or $\rank(\vv V')=3$ 
and $\DD(\vv V)=\DD(\vv W)+\DD(\vv V\setminus \vv W)$; 
 \item\label{it2:lem:chooseinrank3} or $\rank(\vv V')=3$ and there is a subconfiguration $\vv T\subset \vv V$ of rank $2$ such that $\DD(\vv T)\geq 2$, $\lin(\vv T)\cap \vv V=\vv T$ and $\lin(\vv T)\cap \lin(\vv W)=\veczero$.
\end{enumerate}
\end{lemma}

\begin{proof}
 We may assume that $\rank(\vv V')=3$, otherwise we are in case~\eqref{it0:lem:chooseinrank3}. Two things can happen:
 If $\rank(\vv V'\setminus\vv W')=2$, then $\vv V'$ must be the direct sum of~$\vv V'\setminus \vv W'$ and~$\vv W'$. Thus, $\DD(\vv V')=\DD(\vv W')+\DD(\vv V'\setminus\vv W')$, which implies that $\DD(\vv V)=\DD(\vv W)+\DD(\vv V\setminus\vv W)$ because adding pairs of antipodal vectors does not change the covector discrepancy  by Lemma~\ref{lem:removecs}, and we are in case~\eqref{it1:lem:chooseinrank3}.

Otherwise, $\rank(\vv V'\setminus\vv W')=3$. Let $\vv X$ be an affine diagram of $\vv V'$, where the point~$\vv x$ represents some point in~$\vv W'$.

Pretend for a moment that there are no repeated points and that the points are not colored. Since $\rank(\vv V'\setminus\vv W')=3$, the points in $\vv X\setminus \vv x$ are not contained in a common line. Therefore, by Lemma~\ref{lem:strongSylGal}, there is a line~$\vv L$ containing exactly two points of $\vv X$ and that avoids~$\vv x$. 

If we stop pretending, $\vv L\cap \vv X$ represents a subconfiguration $\vv T'\subset \vv V'$ such that $\lin(\vv T')\cap\lin(\vv W)=\veczero$ (because $\vv x$ was not in $\vv L$) and such that $\vv T'$ is the direct sum of two configurations $\vv T_1$ and $\vv T_2$ of rank~$1$. Since $\vv V'$~did not contain antipodal pairs, $\DD(\vv T_i)\geq1$ for $i\in \{1,2\}$. 
Hence $\DD(\vv T')=\DD(\vv T_1)+\DD(\vv T_2)\geq2$ (where the first equality is because $\vv W'=\vv T_1\oplus\vv T_2$). Our desired configuration is $\vv T=\lin(\vv T')\cap \vv V$, which still fulfills $\DD(\vv T)\geq 2$ by Lemma~\ref{lem:removecs}. 
\end{proof}

Before proving the theorem, we still need an easy lemma that follows from the fact that every covector of $\vv V/\vv W$ is also a covector of $\vv V\setminus \vv W$. 

\begin{lemma}\label{lem:DDsubtractcontract}
For any configuration $\vv V$ and any subconfiguration $\vv W\subseteq \vv V$, $\DD(\vv V\setminus \vv W)\geq \DD(\vv V/\vv W)$.\qed
\end{lemma}

\begin{theorem}\label{thm:DD=2}
Any irreducible vector configuration $\vv V$ of rank $r$ and discrepancy $\DD(\vv V)\leq 2$ admits a \codegreeG decomposition of length $\ge r-\DD(\vv V)$. 
\end{theorem}

Our proof is built on top of the proof of Proposition~\ref{prop:2DD}. Instead of repeating the (several) coincident parts, we present a series of ``patches'' on the details that differ among the two proofs.

\begin{proof}[Proof sketch.]
By Corollary~\ref{cor:minimal+1}, it suffices to consider $\DD(\vv V)=2$. 
We will show that in this case, $\vv V$ is centrally symmetric except for either:
\begin{enumerate}[(A)]
 \item\label{it:1Wleq3} a subconfiguration of rank~$\leq 3$. Then $\vv V$ admits a decomposition into one piece $\vv W$ of rank~$\leq 3$ and at least $r-3$ centrally symmetric pieces (by Corollary~\ref{cor:Lawrencedecomposition} and because $\rank(\vv W\setminus \vv V)\geq r-3$);
 \item\label{it:2Wleq2} or two skew configurations of rank~$\leq2$. Then $\vv V$ admits a decomposition into these two pieces of rank~$\leq 2$ and at least $r-4$ centrally symmetric pieces  analogously.
\end{enumerate}

Since most of the proof is exactly the same as that of Proposition~\ref{prop:2DD}, we just mark those points where the proofs are different. 

Because of Lemma~\ref{lem:nocs}, there is a subconfiguration $\vv W_0$ such that
\begin{eqnarray}
 \rank(\vv W_0) &=& 1,\label{eq:rW0=1}\\
 \vv W_0 &=& \lin(\vv W_0)\cap \vv V,\\
 \DD(\vv V)&=& \DD(\vv W_0)+\DD(\vv V/\vv W_0) \quad\text{ and } \label{eq:DDVDDW=DDVmodW0}\\
 \DD(\vv W_0) &\geq& 1.  
\end{eqnarray}
Comparing to Proposition~\ref{prop:2DD}, here we insist on having $\rank(\vv W_0)=1$, while we allow $\DD(\vv W_0)<\DD(\vv V)=2$ as long as $\DD(\vv W_0)>0$.

We distinguish two cases depending on the value of $\DD(\vv W_0)$, which by construction is either $1$ or $2$:

\smallskip

\begin{compactenum}[1)]
 \item If $\DD(\vv W_0)=2$, then the proof exactly parallels that of Proposition~\ref{prop:2DD}. Indeed, $\DD(\vv V/\vv W_0)=0$ and the proof that Equation~\eqref{eq:rankcondition} holds is exactly the same. Then we reach~\eqref{eq:rankinequality} and get
\begin{equation}
  \rank(\vv W) 
  \ \leq\ 
  \rank(\vv W_0)+\rank\bigcup_{i\in I}\vv U_i
  \stackrel{\eqref{eq:rW0=1}}{=} 1+\rank\bigcup_{i\in I}\vv U_i\stackrel{\eqref{eq:rankcondition}}{\leq}3.
\end{equation}
This proves that $\vv V$ is centrally symmetric except for perhaps a subconfiguration $\vv W$ of rank~$\leq3$; hence $\vv V$ fulfills \eqref{it:1Wleq3}.\\

\item If $\DD(\vv W_0)=1$, then $\DD(\vv V/\vv W_0)=1$ by \eqref{eq:DDVDDW=DDVmodW0}. Applying Proposition~\ref{prop:2DD}, we deduce that $\vv V/\vv W_0$ is centrally symmetric except for a subconfiguration $\vv U_1$ of rank~$\leq 2$ with $\DD(\vv U_1)=1$.

Substituting $\vv U_1$ by $\lin(\vv U_1)\cap (\vv V/\vv W_0)$ if necessary, we can assume that $\lin(\vv U_1)\cap(\vv V/\vv W_0)=\vv U_1$. Then we can complete $\vv U_1$ to a decomposition of $\vv V/\vv W_0$ into factors $\vv U_1,\dots,\vv U_m$ fulfilling 
\begin{align*}
 \vv U_i&=\lin(\vv U_i)\cap (\vv V/\vv W_0) &\text{ for } 1\leq i\leq m,\\
 \DD(\vv V/\vv W_0)&=\DD(\vv U_i)+\DD((\vv V/\vv U_i)/\vv W_0) &\text{ for } 1\leq i\leq m,\\
 \rank(\vv U_i)&=1, \quad \DD(\vv U_i)=0 &\text{ for }2\leq i\leq m, \\
\rank(\vv U_1) &\leq2, \quad \DD(\vv U_1)=1.
\end{align*}
This is achieved by grouping together collinear pairs among the centrally symmetric vectors. 

\smallskip

Next, for $1\leq i\leq m$ we define $\vv W_i\subset \vv V$ to be the subconfiguration of~$\vv V$ that fulfills $\vv W_0\subset\vv W_i$ and $\vv W_i/\vv W_0=\vv U_i$. By Lemma~\ref{lem:iterativeextremal} (with $\vv U=\vv U_i$, $\vv W=\vv W_0$ and $\vv U\cup \vv W=\vv W_i$) we know that
\begin{equation}\label{eq:WiUiW0}
\DD(\vv W_i)=\DD(\vv U_i)+\DD(\vv W_0) \qquad\text{for }1\le i\le m.
\end{equation}
In addition, Lemma~\ref{lem:DDsubtractcontract} implies that
\begin{equation}\label{eq:Wi-W0}
\DD(\vv W_i\setminus\vv W_0)\geq \DD(\vv W_i/\vv W_0)=\DD(\vv U_i) \qquad\text{for }1\le i\le m.
\end{equation}

We define a set of indices $I=\{i \mid 1\le i\le m, \; \DD(\vv W_i\setminus \vv W_0)>\DD(\vv U_i)\}$ like in the proof of Proposition~\ref{prop:2DD} (observe that here we require $\DD(\vv W_i\setminus \vv W_0)>\DD(\vv U_i)$ instead of just $\DD(\vv W_i\setminus \vv W_0)>0$). Again, the key for our success will be to bound the size of $I$.

\smallskip
 
From now on, we distinguish two cases according to $\rank(\vv U_1)$. Here, the case $\rank(\vv U_1)=2$ means that $\vv U_1$ cannot be further subdivided into a subconfiguration of rank~$1$ and centrally symmetric pairs.

\medskip

\begin{compactenum}[{2.}1)]
\item\label{it:rankW1=2} If $\rank(\vv U_1)=2$, so that $\vv W_1$ has rank~$3$, we claim that $|I|\le1$. If this is true, this concludes the proof of this case.  
Indeed, if $I\subseteq\{1\}$, then $\vv V\setminus \vv W_1=\bigcup_{j=2}^m(\vv W_j\setminus \vv W_1)$~is centrally symmetric, because each of these pieces is centrally symmetric: by the definition of~$I$, for $j\neq 1$ we have
\begin{equation}\label{eq:DDWj-W0=0}
  0 \le \DD(\vv W_j\setminus\vv W_0) \le \DD(\vv U_j) = 0.
\end{equation}

Thus, if $I\subseteq\{1\}$ then \eqref{it:1Wleq3} holds. If $I=\{i\}$ with $i\neq 1$, then \eqref{eq:DDWj-W0=0} still holds for $j\notin\{1,i\}$. Therefore, $\vv V$ is centrally symmetric except for $\vv W_1$ and $\vv W_i$. Moreover, $\DD(\vv W_1\setminus \vv W_0)=1$ because
\begin{align*}
 1=\DD(\vv U_1)\stackrel{\eqref{eq:Wi-W0}}{\leq}\DD(\vv W_1\setminus \vv W_0)&\stackrel{1\notin I}{\leq}\DD(\vv U_1)=1,
\end{align*}
and $\DD(\vv W_i)=1$ by
\begin{align*}
 \DD(\vv W_i)&\stackrel{\eqref{eq:WiUiW0}}{=}\DD(\vv U_1)+\DD(\vv  W_0)=1 \qquad\text{for } 2\le i\le m.
\end{align*}

Notice that $\rank(\vv W_i)=2$ because $\rank(\vv U_i)=\rank(\vv W_0)=1$. Moreover, since $\DD(\vv W_1\setminus \vv W_0)=1$, by Proposition~\ref{prop:2DD} $\vv W_1\setminus \vv W_0$ is centrally symmetric except for a subconfiguration~$\vv W'$ of rank~$\leq 2$. Hence, either $\vv W_i$ and $\vv W'$ are skew (\eqref{it:2Wleq2} holds) or they belong to the same subspace of rank~$3$ (hence~\eqref{it:1Wleq3} holds). \\

Therefore, we only need to prove that $|I|\leq 1$. To do so, we assume that there are at least two indices in $I$ in order to reach a contradiction. As in the proof of Proposition~\ref{prop:2DD}, we define some subconfigurations $\vv T_j$ that guide us towards the construction of a hyperplane $\vvh H$ that contradicts the assumption~$|I|>1$. 
\\

The discussion is slightly different depending on whether $1\in I$ or not. 
Since the only distinguished index is $1$, there is no loss of generality assuming that one case is $\{1,2\}\subseteq I$ and the other is $\{2,3\}\subseteq I$.
We temporarily split the proof according to these cases to define some convenient subconfigurations $\vv T_1$, $\vv T_2 $ and $\vv T_3$. The proofs are then joined again to find $\vvh H$.

\medskip

 \begin{compactenum}[{2.1}.1)]
 \item $\{1,2\}\subseteq I$. Set $\vv T_1\subset\vv W_1$ to be such that 

\begin{align}
 \DD(\vv T_1)&= 2,\label{eq:condT1_DD=2}\\
 \DD(\vv T_1)&=\DD(\vv W_1)-\DD(\vv W_1/\vv T_1),\label{eq:condT1_DDdec}\\
 \rank(\vv T_1)&= 2,\quad\text{ and }\label{eq:condT1_rank}\\
 \lin(\vv W_0)\cap\lin(\vv T_1)&=\veczero.\label{eq:condT1_lin}
\end{align}

This subconfiguration can be found by Lemma~\ref{lem:chooseinrank3}. Indeed, setting $\vv V=\vv W_1$ and $\vv W=\vv W_0$ we can see that
$\vv W_1$  fulfills neither condition~\eqref{it0:lem:chooseinrank3} nor \eqref{it1:lem:chooseinrank3} of that lemma:

\smallskip
 
\begin{itemize}
 \item If $\vv W_1$ is centrally symmetric except for a configuration of rank~$\leq 2$ (which must contain $\vv W_0$, because $\vv W_0$ is not centrally symmetric since $\DD(\vv W_0)>0$), then $\vv W_1/\vv W_0$ is centrally symmetric except for a configuration of rank~$1$. But we assumed that this could not happen since otherwise we would be in the case~\ref{it:rankW1=1} ($\rank(\vv U_1)=1$).
 \item If $\DD(\vv W_1\setminus\vv W_0)=\DD(\vv W_1)-\DD(\vv W_0)$, then $1\notin I$ by the definition of $I$, which contradicts $\{1,2\}\subseteq I$.
\end{itemize}

\smallskip

Hence, condition \eqref{it2:lem:chooseinrank3} holds and we can find some rank~2 subconfiguration $\vv T_1\subset\vv W_1$ (called~$\vv T$ in the lemma) fulfilling \eqref{eq:condT1_rank}, \eqref{eq:condT1_lin} and $\DD(\vv T_1)\geq 2$. 
This implies \eqref{eq:condT1_DD=2} and \eqref{eq:condT1_DDdec} just by observing that $\DD(\vv W_1/\vv T_1)\geq 0$ by Corollary~\ref{cor:easyboundDD} and that $\DD(\vv T_1)+\DD(\vv W_1/\vv T_1)\leq \DD(\vv W_1)\leq \DD(\vv V)\leq 2$ by Corollary~\ref{cor:subspacequotient}.

\smallskip

Next, let $\vv T_2\subseteq \vv W_2$ be such that 
\begin{align}
\DD(\vv T_2)&\geq 1,\notag\\ 
\DD(\vv T_2)&=\DD(\vv W_2)-\DD(\vv W_2/\vv T_2),\notag\\
\rank(\vv T_2) &= 1,\quad\text{ and }\notag\\
 \lin(\vv W_0)\cap\lin(\vv T_2)&=\veczero\label{eq:condT2_lin}. 
\end{align}

We know that such a subconfiguration exists because of Lemma~\ref{lem:strongnocs} (with $\vv V=\vv W_2$ and $\vv W=\vv W_0$).

It is not hard to see that $\vv T_1$, $\vv T_2$ and $\vv W_0$ are linearly independent (\ie $\lin(\vv T_1\cup \vv T_2\cup \vv W_0)=\lin(\vv T_1)\oplus\lin(\vv T_2)\oplus\lin(\vv W_0)$). Indeed, choose respective bases $\vv B_0$, $\vv B_1$ and $\vv B_2$ for $\vv W_0$, $\vv T_1$ and~$\vv T_2$. For $j=\{1,2\}$, the elements of the projections $\tilde{\vv B}_j$ of $\vv B_j$ in $\vv V/\vv W_0$ are still linearly independent, because by construction $\lin(\vv T_j)\cap \lin(\vv W_0)=\veczero$. Moreover,
$\lin(\tilde{\vv B}_1)\cap\lin(\tilde{\vv B}_2)= \veczero$ because $\tilde{\vv B}_j\subset \vv U_j$ and $\lin(\vv U_1)\cap\lin(\vv U_2)=\veczero$. This already proves that $\vv B_0$, $\vv B_1$ and $\vv B_2$ are linearly independent, because if they had a linear dependence, this would create a dependence between $\tilde{\vv B}_1$ and $\tilde{\vv B}_2$ in $\vv V/\vv W_0$.

To be able to continue the proof of the two cases together, it is convenient to set $\vv T_3=\emptyset$ when $\{1,2\}\subseteq I$. With $\vv T_1$, $\vv T_2$ and $\vv T_3$ we will be able to find two hyperplanes $\vvh H'$ and $\vvh H$ that will lead to a contradiction. This is done below, after defining analogue subconfigurations $\vv T_1$, $\vv T_2$ and $\vv T_3$ for the case $\{2,3\}\subseteq I$. 

\vspace{.5cm}

 \item Suppose that $\{2,3\}\subseteq I$, and for $i\in\{2,3\}$, let Lemma~\ref{lem:strongnocs} choose as before a subconfiguration $\vv T_i\subset\vv W_i$ such that 
\begin{align*}
 \DD(\vv T_i)&\geq 1,\\ \DD(\vv T_i)&=\DD(\vv W_i)-\DD(\vv W_i/\vv T_i),\\ \rank(\vv T_i) &= 1,\quad\text{ and }\\
\lin(\vv W_0)\cap\lin(\vv T_i)&=\veczero. 
\end{align*}

Notice that $\vv T_2$ and $\vv T_3$ are linearly independent because they project to different $\vv U_i$'s in the contraction $\vv V/\vv W_0$.
 Moreover, choose $\vv T_1\subset\vv W_1$ such that
\begin{align}
 \DD(\vv T_1)&\geq 1,\notag\\
 \DD(\vv T_1)+\DD(\vv W_1/\vv T_1)&=\DD(\vv W_1),\notag\\
 \rank(\vv T_1)&=1,\notag\\
 \lin(\vv W_0\cup\vv T_2\cup\vv T_3)\cap\lin(\vv T_1)&=\veczero.\label{eq:condlinT1linT2T3}
\end{align}

To find this configuration, we use Lemma~\ref{lem:strongnocs} with $\vv V=\vv W_1$ and $\vv W=\vv W_1\cap\lin(\vv W_0\cup\vv T_2\cup\vv T_3)$. Observe that the rank of $\vv W_1\cap\lin(\vv W_0\cup\vv T_2\cup\vv T_3)$ is at most $2$, since $\vv W_0\cup\vv T_2\cup\vv T_3$ is a configuration of rank~$3$ that contains vectors that do not belong to $\lin(\vv W_1)$. 
Moreover, $\vv W_1\setminus \lin(\vv W_0\cup\vv T_2\cup\vv T_3)$ cannot be centrally symmetric because as we already discussed, if~$\vv W_1$ were centrally symmetric except for a configuration of rank~$\leq 2$, we would be in the case~\ref{it:rankW1=1} ($\rank(\vv U_1)=1$). Since $\vv W_1\setminus\lin(\vv W_0\cup\vv T_2\cup\vv T_3)$ is not centrally symmetric, Lemma~\ref{lem:strongnocs} states that we can find a suitable $\vv T_1$ avoiding $\lin(\vv W_1)\cap\lin(\vv W_0\cup\vv T_2\cup\vv T_3)$.

Again, we can see that $\vv T_1$, $\vv T_2$ and $\vv T_3$ are linearly independent (\ie $\rank(\vv T_1\cup\vv T_2\cup\vv T_3)=\rank(\vv T_1)+\rank(\vv T_2)+\rank(\vv T_3)$). Indeed, we have seen that $\rank(\vv T_2\cup\vv T_3)=\rank(\vv T_2)+\rank(\vv T_3)$, and $\vv T_1$ does not belong to $\lin(\vv T_2\cup\vv T_3)$ because of \eqref{eq:condlinT1linT2T3}. 

Similarly, we can prove that $\lin(\vv W_0)\cap\lin(\vv T_1\cup \vv T_2\cup \vv T_3)=\veczero$. One sees first that $\vv W_0\not\subset\lin(\vv T_2\cup \vv T_3)$, because these spaces are independent in the contraction; then the claim follows from \eqref{eq:condlinT1linT2T3}.

\end{compactenum}

\medskip

We rejoin now the proofs for cases 2.1.1 and 2.1.2. We define hyperplanes $\vvh H'$ and $\vvh H$ exactly like in Proposition~\ref{prop:2DD} and follow that proof to reach the same contradiction in~\eqref{eq:lasteq}. Indeed, observe that in both cases we have that
$\DD(\vv T_1)+\DD(\vv T_2)+\DD(\vv T_3)\geq 3$; with the convention $\vv T_3=\emptyset$ (with $\DD(\vv T_3)=0$) in the case $\{1,2\}\subseteq I$. Skipping some intermediate calculations that are equivalent to those in Proposition~\ref{prop:2DD}, one reaches

\begin{eqnarray*}
  \big|{\vvh H}^+\cap \vv V\big|
  &=&
  \big|{\vvh H}^+\cap \vv W_0\big|
   + \sum_{i=1}^{\l+1} \big|{\vvh H}^+\cap \vv W_i\setminus \vv W_0\big|+ 
  \\ &&\qquad
  {}+\sum_{i=\l+2}^{m} |{\vvh H}^+\cap \vv W_i\setminus \vv W_0|\notag\\
  &\geq& \codegG(\vv W_0)+\sum_{j=1}^3\DD(\vv T_j)+\sum_{i=1}^{m} \codegG(\vv U_i)\\
  &=&  \codegG(\vv V)+\sum_{j=1}^3\DD(\vv T_j)\\
  &=& \codegG(\vv V)+3= \codegG(\vv V)+\DD(\vv V)+1\\
 &=& r+\degG(\vv V)+1.
\end{eqnarray*}

This contradicts the definition of $\degG(\vv V)=\max_{\vvh H}|{\vvh H}^+\cap \vv V\big|-r$, and shows that $|I|\leq 1$, concluding the proof of the case~\ref{it:rankW1=2}.
\medskip

\item\label{it:rankW1=1} If $\rank(\vv U_1)=1$ and $I=\emptyset$ or $I=\{1\}$, we are trivially done, since $\vv V\setminus \vv W_1$ is then centrally symmetric by the arguments leading up to~\eqref{eq:DDWj-W0=0}, and we have decomposed $\vv V$ into a configuration $\vv W_1$ of rank~$\leq 2$ and a bunch of centrally symmetric vectors (therefore~$\vv V$ fulfills \eqref{it:1Wleq3}).

Otherwise, if there is some $i>1$ such that $i\in I$, we merge all the sets $\vv U_j$ with $\vv U_j\subset\lin(\vv U_1\cup \vv U_i)$ into a unique set $\vv U_1'$ of rank exactly~$2$. We will now see that $\vv W_1'$, the subconfiguration of $\vv V$ containing $\vv W_0$ and such that $\vv W_1'/\vv W_0=\vv U_1'$,  fulfills condition~\eqref{it2:lem:chooseinrank3} of Lemma~\ref{lem:chooseinrank3} (with $\vv V=\vv W_1'$ and $\vv W=\vv W_0$): 
\begin{itemize}
 \item $\vv W_1'$ cannot fulfill condition~\eqref{it0:lem:chooseinrank3} because there are at least three linearly independent subspaces of rank~$1$ that are not centrally symmetric. 
 Namely, $\vv W_0$ and some $\vv T_1\subset \vv W_1\setminus \vv W_0$ and $\vv T_i \subset \vv W_i\setminus \vv W_0$ where $i\in I$ and $\vv T_1$ and $\vv T_i$ are chosen using Lemma~\ref{lem:nocs}. Observe that $\vv T_1$, $\vv T_i$ and $\vv W_0$ are linearly independent because $\vv U_1$ and $\vv U_i$ are linearly independent in $\vv V/\vv W_0$.
The configurations $\vv T_1$~and~$\vv T_i$ avoid $\lin(\vv W_0)$ because $\lin(\vv W_0)\cap\vv V=\vv W_0$.
 \item Neither can $\vv W_1'$ fulfill condition~\eqref{it1:lem:chooseinrank3}. If it did, this would mean that $\vv W_1'\setminus \vv W_0$~is centrally symmetric except for a subconfiguration of rank~$2$. Since neither $\vv T_1$ nor $\vv T_i$ (from the previous point) are centrally symmetric, this means that this subconfiguration of rank~$2$ must be $\lin(\vv T_1\cup \vv T_i)\cap \vv W_1'$. 
 However, with Lemma~\ref{lem:strongnocs} we can find a configuration $\vv T_i'$ in $\vv W_i\setminus \vv W_0$ with $\DD(\vv T_i')>0$ (hence, not centrally symmetric) such that $\lin(\vv T_i')\cap \lin(\vv T_1\cup \vv T_j)=\veczero$. This would contradict the assumption that $\vv W_1'\setminus \vv W_0$ is centrally symmetric except for $\lin(\vv T_1\cup \vv T_j)\cap \vv W_1'$.
\end{itemize}

Since the proof of the previous case, $\rank(\vv U_1)=2$, only uses that $\vv W_1$~and~$\vv W_0$ fulfill condition~\eqref{it2:lem:chooseinrank3} of Lemma~\ref{lem:chooseinrank3}, 
the proof for the case $\rank(\vv U_1)=1$ and $\{i\}\subseteq I$ with $i>1$ is analogous, and proceeds as follows. 

First, we complete $\vv U_1'$ to a \codegreeG decomposition $\vv U_1',\dots, \vv U'_{m'}$ of~$\vv V/\vv W_0$ by aggregating collinear antipodal vectors. Then we define sets $\vv W_j'\supset\vv W_0$ such that $\vv W_j'/\vv W_0=\vv U_j'$, and the family of indices $I'=\{j \mid 1\le j\le m', \; \DD(\vv W_j'\setminus \vv W_0)>\DD(\vv U_j')\}$. 
Copying the proof for the case $\rank(\vv U_1)=2$ yields that $|I'|\leq 1$, which implies that either $\vv V$ is centrally symmetric except for $\vv W_1'$ (hence $\vv V$ fulfills~\eqref{it:1Wleq3}) or $\vv V$ is centrally symmetric except for some $\vv W_k'$ of rank~$2$ with $k>1$ and a subconfiguration of $\vv W_1'$ of rank~$\leq 2$. This implies that $\vv V$ fulfills \eqref{it:2Wleq2} or \eqref{it:1Wleq3}, depending on whether these configurations are skew or not, respectively.
\end{compactenum}
\end{compactenum}
This concludes the proof of Theorem \ref{thm:DD=2}.
\end{proof}

As a consequence, we have settled the conjecture for vector configurations in $\RR^r$ for $r\leq 4$.

\begin{corollary}~\label{cor:rank4}
Let $\vv V$ be a vector configuration in $\RR^r$ with $n=r+d+1$ elements and dual degree $\degG(V)=\dd$. If $r\leq 4$, then $\vv V$ admits a \codegreeG decomposition of length at least $r-\DD(\vv V)= d+1-2\dd$.
\end{corollary}
\begin{proof}
We can assume that $\vv V$ is irreducible by Observation~\ref{obs:pyramids}. 
The result is not trivial only if $r-\DD(\vv V) \geq 2$. Hence, $0\leq \DD(\vv V) \leq r-2\leq 2$. If $\DD(\vv V)=0$ we apply Corollary~\ref{cor:Lawrencedecomposition}, if $\DD(\vv V)=1$ we apply Corollary~\ref{cor:minimal+1} and if $\DD(\vv V)=2$ we apply Theorem~\ref{thm:DD=2}.
\end{proof}

It seems that we should be able to adapt the proofs of Proposition~\ref{prop:2DD} and Theorem~\ref{thm:DD=2} to obtain an inductive proof for Conjecture~\ref{conj:d+1-2dd} that would mimic the proof of Proposition~\ref{prop:deg1}. However, there are several spots where we still need more understanding to turn this to a general proof.
For example, these proofs rely on the fact that we understand the degrees of $\vv V/\vv W_i$, and that all the $\vv W_i$ fulfill $\lin(\vv W_i)\cap (\vv V/\vv W)=\vv W_i$. We do not know whether these properties hold in general. 

Moreover, in the general case we would have more constraints on the choice of the hyperplane~${\vvh H}$ in order to reach the contradiction. In the proof of Proposition~\ref{prop:deg1} we could use Proposition~\ref{prop:nosmallcircuits} and its corollaries to certify linear independence of certain subsets. 
Lemma~\ref{lem:strongnocs} sufficed for Proposition~\ref{prop:2DD}. Furthermore, for Theorem~\ref{thm:DD=2}, we had to use Lemma~\ref{lem:chooseinrank3}, whose proof uses the Sylvester-Gallai Theorem, which is a result much stronger than what we need (and also stronger than what we can prove in higher dimensions). Finding a good analogue of Lemma~\ref{lem:chooseinrank3} in higher dimensions might be the first ingredient for a definitive proof of Conjecture~\ref{conj:strongd+1-2dd}.

\section{One conjecture to prove them all}
\label{sec:consequences}

In this last section, we explore how some of the main results of previous chapters would follow from Conjecture~\ref{conj:strongd+1-2dd}. In particular, we derive proofs for Proposition~\ref{prop:dd=0}, Theorem \ref{thm:dd=1}, Theorem~\ref{thm:d+1-3dd} and Corollary~\ref{cor:easybound} that depend on Conjecture~\ref{conj:strongd+1-2dd} being true.

Conjecture~\ref{conj:strongd+1-2dd} states that every $d$-dimensional point configuration of degree~$\dd$ admits a codegree decomposition of length $\geq d+1-2\dd$. 
Therefore, to assume that Conjecture~\ref{conj:strongd+1-2dd} holds is equivalent to assume that every point configuration is combintorially equivalent to a configuration $\vv A$ that fulfills the following Assumption~\ref{ass:affdec}. It is stated in terms of affine decompositions to simplify the proofs, and this can be done without loss of generality because of Corollary~\ref{cor:definitionscombdecomposable}. 

\begin{assumption}\label{ass:affdec}
$\vv A$~is a $d$-dimensional configuration of $n$ points and degree $\degc(\vv A)=\dd$ 
such that there is a subset $\vv A_0\subset \vv A$ and a codegree preserving projection~$\pi:\vv A/\vv A_0\to \vv B=\vv B_1\join \vv B_2 \join \dots \join \vv B_m$ with $m\geq d+1-2\dd$ and $0\le \dim(\vv B_i)\leq 2\degc(\vv B_i)$.
\end{assumption}

The last assumption, $\dim(\vv B_i)\leq 2\degc(\vv B_i)$, is also a conclusion of the conjecture, since otherwise we could apply induction to decompose further. Moreover, it lets us understand the number of factors a little better and shows that the assumption $m\geq d+1-2\dd$ is redundant. Here and throughout, we use the convention $\dim(\emptyset)=-1$.

\begin{proposition}\label{prop:boundm}
 If $\vv A$ fulfills Assumption~\ref{ass:affdec}, then
\begin{eqnarray}
  m &\geq& 2\dim(\vv A)+1-2\degc(\vv A)-\dim(\vv B)\notag \\ 
    &=& \phantom{2}\dim(\vv A) + 1 - 2\degc(\vv A) + \dim(\ker\pi)+\dim(\vv A_0)+1\notag\\
    &\geq&\phantom{2}\dim(\vv A) + 1 - 2\degc(\vv A) + \dim(\ker\pi)\label{eq:firstboundm}\\
    &\geq& \phantom{2}\dim(\vv A) + 1 - 2\degc(\vv A). \notag 
\end{eqnarray}
\end{proposition}

\begin{proof}
Throughout, we abbreviate $\dim(\vv A)=d$, $\degc(\vv A)=\dd$, $\dim(\vv A_0)=s$, $\dim(\vv B_i) = d_i$, $\degc( \vv B_i) = \dd_i$, and $e=\dim(\ker\pi)$. With this notation, we derived in Section~\ref{sec:operations} that $\dim(\vv B)+1=\sum_{i=1}^m(d_i+1)$ and $\dim(\vv A/\vv A_0)=d-s-1$. Using these identities, the relation $\dim(\vv A/\vv A_0)=\dim\vv B + \dim(\ker\pi)$ can be phrased as
\begin{equation}\label{eq:d+1}
d+1 = e+s+1+\sum_{i=1}^m{(d_i+1)}.
\end{equation}
The assumption that $\pi$ preserves codegrees translates into
\begin{equation}\label{eq:d+1-dd}
  d+1-\dd=\sum_{i=1}^m{(d_i+1-\dd_i)}.
\end{equation}
Finally, taking the linear combination $2\eqref{eq:d+1-dd}-\eqref{eq:d+1}$ yields
\[
 d+1-2\dd \ = \ \sum_{i=1}^m(d_i+1-2\dd_i)-e-s-1 \ \leq\ m-e-s-1,
\]
where in the second step we used the assumption that $d_i\leq 2\dd_i$. Therefore,
\begin{align*}
  m&\geq (d+1-2\dd)+e+s+1,
\end{align*}
from where the expressions in the statement can be easily recovered.
\end{proof}

Next, we show how our previous results prove some special cases of Conjecture~\ref{conj:strongd+1-2dd}, and how it in turn implies several of our previous results. 

\begin{proposition}[Conjecture~\ref{conj:strongd+1-2dd} holds for $\dd=0$]\label{prop:conjdd=0}
If $\degG(\vv A)=0$, then~$\vv A$ fulfills Assumption~\ref{ass:affdec}.
\end{proposition}

\begin{proof}
  Indeed, by Proposition~\ref{prop:dd=0}, $\vv A$ is the vertex set of a simplex~$\simp{d}$, possibly with repetitions. For $1\le i\le d+1$, let $\vv B_i\subset\vv A$ be the set of all copies of the $i$-th vertex of~$\simp{d}$. To see that the $\vv B_i$ form a codegree decomposition of~$\vv A$, note that $\dim\vv B_i=0$, so that $\codegc(\vv B_i)=1$; there are $m=d+1$ such sets; and $\codegc(\vv A)=d+1=\sum_{i=1}^m\codegc(\vv B_i)$.
\end{proof}

The converse also holds:

\begin{proposition}[Conjecture~\ref{conj:strongd+1-2dd} implies Proposition~\ref{prop:dd=0}]
If $\vv A$ fulfills Assumption~\ref{ass:affdec} and $\dd=0$, then $\vv A$ is the set of vertices of a $d$-simplex (possibly with repetitions).
\end{proposition}
\begin{proof}
If $m\geq d+1$, then $\vv B_1\join \vv B_2 \join \dots \join \vv B_m$ must be $d$-dimensional and each $\vv B_i$ $0$-dimensional. Thus $\vv A=\vv B_1\join \vv B_2 \join \dots \join \vv B_m$, which is the set of vertices of a $d$-simplex.
\end{proof}

\begin{proposition}[Conjecture~\ref{conj:strongd+1-2dd} holds for $\dd=1$]\label{prop:conjdd=1}
If $\degG(\vv A)=1$, then~$\vv A$ is combinatorially equivalent to a configuration that fulfills Assumption~\ref{ass:affdec}.
\end{proposition}

\begin{proof}
 This is a consequence of Theorem~\ref{thm:dd=1}. Indeed, if $\vv A$~is a $k$-fold pyramid over a polygon, it is a join of a $2$-dimensional configuration of codegree~$2$ with $(d-1)$ point configurations of dimension~$0$ and codegree~$1$. 

If on the other hand ${\vv A}$ is a weak Cayley configuration of length~$d$, then a contraction of~$\vv A$ projects onto the vertex set of a simplex, a join of $d$~point configurations of dimension~$0$ and codegree~$1$. 
\end{proof}

The reciprocal is also easy.

\begin{proposition}[Conjecture~\ref{conj:strongd+1-2dd} implies Theorem~\ref{thm:dd=1}]
If $\vv A$ fulfills Assumption~\ref{ass:affdec} and $\dd=1$, then 
\begin{enumerate}
 \item ${\vv A}$ is a $k$-fold pyramid over a $2$-dimensional point configuration without interior points; or
 \item ${\vv A}$ is a weak Cayley configuration of length $d$.
\end{enumerate}
\end{proposition}
\begin{proof}
As before, let $d_i$ be the dimension of $\vv B_i$, $d_i+1+r_i$ its number of elements and $\dd_i$ its degree. We can assume $d_i\leq 2\dd_i$. Let $s=\dim(\vv A_0)$ and $e=\dim(\vv A/\vv A_0)-\dim(\vv B)$. 

Observe that $m\ge d-1$ by Proposition~\ref{prop:boundm}. Moreover, $m\le d$. Indeed, by Lemma~\ref{lem:decompiscayley}, if $m\geq d+1$, then $\vv A$ is a weak Cayley configuration of length $m\geq d+1$. Then, by Proposition~\ref{prop:deg-cayley}, $\degc(\vv A)\leq d+1-m\leq 0$, which would contradict the hypothesis that the degree of $\vv A$ is $1$. 

Hence, we have two possible values for $m$:
\begin{itemize}

\item If $m=d$, then $\vv A$ is a weak Cayley configuration of length $d$ by Lemma~\ref{lem:decompiscayley}.

\item Finally, suppose $m=d-1$. Since $d=d+1-\dd=\sum_{i=1}^m (d_i+1-\dd_i)$ by hypothesis, and $d_i+1-\dd_i=\codegG(\vv B_i)\geq 1$ for all $i$ because the codegree of a point configuration is always at least $1$, we can assume that
\begin{eqnarray}
 d_1+1-\dd_1 &=& 2 \qquad\text{and}\label{eq:case1} \\
 d_i+1-\dd_i &=& 1 \qquad\text{for } i\ge2.\notag
\end{eqnarray}
Now $d_1=\dd_1+1$ by~\eqref{eq:case1}, and the assumption $d_1\leq 2\dd_1$ implies that 
\[
  2 = d_1+1-\dd_1 \leq \dd_1+1 = d_1.
\]
Combining this with the identity
\[
  d+1=e+s+1+\sum_{i=1}^m{(d_i+1)},
\]
from \eqref{eq:d+1} implies that $e=0$, $s=-1$, $d_1=2$ and $d_i=0$ for $i>1$, because by definition $e\geq0$, $s\geq-1$ and $d_i\geq0$. Thus, $\vv A$ is, in effect, a $k$-fold pyramid over a $2$-dimensional point configuration of codegree~$2$, which cannot have interior points. 
\end{itemize}
This proves our claim.
\end{proof}

Corollary~\ref{cor:easybound} also follows from  Conjecture~\ref{conj:strongd+1-2dd}:

\begin{proposition}[Conjecture~\ref{conj:strongd+1-2dd} implies Corollary~\ref{cor:easybound}]
If $\vv A$ fulfills Assumption~\ref{ass:affdec}, $r:=n-d-1$ and 
\[d \ge r+2 \dd,\]
then $\vv A$ is a pyramid.
\end{proposition}
\begin{proof}
Observe first that if a contraction of $\vv A$ is a pyramid, then so is~$\vv A$; if a projection of $\vv A/\vv A_0$ is a pyramid, then so is~$\vv A/\vv A_0$; and if some $\vv B_i$ is a simplex, then $\vv B_1\join \dots \join \vv B_m$ is a pyramid.

Set $d_i=\dim(\vv B_i)$ and let $|\vv B_i| = d_i+1+r_i$. We prove that some $\vv B_i$ must be a simplex by showing that some $r_i=0$. Set also $s=\dim(\vv A_0)\leq |\vv A_0|-1$ and $e=\dim(\vv A/\vv A_0)-\dim(\vv B)$.

By counting the elements in $\vv A$ and $\vv B$, we get
\begin{align}
r+d+1=|\vv A|&=|\vv A_0|+\sum_{i=1}^m|\vv B_i|\notag\\
&\geq s+1+\sum_{i=1}^m(r_i+d_i+1).\label{eq:r+d+1}
\end{align}
The linear combination \eqref{eq:r+d+1}-\eqref{eq:d+1} yields
\begin{equation*}\label{eq:r}
r+e\ \geq\ \sum_{i=1}^m{r_i}.
\end{equation*}
Hence, if we assume that each $r_i\geq 1$, this implies
\begin{equation*}\label{eq:rbound}
  r\ \geq \ -e+\sum_{i=1}^m{r_i}
  \ \ge\ 
  m-e
  \ \stackrel{\eqref{eq:firstboundm}}{\geq}\ 
  d+1-2\dd,
\end{equation*}
which in turn implies $d<r+2\dd$. Thus, $d\geq r+2\dd$ forces some $r_i=0$, and we conclude that $\vv A$ is a pyramid.
\end{proof}

Finally, we have already commented that Conjecture~\ref{conj:strongd+1-2dd} is stronger than Conjecture~\ref{conj:d+1-2dd}, which is stronger than Theorem~\ref{thm:d+1-3dd}.

\begin{proposition}[Conjecture~\ref{conj:strongd+1-2dd} implies Conjecture~\ref{conj:d+1-2dd}]
 If $\vv A$ fulfills Assumption~\ref{ass:affdec} and $d > 2\dd$, then $\vv A$ is a weak Cayley configuration of length at least $d+1-2\dd$.
\end{proposition}
\begin{proof}
This is a direct consequence of Lemma~\ref{lem:decompiscayley}.
\end{proof}

\newpage \thispagestyle{empty} \phantom{ a }\newpage
\backmatter
\renewcommand{\namepart}{Bibliography}
\setcounter{compteurlof}{0}
\bibliographystyle{alpha}
\bibliography{Thesis}

\renewcommand{\namepart}{Index}
\printindex
\end{document}